\documentclass[10pt]{article}
\usepackage{latexsym}
\usepackage[all]{xy}
\usepackage{amsmath, amsthm, amssymb, amsfonts}
\usepackage{enumerate}
\usepackage[latin1]{inputenc}
\usepackage{multirow}
\usepackage{hyperref}
\usepackage{abstract, indentfirst}
\usepackage{amscd, stmaryrd, setspace,color}
\usepackage{graphicx}
\usepackage{turnstile}
\usepackage{wasysym}

\newcommand{\ie}{{i.e. }}
\newcommand{\et}{\qquad\textrm{and}\qquad}
\newcommand{\bs}{\backslash}

\newcommand{\tto}{\longrightarrow}
\newcommand{\mto}{\longmapsto}
\newcommand{\lrhup}{-\!\!\!\rhup}
\newcommand{\Hom}{{\textrm{Hom}}}
\newcommand{\HOM}{\textrm{\sc Hom}}
\newcommand{\END}{\textrm{\sc End}}
\newcommand{\un}{\mathbf{1}}

\newcommand{\FF}{\mathbb{F}}

\newcommand{\ZZ}{\mathbb{Z}}

\newcommand{\cM}{\mathcal{M}}

\newcommand{\cT}{\mathcal{T}}
\newcommand{\cU}{\mathcal{U}}
\newcommand{\cV}{\mathcal{V}}

\newcommand{\bA}{\mathbf{A}}
\newcommand{\bB}{\mathbf{B}}
\newcommand{\bC}{\mathbf{C}}
\newcommand{\bD}{\mathbf{D}}
\newcommand{\bE}{\mathbf{E}}

\newcommand{\bV}{\mathbf{V}}

\newcommand{\bCat}{\mathbf{Cat}}

\renewcommand{\d}{{\sf d}}

\newcommand{\colim}{\mathrm{co}\!\lim}
\newcommand{\ccolim}{\mathrm{colim}}

\newcommand{\Mon}{\mathsf{Mon}}
\newcommand{\coMon}{\mathsf{coMon}}

\newcommand{\rhup}{\rightharpoonup} 

\newcommand{\Set}{\mathsf{Set}}
\newcommand{\gSet}{\mathsf{gSet}}

\newcommand{\Vect}{\mathsf{Vect}}

\newcommand{\gVect}{\mathsf{gVect}}

\newcommand{\dgVectfin}{\mathsf{dgVect^{fin}}}
\newcommand{\dgVectgrfin}{\mathsf{dgVect^{gr.fin}}}
\newcommand{\dgVecta}{\mathsf{dgVect^{a}}}
\newcommand{\dgVectb}{\mathsf{dgVect^{b}}}
\newcommand{\dgVectgrfina}{\mathsf{dgVect^{gr.fin,a}}}
\newcommand{\dgVectgrfinb}{\mathsf{dgVect^{gr.fin,b}}}

\newcommand{\dgVect}{\mathsf{dgVect}}
\newcommand{\ddgVect}{\mathsf{d^2gVect}}

\newcommand{\Mod}{\mathsf{Mod}}
\newcommand{\Bimod}{\mathsf{Bimod}}
\newcommand{\Comod}{\mathsf{Comod}}
\newcommand{\Bicomod}{\mathsf{Bicomod}}

\newcommand{\gAlg}{\mathsf{gAlg}}
\newcommand{\dgAlg}{\mathsf{dgAlg}}
\newcommand{\dgAlgfin}{\mathsf{dgAlg^{fin}}}
\newcommand{\dgAlgfincirc}{\mathsf{dgAlg^{fin}_\circ}}

\newcommand{\dgAlgfnilcirc}{\mathsf{dgAlg^{fin,nil}_\circ}}
\newcommand{\dgAlgcom}{\mathsf{dgAlg^{com}}}

\newcommand{\EAlg}{\mathsf{gAlg^\$}}

\newcommand{\gfAlg}{\mathsf{gfAlg}}
\newcommand{\fAlg}{\mathsf{fAlg}}

\newcommand{\nfAlg}{\mathsf{fAlg^{nil}}}

\newcommand{\Coalg}{\mathsf{Coalg}}
\newcommand{\gCoalg}{\mathsf{gCoalg}}
\newcommand{\dgCoalg}{\mathsf{dgCoalg}}
\newcommand{\dgCoalgfin}{\mathsf{dgCoalg^{fin}}}

\newcommand{\dgCoalgnil}{\mathsf{dgCoalg^{conil}}}
\newcommand{\dgCoalgnilcirc}{\mathsf{dgCoalg^{conil}_\circ}}
\newcommand{\dgCoalgnilbullet}{\mathsf{dgCoalg^{conil}_\bullet}}

\newcommand{\dgCoalgcom}{\mathsf{dgCoalg^{com}}}

\newcommand{\gfCoalg}{\mathsf{gfCoalg}}
\newcommand{\fCoalg}{\mathsf{fCoalg}}

\newcommand{\dgBialg}{\mathsf{dgBialg}}

\newcommand{\Lie}{\mathsf{Lie}}

\newcommand{\dgLie}{\mathsf{dgLie}}

\newcommand{\OpMod}{\mathbf{OpMod}}
\newcommand{\MMod}{\mathbf{Mod}}

\newcommand{\Der}{{\sf Der}}
\newcommand{\Coder}{{\sf Coder}}
\newcommand{\Prim}{{\sf Prim}}

\newcommand{\mc}{\textrm{\sc mc}}
\newcommand{\bOmega}{{\bf \Omega}}
\newcommand{\Beta}{\mathrm{B}}
\newcommand{\Betaext}{\mathrm{B^{ext}}}

\definecolor{mat}{rgb}{0.57,0.75,1}

\newcommand{\oddbinom}[2]{\left[\genfrac{}{}{0pt}{}{#1}{#2}\right]}

\newcommand{\Dfrac}[2]{\frac{#1}{#2}}
\renewcommand{\dfrac}[2]{\frac{#1}{#2}}

\newtheorem{thm}{Theorem}[subsection]
\newtheorem{prop}[thm]{Proposition}
\newtheorem{cor}[thm]{Corollary}
\newtheorem{lemma}[thm]{Lemma}

\theoremstyle{definition}
\newtheorem{defi}[thm]{Definition}
\newtheorem{ex}[thm]{Example}

\newtheorem{rem}[thm]{Remark}
\newtheorem{nota}[thm]{Notation}

\begin{document}

\title{{\bf Sweedler theory of (co)algebras \\ and the bar-cobar constructions}
\vspace{1cm}}

\author{
\sc M. Anel \footnote{
\href{http://www.math.ethz.ch/}{\sc eth Zürich}, 
\href{mailto:matthieu.anel@math.ethz.ch}{matthieu.anel@math.ethz.ch} 
} \hspace{1cm} 
A. Joyal \footnote{
\href{http://www.cirget.uqam.ca/}{\sc cirget, UQ\`AM}, 
\href{mailto:joyal.andre@uqam.ca}{joyal.andre@uqam.ca} 
}\\ 
}

\maketitle

\vspace{1cm}

\begin{abstract}

\noindent We prove that the category of dg-coalgebras $(\dgCoalg,\otimes, \HOM)$ is symmetric monoidal closed and that the category of dg-algebras $(\dgAlg, \{-,-\},\rhd,[-,-],\otimes)$ is enriched, tensored, cotensored and strongly monoidal over $\sf dgCoalg$.
For $A$ and $B$ two dg-algebras, the enriched hom $\{A,B\}$ is \cite{Sw}'s universal measuring coalgebra, the cotensor of an algebra $A$ by a coalgebra $C$ is the convolution algebra $[C,A]$ and the tensor of $A$ by $C$ is a new operation $C\rhd A$ called the Sweedler product. We call the resulting structure Sweedler theory. Sweedler theory exists in many contexts, we detail also the case of (co)augmented (co)algebras.

\medskip
\noindent Sweedler operations can be used to produce various adjunctions between categories of dg-algebras and dg-coalgebras:
\begin{center}
\begin{tabular}{c}
$\xymatrix{C\rhd - :\dgAlg \ar@<.6ex>[r]&\ar@<.6ex>[l] \dgAlg :[C,-]}$\\
$\xymatrix{[-,A] :\dgCoalg \ar@<.6ex>[r]&\ar@<.6ex>[l] \dgAlg^{op} :\{-,A\}}$\\
$\xymatrix{-\rhd A :\dgCoalg \ar@<.6ex>[r]&\ar@<.6ex>[l] \dgAlg :\{A,-\}}$
\end{tabular}
\end{center}
We apply this formalism to reconstruct several known adjunctions, particularly the bar-cobar adjunction.

\medskip
\noindent The present paper is partly expository and purely algebraic, we do not investigate the homotopical aspects of the theory.
In a second paper, we will extend our theory to operads and cooperads.

\end{abstract}

\newpage
\setcounter{tocdepth}{3}
\tableofcontents

\newpage

\section*{Introduction}
\addcontentsline{toc}{section}{\bfseries Introduction}

This paper is partly expository and self-contained.
The initial motivation was our desire to understand conceptually Koszul duality for (co)algebras
and (co)operads and for this,  it was natural to first study the bar-cobar duality. It turns out that the bar-cobar  duality can be understood in terms of very general constructions on dg-algebras and dg-coalgebra and that our theory has applications beyond  the duality. A central notion in this paper 
is that of measuring  introduced  by E. M. Sweedler in his book \cite{Sw}, 
and for this reason, we call our theory the {\em Sweedler theory} of algebras and coalgebras.
Many of the results presented here can be further extended
to algebras and coalgebras in a general locally presentable symmetric monoidal closed category,  
but we have limited the exposition to the case of dg-algebras and coalgebras over a field.

\medskip

The Sweedler theory of operads and cooperads was sketched in \cite{AJams} and it will be developed  in a second paper.

\paragraph{Contents}

Chapter \ref{basicalgebra} is a recollection on associative algebras, coalgebras, bialgebras and Lie algebras. It may be skipped at first reading. Chapters \ref{coalgebras} and \ref{algebras} is the core of the paper.
We show that the category of algebras is enriched and bicomplete over the category of coalgebras,
and in particular that it is tensored and cotensored over coalgebras.
Chapter \ref{applications} contains applications; we reconstruct many known adjunctions between algebras and coalgebras, in particular the bar-cobar adjunction.
Chapter \ref{Sweedlercontexts} is a remake (without proofs) of chapters \ref{coalgebras} and \ref{algebras} for non-(co)unital (co)algebras and pointed coalgebras.
Appendix \ref{barcobar} is a recollection of the classical bar and cobar construction following \cite{LV}. 
Appendix \ref{AppendixCattheory} is a collection of categorical results 
 used in the paper: Gray trialities, monoidal structures and functors, enriched category theory...

\subsection*{Main results}

\paragraph{The extended bar-cobar adjunction}

Recall from \cite{EB, Proute} that if $A$ is a dg-algebra and $C$ is a dg-coalgebra (over a field $\mathbb{F}$),
then a twisting cochain $\alpha:C\to A$ is a linear map of degree $-1$ satisfying the Maurer-Cartan equation $d\alpha +\alpha^2=0$ in the convolution dg-algebra $[C,A]$. 
If $A$ is pointed (= coaugmented) and $C$ is pointed  (=augmented),
the cochain $\alpha$  is said to be pointed (=admissible) if $\epsilon \alpha=0$ and $\alpha e=0$ where 
$\epsilon:A\to \FF$ is the augmentation of $A$ and $e:\FF\to C$ is the coaugmentation of $C$.
The cobar construction \cite{Ad1} takes a pointed dg-coalgebra $C$ to a pointed dg-algebra $\Omega C$ and we have a natural
isomorphism
$$
 Tw_\bullet(C,A) \simeq \dgAlg_\bullet(\Omega C,A)$$
where $\dgAlg_\bullet$ is the category of pointed dg-algebras.
One application of this work is to show that the functor $Tw_\bullet(-,A)$
is representable by a pointed dg-coalgebra $\Betaext A$ 
for any pointed dg-algebra $A$.  From the natural bijections 
$$
\dgAlg_\bullet(\Omega C,A)\simeq Tw_\bullet(C,A) \simeq \dgCoalg_\bullet(C,\Betaext  A).
$$
we obtain an adjunction
$$\xymatrix{
\Omega:  \dgCoalg_\bullet \ar@<.6ex>[r]&\ar@<.6ex>[l] \dgAlg_\bullet:\Betaext
}$$
where $\dgCoalg_{\bullet}$ is the category of pointed dg-coalgebras.
We call $\Betaext A$ the {\it extended} bar construction of $A$.
Recall that the classical  bar construction \cite{EM} takes a pointed dg-algebra $A$  to a pointed
{\it conilpotent} dg-coalgebra denoted $\Beta A$.
The category $\dgCoalgnilbullet$ of pointed conilpotent dg-coalgebra
 is a full coreflexive subcategory of the category $\dgCoalg_{\bullet}$
 and the coreflexion  functor takes a pointed dg-coalgebra $C$ to its {\it conilpotent radical} 
$R^c(C)\subseteq C$. We have $\Beta A=R^c(\Betaext  A)$ for any dg-algebra $A$.
This yields the classical bar-cobar adjunction (theorems \ref{classicbarcobar} and \ref{represtw})
$$\xymatrix{
\Omega: \dgCoalgnilbullet \ar@<.6ex>[r]&\ar@<.6ex>[l] \dgAlg_\bullet:\Beta.
}$$

In theorem \ref{represtw}, we construct 
 the {\em unpointed bar-cobar adjunction}
 $$\xymatrix{
\bOmega: \dgCoalg \ar@<.6ex>[r]&\ar@<.6ex>[l] \dgAlg:\bB,
}$$
between the categories of unpointed dg-coalgebras and dg-algebras.
The adjunction is obtained by composing two natural bijections
$$
\dgAlg(\bOmega C,A)\simeq Tw(C,A) \simeq \dgCoalg(C,\bB A),
$$
where $Tw(C,A)$ is the set of all twisting cochains between an unpointed dg-coalgebra
and an unpointed dg-algebra.

\paragraph{Sweedler theory}

Sweedler theory is the name we give to the following set of functors relating the categories $\dgCoalg$ of coalgebras and $\dgAlg$ algebras:
\begin{center}
\begin{tabular}{lrl}
\rule[-2ex]{0pt}{4ex} the tensor product of coalgebras & $\otimes$&$\!\!\!\!: \dgCoalg\times \dgCoalg \to \dgCoalg$, \\
\rule[-2ex]{0pt}{4ex} the coalgebra internal hom & $\HOM$&$\!\!\!\!:\dgCoalg^{op}\times \dgCoalg \to \dgCoalg$, \\
\rule[-2ex]{0pt}{4ex} Sweedler hom & $\{-,-\}$&$\!\!\!\!:\dgAlg^{op}\times \dgAlg \to \dgCoalg$, \\
\rule[-2ex]{0pt}{4ex} Sweedler product & $\rhd$&$\!\!\!\!:\dgCoalg\times \dgAlg \to \dgAlg$, \\
\rule[-2ex]{0pt}{4ex} the convolution product & $[-,-]$&$\!\!\!\!:\dgCoalg^{op}\times \dgAlg \to \dgAlg$, \\
\rule[-2ex]{0pt}{4ex} and the tensor product of algebras & $\otimes$&$\!\!\!\!:\dgAlg\times \dgAlg \to \dgAlg$.
\end{tabular}
\end{center}
The two tensor products $\otimes$ and the convolution functor $[-,-]$ are classical, the functor $\{-,-\}$ is Sweedler's universal measuring coalgebra from \cite{Sw} and \cite{Vasilakopoulou}, $\HOM$ is constructed in \cite{Porst}, but the functor $\rhd$ seems new.

\medskip

These functors satisfy the following structure, which is our main result.
\begin{thm}[theorems \ref{homcoalg} and \ref{enrichmentalgcoal}]\ 
\begin{enumerate}
\item The category $(\dgCoalg,\otimes, \HOM)$ is symmetric monoidal closed.
\item The category $(\dgAlg,\{-,-\},\rhd,[-,-],\otimes )$ is enriched, tensored, cotensored and symmetric monoidal over $\dgCoalg$.
\end{enumerate}
\end{thm}

\medskip

The first result says that, for any two dg-coalgebras $C$ and $D$, there exists a coalgebra $\HOM(C,D)$ such that, for any dg-coalgebra $E$, we have natural bijections between 
\begin{center}
\begin{tabular}{lc}
\rule[-2ex]{0pt}{4ex} coalgebra maps & $E\otimes C\to D$, \\
\rule[-2ex]{0pt}{4ex} coalgebra maps & $E\to \HOM(C,D)$, \\
\rule[-2ex]{0pt}{4ex} and coalgebra maps & $C\to \HOM(E,D)$.
\end{tabular}
\end{center}
This result was proven by H.-E. Porst \cite{Porst} in the general setting of locally presentable categories. 
We proceed the same way, but we prove by hand the $\omega$-presentability of $\dgCoalg$ (theorem \ref{finprescoalg}).
Then, the existence of the cofree coalgebra functor $T^\vee$  (theorem \ref{cofreecoalg}),
the comonadicity of $\dgCoalg$ over $\d\Vect$ (theorem \ref{comonadic1})
and the $\HOM$ functor (theorem \ref{homcoalg})
are all consequences of the special adjoint theorem.

\bigskip

The result about algebras seems new stated in this form, although the enrichment of $\dgAlg$ over $\dgCoalg$ has been previously proved in \cite{Vasilakopoulou}. This enrichment is not obvious and requires some definitions that we shall give now.
Recall that if $C$ is a dg-coalgebra and $A$ a dg-algebra, the complex of graded morphisms $[C,A]$ has the structure of a dg-algebra called the {\em convolution algebra}.
If $C$ is a dg-coalgebra and $A$ and $B$ are two dg-algebras, a map $C\otimes A\to B$ is called a {\em measuring} if the corresponding map $A\to [C,B]$ is a map of algebras \cite[ch. VII]{Sw}.
If ${\cal M}(C,A;B)$ is the set of measurings $C\otimes A\to B$, these sets define a functor
$$\xymatrix{
{\cal M} : \dgCoalg^{op}\times\dgAlg^{op}\times \dgAlg\ar[r]& \Set}.
$$
By definition of measurings, this functor is representable in the second variable by the convolution algebra. This means that, when $C$ and $B$ are fixed, there is a natural bijection between measurings $C\otimes A\to B$ and dg-algebra maps $A\to [C,B]$. In his book \cite[ch. VII]{Sw}, Sweedler proves that this functor is also representable in the first variable, \ie that, when $A$ and $B$ are fixed, there exists a coalgebra $\{A,B\}$ (denoted $M(A,B)$ in \cite{Sw}) and a natural bijection between measurings $C\otimes A\to B$ and dg-coalgebra maps $C\to \{A,B\}$ (theorem \ref{Coalghomalg1}). 
One of our main results (theorem \ref{Sproduct}) is that this functor is also representable in the third variable, \ie that, when $C$ and $A$ are fixed, there exists a dg-algebra $C\rhd A$ and a natural bijection between measurings $C\otimes A\to B$ and dg-algebras maps $C\rhd A\to B$. 

Altogether, we have the fundamental bijections between
\begin{center}
\hfill
\begin{tabular}{lc}
\rule[-2ex]{0pt}{4ex} measurings & $C\otimes A\to B$, \\
\rule[-2ex]{0pt}{4ex} dg-coalgebra maps & $C\to \{A,B\}$, \\
\rule[-2ex]{0pt}{4ex} dg-algebra maps & $C\rhd A\to B$, \\
\rule[-2ex]{0pt}{4ex} and dg-algebra maps & $A\to [C,B]$.
\end{tabular}
\hfill (ST)
\end{center}
We shall call $\{A,B\}$ the {\em Sweedler hom} of $A$ and $B$ and $C\rhd A$ the {\em Sweedler product} of $A$ by $C$.
We shall also refer to $[C,A]$ as the {\em convolution product} of $A$ by $C$.

\medskip

We prove in theorem \ref{enrichmentalgcoal} that the coalgebras $\{A,B\}$ are equipped with a composition law
$$\xymatrix{
{\bf c}:\{B,C\}\otimes \{A,B\}\ar[r]& \{A,C\}
}$$
which turns them into an enrichment of $\dgAlg$ over $\dgCoalg$.
Moreover the adjunctions (ST) says that this enrichment is tensored:
the Sweedler product 
$$\xymatrix{
\rhd : \dgCoalg\times \dgAlg\ar[r]& \dgAlg
}$$
defines is a left action of $(\dgCoalg,\otimes)$ on $\dgAlg$, and cotensored: the convolution product
$$\xymatrix{
[-,-] : \dgCoalg^{op}\times \dgAlg\ar[r]& \dgAlg
}$$
defines is a right action of $(\dgCoalg,\otimes)$ on $\dgAlg$.
The category $\dgAlg$ has all limits and colimits, then, the existence of tensor and cotensor products says that it has also all weighted (co)limits over $\dgCoalg$.
(We refer the reader to Appendix \ref{AppendixCattheory} for details on enriched category theory.)

\bigskip

Finally, we prove that the tensor product of algebras is compatible with the enrichment.
\begin{thm}[theorem \ref{tensorenrichmentof}]
The tensor product of algebras 
$$\xymatrix{
\otimes : \dgAlg\times \dgAlg\ar[r]& \dgAlg
}$$
is a strong symmetric monoidal structure for the enrichment of $\dgAlg$ and $\dgAlg\times \dgAlg$ over $\dgCoalg$.
Moreover, the functors $[-,-]$ and $\{-,-\}$ are strong lax functors and $\rhd$ is a strong colax functor.
\end{thm}

This last fact implies the existence of various functors between enriched categories of (co)commutative (co)algebras and bialgebras (see section \ref{monoidalalg}). This simple fact also implies the well known formulas for iterations of the bar and cobar constructions (see section \ref{iteratedbarcobar}).

\paragraph{From Sweedler theory to the bar-cobar adjunctions}

The bijections (ST) say in particular that, for any dg-coalgebra $C$ and any dg-algebra $A$, we have adjunctions
\begin{center}
\hfill
\begin{tabular}{c}
\rule[-2ex]{0pt}{4ex} $\xymatrix{C\rhd - :\dgAlg \ar@<.6ex>[r]&\ar@<.6ex>[l] \dgAlg :[C,-]}$,\\
\rule[-2ex]{0pt}{4ex} $\xymatrix{[-,A] :\dgCoalg \ar@<.6ex>[r]&\ar@<.6ex>[l] \dgAlg^{op} :\{-,A\}}$,\\
\rule[-2ex]{0pt}{4ex} $\xymatrix{-\rhd A :\dgCoalg \ar@<.6ex>[r]&\ar@<.6ex>[l] \dgAlg :\{A,-\}}$.
\end{tabular}
\hfill
\begin{tabular}{c}
\rule[-2ex]{0pt}{4ex} (I)\\
\rule[-2ex]{0pt}{4ex} (II)\\
\rule[-2ex]{0pt}{4ex} (III)
\end{tabular}
\end{center}
These adjunctions encompass several known constructions on algebras and coalgebras. They are all detailed in chapter \ref{applications}
(\footnote{We have abstracted the structure of such families of adjunctions and their relation a measuring-like functor under the name {\em Gray triality}, this theory is develop in appendix \ref{triality}.}).

\medskip

Examples of type I includes:
products and coproduct of algebras (section \ref{productcoproduct});
non-commutative analogs of Weil restriction of scalars (section  \ref{Weilrestriction}) and
de Rham complexes (section  \ref{derhamexampleadj});
the contruction of matrix algebras (section \ref{matrixexampleadj}) and jet algebras (sections \ref{jetexample1adj} and \ref{jetexample2adj}).

Examples of the other types are rarer. To our knowledge, the only example of type II is Sweedler's duality (developped in section \ref{sweedlerduality}) and the main example of type III is the bar-cobar adjunction that we shall explain now.

\bigskip

Let $\mc$ be the free dg-algebra on one generator $u$ of degree $-1$ and with differential defined by $du+u^2=0$. $\mc$ is called the {\em Maurer-Cartan algebra}. For any dg-algebra $A$, there exists a natural bijection between
\begin{center}
\begin{tabular}{l}
\rule[-2ex]{0pt}{4ex} Maurer-Cartan elements of $A$ \\
\rule[-2ex]{0pt}{4ex} and maps of dg-algebras $\mc\to A$.
\end{tabular}
\end{center}
The bijections (ST) give, for any dg-coalgebra $C$ and any dg-algebra $A$ bijections
\begin{center}
\hfill
$\dgAlg(C\rhd \mc, A) = \dgAlg(\mc, [C,A]) = \dgCoalg(C, \{\mc,A\}).$
\hfill
(BCB)
\end{center}
Because $\dgAlg(\mc, [C,A]) = Tw(C,A)$, these bijections are close to the usual characterization of the bar-cobar adjunction, only the pointing condition on the twisting cochains is missing. This adjunction is the unpointed bar-cobar adjunction $\bOmega\dashv\bB$ mentionned above.

\medskip

To obtain the classical bar-cobar adjunction, we need to adapt Sweedler operations to the setting of pointed (co)algebras, which is done in section \ref{pointedsweedlertheory}. We have found convenient to have results in both the equivalent languages of pointed (co)algebras and non-(co)unital (co)algebras; depending on the context, it is nicer to work with one or the other.
In the end of chapter \ref{Sweedlercontexts}, we study briefly Sweedler theory of other types of (co)algebras.

If $\dgCoalg_\bullet$ and $\dgAlg_\bullet$ are the categories of pointed dg-coalgebras and dg-algebras, there exists pointed Sweedler operations:
\begin{center}
\begin{tabular}{lrl}
\rule[-2ex]{0pt}{4ex} the smash product of coalgebras & $\wedge$&$\!\!\!\!: \dgCoalg_\bullet\times \dgCoalg_\bullet \to \dgCoalg_\bullet$, \\
\rule[-2ex]{0pt}{4ex} the internal hom & $\HOM_\bullet$&$\!\!\!\!:(\dgCoalg_\bullet)^{op}\times \dgCoalg_\bullet \to \dgCoalg_\bullet$, \\
\rule[-2ex]{0pt}{4ex} the pointed Sweedler hom & $\{-,-\}_\bullet$&$\!\!\!\!:\dgAlg_\bullet^{op}\times \dgAlg_\bullet \to \dgCoalg_\bullet$, \\
\rule[-2ex]{0pt}{4ex} the pointed Sweedler product & $\rhd_\bullet$&$\!\!\!\!:\dgCoalg_\bullet\times \dgAlg_\bullet \to \dgAlg_\bullet$, \\
\rule[-2ex]{0pt}{4ex} the pointed convolution product & $[-,-]_\bullet$&$\!\!\!\!:(\dgCoalg_\bullet)^{op}\times \dgAlg_\bullet \to \dgAlg_\bullet$, \\
\rule[-2ex]{0pt}{4ex} and the smash product of algebras & $\wedge$&$\!\!\!\!:\dgAlg_\bullet\times \dgAlg_\bullet \to \dgAlg_\bullet$,
\end{tabular}
\end{center}
such that, for any pointed dg-coalgebra $C$, and any pointed algebras $A$ and $B$ the exists natural bijections between
\begin{center}
\begin{tabular}{lc}
\rule[-2ex]{0pt}{4ex} pointed dg-coalgebra maps & $C\wedge D\to E$, \\
\rule[-2ex]{0pt}{4ex} pointed dg-algebra maps & $C\to \HOM_\bullet(D,E)$, \\
\rule[-2ex]{0pt}{4ex} pointed dg-algebra maps & $D\to \HOM_\bullet(C,E)$,
\end{tabular}
\end{center}
and between
\begin{center}
\begin{tabular}{lc}
\rule[-2ex]{0pt}{4ex} pointed dg-coalgebra maps & $C\to \{A,B\}_\bullet$, \\
\rule[-2ex]{0pt}{4ex} pointed dg-algebra maps & $A\to [C,A]_\bullet$, \\
\rule[-2ex]{0pt}{4ex} pointed dg-algebra maps & $C\rhd_\bullet A\to B$.
\end{tabular}
\end{center}
These last three notions can be constructed similarly to the unpointed ones around a notion a {\em pointed measuring} (see definition \ref{pointedmeasuring}).

All these functors satisfy the same structure as before.
\begin{thm}[theorems \ref{closedmonoidpointedcoalg} and \ref{enrichmentpointedalgcoal}]\ 
\begin{enumerate}
\item The category $(\dgCoalg_\bullet,\wedge, \HOM_\bullet)$ is symmetric monoidal closed.
\item The category $(\dgAlg_\bullet,\{-,-\}_\bullet,\rhd_\bullet,[-,-],\wedge )$ is enriched, tensored, cotensored and symmetric monoidal over $\dgCoalg_\bullet$.
\end{enumerate}
\end{thm}

\medskip
The algebra $\mc$ is naturally pointed by the augmentation sending $u$ to 0 and we have $\dgAlg_\bullet(\mc, [C,A]_\bullet) = Tw_\bullet(C,A)$. Our main application of Sweedler theory is the following result. Recall that $R^c:\dgCoalg_\bullet\to \dgCoalgnilbullet$ is the functor associating to a coalgebra its sub-coalgebra of conilpotent elements.

\begin{thm}[theorem \ref{represtw}]
The adjunction
$$\xymatrix{
-\rhd_\bullet \mc :\dgCoalgnilbullet \ar@<.6ex>[r]&\ar@<.6ex>[l] \dgAlg_\bullet : R^c\{\mc,-\}_\bullet
}$$
coincides with the classical bar-cobar adjunction $\Omega \dashv \Beta $.
\end{thm}

The proof of this theorem is obvious since both adjunctions represent the same bifunctor $Tw_\bullet$ of pointed twisting cochains, but the unravelling of the isomorphisms $\Omega\simeq-\rhd_\bullet \mc$ and $\Beta \simeq R^c\{\mc,-\}_\bullet$ does enlight the classical constructions. For example, the external part of the differentials of the bar and cobar constructions are exactly the differential induced by that of $\mc$ (see section \ref{twistingcochainsection}). This computation also enlight the necessity of some minus signs (see appendix \ref{signtw} for a discussion).

\medskip

From there, it is clear how to define the two variations of the bar-cobar adjunction mentionned above.
The {\em extended bar-cobar adjunction} is the adjunction
$$\xymatrix{
\Omega=-\rhd_\bullet \mc :\dgCoalg_\bullet \ar@<.6ex>[r]&\ar@<.6ex>[l] \dgAlg_\bullet :\{\mc,-\}_\bullet=\Betaext
}$$
representing the bifunctor $Tw_\bullet(C,A)=\dgAlg_\bullet(\mc,[C,A]_\bullet)$,
and the {\em unpointed bar-cobar adjunction} is, as mentionned already, the adjunction
$$\xymatrix{
\bOmega=-\rhd \mc :\dgCoalg \ar@<.6ex>[r]&\ar@<.6ex>[l] \dgAlg :\{\mc,-\}=\bB.
}$$
representing the bifunctor $Tw(C,A)=\dgAlg(\mc,[C,A])$.
These results are also stated in theorem \ref{represtw}.

\medskip
We shall give a few properties of these new adjunctions, but a full study will not be done here.

\subsection*{Other results}

Beside the application to the bar and cobar adjunction we were motivated by the following problems:
\begin{enumerate}
\item understand better the structure of the coalgebras $\HOM(C,D)$ and $\{A,B\}$;
\item extract the classical constructions of $\Beta A$ and $\Omega C$ from $\{\mc,A\}_\bullet$ and $C\rhd_\bullet \mc$;
\item and particularly understand the internal and external parts in the differential of $\Beta A$ and  $\Omega C$ 
\item as well as all the minus signs involved in their definitions and that of universal twisting cochains.
\end{enumerate}
Let us explain our answers.

\paragraph{Distinguished isomorphisms}
Working with Sweedler operations quickly lead us to the following distinguished isomorphisms
\begin{center}
\begin{tabular}{cl}
\rule[-2ex]{0pt}{4ex} $\HOM(C,T^\vee(X)) = T^\vee([C,X])$ &  (proposition \ref{exemplehomcoalg})\\
\rule[-2ex]{0pt}{4ex} $C\rhd T(X) = T(C\otimes X)$ & (proposition \ref{Sweedleroftensor})\\
\rule[-2ex]{0pt}{4ex} $\{T(X),B\} = T^\vee([X,B])$ & (proposition \ref{exSweedlerhom}).
\end{tabular}
\end{center}
They are useful tools for computations but we realized later that they also had a nice interpretation in terms of strong adjunctions.
For example, the isomorphism $\{T(X),B\} = T^\vee([X,B])$ states that the adjunction $T:\dgAlg \rightleftarrows \dgVect: U$ can be enriched over $\dgCoalg$. The enrichment of $\dgAlg$ over $\dgCoalg$ was described above, let us say a word on that of $\dgVect$.
For $X$ and $Y$ two dg-vector spaces, their hom coalgebra is by definition $T^\vee([X,Y])$ where $T^\vee$ is the cofree coalgebra functor and $[X,Y]$ is the dg-vector space of morphisms from $X$ to $Y$. 
Then it should be clear now that the isomorphisms
$$
\HOM(C,T^\vee(X))=T^\vee([UC,X]) \et \{T(X),B\} = T^\vee([X,UB])
$$
expresses the adjunctions $U\dashv T^\vee$ and $T\dashv U$ are enriched over $\dgCoalg$.

\begin{thm}[theorems \ref{strongcomonadicitycoalg} and \ref{strongmonadicityalg}]
The adjunctions
$$\xymatrix{
U : \dgCoalg \ar@<.6ex>[r]& \dgVect:T^\vee \ar@<.6ex>[l]
}
\et
\xymatrix{
T : \dgVect \ar@<.6ex>[r]& \dgAlg:U \ar@<.6ex>[l].
}$$
are enriched and (co)monadic over $\dgCoalg$.
\end{thm}

\bigskip
These isomorphisms together with the (co)monadicity of the categories $\dgCoalg$ and $\dgAlg$ over $\dgVect$ allows to 
give a copresentation of $\HOM(C,D)$ and $\{A,B\}$ in terms of cofree coalgebras.
The reader is refered to corollary \ref{HOMsubcoalg} and remark \ref{calculHOM} for $\HOM(C,D)$ 
and to corollary \ref{betterShom} and remark \ref{calculSHOM} for $\{A,B\}$.

\paragraph{Structure of the hom coalgebras}
Beside the distinguished isomorphisms, 
the following result helps to understand the structure of $\HOM(C,D)$ and $\{A,B\}$.

\begin{thm}[lemmas \ref{underlyingCoalg} and \ref{underlyingAlg}, corollaries \ref{corprimicoderiv} and \ref{corprimideriv}, theorems \ref{primiLiecoderiv} and \ref{primiLiederiv}]\label{thmstructurehomcoalg}
$ $
\begin{enumerate}
\item The atoms of $\HOM(C,D)$ are in bijection with the coalgebra maps $C\to D$.
\item The atoms of $\{A,B\}$ are in bijection with the algebra maps $A\to B$.
\item If $f\in \HOM(C,D)$ is a atom, the dg-vector space $\Prim_f(\HOM(C,D))$ of $f$-primitive elements of $\HOM(C,D)$ is isomorphic to the the dg-vector space $\Coder(f)$ of $f$-coderivations $C\to D$.
\item If $f\in \{A,B\}$ is a atom, the dg-vector space $\Prim_f(\{A,B\})$ of $f$-primitive elements of $\{A,B\}$ is isomorphic to the the dg-vector space $\Der(f)$ of $f$-derivations $A\to B$.
\item If $C=D$, $A=B$ and the atoms $f$ are the identity maps, the previous isomorphisms enhanced into Lie algebra isomorphisms
$$
\Prim(\END(C)) = \Coder(C)
\et 
\Prim(\{A,A\}) = \Der(A).
$$
Moreover these Lie algebras are naturally equipped with a square operation for odd elements and these isomorphisms commutes to these operations.
\label{thmstructurehomcoalgLie}
\end{enumerate}
\end{thm}

The coalgebras $\HOM(C,D)$ and $\{A,B\}$ are by definition equipped with evaluation maps
${\bf ev}:\HOM(C,D)\otimes C\to D$ and ${\bf ev}:\{A,B\}\rhd A\to B$ which have a certain universal property.
These maps gives maps $\Psi:\HOM(C,D)\to [C,D]$ and $\Psi:\{A,B\}\to [A,B]$ that we call {\em reduction maps}.
The maps $\Psi$ inherit the universal properties of the $\bf ev$s (see definitions \ref{remreductioncoalg} and \ref{defcomeasuring} of {\em comorphism} and {\em comeasuring}) and are useful in many proofs (see section \ref{proofbyreduction}). They also have a nice interpretation as part of enriched functors (propositions \ref{strongforgetfulVectcoalg} and \ref{strongforgetfulVectalg}). 

In particular, we have maps of algebras
$$
\Psi:\END(C)\to [C,C] \et \Psi:\{A,A\}\to [A,A].
$$
These maps sends atoms to (co)algebra maps and primitives to (co)derivations, they are the main tool to prove the above isomorphisms.
They are also useful to deal with actions of bialgebras (propositions \ref{hombialg} and \ref{hombialgalg} which are the key to theorems 
\ref{Qhomcoalg} and \ref{Qenrichmentalgcoal}).

All these results have non-(co)unital and pointed analogs developped in sections \ref{nonunitalsweedlertheory} and \ref{pointedsweedlertheory}.

\paragraph{Meta-morphisms}
The fact that (co)derivations are primitive elements of the enriched homs $\HOM(C,D)$ and $\{A,B\}$ lead us to the remark that they could be transported by the Sweedler operations and that they could explain the decomposition of the differential in the bar and cobar constructions. 

This gave us the motivation to introduce and study what we called {\em meta-morphisms}.
Meta-morphisms are simply the elements of $\HOM(C,D)$ and $\{A,B\}$, they can be composed and evaluated at some element of the domain.
The composition is done using the strong compositions maps
${\bf c}:\HOM(D,E)\otimes \HOM(C,D)\to \HOM(C,E)$, ${\bf c}:\{B,E\}\otimes \{A,B\}\to \{A,E\}$ 
and the evaluations are done using the strong evaluation maps ${\bf ev}:\HOM(C,D)\otimes C\to D$, ${\bf ev}:\{A,B\}\rhd A\to B$.
Moreover all Sweedler operations are strong and meta-morphisms can be passed through any of them: they can be tensored, $\HOM$-ed, $\{-,-\}$-ed, $\rhd$-ed and $[-,-]$-ed. For example, the Sweedler product produce a map of coalgebras
$$\xymatrix{
\rhd:\HOM(C,D)\otimes \{A,B\} \ar[r] & \{C\rhd A,D\rhd B\}.
}$$
Most computations with meta-morphisms are easy, using them is no more difficult than taking seriously that they are indeed morphisms.

\medskip
To understand our use of meta-morphisms, let us explain a bit of context. 
We developped the Sweedler theory of dg-(co)algebras but it should be clear to the reader that we can drop the differential and develop the same theory for graded (co)algebras. 
Then, it is a classical method to study dg-objects by first building a graded object and then enhancing it with a natural differential.
Having both differential and graded Sweedler theories, it is natural to think that if one forgets the differential of the Sweedler dg-operations, one obtains the corresponding Sweedler graded operation. 

\begin{thm}[theorems \ref{dgtoghomcoalg} and \ref{dgtoghomalg}]
The forgetful functors
$$\xymatrix{
\dgCoalg\ar[r]& \gCoalg 
}\et
\xymatrix{
\dgAlg\ar[r]& \gAlg 
}$$
preserve all Sweedler operations.
\end{thm}

Moreover we proved that the graded Sweedler operation applied to dg-(co)algebras is automatically enhanced with a differential that can be computed using the calculus of meta-morphisms. This computation is exactly what is needed to understand the two parts of the differentials
of the bar and cobar constructions. Let us explain this.

\paragraph{The explicit description of the bar-cobar constructions}
As a graded object the algebra $\mc=T(u)$ is free on one generator of degree $-1$ and
the differential $d_\mc$ is defined by $d_\mc u=-u^2$.
We can use the pointed analogs of the above distinguished isomorphisms
$$
C\rhd_\bullet T_\bullet(u) = T_\bullet(C_-\otimes u)
\et
R^c\{T(u),B\}_\bullet = T^c_\bullet([u,B_-])
$$
to compute the underlying graded object of $\Beta A=R^c\{\mc,A\}_\bullet$ and $\Omega C=C\rhd_\bullet \mc$.
With the identification $u=s^{-1}$ and $u^\star=s$ we found the classical formulas for those objects as written for example in \cite{LV}.

Then, the differentials $d_C$, $d_A$ and $d_\mc$ of $C$, $A$ and $\mc$ can be seen as primitive elements in $\END_\bullet(C)$, $\{A,A\}_\bullet$ and $\{\mc,\mc\}_\bullet$ respectively.
The functors $\rhd_\bullet$ and $\{-,-\}_\bullet$ are enriched over $\dgCoalg$ and gives maps of bialgebras
\begin{eqnarray*}
\rhd_\bullet:\END_\bullet(C)\otimes \{\mc,\mc\}_\bullet &\tto &  \{C\rhd_\bullet \mc,C\rhd_\bullet \mc\}_\bullet\\
	f\otimes g &\mto & f\rhd_\bullet g\\
R^c\{-,-\}_\bullet:\{\mc,\mc\}_\bullet^o \otimes \{A,A\}_\bullet &\tto &  \END_\bullet(\{\mc,A\}_\bullet)\\
	f\otimes g &\mto & \{f, g\}_\bullet
\end{eqnarray*}
which induce Lie algebra maps between primitive elements of pointed (co)derivations
\begin{eqnarray*}
\rhd'_\bullet:\Coder_\bullet(C)\times \Der_\bullet(\mc) &\tto & \Der_\bullet(C\rhd_\bullet \mc)\\
	(d_C,d_\mc) &\mto & d_C\rhd_\bullet 1_\mc + 1_C\rhd_\bullet d_\mc \\
R^c\{-,-\}'_\bullet:\Der_\bullet(\mc) \times \Der_\bullet(A) &\tto & \Coder_\bullet(\{\mc,A\}_\bullet)\\
	(d_\mc,d_A) &\mto & \{1_\mc, d_A\}_\bullet - \{d_\mc, 1_A\}_\bullet
\end{eqnarray*}
The images of $d_C$ and $d_A$ are exactly the two internal parts of the bar and cobar differentials, and the two images of $d_\mc$
are exactly the two external parts, the minus sign in the second formula comes from the contravariance in the first variable
(see the computation after theorem \ref{represtw}).
Moreover, from the calculus of meta-morphisms it is trivial to check that the image of a square zero (co)derivation is still of square zero.

\medskip
We study a few consequences of our abstract formalism in sections \ref{csqsweedlerformalism} and \ref{iteratedbarcobar}.
Let us only mention that all Sweedler operations have a lax or colax structure with respect to the tensor products of algebras and coalgebras 
(proposition \ref{strength for otimes and hom of coalg} and corollaries \ref{convolutionlax}, \ref{Sproductcolax} and \ref{Shomlax})
and it is possible to extract from this the classical (co)shuffle bialgebra structure on bar and cobar constructions applied to (co)commutative (co)algebras (theorem \ref{barcobarcommutative}).

\paragraph{Signs issues}
Signs problems have been our worst enemy studying the bar and cobar constructions. Multiple conventions exists for the signs of the differentials and the universal twisting cochains, we have found helpful to write a comparative table of the notations from different references which we have included in appendix \ref{signtw}. Whatever conventions are chosen, signs cannot totally disappear and the appendix details why. It also argue for our conventions, which are different from everybody else.

\subsection*{Acknowledgments}

The work of J.-L. Loday and B. Vallette on Koszul duality was an inspiring source for the present work, especially their book \cite{LV}.
\medskip

The first author was supported by the \href{http://www.cirget.uqam.ca/}{CIRGET} at UQ\`AM and at the \href{http://www.math.ethz.ch/}{ETH} in Zürich. Both institutions are gracefully acknowledge for their support and funding.

\newpage
\section{Elementary dg-algebra}\label{conventions}\label{basicalgebra}

\subsection{dg-Vector spaces}

\paragraph{Graded sets}

Recall that  a  {\it $\mathbb{Z}$-graded set}  is a set $E$ equipped with a map $|-|:E\to \mathbb{Z}$ called the {\it degree map}. We shall say that an element $x\in E$ has {\it degree} $|x|$ and we shall denote by $E_n$ the set of elements $x\in E$ of degree $|x|=n\in \mathbb{Z}$.
We shall say that $x\in E$ is {\it odd} (resp. {\it even}) if $|x|$ is odd (resp. even).
A graded set $(E,|-|)$ is said to be {\it graded finite} (resp {\it finite}) if the set $E_n$ is finite for every $n\in \mathbb{Z}$ (resp. if $E$ is  finite).
A {\it map} of graded sets is a map $f:E\to F$ preserving the degrees. 
We shall denote the category of graded sets by $\gSet$.

\medskip
$\gSet$ is a monoidal category when equipped with the graded product
$$
(E\times F)_n = \coprod_{i+j=n} E_i\times F_j.
$$

\paragraph{Graded vector spaces}

Recall that a {\it $\mathbb{Z}$-graded vector space} is a family of pairwise disjoint vector spaces $X=(X_n:n\in \mathbb{Z})$ indexed by the set of integers $\mathbb{Z}$. 
We shall say more simply that $X$ is a {\it graded vector space} if the context is clear. 
We say that an element $x\in X_n$ has {\it degree} $n$ and we write $|x|=n$. 
We shall often identify $X$ with its {\it total space}
$$
X=\bigoplus_{n\in \mathbb{Z} } X_n.
$$

A {\it map} of graded vector spaces $f:X\to Y$ is a family of linear maps $f_n:X_n\to Y_n$. 
We shall denote the category of graded vector spaces  by $\gVect$.
The obvious forgetful functor $\gVect\to \gSet$ has a left adjoint which associates to a graded
set $E=(E_n:n\in \mathbb{Z})$ the graded vector space $\FF E =(\FF E_n:n\in \mathbb{Z})$.
We shall sometimes call the elements of $E$ {\em graded variables} and says that a graded vector space $X$ is generated by the set of graded variable $E$ if $X=\FF E$.

\begin{defi}\label{deffinite}

We shall say that a graded vector space $X$ is {\it graded finite} (resp. {\it finite}) if the vector space $X_n$ is finite dimensional for every $n\in \mathbb{Z}$ (resp. if the total space of $X$ is finite dimensional).
We shall say that a graded vector space $X$ is {\em positive} (resp. {\em negative}) if $X_i=0$ for $i<0$ (resp. for $i>0$).
We shall say that $X$ is {\em strictly positive} (resp. {\em strictly negative}) if $X_i=0$ for $i\leq 0$ (resp. for $i\geq0$).
We shall say that $X$ is {\em bounded below} (resp. {\em bounded above}) if $X_i=0$ for $i\ll0$ (resp. for $i\gg0$).

\end{defi}

We shall say that a family $(X(i):i\in I)$ of graded vector spaces is {\it locally finite} if for every $n\in \mathbb{Z}$
we have $X(i)_n=0$ except for a finite number of elements $i\in I$.

\begin{lemma}\label{convergence}
Let $(X(i):i\in I)$ be a locally finite family of graded vector spaces.
Then the canonical map
$$\xymatrix{
\bigoplus_{i\in I} X(i)\ar[r]& \prod_{i\in I} X(i)
}$$
is an isomorphism.
\end{lemma}

\begin{proof} 
It suffices to show that the canonical map
$$\xymatrix{
\bigoplus_{i\in I} X(i)_n \ar[r]& \prod_{i\in I} X(i)_n
}$$
is an isomorphism for every $n\in \mathbb{Z}$.
This reduces the problem to the case where $I$ is finite,
since the family $(X(i):i\in I)$ is locally finite.
But the result is true in this case, since
the category of vector spaces is additive.
\end{proof}

In particular, a graded vector space $X$ is both the sum and the product of the family of graded vector spaces $X(i)$ where $X(i)_i=X_i$ and $X(i)_n=0$ otherwise.

\bigskip

\paragraph{Graded tensor product}

The {\it tensor product} of two graded vector spaces $X$ and $Y$ is the graded vector space $X\otimes Y$ defined by putting
$$
(X\otimes Y)_n = \bigoplus_{i+j=n } X_i\otimes Y_j
$$
for every $n\in\mathbb{Z}$.
This defines a monoidal structure on the category $\gVect$ for which the unit object is the field $\FF$ viewed as a graded vector space concentrated in dimension $0$.  If $E$ is a graded set, we shall put $E\otimes X= \FF E\otimes X$
and $X\otimes E= X\otimes \FF E$.
The category $\gVect$ is symmetric monoidal (Appendix \ref{Catmonoidales}) with a symmetry $\sigma=\sigma(X,Y):X\otimes Y \to Y\otimes X$
given by the  {\it Koszul rule} 
$$
\sigma(x\otimes y)=y\otimes x\ (-1)^{|x||y|}.
$$
Remark that the forgetful functor $U:\gVect\to \gSet$ is monoidal but not symmetric monoidal: the sign rule cannot be given a sense in $\gSet$.

\paragraph{Graded hom}\label{lambdatransform}
 
The monoidal category $\gVect$ is closed;
the internal hom is the graded vector space $[X,Y]$
obtained by putting
$$
[X,Y]_n=\prod_{i}[X_i,Y_{i+n}]
$$
for every $n\in \mathbb{Z}$.
We shall say that an element $f\in [X,Y]_n$ is a {\it graded morphism} 
and  write $f:X\rhup Y$ or $f:X\rhup_n Y$.
The adjunction
$$\lambda^2:[X\otimes Y,Z]\simeq [X,[Y,Z]]$$
takes a graded morphism $f:X\otimes Y \rhup Z$ to the graded morphism $\lambda^2(f):X\rhup [Y,Z]$
defined by putting $$\lambda^2(f)(x)(y)= f(x\otimes y).$$
The other adjunction
$$\lambda^1: [X\otimes Y,Z]\simeq [Y,[X,Z]]$$
takes a graded morphism $f:X\otimes Y \rhup Z$ to the graded morphism $\lambda^1(f):Y\rhup [X,Z]$
defined by putting $$\lambda^1(f)(y)(x)= f(x\otimes y)(-1)^{|x||y|}.$$

\medskip

The counit of the adjunction $\lambda^2$ is the evaluation $ev:[Y,Z]\otimes Y\to Z$
defined by putting $ev(f \otimes y )=f(y)$
and its unit is the map $\eta:X\to [Y,X\otimes Y]$
defined by putting $\eta(x)(y)=x\otimes y$.
The counit of $\lambda^1$
is the reversed evaluation $rev:X\otimes [X,Z]\to Z$
defined by putting $rev(x\otimes f ) =f(x)(-1)^{|f||x|}$.

\medskip

The composition law
$$\xymatrix{
c:[Y,Z]\otimes [X,Y]\ar[r]& [X,Z]
}$$
takes a pair of graded morphisms $(g,f)\in [Y,Z]_m\otimes [X,Y]_n$
to their composite $gf=g\circ f \in [X,Z]_{m+n}$.
This defined a category $\gVect^\#$ enrichied over $\gSet$
whose objects are the graded vector spaces and whose
morphisms are the graded morphisms.
The category $\gVect$ is the subcategory of graded morphisms
of degree 0 of the category $\gVect^\#$.

\paragraph{Monoidal strength}

The tensor product functor
$\otimes: \gVect \times \gVect\to \gVect$
is strong (=enriched) by a general result of category theory \cite{AM}. 
Its strength is given by the map
$$\xymatrix{
\theta_\otimes :[X,X']\otimes [Y,Y']\ar[r]& [X\otimes Y, X'\otimes Y']
}$$
defined by putting
$(f\otimes g)(x\otimes y)=f(x)\otimes g(y)\ (-1)^{|g| |x|}$
for $f:X\rhup X'$, $g:Y\rhup Y'$, $x\in X$ and $y\in Y$.
Moreover, if $k:X'\to X''$ and $r:Y'\to Y''$, then
$$
(k\otimes r)(f\otimes g)=kf\otimes rg\ (-1)^{|r| |f|}.
$$
In particular $\theta_\otimes :[X,X]\otimes [Y,Y]\to [X\otimes Y, X\otimes Y]$ is an algebra map.

\medskip

We should be careful in taking the tensor product of commutative diagrams of graded morphisms.
For example, the tensor product of two commutative triangles of graded morphisms
$$\xymatrix{
X  \ar@^{>}[rrd]_{h} \ar@^{>}[rr]^-{f}  && Y  \ar@^{>}[d]^{g} \\
&& Z 
}
\quad\qquad
\xymatrix{
X'  \ar@^{>}[rrd]_{h' } \ar@^{>}[rr]^-{f'}  && Y'  \ar@^{>}[d]^{g'} \\
&& Z' 
}$$
is a triangle of graded morphisms
$$\xymatrix{
X\otimes X'  \ar@^{>}[rrd]_{h\otimes h' } \ar@^{>}[rr]^-{f\otimes f'}  && Y\otimes Y'  \ar@^{>}[d]^{g\otimes g'} \\
&& Z \otimes Z'
}$$
which commutes up only up to the sign, since we have
$$
h\otimes h' =( g\circ f)\otimes (g'\circ f')= (g\otimes g')\circ (f\otimes f')(-1)^{|g'||f|}.
$$
Hence the category $\gVect^\#$ is {\bf not} monoidal (but the sub-category of graded morphisms of even degree is).
This is related to the fact that the forgetful functor $U:\gVect\to \gSet$ is not symmetric monoidal.

\medskip

The hom functor $[-,-]: \gVect^{op}\times \gVect\to \gVect$ is strong by a general result of category theory (see appendix \ref{appclosedcategories}). 
Its strength is given by the map
$$\xymatrix{
\theta_{[-,-]}:[X',X]\otimes [Y,Y']\ar[r]& [[X,Y], [X',Y']], 
}$$
obtained by putting $\theta_{[-,-]}(f\otimes g)(u)=g u f (-1)^{|f|(|u|+|g|)}$ for $f:X'\rhup X$, $g:Y\rhup Y'$ and $u:X'\rhup Y$. 
We will note for short $\hom(f,g)$ or $[f,g]$ instead $\theta_{[-,-]}(f\otimes g)$. 
There is a risk of confusion between the notation $[f,g]$ and that of a Lie bracket, however the context should prevent any mistake.
We will also put $\hom(X,g)$ or $[X,g]$ and $\hom(f,Y)$ or $[f,Y]$ instead $\theta_{[-,-]}(1_X\otimes g)$ and  $\theta_{[-,-]}(f\otimes 1_Y)$ respectively.
If $k:X''\rhup X'$ and $r:Y'\rhup Y''$, we have
$$
[fk,rg] =[k,r][f,g](-1)^{|f|(|k|+|r|)}.
$$

In particular $\theta_{[-,-]}:[X,X]^{op}\otimes [Y,Y]\to [[X,Y], [X,Y]]$ is an algebra map.

\paragraph{Suspension}\label{suspension}

We shall denote by $S^n$ the graded vector space freely generated by one element $s^n$ of degree $n\in\mathbb{Z}$.
By definition, $S^0=\FF$ if we put $s^0=1$.
The  {\it $n$-fold suspension} or the {\em $(-n)$-fold desuspension} $S^n(X)$ of a graded vector space $X$ is defined by putting $S^n(X)=S^n\otimes X$.
By definition,  $S^n(X)_i =  s^n\otimes X_{i-n}$ for every $i\in\mathbb{Z}$. 

We shall be careful with the notation $s^n$ which suggests that we have formulas $s^m\otimes s^n = s^{m+n}$ for any $m,n\in \ZZ$.
Such formulas are false because of the degree of the elements: if $n\geq 0$ is odd, we have $s^n\otimes s^{-n} = -s^{-n}\otimes s^n$ so both cannot be equal to $s^0=1$. In particular, $s$ and $s^{-1}$ are not inverse to each other.
Rather, it is convenient to think $s$ and $s^{-1}$ as dual to each other: if we set $s\otimes s^{-1} = 1$ (and thus $s^{-1}\otimes s=-1$), we can identify $s=(s^{-1})^\star$.

\paragraph{Differential graded vector spaces}\label{dgvectorspaces}
A graded endomorphism $d:X\rhup X$ of degree $-1$ is called a {\em differential} or a {\em boundary operator} if $d^2=0$.
A {\em dg-vector space}, or a {\em complex of vector spaces} is defined as a graded vector space $X$ together with the choice of a differential $d_X$.
A map of dg-vector spaces $(X,d_X)\to (Y,d_Y)$ is a map of graded vector spaces $f:X\to Y$ such that $d_Yf=fd_X$ as graded morphisms of degree $-1$.
The category of dg-vector spaces is noted $\dgVect$ (\footnote{We have prefered the name {\em dg-vector space} rather than {\em complex} to emphasize the relative nature of Sweedler theory with respect to its basic objects. In particular the theory could be developped for actual vector spaces, or graded vector spaces, or any objects which behave like vector spaces.}).

If $X$ is a dg-vector space, we shall sometimes use the notation $|X|$ to refer to the underlying graded vector space and note $|-|$ or $U_d$ the functor $\dgVect \to \gVect$ forgetting the differential.

\medskip

Any graded vector space can be enhanced to a dg-vector space with zero differential. This defines a fully faithful embedding
$\gVect\to \dgVect$. We will say that a dg-vector space is generated by a graded $E$ if it is $\FF E$ viewed with zero differential.
In particular, the graded vector space $S^n=\FF s^n$ will be seen as such.

If $X$ is a dg-vector space, a {\em dg-vector subspace} of $X$ is a graded subspace stable by the differential.

\medskip
The category $\dgVect$ inherits of all the structures of $\gVect$.
The graded tensor product is enhanced in a tensor of dg-vector space with differential
$$
d_{X\otimes Y} = d_X\otimes Y + X\otimes d_Y.
$$
which is explicitely given on elements by $d_{X\otimes Y}(x\otimes y) = (d_Xx)\otimes y + x\otimes d_Yy\ (-1)^{|x|}$.
This monoidal structure is still symmetric and closed. 
The unit object is given by the vector space $\FF$ concentrated in degree 0 with 0 differential.
The internal hom between $(X,d_X)$ and $(Y,d_Y)$ is given by the graded vector space $[X,Y]$ equipped with the differential
$$
d_{[X,Y]} = [X,d_Y] - [d_X,Y]
$$
which is explicitely given on elements by $d_{[X,Y]}(f) = d_Yf - fd_X\ (-1)^{|f|}$.
In particular, we see that maps of dg-vector spaces $X\to Y$ are in bijection with the subspace $Z_0[X,Y]$ of 0-cycles of $[X,Y]$.

\medskip

It is useful to regard a complex $(X,d_X)$ as a module over the graded commutative algebra $\d=\FF[\delta]/(\delta^2)$, where  $\delta$ is a symbol of degree $-1$. $\d$ is actually a cocommutative Hopf algebra with the coproduct $\Delta:\d\to \d\otimes \d$
defined by $\Delta(\delta)=\delta\otimes 1+1\otimes \delta $ and the antipode defined as the unique map of algebras $S:\d\to \d$ such that $S(\delta)=-\delta$ ($\d$ is the $\d_{-1}$ of example \ref{primitiveHopf}).
The previous symmetric monoidal closed structure of $\dgVect$ is the one associated to the Hopf algebra $\d$.

\begin{prop}\label{dVectvsVect}
The category $\dgVect$ is symmetric closed monoidal and the forgetful functor $|-|:\dgVect \to \gVect$ is monoidal and preserve internal hom.
Moreover, this functor is faithful and admits a left adjoint $X\to \d\otimes X$ and a right adjoint $X\mapsto [\d,X]$.
Hence limits and colimits exist in $\dgVect$ and they can be computed in $\gVect$.
\end{prop}

\begin{proof}
The first part is obvious by construction, but the whole result is a general fact for Hopf algebras and will be proven in proposition \ref{freeandcofreeQmodules}.
\end{proof}

\paragraph{Graded dual}

The graded dual $X^{\star}=[X,\FF]$ of a dg-vector space $X$ is obtained by putting $(X^{\star})_n=(X_{-n})^{\star}$ and the evalutation map $ev:X^\star \otimes X\to \FF$ by putting $ev(\phi\otimes x)=\phi(x)$.
The {\it transpose} of a graded morphism $f:X\rhup_n Y$ is the graded morphism ${}^t\! f:=[f,\FF]:Y^\star \rhup_{-n} X^\star$ defined by putting ${}^t\! f(\phi)=\phi\circ  f(-1)^{|\phi||f|}$ for every $\phi\in X^\star$. 
If $f:X\rhup Y$ and $g:Y\rhup Z$, then 
$$
{}^t\!(g\circ f)={}^t\!f\circ {}^t\!g\ (-1)^{|f||g|}.
$$
The contravariant functor $(-)^\star: \dgVect^{op} \to \dgVect$ is lax monoidal, with a  lax structure is given by the map
$$\xymatrix{
X^\star \otimes Y^\star \ar[r]& (X\otimes Y)^\star
}$$
defined by putting $(\phi\otimes \psi)(x\otimes y)=\phi(x)\psi(y)(-1)^{|x||\psi}$.

\medskip
If $X=\FF E$ is generated by a graded finite set, $X_n$ is generated by the finite set $E_n$ and thus $(X^\star)_n=(X_{-n})^\star$ is generated by the dual of the basis $E_{-n}$.
Precisely, each element $e\in E_{-n}$ defines a dual map $e^\star:X=\FF E_{-n}\to \FF$ which sends $e$ to 1 and all other elements of $E$ to 0. $e^\star$ is called the {\em dual} of the element $e$. 

Let $E^\star$ be the graded set such that $E^\star_n = E_{-n}$, the canonical isomorphisms $(X^\star)_n =\FF E_{-n}$ give a canonical isomorphism $X^\star = \FF (E^\star)$. 
$X^\star \otimes X$ is generated by $E^\star \times E$ and the evaluation map $ev:X^\star \otimes X \to \FF$ is given $e^\star\otimes f\mapsto e^\star(f)$ if $|e|+|f|=0$ and $e^\star\otimes f\mapsto 0$ if not.

\medskip
The canonical map $i_X:=\lambda^1(ev):X\to X^{\star\star}$ takes an element $x\in X_n$ to the element $(-1)^ni_n(x)\in (X_n)^{\star\star}$, where $i_n$ is the canonical map $X_n\to (X_n)^{\star\star}$.
The map $i_X:X\to X^{\star \star}$ is invertible if and only if $X$ is graded finite.

\medskip

If $X$ and $Y$ are graded vector spaces, we shall use the canonical maps
$$\xymatrix{
Y\otimes X^\star \ar[r]& [X,Y]
} \et 
\xymatrix{
X^\star \otimes Y \ar[r]& [X,Y]
}$$
defined by putting $(y\otimes \phi)(x)=y\phi(x) =\phi(x)y$ and $(\phi \otimes y)(x)=\phi(x)y\ (-1)^{|x||y|}$.
The two maps are invertible when $X$ or $Y$ is finite.

\medskip

A dg-vector space is said to be {\em finite} (resp. {\em graded finite}) if its underlying graded space is.
We shall denote  by $\dgVectfin$ the full subcategory of $\dgVect$ spanned the finite dg-vector spaces. 
The category $\dgVectfin$ is symmetric monoidal closed and the duality functor  $(-)^\star:\dgVectfin^{op}\to \dgVectfin$ is a contravariant equivalence of symmetric monoidal categories.

\medskip

A dg-vector space is said to be {\em bounded above} (resp. {\em bounded below}) if its underlying graded space is.
We denote by $\dgVectb$ (resp. $\dgVecta$) the full subcategory of $\dgVect$ whose objects are the dg-vector spaces bounded below (resp. bounded above).
The tensor product of two dg-vector spaces bounded below (resp. above) is bounded below (resp. above).
Hence the subcategories $\dgVectb$ and $\dgVecta$  are monoidal.
  
\medskip

A dg-vector space is said to be {\em graded finite}) if its underlying graded space is.
We denote by $\dgVectgrfina$ (resp. $\dgVectgrfinb$) the full subcategory of $\dgVectgrfin$ whose objects are the graded finite dg-vector spaces bounded above (resp. bounded below).
Because of the bound condition, the subcategories $\dgVectgrfina$ and $\dgVectgrfinb$ are stable by the tensor product.
In particular we have the important result that the duality functor $(-)^\star: \dgVect^{op} \to \dgVect$ induces an equivalence of categories 
$(\dgVectgrfina)^{op} \simeq \dgVectgrfinb$.

\begin{lemma}\label{lemmagradedfinite}
The duality functors
$$\xymatrix{
(-)^\star:  (\dgVectgrfina)^{op} \ar@<.6ex>[r]& \dgVectgrfinb : (-)^\star \ar@<.6ex>[l]
}$$
are inverse equivalence of monoidal categories.
\end{lemma}

\paragraph{Triality}\label{trialitydgVect}

In this last section, we present the monoidal structure of $\dgVect$ using the language of triality of appendix \ref{triality}. 
This is useful to understand the measuring functor of section \ref{measuring}.

For $X$, $Y$ and $Z$ three dg-vector spaces, we shall say that a map of graded sets $X\times Y\to Z$ is {\em bilinear} if it is linear and differential in each variable.
There are canonical bijections between
\begin{center}
\begin{tabular}{lc}
\rule[-2ex]{0pt}{4ex} bilinear maps & $X\times Y\to Z$, \\
\rule[-2ex]{0pt}{4ex} linear maps & $f:X\otimes Y\to Z$, \\
\rule[-2ex]{0pt}{4ex} linear maps & $\lambda^1f:Y\to [X,Z]$, \\
\rule[-2ex]{0pt}{4ex} linear maps & $\lambda^2f:X\to [Y,Z]$, \\
\end{tabular}
\end{center}
Let $\cT(X,Y;Z)$ be the set of bilinear maps $X\times Y\to Z$, we have isomorphisms $\cT(X,Y;Z)=Z_0[X\otimes Y,Z]=Z_0[Y,[X,Z]]$.

The sets $\cT(X,Y;Z)$ define a triality
$$\xymatrix{
\cT:\dgVect^{op}\times \dgVect^{op}\times \dgVect \ar[r]& \Set
}$$
and the above bijections proves that this triality is representable in each on its variable.

The functor $\cT$ is lax monoidal when $\dgVect$ is equipped with the graded tensor product and $\Set$ with the cartesian product.
The lax monoidal structure is given by the pair $(\alpha,\alpha_0)$
where $\alpha_0$ is the surjection $\cT(\FF,\FF;\FF)\to \{*\}$ and
$$\xymatrix{
\alpha:\cT(X_1,Y_1;Z_1)\times \cT(X_2,Y_2;Z_2)\ar[r]& \cT(X_1\otimes X_2,Y_1\otimes Y_2;Z_1\otimes Z_2)
}$$
is the map
$$\xymatrix{
Z_0[X_1\otimes Y_1,Z_1]\times Z_0[X_2\otimes Y_2,Z_2]
\ar[r]& Z_0[X_1\otimes Y_1\otimes X_2\otimes Y_2,Z_1\otimes Z_2]
\ar[r]& Z_0[X_1\otimes X_2\otimes Y_1\otimes Y_2,Z_1\otimes Z_2]
}$$
where the first map is the functoriality of $\otimes$ and the second the permutation $Y_1\otimes X_2\simeq X_2\otimes Y_1$.

The lax condition on $\alpha_0$ is immediate and on $\alpha$, it amounts to the commutativity of the square
$$\xymatrix{
Z_0[X_1\otimes Y_1,Z_1]\times Z_0[X_2\otimes Y_2,Z_2]\times Z_0[X_3\otimes Y_3,Z_3]\ar[r]\ar[d]&  Z_0[X_1\otimes X_2\otimes Y_1\otimes Y_2,Z_1\otimes Z_2]\times Z_0[X_3\otimes Y_3,Z_3]\ar[d]\\
Z_0[X_1\otimes Y_1,Z_1]\times Z_0[X_2\otimes X_3\otimes Y_2\otimes Y_3,Z_2\otimes Z_3]\ar[r]&  Z_0[X_1\otimes X_2\otimes X_3\otimes Y_1\otimes Y_2\otimes Y_3,Z_1\otimes Z_2\otimes Z_3]
}$$
which is easily deduced from the universal property of $\otimes$ as representing bilinear maps.

It will be useful in the proof of proposition \ref{laxmeasuring} to reformulate this condition as the commutation of the diagram
$$\xymatrix{
Z_0[Y_1,[X_1,Z_1]]\times Z_0[Y_2,[X_2,Z_2]]\times Z_0[Y_3,[X_3,Z_3]]\ar[r]\ar[d]&  Z_0[Y_1\otimes Y_2,[X_1\otimes X_2,Z_1\otimes Z_2]]_0\times [Y_3,[X_3,Z_3]]\ar[d]\\
Z_0[Y_1,[X_1,Z_1]]\times Z_0[Y_2\otimes Y_3,[X_2\otimes X_3,Z_2\otimes Z_3]]\ar[r]&  Z_0[Y_1\otimes Y_2\otimes Y_3,[X_1\otimes X_2\otimes X_3,Z_1\otimes Z_2\otimes Z_3]]
}$$

We leave to the reader the proof that the lax structure $(\alpha,\alpha_0)$ is symmetric.

\subsubsection{Pointed vector spaces}\label{pointedvectorspaces}

A dg-vector space $X$ will be called {\em pointed} if it is equipped with two maps $e:\FF\to X$ and $\epsilon:X\to \FF$ such that $\epsilon e= id_\FF$. We shall call $e$ the {\em unit} and $\epsilon$ the {\em counit}. For simplicity, we shall also denote by $e$ the element $e(1)\in X$. We have a canonical isomorphism $\ker \epsilon \simeq X/\FF e$. 

A {\em map} of pointed dg-vector spaces $(X,e_X,\epsilon_X)\to (X,e_Y,\epsilon_Y)$ is a map of dg-vector spaces $f:X\to Y$ commuting with the units and counits, \ie such that $fe_X=e_Y$ and $\epsilon_Y f=\epsilon_X$.
We shall denote by $\dgVect_\bullet$ the category of pointed dg-vector spaces.

\medskip
If $X$ is a pointed dg-vector space, we shall denote $X_-= \ker \epsilon \simeq X/\FF e$.
Let $i:X_-\to X$ and $p:X\to X_-$ be the canonical inclusion and projection maps.
The natural maps $p\oplus \epsilon:X\to X_-\oplus \FF$ and $i\otimes e:X_-\otimes \FF \to X$ are inverse isomorphisms

If $X$ is a dg-vector space, the space $X_+:=X\oplus \FF$ is canonically pointed.
These constructions define an equivalence of categories
$$\xymatrix{
(-)_+:\dgVect \ar@<.6ex>[r]&\dgVect_\bullet:(-)_-. \ar@<.6ex>[l]
}$$
We define the {\em wedge} product $X\vee Y$ of two pointed dg-vector spaces $X$ and $Y$ by the formula
$$
X\vee Y := (X_-\oplus Y_-)_+.
$$
This operation corresponds by the equivalence to the direct sum $\oplus$, it is therefore a sum and a product in $\dgVect_\bullet$.

Using the equivalence we can also transfer the monoidal structure of $\dgVect$ to $\dgVect_\bullet$.
We define the {\em smash} product $X\wedge Y$ of two pointed dg-vector spaces $X$ and $Y$ and the corresponding internal hom $[X,Y]_\bullet$ by the formulas
$$
X\wedge Y := (X_-\otimes Y_-)_+
\et 
[X,Y]_\bullet := [X_-,Y_-]_+.
$$

There exist a canonical decomposition 
$$
X\otimes Y = X_-\otimes Y_-\oplus X_-\oplus Y_-\oplus \FF.
$$
From which we see that $X\wedge Y = X_-\otimes X_-\oplus \FF$ is a retract of $X\otimes Y$.
Let $$
\alpha: X\wedge Y \to X\otimes Y
\et
\beta: X\otimes Y\to X\wedge Y
$$
be the canonical inclusion and projection.
The maps $\alpha$ and $\beta$ define a bilax monoidal structure on the forgetful functor
$$\xymatrix{
U:\dgVect_\bullet \ar[r]&\dgVect.
}$$

\begin{lemma}\label{pointedbicart}
Let $X$ and $Y$ be two pointed vector spaces, then the squares
$$
\vcenter{\xymatrix{
X\wedge Y\ar[r]\ar[d] &\FF\ar[d]\\
X\otimes Y\ar[r] & X\oplus Y
}}
\et
\vcenter{\xymatrix{
X\oplus Y\ar[r]\ar[d] &\FF\ar[d]\\
X\otimes Y\ar[r] & X\wedge Y
}}$$
are respectively cartesian and cocartesian.
\end{lemma}
\begin{proof}
Direct from the above decomposition of $X\otimes Y$.
\end{proof}

\begin{lemma}\label{pointedhomcart}
Let $X$ and $Y$ be two pointed vector spaces, then the square
$$\xymatrix{
[X,Y]_\bullet \ar[rr]\ar[d]_\epsilon && [X,Y] \ar[d]^{([X,\epsilon_Y],[e_X,Y])}\\
\FF \ar[rr]^-{(\epsilon_X,e_Y)} &&[X,\FF]\times [\FF,Y]
}$$
is cartesian.
\end{lemma}
\begin{proof}
Direct using that $[X,Y]=[\FF,\FF]\oplus [X_-,\FF]\oplus [\FF,Y_-]\oplus [X_-,Y_-]$.
\end{proof}

\medskip

For three pointed vector spaces $X$, $Y$ and $Z$, we shall say that a map $f:X\otimes Y\to Z$ in $\dgVect$ is an {\em expanded map of pointed vector spaces} if it factors (necessarily uniquely) through the projection $X\otimes Y\to X\wedge Y$.

Recall that $X\otimes Y$ is pointed by $(e_X\otimes e_Y,\epsilon_X\otimes \epsilon_Y)$.
A map $F:X\otimes Y\to Z$ is expanded pointed iff 
it is pointed and the following diagram commute
$$\xymatrix{
X\oplus Y\ar[rr]^-{\epsilon_X+\epsilon_Y}\ar[d]_{X\otimes e_Y\oplus e_X\otimes Y} &&\FF\ar[d]^{e_Z}\\
X\otimes Y\ar[rr]^-f && Z.
}$$
This last condition is also equivalent to the commutation of the following diagram
$$\xymatrix{
\FF\ar[rr]^-{e_X} \ar[d]_{e_Z}&& X\ar[d]^{\lambda^2 f} \ar[rr]^-{\epsilon_X} && \FF\ar[d]^{e_Z}\\
Z\ar[rr]^-{[\epsilon_Y,Z]} && [Y,Z] \ar[rr]^-{[e_Y,Z]} && Z.
}$$
In terms of elements, $f:X\otimes Y\to Z$ is an expanded map of pointed vector spaces iff 
$$
\epsilon_Z(f(x,y)) = \epsilon_X(x)\epsilon_Y(y)\ , \qquad
f(e_X,e_Y) = e_Z\ ,
$$
$$
f(e_X,y) = \epsilon_Y(y)e_Z
\et 
f(x,e_Y) = \epsilon_X(x)e_Z
$$
for every $x\in X$ and $y\in Y$.
Notice that the condition $f(e_X,e_Y) = e_Z$ is implied by either of the last two.

\subsection{dg-Algebras}

\subsubsection{Categories of algebras}

We shall be concerned in this work with several types of associative algebras in $\gVect$: unital or not, augmented or not, differential or not.
We will define and study mainly differential graded unital associative algebras, which are the most structured, all results concerning non-unital and non-differential algebras can then be deduced by forgetting the extra structure. As for augmented algebras, we shall use the classical fact that their theory is equivalent to that of non-unital algebras.

\bigskip

A {\it differential graded unital associative algebra} (or a {\em dg-algebra} for short)  is a monoid object in the monoidal category $\dgVect$.
Explicitely it is 
a graded vector space $A$ equipped with 
a linear map $m:A\otimes A \to A$, called the {\it multiplication} or the {\it product}, 
a linear map $e:\FF\to A$, called the {\it unit}
and a graded morphism $d:A\rhup A$ of degree $-1$, called the {\it differential}
satisfying the following conditions:
\begin{itemize}
\item {\it associativity}
$$\xymatrix{
 A\otimes A\otimes A  \ar[rr]^-{m\otimes A} \ar[d]_{A\otimes m} && A\otimes A \ar[d]^-{m} \\
A\otimes A\ar[rr]^-{m} && A
}$$
\item {\it left and right unit}
$$\xymatrix{
A \ar[r]^-{e \otimes A}  \ar@{=}[dr] & A\otimes A \ar[d]^-{m} \\
 & A
}\qquad
\xymatrix{
A \ar[r]^-{A\otimes e}  \ar@{=}[dr] & A\otimes A \ar[d]^-{m} \\
 & A
}$$
\item {\it derivation}
$$\xymatrix{
 A\otimes A  \ar@^{>}[rr]^-{d\otimes A+A\otimes d} \ar[d]_{m} && A\otimes A \ar[d]^-{m} \\
A\ar@^{>}[rr]^-{d} && A.
}$$
\end{itemize}

A {\em graded unital associative algebra} (we shall say a {\em graded algebra}) has only an associative multiplication and a unit, hence it is a monoid object in $\gVect$.
A {\em non-unital dg-algebra} has only an associative multiplication and a differential.
The term {\em algebra} will be used as a synonym of to refer indifferently to any of the previous types of algebras.

We say that a graded algebra $A$ is {\it graded finite} (resp. finite) if its underlying graded vector space is graded finite (resp. finite).

\medskip
 
A {\it map} of differential graded algebras $f:A\to B$ is a homomorphism of monoids in $\dgVect$.
Explicitely, it is a map of graded vector spaces $A\to B$ such that the following diagrams commute,
$$\xymatrix{
A\otimes A \ar[d]_{m_A}  \ar[r]^-{f\otimes f} & \ar[d]^{m_B}B\otimes B \\
A   \ar[r]^-{f} &  B
}
\quad\quad\quad
\xymatrix{
\FF   \ar@{=}[r]  \ar[d]_{e_A}  & \FF \ar[d]^{e_B}  \\
A \ar[r]^f & B
}
\quad\quad\quad
\xymatrix{
A \ar@^{>}[d]_{d_A}  \ar[r]^-{f} & B \ar@^{>}[d]^{d_B} \\
A \ar[r]^-{f} &  B.
}$$
In particular, it is a map of dg-vector spaces.

A {\it map} of graded algebras $f:A\to B$ is a map of graded vector spaces $A\to B$ commuting with the multiplication and the unit only.

A {\it map} of non-unital dg-algebras $f:A\to B$ is a map of graded vector spaces $A\to B$ commuting with the multiplication and the differential only.

\medskip
We shall denote the category of differential graded algebras by $\dgAlg$, that of graded algebras by $\gAlg$ and that of non-unital dg-algebras by $\dgAlg_\circ$.

\begin{rem}
The unit element $e$ of an algebra (differential or not) $A=(A,m,e)$ is uniquely determined by the product $m:A\otimes A\to A$ when it exists. Hence the algebra $(A,m,e)$ is entirely described by the non-unital algebra $(A,m)$. However maps of non-unital algebras may not preserve the units when they exist.
\end{rem}

\medskip
We say that a unital algebra (differential or not) $(A,m,e)$ is {\em augmented} or {\em pointed} if it it is equipped with an algebra map $\epsilon:A\to \FF$, called the {\em augmentation}, such that $\epsilon(1)=1$, \ie such that $\epsilon e=id_\FF$.
A map of augmented algebras is an algebra map $f:A\to B$ commuting with the augmentations
$$\xymatrix{
A   \ar[r]^f  \ar[d]_{\epsilon_A}  & B \ar[d]^{\epsilon_B}  \\
\FF \ar@{=}[r] & \FF.
}$$
The category of augmented dg-algebras is noted $\dgAlg_\bullet$. By definition, we have $\dgAlg_\bullet = \dgAlg/\FF$.

\subsubsection{Examples}

\begin{ex}[Tensor algebra]\label{tensoralg}
The forgetful functor $U:\dgAlg\to \dgVect$ has a left adjoint which associates to a dg-vector space $(X,d)$ the tensor algebra
$$
T(X)=\bigoplus_{n\geq 0} X^{\otimes n}
$$
with differential 
$$
d(x_1\cdots x_n) = \sum_i x_1\cdots (dx_i)\cdots x_n.
$$
We shall say that an element $x=x_1\cdots x_n\in X^{\otimes n}$ has  {\it length} $n$ and {\it degree} $|x|=|x_1|+\cdots +|x_n|$.
If $E$ is a graded set, we shall write $T(E):=T(\FF E)$ and call it the {\em polynomial algebra} with coefficient in $\FF$ in the set of variables $E$.
If $A$ is a dg-algebra, we define the polynomial algebra $A(E)$ as $A\otimes T(E)$.
If $E=\{x\}$, we shall write $T(x)$ or $\FF[x]$ for $T(\{x\})$ and $A[x]$ for $A(\{x\})$.
\end{ex}

\begin{ex}[Pointed tensor algebra]\label{tensoralgaugmented}
The forgetful functor $U:\dgAlg_\bullet\to \dgVect$ has a left adjoint which associates to a dg-vector space $(X,d)$ the tensor algebra
$$
T_\bullet(X)=\bigoplus_{n\geq 0} X^{\otimes n}
$$
with differential 
$$
d(x_1\cdots x_n) = \sum_i x_1\cdots (dx_i)\cdots x_n
$$
and augmentation the projection $\epsilon:T_\bullet \to \FF$ to the $n=0$ factor.

As in the non-augmented case we can define $T_\bullet(E)$ and for any graded set $E$.
\end{ex}

\begin{ex}[Non-unital tensor algebra]\label{tensoralgnu}
The forgetful functor $U:\dgAlg_\circ\to \dgVect$ has a left adjoint which associates to a dg-vector space $(X,d)$ the non-unital tensor algebra
$$
T_\circ(X)=\bigoplus_{n> 0} X^{\otimes n}
$$
with differential 
$$
d(x_1\cdots x_n) = \sum_i x_1\cdots (dx_i)\cdots x_n.
$$
As before we can define $T_\circ(E)$ and for any graded set $E$.
\end{ex}

\begin{ex}[Square zero algebra]\label{squ0alg} If $(X,d)$ is a dg-vector space, then the dg-vector space $T_1(X)=\FF\oplus X$ has the structure of a dg-algebra if the product $T_1(X)\otimes  T_1(X)\to T_1(X)$ is defined by putting 
$$
( \lambda +x)(\lambda' +x')=\lambda\lambda'+\lambda x'+\lambda' x
$$
and the differential is $d( \lambda +x) = dx$.
We shall say that $T_1(X)$ is the {\it square 0 extension} of $\FF$ by $X$.

The square zero algebra is naturally a pointed dg-algebra: the augmentation is given by $p_1:\FF\oplus X\to \FF$.

The square zero algebra has a non-unital analog, it is the vector space $X$ with product given by the zero map $X\otimes X\to X$.

\end{ex}

\begin{ex}[Divided powers algebra]\label{dividedpower}
If $A$ is an algebra, then a finite linear combination 
$$
p(x)=\sum_{n} f_n \Dfrac{x^n}{n!}
$$
of divided powers $\dfrac{x^n}{n!}$ with coefficients $f_n\in A$ is called a {\it divided powers polynomial} with coefficients in $A$.
The divided powers $\dfrac{x^n}{n!}$ are symbols of degree 0 for every $n\geq 0$ (we use a double fraction to recall that they are symbols rather than actual fractions). The product of two divided powers polynomials is the $A$-linear extension of the formulas
$$
\Dfrac{x^n}{n!}\ \ \Dfrac{x^m}{m!}\ \ = {n+m\choose m}\Dfrac{x^{m+n}}{(m+n)!}.
$$
Notice that $A\{x\}=A\otimes \FF\{x\}$.

\end{ex} 

\begin{ex}[Formal power series algebra] \label{formalpower}

Let $A$ be a dg-algebra and $x$ be a variable of degree $d$, the {\em formal power serie algebra}  in the variable $x$ with coefficients in $A$ is defined as the dg-vector space $\FF[[x]]$ such that an element $f\in \FF[[x]]_n$ is a formal power series 
$$
f=f(x)=\sum_{i\geq 0} f_i x^i
$$
with $f_i\in A_{n-di}$ for every $i\geq 0$.

For more on this example see example \ref{convolutionformalpowerseries}.
\end{ex} 

\begin{ex}[Formal divided power series algebra]\label{Formal divided power series algebra}
Let $A$ be a dg-algebra and $x$ be a variable of degree $0$, the {\em formal divided power serie algebra} in the variable $x$ with coefficients in $A$ is defined as the dg-vector space $\FF\{\{x\}\}$ such that 
an element $f\in \FF\{\{x\}\}_n$ is a formal power series 
$$
f=f(x)=\sum_{i\geq 0} f_i \Dfrac{x^i}{i!}
$$
with $f_i\in A_n$ for every $i\geq 0$ and
where the $\dfrac{x^i}{i!}$ are the same formal symbols as in the example \ref{dividedpower}.

For more on this example see example \ref{convolutionformaldividedpowerseries}.
\end{ex}

\begin{ex}[Commutative algebras]
Recall that a dg-algebra $A$ is said to be {\it graded commutative}
if we have $ab=ba(-1)^{|a||b|}$ for every $a,b\in A$;
in which case we have $2a^2=0$ when $|a|$ is odd.
The full subcategory of $\dgAlg$ generated by commutative dg-algebras is noted $\dgAlgcom$.
The inclusion functor $\dgAlgcom\to \dgAlg$ has a left adjoint sending $A$ to $A_{ab}:=A/[A,A]$ where $[A,A]$ is the bilateral ideal generated by commutators $[a,b]:=ab-ba(-1)^{|a|||b}$.
\end{ex}

\subsubsection{Comparison functors}\label{comparisonfunctoralgebras}

There exists several functors between the categories of algebras.
For a dg-algebra, we can forget the differential
$$\xymatrix{
U_d: \dgAlg \ar[r] & \gAlg
}$$
or the unit
$$\xymatrix{
U_e: \dgAlg \ar[r] & \dgAlg_\circ.
}$$
And for pointed dg-algebra, we can forget the augmentation
$$\xymatrix{
U_\epsilon:\dgAlg_\bullet \ar[r] & \dgAlg
}$$
For a pointed dg-algebra $A$, the kernel of the augmentation $A_-:=\ker \epsilon$ is naturally a non-unital algebra.
This defines a functor
$$\xymatrix{
(-)_-:\dgAlg_\bullet \ar[r] & \dgAlg_\circ
}$$

Reciprocally, if $A$ is a non-unital algebra, $A_+:=\FF\oplus A$ can be equipped with a unital algebra structure by putting 
$$
(\alpha,a)(\beta,b)=(\alpha \beta,\alpha b+\beta a+ m(a,b))
$$
and the unit element $e=(1,0)$. Moreover $A_+$ is naturally pointed by the projection to $\FF$.
This defines a functor
$$\xymatrix{
(-)_+:\dgAlg_\circ \ar[r] & \dgAlg_\bullet
}$$

\begin{prop}\label{equivnonunitalpointedalg}
The two functors
$$\xymatrix{
\dgAlg_\bullet \ar@<.6ex>[rr]^-{(-)_-} && \dgAlg_\circ \ar@<.6ex>[ll]^-{(-)_+}
}$$
are inverse equivalences of categories.
\end{prop}
\begin{proof}
A direct computation proves that if $A$ is a non-unital algebra and $B$ a pointed dg-algebra, 
then $(A_+)_-$ is naturally isomorphic to $A$ and $(B_-)_+$ is naturally isomorphic to $B$.
\end{proof}

\begin{cor}\label{adjunctionunitalnonunital}
\begin{enumerate}
\item The composition $U_\epsilon(-)_+:\dgAlg_\circ\to \dgAlg_\bullet \to \dgAlg$ is left adjoint to $U_e$ 
$$\xymatrix{
U_\epsilon(-)_+ :\dgAlg_\circ \ar@<.6ex>[rr] && \dgAlg \ar@<.6ex>[ll]{:U_e}.
}$$
The unit of the adjunction is the inclusion $A\to \FF \oplus A$.

\item The composition $(U_e-)_+:\dgAlg\to \dgAlg_\circ \to \dgAlg_\bullet$ is right adjoint to $U_\epsilon$ 
$$\xymatrix{
U_\epsilon :\dgAlg_\bullet \ar@<.6ex>[rr] && \dgAlg \ar@<.6ex>[ll]{:(U_e-)_+}.
}$$
The counit of the adjunction is the projection $\FF \oplus A\to A$.
\end{enumerate}
\end{cor}

\begin{rem}
We can summarize this corollary by saying that, up to the equivalence $\dgAlg_\bullet\simeq \dgAlg_\circ$, the functor $U_\epsilon$ is left adjoint to $U_e$.
\end{rem}

\bigskip
We shall prove in theorem \ref{dgtoghomalg} that the forgetful functor $U_d:\dgAlg\to \gAlg$ has both a left and a right adjoint.
With the consequence that limits and colimits exist in $\dgAlg$ and they can be computed in $\gAlg$.

\subsubsection{Opposite algebra and anti-homomorphisms}\label{anti-homo-def}

If $A$ and $B$ are graded algebras, we shall say that a 
linear map $f:A\to B$ is an {\it anti-homomorphism}, if $f(1)=1$ and $f(xy)=f(y)f(x)(-1)^{|x||y|}$
for every $x,y\in A$.
The {opposite} of a graded algebra $A$
is a graded algebra $A^o$ anti-isomorphic to $A$. 
By definition, the map $x\mapsto x^o$ is a bijection 
$A\to A^o$ and for every $a,b\in A$ we have
$a^ob^o=(ba)^o (-1)^{|a||b|}.$
There is an obvious bijection between the anti-homomorphisms
of algebras $A\to B$ and the homomorphisms
of algebras $A^o\to B$ and $A\to B^o$.

\medskip

\begin{ex} The inclusion $i:V\to T(V)$ can be extended uniquely as an anti-automorphism
of the tensor algebra $(-)^o:T(V)\to T(V)$. By definition,
$$(x_1\otimes \cdots \otimes x_n)^{o}=x_n \otimes \cdots \otimes x_1\ (-1)^{\sum_{i<j}|x_i||x_j|}.$$
Thus, $T(V)^o=T(V)$.
\end{ex}

\begin{ex} If $X$ is a graded vector space, then the graded vector space $End(X)=[X,X]$
has the structure of a graded algebra in which the product  is the composition 
of graded morphisms.
The transposition operation
$${}^t(-):End(X)\to End(X^\star)$$
is an anti-homomorphism of algebras.
The algebra $End(X)$ is finite when $X$ is finite, in which case
the transposition operation is invertible.
\end{ex}

\subsubsection{Monoidal structures}\label{monoidalstructurealgebras}

The tensor product $A\otimes B$ of two algebras has the structure of an algebra with multiplication define as the composite of the maps
$$\xymatrix{
A\otimes B\otimes A\otimes B \ar[rr]^-{A\otimes \sigma \otimes B}  &&  A\otimes A\otimes B\otimes B \ar[rr]^-{m_A \otimes m_B} & & A\otimes B.
}$$
Explicitely, this gives
$$
(a\otimes b) (a'\otimes b')=(a a')\otimes (b b')\ (-1)^{|b||a'|},
$$
for $a,a'\in A$ and $b,b'\in B$.
If $A$ and $B$ are unital, so is $A\otimes B$ with unit defined by $e_A\otimes e_B:\FF\simeq \FF\otimes \FF \to A\otimes B$.
If $A$ and $B$ are differential, so is $A\otimes B$ with differential defined by $d_{A\otimes B}=d_A\otimes B+A\otimes d_B$, as for dg-vector spaces.

\medskip
\begin{lemma}\label{tensorfromsum}
Let $A\cup B$ be the sum of $A$ and $B$ as algebras, then 
$A\otimes B = A\cup B/(ab-ba (-1)^{|a||b|} ; a\in A,b\in B)$.
\end{lemma}
\begin{proof}
A straightforward computation proves that $A\cup B = T(A\oplus B)/(a\otimes a' = aa', b\otimes b' = bb')$.
The map $A\otimes 1_B: A\to A\otimes B$ and $1_A\otimes B: B\to A\otimes B$ induce an algebra map $A\cup B\to A\otimes B$
which send an element $a_1b_1\dots a_nb_n$ to $\pm (a_1\dots a_n)\otimes (b_1\dots b_n)$.
Hence the map factors through the quotient $A\cup B/(ab-ba (-1)^{|a||b|} ; a\in A,b\in B)$ and it is easy to see that it becomes an isomorphism.
\end{proof}

\medskip
The category $\dgAlg$ is symmetric monoidal, the unit object is the algebra $\FF$. 
The category $\gAlg$ is symmetric monoidal, the unit object is the algebra $\FF$. 
The category $\dgAlg_\circ$ is symmetric monoidal, the unit object is the algebra $\FF e$ where $e$ is an idempotent. 
The three forgetful functors $\dgAlg\to \dgVect$, $\gAlg\to \gVect$ and $\dgAlg_\circ\to \dgVect$ are clearly symmetric monoidal.

\medskip
Using the equivalence of proposition \ref{equivnonunitalpointedalg}, we can transport the monoidal structure of $\dgAlg_\circ$ to a monoidal $\dgAlg_\bullet$. We define the {\it smash product} $A\wedge B$ of two pointed algebras as the smash product of the underlying pointed vector spaces
$$
A\wedge B = ( A_{-} \otimes  B_{-})_{+}.
$$
The unit object for the smash product is the algebra $\FF_{+}=\FF\oplus \FF e$ generated by one idempotent $e$;
the algebra  $\FF_{+}$ is pointed by the first projection $ \FF\oplus \FF e\to  \FF$.
Notice the isomorphism $\FF_{+}\simeq \FF e\oplus  \FF(1-e)\simeq \FF\times \FF$.
Hence the algebra $\FF_{+}$ is isomorphic to the cartesian product of the algebra $\FF$ with itself.

$(\dgAlg_\bullet,\wedge , \FF_{+})$ is a symmetric monoidal category and the forgetful functor $\dgAlg_\bullet\to \dgVect_\bullet$ is monoidal.

\medskip
The pointed algebra $A\wedge B$ has also the following description.

\begin{prop}\label{smashproductofalgebras}
We have a pullback square of algebras
$$\xymatrix{
A\wedge B \ar[rr] \ar[d]_\epsilon &&  \ar[d]^{(A\otimes  \epsilon_B,\epsilon_A\otimes B) }  A\otimes B \\
\FF \ar[rr]&& A \times  B.
}$$
\end{prop}

\begin{proof} 
Because $A\times B = A\oplus B$ as vector spaces, we know the square is cartesian by lemma \ref{pointedbicart}.
We need only to prove that the maps are algebra maps. This is a straightforward computation left to the reader.
\end{proof}

\begin{prop}\label{monoidalfunctorsalg}
\begin{enumerate}
\item The forgetful functor
$$\xymatrix{
U_d: \dgAlg \ar[r] & \gAlg
}$$
is monoidal. 

\item The forgetful functor
$$\xymatrix{
U_e: \dgAlg \ar[r] & \dgAlg_\circ.
}$$
is monoidal. Hence its left adjoint $U_\epsilon(-)_+$ is colax monoidal.

\item The forgetful functor
$$\xymatrix{
U_\epsilon:(\dgAlg_\bullet,\wedge,\FF_+) \ar[r] & (\dgAlg,\otimes,\FF)
}$$
is colax monoidal and its right adjoint is monoidal.
The colax structure is given by the map $A\wedge B\to A\otimes B$ from proposition \ref{smashproductofalgebras} and by the map $\FF_+\simeq \FF\oplus \FF e\to \FF$ sending $e$ to $1$.
\end{enumerate}
\end{prop}
\begin{proof}
The monoidal structure of 1. and 2. are obvious. 
The consequences on the colax structure on the left adjoint is a general fact about adjoint to monoidal functors.
Then 3. derives from 2. and the fact that up to the equivalence $\dgAlg_\bullet \simeq \dgAlg_\circ$ the functor $U_\epsilon$ corresponds to the left adjoint to $U_e$.
\end{proof}

\subsubsection{Generation and separation}

We introduce two definitions that will useful in some proofs later.

\begin{defi}\label{separation}
Let $A$ be an algebra and $X$ a vector space
\begin{enumerate}
\item we shall say that a linear map $f:X\to A$ is {\it generating} if the associated algebra map $g:T(X)\to A$ is surjective.
\item we shall say that a linear map $f:X\to A$ is {\it separating} if the implication $uf=vf\Rightarrow u=v$ is true for any pair of algebra maps $u,v:A\to B$. 
\end{enumerate}
\end{defi}

Every linear map $f:X\to A$ can be extended as an algebra map $g:T(X)\to A$,
the map $f:X\to A$ is then separating if and only if the map $g:T(X)\to A$ is a epimorphism of algebras.
In particular, any generating map is separating.
For example, the free map $p:V\to T(X)$ is generating, hence separating.

\subsubsection{Modules}

We shall only develop the theory of unital modules, that of non-unital modules is analog and all results are valid also in this context.

\bigskip

Recall that if $A$ is a dg-algebra, then a {\it left $A$-module} is a dg-vector space $X$ equipped with a left action $\alpha:A\otimes X\to X$
which is associative and unital,
$$\xymatrix{
A\otimes A\otimes X  \ar[rr]^-{m\otimes X} \ar[d]_{A\otimes \alpha} && A\otimes X \ar[d]^-{\lambda} \\
A\otimes X\ar[rr]^-{\alpha} && X
}
\qquad
\xymatrix{
X \ar[r]^-{e \otimes X}  \ar@{=}[dr] & A\otimes X \ar[d]^-{\lambda} \\
& X
}$$
Let us put $ax=\alpha(a \otimes x)$. The conditions means that we have $a(bx)=(ab)x$ and $ex=x$ for every $a,b\in A$ and $x\in X$.
If $X$ and $Y$ are (left) $A$-modules, then a {\it module map} $f:X\to Y$ is linear map such that $f(ax)=af(x)$ for every (homogenous) $a\in A$ and $x\in X$. 
We shall denote by ${\Mod}(A)$ the category whose objects are the left $A$-modules and whose arrows are the module maps.

\medskip

The category ${\Mod}(A)$ is also enriched over the category $\dgVect$. 
If $X$ and $Y$ are two $A$-modules, the dg-vector space of morphisms of $A$-modules between them is the object $\Hom_A(X,Y)$ defined by the equalizer
$$\xymatrix{
\Hom_A(X,Y)\ar[r]& [X,Y] \ar@<.6ex>[rr]^-{[a_X,Y]}\ar@<-.6ex>[rr]_-{[A\otimes -, a_Y]} && [A\otimes X,Y]
}$$
where $a_X:A\otimes X\to X$ and $a_Y:A\otimes Y\to Y$ are the $A$-module structures.
Concretly, a graded morphism $f\in [X,Y]_n$ belongs to $\Hom_{A}(X,Y)_n$ iff we have 
$$
f(ax)=af(x)\ (-1)^{n|a|}
$$
for every $a\in A$ and $x\in X$.

\begin{prop} \label{enrichmentofmodules}
The enriched category ${\Mod}(A)$ admits tensor and cotensor products over $\dgVect$.
The tensor product of an $A$-module $X$ by a dg-vector space $V$ is the dg-vector space $V\otimes X$ equipped with the $A$-module structure defined by putting $a(v\otimes x)=(-1)^{|a||v|}v\otimes ax$ for $a\in A$, $v\in V$ and $x\in X$.
\end{prop}

\begin{proof} 
A straightforward computation shows the bijections between
\begin{center}
\begin{tabular}{lc}
\rule[-2ex]{0pt}{4ex}  linear maps & $V\to \Hom_A(X,Y)$, \\
\rule[-2ex]{0pt}{4ex}  $A$-module maps & $V\otimes X\to Y$, \\
\rule[-2ex]{0pt}{4ex}  and $A$-module maps & $X\to [V,Y]$.
\end{tabular}
\end{center}
\end{proof}

In particular, the $n$-fold supension $S^n(X)=S^n\otimes X$ of an $A$-module $X$ has the structure of an $A$-module if we put $a(s^nx)=(-1)^{|a|n}s^nax$ for $a\in A$ and $x\in X$.

\bigskip

The functor $A\otimes(-):\dgVect \to \dgVect$ has the structure of a monad with multiplication $m\otimes X:A\otimes A\otimes X \to A\otimes X$
and unit $e\otimes X:X \to A\otimes X$. An algebra for this monad is a (left) $A$-module.
Dually, the functor $[A,-]:\dgVect\to \dgVect$ has the structure of a comonad, since it is right adjoint to the functor $A\otimes(-)$; its comultiplication is given by the map $[m,X]:[A,X] \to [A\otimes A,X]=[A,[A,X]]$ and its counit by the map $[e, X]:[A,X]\to [\FF,X]=X$. 
If $\alpha:A\otimes X \to X$ is an algebra for the monad $A\otimes (-)$, then $\lambda^1(\alpha):X\to [A,X]$ is a coalgebra for the comonad $[A,-]$.

\begin{prop}\label{moduleandmonad}
If $X$ is a graded vector space, then the following data are equivalent:
\begin{enumerate}
\item An action $a:A\otimes X\to X$ of the algebra $A$ on $X$;
\item  An action $A\otimes X\to X$ of the monad $A\otimes(-)$ on $X$;
\item  A coaction $X\to [A,X]$ of the comonad $[A,-]$ on $X$.
\item  A representation $\pi:A\to [X,X]$ of the algebra $A$;
\end{enumerate}
\end{prop}

\begin{proof} Left to the reader.
\end{proof}
 
\begin{prop} \label{modulesoveran algebra}
The forgetful functor ${\Mod}(A) \to\dgVect$  is faithful and admits a left adjoint $X\to A\otimes X$ and a right adjoint $X\mapsto [A,X]$.
Limits and colimits exist in ${\Mod}(A)$ and they can be computed in $\dgVect$.
\end{prop}

\begin{proof} This follows from the general theory of algebras over a monad and of coalgebras over a comonad.
\end{proof}

A $A$-module $X$ is said to be {\em freely generated} by a map $i:V\to X$ if the corresponding map $A\otimes V\to X$ is an isomorphism of $A$-modules. Through this isomorphism, the map $i$ identifies with the map $e \otimes V:V\to A\otimes V$.
A module $X$ is said to be {\em free} if there exists a map $i:V\to X$ which freely generates it.

\bigskip

There is a dual notion of a {\it right} $A$-module
defined by a right action, $\beta:X\otimes A\to X$
$$\xymatrix{
 X\otimes A\otimes A  \ar[rr]^-{X\otimes m} \ar[d]_{\beta\otimes A} && X\otimes A \ar[d]^-{\beta} \\
X\otimes A\ar[rr]^-{\beta} && X
}\qquad \xymatrix{
X \ar[r]^-{X \otimes e}  \ar@{=}[dr] & X\otimes A \ar[d]^-{\beta} \\
 & X
}$$
Let us put $xa=\beta(x\otimes a)$. The conditions means 
that we have $(xa)b=x(ab)$ and $xe=x$ for every $a,b\in A$ and $x\in X$.
A right $A$-module $X$ is the same thing as a left $A^o$-module
if we put $a^ox=xa(-1)^{|a||x|}$ for every $a\in A$ and $x\in X$.
Because of this we define the category of right $A$-module to be $\Mod(A^o)$.

If $A$ and $B$ are graded algebras, then a left $A$-module structure $\alpha:A\otimes X\to X$
is said to {\it commute} with a right $B$-module structure $\beta:X\otimes B\to X$ 
if following square commutes,
$$\xymatrix{
 A\otimes X\otimes B  \ar[rr]^-{\alpha\otimes B} \ar[d]_{A\otimes \beta} && X\otimes B \ar[d]^-{\beta} \\
A\otimes X\ar[rr]^-{\alpha} && X.
}$$
In the notation introduced above, this means that we have
$(ax)b=a(xb)$  for every $a\in A$, $x\in X$ and $b\in B$.
In which case $X$ is said to be an $(A,B)$-{\it bimodule}.
An $(A,B)$-bimodule is the same thing
as a left module over $A\otimes B^{o}$, or as a right module over $A^o\otimes B$.

Morphisms of bimodules are equivalently defined as morphisms of left modules over $A\otimes B^{o}$ or of right modules over $A^o\otimes B$. The category of $(A,B)$-bimodules is noted $\Bimod(A,B)$. We have an equivalence $\Bimod(A,B)=\Mod(A\otimes B^{o})$.
If $A=B$, we shall put $\Bimod(A):=\Bimod(A,A)$.

As a particular case of a module category, $\Bimod(A,B)$ is enriched, tensored and cotensored over $\dgVect$.

\bigskip
If $f:A\to B$ is an algebra map, then $B$ has a canonical structure of a $A$-bicomodule.
given by $a\cdot b=f(a)b$ and $b\cdot a = bf(a)$. This produces a functor $A\bs \dgAlg\to \Bimod(A)$. 
If $X$ is a $A$-bimodule, the tensor $A$-algebra of $X$ is defined as 
$$
T_A(X) : =  \bigoplus_{n\geq 0} X^{\otimes_A n}= A\oplus X\oplus X\otimes_AX \oplus X\otimes_AX\otimes_AX \oplus \dots 
$$

\begin{prop}\label{tensoralgebraforbimodules}
The functor $T_A$ is left adjoint to the functor $A\bs \dgAlg \to \Bimod(A)$.
\end{prop}
\begin{proof}
This can be checked directly on elements but we are going to give a more conceptual proof to generlize the statement to coalgebras in proposition \ref{cofreecoalgebraforbicomodules}.
Let $X$ be a $A$-bicomodule and let $\lambda_X:A\otimes X\to X$ and $\rho_X:X\otimes A\to X$ be the left and right action of $A$ on $X$, we defined $\mu =(m_A,\lambda_X,\rho_X)p:(A\otimes X)^2\to A\otimes A \oplus A\otimes X\oplus X\otimes A \to A\oplus X$ 
where the map $p$ is the canonical projection.
Then, $T_A(X)$ can be described as the coequalizer in $\dgAlg$ of
$$\xymatrix{
T((A\oplus X)^2) \ar[rr]^-{T(\mu)} \ar[rd]_{\alpha}&& T(A\oplus X)\\
&T(A\oplus X)\otimes T(A\oplus X) \ar[ru]_-{m_T}
}$$
where $\alpha:T(X\otimes X)\to T(X)\otimes T(X)$ is the colax structure of $T$ given by $\alpha(\prod (x_i\otimes y_i)) = \pm (\prod x_i )\otimes (\prod y_i )$ (the sign depends on the degrees of the $x_i$ and $y_i$)
We leave to the reader the proof that a cocone from this diagram to $B$ is equivalent to the data of an algebra map $f:A\to B$ and an $A$-bimodule map $g:X\to B$.
\end{proof}

\bigskip

The tensor product of a right $A$-module $X$ and a left $A$-module $Y$ is the object $X\otimes_AY$ defined as the coequalizer
$$\xymatrix{
X\otimes A\otimes Y \ar@<.6ex>[rr]^-{X\otimes a_Y}\ar@<-.6ex>[rr]_-{a_X\otimes Y}&& X\otimes Y \ar[r]& X\otimes_AY
}$$
where $a_X:A\otimes X\to X$ and $a_Y:Y\otimes A\to Y$ are the $A$-module structures.

Let $f:A\to B$ be a dg-algebra map, any $B$-module $X$ can be seen as a $A$-module by the formula $a\cdot x = f(a)x$ for any $a\in A$ and $x\in X$. This produces a functor $\Mod(B)\to \Mod(A)$ called the {\em restriction functor}. In particular $B$ can be viewed as a $A$-module.

Let $X$ be a left $A$-module, then the dg-vector spaces $B\otimes_AX$ and $\Hom_A(B,X)$ are $B$-modules for the actions
$b'\cdot (b\otimes x) = (b'b)\otimes x$ and $(b\cdot f)(x) = f(bx) \ (-1)^{|f||b|}$.

\begin{prop}
The functor $B\otimes_A-$  and $\Hom_A(B,-) : \Mod(A)\to \Mod(B)$
are respectively left and right adjoint to the restriction functor.
\end{prop}
\begin{proof}
A straightforward computation left to the reader
\end{proof}

\subsubsection{Derivations}\label{derivationalgebra}

\begin{defi} \label{defderivmodule}
Let $A$ be a dg-algebra or a non-unital dg-algebra and $M$ an $A$-bimodule.
We shall say that a map of dg-vector spaces $D:A\to M$ is a {\it derivation of $A$ with values in $M$} (or simply a {\em derivation})
if we have the Leibniz rule
$$
D(ab)=D(a)b+aD(b)
$$
for every $a,b\in A$.
We shall note $der(A,M)$ the set of derivations of $A$ with values in $M$.

\medskip
More generally, we shall say that a graded morphism of degree $n$ of graded vector spaces $D:A\rhup M$ is a {\it graded derivation of $A$ with values in $M$ of degree $n$} (or simply a {\em graded derivation})
if we have the graded Leibniz rule
$$
D(ab)=D(a)b+aD(b)\ (-1)^{n|a|}
$$
for every $a,b\in A$.
We shall note $\Der(A,M)_n$ the vector space of graded derivations of $A$ with values in $M$ of degree $n$.
We have an inclusion $\Der(A,M)_n\subset [A,M]_n$ and the graded vector subspace $\Der(A,M)=\Der(A,M)_*\subset [A,M]$ is stable by the differential of $[A,M]$ which enhance it into a dg-vector space. The differential of a graded derivation $D:A\rhup M$ of degree $n$ is the map  $d_MD-Dd_A \ (-1)^{n}$ of degree $n-1$.
Derivations correspond exactly to graded derivations of degree 0 such that the morphism $D:A\to M$ commute with the differentials, equivalently they correspond to cycles in $\Der(A,M)_0$.
\end{defi}

\begin{rem}\label{remdersuspension}
Let $S^n$ the graded vector space freely generated by one element $s^n$ of degree $n$.
A morphism $D:A\rhup M$ of degree $n$ is a graded derivation iff the map $s^{-n}D:A\to S^{-n}M=S^{-n}\otimes M$ defined by putting $s^{-n}D(x)=s^{-n}\otimes D(x)$ is a graded derivation (of degree 0) with values in the bimodule $S^{-n}M$.
\end{rem}

\begin{ex}\label{innerderivationsofsemialg}
The {\it commutator} $[a,b]$ of two elements in a non-unital algebra $A$
is defined by putting
$[a,b]=ab-ba(-1)^{|a||b|}.$
The morphism $[a,-]:A\rhup A$ is a derivation of degree $|a|$  for every $a\in A$
since we have
$$
[a,x]y+x[a,y](-1)^{|a||x|}=axy-xay(-1)^{|a||x|}+xay^{|a||x|}-xya^{|a||x|+|a||y|}=[a,xy]
$$
\end{ex}

\begin{lemma}\label{smalllemmaforderivation}
If $D:A\rhup M$ is a graded derivation and $e:\FF\to A$ is the unit map of $A$, then $De=0$.
\end{lemma}

\begin{proof}
This is equivalent to $D(1)=0$ but $D(1)=D(1\cdot 1)=D(1)\cdot 1+1\cdot D(1)=D(1)+D(1)$.
\end{proof}

\begin{rem}
If $X$ be a graded vector space equipped with a set of operations $\Phi_i:X^{\otimes k_i} \to X$ of arity $k_i$, let us call a graded morphism $D:X\rhup X$ of degree $n$ a {\it derivation} if we have, for all $i$,
$$
D\Phi_i(x_1,\ldots, x_{k_i})=\sum_{j=1}^{k_i} \Phi_i( x_1,\ldots, x_{j-1},D(x_j),x_{j+1}, \ldots, x_{k_i}) 
(-1)^{n|x_1|+\cdots+n|x_{j-1}| }
$$
for every $x_1,\ldots, x_{k_i}\in X$.
This definition applies in particular to non-unital algebras $(A,m)$ and unital algebras $(A,m,e)$. 
Lemma \ref{smalllemmaforderivation} says that the extra condition of the latter case is superfluous.
In consequence, for a unital algebra $A$, we shall use the same notation $\Der(A)$ to refer to derivations of $A$ viewed as unital or not.
\end{rem}

\bigskip

If $D:A\to M$ is a derivation and $M\to N$ is a map of $A$-bimodules, the composite $fD:A\to N$ is still a derivation.
This defines a functor
$$\xymatrix{
der(A,-):\Bimod(A)\ar[r] & \Set.
}$$
We shall say that a derivation $d:A\to \Omega_A$ is universal if it represents the functor $der(A,-)$.
Equivalently, $d:A\to \Omega_A$ is universal if for any derivation $D:A\to M$, there exists a unique bimodule map $f:\Omega_A\to M$ such that $fd=D$
$$\xymatrix{
\Omega_A\ar@{-->}[rrd]^f \\
A\ar[u]^d\ar[rr]^-D&& M
}$$

\begin{prop}\label{univderiv}
For any dg-algebra $A$ (unital or not), there exists a universal derivation $d:A\to \Omega_A$.
\end{prop}

\begin{proof} 
Let $\Omega_A$ is the $A$-bimodule generated by elements $dx$ for every $x\in A$ and relations
$$
d(xy) = (dx)y+ xdy.
$$
There is a canonical map $d:A\to \Omega_A$ sending $x$ to $dx$. If $d':A\to M$ is any derivation, the map $f:\Omega_A\to M$ is defined by sending $dx$ to $d'(x)$. Relations $f(d(xy)) = f(dx)y+ xf(dy)$ are satisfied because $d'$ is a derivation.
\end{proof}

We shall call the bimodule $\Omega_A$ the {\em bimodule of differentials of $A$} or the {\em cotangent space of $A$}.

\begin{rem}\label{gradedunivder}
An analogous functor can be defined in the graded context if $der(A,M)$ is replaced by $\Der(A,M)_0$.
If $A$ is a dg-algebra (unital or not), the same proof shows that, $d:A\to \Omega_A$ viewed as a graded map, is still the universal graded derivation.
In other words, we have a natural isomorphism $\Der(A,M)_0=\Hom_{A,A}(\Omega_A,M)_0$, 
where, for two $A$-bimodules $M$ and $N$, $\Hom_{A,A}(N,M)$ is the enriched dg-vector space hom in $\Bimod(A)$.
\end{rem}

\begin{rem}\label{remderuniv}
In diagram terms, $\Omega_A$ is the cokernel of the following diagram in the category of $A$-bimodules:
$$\xymatrix{
A\otimes A \otimes A\otimes A \ar@<.6ex>[rrrr]^-{A\otimes m_A\otimes A} \ar@<-.6ex>[rrrr]_-{m_A\otimes A\otimes A + A\otimes A\otimes m_A} &&&& A\otimes A\otimes A.
}$$
And the map $d:A\to \Omega_A$ is induced by $1_A\otimes A\otimes 1_A:A\to A\otimes A\otimes A$.

This diagram is part of a full simplicial diagram of $A$-bimodules:
$$\xymatrix{
\dots &A\otimes A \otimes A\otimes A \ar@<1ex>[r] \ar[r] \ar@<-1ex>[r] & A\otimes A\otimes A \ar@<.6ex>[r]\ar@<-.6ex>[r]& A\otimes A \ar[r] & A
}$$
whose associated complex is the Hochschild complex of $A$.
This simplicial diagram admits a contracting homotopy, making the Hochschild complex an exact complex \cite{CE}.
In consequence, $\Omega_A$ can be seen to be isomorphic as an $A$-bimodule to the kernel of $m:A\otimes A\to A$
and the derivation $d:A\to \Omega_A$ is obtained by putting $d(x)=1\otimes x-x\otimes 1$.

The same argument is valid in the non-differential setting.
\end{rem}

We can proceed similarly with graded derivations. They define a functor
$$\xymatrix{
\Der(A,-):\Bimod(A)\ar[r] & \dgVect.
}$$
We shall say that a graded derivation is universal if it represents the functor $\Der(A,-)$.
We have the following stronger form of proposition \ref{univderiv}.

\begin{prop}\label{univderivgraded}
For any dg-algebra $A$ (unital or not), $d:A\to \Omega_A$, viewed as a graded derivation, represents the functor $\Der(A,-)$.
Equivalently, we have a natural isomorphism $\Der(A,M)=\Hom_{A,A}(\Omega_A,M)$. 
\end{prop}
\begin{proof}
Let us show that $\Der(A,M)$ is naturally isomorphic to $\Hom_{A,A}(\Omega_A,M)$. 
Because of remark \ref{gradedunivder}, $d:A\to \Omega_A$ is also the universal derivation for $A$ as a graded algebra.
Hence graded derivations of degree $0$ are in bijection with elements of $\Hom_{A,A}(\Omega_A,M)_0$.
From remark \ref{remdersuspension} we deduce that graded derivations of degree $n$ are in bijection with elements of $\Hom_{A,A}(\Omega_A,S^{-n}M)_0=\Hom_{A,A}(\Omega_A,M)_n$. This construct a natural isomorphism $\Der(A,M)\simeq \Hom_{A,A}(\Omega_A,M)$ as graded vector spaces. Let us see that it commutes to the differential. If $f:\Omega_A\rhup M$ is a graded morphism, the corresponding derivation is $fd:A\to \Omega_A\rhup M$.
The differential of $f$ is $d_Mf-fd_A (-1)^{|f|}$, it is send to $d_Mfd-fd_Ad (-1)^{|f|}$. Because $d:A\to \Omega_A$ commute with the differential, we have $d_Mfd-fd_Ad (-1)^{|f|}=d_Mfd-fdd_A (-1)^{|f|}$ which is the differential of $fd$ in $[A,M]$.
\end{proof}

\begin{rem}\label{remdersuspension2}
The universal derivation $d:A\to \Omega_A$ can be viewed as a graded morphism $d:A\rhup S^n\Omega_A$.
A consequence of the previous result is that $d:A\rhup S^n\Omega_A$ is the universal derivation of degree $n$ in the sense that 
any graded derivation $A\rhup M$ of degree $n$ factor uniquely through a map $S^n\Omega_A\to M$ of $A$-bimodules.
\end{rem}

\bigskip

Let $f:A\to B$ be a map of dg-algebras (unital or not) and $M$ a $B$-bimodule, then $M$ can be viewed as a $A$-bimodule if we put $a\cdot x=f(a)x$ and $x\cdot a=xf(a)$ for $a\in A$ and $x\in M$. We shall note this bimodule $M_f$.

\begin{defi}\label{fderivation} \label{defderivationalg} 

If $f:A\to B$ is a map of dg-algebras  (unital or not) and $M$ a $B$-bimodule, 
we shall say that a morphism $D:A\rhup M$ of degree $n$ is an $f$-{\it graded derivation with values in $M$}  (or a {\em graded $f$-derivation}) if it is a graded derivation of $A$ with values in the bimodule $M_f$. The $f$-derivations form a dg-vector space $\Der(A,f;M):=\Der(A,M_f)$.

In case $M=B$ with the canonical bimodule structure, we shall put $\Der(f):=\Der(A,f;B)$.
In case $f=id_A$ and $M=A$, we shall put $\Der(A):=\Der(id_A)$.
\end{defi}

Concretely, a morphism $D:A\rhup M$ of degree $n$ is a $f$-derivation if we have
$$
D(ab)=D(a) f(b)+f(a) D(b)(-1)^{n|a|}
$$
for every $x,y\in A$.
In particular, a morphism $D:A\rhup A$ of degree $n$ is a derivation if we have
$$
D(ab)=D(a) b+a D(b)(-1)^{n|a|}
$$
for every $a,b\in A$.

\bigskip

$f$-derivations define a functor
$$\xymatrix{
\Der(A,f;-):\Bimod(B)\ar[r] & \dgVect.
}$$
We shall say that a graded $f$-derivation $d:A\rhup M$ is universal if it represents the functor $\Der(A,f;-)$.

Let $\Omega_A$ be the bimodule of differential of $A$, and let $\Omega_{A,f}:=\Omega_A\otimes_{A^o\otimes A}B^o\otimes B$ be the $B$-bimodule obtained by base change along $f$. There is a canonical $A$-bimodule map $\Omega_A\to \Omega_{A,f}$ and an $f$-derivation $A\to \Omega_{A,f}$.
We have the following variation of proposition \ref{univderivgraded}.

\begin{prop}\label{univfderiv}
The $f$-derivation $A\to \Omega_{A,f}$ is universal. 
Equivalently, we have a natural isomorphism $\Der(A,f;M)=\Hom_{B,B}(\Omega_{A,f},M)$.
\end{prop}
\begin{proof}
Let $A\to M$ be an $f$-derivation, by proposition \ref{univderivgraded} we have $\Der(A,f;M)=\Der(A,M_f)=\Hom_{A,A}(\Omega_A,M_f)$
and by base change, we have $\Hom_{A,A}(\Omega_A,M_f) = \Hom_{B,B}(\Omega_A\otimes_{A^o\otimes A}B^o\otimes B,M)$.
\end{proof}

We shall call $\Omega_{A,f}$ the {\em cotangent space of $A$ at $f$}.

\bigskip

If $A$ is a dg-algebra (unital or not) and $M$ is a bimodule, then the direct sum $A\oplus M$ has the structure of an dg-algebra with the product defined by putting 
$$
(a,x)(a',x') = \left(aa', ax' +xa' \right)
$$
and unit the element $(1,0)$.
We shall say that the algebra $ A\oplus M$ is the {\it square 0 extension} of $A$ by $M$.
The inclusion $i_1:A\to A\oplus M$ and the projection $p_1:A\oplus M \to A$ are maps of algebras,
and their composite is the identity of $A$, 
$$\xymatrix{
A\ar[r]^-{i_1} & A\oplus M \ar[r]^-{p_1} & A.
}$$

\begin{prop}\label{fderiv}
Let $f:A\to B$ be a map of dg-algebras (unital or not) and $M$ be a $B$-bimodule.
Then a map $d:A\to M$ is a graded $f$-derivation (of degree $0$) iff the map $g:A\to B\oplus M$ defined by putting $g(a)=f(a)+D(a):=(f(a),d(a))$ is a map of graded algebras. $d$ is a $f$-derivation iff the map $g$ is a map of dg-algebras.

In other words, we have a cartesian square in $\gSet$
$$\xymatrix{
\Der(A,f;M) \ar[r]^-g\ar[d]& [A,B\oplus M]\ar[d]^{[A,p_1]}\\
\{*\}\ar[r]^-f & [A,B].
}$$
\end{prop}

\begin{proof}
The condition $g(ab)=g(a)g(b)$ means that we have
$$
ab+d(ab)=(a+d(a))(b+d(b))=ab+d(a)b+ad(b)
$$
and equivalently that $d(ab)=d(a)b+ad(b)$.
Thus, $d$ is a derivation iff $g$ is a map of non-unital algebras.
But if $d$ is a derivation, then we have $g(1)=1$ by lemma \ref{smalllemmaforderivation}.
Thus, $g$ is a map of graded algebras in this case.
For the last assertion, $g$ commutes with the differential iff $f$ commutes with the differential.
\end{proof}

\bigskip

\begin{defi}\label{pointed derivationdef}
If $A=(A,\epsilon_A)$ and $(B,\epsilon_B)$ are two pointed algebras, and $f:A\to B$ is a map of pointed algebras,
we shall say that an $f$-derivation $D:A\rhup B$ is {\it pointed} if $\epsilon_B D=0$.

Pointed $f$-derivations form a dg-vector subspace $\Der_\bullet(f)\subset \Der(f)$. 
If $A=B$ and $f=id_A$ we shall put $\Der_\bullet(A)$ instead of $\Der_\bullet(id_A)$.
\end{defi}

Every pointed derivation $D:A\rhup B$ induces a $f_-$-derivation $D_{-}:A_-\rhup B_-$ of the corresponding non-unital algebras.
Using lemma \ref{smalllemmaforderivation}, it is easy to verify that the map $D\mapsto D_{-}$
is an isomorphism of vector spaces $\Der_\bullet(f)\simeq \Der(f_{-})$.

There exists an analog of proposition \ref{fderiv} for pointed derivations where $A\oplus B$ need to be replaced by the wedge $A\vee B=(A_-\oplus B_-)_+$. We leave the statement to the reader.

\bigskip

If $X$ is a graded vector space and  $D:T(X)\rhup M$ is a derivation of degree $0$
with values in a $T(X)$-bimodule $M$, then
for every $x_1\otimes \cdots \otimes  x_n\in X^{\otimes n}$
we have
$$D(x_1\otimes \cdots \otimes  x_n)=\sum_{i=1}^{n} (x_1\otimes \cdots \otimes  x_{i-1})\cdot Dx_i \cdot (x_{i+1}\otimes  \cdots \otimes  x_n).
$$
This formula shows that  $D$ is determined by its values on $X$.
The following lemma is showing that these values can
be chosen quite arbitrarly.

\begin{lemma}\label{lemmaextensionderiv}
If $X$ is a graded vector space, and $M$ is a $T(X)$-bimodule,
then every morphism $\phi:X\rhup M$ of degree $n$
can be extended uniquely as a derivation $D:T(X)\rhup M$ of degree $n$.
\end{lemma}

\begin{proof} Let us first prove the result in the case where $n=0$.
The linear map $f:X\to  T(X)\oplus M $
defined by putting $f(v)=(v,\phi(v))$
can be extended uniquely
as a map of algebras $g: T(X)\to T(X)\oplus M $.
We have $gi=f$, where  $i:X\to T(X)$ is the inclusion.
Thus, $p_1gi=p_1f=i$,
where $p_1:   T(X)\oplus M  \to T(X)$ is the
projection. But $p_1 g: T(X)\to T(X)$ is an algebra map,
since $p_1$ is an algebra map.
Thus, $p_1g=id$ since the map  $i:X\to T(X)$
is generating.
Thus, $g=(id, D)$, where $D: T(X)\to M$
is a derivation of degree $0$ by proposition \ref{fderiv}.
Moreover, we have $Di=\phi$, since we have $gi=f$.
This proves the existence of $D$ in the case $n=0$.
The uniqueness is left to the reader.
Let us now consider the case of a map $\phi:X\to M$ of general degree $n$.
In this case, $s^{-n}\phi:X\to S^{-n}M$
is a map of degree 0.  It can thus be extended uniquely as a 
derivation $s^{-n}D:T(X)\to S^{-n}M$ of degree 0.
The resulting morphism $D:T(X)\to M$ is a derivation of degree $n$
which is extending $\phi$. Its uniqueness is clear.
\end{proof}

\begin{prop}\label{extensionderivation}
If $X$ is a graded vector space, then
every morphism $\phi:X\rhup T(X)$ of degree $n$ 
can be extended uniquely as a derivation $D:T(X)\rhup T(X)$
of degree $n$.
\end{prop}

\begin{proof} This follows from proposition \ref{lemmaextensionderiv}.
\end{proof}

The result above has a pointed version, which we state without proof:

\begin{prop}\label{extensionderivationpointed}
If $X$ is a graded vector space, then
every morphism $\phi:X\rhup T_\circ(X)$ of degree $n$ 
can be extended uniquely as a (non-unital) derivation $D:T_\circ(X)\rhup T_\circ(X)$ of degree $n$.
or equivalently as a pointed derivation $D:T_\bullet(X)\rhup T_\bullet(X)$ of degree $n$.
\end{prop}

The vector space $ T(X)\otimes X\otimes T(X)$ has the structure of a bimodule
over the algebra $T(X)$. The bimodule is freely generated
by the map  $\phi=1\otimes X\otimes 1:X\to T(X)\otimes X\otimes T(X).$.
The map $\phi$ can be extended uniquely as a derivation $d:T(X)\to T(X)\otimes X\otimes T(X)$
by proposition \ref{lemmaextensionderiv}.

\begin{prop}\label{univderfree}
The derivation $d:T(X)\to T(X)\otimes X\otimes T(X)$ defined above is universal.
Hence we have $$ \Omega_ {T(X)}=T(X)\otimes X\otimes T(X).$$
\end{prop}

\begin{proof} Let $D:T(X)\to M$ be a derivation with values in a $T(X)$-bimodules $M$.
There is then a unique map of bimodules $f:T(X)\otimes X\otimes T(X)\to M$
such that $f\phi=Di(x)$ for every $x\in X$, since the bimodule $T(X)\otimes X\otimes T(X)$
is freely generated by $\phi$. Let us show that $fd=D$.
The map $fd:T(X)\to M$ is a derivation, since $d$ is a derivation and $f$ is a bimodule map.
If $i:X\to T(X)$ is the inclusion, then we have $fdi=f\phi=Di$ by definition of $f$.
It then follows that  $fd=D$ by the uniqueness part of proposition \ref{lemmaextensionderiv}.
\end{proof}

Using propositions \ref{univderfree} and \ref{univfderiv}, we can improve lemma \ref{lemmaextensionderiv}.
\begin{cor}\label{caracderfree}
Let $i:X\to T(X)X$ be the cogenerating map, then the map $D\mapsto Di$ induces an isomorphism of vector spaces
$$
\Der(T(X)) = [X,T(X)].
$$
In particular, a derivation $D$ is zero iff $Di=0$.
\end{cor}

\subsection{dg-Coalgebras}

\subsubsection{Categories of coalgebras}

As for algebras, we shall be concerned with several types of coassociative coalgebras (unital or not, coaugmented or not, differential or not)and 
we will define differential graded counital coassociative coalgebras and then deduce all results concerning non-counital and non-differential coalgebras by forgetting the extra structure. Coaugmented coalgebras are dealt through an equivalence with non-counital coalgebras.

\bigskip

A {\it differential graded counital coassociative coalgebra} (or a {\em dg-coalgebra} for short) is a comonoid object in the monoidal category $\dgVect$.
Explicitely it is 
a graded vector space $C$ equipped with 
a linear map $\Delta:C\to C\otimes C$ (of degree 0) called the {\it comultiplication} or the the {\it coproduct},
a linear map $\epsilon:C\to \FF$ (of degree 0) called the {\it counit} 
and a graded morphism $d:C\rhup C$ of degree $-1$ called the {\it differential}
satisfying the following conditions:
\begin{itemize}
\item {\it coassociativity}
$$\xymatrix{
C \ar[rr]^-{\Delta} \ar[d]_{\Delta} && C\otimes C \ar[d]^-{\Delta\otimes C} \\
C\otimes C\ar[rr]^-{C\otimes \Delta} && C\otimes C\otimes C
}$$
\item {\it left and right counit}
$$\xymatrix{
C \ar[r]^-{\Delta}  \ar@{=}[dr] & C\otimes C \ar[d]^-{\epsilon\otimes C} \\
 & C
}\qquad
\xymatrix{
C \ar[r]^-{\Delta}  \ar@{=}[dr] & C\otimes C \ar[d]^-{C\otimes \epsilon} \\
 & C
}$$
\item {\it coderivation}
$$\xymatrix{
C \ar@^{>}[rr]^-d \ar[d]_{\Delta} && C \ar[d]^-{\Delta} \\
C\otimes C\ar@^{>}[rr]^-{d\otimes C+C\otimes d} && C\otimes C.
}$$
\end{itemize}

A {\em graded counital coalgebra} (we shall say a {\em graded coalgebra}) has only a coassociative comultiplication and a counit, hence it is a comonoid object in $\gVect$.
A {\em non-counital dg-coalgebra} has only a coassociative comultiplication and a differential.
The term {\em coalgebra} will be used to refer indifferently to any of the previous types of coalgebras.

We say that a graded coalgebra $C$ is {\it graded finite} (resp. finite) if its underlying graded vector space is graded finite (resp. finite).

\medskip

We shall often use Sweedler's convention of denoting the sum
$$
\Delta(c)=\sum_{i\in I}c^{(1)}_i \otimes c^{(2)}_i
$$
by $\Delta(c)=c^{(1)}\otimes c^{(2)}$.
The iterated coproduct $\Delta^{(n)}:C\to C^{\otimes n}$ is defined recursively
for $n\geq 2$ by putting $\Delta^{(2)}=\Delta$ and $\Delta^{(n+1)}=(\Delta^{(n)}\otimes C)\Delta$.
By convention, $\Delta^{(1)}=id$ and, if the coalgebra is unital, $\Delta^{0}=\epsilon$.
We shall use Sweedler's convention,
$$
\Delta^{(n)}(c)=c^{(1)}\otimes \dots \otimes c^{(n)}.
$$

A {\it map} of differential graded algebras $f:C\to D$ is a homomorphism of comonoids in $\dgVect$.
Explicitely, it is a map of graded vector spaces $A\to B$ such that the following diagrams commute,
$$\xymatrix{
C \ar[d]_{\Delta_C}  \ar[r]^f & \ar[d]^{\Delta_D}D \\
C\otimes C   \ar[r]^-{f\otimes f} &  D\otimes D
}
\quad\quad\quad
\xymatrix{
C \ar[d]_{\epsilon_C}  \ar[r]^f & C \ar[d]^{\epsilon_D} \\
\FF   \ar@{=}[r] & \FF 
}
\quad\quad\quad
\xymatrix{
C \ar@^{>}[d]_{d_C}  \ar[r]^-{f} & D \ar@^{>}[d]^{d_D} \\
C \ar[r]^-{f} &  D.
}$$
In particular, it is a map of dg-vector spaces.

This means that we have 
$$
f\Delta(c)=f(c^{(1)})\otimes f(c^{(2)}),
\qquad\qquad \epsilon f(c)=\epsilon(c)
\et d_Df(c) = f(d_Cc)
$$
for every $c\in C$.

A {\it map} of graded coalgebras $f:C\to D$ is a map of graded vector spaces $C\to D$ commuting with the comultiplication and the counit only.

A {\it map} of non-counital dg-coalgebras $f:C\to D$ is a map of graded vector spaces $C\to D$ commuting with the comultiplication and the differential only.

\medskip
We shall denote the category of differential graded algebras by $\dgCoalg$, that of graded algebras by $\gCoalg$ and that of non-unital dg-algebras by $\dgCoalg_\circ$.

\begin{rem}
The counit $\epsilon$ of a coalgebra $C=(C,\Delta,\epsilon)$ is uniquely determined by the coproduct $\Delta:C\to C\otimes C$ when it exists. Hence a coalgebra $(C,\Delta,\epsilon)$ is entirely described by the pair $(C,\Delta)$.
However maps of non-counital coalgebras may not preserve the counits when they exist.
\end{rem}

\subsubsection{Atoms and coaugmentations}\label{atoms}

The graded vector space $\FF$ (concentrated in degree 0) has the structure of a dg-coalgebra in which the coproduct $\Delta: \FF\to \FF\otimes \FF$ is the canonical isomorphism and the counit $\epsilon:\FF\to \FF$ is the identity map.
Thus, $\Delta(1)=1\otimes 1$ and $\epsilon(1)=1$.
For this coalgebra structure on $\FF$, the augmentation map $\epsilon:C\to \FF$ of any counital coalgebra is a coalgebra map.

\medskip

\begin{defi}\label{defatom}
We shall say that an element $e$ of a dg-coalgebra $C=(C,\Delta,\epsilon)$ is an {\it atom} or a {\em group like element} 
if $\Delta(e)=e\otimes e$, $\epsilon(e)=1$ and $de=0$.
\end{defi}

The element $1\in \FF$ is an atom in the dg-coalgebra $\FF$.
The image of an atom $e\in C$ by a map of dg-coalgebras $f:C\to D$ is an atom $f(e)\in D$.

\begin{prop}\label{atomsandmap}
If $\phi:\FF\to C$ is a map of dg-coalgebras, then $\phi(1)\in C$ is an atom of $C$. 
Moreover, the map $\phi\mapsto \phi(1)$ is a bijection between the maps of coalgebras  $\FF\to C$ and the atoms of $C$.
\end{prop}
\begin{proof}
Left to the reader.
\end{proof}

We say that a counital coalgebra (differential or not) $(C,\Delta,\epsilon)$ is {\em coaugmented} or {\em pointed} if it it is equipped with an atom $e$, \ie with a coalgebra map $e:\FF\to C$, called the {\em coaugmentation}.
A map of coaugmented coalgebras is a coalgebra map $f:C\to D$ commuting with the coaugmentations
$$\xymatrix{
A   \ar[r]^f  \ar[d]_{\epsilon_A}  & B \ar[d]^{\epsilon_B}  \\
\FF \ar@{=}[r] & \FF.
}$$
The category of coaugmented dg-coalgebras is noted $\dgCoalg_\bullet$. By definition, we have $\dgCoalg_\bullet = \FF\bs\dgCoalg$.

The coalgebra $\FF$ is naturally pointed and it is both initial and terminal in category $\dgCoalg_\bullet $.
Hence the category $\dgCoalg_\bullet $ is {\it pointed}. 
By definition, the {\it  nul map} $\theta_{CD}:C\to D$  between two pointed coalgebras $C$ and $D$ is obtained by composing the counit $\epsilon_C:C\to \FF$ with the coaugmentation $e_D:\FF \to D$.

\bigskip

\begin{defi} \label{nuatomdef} If $C=(C,\Delta)$ is a non-counital (dg-)coalgebra, 
we shall say that an element $e\in C$ is a {\it non-unital atom} if $\Delta(e)=e\otimes e$ and $de=0$.
For example, $0\in C$ is a non-unital atom.
\end{defi}

\begin{prop}\label{seiatomsandmap} If $\phi:\FF\to C$ is a map of non-counital coalgebras, then $\phi(1)\in C$
is a non-unital atom. Moreover, the map $\phi\mapsto \phi(1)$ induces a bijection
between the maps of non-counital coalgebras  $\FF\to C$
and the non-unital atoms of $C$. 
\end{prop}

\begin{proof} Left to the reader.\end{proof}

\subsubsection{Examples}

\begin{prop}\label{C0isacoalgebra} If $C$ is a dg-coalgebra, then $C_0$ has the structure of a coalgebra 
with the coproduct $\Delta_0:C_0\to C_0\otimes C_0$ defined by putting $\Delta_0(x)=(p_0\otimes p_0)\Delta(x)$ where $p_0:C\to C_0$ is the projection.
If $C$ is unital so is $C_0$: $\epsilon:C\to \FF$ factors through $p_0:C\to C_0$ for degree reason and is the induced map $C_0\to \FF$ is a counit for $C_0$.
Moreover, the projection $p_0:C\to C_0$ is a coalgebra map.
\end{prop}

\begin{proof} Left to the reader.
\end{proof}

\begin{ex}[Diagonal coalgebra]\label{diagonalcoalg}
If $I$ is a set viewed as a graded set concentrated in degree 0, then the dg-vector space $\FF I$ freely generated by $I$ has the structure of a dg-coalgebra (concentrated in degree $0$). 
The comultiplication $\Delta:\FF I\to \FF I\otimes \FF I$ is the linear extension of the diagonal $I\to I\times I$ and the counit $\epsilon:\FF I\to \FF$ is the linear extension of the map $I\to 1$. 
We shall say that $\FF I$  is the {\it diagonal coalgebra} of the set $I$. 
In particular, the coalgebra structure on $\FF=\FF1$ is the one defined earlier.

\begin{prop}\label{atomsrightadjoint}
The functor $At:\dgCoalg\to \Set$ which associates to a coalgebra $C$ its set of atoms is right adjoint to the functor which associates to a set $I$ the diagonal coalgebra $\FF I$.
\end{prop}
\begin{proof}
Left to the reader.
\end{proof}

\end{ex}

\begin{ex}[Tensor coalgebra]\label{tensorcoalg}
The tensor algebra of a graded vector space $X$ has the structure of a coalgebra $T^c(X)=(T(X),\Delta,\epsilon)$ called the {\it tensor coalgebra} of $X$. 
By definition, the coproduct $\Delta:T(X)\to T(X)\otimes T(X)$ is defined by putting
$$
\Delta(x_1\otimes \cdots \otimes x_n)=\sum_{i=0}^n (x_1\otimes \cdots \otimes x_i)\otimes (x_{i+1}\otimes \cdots \otimes x_n)
$$
and the counit $\epsilon:T(X)\to \FF$ is the canonical projection.
\end{ex}

\begin{rem}
In opposition to the situation of algebras, the functor $T^c$ is not adjoint to the forgetful functor $U:\dgCoalg\to \dgVect$.
We shall see in proposition \ref{structureconil} that $T^c$ is right adjoint to the forgetful functor when restricted to the subcategory of conilpotent coalgebras. 
We shall see also in theorem \ref{cofreecoalg} that the forgetful functor $U:\dgCoalg\to \dgVect$ has a right adjoint $T^\vee$. 
\end{rem}

\begin{ex}[Pointed tensor coalgebra]\label{ptensorcoalg}
The tensor coalgebra $T^c(X)$ is naturally coaugmented.
The coaugmentation is given by the inclusion of the $\FF\to T^c(X)$
We shall denote $T^c_\bullet(X)$ the tensor coalgebra viewed as a pointed coalgebra.
\end{ex}

\begin{ex}[Non-counital tensor coalgebra]\label{ncutensorcoalg}
The non-unital tensor algebra $T_\circ(X)$ of a dg-vector space $X$ (exemple \ref{tensoralg}) has the structure of a coalgebra noted $T_\circ^c(X)$ and called the {\it non-counital tensor coalgebra} of $X$ when equipped with the {\it reduced coproduct}
$\Delta_{\circ}:T_{\circ}(X)\to T_{\circ}(X) \otimes T_{\circ}(X)$ defined by putting 
$$
\Delta_{\circ}(x_1\otimes \cdots \otimes x_n)=\sum_{0<i<n}(x_1\otimes \cdots \otimes x_i)\otimes (x_{i+1}\otimes \cdots \otimes x_n).
$$
We have $T^c_\bullet(X)_- = T^c_\circ(X)$ and $T^c_\circ(X)_+ = T^c_\bullet(X)$.
\end{ex}

\begin{ex}[Primitive coalgebra]\label{primitivecoalgebra}
If $(X,d)$ is a dg-vector space, then the dg-vector space $X_+:=T^c_1(X):=\FF\oplus X$ has the structure of a dg-coalgebra if the coproduct $\Delta:T^c_1(X)\to T^c_1(X)\otimes  T^c_1(X)$ is defined by putting 
$$
\Delta(1) = 1\otimes 1,\qquad \Delta(x) = x\otimes 1+1\otimes x
$$
for any $x\in X$, and the differential by $d1=0$ and $dx = d_Xx$.
We shall say that $T^c_1(X)$ is the {\it primitive coalgebra} generated by $X$.

The primitive coalgebra is naturally a pointed dg-coalgebra: the coaugmentation is given by the inclusion $\FF\to \FF\oplus X$. We note $T^c_{\bullet,1}(X)$ this pointed coalgebra.

The primitive coalgebra has a non-counital analog $T^c_{\circ,1}(X)$: it is the vector space $X$ with coproduct given by the zero map $X\otimes X\to X$.
We have $T^c_{\bullet,1}(X)_- = T^c_{\circ,1}(X)$ and $T^c_{\circ,1}(X)_+ = T^c_{\bullet,1}(X)$.
\end{ex}

\begin{ex}[Coshuffle coalgebra] \label{coshufflecoalgebra}

Let $X$ be a dg-vector space, the tensor algebra $T(X)$ can be endowed with another counital dg-coalgebra structure called the {\em coshuffle coalgebra} and noted $T^{csh}(X)$. The counit and differential are the same as for $T^c(X)$ but the {\em shuffle coproduct} is defined by 
$$
\Delta (x_1\otimes \dots \otimes x_n) = \sum_{i+j=n} \sum_{\sigma\in Sh(i,j)} (x_{\sigma(1)}\otimes \dots \otimes x_{\sigma(i)})\otimes (x_{\sigma(i+1)}\otimes \dots \otimes x_{\sigma(i+j)}) \ (-1)^{|\sigma|}
$$
where $Sh(i,j)$ is the set of $(i,j)$-shuffle permutations of $n$ elements (see \cite[ch. 1.3.2]{LV} for more details) and $(-1)^{|\sigma|}$ is the sign associated to the permutation $\sigma$ of the graded variables $x_i$ 
(In example \ref{coshuffle bialgebra} will see how this coalgebra structure is part of a bialgebra structure).

\medskip
If case $X$ is generated by a single variable $x$ of even degree, the coproduct on $T^{csh}(x)=\FF[x]$ is simply given by 
$$
\Delta(x^n)=\sum_{k=0}^n {n\choose k} x^k\otimes x^{n-k}.
$$ 
where the ${n\choose k}$ are the binomial coefficients.
In case the variable $x$ is of odd degree, the formula is
$$
\Delta(x^n)=\sum_{k=0}^n \oddbinom{n}{k} x^k\otimes x^{n-k}
$$ 
where the $\oddbinom{n}{k}$ are the odd binomial coefficients.
They can be defined as the coefficients of the expansion of $(x+y)^n$ when $x$ and $y$ are commuting odd variables.
The first coefficients are given by the table
$$\begin{array}{c|cccccccccc}
\rule[-2ex]{0pt}{5ex} \oddbinom{n}{k}& k=0&1&2&3&4&5&6&7&8&\dots\\
\hline
n=0&1\\
1&1&1\\
2&1&0&1\\
3&1&1&1&1\\
4&1&0&2&0&1\\
5&1&1&2&2&1&1\\
6&1&0&3&0&3&0&1\\
7&1&1&3&3&3&3&1&1\\
8&1&0&4&0&6&0&4&0&1\\
\vdots
\end{array}$$
The rule of computation is the same as Pascal triangle $\oddbinom{n+1}{k+1} = \oddbinom{n}{k+1} + \oddbinom{n}{k}$
but with the constraint that $\oddbinom{n}{k}=0$ whenever $n$ is even and $k$ is odd. 
In particular, we have $\oddbinom{2n}{2k}=\binom{n}{k}$.

\end{ex}

\begin{ex}[Dual of finite algebra]\label{dualfinitealgex}
If $(A,m,e)$ is a finite dg-algebra, its dual $A^\star$ is a dg-coalgebra. Because $A$ is finite, the canonical map $[A,\FF]\otimes [A,\FF]\to [A\otimes A,\FF]$ is an isomorphism, then the coproduct is defined as 
$m^\star =[m,\FF]: [A,\FF]\to [A\otimes A,\FF]\simeq [A,\FF]\otimes [A,\FF]$,
the counit by 
$e^\star =[e,\FF]: [A,\FF]\to [\FF,\FF]\simeq \FF$ 
and the differential is the canonical one on $[A,\FF]$.
More conceptually, this is a consequence of lemma \ref{lemmagradedfinite}.
This gives an equivalence of categories
$$
(\dgCoalg^{\sf fin})^{op} \simeq \dgAlg^{\sf fin}.
$$
\end{ex}

\begin{ex}[Dual of bounded graded finite algebra]\label{dualboundedgradedfinitealg}
We have the following strengthening of the previous construction.
If $(A,m,e)$ is a graded finite dg-algebra, bounded above or below, then its dual $A^\star$ is a dg-coalgebra.
This is a consequence of the equivalence of monoidal categories of lemma \ref{lemmagradedfinite}.

Let $\dgCoalg^{\sf gr.fin.a}$ and $\dgAlg^{\sf gr.fin.a}$ (resp. $\dgCoalg^{\sf gr.fin.b}$ and $\dgAlg^{\sf gr.fin.b}$) be the categories of graded finite bounded above (resp. below) dg-coalgebra and dg-algebras.

\begin{prop}\label{equivgrfinalgcoalg}
The contravariant equivalence of lemma \ref{lemmagradedfinite} induces equivalences
$$
(\dgCoalg^{\sf gr.fin.a})^{op} \simeq \dgAlg^{\sf gr.fin.b}
\et 
(\dgCoalg^{\sf gr.fin.b})^{op} \simeq \dgAlg^{\sf gr.fin.a}.
$$
\end{prop}
\end{ex}

\begin{ex}[Endomorphism coalgebra]\label{coendomorphisms}
If $X$ is a finite dg-vector space, then the algebra $End(X)=[X,X]=X^\star \otimes X$ is canonically self dual.
It follows that it has the structure of a finite dg-coalgebra denoted $End^c(X)$.
The coalgebra structure can be described explicitly as follows.  
The object $X$ is perfect, since it is finite.
Let $\eta: \FF \to X\otimes X^\star$ be the unit of the duality $X^\star \dashv X$
and $\epsilon=ev:X^\star \otimes X\to \FF$
be the counit.
Then the coproduct $\Delta:End(X) \to End(X)\otimes End(X)$
is the composite
$$\xymatrix{
X^\star\otimes X \ar[rr]^-{X^\star \otimes \eta \otimes X}&& X^\star \otimes X  \otimes X^\star \otimes X.
}
$$
and the counit $\epsilon:End(X)\to \FF$
is the counit $\epsilon:X^\star \otimes X\to \FF$
\end{ex}

\begin{ex}
Recall that a a {\it Frobenius-Poincar\'e algebra} of level $n\in \mathbb{Z}$ is a finite dimensional dg-algebra $A=(A,\mu,e)$
equipped with a linear map $t: A\to S^n$ (called the {\it orientation}) for which the pairing $t\mu:A\otimes A\to A\to S^n$
is symmetric and perfect. If we transport the algebra structure of $A$ along the isomorphism $ \lambda^1(t\mu ):A \to [A,S^n]$
we obtain an algebra structure on $[A,S^n]=A^\star\otimes  S^n$, and hence a coalgebra structure
on $A\otimes S^{-n}$.
\end{ex}

\begin{ex}[Cocommutative coalgebras]
Recall that a dg-coalgebra $C$ is said to be {\it graded cocommutative}
if we have $x^{(1)}\otimes x^{(2)}=x^{(1)}\otimes x^{(2)}(-1)^{|x^{(1)}||x^{(2)}|}$ for every $x\in C$.
The full subcategory of $\dgCoalg$ generated by cocommutative dg-coalgebras is noted $\dgCoalgcom$.
The inclusion functor $\dgCoalgcom\to \dgCoalg$ has a right adjoint sending $C$ to the subcoalgebra $C^{ab}\subset C$ of elements with symmetric coproducts.
\end{ex}

\subsubsection{Comparison functors}\label{comparisonfunctorcoalgebras}

There exists several functors between the categories of coalgebras.
For a dg-coalgebra, we can forget the differential
$$\xymatrix{
U_d: \dgCoalg \ar[r] & \gCoalg
}$$
or the counit
$$\xymatrix{
U_\epsilon: \dgCoalg \ar[r] & \dgCoalg_\circ.
}$$
And for pointed dg-algebra, we can forget the coaugmentation
$$\xymatrix{
U_e:\dgCoalg_\bullet \ar[r] & \dgCoalg
}$$

If $C=(C,\Delta,\epsilon, e)$ is a pointed coalgebra, then the dg-vector space $\overline{C}=Coker(e)=C/\FF e$ has the structure of a non-unital coalgebra, where the coproduct $\overline{\Delta}:\overline{C} \to \overline{C}\otimes \overline{C}$ is induced by $\Delta$.
Let us denote by $C_{-}$ the kernel of the counit $\epsilon:C\to \FF$.
We have a decomposition $C=\FF e \oplus C_{-}$ since we have $\epsilon(e)=1$. 
Hence the quotient map $C\to \overline{C}$ induces an isomorphism $i:C_{-}\simeq \overline{C}$.
We shall denote by $\Delta_{-}$ the coproduct on $C_{-}$ obtained by transporting  the coproduct $\overline{\Delta}$ along the isomorphism $i:C_{-}\simeq \overline{C}$.
It is easy to see that we have
$$
\Delta_{-}(x)=\Delta(x)-x\otimes e-e\otimes x
$$
for every $x\in C_{-}$.
This defines a non-unital coalgebra $C_{-}=(C_{-}, \Delta_{-})$ and a functor 
$$\xymatrix{
(-)_-:\dgCoalg_\bullet \ar[r] & \dgCoalg_\circ
}$$
Conversely, if $(C,\Delta)$ is a non-counital coalgebra, the dg-vector space $C_+:=\FF\oplus C$ can be equipped with a counital coalgebra structure by putting 
$\Delta_{+}(e)=e\otimes e$ for $e=(1,0)$ and 
$$\xymatrix{
\Delta_{+}(x)=e\otimes x+x\otimes e+ \Delta(x)
}$$
for $x=(0, x)\in \FF\oplus C$. 
The counit $\epsilon:C_+ \to \FF$ is the projection $\FF \oplus C \to \FF$ and 
the element $e=(1,0)\in \FF\oplus C$ is an atom, thus $C_+$ is naturally pointed.
This defines a functor 
$$\xymatrix{
(-)_+:\dgCoalg_\circ \ar[r] & \dgCoalg_\bullet
}$$
The following result is elementary but useful.

\begin{prop}\label{equivnonunitalpointedcoalg}
The two functors
$$\xymatrix{
\dgCoalg_\bullet \ar@<.6ex>[rr]^-{(-)_-} && \dgCoalg_\circ \ar@<.6ex>[ll]^-{(-)_+}
}$$
are inverse equivalences of categories.
\end{prop}

\begin{cor}\label{adjunctioncounitalnoncounital}
\begin{enumerate}
\item The composition $U_e(-)_+:\dgCoalg_\circ\to \dgCoalg_\bullet \to \dgCoalg$ is right adjoint to $U_\epsilon$ 
$$\xymatrix{
U_\epsilon :\dgCoalg \ar@<.6ex>[rr] && \dgCoalg_\circ \ar@<.6ex>[ll]{:U_e(-)_+}.
}$$
The counit of the adjunction is the projection $C\oplus \FF \to C$.

\item The composition $(U_\epsilon-)_+:\dgCoalg\to \dgCoalg_\circ \to \dgCoalg_\bullet$ is left adjoint to $U_e$ 
$$\xymatrix{
(U_\epsilon-)_+ :\dgCoalg \ar@<.6ex>[rr] && \dgCoalg_\bullet \ar@<.6ex>[ll]{:U_e}.
}$$
The unit of the adjunction is the inclusion $C\to C\oplus \FF$.
\end{enumerate}
\end{cor}

\begin{rem}
We can summarize this corollary by saying that, up to the equivalence $\dgCoalg_\bullet\simeq \dgCoalg_\circ$, the functor $U_\epsilon$ is left adjoint to $U_e$.
\end{rem}

\bigskip

We shall prove in section \ref{dgtoghomcoalg} that the forgetful functor $U_d:\dgCoalg\to \gCoalg$ has both a left and a right adjoint.
With the consequence that limits and colimits exist in $\dgCoalg$ and they can be computed in $\gCoalg$.

\subsubsection{Opposite coalgebra and anti-homomorphism}

If $C$ and $D$ are graded coalgebras, we shall say that a 
linear map $f:C\to D$ is an {\it anti-homomorphism}
if it is an anti-homorphism of comonoids.
This means that we have $$f\Delta(c)=f(c^{(2)})\otimes f(c^{(1)})(-1)^{|c^{(1)}||c^{(2)}|}$$
and $\epsilon f(c)=\epsilon(c)$ for every $c\in C$.
The {opposite} of a graded coalgebra $C$
is a graded algebra $C^o$ anti-isomorphic to $C$. 
By definition, the map $x\mapsto x^o$ is a bijection 
$C\to C^o$ and  for every $c\in C$ we have
$$\Delta(c^o)=(c^{(2)})^o\otimes (c^{(1)})^o(-1)^{|c^{(1)}||c^{(2)}|}\et \epsilon c^o=\epsilon(c).$$
There is an obvious bijection between the anti-homomorphisms
of coalgebras $C\to D$ and the homomorphisms
of coalgebras $C^o\to D$ and $C\to D^o$.

\subsubsection{Monoidal structures}\label{monoidalstructurecoalgebras}

The tensor product  $C\otimes D$ of (graded) coalgebras $C$ and $D$ has the structure of a coalgebra, with the comultiplication equal to the composite of the maps
$$\xymatrix{
C\otimes D \ar[rr]^-{\Delta \otimes \Delta }  &&  C\otimes C\otimes D\otimes D 
	\ar[rr]^-{C\otimes \sigma \otimes D} & & (C\otimes D)\otimes (C\otimes D) 
}$$
By definition,
$$
\Delta(c\otimes d)=(c^{(1)}\otimes d^{(1)}) \otimes (c^{(2)}\otimes d^{(2)})  (-1)^{|d^{(1)}||c^{(2)}|}
$$
for $c\in C$ and $d\in D$.

If $C$ and $D$ are counital, so is $C\otimes D$ with counit defined by $\epsilon_A\otimes \epsilon_B: A\otimes B\to\FF\otimes \FF \simeq \FF$.
If $C$ and $D$ are differential, so is $C\otimes D$ with differential defined by $d_{C\otimes D}=d_C\otimes D+C\otimes d_D$, as for dg-vector spaces.

\medskip
The category $\dgCoalg$ is symmetric monoidal, the unit object is the coalgebra $\FF$. 
The category $\gCoalg$ is symmetric monoidal, the unit object is the coalgebra $\FF$. 
The category $\dgCoalg_\circ$ is symmetric monoidal, the unit object is the coalgebra $\FF e$ where $e$ is a non-unital atom. 
The three forgetful functors $\dgCoalg\to \dgVect$, $\gCoalg\to \gVect$ and $\dgCoalg_\circ\to \dgVect$ are clearly symmetric monoidal.

\medskip
Using the equivalence of proposition \ref{equivnonunitalpointedcoalg}, we can transport the monoidal structure of $\dgCoalg_\circ$ to a monoidal $\dgCoalg_\bullet$. We define the {\it smash product} $C\wedge D$ of two pointed coalgebras as the smash product of the underlying pointed vector spaces
$$
C\wedge D = ( C_{-} \otimes  D_{-})_{+}.
$$
The unit object for the smash product is the pointed coalgebra  $\FF_{+}=\FF\oplus \FF e$ where 
$$
\Delta(1)=1\otimes 1
\et
\Delta(e)=1\otimes e+e\otimes 1+e\otimes e
$$
with the base point $1=(1,0)$. The counit $\epsilon:\FF_{+}\to \FF$ is the first projection $\FF\oplus \FF e\to \FF$.
Notice that $\Delta(1+e)=(1+e)\otimes (1+e)$. Thus, $\FF_{+}=\FF1\oplus  \FF(1+e)$ is isomorphic to the coproduct of the coalgebra $\FF$ with itself.

$(\dgCoalg_\bullet,\wedge , \FF_{+})$ is a symmetric monoidal category and the forgetful functor $\dgCoalg_\bullet\to \dgVect_\bullet$ is monoidal.

The pointed algebra $C\wedge D$ has also the following description.

\begin{prop}\label{smashproductofcoalgebras}
We have a pushout square of coalgebras
$$\xymatrix{
C\oplus D \ar[r]^-{(\epsilon_C,\epsilon_D)} \ar[d]_{(C \otimes e_D, e_C\otimes  D)} & \ar[d]^{e}  \FF \\
C\otimes D \ar[r] & C\wedge D.
}$$ 
\end{prop}

\begin{proof} 
We know that the square is a pushout by lemma \ref{pointedbicart}.
The rest is straightforward computation dual to that of proposition \ref{smashproductofalgebras}.
\end{proof}

\medskip

\begin{prop}\label{monoidalfunctorscoalg}
\begin{enumerate}
\item The forgetful functor
$$\xymatrix{
U_d: \dgCoalg \ar[r] & \gCoalg
}$$
is monoidal. 

\item The forgetful functor
$$\xymatrix{
U_\epsilon: \dgCoalg \ar[r] & \dgCoalg_\circ.
}$$
is monoidal. Hence its right adjoint $U_e(-)_+$ is lax monoidal.

\item The forgetful functor
$$\xymatrix{
U_e:(\dgCoalg_\bullet,\wedge,\FF_+) \ar[r] & (\dgCoalg,\otimes,\FF)
}$$
is lax monoidal and its left adjoint is monoidal.
The lax structure on $U_e$ is given by the map $C\otimes D\to C\wedge D$ from proposition \ref{smashproductofcoalgebras} and by the map $\FF\to \FF\oplus \FF e\simeq \FF_+$ pointing the atom $e+1$.
\end{enumerate}
\end{prop}
\begin{proof} 
Dual to that of proposition \ref{monoidalfunctorsalg}.
\end{proof}

\subsubsection{Conilpotent coalgebras}

If $(A,m)$ is a non-unital algebra, and $n\geq 1$ let us denote by  $A^{(n)}$
the image of the $n$-fold product map  $m^{(n)}:A^{\otimes n}\to A$.
Recall that a non-unital algebra $A$ is said to be {\em nilpotent} if $A^{(n)}=0$
for some $n\geq 1$.
We shall say that a pointed algebra $A$ is {\em nilpotent} if the corresponding non-unital algebra $A_-$ is nilpotent.

\medskip

Let us denote the category of {\em finite} non-unital algebras by $\dgAlgfincirc$  and 
its subcategory of finite {\em nilpotent} non-unital algebras by $\dgAlgfnilcirc$.

\begin{lemma} \label{reflectionnilpotent} 
The inclusion functor $\dgAlgfnilcirc \subset \dgAlgfincirc$
has a left adjoint
$$R: \dgAlgfincirc \to \dgAlgfnilcirc.$$
\end{lemma}

\begin{proof} If $A$ is a finite non-unital algebra,
then the chain of ideals
$$A^{(1)}\supseteq A^{(2)}\supseteq A^{(3)}\supseteq \cdots $$
is stationnary, since $A$ is finite.
Let $n_{0}> 0$ be an integer such that
$A^{(n_0+1)}=A^{(n_0)}$.
The non-unital algebra
$R(A)=A/A^{n_0}$ is nilpotent
since we have $R(A)^{(n_0)}=A^{(n_0)}/A^{(n_0)}=0$.
Let us show that the quotient map $q:A\to R(A)$
is reflecting $A$ in the subcategory $ \dgAlgfnilcirc$.
If $B$ is a nilpotent algebra,
then we have $ B^n=0$ for some $n>0$ and 
we may suppose that $n\geq n_0$.
If $f:A\to B$ is a map of non-unital algebras,
then we have 
$f(A^{n_0})=f(A^{n})\subset B^{n}=0$.
It follows that there is a unique map of non-unital algebras
$f':R(A)\to B$ such that $f'q=f$.
\end{proof}

If $(C,\Delta)$ is a non-counital coalgebra, and $n\geq 1$ let us denote $\Delta^{(n)}$ the map $C\to C^{\otimes n}$ defined inductively by putting 
$\Delta^{(1)}(x)=x$ and 
$$
\Delta^{(n)}(x)=(\Delta^{(n-1)}\otimes C)\Delta(x)
$$
for $n>1$. We may use Sweedler notation and write that $\Delta^{(n)}(x)=x^{(1)}\otimes \cdots \otimes x^{(n)}$.
The coassociativity of $\Delta$ implies that we have $\Delta^{(p+q)}=(\Delta^{(p)}\otimes \Delta^{(q)})\Delta$ for every $p,q\geq 1$.
Thus,
$$x^{(1)}\otimes \cdots \otimes x^{(p+q)}=\Delta^{(p)}(x^{(1)}) \otimes \Delta^{(q)}(x^{(2)}).$$

\medskip

\begin{defi} \label{defconilpotent} 
If $(C,\Delta)$ is a non-counital coalgebra, we shall say that an element $x\in C$ is {\it conilpotent} if we have $\Delta^{(n)}(x)=0$ for $n\gg0$.
We shall say that a non-counital coalgebra $(C,\Delta)$ is {\it conilpotent} if every element $x\in C$ is conilpotent.
We shall say that a pointed coalgebra $(C,\Delta, e)$ is {\it conilpotent} if the corresponding non-counital coalgebra $C_{-}$ is conilpotent.
\end{defi}

\begin{ex}\label{tensorconilpotent}
In the tensor non-counital coalgebra $T^c_{\circ}(X) = (T_{\circ}(X),\Delta_{-})$ of example \ref{ncutensorcoalg}, we have $\Delta^{(N)}(x_1\otimes \dots \otimes x_n)=0$ as soon as $N>n$. This proves that $T^c_{\circ}(X)$ is conilpotent,
hence also the pointed coalgebra $T^c_\bullet (X)=(T_\bullet(X),\Delta,\epsilon,1)$.
\end{ex}

\begin{ex}\label{primitiveconilpotent}
The non-counital primitive coalgebra $T^c_{\circ,1}(X)=(X,\Delta=0)$ of example \ref{primitivecoalgebra} is clearly conilpotent since $\Delta_-(x)=0$ for every $x\in X$.
Hence the pointed primitive coalgebra $T^c_1 (X)=(T^c_{\circ,1}(X))_+$ is also conilpotent.
\end{ex}

\begin{ex}\label{coshuffleconilpotent}
The coshuffle coproduct on $T(X)$ of example \ref{coshufflecoalgebra} is pointed by 1 and conilpotent.
Each iteration of the coproduct $\Delta_-$ reduces strictly the size of the factors and $\Delta^{(N)}_-(x_1\otimes \dots \otimes x_n)=0$ as soon as $N>n$.
\end{ex}

\medskip
\begin{prop}\label{singleatom}
If a pointed coalgebra $(C,e)$ is conilpotent, then $e$
is the unique atom of $C$.
\end{prop}

\begin{proof} We shall use the decomposition $C=\FF e\oplus C_{-}$.
If $f=\lambda e+x$ is an atom, then $\lambda =1$, since $\epsilon(f)=1$.
Thus, 
$$\Delta(f)=\Delta(e)+\Delta(x) =e\otimes e+e\otimes x+x\otimes e +\Delta_{-}(x).$$
Hence the condition
$$\Delta(f)=f\otimes f=( e+x)\otimes (e+x)=e\otimes e+e\otimes x+x\otimes e+x\otimes x$$
is equivalent to the condition $\Delta_{-}(x)=x\otimes x$.
But we have $\Delta_{-}^{(n)}(x)=0$ for some $n>0$, since the non-counital coalgebra $C_{-}$
is conilpotent. Thus, $x^n=0$ in $C^{\otimes n}$.
It follows that $x=0$.
\end{proof}

Let $C$ be a counital coalgebra, because of proposition \ref{singleatom} we can define $C$ to be conilpotent iff it has a single atom and it is conilpotent when pointed by this atom. Conilpotency becomes a property for coalgebras. This distinguishes a full subcategory $\dgCoalgnil$ of $\dgCoalg$.

We shall denote the category of conilpotent non-counital dg-coalgebras by $\dgCoalgnilcirc$ and the category of conilpotent pointed coalgebras by $\dgCoalgnilbullet$.
The functors
$$\xymatrix{
(-)_+:\dgCoalgnilcirc\ar@<.6ex>[r]&\ar@<.6ex>[l] \dgCoalgnilbullet:(-)_{-}
}$$
are inverse equivalences of categories. 
The functor forgetting the coaugmentation $\dgCoalgnilbullet\to \dgCoalg$ is fully faithful and realizes an equivalence
$\dgCoalgnilbullet\simeq \dgCoalgnil$.

\begin{lemma}\label{finitenilcoalg}
A finite non-counital coalgebra $C$ is conilpotent if and only if the dual non-unital algebra $C^\star$ is nilpotent.
\end{lemma}

\begin{proof}
Obviously, a finite non-counital coalgebra is conilpotent if and only if $\Delta^{(n)}=0$ for some $n>0$.
This proves the result, since the $n$-fold product $m^{(n)}:(C^\star)^{\otimes n}\tto C^\star$ is obtained by
transposing the $n$-fold coproduct $\Delta^{(n)}):C\to C^{\otimes n}$.
\end{proof}

The {\it conilpotent radical} (or simply the {\em radical}) $R^c(C)$ of a non-counital coalgebra $C$ is defined to be the dg-vector space of conilpotent elements $x\in C$.
A non-counital coalgebra $C$ is conilpotent if and only if $R^c(C)=C$.
Similarly, the {\it radical} $R^c(C)$ of a pointed coalgebra $C$ is defined by putting $R^c(C)=R^c(C_{-})_+$.

\begin{prop}  \label{radical} 
The radical $R^c(C)$ of a non-counital coalgebra $C$ is the largest conilpotent sub-coalgebra of $C$.
The functor $R^c:\dgCoalg_{\circ}\to \dgCoalgnilcirc$ is right adjoint to the inclusion $\dgCoalgnilcirc \subset \dgCoalg_{\circ}$.
Similarly, the functor
$R^c:\dgCoalg_{\bullet}\to  \dgCoalgnilbullet$
is right adjoint to the inclusion $\dgCoalgnilbullet \subset \dgCoalg_\bullet$.
\end{prop}

\begin{proof} 
Let us show that $R^c(C)$ is a sub-coalgebra of $(C,\Delta)$.
We have $R^c(C)\otimes R^c(C)=(R^c(C)\otimes C )\cap (C\otimes  R^c(C))$.
Hence it suffices to show that $\Delta(R^c(C))\subseteq C\otimes  R^c(C)$, since the inclusion $\Delta(R^c(C))\subseteq R^c(C)\otimes  C$ follows by symmetry.
Let $(u_i:i\in I)$ be a graded basis of $C$. For every $x\in C$ we have
$$
\Delta(x)=\sum_{i\in I} u_i\otimes x_i
$$
where $x_i=0$ except for a finite number of elements $i\in I$.
If $x\in R^c(C)$, let us show that $x_i\in R^c(C)$.
But we have
$$
\Delta^{(n)}(x)=\sum_{i\in I} u_i\otimes \Delta^{(n-1)}(x_i)
$$
for every $n>1$. 
Hence we have $\Delta^{(n-1)}(x_i)=0$ if $\Delta^{(n)}(x)=0$.
This proves that $\Delta(R^c(C))\subseteq C\otimes  R^c(C)$ and hence that $R^c(C)$  is a sub-coalgebra of $C$.
It is clear by construction of $R^c(C)$ that it is the largest conilpotent sub-coalgebra of $C$.
If $D$ is a  nilpotent non-counital coalgebra, then every map $f:D\to C$
factors (uniquely) through the inclusion $R^c(C)\subseteq C$,
since the image by $f$ of a conilpotent element is conilpotent.
This shows that the functor $R^c:\dgCoalg_{\circ}\to \dgCoalgnilcirc$ is right adjoint to the inclusion $\dgCoalgnilcirc \subset \dgCoalg_{\circ}$.
The second statement is deduced by equivalence.
\end{proof}

\medskip

If $X$ is a dg-vector space, we shall say that a conilpotent non-counital coalgebra $C$ equipped with a linear map $p: C\to X$ is {\it cofree conilpotent} if for any conilpotent non-counital coalgebra $E$ and any linear map $f:E\to X$, there exists a unique map of non-counital coalgebras $g:E\to C$ such that $pg=f$.
$$\xymatrix{
E \ar[drr]_-{f} \ar@{-->}[rr]^-{g}&&C \ar[d]^{p} \\
&& X.
}$$
Similarly, we shall say that a conilpotent pointed coalgebra $C=(C,e)$ equipped with a linear map $p: C_{-}\to X$ such that is {\it cofree conilpotent} if for any conilpotent pointed coalgebra $E$ and any linear map $f:E_{-}\to X$, there exists a unique map of pointed coalgebras $g:E\to C$ such that $pg_{-}=f$.

\medskip

\begin{lemma}\label{coniltensorcolag}
If $X$ is a graded vector space, and $p:T^c_\circ (X)\to X$  is the projection then we have
$$
x=\sum_{n\geq 1}p^{\otimes n} \Delta^{(n)}(x)
$$
for every $x\in T^c_\circ (X)$, where $\Delta$ is the reduced deconcatenation.
\end{lemma}

\begin{proof} If $x=x_1\otimes \cdots \otimes x_k\in X^{\otimes k}$,
then 
$p^{\otimes n} \Delta^{(n)}(x)=0$ unless $n=k$ in which case
$p^{\otimes n} \Delta^{(n)}(x)=x$.
\end{proof}

\begin{prop}\label{structureconil} 
If a non-unital coalgebra $C$ is conilpotent then for any linear map $f:C\to X$,
there exists a unique map of non-counital coalgebras $g:C\to T^c_\circ(X)$ such that $pg=f$.
Moreover, we have
$$g(x)=\sum_{n\geq 1}f^{\otimes n} \Delta^{(n)}(x)=\sum_{n\geq 1} f(x^{(1)}) \otimes \cdots \otimes f(x^{(n)})$$
for every $x\in C$.

Hence the non-counital coalgebra $T^c_\circ(X)$ equipped 
with the projection $p:T^c_\circ (X)\to X$ is cofree conilpotent.
Similarly, the pointed coalgebra $T^c_\bullet (X)$ equipped with the projection $p:T^c(X)\to X$ is cofree conilpotent.
\end{prop}

\begin{proof} 
Let us show that if $C$ is a conilpotent non-unital coalgebra and  $f:C\to X$ is a linear map,
then there exists a unique map of non-counital coalgebras $g:C\to T^c_\circ(X)$ such that $pg=f$.
The formula in the proposition is defining a map $g:C\to T^c_\circ(X)$,
since the sum  is finite by the conilpotency of $x$.
Let us show $g$ is a map of non-counital coalgebras. 
For every $x\in C$, we have
\begin{eqnarray*}
\Delta g(x) &=&\sum_{n\geq 1} \Delta \bigl( f(x^{(1)}) \otimes \cdots \otimes f(x^{(n)})\bigr)\\
&=& \sum_{n\geq 1} \sum_{1<i<n} \bigl( f(x^{(1)}) \otimes \cdots \otimes f(x^{(i)})\bigr)\otimes \bigl( f(x^{(i+1)}) \otimes \cdots \otimes  f(x^{(n)}) \bigr)\\
&=& \sum_{i.j \geq 1} f^{\otimes i}\Delta^{(i)}(x^{(1)})\otimes  f^{\otimes j}\Delta^{(j)}(x^{(2)})\\
&=& \sum_{i\geq 1} f^{\otimes i}\Delta^{(i)}(x^{(1)})\otimes  \sum_{j\geq 1} f^{\otimes j}\Delta^{(j)}(x^{(2)})\\
&=& g(x^{(1)})\otimes g(x^{(2)})\\
&=& (g\otimes g)\Delta(x).
\end{eqnarray*}
by coassociativity.
This proves $g$ is a map of non-counital coalgebras. 
Obviously, we have $pg =f$. It remains to prove the uniquness of $g$.
If $h:C\to T^c_\circ (X)$ is a map of non-counital coalgebras such that
$ph=f$, let us show that $h=g$. 
For every $x\in C$ we have
$$h(x)=\sum_{n\geq 1}p^{\otimes n} \Delta^{(n)}(h(x))$$
by lemma \ref{coniltensorcolag}.
But we have
$\Delta^{(n)}(h(x))=h^{\otimes n}\Delta^{(n)}(x)$,
since $h$ is a map of non-counital coalgebras.
Thus
\begin{eqnarray*}
h(x) &=& \sum_{n\geq 1}p^{\otimes n}h^{\otimes n} \Delta^{(n)}(x)\\
 &=& \sum_{n\geq 1}(ph)^{\otimes n} \Delta^{(n)}(x)\\
  &=& \sum_{n\geq 1}f^{\otimes n} \Delta^{(n)}(x)\\
  &=& g(x).
\end{eqnarray*}
We have proved that the conilpotent coalgebra $T^c_\circ(X)$ is cofreely cogenerated 
by the projection $p:T^c_\circ (X)\to X$.
\end{proof}

\begin{lemma}  \label{radical2} 
The tensor product $C\otimes D$ of two non-counital coalgebras is conilpotent if  $C$ or $D$ is conilpotent.
Similarly, the smash product $C\wedge D$ of two pointed coalgebras is conilpotent if $C$ or $D$ is conilpotent.
\end{lemma}

\begin{proof} 
If $C$ and $D$ are non-counital coalgebras, then we have
$$
\Delta^{(n)}(x\otimes y)= \sigma^{(n)} (\Delta^{(n)}(x)\otimes \Delta^{(n)}(y))
$$
for every $x\in C$ and $y\in D$, where $\sigma^{(n)}$ is the canonical isomorphism $C^{\otimes n}\otimes D^{\otimes n} \simeq  (C\otimes D)^{\otimes n}$.
It follows that $\Delta^{(n)}(x\otimes y)=0$ if $\Delta^{(n)}(x)=0$.
\end{proof}

\subsubsection{Comodules}

We shall only develop the theory of counital comodules, that of non-counital comodules is analog and all results are valid also in this context.

\bigskip

Recall that if $C=(C,\Delta_C,\epsilon)$ is a dg-coalgebra, then a left $C$-{\it comodule} is a dg-vector space $X$ equipped with a linear map $\alpha: X\to C\otimes X$ which is coassociative and counital,
$$\xymatrix{
X \ar[rr]^-{\Delta_X} \ar[d]_{\Delta_X} && C\otimes X \ar[d]^-{\Delta_C \otimes X} \\
C\otimes X\ar[rr]^-{C\otimes \Delta_X} && C\otimes C\otimes X
}
\quad \quad \quad 
\xymatrix{
X \ar[r]^-{\Delta_X}  \ar@{=}[dr] & C\otimes X \ar[d]^-{\epsilon\otimes X} \\
 & X
}$$
We shall occasionnally follow Sweedler's convention of denoting the sum
$$
\Delta_X(x)=\sum_{i\in I}x_i^{(1)}\otimes x_i^{(0)}
$$
by $\Delta_X(x)=x^{(1)}\otimes x^{(0)}$.
If $X=(X,\Delta_X)$ and $Y=(Y,\Delta_Y)$ are left $C$-comodules, 
then a  {\it morphism}  $f:X\to Y$ is linear map 
such that $\Delta_Y f(x)=x^{(1)}\otimes f(x^{(0)})$ for every $x\in X$.
We shall denote the category of (left) $C$-comodules
by ${\Comod}(C)$. 

If $X$ and $Y$ are two $C$-comodules, the dg-vector space of morphisms of $C$-modules between them is the object $\Hom_C(X,Y)$ defined by the equalizer
$$\xymatrix{
\Hom_C(X,Y)\ar[r]& [X,Y] \ar@<.6ex>[rr]^-{[X,\Delta_Y]}\ar@<-.6ex>[rr]_-{[\Delta_X,C\otimes -]} && [X,C\otimes Y]
}$$
where $\Delta_X: X\to C\otimes X$ and $\Delta_Y:Y\to C\otimes Y$ are the $C$-comodule structures.
Concretly, a graded morphism $f\in [X,Y]_n$ belongs to $\Hom_{A}(X,Y)_n$ iff we have 
$$
\Delta_Y f(x)=x^{(1)}\otimes f(x^{(0)}) \ (-1)^{n|x^{(1)}|}
$$
for every $a\in A$ and $x\in X$.

\bigskip

If $C=(C,\Delta,\epsilon)$ is a dg-coalgebra, then the functor $C\otimes(-):\dgVect \to \dgVect$ has the structure of a comonad with multiplication $\Delta \otimes X:C\otimes X \to C\otimes C\otimes X$ and counit $\epsilon \otimes X:C\otimes X \to X$. 
A coalgebra for this comonad is a (left) $C$-comodule.

\begin{prop} \label{gammavectforget}
The functors $X\to C\otimes X$ is right adjoint to the forgetful functor ${\Comod}(C) \to\dgVect$.
Colimits exists in ${\Comod}(C)$ and they can be computed in $\dgVect$.
\end{prop}
\begin{proof}
This follows from the general theory of of coalgebras over a comonad.
\end{proof}

A $C$-comodule $N$ is said to be {\em cofreely generated} by a map $p:N\to X$ if the corresponding map $N\to C\otimes X$ is an isomorphism of $C$-comodules. Through this isomorphism, the map $p$ identifies with the map $\epsilon \otimes X:C\otimes X \to X$.
A comodule $N$ is said to be {\em cofree} if there exists a map $p:N\to X$ which cofreely generates it.

\medskip

The category ${\Comod}(C)$ is enriched, tensored and cotensored over the category $\dgVect$.
If $(X,\alpha)$ is a $C$-comodule and $V$ is a dg-vector space, then the composite 
$$\xymatrix{
V\otimes X\ar[rr]^-{V\otimes \alpha} && V\otimes C\otimes X \ar[r]^-{\sigma\otimes X} & C\otimes V\otimes X
}$$
is a coaction (again denoted $\alpha$) of $C$ on $V\otimes X$. 
By definition, $\alpha(v\otimes x)=x^{(1)}\otimes (v\otimes  x^{(0)}) (-1)^{|x^{(1)} ||v|}$ for $v\in V$ and $x\in X$, where $\alpha(x)=x^{(1)}\otimes x^{(0)}$.
The resulting $C$-comodule $V\otimes X$ is the tensor product of $(X,\alpha)$ by $V$.
In particular, the $n$-fold {\it supension}  $S^n(X)=S^n\otimes X$ of $X=(X,\alpha)$ is the $C$-comodule defined by putting
$\alpha(s^nx)=x^{(1)}\otimes (s^n x^{(0)}) (-1)^{|x^{(1)} | n}$.
The dg-vector spaces $\hom^C(X,Y)$ of morphisms between two $C$-comodules $X=(X,\alpha)$ and $Y=(Y,\alpha)$ is a subspace of $[X,Y]$. 
By construction,  a morphism $f:X\rhup Y$  of degree $n$ belongs to $\hom^C(X,Y)$ iff we have
$$
\alpha f(x)=x^{(1)} \otimes f(x^{(0)})\ (-1)^{n|x^{(1)}|}
$$
for every $x\in X$, where  $\alpha(x)=x^{(1)}\otimes x^{(0)}$.

\medskip

Dually, a  {\it right} $C$-comodule is a dg-vector space $X$ equipped with a linear map $\beta:X\to X\otimes C$ which is coassociative and counital,
$$\xymatrix{
X \ar[rr]^-{\beta} \ar[d]_{\beta} && X\otimes C \ar[d]^-{X\otimes \Delta} \\
X\otimes C \ar[rr]^-{\beta \otimes C} && X\otimes C\otimes C
}
\quad \quad \quad 
\xymatrix{
X \ar[r]^-{\beta}  \ar@{=}[dr] & X\otimes C \ar[d]^-{X\otimes \epsilon} \\
 & X
}$$
A right comodule $X$ over $C$ is the same thing as a left comodule over the opposite coalgebra $C^o$.
If $C$ and $D$ are graded coalgebras, then a $(C,D)$-{\it bicomodule} $X$ is a graded vector space $X$ equipped with a left coaction $\alpha:X\to C\otimes X$ and a right coaction $\beta:X\to X\otimes D$ which commutes, that is, for which the following square commutes, 
$$\xymatrix{
X \ar[rr]^-{\alpha} \ar[d]_{\beta} && C\otimes X \ar[d]^{C\otimes \beta } \\
X\otimes D \ar[rr]^-{\alpha \otimes D} && C\otimes X\otimes D.
}$$
When $C=D$ we shal say a $C$-bicomodule instead of a $C$-bicomodule.
A $(C,D)$-bicomodule is the same thing as a left comodule over $C\otimes D^o$, or as a right comodule over $C^o\otimes D$.
The coalgebra $C$ has the structure of a $C$-bicomodule.  Equivalently, it has the structure of a left comodule over the coalgebra $C\otimes C^o$.

Morphisms of bicomodules are equivalently defined as morphisms of left comodules over $C\otimes D^{o}$ or of right comodules over $C^o\otimes D$. The category of $(C,D)$-bicomodules is noted $\Bicomod(C,D)$. We have an equivalence $\Bicomod(C,D)=\Comod(C\otimes D^{o})$. If $C=D$, we shall put $\Bicomod(C):=\Bicomod(C,C)$.
As a particular case of a comodule category, $\Bicomod(C,D)$ is enriched, tensored and cotensored over $\dgVect$.

\medskip

\begin{defi}  If $X=(X,\alpha)$ is a left comodule over a dg-coalgebra $C$,
we shall say that a graded subspace $Y\subseteq X$ is a 
 {\it sub-comodule} if we have $\alpha (Y)\subseteq C\otimes Y$.
 The notion of sub-comodule of a right comodule (resp. bicomodule) is defined similarly.
\end{defi}

If $Y$ is a  sub-comodule of a (left) $C$-comodule $X=(X,\alpha)$, then the vector space $Y$ has the structure of a
(left) $C$-comodule with
 the coaction $Y\to C\otimes Y$ induced by the coaction
$\alpha:X\to C\otimes X$.
This turns the inclusion $Y\subseteq X$ 
into a morphism of comodules $Y\to X$.

\medskip

\begin{prop} [\cite{Sw}] \label{gradedCatComodfinitude0}
Let $X=(X,\alpha)$ be a left comodule over a dg-coalgebra $C$.
If $ev:C^\star \otimes  C\to \mathbb{F}$ is the evaluation,
then the composite
$$\xymatrix{
X\otimes C^\star \ar[r]^-{\sigma}   & C^\star \otimes X \ar[rr]^-{C^\star \otimes \alpha}  &&C^\star \otimes  C\otimes X \ar[rr]^-{ev \otimes X} && X
}$$
is a right action $\lfloor$ of the algebra $C^\star $ on $X$.
Thus, $X$ has the structure of a right module over $C^\star $.
\end{prop}

\begin{proof} By definition, we have
$$x\lfloor \phi =\phi( x^{(0)}) \otimes  x^{(1)}\ (-1)^{| x||\phi|} $$
for $\phi \in C^\star $ and $x\in X$.
The verification that this defines a right action of $C^\star $ on $X$ is 
straightforward.
\end{proof}

\begin{lemma} [\cite{Sw}] \label{subcomodule}
Let $X=(X,\alpha)$ be a left comodule over a dg-coalgebra $C$.
Then a dg-vector subspace $Y$ of $X$ is a sub-comodule iff $Y$ is a sub-module with respect to the right action $C^\star$.
Every dg-vector subspace $V\subseteq X$ is contained in a smallest sub-comodule $V\lfloor C^\star$.
The sub-comodule $V\lfloor C^\star$ is finite dimensional when $V$ is finite dimensional (in particular if $V$ is generated by a finite graded set).
\end{lemma}

\begin{proof} Every dg-sub-comodule of $X$ 
 is obviously closed under the (right) action of $C^\star$.
Conversely, if $Y$ is a dg-vector subspace of $X$ closed 
under the action of $C^\star$, let us show that $Y$ is a sub-comodule.
It is enough to check it at the level of the underlying graded objects.
Let us choose a graded linear basis $(u_i:i\in I)$ of the coalgebra $C$
and let $u^i:C\to \mathbb{F}$ be the linear form such that $u^i(u_j)=\delta^i_j$
(Kronecker delta) for every $j\in J$.
We shall use the decomposition
$$C\otimes X\simeq \FF I\otimes X \simeq \bigoplus_{i\in J} u_i\otimes X.$$
If $x\in X$, then $\alpha(x)=\sum_{i\in I}u_i\otimes x_i$,
where $x_i=0$ except for a finite number of elements $i\in I$.
Notice that $x\lfloor u^i= x_i(-1)^{|u^i||x|}$ for every $i\in I$.
If $x\in Y$ then $x_i=x\lfloor u^i(-1)^{|u^i||x|} \in Y$, since $Y$ is closed 
under the action of $C^\star$; thus $\alpha(x)\in C\otimes Y$.
This shows that $Y$ is a sub-comodule.
If $V\subseteq X$ is a graded subspace, then $V\lfloor C^\star$
is the sub-module of $X$ generated by $Y$.
It is thus the smallest sub-comodule of $X$ containing $V$.
If $V$ has finite dimension, let us show
that $V\lfloor C^\star$ has finite dimension.
 If $V=\mathbb{F}v_1+\cdots+\mathbb{F}v_n$, then $V \lfloor C^\star=v_1\lfloor C^\star+\cdots +v_n\lfloor C^\star$.
Hence it suffices to show that $v\lfloor C^\star$
has finite dimension for every $v\in X$.
But we have $\alpha(v)=\sum_{i\in F}u_i\otimes v_i$
for a finite set $F\subseteq I$.
Hence we have $v\lfloor \phi=\sum_{i\in F}\phi(u_i)v_i(-1)^{|v||\phi|}$
for every $\phi \in C^\star$.
This shows $v\lfloor C^\star$ is included in the linear span of the set $\{v_i\mid i\in F\}$. 
Thus, $v\lfloor C^\star$ is finite dimensional.

Finally, if $Y$ is a finite graded set of a dg-vector space $V$, the sub-dg-vector space generated by $Y$ is finite dimensional: it is generated as a graded vector space by the elements of $Y$ and their differentials.
\end{proof}

\begin{prop}[\cite{Sw}] \label{Comodunion}
A graded comodule is the directed union of its finite dimensional sub-comodules. 
\end{prop}

\begin{proof} Let $X$ be a left comodule over the coalgebra $C$.
If $S\subseteq X$ is a graded subset, let us denote by $\overline S$  the sub-comodule of $X$
generated by $S$. The comodule $X$ is obviously the directed
union of the sub-comodules $\overline S$, where $S$ runs in the finite
graded subsets of $X$.  This proves the result, since the sub-comodule
$\overline S$ is finite dimensional when $S$ is finite by lemma \ref{subcomodule}.
\end{proof}

Recall that an object $X$ in a cocomplete category $\cal C$ is said to be $\omega$-{\it compact} if the functor ${\cal C}(X,-)$ commute with directed limits. Recall also that a cocomplete category $\cal C$ is said to be $\omega$-{\it presentable} if  the category of $\omega$-compact objects of $\cal C$  is essentially small and every objet of $\cal C$ is a colimit of $\omega$-compact objects.

\bigskip

The {\em cotensor product} of a right $C$-comodule $X$ and a left $C$-comodule $Y$ is the object $X\otimes^CY$ defined as the coequalizer
$$\xymatrix{
X\otimes^CY \ar[r]& X\otimes Y \ar@<.6ex>[rr]^-{X\otimes \Delta_Y}\ar@<-.6ex>[rr]_-{\Delta_X\otimes Y}&& X\otimes C\otimes Y 
}$$
where $\Delta_X:A\otimes X\to X$ and $\Delta_Y:Y\otimes A\to Y$ are the $C$-comodule structures.

Let $f:D\to C$ be a dg-coalgebra map, any $D$-comodule $X$ can be seen as a $C$-comodule by the formula 
$\Delta(x) = x^{(1)}\otimes f(x^{(0)})$ for any $x\in X$. 
This produces a functor $\Comod(D)\to \Comod(C)$ called the {\em corestriction functor}. In particular $D$ can be viewed as a $C$-comodule.

Let $X$ be a left $C$-comodule, then the dg-vector spaces $D\otimes^CX$  is a $D$-comodules for the coaction
$(d\otimes x)^{(1)} \otimes (d\otimes x)^{(0)}:= d^{(1)} \otimes (d^{(1)}\otimes x)$.

\begin{prop}
The functor $D\otimes^C- : \Comod(C)\to \Comod(D)$ is right adjoint to the restriction functor.
\end{prop}
\begin{proof}
A straightforward computation left to the reader
\end{proof}

\subsubsection{Coderivations}

\begin{defi} \label{defcoderivcomodule}
Let $C$ be a dg-coalgebra or a non-counital dg-coalgebra and $N$ an $C$-bicomodule.
We shall say that a map of dg-vector spaces $D:N\to C$ is a {\it coderivation of $C$ with domain $N$} (or simply a {\em coderivation})
if we have the co-Leibniz rule
$$
\Delta D(x) = (D \otimes id)\rho(x) + (id\otimes D)\lambda(x)
$$
for every $x\in N$, were $\lambda:N\to C\otimes N$ and $\rho:N\to N\otimes C$ are respectively the left and right coaction of $C$ on $N$.
With Sweedler notation this condition reads
$$
D(x)^{(1)}\otimes D(x)^{(2)} = D(x_\rho^{(0)})\otimes x_\rho^{(1)} + x_\lambda^{(1)}\otimes D(x_\lambda^{(0)})
$$
where 
$x_\rho^{(0)}\otimes x_\rho^{(1)}\in N\otimes C$ 
and 
$x_\lambda^{(1)}\otimes x_\lambda^{(0)}\in C\otimes N $ 
are respectively the right and left coaction of $x$.
We shall note $coder(N,C)$ the set of derivations of $C$ with domain $N$.

\medskip
More generally, we shall say that a graded morphism of degree $n$ of graded vector spaces $D:N\rhup C$ is a {\it graded coderivation of $C$ with domain $N$ of degree $n$} (or simply a {\em graded coderivation})
if we have the graded co-Leibniz rule
$$
\Delta D(x) = (D \otimes id)\rho(x) + (id\otimes D)\lambda(x)
$$
for every $x\in N$.
With Sweedler notation this condition reads
$$
D(x)^{(1)}\otimes D(x)^{(2)} = D(x_\rho^{(0)})\otimes x_\rho^{(1)} + x_\lambda^{(1)}\otimes D(x_\lambda^{(0)})\ (-1)^{n|x_\lambda^{(1)}|}.
$$
We shall note $\Coder(N,C)_n$ the vector space of graded coderivations of $C$ with domain $N$ of degree $n$.
We have an inclusion $\Coder(N,C)_n\subset [N,C]_n$ and the graded vector subspace $\Coder(N,C)=\Coder(N,C)_*\subset [N,C]$ is stable by the differential of $[N,C]$ which enhance it into a dg-vector space. The differential of a graded coderivation $D:N\rhup C$ of degree $n$ is the map  $d_CD-Dd_N \ (-1)^{n}$ of degree $n-1$.
Coderivations correspond exactly to graded coderivations of degree 0 such that the morphism $d:N\to C$ commute with the differentials, equivalently they correspond to cycles in $\Coder(N,C)_0$.
\end{defi}

\begin{ex}\label{internalcoderivation}
If $C$ is a coalgebra, then the vector space $C$ has the structure of a bimodule over
the dual algebra $C^\star$ by proposition \ref{gradedCatComodfinitude0}.
The right action
$\lfloor :C \otimes C^\star \to C$ is defined by putting $x \lfloor \phi =\phi(x^{(1)}) x^{(2)}(-1)^{|\phi||x|}$
and the left action $\rfloor :C^\star \otimes C\to C$ by putting
$\phi \rfloor x =x^{(1)} \phi(x^{(2)})(-1)^{|\phi||x^{(1)}|}$.
Then the map $ [\phi,-]=C\to C$
defined by putting $[\phi,x]=\phi \rfloor x -x \lfloor \phi (-1)^{|x||\phi|}$
is a coderivation of degree $|\phi|$. The coderivation $[\phi,-]$ is said to be {\it inner}.
\end{ex}

\begin{rem}\label{remcodersuspension}
Let $S^n$ the graded vector space freely generated by one element $s^n$ of degree $n$.
If $N$ is a $C$-bicomodule, then $S^nN=S^n\otimes N$ is a $C$-bicomodule for the coactions
$$
\lambda(s^nx) = x_\lambda^{(1)}\otimes s^nx_\lambda^{(0)}\ (-1)^{n|x_\lambda^{(0)}|} \in C\otimes S^nN
$$
$$
\rho(s^nx) = s^nx_\rho^{(0)}\otimes x_\rho^{(1)} \in S^nN\otimes C.
$$
A morphism $D:N\rhup C$ of degree $n$ is a graded derivation iff the map $D':S^nN\to C$ 
defined by putting $D'(s^nx)=D(x)$ is a graded derivation (of degree 0) with domain $S^nN$.
\end{rem}

\begin{lemma} \label{smalllemmaforcoderivation2}
If $D:N\to C$ is a coderivation, then $\epsilon D=0$.
\end{lemma}

\begin{proof}
The counit  $\epsilon:N\to \FF$ is a map of coalgebras. 
Hence the morphism  $\epsilon D:N\rhup \FF$ is a coderivation $d:N\to \FF$ with values in the coalgebra $\FF$. 
The left coaction $\lambda:N\to \FF\otimes N$ is the canonical isomorphism $\lambda(x)=1\otimes x$ and the right coaction  $\rho:N\to N\otimes \FF$ is the canonical isomorphism $\rho(x)=x\otimes 1$.
For every $x\in N$ we have
$$
\Delta d(x) =1\otimes d(x)+d(x)\otimes 1.
$$
This proves that $d(x)=0$, since $\FF\otimes \FF=\FF$, $\Delta d(x)=d(x)=1\otimes d(x)=d(x)\otimes 1$.
\end{proof}

\begin{rem}\label{nucoder}
If $X$ be a graded vector space equipped with a set of co-operations $\Phi_i:X \to X^{\otimes k_i}$ of arity $k_i$, let us call a graded morphism $D:X\rhup X$ of degree $n$ a {\it coderivation} if, for every $x\in X$ and $i$, we have
$$
\Phi_i(Dx)=\sum_{j=1}^k  x^{(1)}\otimes \cdots \otimes  x^{(j-1)} \otimes Dx^{(j)} \otimes x^{(j+1)} \otimes \cdots \otimes x^{(k_i)}
(-1)^{n|x_1|+\cdots+n|x_{j-1}| },
$$
where $\Phi_i(x)= x^{(1)}\otimes \cdots \otimes x^{(k_i)}$ in Sweedler's notation.
This definition applies in particular to non-counital coalgebras $(C,\Delta)$ and coalgebras $(C,\Delta,\epsilon)$. 
Lemma \ref{smalllemmaforcoderivation2} says that the extra condition of the latter case is superfluous.
In consequence, for a counital coalgebra $C$, we shall use the same notation $\Coder(C)$ to refer to coderivations of $C$ viewed as counital or not.
\end{rem}

\bigskip

If $D:N\to C$ is a derivation and $M\to N$ is a map of $C$-bicomodules, the composite $Df:M\to C$ is still a coderivation.
This defines a functor
$$\xymatrix{
\Coder(-,C):\Bicomod(C)^{op}\ar[r] & \Set.
}$$
We shall say that a coderivation $d:\Omega^C\to C$ is couniversal if it represents the functor $coder(-,C)$.
Equivalently, $d:\Omega^C\to C$ is couniversal if for any coderivation $D:N\to C$, there exists a unique bicomodule map $f:N\to \Omega^C$ such that $df=D$
$$\xymatrix{
&&\Omega^C\ar[d]^{d} \\
N\ar[rr]^-D\ar@{-->}[rru]^f&& C.
}$$

\begin{prop}\label{univcoderiv}
For any dg-coalgebra $C$ (counital or not), there exists a couniversal coderivation $d:\Omega^C\to C$.
\end{prop}

\begin{proof} 
The construction is dual to proposition \ref{univderiv}. Using remark \ref{remderuniv}, $\Omega^C$ is the equalizer in $C$-bicomodules of the diagram
$$\xymatrix{
C\otimes C\otimes C \ar@<.6ex>[rrrr]^-{C \otimes \Delta_C \otimes C } \ar@<-.6ex>[rrrr]_-{\Delta_C \otimes C\otimes C + C\otimes C\otimes \Delta_C} &&&& C\otimes C\otimes C\otimes C.
}$$
The map to $d:\Omega^C\to C$ is induced by $\epsilon_C\otimes C\otimes \epsilon_C: C\otimes C\otimes C\to C$.
\end{proof}

We shall call the bicomodule $\Omega^C$ the {\em bicomodule of codifferentials of $C$} or the {\em tangent space of $C$}.

\begin{rem}\label{gradedunivcoder}
An analogous functor can be defined in the graded context if $\Coder(N,C)$ is replaced by $\Coder(N,C)_0$.
If $C$ is a dg-coalgebra (counital or not), the same proof shows that, $d:\Omega^C\to C$ viewed as a graded map, is still the universal graded coderivation.
In other words, we have a natural isomorphism $\Coder(N,C)_0=\Hom_{C,C}(N,\Omega^C)_0$, 
where, for two $C$-bicomodules $M$ and $N$, $\Hom_{C,C}(N,M)$ is the enriched dg-vector space hom in $\Bicomod(C)$.
\end{rem}

\begin{rem}\label{remcoderuniv}
As in remark \ref{remderuniv}, the previous diagram is part of a full cosimplicial diagram of $A$-bimodules:
$$\xymatrix{
C\ar[r] &C\otimes C  \ar@<.6ex>[r] \ar@<-.6ex>[r] & C\otimes C\otimes C \ar@<1ex>[r] \ar[r]\ar@<-1ex>[r] & C\otimes C\otimes C\otimes C & \dots
}$$
whose associated complex is the Hochschild complex of $C$.
Similarly to the algebra case, this cosimplicial diagram admits a "contracting homotopy", making the Hochschild complex of $C$ an exact complex.
In consequence, $\Omega^C$ can also be constructed as the cokernel of $\Delta:C\to C\otimes C$ in the category of $C$-bicomodules
and the universal coderivation $d:\Omega^C\to C$ is induced by the map $\epsilon\otimes C-C\otimes \epsilon : C\otimes C\to C$.
\end{rem}

We can proceed similarly with graded coderivations. They define a functor
$$\xymatrix{
\Coder(-,C):\Bicomod(C)^{op}\ar[r] & \dgVect.
}$$
We shall say that a graded coderivation is universal if it represents the functor $\Coder(-,C)$.
We have the following stronger form of proposition \ref{univcoderiv}.

\begin{prop}\label{univcoderivgraded}
For any dg-coalgebra $C$ (counital or not), $d:\Omega^C\to C$, viewed as a graded coderivation, represents the functor $\Coder(-,C)$.
Equivalently, we have a natural isomorphism $\Coder(N,C)=\Hom_{C,C}(N,\Omega^C)$.
\end{prop}
\begin{proof}
Dual to proposition \ref{univderivgraded}.
\end{proof}

\bigskip

Let $f:D\to C$ be a map of dg-coalgebras (counital or not) and $N$ a $D$-bicomodule, then $M$ can be viewed as a $C$-bicomodule with coactions 
$$\xymatrix{
N\ar[r]^-\lambda & D\otimes N\ar[rr]^-{f\otimes N} && C\otimes N
}$$
$$\xymatrix{
N\ar[r]^-\rho & N\otimes D\ar[rr]^-{N\otimes f} && N\otimes C.
}$$
We shall note this bicomodule $N^f$.

\begin{defi}\label{fcoderivation} \label{defcoderivationcoalg}

If $f:D\to C$ is a map of dg-coalgebras  (counital or not) and $N$ a $D$-bicomodule, 
we shall say that a morphism $D:N\rhup C$ of degree $n$ is an $f$-{\it graded coderivation with domain $N$}  (or a {\em graded $f$-coderivation}) if it is a graded coderivation of $C$ with values in the bicomodule $N^f$. 
The $f$-coderivations form a dg-vector space $\Coder(N;C,f):=\Coder(N^f,C)$.

In case $N=D$ with the canonical bimodule structure, we shall put $\Coder(f):=\Coder(D;C,f)$.
In case $f=id_C$ and $N=C$, we shall put $\Coder(C):=\Coder(id_C)$.
\end{defi}

Concretely, a morphism $D:N\rhup C$ of degree $n$ is a $f$-coderivation if we have
$$
D(x)^{(1)}\otimes D(x)^{(2)} = D(x_\rho^{(0)})\otimes f(x_\rho^{(1)}) + f(x_\lambda^{(1)})\otimes D(x_\lambda^{(0)}) \ (-1)^{n|x_\lambda^{(1)}|}
$$
where 
$x_\rho^{(0)}\otimes x_\rho^{(1)}\in N\otimes C$ 
and 
$x_\lambda^{(1)}\otimes x_\lambda^{(0)}\in C\otimes N $ 
are respectively the right and left coaction of $x$.

\bigskip

$f$-coderivations define a functor
$$\xymatrix{
\Coder(-;C,f):\Bimod(B)\ar[r] & \dgVect.
}$$
We shall say that a graded $f$-coderivation $d:N\rhup C$ is couniversal if it represents the functor $\Coder(-;C,f)$.

Let $\Omega^C$ be the bicomodule of codifferentials of $C$, and let $\Omega^{C,f}:=\Omega^C\otimes^{C^o\otimes C}D^o\otimes D$ be the $D$-bicomodule obtained by base change along $f$. There is a canonical $C$-bicomodule map $\Omega^{C,f}\to \Omega^C$ and an $f$-coderivation $\Omega^{C,f}\to C$.
We have the following variation of proposition \ref{univderivgraded}.

\begin{prop}\label{univfcoderiv}
The $f$-coderivation $\Omega^{C,f}$ is couniversal. 
Equivalently, we have a natural isomorphism $\Coder(N;C,f)=\Hom_{D,D}(N,\Omega^{C,f})$.
\end{prop}
\begin{proof}
Let $N\to C$ be an $f$-derivation, by proposition \ref{univderivgraded} we have $\Coder(N;C,f)=\Coder(N^f,C)=\Hom_{C,C}(N^f,\Omega^C)$
and by base change, we have $\Hom_{C,D}(N^f,\Omega^C) = \Hom_{D,D}(N,\Omega^C\otimes^{C^o\otimes C}D^o\otimes D)$.
\end{proof}

We shall call $\Omega^{C,f}$ the {\em tangent space of $C$ at $f$}. We shall see in proposition \ref{prim=coder} another characterization in the case $D=\FF$.

\bigskip

If $C$ is a dg-coalgebra (counital or not) and $N$ is a bicomodule, then the direct sum $C\oplus N$ has the structure of a dg-coalgebra with the coproduct defined by putting 
$$
\Delta(x,y)=\Delta(x)+\rho(y)+\lambda(y)
$$
for $(x,y)\in C\oplus N$,
where $\lambda:N\to C\otimes N$ and $\rho:N\to N\otimes C$ denote respectively the left and the right coaction of $C$ on $N$.
The counit $\epsilon: C\oplus N \to \FF$ is defined by putting $\epsilon(x,y)=\epsilon(x)$.
The inclusion $i_1:C\to C\oplus N$ and projection $p_1:C\oplus N \to C$ are maps of coalgebras and we have $p_1i_1=1_C$,
$$\xymatrix{
C\ar[r]^-{i_1} & C\oplus N  \ar[r]^-{p_1} & C.
}$$

\begin{prop}\label{fcoderiv}
Let $f:D\to C$ be a map of dg-coalgebras (counital or not) and $N$ be a $D$-bicomodule.
Then a map $D:N\to C$ is a graded $f$-coderivation (of degree $0$) iff the map $g: C\oplus N \to C$
defined by putting $g(x,y)=f(x)+D(y)$ is a map of graded coalgebras.
$D$ is a $f$-derivation iff the map $g$ is a map of dg-coalgebras.

In other words, we have a cartesian square in $\gSet$
$$\xymatrix{
\Coder(N;C,f) \ar[r]^-g\ar[d]& [D\oplus N,C]\ar[d]^{[i_1,C]}\\
\{*\}\ar[r]^-f & [D,C].
}$$
\end{prop}

\begin{proof}
Dual to proposition \ref{fderiv}.
\end{proof}

\begin{defi}\label{pointedcoderivationdef}
If $C=(C,e_C)$ and $(D,e_D)$ are two pointed coalgebras, and $f:D\to C$ is a map of pointed coalgebras,
we shall say that an $f$-coderivation $D:D\rhup C$ is {\it pointed} if $De_C=0$.
Pointed $f$-derivations form a dg-vector subspace $\Coder_\bullet(f)\subset \Coder(f)$. 
If $C=D$ and $f=id_C$ we shall put $\Coder_\bullet(C)$ instead of $\Coder_\bullet(id_C)$.
\end{defi}

Every pointed derivation $D:D\rhup C$ induces a $f_-$-coderivation $D_{-}:D_-\rhup C_-$ of the corresponding non-unital coalgebras.
Using lemma \ref{smalllemmaforcoderivation2}, it is easy to verify that the map $D\mapsto D_{-}$
is an isomorphism of vector spaces $\Coder_\bullet(f)\simeq \Coder(f_{-})$.

\bigskip

If $N$ is a bicomodule over a pointed coalgebra $C=(C,e)$, then the reduced coactions

$$\overline{\lambda}: N\to C_{-}\otimes N \et  \overline{\rho} : N\to N \otimes C_{-}$$
are defined by putting 
$\overline{\lambda}(y)=\lambda(y)-e\otimes y$ and $\overline{\rho}(y)=\rho(y)-y\otimes e$.
The coalgebra $ C\oplus N $ is pointed by $(e,0)\in  C\oplus N$
and the reduced coproduct
$\overline{\Delta}$ of the non-counital coalgebra $  C\oplus N_{-}=C_{-}\oplus N$
is defined by putting
\begin{eqnarray*}
\overline{\Delta}(x+y) &= & \Delta (x) +\lambda(y)+\rho(y)  - e\otimes x-x\otimes e-e\otimes y-y\otimes e \\
&=&\overline{\Delta}(x) +\overline{\lambda}(y)+\overline{\rho}(y) 
\end{eqnarray*}
for $x+y=(x,y)\in C_{-}\oplus N$, 
where $\overline{\Delta}(x)=\Delta(x)-e\otimes x-x\otimes e$.
For every $(x,y)\in C_{-}\oplus N$ we have
$$\overline{\Delta}^{n}(x,y)=\overline{\Delta}^{n}(x)+\sum_{i+1+j=n}\overline{\Delta}^{(i,j)}(y),$$
where the map $\overline{\Delta}^{(i,j)}:N\to C^{\otimes i}_{-}\otimes N\otimes C^{\otimes j}_{-}$
is the diagonal of the following commutative square
$$\xymatrix{
N\ar[rr]^-{\overline{\lambda}^i} \ar[d]_{\overline{\rho}^{j}}  &&  C^{\otimes i}_{-} \otimes N \ar[d]^{C^{\otimes i}_{-}\otimes \overline{\rho}^{j}}\\
N \otimes C^{\otimes j}_{-}  \ar[rr]^-{\overline{\lambda}^{i}\otimes  C^{\otimes j} }  && C^{\otimes i}_{-} \otimes N \otimes C^{\otimes j}_{-}.
}
$$

\begin{lemma}\label{nilpotentextension} If the pointed coalgebra $C$ is conilpotent ,
then so is the pointed coalgebra $ C\oplus N $. 
 \end{lemma}

\begin{proof}  Let us show that for every $(x,y)\in C_{-}\oplus N$ we have
$\overline{\Delta}^{n}(x,y)=0$ for $n\gg 0$.
But if $\overline{\lambda}(y)=y^{(1)}\otimes y^{(0)}$, then we have $\overline{\lambda}^{i}(y)=\overline{\Delta}^{i-1}(y^{(1)})\otimes y^{(0)}$
for every $i>1$ by coassociativity.  Hence we have $\overline{\lambda}^{i}(y)=0$ when $i\gg 0$, since the pointed coalgebra $C$ is conilpotent.
Similarly, we have $\overline{\rho}^{i}(y)=0$ when $i\gg 0$.
 It then follows from the formula above that we have $\overline{\Delta}^{n}(x,y)=0$ for $n\gg 0$.
\end{proof}

\medskip

\begin{prop}\label{nilradicalcoder} The
radical $R^cC$
of a pointed coalgebra $(C,e)$ is preserved 
by every graded coderivation $D:C\rhup C$, pointed or not.
\end{prop}

\begin{proof}  If $D:C\rhup C$ is a coderivation of degree $n$, then
the map $f: C\oplus s^n C\to C$ defined by putting $f(x,s^ny)=x+D(y)$
is a map of pointed coalgebras by lemma \ref{fcoderiv} since $f(e,0)=e$.
Let us put $E= R^cC$.
The coalgebra $ E\oplus s^n E$ is conilpotent
by lemma \ref{nilpotentextension} and it is a 
subcoalgebra of $  C\oplus s^n C$,
since $E$ is a subcoalgebra of $C$.
Hence we have $f(E)\subseteq R^cC=E$ by proposition \ref{radical}.
Thus $D(y)=f(0,s^ny)\in E$ for every $y\in E$.
\end{proof}

By the equivalence between pointed coalgebras and non-counital coalgebras, we have the following consequence.
\begin{cor}\label{nilradicalcodernu}
The radical $R^cC$ of a non-counital coalgebra $C$ is preserved by every graded coderivation $D:C\rhup C$.
\end{cor}

\medskip

Let us denote by  $p:T^c(X)\to X$ the cogenerating map
of the cofree conilpotent coalgebra on $X$.

\begin{lemma} \label{coextensioncoderconillemma}
 If $X$ is a graded vector space and $N$ is a $T^c(X)$
 bicomodule, then every morphism $\phi:N\rhup X$ of degree $n$
 can be coxtended uniquely as a coderivation  $D:N \rhup T^c(X)$
 of degree $n$.
 \end{lemma}

\begin{proof} 
 Let us first suppose that $n=0$.
Let $f: T^c(X) \oplus N  \to X$ be the map
defined by putting $f(x,y)=p(x)+\phi(y)$ for every $(x,y)\in T^c(X) \oplus N$.
The coalgebra $  T^c(X) \oplus N $ is conilpotent by lemma \ref{nilpotentextension}, since the coalgebra
$T^c(X) $ is conilpotent.  Moreover, we have $f(1,0)=p(1)=0$.
Hence there exists
a unique map of pointed coalgebras $g:T^c(X) \oplus N \to T^c(X)$
such that $pg=f$ by proposition \ref{structureconil}.
If  $i_1:T^c(X) \to  T^c(X) \oplus N $
is the inclusion, then we have $pgi_1=fi_1=p=p(id)$.
But $gi_1: T^c(X) \to T^c(X)$ is a coalgebra map,
since $i_1$ is a coalgebra map. Thus $gi_1=id$, since $p$
is cogenerating. It then follows from lemma
\ref{fcoderiv} that we have $g(x,y)=x+D(y)$, where $D:N\to T^c(X)$ is a coderivation of degree 0.
We have $pD=\phi$, since we have $pg=f$. The uniquness of $D$
is left to the reader. 
Let us now consider the case of a morphism $\phi:N\rhup X$ of general degree $n$.
In this case, the morphism $\phi s^{-n}:S^{n}N\rhup X$ defined by putting
$\phi s^{-n}(s^n x)=\phi(x)$ for $x\in N$ has degree 0.  It can thus be coextended uniquely as a 
coderivation $Ds^{-n}:S^{n}N \to T^c(X) $ of degree 0.
The resulting morphism $D:N\rhup T^c(X) $ is a coderivation of degree $n$
which is extending $\phi$. The uniqueness of $D$ is clear.
 \end{proof}

\begin{ex}\label{unpointedcoder} If $ X$ is a graded vector space and $v\in X_n$, then the map  $D:T^c(X) \rhup T^c(X)$ defined by putting 
$$
D(x_1\otimes \cdots \otimes x_k)=\sum_{i=0}^{k} x_1\otimes \cdots \otimes x_{i} \otimes v \otimes x_{i+1}\otimes \cdots \otimes x_k(-1)^{n(|x_1|+\cdots + |x_{i}|)}
$$
for every $x_1\otimes \cdots \otimes x_k\in X^{\otimes k}$
is a coderivation of degree $n$.
\end{ex}

The vector space $ T^c(X)\otimes X\otimes T^c(X)$ has the structure of a bicomodule
over the coalgebra $T^c(X)$. The bicomodule is cofreely cogenerated
by the map  $\phi=\epsilon \otimes X\otimes \epsilon: T^c(X)\otimes X\otimes T^cX)\to X$.
The map $\phi$ can be coextended uniquely as a coderivation $d:T^c(X)\otimes X\otimes T^c(X)\to T^c(X)$
by lemma \ref{coextensioncoderconillemma}.

\begin{prop} \label{univcodercofreeconil} 
The coderivation $d:T^c(X)\otimes X\otimes T^c(X)\to T^c(X)$  defined above is couniversal.
Hence we have $$ \Omega^{T^c(X)}= T^c(X)\otimes X\otimes T^c(X).$$
\end{prop}

\begin{proof} 
Let $N$ be a bicomodule over the coalgebra $T^c(X)$
and  $D:N\to T^c(X)$ be a coderivation of degree 0.
There is then a unique map of bicomodules $f: N\to T^c(X)\otimes X\otimes T^c(X)$
such that $\phi f=pD$, since the bicomodule $T^c(X)\otimes X\otimes T^c(X)$
is cofreely cogenerated by $\phi$.
Let us show that $df=D$.
But $df:N\to T^c(X)$ is a coderivation, since $d$ is a coderivation and $f$ is a bicomodule map.
Moreover, we have $pdf=\phi f=pD$ by definition of $f$.
It then follows by lemma \ref{coextensioncoderconillemma} that $df=D$.
\end{proof}

Using propositions \ref{univcodercofreeconil} and \ref{univfcoderiv}, we can improve lemma \ref{coextensioncoderconillemma}.
\begin{cor}\label{caraccoderconil}
Let $p:T^c(X)\to X$ be the cogenerating map, then the map $D\mapsto pD$ induces an isomorphism of vector spaces
$$
\Coder(T^c(X)) = [T^c(X),X].
$$
In particular, a derivation $D$ is zero iff $pD=0$.
\end{cor}

\bigskip

\begin{prop} \label{coextensioncoderconil}
 If $X$ is a graded vector space, then every morphism $\phi:T^c(X)\rhup X$ 
 of degree $n$ can be coextended uniquely as a 
 coderivation  $D:T^c(X) \rhup T^c(X)$
 of degree $n$. Moreover, we have
  $$D(x_1\otimes \cdots \otimes x_k)=\sum_{0\leq i\leq j\leq k} x_1\otimes \cdots \otimes x_{i} \otimes \phi(x_{i+1}\otimes \cdots \otimes x_{j})\otimes x_{j+1}\otimes \cdots \otimes x_k(-1)^{n(|x_1|+\cdots + |x_{i}|)}$$
  for every $x_1\otimes \cdots \otimes x_k\in X^{\otimes k}$.
   \end{prop}

\begin{proof} The existence and uniqueness follows from lemma \ref{coextensioncoderconillemma}.
Recall that the reduced right coaction of the non-counital coalgebra $T^c_\circ(X)=T^c(X)_{-}$ on the comodule $N=S^n T^c(X)$
is given by
the formula
$${\rho}_{-} (s^nx_1\otimes \cdots \otimes x_k)=\sum_{i=0}^{k-1}s^n(x_1\otimes \cdots \otimes x_i)\otimes (x_{i+1}\otimes \cdots \otimes x_k)$$
and the reduced left coaction  by
the formula
$${\lambda}_{-} (s^nx_1\otimes \cdots \otimes x_k)=\sum_{i=1}^k (x_1\otimes \cdots \otimes x_i)\otimes s^n ( x_{i+1}\otimes \cdots \otimes x_k)(-1)^{n(|x_1|+\cdots + |x_{i}|)}.$$
The reduced coproduct ${\Delta}_{-} $ of the non-counital coalgebra $  C\oplus N_{-}=C_{-}\oplus N$
is given by the formula
$${\Delta}_{-} (y+s^nx) ={\Delta}_{-}(y) +{\lambda_{-} }(x)+{\rho}_{-} (x) $$
for $y+s^nx=(y,s^nx)\in C_{-}\oplus N$, and the iterated reduced coproduct ${\Delta}^r$ 
by the formula 
$${\Delta}^{r}_{-} (y+s^nx)={\Delta}^{r}_{-} (y) +\sum_{i+1+j=r}{\Delta}_{-}^{(i,j)}(s^nx)$$
where ${\Delta}^{(i,j)}_{-} :N\to C_{-}^{\otimes i}\otimes N\otimes  C_{-}^{\otimes j}$
is defined by putting
$${\Delta}^{(i,j)}_{-} (s^nx)=({\lambda}_{-}^i \otimes id){\rho}_{-}^j(s^n x).$$
The map $f:  T^c(X) \oplus N  \to T^c(X)$
 defined by putting $f(y,s^nx)=y+D(x)$ is a coalgebra map by lemma \ref{fcoderiv},
 since the map $Ds^{-n}:N\to  T^c(X)$ is a coderivation of degree 0.
The coalgebra  $T^c(X) \oplus N$
is pointed and conilpotent by lemma \ref{nilpotentextension}.
Moreover, the map $f$ is pointed,
since $f(1,0)=1$. It then follows from lemma \ref{structureconil}
that  we have
 $$f(y+s^nx)=\sum_{r\geq 1} (pf)^{\otimes r}{\Delta}_{-}^{(r)}(y,s^nx)$$
for every $y+s^nx=(y,s^nx)\in T^c(X)_{-} \oplus N$.
Hence we have
$$D(x)=\sum_{r\geq 1}(pf)^{\otimes r} {\Delta}_{-}^{(r)}(0,s^nx)$$
It follows that
\begin{eqnarray*}
D(x) &= & \sum_{i,j\geq 0} ((pf)^{\otimes i}\otimes pf \otimes (pf)^{\otimes i}){\Delta}_{-}^{(i,j)}(s^nx) \\
&=&\sum_{i,j\geq 0} (p^{\otimes i}\otimes \phi s^{-n} \otimes p^{\otimes i}){\Delta}_{-}^{(i,j)}(s^nx)
\end{eqnarray*}
since $pf(y,0)=p(y)$ and $pf(0, s^n x)=pD(x)=\phi s^{-n}(s^n x)$.
But we have $(p^{\otimes i}\otimes \phi s^{-n}  \otimes p^{\otimes j}){\Delta}_{-}^{(i,j)}(s^nx)=0$,
unless $i+j\leq k$, in which case 
$$(p^{\otimes i}\otimes \phi s^{-n}  \otimes p^{\otimes j}){\Delta}_{-}^{(i,j)}(s^nx)=
x_1\otimes \cdots \otimes x_i \otimes \phi(x_{i+1}\otimes \cdots \otimes x_{k-j}) \otimes x_{k-j+1}\otimes \cdots \otimes x_k
(-1)^{n(|x_1|+\cdots + |x_{i}|)}.$$
Thus, 
$$D(x_1\otimes \cdots \otimes x_k)= \sum_{i+j\leq k} x_1\otimes \cdots \otimes x_{i} \otimes \phi(x_{i+1}\otimes \cdots \otimes x_{k-j})\otimes x_{k-j+1}\otimes \cdots \otimes x_k(-1)^{n(|x_1|+\cdots + |x_{i}|)}.$$
The formula to be proved is  obtained by reindexing this sum.
\end{proof}

The result above has a pointed version, which we state without proof:

\begin{prop} \label{coextensioncoderconilpointed}
If $X$ is a graded vector space, then every morphism $\phi:T^c_\circ(X)\rhup X$ of degree $n$ can be coextended uniquely as a pointed coderivation $D:T^c_\bullet(X) \rhup T^c_\bullet(X)$ of degree $n$ given by the same formula.
\end{prop}

\subsubsection{Primitive elements of coalgebras}

\begin{defi}\label{defprim}
For $C$ a dg-coalgebra of a graded coalgebra,
we shall say that an element $x$ of a coalgebra $C$ is {\it primitive with respect to an atom} $e\in C$, 
if we have $\Delta(x)=e \otimes x + x \otimes e $.
We shall say that an element of a pointed coalgebra $(C,e)$ is {\it primitive} of it is primitive with respect to $e$. 
\end{defi}

\begin{lemma}\label{primaugmentation} 
If $x\in C=(C,e)$ is primitive, then $\epsilon(x)=0$.
\end{lemma}

\begin{proof} 
If $x$ is primitive with respect to $e\in C$, then we have
$$
x=(\epsilon\otimes C)\Delta(x)=(\epsilon\otimes C)(e\otimes x+x\otimes e)=x+\epsilon(x)e
$$
since $\epsilon(e)=1$.
Thus, $\epsilon(x)e=0$ and hence $\epsilon(x)=0$, since $\epsilon(e)=1$.
\end{proof}

\begin{lemma}\label{lemmaprim}
The graded vector space of primitive elements of a pointed coalgebra $C=(C,e)$ is the kernel of the reduced coproduct $\Delta_{-}:C_{-} \to C_{-}\otimes C_{-}$.
Moreover, if $C$ is a dg-coalgebra, it is stable by the differential.
\end{lemma}

\begin{proof}
If $x\in \Prim(C,e)$, then $x\in C_{-}$, since $\epsilon(x)=0$ by  \ref{primaugmentation}.
But we have 
$$
\Delta(x)=e\otimes x+x\otimes e+\Delta_{-}(x),
$$
since $x\in C_{-}$.
Thus $x$ is primitive if and only if $\Delta_{-}(x)=0$.

Then the stability if $\ker \Delta_-$ by the differential is a consequence of the commutation of $\Delta_-$ with the differential.
\end{proof}

We shall say that an element of a non-counital coalgebra $(C,\Delta)$ is {\em primitive} if $\Delta(x)=0$.
It is clear from the lemma that if $D$ is a pointed coalgebra, the primitive elements of $D$ are in bijection with those of $D_-$.
In particular if $C$ is a non-counital coalgebra, the primitive elements of $C_+$ are in bijection with those of $C$.

\begin{ex}\label{primfreecomil}
For $X$ a dg-vector space, let $(T^c_\bullet(X),1)$ be the tensor coalgebra of example \ref{tensorcoalg}, 
then every element of $X$ is primitive. The dg-vector space of primitive elements is isomorphic to $X$.
If $T^c_\circ(X)$ be the non-counital tensor coalgebra of example \ref{ncutensorcoalg}, every element of $X$ is also primitive and
the dg-vector space of primitive elements is again isomorphic to $X$.
\end{ex}

\begin{ex}\label{primcoalg}
For $X$ a dg-vector space, let $T^c_{\circ,1}(X)$ be the non-counital primitive coalgebra of example \ref{primitivecoalgebra}, then every element of $X=T^c_{\circ,1}(X)$ is primitive (hence the name chosen for $T^c_{\circ,1}(X)$). If $(T^c_{\bullet,1}(X),1)$ is the pointed primitive coalgebra, every element of $X=(T^c_{\bullet,1}(X))_-$ is also primitive. In both case the dg-vector space of primitive elements is isomorphic to $X$.
\end{ex}

\begin{ex}\label{freelieasprimitive}
For $X$ a dg-vector space, let $T^{csh}(X)$ be the coshuffle coalgebra of example \ref{coshufflecoalgebra},
then the dg-vector subspace of primitive elements of $T^{csh}(X)$ is the subspace $L(X)$ generated by $X$ under the Lie bracket.
$L(X)$ is the free dg-Lie algebra generated by $X$ \cite{Reutenauer}.
\end{ex}

\medskip

We shall denote the graded vector space of primitive elements of $C$ with respect to an atom $e\in C$  by $\Prim(C,e)$ or more simply by $\Prim(C)$ if the pointing is clear or if $C$ is non-counital. 
If $C$ is a dg-coalgebra, the restriction of $d_C$ makes $\Prim(C,e)$ into a dg-vector space by lemma \ref{lemmaprim}.
In particular $\Prim$ defines a functor
$$\xymatrix{
\Prim: \dgCoalg_\bullet \ar[r] &\dgVect.
}$$
Let us also consider the functor $T^c_{\bullet,1}:\dgVect \to \dgCoalg_\bullet$ of example \ref{primcoalg}.

\begin{prop}\label{adjointPrim}
There is an adjunction
$$\xymatrix{
T^c_{\bullet,1}:\dgVect \ar@<.6ex>[r]& \dgCoalg_\bullet :\Prim. \ar@<.6ex>[l]
}$$
\end{prop}

\begin{proof}
Coalgebra maps $X\to C_-$ sends the elements of $X$ in $\Prim(C) =\ker\Delta_-$.
This produces a bijection between non-counital coalgebra maps $X=T^c_{\circ,1}(X)\to C_-$ and linear maps $X\to \Prim(C)$.
The conclusion follows form the natural bijection between non-counital coalgebra maps $X=T^c_{\circ,1}(X)\to C_-$ and pointed coalgebra maps $T^c_{\bullet,1}(X)\to C$.
\end{proof}

\begin{cor}\label{univpropofdeltanew}
If $(C,e)$ is a pointed dg-coalgebra and $X$ a dg-vector space, then there are natural bijections between 
\begin{enumerate}
\item maps of dg-vector spaces $X\to \Prim(C,e)$,
\item maps of pointed dg-coalgebras $T^c_{\bullet,1}(X)\to (C,e)$,
\item maps of $C$-bicomodules $X\to \Omega^C$ (where $X$ is viewed as a $C$-bicomodule through $e:\FF\to C$),
\item and maps of dg-vector spaces $X\to \Omega^{C,e}$.
\end{enumerate}
\end{cor}

\begin{proof}
The bijection $1\leftrightarrow 2$ is given by proposition \ref{adjointPrim}.
The bijection $2\leftrightarrow 3$ is given by the universal property of $\Omega^C$ applied to $D\oplus N = \FF\oplus X$.
The bijection $3\leftrightarrow 4$ is proposition \ref{univfcoderiv}.
\end{proof}

We now prove that the primitive dg-vector space is the tangent complex.

\begin{cor}\label{prim=coder}
Let $(C,e)$ be a pointed dg-coalgebra, then we have a canonical isomorphism of dg-vector spaces
$$
\Prim(C,e) = \Omega^{C,e}.
$$
In other words, primitive elements of degree $n$ are the same thing as graded derivations of degree $n$ with values in the bicomodule given by the atom $e:\FF\to C$.
\end{cor}
\begin{proof}
From the equivalence $1\leftrightarrow 4$ of corollary \ref{univpropofdeltanew} both objects have the same functor of points.
They are isomorphic by Yoneda's lemma.
The last statement is a consequence but can be proven independently from corollary \ref{univpropofdeltanew} by taking $X=\FF\delta\to \FF d\delta$ where $\delta$ is in degree $n$ and $d\delta$ is the differential of $\delta$.
\end{proof}

\medskip

We finish by computing the primitive elements of a tensor product of coalgebras.

\begin{prop}\label{primitiveoftensor}
Let $C=(C,e)$ and $D=(D,u)$ be two pointed dg-coalgebras.
If $x\in C$ and $y\in D$ are primitive elements (of the same degree),
then the element $i(x,y)=x\otimes u+e\otimes y$ is primitive in $C\otimes D=(C\otimes D,e\otimes u)$.
Moreover, the map 
$$
i:\Prim(C)\oplus \Prim(D)\to \Prim(C\otimes D)
$$
is an isomorphism of dg-vector spaces.
\end{prop}

\begin{proof}
If $x\in C$ is primitive, then $i_1(x)=x\otimes u\in C\otimes D$
is primitive, since $C\otimes u:C\to C\otimes D$ is a map of pointed coalgebras.
Similarly,  if $y\in D$ is primitive, then $i_2(y)=e\otimes y\in C\otimes D$
is primitive.  It follows that $i(x,y)=x\otimes u+e\otimes y$  is primitive.
If $z\in C\otimes D$ is primitive, then $p_1(z)=(C\otimes \epsilon)(z)\in C\otimes D$
is primitive, since $C\otimes \epsilon: C\otimes D\to C$ is a map of pointed coalgebras.
Similarly,  $p_2(z)=(\epsilon \otimes D)(z)\in C\otimes D$ is primitive.
Obviously, $p_1i_1(x)=x$ and $p_2i_2(y)=y$.
Moreover, $p_2i_1(x)=0$ and $p_1i_2(y)=0$, since $x\in C_{-}$ and $y\in D_{-}$.
Hence the map $(p_1,p_2) : \Prim(C\otimes D)\to \Prim(C)\oplus \Prim(D)$
is a retraction of the map $i : \Prim(C)\oplus \Prim(D)\to \Prim(C\otimes D)$.
The result will be proved if we show that the map $(p_1,p_2)$
is monic. We have
$$\ker(p_1,p_2)\subseteq \ker(C\otimes \epsilon)\cap \ker(\epsilon \otimes D)=(C\otimes D_{-})\cap (C_{-}\otimes D)=C_{-}\otimes D_{-}.$$
Hence it suffices to show that $C_{-}\otimes D_{-}\cap \Prim(C\otimes D)=0$.
We saw above that $\Prim(C\otimes D)=\ker (\Delta_{-})$.
For every $x\otimes y\in  C_{-}\otimes D_{-}$ we have
$$\Delta_{-}(x\otimes y)=(x^{(1)}\otimes y^{(1)})\otimes (x^{(2)}\otimes y^{(2)}) -(e\otimes u)\otimes (x\otimes y)-(x\otimes y)
\otimes(e\otimes u).$$
Thus, 
$$(C\otimes \epsilon\otimes \epsilon \otimes D)\Delta_{-}(x\otimes y)=
\epsilon(x^{(1)}) y^{(1)})\otimes x^{(2)}\epsilon(y^{(2)}))= y^{(1)}) \epsilon(y^{(2)})) \otimes  \epsilon(x^{(1)}) x^{(2)}=x\otimes y.$$
It follows that 
$(C\otimes \epsilon\otimes \epsilon \otimes D)\Delta_{-}(z)=z$
for every $z\in C_{-}\otimes D_{-}$. Hence the map $\Delta_{-}$
is monic in $C_{-}\otimes D_{-}$. It follows that
$$
(C_{-}\otimes D_{-})\cap \Prim(C\otimes D)= (C_{-}\otimes D_{-})\cap \ker(\Delta_{-})=0.
$$
\end{proof}

\subsection{dg-Bialgebras and Hopf dg-algebras}\label{Hopf}

\begin{defi} 
Recall that a {\em differential graded bialgebra} (or a {\em dg-bialgebra} or simplys a {\em bialgebra}) is a monoid object in the category $\dgCoalg$. It is also a comonoid object in the category $\dgAlg$.
More concretely, it is a dg-vector space $B$ equipped with an algebra structure $(B,\mu,e)$ and a coalgebra structure $(B,\Delta,\epsilon)$ satisfying the compatiblity conditions expressed by the following four commutative diagrams
$$\xymatrix{
B\otimes B \ar[d]_{\Delta\otimes \Delta}  \ar[r]^\mu & B \ar[r]^\Delta &  B\otimes B  \\
B\otimes B \otimes B\otimes B   \ar[rr]^-{B\otimes \sigma \otimes B}& & B\otimes B\otimes B\otimes B\ar[u]^{\mu\otimes \mu}
}$$
$$\xymatrix{
B\otimes B \ar[r]^-{\mu}  \ar[dr]_-{\epsilon\otimes \epsilon} & B \ar[d]^-{\epsilon} \\
	& \FF
}\qquad
\xymatrix{
\FF \ar[r]^-{e}  \ar@{=}[dr] & B \ar[d]^-{\epsilon} \\
	& B
}\qquad
\xymatrix{
\FF \ar[dr]^-{e\otimes e}  \ar[d]_e & \\
B \ar[r]^-{\Delta} & B\otimes B 
}$$
\end{defi}
The unit element $e$ of a bialgebra $B=(B,\mu,e,\Delta,\epsilon)$
is uniquely determined by the product $\mu:B\otimes B\to B$
and the counit is uniquely determined by the coproduct $\Delta:B\to B\otimes B$.
Hence we may describe the bialgebra as a triple $B=(B,\mu,\Delta)$.

\medskip

For example, the field $\FF$ has the structure of a bialgebra.
A {\it map of bialgebras} $B\to B'$ is a linear map which is both an algebra map and a coalgebra map.
We shall denote the category of bialgebras by $\dgBialg$.
The bialgebra $\FF$ is both a terminal object and an initial object of the category $\dgBialg$.
The tensor product of two bialgebras has the structure of a bialgebra.
This gives the category $\dgBialg$ a symmetric monoidal structure in which the unit object is the bialgebra $\FF$.

\medskip
 
If $B=(B,m,e,\Delta,\epsilon)$ is a bialgebra, then the algebra $(B,m, e)$ is augmented by the map $\epsilon:B\to \FF$, and the coalgebra $(B,\Delta, \epsilon)$ is coaugmented by the map $e:\FF \to B$. 
Let us put $U_\Delta(B)=(B,m, e,\epsilon)$ and $U_m(B)=(B,\Delta, \epsilon,e)$ the functors forgetting respectively the coproduct and the product. 
This defines two forgetful functors
$$\xymatrix{
\dgAlg_\bullet  &\ar[l]_-{U_\Delta}\dgBialg \ar[r]^-{U_m} & \dgCoalg_\bullet .
}$$

We shall say that a bialgebra is {\em commutative} (or {\em cocommutative}) if its underlying algebra (or coalgebra) is.
We shall say that a bialgebra is {\em conilpotent} if its underlying pointed coalgebra is.

Let $\dgBialg^{\sf conil}$ be the full subcategory of $\dgBialg$ generated by conilpotent bialgebras.

\begin{prop}\label{laxcorad}
The radical functor $R^c$ provides a right adjoint to the inclusion $\dgBialg^{\sf conil}\subset \dgBialg$.
\end{prop}
\begin{proof}
The inclusion $\iota:\dgCoalgnil\subset \dgCoalg$ is monoidal, hence $R^c$ is a right monoidal functor. In particular the adjunction $\iota \dashv R^c$ passes to the category of monoids to give an adjunction $\iota:\dgBialg^{\sf conil}\rightleftarrows \dgBialg:R^c$ where $\iota$ is the canonical inclusion.
\end{proof}

\medskip

Recall that an algebra $A=(A,\mu)$ has an opposite $(A^o,\mu^{o})$
where $\mu^{o}(x^o,y^o)=\mu(y,x)^o (-1)^{|x||y|}$. 
 
\begin{defi}\label{defoppositebialgebra}
The {\it opposite} 
of a bialgebra $Q=(Q,\mu,\Delta)$ is defined to be the bialgebra
$Q^o=(Q^o,\mu^o,\Delta)$, where 
$$
\mu^{o}(a^o,b^o)=\mu(a,b)^o (-1)^{|a||b|} \et \Delta(x^o)=x^{(1)o}\otimes x^{(2)o}.
$$
The map $x\mapsto x^o$ is an anti-isomorphism of algebras $Q\to Q^o$ and an isomorphism of coalgebras. 
\end{defi}

We shall not need the notion of co-opposite bialgebra where the coproduct is reversed.

\subsubsection{Examples}\label{examplesbialg}

\begin{ex}[Coshuffle bialgebra]\label{coshuffle bialgebra}
If $X$ is a  dg-vector space, the {\it coshuffle coproduct} on the tensor algebra $T(X)$ defined in example \ref{coshufflecoalgebra} 
can be characterized as the unique map of algebras $\Delta:T(X)\to T(X)\otimes T(X)$
such that the following square commute
$$\xymatrix{
X\ar[rr]^-{\delta}\ar[d]_{i_X} && X \oplus X\ar[d]^{i_X\otimes 1 \oplus 1\otimes i_X}\\
T(X)\ar[rr]^-{\Delta} && T(X)\otimes T(X)
}$$
where $\delta:X\to X\oplus X$ is the diagonal 
and $i_X:x\to T(X)$ is the canonical inclusion.
The {\em shuffle counit} $\epsilon$ is the unique algebra map such that the following square commute
$$\xymatrix{
X\ar[rr]\ar[d]_{i_X} && 0\ar[d]\\
T(X)\ar[rr]^-{\epsilon} && T(0)\simeq \FF.
}$$
Equivalently, $\Delta$ and $\epsilon$ are the unique algebras map such that $\Delta(v)=1\otimes v+v\otimes 1$ and $\epsilon(v)=0$ for every $v\in X$.

This coalgebra structure gives the algebra $T(X)$ the structure of a bialgebra that we note $T^{csh}(X)$ and call the {\em coshuffle bialgebra}. It is a conilpotent bialgebra by example \ref{coshuffleconilpotent}.
\end{ex}

\begin{ex}[Shuffle bialgebra]\label{shuffle bialgebra2}
If $X$ is a graded vector space, the {\it shuffle product} on the tensor coalgebra $T^c(X)$ is defined as the unique map of algebras $\mu:T^c(X)\otimes T^c(X)\to T^c(X)$ such that the following square commute
$$\xymatrix{
T^c(X)\otimes T^c(X) \ar[rr]^-{\mu} \ar[d]_{p_X\otimes \epsilon \oplus \epsilon\otimes p_X}&& T^c(X)\ar[d]^{p_X}\\
X\oplus X \ar[rr]^-{+} && X
}$$
where $+:X\oplus X\to X$ is the sum,
$\epsilon$ is the canonical projection $T(X)\to \FF$
 and $p_X:T^c(X)\to X$ is the canonical projection.
The {\em shuffle unit} $e$ is the unique coalgebra map such that the following square commute
$$\xymatrix{
\FF\ar[rr]^-{e} \ar[d]&& T^c(X)\ar[d]^{p_X}\\
0\ar[rr] && X.
}$$
Equivalently, $\Delta$ and $\epsilon$ are the unique coalgebra maps such that 
$\mu(x\otimes y)=x\epsilon(y)+\epsilon(x)y$ for every $x,y\in T(X)$, and $e(1)=1$.

This algebra structure gives the coalgebra $T^c(X)$ the structure of a bialgebra that we note $T^{sh}(X)$ and call the {\em shuffle bialgebra}.
A variation using the cofree coalgebra functor rather is constructed in corollary \ref{shuffleproductTv}.
\end{ex}

\begin{rem}
Recall that $T(X)= \oplus_{n\geq 0}X^{\otimes n}$.
Using the isomorphism $X\otimes \FF\oplus \FF\otimes Y=X\oplus Y$, we define
the map $i':X\oplus Y\to T(X)\otimes T(Y)$ to be the canonical inclusion and 
the map $p':T(X)\otimes T(Y)\to X\oplus Y$ to be the canonical projection.
The map $i'$ extends to a canonical map of algebras $i:T(X\oplus Y)\to T(X)\otimes T(Y)$ and
the map $p'$ extends to a canonical map of conilpotent coalgebras $p:T^c(X)\otimes T^c(Y)\to T^c(X\oplus Y)$.
Moreover $i$ is a colax structure on the functor $T:(\dgVect, \oplus)\to (\dgAlg,\otimes)$
and $p$ is a lax structure on the functor $T^c:(\dgVect, \oplus)\to (\dgCoalgnil,\otimes)$.

The maps $X\oplus X\to X$ and $X\to 0$ used in example \ref{coshuffle bialgebra} are the structure maps of $X$ as a algebra in $(\dgVect,\oplus)$. Dually, the maps $\delta: X\to X\oplus X$ and $0\to X$ used in example \ref{shuffle bialgebra2} are the structure maps of $X$ as a coalgebra in $(\dgVect,\oplus)$. Together this equips $X$ with a bialgebra structure in $(\dgVect,\oplus)$.
The coshuffle coproduct and shuffle product are the image of these structures through the (co)lax structure of $T$ and $T^c$. 

An analog of the shuffle product can be defined on the cofree coalgebra functor $T^\vee(X)$ (see corollary \ref{shuffleproductTv}).
\end{rem}

\bigskip

We recall that a bialgebra $H$ is said to be a {\it Hopf algebra} if the identity map $id:H\to H$ has an inverse in the convolution algebra $[H,H]$. Such an inverse is a map $S:H\to H$ such that, for $x\in H$
$$
S(x^{(1)})x^{(2)} = x^{(1)}S(x^{(2)}) = \epsilon(x)1_H.
$$
$S$ is called the {\it antipode} and it is unique when it exists.
The antipode is an anti-homomorphism of algebras $S:H\to H$
and an anti-homomorphism of coalgebras.
If $H$ is cocommutative, the antipode defines a bialgebra map $S:H\to H^{op}$.

\begin{ex}[Coshuffle Hopf algebra] \label{shufflehopf} 
The coshuffle bialgebra $T^{csh}(X)$ of example \ref{coshuffle bialgebra} is a Hopf algebra.
The antipode $S:T(Z)\to T(Z)$ is the unique anti-homomorphism of algebras such that $S(x)=-x$ for every $x\in X$.

In particular, if $\delta$ is a variable of degree $n$, then the bialgebra $T^{csh}(\FF\delta)$
has the structure of a cocommutative Hopf algebra with the coproduct
$\Delta:\FF[\delta] \to\FF[\delta] \otimes \FF[\delta]$ defined by putting
$\Delta(\delta)=\delta \otimes 1+1\otimes \delta.$
The antipode $S:\FF[\delta]\to \FF[\delta]$ is the unique anti-homomorphism of algebras such that $S(\delta)=-\delta$.
If $|\delta|$ is even we have $S(\delta^n)=\delta^n(-1)^n$ for every $n$, and
if $|\delta|$ is odd we have $S(\delta^n)=\delta^n(-1)^{{n\choose 2}+n}$ for every $n$.
\end{ex}

\begin{ex}[Primitive Hopf algebra] \label{primitiveHopf} 
For $n$ odd (!), let us consider the coalgebra $\d_n=T^{csh}_1(\FF\delta) = \FF\delta_+=\FF\oplus \FF\delta $ where $\delta$ is a variable of degree $n$.
It is also an algebra where for the product defined by $\delta^2=0$.
Because $|\delta|$ is odd, we have
$$
\Delta(\delta)^2 = (\delta\otimes 1+1\otimes \delta)^2 = \delta\otimes \delta +\delta \otimes \delta \ (-1)^{|\delta|^2} = 0.
$$
This proves that $T^{csh}_1(\FF\delta) = \FF\delta_+$ is a bialgebra. It is also a Hopf algebra with antipode the unique anti-homomorphism of algebras such that $S(\delta)=-\delta$.
\end{ex}

\subsubsection{Non-biunital bialgebras}\label{nonbiunitalbialgebras}

If $B=(B,m,e,\Delta,\epsilon)$ is a bialgebra, then the vector space $B_{-}=\ker(\epsilon)$
has the structure of a non-unital algebra $(B_{-},m_{-})$
and of a non-counital coalgebra $(B_{-},\Delta_{-})$. 
Moreover, for every $x,y\in B_{-}$ we have
$$
\Delta_{-}(xy)=x\otimes y+y\otimes x +\Delta_{-}(x)(1\otimes y+ y\otimes 1)+ (1\otimes x+ x\otimes 1)\Delta_{-}(y) + \Delta_{-}(x)\Delta_{-}(y).
$$

\begin{defi} 
Let $P$ be a vector space equipped with a non-unital algebra structure $\mu:P\otimes P\to P$, $xy=\mu(x,y)$ together with a non-counital 
coalgebra structure $\Delta:P\to P\otimes P$. 
We shall say that the triple $(P,\mu, \delta)$ is a {\em non-biunital bialgebra} if the following compatibility condition is satisfied,
$$
\Delta(xy)=x\otimes y+y\otimes x +\Delta(x)(1\otimes y+ y\otimes 1)+ (1\otimes x+ x\otimes 1)\Delta(y)
+ \Delta(x)\Delta(y).
$$
\end{defi}

Beware that in this  formula
$$
\Delta(x)(1\otimes y+ y\otimes 1):= x_{(1)}\otimes x_{(2)}y+x_{(1)}y\otimes x_{(2)}(-1)^{|y||x_{(2)}|}
$$
$$
(x\otimes 1+1\otimes x)\Delta(y):=xy_{(1)}\otimes y_{(2)}+y_{(1)}\otimes xy_{(2)} (-1)^{|x||y_{(1)}|}
$$

\medskip

A {\it map of non-biunital bialgebras} is both a map of non-unital algebras and a map of non-counital coalgebras.
We shall denote the category of non-biunital bialgebras by $\dgBialg_\circ$.
If $B=(B,m,\Delta)$ is a non-biunital bialgebra, 
let us put $U_\Delta(B)= (B,m)$ and $U_m(B)= (B,\Delta)$
the functors forgetting respectively the coproduct and the product. 
This defines two forgetful functors
$$\xymatrix{
\dgAlg_\circ  &\ar[l]_-{U_\Delta}\dgBialg_\circ \ar[r]^-{U_m} & \dgCoalg_\circ .
}$$

\bigskip

If $B=(B,m,e,\Delta,\epsilon)$ is a bialgebra, then the vector space $B_{-}=\ker(\epsilon)$
has the structure of a non-biunital bialgebra.
This defines a functor $(-)_{-}:\dgBialg \to \dgBialg_\circ.$
Conversely, $P=(P,m,\Delta)$ is a non-biunital bialgebra, then the vector space $P_{+}=\FF\oplus P$
has the structure of a  bialgebra with unit $(1,0)$, with counit the projection $\epsilon:P_{+}\to \FF$
and with coproduct $\Delta_{+}$. This defines a functor $(-)_{+}:\dgBialg_\circ \to \dgBialg$.

\begin{prop}\label{equivBilagnonunital} The functors 
$$\xymatrix{
(-)_{-}:\dgBialg \ar@<.6ex>[r] &\dgBialg_\circ:(-)_{+} \ar@<.6ex>[l]
}$$
are inverse equivalences of categories.
\end{prop}

We have in fact a commutative square
$$\xymatrix{
\dgAlg_\bullet \ar@<-.6ex>[dd]_{(-)_-}&& \dgBialg \ar@<-.6ex>[dd]_{(-)_-} \ar[ll]_-{U_\Delta}\ar[rr]^-{U_m} && \dgCoalg_\bullet \ar@<-.6ex>[dd]_{(-)_-}\\
&&\\
\dgAlg_\circ\ar@<-.6ex>[uu]_{(-)_+} && \dgBialg_\circ \ar@<-.6ex>[uu]_{(-)_+}\ar[ll]_-{U_\Delta} \ar[rr]^-{U_m}&& \dgCoalg_\circ\ar@<-.6ex>[uu]_{(-)_+}
}$$
where the horizontal functors are equivalences.
We see that $\dgBialg_\bullet$ can be described also the category of non-unital monoids in $(\dgCoalg_\bullet,\wedge)$
or of non-counital comonoids in $(\dgAlg_\bullet,\wedge)$.
In particular, we deduce that $U_m:\dgBialg\to \dgCoalg_\bullet$ have a left adjoint. 
For a non-counital coalgebra $(C,\delta)$ it is given by $T(C_-)$. The product and the unit are the tensor product and the tensor unit,
the coproduct and counit are the unique algebra map $\Delta:T(C_-)\to T(C_-)\otimes T(C_-)$ and $\epsilon:T(C_-)\to \FF$
such that
$$
\Delta(x) = x\otimes 1 +1\otimes x + \delta(x)
\et 
\epsilon(x)=0
$$
for any $x\in C_-$.
We shall call {\em quasi-coshuffle bialgebra} this bialgebra and note it by $T^{qcsh}_\bullet(C)$ when $C$ is a pointed coalgebra (and $T^{qcsh}(C)$ if $C$ is a non-counital coalgebra).
The name is justified by the fact that, if $C$ is a primitive coalgebra $T_{\bullet,1}(X)$, the bialgebra $T^{qcsh}_\bullet(C)$ is the coshuffle bialgebra $T^{csh}(X)$.

\begin{prop}\label{freebialgfromcoalg}
The functor $U_m:\dgBialg\to \dgCoalg_\bullet$ has a left adjoint given by the quasi-coshuffle bialgebra functor $T^{qcsh}_\bullet$.
\end{prop}

\medskip

We shall prove in proposition \ref{adjointforgetbialgalg} that $U_\Delta:\dgBialg\to \dgAlg_\bullet$ has a right adjoint.
For now, we recall the construction of the adjoint of the restriction of $U_\Delta$ to the categorie $\dgBialg^{\sf conil}$ of conilpotent bialgebras.

Let  $A=(A,m)$ is a non-unital algebra, let $T^c_\bullet(A)=(T^c(A),\Delta, \epsilon, e)$ be the cofree conilpotent pointed coalgebra on the vector space $A$ and let  $p:T^c(A)\to A$ be the cogenerating map. 

\begin{lemma}\label{cofreebialgebra1} 
There is a unique map of coalgebras $\mu:T^c(A)\otimes T^c(A) \to T^c(A)$ such that
$$
p\mu(x,y)=p(x)p(y)+\epsilon(x)p(y)+p(x)\epsilon(y)
$$
for every $x, y \in T^c(A)$.
The product $\mu$ is associative and the pair $(\mu,e)$ is defining a bialgebra structure on the coalgebra $T^c(A)=(T^c(A),\Delta, \epsilon)$.
\end{lemma} 

If $A$ is non-unital algebra, then the product $\mu:T^c(A)\otimes T^c(A) \to T^c(A)$ defined above is called the {\it quasi-shuffle product} and the algebra $T^{qsh}_\bullet(A)=(T^c(A),\mu)$ the {\it quasi-shuffle algebra} of $A$.
The deconcatenation coproduct $\Delta$ gives the quasi-shuffle algebra $T^c(A)$ the structure of a bialgebra  $(T^c(A),\mu, \Delta)$ called
the {\it quasi-coshuffle bialgebra} of $A$.
In particular, if the product of $A$ is the zero map, the quasi-shuffle product on $T^c(A)$ coincides with the shuffle product defined in example \ref{shuffle bialgebra2}.

\medskip
If $A=(A,m,e_A, \epsilon_A)$ is a pointed algebra and $p:T^c(A_{-})\to A_{-}$ is the cogenerating map, let us denote by $q:T^c(A_{-})\to A$ the map defined by putting $q(x)=\epsilon_A(x)e_A+p(x)$.

\begin{prop}\label{cofreebialgebra4} 
The forgetful functor $U_\Delta:\dgBialg^{\sf conil} \to \gAlg_\bullet$ admits a right adjoint which associates to a pointed algebra $A$
the quasi-coshuffle bialgebra $(T^c(A_{-}),\mu, \Delta)$.
Moreover, the canonical map  $q:T^c(A_{-}) \to A$ is the counit of the adjunction.
\end{prop} 

The proof is the same as that proposition \ref{adjointforgetbialgalg}.

\subsubsection{Pointed and non-counital bialgebras}\label{pointed bialgebras}

 A {\it non-counital bialgebra} is defined to be a monoid object in the category of non-counital coalgebras $\dgCoalg_\circ$.
Let us denote the category of non-unital bialgebras by $\dgBialg_{\!\not\,\epsilon}$.

\medskip

If $E=(E,\mu,1,\Delta)$ is a non-unital bialgebra, let us show that $\FF\times E$
The algebra structure of $\FF\times E$ is the product
$(\alpha,x)(\beta, y)=(\alpha\beta, xy)$.
The coalgebra structure of $\FF\times E$ is identical to the  coalgebra structure of $E_{+}=\FF\oplus E$.
Let us put $\theta =(1,0)$. Then $\Delta_{+}(\theta)=\theta\otimes \theta$
and we have $\Delta_{+}(x)=\theta \otimes x+x\otimes \theta+\Delta(x)$ for every $x\in E$.

Observe that $\theta \theta=\theta$ and that $\theta x=x\theta=0$ for every $x\in E$.
Let us verify that $\Delta_{+}:\FF\times E \to (\FF\times E)\otimes (\FF\times E)$
is an algebra map.  For every $(\alpha,x)$ and $(\beta,y) \in \FF\times E$
we have
\begin{eqnarray*}
\Delta_{+}(\alpha e+x)\Delta_{+}(\beta e+ y) &=& (\alpha e\otimes e+e\otimes x+x\otimes e+x_{1}\otimes x_{(2)})
(\beta e\otimes e+e\otimes y+y\otimes e+y_{1}\otimes y_{(2)})\\
&=& \alpha \beta e\otimes e + e\otimes xy +xy\otimes e+x_{1}y_{1} \otimes x_{(2)}y_{(2)} (-1)^{|x_{(2)}|||y_{(1)}||}\\
&=&\Delta_{+}(\alpha\beta+xy).
\end{eqnarray*}
The unit of the algebra $\FF\otimes  E$ is the element $(1,1)=\theta+1$.
We have
$$\Delta_{+}(\theta+1)=\theta \otimes \theta +\theta\otimes 1+1\otimes \theta+1\otimes 1=(\theta+1)\otimes (\theta+1).$$
We have proved that $\Delta_{+}$ is a map of algebras.
It is obvious that the projection $\epsilon:\FF\times E\to \FF$
is a map of algebras.

\begin{prop}\label{non-unital bialgebra} The forgetful functor $U_\epsilon:\dgBialg\to \dgBialg_{\!\not\,\epsilon}$
has a right adjoint which associates to a non-unital bialgebra $E$ the bialgebra
$\FF\times E$.
\end{prop}

\begin{defi}
If $B$ is a bialgebra, we shall say that an element  $\theta\in B$
is {\it absorbant}
if it satisfies the following conditions
\begin{itemize}
\item $\Delta(\theta)=\theta \otimes \theta \et \epsilon(\theta)=1$ 
\item $x\theta =\epsilon(x)\theta=\theta x$
for every $x\in B$. 
\end{itemize}
\end{defi}

The first condition means that $\theta$ is an atom of the coalgebra $(B,\Delta,\epsilon)$.
The second condition means that the following two squares commute,
$$\xymatrix{
B\ar[d]_\epsilon \ar[r]^(0.4){B\otimes \theta} & B\otimes B\ar[d]_{\mu}& \ar[l]_(0.35){\theta \otimes B} B\ar[d]^\epsilon \\
\FF \ar[r]^-{\theta} & B& \ar[l]_-{\theta} B.
}
$$
For example, if $E$ is a non-unital bialgebra, then the element $\theta=(1,0)\in \FF\times E$
is absorbant. 

\medskip
Observe that of $\theta_1$ and $\theta_2$
are two absorbant elements of a bialgebra $B$, then

$$\theta_1=\epsilon(\theta_2)\theta_1=\theta_2 \theta_1=\epsilon(\theta_1)\theta_2=\theta_2.$$

Thus, an absorbant element $\theta\in B$ is unique when it exists.
We shall say that a bialgebra $B$ equipped with an absprbant element $\theta\in B$
 is {\it pointed}.  A {\it map of pointed bialgebras} $(B,\theta)\to (B',\theta')$
 is a map of biagebras $f:B\to B'$ such that $f(\theta)=\theta'$.
We shall denote the category of pointed bialgebras by $\dgBialg_\bullet$.

\medskip

It follows from the definition that an aborbant element $\theta\in B$ is a central idempotent of the algebra $B$.
We thus have a product decomposition $B=B\theta \times B(1-\theta)$.
Moreover,   the map $\epsilon:B\to \FF$
induces an isomorphism $B\theta\simeq \FF$ and we have
$B(1-\theta)=\ker(\epsilon)$.
We thus obtain a a product decomposition $B=\FF\theta \times B_{-}$,
where $B_{-}=\ker(\epsilon)$ is an algebra with unit $1-\theta$.
The vector space $B_{-}$ has then the structure of a non-counital coalgebra with the coproduct
 $\Delta_{-}$ defined 
 by putting
 $$\Delta_{-}(x)=\Delta(x)-\theta \otimes x-x\otimes \theta$$
 for $x\in B_{-}$. 
 Let us show that $\Delta_{-}:B_{-}\to B_{-} \otimes B_{-}$
 is a map of algebras.
 For every $x, y\in B_{-}$
 we have 
\begin{eqnarray*}
\Delta(x)(\theta \otimes y) &=&(x_{(1)}  \otimes x_{(2)})(\theta  \otimes y) =x_{(1)} \theta  \otimes x_{(2)}y \\
&=& \epsilon(x_{(1)}) \theta  \otimes  x_{(2)}y = \theta  \otimes  \epsilon(x_{(1)} x_{(2)}y \\
&=& \theta  \otimes xy 
\end{eqnarray*}
Similarly, we have $\Delta(x)(y \otimes \theta)=xy\otimes \theta $, 
$(\theta \otimes x)\Delta(y)=\theta \otimes xy $ and $(x \otimes \theta)\Delta(y)=xy\otimes \theta.$
It follows that
\begin{eqnarray*}
\Delta_{-}(x)\Delta_{-}(y) &=& (\Delta(x)-\theta \otimes x-x\otimes \theta)
(\Delta(y)-\theta \otimes y-y\otimes \theta) \\
&=&\Delta(x)\Delta(y)-\theta \otimes xy-xy\otimes \theta\\
&& -\theta \otimes xy+\theta \otimes xy -xy \otimes \theta +xy\otimes \theta\\
&=& \Delta(xy)-\theta \otimes xy-xy\otimes \theta \\
&=& \Delta_{-}(xy).
\end{eqnarray*}
Moreover,
\begin{eqnarray*}
\Delta_{-}(1-\theta) &=& 1\otimes 1-\theta \otimes \theta -\theta \otimes (1-\theta)- (1-\theta) \otimes \theta \\
&=& 1\otimes 1-\theta \otimes 1- 1 \otimes \theta + \theta \otimes \theta  \\
&=& (1-\theta)\otimes (1-\theta). 
\end{eqnarray*}
This shows that $\Delta_{+}:B_{-}\to B_{-} \otimes B_{-}$
 is a map of algebras.
 Thus, $B_{-}$ has the structure of a non-unital bialgebra.
We have defined a functor $(-)_{-}: \dgBialg_\bullet \to \dgBialg_{\!\not\,\epsilon}$.

\begin{prop}\label{equivBilagpointed}
The functors 
$$\xymatrix{
(-)_{-}: \dgBialg_\bullet \ar@<.6ex>[r]& \dgBialg_{\!\not\,\epsilon}: \FF\times(-) \ar@<.6ex>[l]
}$$
are inverse equivalences of categories.
\end{prop}

We have in fact a commutative square of adjunctions
$$\xymatrix{
\dgBialg_\bullet \ar@<.6ex>[rr]^-{(-)_-} \ar@<.6ex>[dd]^{U}&& \dgBialg_{\!\not\,\epsilon} \ar@<.6ex>[ll]^-{(-)_+}\ar@<.6ex>[dd]^{U}\\
&&\\
\dgCoalg_\bullet \ar@<.6ex>[rr]^-{(-)_-} \ar@<.6ex>[uu]^{T}&& \dgCoalg_\circ \ar@<.6ex>[ll]^-{(-)_+} \ar@<.6ex>[uu]^{T}
}$$
where the horizontal functors are equivalences.
We deduce that $\dgBialg_\bullet$ can be described also the category of monoids in the category $(\dgCoalg_\bullet,\wedge)$.

\subsubsection{Modules over a cocommutative Hopf algebra}\label{moduleoverhopfalgebra}

A {\it module} over a bialgebra $Q=(Q,\mu,e,\Delta,\epsilon)$ is defined to be 
a module over the algebra $(Q,\mu,e)$.
We shall denote the category of (left) $Q$-modules by $\Mod(Q)$.
It follows from proposition \ref{enrichmentofmodules} that
the category $\Mod(Q)$ is enriched and bicomplete over the 
category $\dgVect$. 
If $X$ and $Y$ are $Q$-modules we shall note $\Hom_{Q}(X,Y)$ this enrichment.
$\Hom_{Q}(X,Y)$ is naturally a vector sub-space of $[X,Y]$: a morphism $f:X\rhup Y$ of degree $n$ belongs to $\Hom_{Q}(X,Y)_n$ iff we have 
$$
f(q\cdot x)=q\cdot f(x)\ (-1)^{n|q|}
$$
for every $q\in Q$ and $x\in X$.

\begin{prop} \label{freeandcofreeQmodules}
The forgetful functor 
${\Mod}(Q) \to\dgVect$ admits a left adjoint $X\to Q\otimes X$ 
and a right adjoint $X\mapsto [Q,X]$.
Limits and colimits exist in ${\Mod}(Q)$ and they can be computed in $\dgVect$.
\end{prop}

\begin{proof}
This follows from proposition \ref{moduleandmonad}.
\end{proof}

The tensor product of two $Q$-modules $X$ and $Y$ has the structure of a $Q$-module $X\otimes Y$.
The left action of $Q$ on $X\otimes Y$ is the composite of the maps,
$$\xymatrix{
Q \otimes X\otimes Y\ar[rr]^-{\Delta\otimes X\otimes Y } && Q \otimes Q \otimes X\otimes Y \ar[r]^-{\sigma_{23}} &  Q \otimes X\otimes  Q \otimes Y\ar[rr]^-{a_X\otimes a_Y} && X\otimes Y,
}$$
 In terms of elements,
$
q(x\otimes y) = q^{(1)}x \otimes q^{(2)}y \ (-1)^{|q^{(2)}||x|}
$
for $q\in Q$, $x\in X$ and $y\in Y$.
The unit object for the tensor product is the field $\FF$ equipped the trivial action $\epsilon\otimes \FF:Q\otimes \FF \to \FF$.
Beware that the symmetry $\sigma:X\otimes Y\to Y\otimes X$ is generally not a map
of $Q$-modules, unless $Q$ is cocommutative.

The action of $Q$ on $X\otimes Y$ can also be described by the representation
$$\xymatrix{
Q \ar[r]^-{\Delta} & Q \otimes Q \ar[rr]^-{\pi_X\otimes \pi_Y} && [X,X]\otimes [Y,Y] \ar[r]^-{\theta} &  [X\otimes Y,X\otimes Y]
}$$
where the $\pi$s are the actions of $Q$ on $X$ and $Y$ given by proposition \ref{moduleandmonad} and $\theta$ is 
the strength of the strong functor $\otimes$ in $\dgVect$ (in particular it is an algebra map).
All the maps of this diagram are algebra maps, therefore their composition is indeed a $Q$-module structure on $X\otimes Y$.

\bigskip

Moreover, if $Q$ is a cocommutative Hopf algebra, the symmetric monoidal category $(\Mod(Q),\otimes,\FF)$ is closed.
The internal hom between two $Q$-modules $X$ and $Y$ is given by $[X,Y]$ equipped with the action
$$\xymatrix{
Q\otimes [X,Y] \ar[rr]^-{(S\otimes Q)\circ \Delta}&& Q^o \otimes Q\otimes [X,Y]\ar[r]^-{a}& [X,X]^o \otimes [Y,Y]\otimes [X,Y]\ar[r]^-{c}& [X,Y].
}$$
where $a$ is given the two action maps $Q\to [X,X]$ and $Q\to [Y,Y]$ and $c$ is the strong composition in $\dgVect$.
In terms of elements, the action is given by
$$
(q\cdot f)(x) = q^{(1)}f(S(q^{(2)}) x) \ (-1)^{|q^{(2)}||f|}.
$$
In particular, we have $\Hom_Q(\FF,[X,Y])=\Hom_Q(X,Y)$.
The proof that this is the internal hom will be given in proposition \ref{freeandcofreeQmodules2}.

The action of $Q$ on $[X,Y]$ can also be described as the representation
$$\xymatrix{
Q \ar[r]^-{\Delta} & Q \otimes Q\ar[r]^-{S\otimes Q} & Q^o  \otimes Q \ar[rr]^-{\pi_X\otimes \pi_Y} && [X,X]^o \otimes [Y,Y] \ar[r]^-{\theta} &  [[X, Y],[X, Y]]
}$$
where $S$ is the antipode and $\theta$ is the strength of the internal hom $[-,-]$ of $\dgVect$ (in particular it is an algebra map).
Again all maps are algebra maps and the composition is a $Q$-module structure on $[X,Y]$.

\begin{ex}
If $\delta$ is a variable of degree $n$, then a module over the bialgebra $\FF[\delta]$ of example \ref{shufflehopf} is a dg-vector space $X$ equipped with a morphism $\delta_X:X\rhup X$ of degree $n$. 
The tensor product of two $\FF[\delta]$-modules $(X,\delta_X)$ and $(Y,\delta_Y)$ is the dg-vector space $X\otimes Y$ equipped with the morphism $\delta=\delta_X\otimes Y+X\otimes \delta_Y$.
By definition, we have
$$
\delta(x\otimes y)=\delta_X(x)\otimes y+x\otimes \delta_Y(y)(-1)^{n|x|}
$$
for every $x\in X$ and $y\in Y$.
The monoidal category ${\Mod}(\FF[\delta])$ is symmetric and closed, since the bialgebra $\FF[\delta]$ is a cocommutative Hopf algebra.
The internal hom of $(X,\delta_X)$ and $(Y,\delta_Y)$  if the dg-vector space $([X,Y],\delta)$ where $\delta$ is defined by
$$
\delta(f)(x) = \delta_Y(f(x)) - f(\delta_X(x)) \ (-1)^{|f||\delta|}
$$
for any graded morphism $f$ and any $x\in X$.
\end{ex}

\begin{ex}\label{fromgtodg}
If $\delta$ is of odd degree, we can impose $\delta^2$ in the previous example.
Let $Q=\d=\FF\delta_+$ be the bialgebra of example \ref{primitiveHopf}.
A $Q$-module dg-vector space $X$ equipped with a morphism $\delta_X:X\rhup X$ of degree $n$ and of square zero.
If $n=-1$, $X$ is a graded vector space equipped with two commuting differentials.

Remark that $Q$ is in fact a graded bialgebra rather than a dg-bialgebra, so it make sense to look at $Q$-module in the category $\gVect$ of graded vector spaces. In this context and still with $n=-1$, a $Q$-module is simply a dg-vector space. Moreover the formulas of the previous example show that $\Mod(Q)=\dgVect$ as a symmetric monoidal closed category.
\end{ex}

\begin{prop} \label{freeandcofreeQmodules2}
Let $Q$ be a cocommutative Hopf algebra, then
the category $(\Mod(Q),\otimes, \FF, [-,-])$ is symmetric monoidal closed and 
the forgetful functor ${\Mod}(Q) \to\dgVect$ is symmetric monoidal and preserves the internal hom.
\end{prop}
\begin{proof}
We need only to prove the closeness of the tensor product. 
This is equivalent to prove that the isomorphism $\lambda=\lambda^2:[X\otimes Y,Z]\simeq [X,[Y,Z]]$ induces a bijection between the dg-vector subspaces of $Q$-equivariant maps.
Let $f\in \Hom_Q(X\otimes Y,Z)$, let us prove that $\lambda f\in \Hom_Q(X,[Y,Z])$.
We need to prove that for any $x\in X$, $(\lambda f)(q\cdot x) = (q\cdot \lambda f)(x)$, let $y\in Y$ we have indeed
\begin{eqnarray*}
(q\cdot \lambda f)(x) & = & q^{(1)}(\lambda f(x)(S(q^{(2)})y) \ (-1)^{|q^{(2)}|(|f|+|x|)}\\
	& = & q^{(1)}f(x\otimes S(q^{(2)})y) \ (-1)^{|q^{(2)}|(|f|+|x|)}\\
	& = & f(q^{(1)}x\otimes q^{(2)}S(q^{(3)})y) \ (-1)^{|q^{(3)}|(|f|+|x|)+|q^{(2)}|(|f|+|x|) + |q^{(1)}||f|)}\\
	& = & f(q^{(1)}x\otimes \epsilon(q^{(2)})y) \ (-1)^{|q^{(2)}|(|f|+|x|) + |q^{(1)}||f|)}\\
	& = & f(qx\otimes y) \ (-1)^{|q||f|}\\
	& = & (\lambda f)(q\cdot x)(y).
\end{eqnarray*}
\end{proof}

Let us present a more conceptual proof of this computation.
Remark first that, for two $Q$-modules $X$ and $Y$, $\Hom_Q(X,Y)$ is entirely determined by the $Q$-module structure of $[X,Y]$.
It is defined as the equalizer in $\dgVect$
$$\xymatrix{
\Hom_Q(X,Y) \ar[r]& [X,Y]=[\FF,[X,Y]]\ar@<.6ex>[rr]^-{[\epsilon_Q,[X,Y]]}\ar@<-.6ex>[rr]_-{\lambda^1a} && [Q,[X,Y]]
}$$
where $a:Q\otimes [X,Y]\to [X,Y]$ is the action of $Q$. The closeness will be proven if we show that $\lambda:[X\otimes Y,Z]\simeq [X,[Y,Z]]$ is an isomorphism of $Q$-modules.

The adjunction $\lambda:[X\otimes Y,X]\simeq [X,[Y,Z]]$ is strong for the enrichment of $\dgVect$ over itself. 
This implies the commutativity of the square
$$\xymatrix{
[X,X]^o \otimes [Y,Y]^o \otimes [Z,Z] \ar[r]\ar[d] &[X,X]^o \otimes [[Y,Z],[Y,Z]]\ar[d]\\
[X\otimes Y,X\otimes Y]^o  \otimes [Z,Z]\ar[r]& End([X\otimes Y,Z])= End([X,[Y,Z]]) .
}$$
By composing at the source by the map $(S\otimes S\otimes Q)\circ \Delta^{(3)}:Q\to Q^o \otimes Q^o \otimes Q$
we deduce that $\lambda:[X\otimes Y,Z]\simeq [X,[Y,Z]]$ is an isomorphism of $Q$-modules.

\medskip
Presented this way the proof of proposition \ref{freeandcofreeQmodules2} can be abstracted for any Hopf algebra in any symmetric monoidal closed category. We shall use this in theorem \ref{Qhomcoalg}.

\subsection{Lie dg-algebras and (co)derivations}

Recall that a {\it Lie dg-algebra} (or {\it Lie algebra} for short) is a dg-vector space $X$ equipped with a bilinear operation
$$
[-,-]:X\otimes X\to X
$$
called the {\it bracket} satisfying the following two conditions
\begin{itemize}
\item $[x,y]=-[y,x](-1)^{|x||y|}$ (anti-symmetry)
\item $[x,[y,z]]=[[x,y],z]+ [y,[x,z]] (-1)^{|x||y|}$ (Jacobi identity).
\end{itemize}

If $X$ and $Y$ are (graded) Lie algebras, we shall say that a linear map $f:X\to Y$ is a {\it homomorphism } (resp.  {\it an anti-homomorphism}) if we have 
$$
f[x,y]=[fx,fy] \quad \quad {\rm ( } {\rm resp.} \quad f[x,y]=[fy,fx](-1)^{|x||y|}  {\rm ) }
$$
for every $x,y\in X$. 
A linear map $f:X\to Y$ is an anti-homomorphism if and only if the opposite map $-f:X\to Y$ is a homomorphism.

We shall note $\dgLie$ the category of Lie dg-algebras.

\begin{rem}
All Lie algebras considered in this paper will be Lie sub-algebras of endomorphisms Lie algebras 
and therefore equipped with a quadratic {\it squaring operation} $x\mapsto x^2$ 
defined for elements of odd degree. We have $x^2={1\over 2} [x,x]$
when the characteristic of the field $\FF$ is $\neq 2$.
The squaring operation can be extended to
elements of even degree when the field $\FF$ has characteristic 2.

\end{rem}

\medskip

\begin{prop}\label{semialgtoLie}
If $A$ is a non-unital algebra, then commutator operation $ [x,y]=xy-yx(-1)^{|x||y|}$ is a Lie algebra structure on $A$.
\end{prop}

\begin{proof}
Obviously $ [x,y]=[y,x](-1)^{|x||y|}$.
Let us prove the Jacobi identity.
The map $[x,-]:A\to A$ is a derivation of the non-unital algebra $A$
by example \ref{innerderivationsofsemialg}.
Thus
\begin{eqnarray*} 
[x,[y,z]] &= & [x,yz]-[x,zy](-1)^{|y||z|} \\
&= & [x,y]z+y[x,z](-1)^{|x||y|} -[x,z]y (-1)^{|y||z|}-z[x,y] (-1)^{|y||z|+|x||z|} \\
&= & [[x,y],z]+[y,[x,z]](-1)^{|x||y|}
\end{eqnarray*}
This is equivalent to the Jacobi identity.
\end{proof}

The proposition above shows that a non-unital algebra $A$ has the structure of a Lie algebra  $(A,[-,-])$.
We shall denote this Lie algebra by $Lie(A)$.

\begin{ex}\label{theLiealgebraofendomorphism}
If $X$ is a graded vector space, then the algebra of endomorphisms
of $X$ has the structure of a Lie algebra $Lie(End(X))$ in which
the bracket operation is the commutator 
$[f,g]=fg-gf(-1)^{|f||g|}.$
We shall denote this Lie algebra by $gl(X)$.
\end{ex}

Let $X$ be a (dg-)vector space, recall that the free Lie algebra generated by $X$
can be characterized as the Lie sub-algebra of $Lie(T(X))$ generated by $X$ \cite{Reutenauer}.

\begin{prop}
The category $\dgLie$ is bicomplete and the forgetful functor $\dgLie\to \dgVect$ has a left adjoint which associates to a vector space $V$ the Lie algebra $L(V)$ freely generated by $V$.
\end{prop}

The forgetful functor $\Lie\to \dgVect$ preserves and reflects limits, since it is continuous and conservative. 
In particular, the cartesian product of a family of Lie algebras $(L_i:i\in I)$ has the structure of a Lie algebra $\prod L_i$.
We shall denote the cartesian product of two Lie algebras $L$ and $P$ by $L\times P$.

\bigskip

A homomorphism of algebras $f:A\to B$ induces a homorphism of Lie algebras $f:Lie(A)\to Lie(B)$.
An anti-homomorphism of algebras $f:A\to B$ induces an anti-homorphism of Lie algebras $f:Lie(A)\to Lie(B)$.

\begin{prop}\label{Lieanti}
If $f:A\to B$ is an anti-homomorphism of non-unital algebras, then $-f$ is a homomorphism of Lie algebras $Lie(A)\to Lie (B)$.
\end{prop}

\begin{proof} 
An anti-homomorphism of non-unital algebras $f:A\to B$ induces an anti-homorphism of Lie algebras $f:Lie(A)\to Lie(B)$.
It follows that $-f:Lie(A)\to Lie (B)$ is a homomorphism of Lie algebras.
\end{proof}

\medskip

\begin{defi}\label{defcommutation}
If $P$, $Q$ and $L$ are Lie algebras, we shall say that two homomorphisms $f:P\to L$ and $g:Q\to L$ {\it commutes} if we have $[f(x),g(y)]=0$ for every $x\in P$ and $y\in Q$.
\end{defi}

\begin{lemma} \label{commutehomoLie}
Let $f:P\to L$ and $g:Q\to L$ be a commuting pair of homomorphisms of Lie algebras.
Then the map $h:P\times Q\to L$ defined by putting $h(x,y)=f(x)+g(y)$ is a homomorphism of Lie algebras.
\end{lemma}

\begin{proof} 
If $(x,y)\in P_n\oplus Q_n $ and $(u,v)\in P_m\oplus Q_m $ then
\begin{eqnarray*} 
[h(x,y),h(u,v)]  &= & [f(x)+g(y),f(u)+g(v)] \\
&= &[f(x),f(u)] +[f(x),g(v)] + [g(y),f(u)] + [g(y),g(v)]\\
&= &[f(x),f(u)]  + [g(y),g(v)]\\
&= & f[x,u]+g[y,v] \\
&= &h([x,u],[y,v])
\end{eqnarray*}
\end{proof}

\begin{prop}\label{laxLie}
If $A$ and $B$ are algebras then the map $\rho:A\times  B\to A\otimes B$ defined by putting 
$\rho(a,b)=a\otimes 1+1\otimes b$ is a homomorphism of Lie algebras.
The homomorphism preserves the square of odd elements.
\end{prop}

\begin{proof}
The homomorphisms $i_1:Lie(A)\to Lie(A\otimes B)$ and $i_2:Lie(B)\to Lie(A\otimes B)$ 
induced by the canonical maps $i_1:A\to A\otimes B$ and $i_2:B\to A\otimes B$ commute
since 
$$
i_1(a)i_2(b)=(a\otimes 1)(1\otimes b)=a\otimes b=(1\otimes b)(a\otimes 1)(-1)^{|a||b|}=i_2(b)i_1(a)(-1)^{|a||b|}.
$$
Hence the map $\rho:A\times  B\to A\otimes B$ defined by putting $\rho(a,b)=i_1(a)+i_2(b)$
is a homomorphism of Lie algebras by lemma \ref{commutehomoLie}.
Moreover, if $n$ is odd and $(a,b)\in A_n\times B_n$ then
$$
\rho(a,b)^2 =(a\otimes 1+1\otimes b)^2 = aa\otimes 1+a\otimes b (-1)^{n n} +a\otimes b+1\otimes bb=a^2\otimes 1+1\otimes b^2=\rho(a^2,b^2).
$$
\end{proof}

\bigskip

If $L$ is a Lie algebra, we shall say
that a vector space $Z$ equipped with an
operation $\alpha:L\otimes Z\to Z$ is a (left) {\it Lie module over} $L$ ,
if for every $x,y\in L$ and $z\in Z$ we have
$$[x,y]\cdot z=x\cdot (y\cdot z)-y\cdot (x\cdot z)(-1)^{|x||y|},$$
where $x\cdot z=\alpha(x\otimes z)$.
This condition means that the map $\pi=\lambda^2(\alpha):L \to [Z,Z]$
defined by putting $\pi(x)(z)=\alpha(x\otimes z)$
is a homomorphism of Lie algebras $L\to gl(Z)$.
If $L$ and $P$ are Lie algebras, we shall say
that an action of $L$ on a vector space $Z$
{\it commutes} with an action of $P$ on $Z$
if we have
$$x\cdot (y\cdot z)=y\cdot (x\cdot z)(-1)^{|x||y|}$$
for every $x\in L$, $y\in P$ and $z\in Z$.
In which case the vector space $Z$
has the structure of a Lie module over the
Lie algebra $L\times P$ if we 
put $(x,y)\cdot z=x\cdot z+y\cdot z$.
The map
$\pi:L\times P\to gl(Z)$
defined by putting $\pi(x,y)(z)=x\cdot z +y\cdot z$
is then a homomorphism of Lie algebras.

\medskip

If $X$ and $Y$ are graded vector spaces, then 
the vector space $X\otimes Y$ has the structure
of a module over the algebra $End(X)\otimes End(Y)$ if we
put $(f\otimes g)(x\otimes y)=f(x)\otimes g(y)(-1)^{|g||x|}$
for $f\in End(X)$, $g\in End(Y)$, $x\in X$ and $y\in Y$.

\begin{prop}\label{Lietensor}
If $X$ and $Y$ are graded vector spaces,
then the map
$$\pi: gl(X) \times  gl(Y)  \to gl(X\otimes Y)$$
defined by putting $\pi(f,g)=f\otimes Y+X\otimes g$
is a homomorphism of Lie algebras.
The homomorphism preserves the square of odd elements.
\end{prop}

\begin{proof} 
This follows from definition \ref{defcommutation}, since the homomorphisms
$f\mapsto f\otimes Y$ and $g\mapsto X\otimes g$
commute.
\end{proof}

\medskip

If $L$ is a Lie algebra,
we shall say a vector space $Z$ equipped with an
operation $\beta:Z\otimes L\to Z$ is a 
 {\it right Lie module over}  $L$
if, for every $x,y\in L$ and $z\in Z$, we have
$$z\cdot [x,y]= (z\cdot x)\cdot y -(z\cdot y)\cdot x (-1)^{|x||y|},$$
where $z\cdot x=\beta(z\otimes x)$.
This condition means that the map $\pi=\lambda^1(\beta):L\to [Z,Z]$
defined by putting $\pi(x)(z)=\beta(z\otimes x)(-1)^{|z||x|}$
is an anti-homomorphism of Lie algebras $L\to gl(Z)$.
Every right $L$-module $(Z,\beta)$ has the structure
of a left $L$-module $(Z,\alpha)$ if we put $\alpha(x\otimes z):=-\beta(z\otimes x) (-1)^{|x||x|}$
for $x\in L$ and $z\in Z$.
If $L$ and $P$ are Lie algebras, we shall say
that a vector space $Z$ equipped with a left action by $L$
and a right action by $P$ is a {\it Lie bimodule} over the pair $(L,P)$
if the two actions commutes, meaning that 
we have 
$$(x\cdot z)\cdot y =x\cdot (z\cdot y)$$
for every $x\in L$, $z\in Z$ and $y\in P$.
In which case the vector space $Z$
has the structure of a (left) Lie module over the
Lie algebra $L\times P^o$ if we 
put $(x,y^o)\cdot z=x\cdot z+ z\cdot y(-1)^{|z||y|}$.

\begin{prop}\label{repLiebimodules}
If $X$ is a Lie bimodule over a pair of Lie algebras $(L,P)$, then the map
$\pi:L\times P\to gl(Z)$
defined by putting $\pi(x,y)(z)=x\cdot z -z\cdot y(-1)^{|z||y|}$
is a homomorphism of Lie algebras.
\end{prop}

\begin{proof} The map 
$\pi_0: L \times  P^o \to gl(Z)$
defined by putting $\pi_0(x,y^o)(z)=x\cdot z +z\cdot y(-1)^{|z||y|}$
is a homomorphism of Lie algebras, since
$Z$ is a Lie module over  the Lie algebra $L \times  P^o$.
But the map $y\mapsto -y^o$ is a homomorphism
of Lie algebras $P\to P^o$ by proposition \ref{Lieanti}, since the same map
is an anti-homomorphism of Lie algebras.
Hence the map
$\pi: L \times  P \to gl(Z)$
defined by putting $\pi(x,y)(z)=x\cdot z -z\cdot y(-1)^{|z||y|}$
is a homomorphism of Lie algebras.
\end{proof}

\medskip

If $X$ and $Y$ are graded vector spaces, then 
the vector space $[X,Y]$ has the structure
of bimodule over the pair of algebras $(End(Y),End(X))$.
The left and the right actions are defined respectively by the composition laws
$$[Y,Y]\otimes [X,Y]\to [X,Y] \et  [X,Y]\otimes  [X,X]\to [X,Y].$$

\begin{prop}\label{Liebimodules2}
If $X$ and $Y$ are graded vector spaces, 
then the map
$$\pi: gl(Y) \times  gl(X)  \to gl([X,Y])$$
defined by putting $\pi(g,f)=[X,g]-[f,Y]$
is a homomorphism of Lie algebras.
The homomorphism preserves the square of odd elements.
\end{prop}

\begin{proof} The vector space $[X,Y]$ has the structure of a bimodule over the pair of algebras $(End(Y),End(X))$.
It thus have the structure of a Lie bimodule over the pair of Lie algebras  $(gl(Y), gl(X))$.
The result then follows from proposition \ref{repLiebimodules}, since we have
$$\pi(g,f)(h)=[X,g](h)-[f,Y](h)=gh-hf(-1)^{|f||h|}$$
for every $h\in [X,Y]$.
We leave to the reader the verification that $\pi$
preserves the square of odd pairs $(g,f)$.
\end{proof}

\subsubsection{Primitive elements of bialgebras}

\begin{defi}\label{primitivedefinition2} 
If $Q$ is a bialgebra, we shall say that an element $x\in Q$ is {\it primitive}
if it is primitive with respect to the unit element $1\in Q$.
\end{defi}

Recall that for an algebra $A$, $Lie(A)$ is the Lie algebra $(A,[-,-])$ where $[-,-]$ is the commutator.

\begin{prop}
The vector space $\Prim(Q)$ of primitive elements of $Q$ is a Lie subalgebra of $Lie(Q)$.
If $x\in Q$ is odd and primitive, then $x^2$ is primitive.
\end{prop}
\begin{proof}
Let us show that the commutator $[x,y]$ of two primitive elements of a bialgebra $Q$ is primitive.
We have
\begin{eqnarray*}
\Delta(xy) &=& \Delta(x)\Delta(y)\\
 &=& (x\otimes 1+1\otimes x)(y\otimes 1+1\otimes y)\\
 &=& xy\otimes 1+1\otimes xy +x\otimes y + y\otimes x\ (-1)^{|x||y|}
\end{eqnarray*}
and
\begin{eqnarray*}
\Delta(yx) &=& yx\otimes 1+1\otimes yx +y\otimes x + x\otimes y\ (-1)^{|x||y|}.
\end{eqnarray*}
Hence
\begin{eqnarray*}
\Delta\left(xy-yx\ (-1)^{|x||y|}\right) &=& \Delta(xy)-\Delta(yx)\ (-1)^{|x||y|}\\
 &=& [x,y]\otimes 1+1\otimes [x,y].
\end{eqnarray*}
As for the second statement, the first comptuation shows that
\begin{eqnarray*}
\Delta(x^2) &=& x^2\otimes 1+1\otimes x^2 +x\otimes x + x\otimes x\ (-1)^{|x||x|}\\
&=& x^2\otimes 1+1\otimes x^2
\end{eqnarray*}
since $(-1)^{|x||x|}=-1$ when $x$ is odd.
\end{proof}

\begin{lemma}\label{functoprim}\label{derivationmorbialg}
Let $f:Q_1\to Q_1$ be a map of bialgebras, then $f$ induces by restriction a map of Lie algebras $f':\Prim(Q_1)\to \Prim(Q_2)$ which preserves the square of odd elements.
\end{lemma}
\begin{proof}
The image of a primitive element of $Q_1$ is a primitive element of $Q_2$ since $f(1)=1$. 
This prove that $f$ induces a map $f':\Prim(Q_1)\to \Prim(Q_2)$.
$f$ is a map of algebras so $f([x,y])=[f(x),f(y)]$ and $f(x^2)=f(x)^2$ for every $x,y\in Q_1$. 
This proves that $f'$ is a map of Lie algebras and preserves the square of odd elements.
\end{proof}

We shall call the map $f':\Prim(Q_1)\to \Prim(Q_2)$ the {\em derivative } of $f':Q_1\to Q_2$.

\medskip
In other words $\Prim$ defines a functor
$$
\Prim:\dgBialg \to \dgLie
$$
where $\dgLie $ denotes the category of (graded) Lie algebras. 
The next proposition shows that this functor is monoidal when $\dgLie$ is equipped with the product.

\begin{prop}\label{primitiveoftensorbialg}
If $Q_1$ and $Q_2$ are two bialgebras, then $Q_1\otimes Q_2$ is a bialgebra and the map 
$$
i:\Prim(Q_1)\times \Prim(Q_2)\to \Prim(Q_1\otimes Q_2)
$$
given by $i(x,y)=x\otimes 1+1\otimes y$ is an isomorphism of Lie algebras which preserves the square of odd elements.
\end{prop}

\begin{proof}
$i$ is a Lie algebra map by lemma \ref{functoprim} and an isomorphism by proposition \ref{primitiveoftensor}. The preservation of squares is proposition \ref{laxLie}.
\end{proof}

Recall from example \ref{primitiveHopf} the shuffle bialgebra $T^{sh}(\delta) = \FF[\delta]$.

\begin{prop}\label{univpropofdeltabialg}
If $Q$ is a bialgebra and $X$ a dg-vector space, then there are natural bijections between 
\begin{enumerate}
\item maps of dg-vector spaces $X\to \Prim(Q)$
\item maps of pointed dg-coalgebras $T^c_{\bullet,1}(X)\to Q$,
\item maps of $C$-bicomodules $X\to \Omega^Q$ (where $X$ is viewed as a $C$-bicomodule through $e:\FF\to C$),
\item maps of dg-vector spaces $X\to \Omega^{Q,1}$,
\item and maps of pointed dg-bialgebras $T^{csh}_\bullet(X)\to Q$,
\end{enumerate}
\end{prop}

\begin{proof}
Most of the assertions are from corollary  \ref{univpropofdeltanew}.
The equivalence $2\leftrightarrow 5$ is given by the adjunction $T^{qcsh}_\bullet:\dgCoalg_\bullet\rightleftarrows \dgBialg:U_m$ of proposition \ref{freebialgfromcoalg} and the isomorphism $T^{qcsh}_\bullet(T^c_{\bullet,1}(X))=T^{csh}(X)$.
\end{proof}

\subsubsection{Lie algebras of (co)derivations}\label{liealgebraofderivations}

This section construct the Lie algebra structure of derivations and coderivations and explicit their compatibility with the tensor products of algebras and coalgebras.

\bigskip

Recall from proposition \ref{semialgtoLie} that if $B$ is an algebra, $Lie(B)$ is the Lie algebra $(B,[-,-])$ where $[-,-]$ is the commutator in $B$.
Recall also from example \ref{theLiealgebraofendomorphism} that if $X$ is a vector space, $gl(X)$ is defined to be the Lie algebra $Lie(End(X))$.

\paragraph{Algebras}

Recall from definition \ref{defderivationalg} that $\Der(A)\subset [A,A]$ is the dg-vector space of graded derivations $A\rhup A$.

\begin{prop}\label{commutatorofderiv}
Let $A$ be a unital or non-unital algebra, then the dg-vector space $\Der(A)$ is a Lie sub-dg-algebra of $gl(A)$.
Moreover, the square $d^2$ (computed in $End(A)$) of an odd derivation is a derivation.
\end{prop}

\begin{proof} 
Let us show that if $d_1$ and $d_2$ are derivations of degree $n_1$ and $n_2$, then $[d_1,d_2]$ is a derivation of degree $n_1+n_2$.
To avoid too much signs, it is convenient to work at the level of graded morphisms rather than with elements.
Let $m:A\otimes A\to A$ be the product of $A$, we have
\begin{eqnarray*}
d_1d_2 m &=&d_1 m (d_2\otimes A +A\otimes d_2)  \\
&=&  m (d_1\otimes A +A\otimes d_1) (d_2\otimes A +A\otimes d_2)  \\
&=&m \bigl(d_1d_2\otimes A +d_2\otimes d_1 (-1)^{n_1n_2} +d_1\otimes d_2 +A\otimes d_1d_2\bigr)
\end{eqnarray*}
\begin{eqnarray*}
d_2d_1 m  &=&d_2 (d_1\otimes A +A\otimes d_1)  \\
&=& m (d_2 \otimes A +A\otimes d_2) (d_1\otimes A +A\otimes d_1)\\
&=&m \bigl(d_2d_1\otimes A +d_1\otimes d_2 (-1)^{n_1n_2} +d_2\otimes d_1 +A\otimes d_2d_1\bigr)
\end{eqnarray*}
Thus,
\begin{eqnarray*}
 [d_1,d_2] m &=& d_1d_2m   -d_2d_1 m  (-1)^{n_1n_2} \\
&=& m\bigl(d_1d_2\otimes A - d_2d_1\otimes A (-1)^{n_1n_2}+ A\otimes d_1d_2 -A\otimes d_2d_1  (-1)^{n_1n_2}\bigr)\\
&=& m\bigl([d_1,d_2]\otimes A+ A\otimes [d_1,d_2])
\end{eqnarray*}
This shows that $[d_1,d_2] $ is a derivation of degree $n_1+n_2$.
The dg-vector space $\Der(A)$ is thus a Lie sub-algebra of the Lie dg-algebra $gl(A)$

Let $d$ be a derivation of degree $n$ with $n$ odd, the first computation above with $d_1=d_2=d$ gives
\begin{eqnarray*}
d^2 m &=&m \bigl(d^2\otimes A +d\otimes d (-1)^{n^2} +d\otimes d +A\otimes d^2\bigr)\\
&=&m(d^2\otimes A +A\otimes d^2)
\end{eqnarray*}
where we have used that $(-1)^{n^2}=-1$. This shows that $d^2$ is a derivation.
\end{proof}

Recall from definition \ref{pointed derivationdef} that, for a pointed algebra $(A,e)$, $\Der_\bullet(A)\subset \Der(A)$ is the dg-vector space of pointed graded derivations $A\rhup A$.

\begin{prop}\label{pointedderivaresubliealg}
If $(A,\epsilon)$ is a pointed algebra, the graded vector space $\Der_\bullet(A)$ is a Lie sub-algebra of $\Der(A)$.
\end{prop}

\begin{proof}
$\Der_\bullet(A)\subset \Der(A)\subset End(A)$ by construction. 
Let $d_1$ and $d_2$ be two pointed derivations let us verify that $[d_1,d_2]$ is still pointed.
We have $\epsilon[d_1,d_2] = \epsilon d_1d_2 - \epsilon d_2d_1 \ (-1)^{|d_1||d_2|} = 0$.
This proves the stability for the Lie bracket.
\end{proof}

\bigskip

\begin{prop}\label{actiondertensor}
Let $A$ and $B$ be two algebras.
If $d_1:A\rhup A$ and $d_2:B\rhup B$ are derivation of degree $n$,
then the morphism $\varpi(d_1,d_2)=d_1\otimes B+A\otimes d_2$
is a derivation of degree $n$ of the 
algebra $A\otimes B$. 
Moreover, the map
$$
\varpi:\Der(A) \times \Der(B) \to \Der(A\otimes B)
$$
so defined is a homomorphism of Lie algebras which preserves
the square of odd pairs $(d_1,d_2)$.
\end{prop}

\begin{proof} 
If $d:A\rhup A$ is a derivation of degree $n$,
let us verify that the morphism $D=d\otimes B:A\otimes B\rhup A\otimes B$ is a derivation
of degree $n$. For every $x_1,x_2\in A$ and $y_1,y_2\in B$ 
we have
\begin{eqnarray*}
D((x_1\otimes y_1)(x_2\otimes y_2))&=&D(x_1x_2\otimes y_1y_2)(-1)^{|y_1||x_2|} \\
&=&d(x_1x_2)\otimes y_1y_2 (-1)^{|y_1||x_2|} \\
&=&d(x_1)x_2 \otimes y_1y_2 (-1)^{|y_1||x_2|} +x_1d(x_2) \otimes y_1y_2 (-1)^{n|x_1|+|y_1||x_2|} \\
&=&(d(x_1) \otimes y_1)(x_2\otimes y_2)  +(x_1\otimes y_1)(d(x_2)\otimes y_2) (-1)^{n|x_1|+n|y_1|} \\
&=&D(x_1 \otimes y_1)(x_2\otimes y_2)  +(x_1\otimes y_1)D(x_2\otimes y_2) (-1)^{n(|x_1|+|y_1|} \\
\end{eqnarray*}
This proves that $d\otimes B$ is a derivation of degree $n$.
Similarly, if $d:B\rhup B$ is a derivation of degree $n$,
then the morphism $A\otimes d: A\otimes B\rhup A\otimes B$ is a derivation
of degree $n$.  Thus, if $d_1:A\rhup A$ and $d_2:C\rhup C$ are derivation of degree $n$,
then the morphism $\varpi(d_1,d_2)=d_1\otimes A+A\otimes d_2$
is a derivation of degree $n$.
It then follows from proposition \ref{Lietensor} that the map
$\varpi :\Der(A) \times \Der(B) \to \Der(A\otimes B)$
is a homomorphism of Lie algebras which preserves
the square of odd pairs $(d_1,d_2)$.
\end{proof}

\paragraph{Coalgebras}

\begin{prop}\label{commutatorofcoderiv}
Let $C$ be a counital or non-counital coalgebra, then the graded vector space $\Coder(C)$ is a Lie sub-algebra of $gl(C)$.
Moreover, the square $d^2$ (computed in $End(C)$) of an odd coderivation is a coderivation.
\end{prop}

\begin{proof}
Dual to proposition \ref{commutatorofderiv}.
\end{proof}

\begin{prop}
If $C$ is a pointed coalgebra, the graded vector space $\Coder_\bullet(C)$ is a Lie sub-algebra of $\Coder(C)$.
\end{prop}

\begin{proof}
Dual to proposition \ref{pointedderivaresubliealg}.
\end{proof}

By proposition \ref{nilradicalcoder}, for a pointed coalgebra $C$, if we restrict a coderivation $D:C\rhup C$ to the radical of $C$, we obtain a coderivation $D^c:R^cC\rhup R^c C$.
\begin{lemma}
This defines a homomorphism of Lie algebras
$$
(-)^c:\Coder_\bullet(C)\to \Coder_\bullet(R^cC).
$$
\end{lemma}
\begin{proof}
The restriction commutes with the commutator of endomorphisms.
\end{proof}

If $C$ is a non-counital coalgebra, we have an analog Lie algebra homomorphism $(-)^c:\Coder(C)\to \Coder(R^cC)$.

\bigskip

\begin{prop}\label{actioncodertensor}
Let $C$ and $D$ be two coalgebras.
If $d_1:C\rhup C$ and $d_2:D\rhup D$ are coderivation of degree $n$,
then the morphism $\varpi(d_1,d_2)=d_1\otimes D+C\otimes d_2$
is a coderivation of degree $n$ of the 
coalgebra $C\otimes D$. 
Moreover, the map
$$
\varpi:\Coder(C) \times \Coder(D) \to \Coder(C\otimes D)
$$
so defined is a homomorphism of Lie algebras which preserves
the square of odd pairs $(d_1,d_2)$.
\end{prop}

\begin{proof} 
If $d:C\rhup C$ is a coderivation of degree $n$,
let us verify that the morphism $d\otimes D:C\otimes D\rhup C\otimes D$ is a coderivation
of degree $n$. We have
$\Delta_C d=(d\otimes C+C\otimes d)\Delta_C $,
since $d$ is a coderivation. Moreover, we have
$$\Delta_{C\otimes D} =(C\otimes \sigma \otimes  D) (\Delta_C \otimes \Delta_D)$$
by definition of the coproduct of $C\otimes D$.
Thus
\begin{eqnarray*}
\Delta_{C\otimes D}(d\otimes D) E &=&(C\otimes \sigma \otimes  D) (\Delta_C \otimes \Delta_D)(d\otimes C)\\
&=&(C\otimes \sigma \otimes  D) (\Delta_C d\otimes \Delta_D)\\
&=&(C\otimes \sigma \otimes  D) (d\otimes C\otimes D\otimes D+C\otimes d\otimes D\otimes D) (\Delta_C\otimes \Delta_D)\\
&=& (d\otimes D\otimes C\otimes D+C\otimes D\otimes d\otimes D)(C\otimes \sigma \otimes  D) (\Delta_C\otimes \Delta_D)\\
&=&\bigl( (d\otimes D)\otimes (C\otimes D)+(C\otimes D)\otimes (d\otimes D)\bigr)\Delta_{C\otimes D}.
\end{eqnarray*}
This shows that the morphism $d\otimes D$ is a coderivation of $C\otimes D$.
Similarly, if $d:D\rhup D$ is a coderivation of degree $n$,
then the morphism $C\otimes d: C\otimes D\rhup C\otimes D$ is a coderivation
of degree $n$.  Thus, if $d_1:C\rhup C$ and $d_2:D\rhup D$ are coderivation of degree $n$,
then the morphism $\varpi(d_1,d_2)=d_1\otimes D+C\otimes d_2$
is a coderivation of degree $n$.
It then follows from proposition \ref{Lietensor} that the map
$\varpi :\Coder(C) \oplus \Coder(D) \to \Coder(C\otimes D)$
is a homomorphism of Lie algebras which preserves
the square of odd pairs $(d_1,d_2)$.
\end{proof}

\newpage
\section{The category of coalgebras}\label{coalgebras}

This chapter will deal solely with dg-coalgebras and dg-vector spaces.
In order to simplify the langage we will often call dg-coalgebras and dg-vector spaces simply {\em coalgebras} and {\em vector spaces}.
More generally, we will often remove the "dg" prefix in the name of dg-subcoalgebras, dg-vector subspace... and a map of dg-vector spaces will simply be called a linear map or a map.

The chapter contains the following results.
\begin{itemize}

\item The category of dg-coalgebras $\dgCoalg$ is $\omega$-presentable, hence it is bicomplete (theorem \ref{finprescoalg}).
A dg-coalgebra is of $\omega$-compact if and only if it is finite dimensional (proposition \ref{finprescoalg}).

\item The forgetful functor $\dgCoalg\to\dgVect$ has a right adjoint $T^\vee: \dgVect\to \dgCoalg$ (theorem \ref{cofreecoalg}).

\item For every graded vector space we have $T(X)^\vee =T^\vee(X^\star)$. And we have $T^\vee(X)=T^c(X)$
when $X$ is graded finite and stricly positive or strictly negative (theorem \ref{tc=tv}).

\item The adjunction $U:\dgCoalg\rightleftarrows\dgVect:T^\vee$ is comonadic (theorem \ref{comonadic1}).

\item The category $\dgCoalg$ is symmetric monoidal closed (theorem \ref{homcoalg}). The internal hom will be noted $\HOM$.

\item The adjunction $U:\dgCoalg\rightleftarrows\dgVect:T^\vee$ is in fact strongly comonadic if both categories are enriched over $\dgCoalg$ (theorem \ref{strongcomonadicitycoalg}).

\end{itemize}

These results for coalgebras are sensibly harder to prove than their analog for algebras, see remark \ref{differencealgcoalg}.

\subsection{Presentability}

\begin{defi}
We shall say that a (dg-)subspace $V\subseteq C$ of a coalgebra $C=(C,\Delta, \epsilon)$ is a {\it sub-coalgebra} 
if we have $\Delta(V)\subseteq V\otimes V$.
\end{defi}

Recall that every coalgebra $C$ has the structure of a $C$-bicomodule.

\begin{lemma} \label{intersectionandtensor}
If $A$ and $B$ are respectively (dg-)subspaces of vector spaces $X$ and $Y$, then we have
$$
(A\otimes Y)\cap (X\otimes B)=A\otimes B.
$$
\end{lemma}

\begin{proof}
The intersection of dg-vector subspaces as graded subspaces is stable by the differential.
Thus, it is enough to prove the result at the level of the underlying graded vector spaces.
Let us forget the differential and pick $A'\subseteq X$ a complement of the subspace $A$ and $B'\subseteq Y$ a complement of $B$,
then we have a decomposition
$$
X\otimes Y=(A\oplus A')\otimes (B\oplus B') =(A\otimes B)\oplus (A'\otimes B) \oplus (A\otimes B')\oplus (A'\otimes B').
$$
The result follows, since $A\otimes Y=(A\otimes B)\oplus (A\otimes B')$ and $X\otimes B=(A\otimes B)\oplus (A'\otimes B)$.
\end{proof}

We shall say that a (dg-)subcoalgebra $E$ of a coalgebra $C$ is {\it generated} by a graded subset $S\subseteq E$ if $E$ is the smallest sub-dg-coalgebra of $C$ containing $S$.

\begin{prop}[\cite{Sw}] \label{subcoalgebravssuncomodule} 
If $C$ is a (dg-)coalgebra, then a subspace $E\subseteq C$ is a sub-coalgebra iff it is a sub-comodule of the $(C,C)$-comodule $C$.
Every graded subset $S\subseteq C$ generates a sub-dg-coalgebra $\overline S$ and coalgebra $\overline S$ is finite dimensional if $S$ is finite.
\end{prop}

\begin{proof} 
We have $E\otimes E=(C\otimes E)\cap  (E\otimes C)$ by lemma \ref{intersectionandtensor}.
It follows that $E$ is a sub-coalgebra iff it is closed under the left and right coactions of $C$ on itself.
The first statement is proved.
The second and third statements then follows from lemma \ref{subcomodule}.
\end{proof}

\begin{thm}[\cite{Sw}] \label{Calgunion}
Every coalgebra is the directed union (in $\dgVect$) of its finite dimensional sub-coalgebras.
\end{thm}

\begin{proof}
A graded subset $S\subseteq C$ generates a sub-coalgebra $\overline S$ by proposition \ref{subcoalgebravssuncomodule}.
Obviously, $C$ is the directed union of the sub-coalgebra $\overline S$, where $S$ is running in the finite graded subsets of $X$.
This proves the result, since the sub-coalgebra $\overline S$ is finite dimensional when $S$ is finite by proposition \ref{subcoalgebravssuncomodule}.
\end{proof}

\begin{prop} \label{coalgcocomp}
The category $\dgCoalg$ is cocomplete and the forgetful functor $U:\dgCoalg\to  \dgVect$ preserves and reflects colimits.
A finite colimit of finite coalgebras is finite.
\end{prop}
\begin{proof} 
Let us construct explicitly the equalizer of two maps $f,g:C \to D$ in the category $\dgCoalg$.
Let us denote by $p:C \to Z$ the coequalizer of $f$ and $g$ in the category $\dgVect$.
We shall prove that the dg-vector space $Z$ has the structure of a dg-coalgebra $(Z,\Delta, \epsilon)$. 
For this, consider the following two diagrams in the category $\dgVect$,
$$
\vcenter{
\xymatrix{
C\ar@<.6ex>[rr]^f\ar@<-.6ex>[rr]_g\ar[d]_{\Delta_C}&&D\ar[r]^p\ar[d]^{\Delta_D} & Z\ar@{-->}[d]^{\Delta}\\
C\otimes C \ar@<.6ex>[rr]^-{f\otimes f} \ar@<-.6ex>[rr]_-{g\otimes g} && D\otimes D \ar[r]_-{p\otimes p} &Z\otimes Z
}}
\et
\vcenter{
\xymatrix{
C\ar@<.6ex>[r]^-f\ar@<-.6ex>[r]_-g\ar[d]_{\epsilon_C}&D\ar[r]^p\ar[d]^{\epsilon_D} & Z\ar@{-->}[d]^{\epsilon}\\
\FF \ar@{=}[r] &\FF  \ar@{=}[r]&\FF
}}.$$
We have
$$
(p\otimes p)\Delta_Df=(p\otimes p)(f\otimes f)\Delta_C =(pf\otimes pf)\Delta_C =(pg\otimes pg)\Delta_C =(p\otimes p)(g\otimes g)\Delta_C =(p\otimes p)\Delta_Dg
$$
and
$$
\epsilon_Df=\epsilon_C =\epsilon_Dg
$$
since $f$ and $g$ are maps of coalgebras. 
Hence there exists a unique map $\Delta : Z \to Z\otimes Z$ in $\dgVect$ such that $\Delta p = (p\otimes p)\Delta_D$ 
and a unique map $\epsilon :Z\to \FF$ such that $\epsilon p=\epsilon_D$. 
We leave to the reader the verification that the pair $(\Delta, \epsilon)$ is a dg-coalgebra structure on $Z$. 
With this structure, $p : D \to Z$ becomes a map of dg-coalgebras; we leave to the reader the verification that $p$ is the coequalizer of $f$ and $g$ in the category $\dgCoalg$. 
This proves that the category admits equalizers and that the forgetful functor $U$ preserves coequalizers.
The proof that the category $\dgCoalg$ admits coproducts and that they are preserved by the functor $U$ is similar.
This shows that the category $\dgCoalg$ is cocomplete and the forgetful functor $U:\dgCoalg\to  \dgVect$ preserves colimits.
The forgetful functor also reflects colimits since it reflects isomorphisms.
The first statement of the proposition is proved.
The second statement follows, since a finite colimit of finite dimensional vector spaces is finite dimensional.
\end{proof}

\begin{cor} \label{Calgunion2} 
Every coalgebra is the directed colimit in $\dgCoalg$ of its finite dimensional sub-coalgebras. 
\end{cor}

\begin{proof}
We saw in theorem \ref{Calgunion} that a graded coalgebra is directed union of its finite dimensional sub-coalgebras.
The result follows since the forgetful functor $U:\dgCoalg\to  \dgVect$ reflects colimits by proposition \ref{coalgcocomp}.
\end{proof}

Recall that a (small) category $I$ is said to be {\it sifted} (resp. {\it directed}) if the colimit functor
$$
\underset{I}{\ccolim}: \Set^I\to \Set
$$
preserves finite products (resp. finite colimits). In particular a directed category is sifted.
A category $I$ is sifted if and only if it is non-empty and the diagonal functor $I\to I\times I$ is cofinal \cite{Lair, ARV}.
A category $I$ is said to be {\it cosifted} (resp. codirected) if the opposite category $I^{op}$ is sifted.

\begin{prop}\label{siftedtensor1}
If a category $I$ is sifted, then the colimit functor
$$
\underset{I}{\ccolim}: \dgVect^I\to \dgVect
$$
is monoidal.
If $I$ is cosifted and finite, then the limit functor
$$
\lim_{I}: \dgVect^I\to \dgVect
$$
is monoidal.
\end{prop}

\begin{proof} Let $A:I\to \dgVect$ and $B:I\to \dgVect $
be two diagrams indexed by the category $I$.  The canonical map
$$
\underset{(i,j)\in I\times I}{\ccolim} A_i\otimes B_j \to \underset{i\in I}{\ccolim} A_i \otimes \underset{j\in I}{\ccolim} B_j
$$
is an isomorphism, since the tensor product functor $\otimes$ is cocontinuous in each
variable. And the canonical map 
$$
\underset{i\in I}{\ccolim} A_i\otimes B_i\to \underset{(i,j)\in I\times I}{\ccolim} A_i\otimes B_j
$$
is an isomorphism, since the diagonal $I\to I\times I$ is cofinal when $I$ is sifted. This proves the first statement.

Let us now suppose that $I$ is cosifted and finite.
The canonical map
$$\lim_{i\in I} A_i \otimes \lim_{j\in I} B_j \to  \lim_{(i,j)\in I\times I} A_i\otimes B_j $$
is an isomorphism, since the tensor product functor $\otimes$ is exact in each
variable and the category $I$ is finite.
And the canonical map 
$$  \lim_{(i,j)\in I\times I} A_i\otimes B_j \to \lim_{i\in I} A_i\otimes B_i$$
is an isomorphism, since the diagonal $I\to I\times I$ is coinitial when $I$ is cosifted.
\end{proof}

\begin{rem}
The statement of proposition \ref{siftedtensor1} about cosifted limits uses the fact that $\otimes$ is left exact, \ie that we are working with over a field $\FF$ and not an arbitrary ring. This will be useful in the proof of the comonadicity theorem \ref{comonadic1}.
\end{rem}

\begin{cor}\label{directedtensorpower}
The tensor power functor $(-)^{\otimes n}:\dgVect \to \dgVect$ preserves sifted colimits for every $n\geq 0$.
\end{cor}
\begin{proof}
If $I$ is a sifted category, then the colimit functor $\colim_{I}: \dgVect^I\to \dgVect$ preserves tensor products by proposition \ref{siftedtensor1}.
Hence the colimit functor $\colim_{I}$ preserves tensor powers. 
It follows that the tensor power functor $(-)^{\otimes n}:\dgVect \to \dgVect$ commutes with sifted colimits.
\end{proof}

Recall that an object $X$ in a cocomplete category $\cal C$ is said to be $\omega$-{\it compact} if the functor ${\cal C}(X,-)$ commute with directed limits.
Recall also that a cocomplete category $\cal C$ is said to be $\omega$-{\it presentable} if the category of $\omega$-compact objects of $\cal C$  is essentially small and every objet of $\cal C$ is a colimit of $\omega$-compact objects.

\begin{thm}\label{finprescoalg}
The category $\dgCoalg$ is $\omega$-presentable.
A coalgebra is $\omega$-compact if and only if it is finite.
\end{thm}

\begin{proof}
Let us first show that if an $\omega$-compact coalgebra $C$ is finite.
The coalgebra $C$ is the directed colimits of the poset ${\cal F}(C)$ of its finite sub-coalgebras by corollary \ref{Calgunion2}.
The functor $\dgCoalg(C,-):\dgCoalg\to \Set$ preserves directed colimits, since $C$ is $\omega$-compact by hypothesis.
It follows that the identity element $1_C\in \dgCoalg(C,C)$ is in the image of the map $\dgCoalg(C,E)\to \dgCoalg(C,C)$ induced by the inclusion $i_E:E\to C$ for some finite sub-coalgebra $E\subseteq C$.
But the relation $1_C=i_E u$ implies that the inclusion $i_E$ is surjective; it is thus an isomorphism.
We have proved that an $\omega$-compact coalgebra is finite.
Conversely, let us show a finite coalgebra $E$ is $\omega$-compact.
For this, we need to show that the functor $\dgCoalg(E,-) :\dgCoalg \to \Set$ preserves directed colimits.
We shall use the fact that a finite limit of (set valued) functors preserving directed colimits preserves directed colimits.
If $X$ and $Y$ are dg-vector spaces, let us put $\Hom(X,Y)=\dgVect(X,Y)=Z_0[X,Y]$.
Let $U:\dgCoalg\to \dgVect$ be the forgetful functor.
We shall prove that the functor $\dgCoalg(E,-)$ is a finite limit of functors of the form $\Hom(U(E),U(-)^{\otimes n})$, and that these functors are preserving directed colimits.
Let us start with the latter.
The forgetful functor $U:\dgCoalg\to \dgVect$ is cocontinuous by proposition \ref{coalgcocomp}.
Hence the functor $U(-)^{\otimes n}:\dgCoalg\to \dgVect$ preserves directed colimits for every $n\geq 0$ by corollary \ref{directedtensorpower} 
(since a directed colimit is sifted).
It follows that the functor $\Hom(U(E),U(-)^{\otimes n}):\dgCoalg\to \Set$ preserves directed colimits for every $n\geq 0$,
since the dg-vector space $U(E)$ is finite.
It remains to show that the functor $\dgCoalg(E,-)$ is a finite limit of functors of the form $\Hom(U(E),U(-)^{\otimes n})$.
If $C$ is a coalgebra, then a linear map $f:E\to C$ is a map of coalgebras if and only if $(f\otimes f)\Delta_E=\Delta_Cf$ and $\epsilon_Cf=\epsilon_E$. 
Let us put $\alpha_C(f)=(f\otimes f)\Delta_E$, $\alpha'_C(f)=\Delta_Cf$, $\beta_C(f)=\epsilon_Cf$ and $\beta'_C(f)=\epsilon_E$.
This defines four natural transformations
$$\xymatrix{
\alpha,\alpha':\Hom(U(E),U(-))\ar@<.6ex>[r]\ar@<-.6ex>[r]& \Hom(U(E),U(-)^{\otimes 2})
}$$
and
$$\xymatrix{
\beta,\beta':\Hom(U(E),U(-))\ar@<.6ex>[r]\ar@<-.6ex>[r]& \Hom(U(E),U(-)^{\otimes 0}),
}$$
and $\dgCoalg(E,C)$ is their common equalizer  (see section \ref{commonequalizer}).
This shows that the functor $\dgCoalg(E,-)$ is a finite limit of functors of the form $\Hom(U(E),U(-)^{\otimes n})$.
We have proved that a finite coalgebra is $\omega$-compact.
Let us now prove that the category $\dgCoalg$ is $\omega$-presentable.
The category $\dgCoalg$ is cocomplete by proposition \ref{coalgcocomp}. 
Obviously, the category of finite coalgebras is essentially small.
Hence the category of $\omega$-compact coalgebras is essentially small.
By corollary \ref{Calgunion2}, every coalgebra is a directed colimit of finite dimensional coalgebras, \ie of $\omega$-compact coalgebras. 
\end{proof}

\begin{rem}
It follows from proposition \ref{dualfinitealg} that the category $\dgCoalg^{op}=({\sf Ind}(\dgCoalgfin))^{op}$ is equivalent to the category ${\sf Pro}(\dgAlgfin)$ of (strict) pro-finite algebras.
\end{rem}

\begin{cor}\label{catcoalgcomplete}
The category $\dgCoalg$ is complete. A contravariant functor $\dgCoalg^{op}\to \Set$ is representable if and only if
it is continuous. A functor $U:\dgCoalg \to \cal C$ with codomain a cocomplete category  $\cal C$ has a 
right adjoint $R:  \cal C \to  \dgCoalg$ if and only if it is cocontinuous; the right adjoint preserves directed colimits
if and only if the functor  $U$ takes a $\omega$-compact object to a $\omega$-compact object.
\end{cor}

\begin{proof} The category $\dgCoalg$ is $\omega$-presentable by theorem \ref{finprescoalg}. 
The statements in the proposition are true for any $\omega$-presentable category \cite{AR}.
\end{proof}

\subsection{Cofree coalgebras}

\begin{defi}
If $C$ is a (dg-)coalgebra and $V$ is a (dg-)vector space, we shall that $C$ is {\em cofreely cogenerated by a linear map $p:C\to V$ (or by $V$)} if for any coalgebra $E$ and any linear map $f:E\to V$, there exists a unique map of coalgebras $g:E\to C$ such that $pg=f$.
$$\xymatrix{
E\ar@{-->}[rr]^-g \ar[rrd]_f && C\ar[d]^p\\
&& V
}$$
\end{defi}

\begin{thm}\label{cofreecoalg}
The forgetful functor $U:\dgCoalg\to \dgVect$ has a right adjoint $T^\vee$ which takes a vector space $X$ to the cofree coalgebra $T^\vee(X)$ cogenerated by $X$; the counit of the adjunction $U\dashv T^\vee$ is the cogenerating map $p:T^\vee(X) \to X$.
The functor $T^\vee$ preserves directed colimits.
\end{thm}

\begin{proof}
The functor $U:\dgCoalg\to \dgVect$ is cocontinous by proposition \ref{coalgcocomp}.
Hence it has a right adjoint $T^\vee$ by corollary \ref{catcoalgcomplete}.
It is easy to see that the coalgebra $T^\vee(X)$ is cofreely cogenerated by the counit $UT^\vee(X) \to X$ of the adjunction $U\dashv T^\vee$.
The functor $T^\vee$ preserves directed colimits by corollary \ref{catcoalgcomplete}, since the functor $U$ takes $\omega$-compact coalgebras to $\omega$-compact (=finite dimensional) vector spaces by theorem \ref{finprescoalg}.
\end{proof}

We will not explicit the construction of $T^\vee$ but let us say that $T^\vee(X)$ can be realized as a certain subspace of the completed tensor algebra $T^\wedge (X) = \prod_n X^{\otimes n}$ containing the tensor coalgebra $T^c(X)$ (see \cite{BL,Ha, Smith}).
Also we shall see in theorem \ref{tc=tv} that $ T^\vee(X)=T^c(X)$ when $X$ is strictly positive or negative.

\begin{prop}\label{atomcofree}
The atoms of $T^\vee(X)$ are in bijection with $Z_0(X)$.
\end{prop}
\begin{proof}
An atom if a coalgebra maps $\FF\to T^\vee(X)$ which are in bijection with linear map $\FF\to X$.
\end{proof}

\subsubsection{Generation and separation}

We finish this section by some definitions.

\begin{defi}\label{coseparation}
Let $C$ be a coalgebra and $V$ a vector space
\begin{enumerate}
\item we shall say that a linear map $f:C\to V$ is {\it cogenerating} if the associated coalgebra map $g:C\to T^{\vee}(V)$ is injective.
\item we shall say that a linear map $f:C\to V$ is {\it separating} if the implication $fu=fv\Rightarrow u=v$ is true for any pair of coalgebra maps $u,v:E\to C$. 
\end{enumerate}
\end{defi}

Recall that if $C$ is a coalgebra and $V$ is a vector space,
then every linear map $f:C\to V$ can be coextended
as a coalgebra map $g:C\to T^{\vee}(V)$;
the map  $f:C\to V$ is separating if and only if
the map $g:C\to T^{\vee}(V)$ is a monomorphism of coalgebras. 
In particular, any cogenerating map is separating.

If $T^\vee(V)$ is a cofree coalgebra, we shall call the canonical cogenerating map $p:T^\vee(V)\to V$ the {\em cofree map}.
It is in particular separating.

\subsection{Applications of the cofree functor}

\subsubsection{Non-conilpotent quasi-shuffle}\label{quasishuffleTv}

As promised in section \ref{nonbiunitalbialgebras} we construct an adjoint to the functor $U_\Delta:\dgBialg\to \dgCoalg_\bullet$.
We only have to copy the definition of quasi-shuffle with $T^\vee$ instead of $T^c$.

$T^\vee$ is a right adjoint, then it sends the terminal object $0\in \dgVect$ to the terminal objet $T^\vee(0)=\FF\in \dgCoalg$.
For a vector space $X$ let $e:\FF\to T^\vee(X)$ be the map image of the zero map $0\to X$.
We define a pointed coalgebra $T^\vee_\bullet(X)=(T^\vee(X),\Delta, \epsilon, e)$.

Let  $A=(A,m)$ is a non-unital algebra, we consider $T^\vee_\bullet(A)=(T^\vee(A),\Delta, \epsilon, e)$ and note 
$p:T^\vee(A)\to A$ the cogenerating map.

\begin{prop}\label{cofreebialgebra11} 
There is a unique map of coalgebras $\mu:T^\vee(A)\otimes T^\vee(A) \to T^\vee(A)$ such that
$$
p\mu(x,y)=p(x)p(y)+\epsilon(x)p(y)+p(x)\epsilon(y)
$$
for every $x, y \in T^\vee(A)$.
The product $\mu$ is associative and the pair $(\mu,e)$ is defining a bialgebra structure on the coalgebra $T^\vee(A)=(T^\vee(A),\Delta, \epsilon)$.
\end{prop}

\begin{proof} 
Let us show that the operation $\mu$ is associative. We need to prove that the following square commutes
$$\xymatrix{
T^\vee(A)\otimes T^\vee(A)\otimes T^\vee(A) \ar[rr]^-{\mu \otimes Id} \ar[d]_{Id \otimes \mu} && T^\vee(A)\otimes T^\vee(A)\ar[d]^{\mu}\\
T^\vee(A)\otimes T^\vee(A) \ar[rr]^-{\mu} &&  T^\vee(A).
}$$
It suffices to show that the square commutes after composition with  $p:T^\vee(A)\to A$,
since $p$ is cogenerating and since every map in the square is a coalgebra map.
For this it suffices to show  that the square commutes after composition with
the map  $q:T^\vee(A)\to A_{+}=\FF\oplus A$ defined by putting $q(x)=(\epsilon(x),p(x))$,
since $p=p_2q$.
Let us show that we have 
$q\mu(Id \otimes \mu)=q\mu(\mu \otimes Id)$.
Let us put $xy:=\mu(x\otimes y)$ for every $x,y\in T^\vee(A)$.
We have $\epsilon(xy)=\epsilon(x)\epsilon(y)$,
since $\mu$ is a coalgebra map. Thus,
$$q(xy)=(\epsilon(xy),p(xy))=(\epsilon(x)\epsilon(y), p(x)p(y)+\epsilon(x)p(y)+p(x)\epsilon(y))=(\epsilon(x),p(x))(\epsilon(y),p(y))=q(x)q(y).$$
It follows that we have
$$q(x(yz))  =  q(x)q(yz)= q(x)q(y)q(z)= q(xy)q(z)= q((xy)z)$$
for every $x, y, z \in T^\vee(A)$. 
Thus, $q\mu(Id \otimes \mu)=q\mu(\mu \otimes Id)$
and it follows that $\mu(Id \otimes \mu)=\mu(\mu \otimes Id)$.
Let us now show that $e$ is a unit for the operation $\mu$.
 We need to prove that the following two triangles commute,
$$\xymatrix{
T^\vee(A) \ar[rr]^-{e \otimes Id} \ar@{=}[drr] && T^\vee(A)\otimes T^\vee(A)\ar[d]^{\mu}\\
 &&  T^\vee(A),
}\quad \quad
\xymatrix{
T^\vee(A) \ar[rr]^-{Id \otimes e} \ar@{=}[drr] && T^\vee(A)\otimes T^\vee(A)\ar[d]^{\mu}\\
 &&  T^\vee(A).
}$$
It suffices to show that the triangles are commuting after composition with  $q:T^\vee(A)\to A_{+}$.
For every $x \in T^\vee(A)$ we have
$q(ex)=q(e) q(x)=q(x)$, since $q(e)=(\epsilon(e), p(e))=(1,0)$ is the unit element of the algebra $A_{+}$.
Thus, $q\mu (e\otimes Id)=q$
and it follows that $\mu (e\otimes Id)=Id$.
Similarly, $\mu(Id \otimes e)=Id$.
\end{proof}

If $A=(A,m,e_A, \epsilon_A)$ is a pointed algebra, then $A_{-}=\ker(\epsilon_A)$ is non-unital algebra and we have $A=\FF e_A\oplus A_{-}$.
If $p:T^\vee(A_{-})\to A_{-}$ is the cogenerating map, let us denote by $q:T^\vee(A_{-})\to A$ the map defined by putting $q(x)=\epsilon_A(x)e_A+p(x)$.

\begin{prop}\label{cofreebialgebra3}\label{adjointforgetbialgalg}
The forgetful functor  $U_\Delta:\dgBialg \to \gAlg_\bullet$ admits a right adjoint which associates to a pointed algebra $A$ the bialgebra $(T^\vee(A_{-}),\mu, \Delta)$. Moreover, the canonical map  $q:T^\vee(A_{-}) \to A$ is the counit of the adjunction.
\end{prop}

\begin{proof} By lemma \ref{cofreebialgebra11}, 
the coalgebra $T^\vee(A)$ admits a bialgebra structure $(\mu, e)$
such that $$p\mu(x\otimes y)=p(x)p(y)+\epsilon(x)p(y)+p(x)\epsilon(y)$$
for every $x, y \in T^\vee(A)$ and such that $pe=0$.
It follows from these identities that $q:T^\vee(A)\to A$ is an algebra map.
The map $q$ is respecting the augmentations,
since we have $\epsilon_A q(x)=\epsilon_A (\epsilon(x)e_A+p(x))=\epsilon(x)$
for every $x \in T^\vee(A)$.
If $B$ is a bialgebra, and $f:B\to A$ is a map
of pointed algebras, let us show that there is a unique map of bialgebras
$g:B\to T^\vee(A)$ such that $qg=f$.
Let $f':B\to A_{-}$ be the map defined by putting
$f'(x) =f(x)-\epsilon_B(x)e_A$ for $x\in B$.
There is a unique coalgebra map $g:B\to T^\vee(A_{-})$ such that
$pg=f'$, since the coalgebra $T^\vee(A_{-})$ is cofree.
We then have $qg(x)=\epsilon(g(x))e_A+p(g(x))=\epsilon_B(x)e_A+f(x)-\epsilon_B(x)e_A=f(x)$.
Let us show that $g$ is a map of algebras.
For this we need to show that the following square commutes
$$\xymatrix{
B\otimes B\ar[rr]^(0.4){g\otimes g} \ar[d]_{m_B}&& T^\vee(A_{-})\otimes T^\vee(A_{-})\ar[d]^\mu\\
B\ar[rr]^g && T^\vee(A_{-})
}$$
The diagram commutes after composing it with  $q:T^\vee(A_{-})\to A$,
since we have
$qg(xy) = f(xy)= f(x)f(y)=qg(x)qg(y)$
for  every $x\in T^\vee(A_{-})$.
It thus commutes after composing it with $p:T^\vee(A_{-})\to A_{-}$.
It follows that the diagram commutes, since the map $p$ is cogenerating.
It remains to show that $g(e_B)=1$.
But we have $\epsilon(ge_B)=1$, since $ge_B:\FF \to T^\vee(A_{-})$
is a map of coalgebras.
Thus, $qge_B=\epsilon(ge_B)e_A+pge_B= e_A+f'(e_B)=e_A$,
since $f'(e_B)=f(e_B)-\epsilon_B(e_B)e_A=e_A-e_A=0$.
Thus $pg(e_B)=0=pe$ and it follows that $ge_B=e$,
since the map $p$ is cogenerating.
\end{proof} 

\begin{cor}\label{shuffleproductTv}
For $X$ a (dg-)vector space, the cofree coalgebra $T^\vee(X)$ is naturally a bialgebra.
\end{cor}
\begin{proof}
Consider $X$ as a non-unital algebra with zero multiplication.
\end{proof}

The product of $T^\vee(X)$ can be called the {\em shuffle product} as the algebra $R^cT^\vee(X)$ is $T^c(X)$ equipped with the shuffle product.

\subsubsection{Coderivations}

Let us denote by  $p:T^{\vee}(X)\to X$ the cogenerating map
of the cofree coalgebra on $X$.

\begin{lemma} \label{coextensionofcoderlemma}
If $X$ is a graded vector space and $N$ is a $T^\vee(X)$
bicomodule, then for every linear morphism $\phi:N\rhup X$ of degree $n$
can be coextended uniquely as a coderivation  $D:N \rhup T^{\vee}(X)$
of degree $n$. 
\end{lemma}

\begin{proof} Let us first suppose that $n=0$.
Let $f: T^{\vee}(X) \oplus N  \to X$ be the map
defined by putting $f(x,y)=p(x)+\phi(y)$ for every $(x,y)\in [T^{\vee}(X) \oplus N]$.
There is a unique map of coalgebras $g:[T^{\vee}(X) \oplus N] \to T^\vee(X)$
such that $pg=f$, since $p$ is cogenerating freely the coalgebra $T^{\vee}(X)$.
If  $i_1:T^\vee(X) \to  T^{\vee}(X) \oplus N $
is the inclusion, then we have $pgi_1=fi_1=p=p(id)$.
But $gi_1: T^{\vee}(X) \to T^{\vee}(X)$ is a coalgebra map,
since $i_1$ is a coalgebra map. Thus $gi_1=id$, since the map $p$
is cogenerating. It then follows from lemma
\ref{fcoderiv} that we have $g(x,y)=x+D(y)$, where $D:N\to  T^{\vee}(X)$ is a coderivation of degree 0.
We have $pD=\phi$, since we have $pg=f$. The uniquness of $D$
is left to the reader. 
Let us now consider the case of a morphism $\phi:N\rhup X$ of general degree $n$.
In this case, the morphism $\phi s^{-n}:S^{n}N\rhup X$ defined by putting
$\phi s^{-n}(s^n x)=\phi(x)$ for $x\in N$ has degree 0.  It can thus be coextended uniquely as a 
coderivation $Ds^{-n}:S^{n}N \to T^{\vee}(X) $ of degree 0.
The resulting morphism $D:N\rhup T^{\vee}(X) $ is a coderivation of degree $n$
which is extending $\phi$. The uniqueness of $D$ is clear.
 \end{proof}

\bigskip

The vector space $ T^\vee(X)\otimes X\otimes T^\vee(X)$ has the structure of a bicomodule
over the coalgebra $T^\vee(X)$. The bicomodule is cofreely cogenerated
by the map  $\phi=\epsilon \otimes X\otimes \epsilon: T^\vee(X)\otimes X\otimes T^\vee(X)\to X$.
The map $\phi$ can be coextended uniquely as a coderivation $d:T^\vee(X)\otimes X\otimes T^\vee(X)\to T^\vee(X)$
by lemma \ref{coextensionofcoderlemma}.

\begin{prop} \label{univcodercofree}
The coderivation $d:T^\vee(X)\otimes X\otimes T^\vee(X)\to T^\vee(X)$  defined above is couniversal.
Hence we have 
$$
\Omega^{T^\vee(X)}= T^\vee(X)\otimes X\otimes T^\vee(X).
$$
\end{prop}

\begin{proof} Similar to the proof of proposition \ref{univcodercofreeconil}.
\end{proof}

\begin{cor}\label{primcofree}
Let $p:T^\vee(X)\to X$ be the cogenerating map, then the map $D\mapsto pD$ induces an isomorphism of vector spaces
$$
\Coder(T^\vee(X)) = [T^\vee(X),X].
$$
In particular, a coderivation $D$ is zero iff $pD=0$.
\end{cor}
\begin{proof}
By proposition \ref{univcodercofree}, we have $\Omega^{T^\vee(X)}= T^\vee(X)\otimes X\otimes T^\vee(X)$.
The result follows from $\Omega^{T^\vee(X),f}=\Omega^{T^\vee(X)}\otimes^{(T^\vee(X))^{op}\otimes T^\vee(X)}D^{op}\otimes D=D\otimes X\otimes D$.
The second statement is the case where $D=\FF$. 
The isomorphism $\Prim(T^\vee(X),e)=\Omega^{T^\vee(X),e}$ is corollary \ref{prim=coder}.
\end{proof}

As a consequence, if $e:\FF\to T^\vee(X)$ is an atom of $T^\vee(X)$, we have $\Prim(T^\vee(X),e)=\Omega^{T^\vee(X),e}=X$, \ie the tangent space of $T^\vee(X)$ at $e$ is $X$.

\medskip
Recall the notion of a pointed coderivation from \ref{pointedcoderivationdef}. 
The coalgebra $T^\vee(X)$ is naturally pointed by the atom $0\in X$, we shall denote this pointed coalgebra $T_\bullet^\vee(X)$
and $T_\circ^\vee(X)$ the corresponding non-counital coalgebra.
The result above has a pointed version, which we state without proof:

\begin{prop}\label{primcofreepointed} \label{coextensionofcoderpointed}
Let $p:T_\bullet^\vee(X)\to X$ be the cogenerating map, then the map $D\mapsto pD$ induces an isomorphism of vector spaces
$$
\Coder(T_\bullet^\vee(X)) = [T_\circ^\vee(X),X].
$$
In particular, a poitned coderivation $D$ is zero iff $pD=0$.
\end{prop}

\subsection{Comonadicity}

We shall prove in this section that the category $\dgCoalg$ is comonadic over the category of dg-vector spaces $\dgVect$
(theorem \ref{comonadic1}). This result uses the fact that we are working over a base field $\FF$ and not an arbitray ring.

\bigskip

Let $\Delta(1)$ be the full subcategory of the simplical category $\Delta$
whose object are $[0]$ and $[1]$. 
$$\xymatrix{
 [0]\ar@<.6ex>[rr]^-{d_0} \ar@<-.6ex>[rr]_-{d_1}  && [1]\ar@/_2pc/[ll]_-{s_0}
}$$
We have $s_0d_0=1_{[0]}=s_0d_1$.
Recall that a {\it reflexive graph} in a category $\cal C$ is a truncated simplicial object
$X:\Delta(1)^{op}\to \cal C$. 
$$\xymatrix{
X_1\ar@<.6ex>[rr]^-{\partial_0} \ar@<-.6ex>[rr]_-{\partial_1}  && X_0\ar@/_2pc/[ll]_-{\sigma_0}
}$$
By definition, $\partial_0\sigma_0=id =\partial_1\sigma_0$.
A coequalizer of the pair $(\partial_0,\partial_1)$ is said to be {\it reflexive}.
For example, if $M=(M,\mu, \eta)$ is a monad on a category $\cal C$, and $(A,a)$
is a $M$-algebra, then
$$\xymatrix{
M^2(A)  \ar@<.6ex>[rr]^-{\mu_A} \ar@<-.6ex>[rr]_-{M(a)}  && M(A)\ar@/_2pc/[ll]_-{M(\eta_A)} \ar[r]^-a & A
}$$
is a reflexive coequalizer diagram in the category of $M$-algebras. 
It shows that every $M$-algebra is a reflexive equalizer of free $M$-algebras.

\medskip

Dually, a {\it reflexive cograph} in a category $\cal C$ is a truncated cosimplicial object
$X:\Delta(1)\to \cal C$,
$$\xymatrix{
 X^0\ar@<.6ex>[rr]^-{\partial^0} \ar@<-.6ex>[rr]_-{\partial^1}  && X^1\ar@/_2pc/[ll]_-{\sigma^0}
}$$
By definition, $\sigma^0 \partial^0=id =\sigma^0\partial^1$.
And equalizer of the pair $(\partial^0,\partial^1)$ is said to be {\it reflexive}.
For example,  if $T=(T,\delta, \epsilon)$ is a comonad on a category $\cal C$, and $(C,\rho)$
is a $T$-coalgebra, then
$$\xymatrix{
C\ar[r]^-{\rho} & T(C) \ar@<.6ex>[rr]^-{T(\rho)} \ar@<-.6ex>[rr]_-{\delta_C}  && T^2(C)\ar@/_2pc/[ll]_-{T(\epsilon_C)}
}$$
is a reflexive equalizer diagram in the category of $T$-coalgebras.
It shows that every $T$-coalgebra is a reflexive coequalizer of free $T$-coalgebras.

\medskip

\begin{lemma} The category $\Delta(1)$ is cosifted.
\end{lemma}

\begin{proof} 
The colimit of a functor $X:\Delta(1)^{op}\to \Set$
is the coequalizer of the pair of maps $\partial_0,\partial_1:X_1\to X_0$;
it is thus
the set  $\pi_0 X$ of connected components of the graph $X$.
It is easy to verify that the functor
$$\pi_0:[\Delta(1)^{op},\Set]\to \Set$$
preserves finite cartesian products.
\end{proof}

The colimit of a reflexive graph $X:\Delta(1)^{op}\to \dgVect$
is the cokernel of the boundary map $\partial =\partial_0-\partial_1:X_1\to X_0$.
We shall denote this cokernel by $h_0(X)$.
Dually, the limit of a reflexive cograph $X:\Delta(1)\to \dgVect$
is the kernel of the coboundary map $\partial =\partial^0-\partial^1:X^0\to X^1$.
We shall denote this cokernel by $h^0(X)$.

\begin{prop}\label{siftedtensor2}
The functors
$$h_0: [\Delta(1)^{op}, \dgVect] \to \dgVect \quad \quad {\rm and}\quad \quad  h^0: [\Delta(1), \dgVect] \to \dgVect $$
are monoidal.  
\end{prop}

\begin{proof} This follows from proposition \ref{siftedtensor1}.
\end{proof}

\begin{thm}\label{comonadic1}
The forgetful functor $U:\dgCoalg\to  \dgVect$ preserves and reflects reflexive equalizers.
The adjunction $U:\dgCoalg\rightleftarrows  \dgVect:T^\vee$ is comonadic.
\end{thm}

\begin{proof} Let us first show that the functor $U$ preserves reflexive equalizers.
Let 
$$\xymatrix{
 C\ar@<.6ex>[rr]^-{f} \ar@<-.6ex>[rr]_-{g}  && D\ar@/_2pc/[ll]_-{\sigma^0}
}$$
be a reflexive cograph in the category $\dgCoalg$ and let us denote by $h:E\to C$
the equalizer of $(f,g)$ in the category $\dgVect$,
$$\xymatrix{
E\ar[r]^h & C\ar@<.6ex>[r]^-{f} \ar@<-.6ex>[r]_-{g}  & D.
}$$
The diagram 
$$\xymatrix{
E\otimes E\ar[r]^-{h\otimes h} & C\otimes C\ar@<.6ex>[rr]^-{f\otimes f} \ar@<-.6ex>[rr]_-{g\otimes g} && D\otimes D
}$$
is also an equalizer in $\dgVect$ by lemma \ref{siftedtensor2}.
 In consequence, there exists a linear map $\Delta:E\to E\otimes E$ such that the following diagram commutes,
$$\xymatrix{
E\ar[rr]^h \ar[d]^{\Delta}&& C\ar@<.6ex>[rr]^-f \ar@<-.6ex>[rr]_-g \ar[d]^{\Delta_C}&& D\ar[d]^{\Delta_D}\\
E\otimes E\ar[rr]_-{h\otimes h} && C\otimes C \ar@<.6ex>[rr]^-{f\otimes f} \ar@<-.6ex>[rr]_-{g\otimes g} && D\otimes D
}$$
Also, there exists a unique linear map $\epsilon: E\to \FF$ such that the following diagram commutes,
$$\xymatrix{
E\ar[r]^h \ar[d]^{\epsilon}& C\ar@<.6ex>[r]^-f \ar@<-.6ex>[r]_-g \ar[d]^{\epsilon_C}& D\ar[d]^{\epsilon_D}\\
\FF\ar@{=}[r] & \FF \ar@{=}[r] &\FF.
}$$
We leave to the reader the proof that $(E,\Delta,\epsilon)$ is a coalgebra structure making $E$ into the equalizer of $f$ and $g$ in 
$\dgCoalg$.
This shows that the functor $U$ preserves reflexive coequalizers.
It also reflects reflexive equalizers, since it is conservative.
It then follows from the "crude" comonadicity theorem \cite{Barr-Wells}
that the adjunction $U\dashv T^\vee$ is comonadic.
\end{proof}

\begin{rem}
The proof of theorem \ref{comonadic1} is the only place where we used explicitely the field nature of the base ring (through the exactness of the tensor product). However, the result is proven without it in \cite[prop. 2.7]{Porst}.
\end{rem}

\subsection{Internal hom}\label{closedmonoidalstructurecoalg}

Recall from section \ref{monoidalstructurecoalgebras} that the tensor product of coalgebras gives the category $\dgCoalg$ a symmetric monoidal structure in which the unit object is the coalgebra $\mathbb{F}$.

\begin{thm}\label{homcoalg}
The symmetric monoidal category $(\dgCoalg,\otimes,\FF)$ is closed.
\end{thm}

\begin{proof}
We have to show that the functor $(-)\otimes C: \dgCoalg \to \dgCoalg$ admits a right adjoint for any coalgebra $C$.
By proposition \ref{catcoalgcomplete}, it suffices to show that the functor $(-)\otimes C$ is cocontinuous, since the category  $\dgCoalg$ is
$\omega$-presentable by proposition \ref{finprescoalg}. 
The forgetful functor $U:\dgCoalg \to \dgVect$ preserves and reflects colimits by proposition \ref{coalgcocomp}.
Hence it suffices to show that the functor $U(-\otimes C):\dgCoalg \to \dgVect$ is cocontinous. 
But we have $U(D\otimes C)=U(D)\otimes U(C)$ for every coalgebra $C$ and $D$ by definition of the tensor product of coalgebras.
Thus, $U(-\otimes C)=U(-)\otimes U(C)$.
But the functor $(-)\otimes  U(C):\dgVect \to \dgVect$ is cocontinuous, since the category $\dgVect$ is closed.
It follows that the functor $U(-)\otimes U(C)$ is cocontinuous, since the functor $U$ is cocontinuous.
\end{proof}

We shall denote the hom object between two coalgebras by $\HOM(C,D)$ and the endomorphism object of a coalgebra $C$ by $\END(C)$.
By definition, the coalgebra $\HOM(C,D)$ is equipped with a map of coalgebras, called the {\em strong evaluation},  
${\bf ev}:\HOM(C,D)\otimes C\to D$ which is couniversal in the following sense:
for any coalgebra $E$ and any map of coalgebras $f:E\otimes C\to D$ there exists a unique map of coalgebras $g:E\to \HOM(C,D)$ such that ${\bf ev}(g\otimes C)=f$.
$$
\xymatrix{
&& \HOM(C, D) \otimes C  \ar[d]^-{{\bf ev}} \\
E \otimes C \ar@{-->}[urr]^-{g\otimes C}  \ar[rr]^-{f}    && D
}$$
We shall denote the morphism $g$ by ${\Lambda}^2(f)$.
In addition, we shall put ${\Lambda}^1(f)={\Lambda}^2(f\sigma):C\to \HOM(E,D)$,
where  $\sigma:C\otimes E\to E\otimes C$ is the symmetry.
More generally, to any map of coalgebras $f:E_1\otimes \cdots \otimes E_n \to C$ 
we can associate a canonical map
$${\Lambda}^k(f):E_1\otimes \cdots \otimes \Hat{E_{k}} \otimes \cdots \otimes E_n \to \HOM(E_k,C)$$
 for each $1\leq k\leq n$. We may also use a canonical map
 $${\Lambda}^{k,r}(f)= E_1\otimes \cdots \otimes \Hat{E_{k}} \otimes \cdots  \otimes \Hat{E_{r}} \otimes \cdots  \otimes E_n \to \HOM(E_k\otimes E_r,C)$$
for $1\leq k<r\leq n$.

The notion of strong limit is given in appendix \ref{stronglimit}.
\begin{prop}\label{stronglimcoalg}
All ordinary (co)limits in $\dgCoalg$ are strong.
\end{prop}
\begin{proof}
This is a formal consequence of $\dgCoalg$ being bicomplete.
We shall prove only the result for limits, the proof for colimits is similar.
Let $C:I\to \dgCoalg$ be a diagram with limit $D$, then $D$ can be writtent as the equalizer in $\dgCoalg$
$$\xymatrix{
D\ar[r]&\prod_iC_i\ar@<.6ex>[r]\ar@<-.6ex>[r] & \prod_{i\to j}C_i.
}$$
$D$ is a strong limit, if for any coalgebra $E$, we have an equalizer in $\dgCoalg$
$$\xymatrix{
\HOM(E,D)\ar[r]&\prod_i\HOM(E,C_i)\ar@<.6ex>[r]\ar@<-.6ex>[r] & \prod_{i\to j}\HOM(E,C_i).
}$$
But this can be deduced form the previous equalizer and the fact that $\HOM(E,-)$ being right adjoint to $E\otimes-$, it commutes with equalizers.
\end{proof} 

\bigskip

The universal property of $\HOM(C,D)$ can be stated another way, using a notion analog to the comeasurings that will be used to define the Sweedler hom $\{-,-\}$ in section \ref{sweedlerhom}.

For $C$, $D$ and $E$ three coalgebras, a linear map $f:E\to [C,D]$ corresponds to a coalgebra map $g:E\otimes C\to D$ if and only if, for any $e\in E$ and $c\in C$,
$$
\Delta \left( f(e)(c)\right) = f(e^{(1)})(c^{(1)})\otimes f(e^{(2)})(c^{(2)}) \ (-1)^{|c^{(1)}||e^{(2)}|} \et \epsilon(f(e)(c))=\epsilon(e)\epsilon(c).
$$
Then, if $g(e,c):=f(e)(c)$, $g$ is a map of coalgebras.

\begin{defi}\label{remreductioncoalg}
We shall call a map of (dg-)vector spaces $E\to [C,D]$ a {\em comorphism of (dg-)coalgebras} if it satisfies the above equations, \ie if the corresponding map $E\otimes C\to D$ is a (dg-)coalgebra map. 
The coalgebra $\HOM(C,D)$ is equipped with a canonical comorphism $\Psi :\HOM(C,D)\to [C,D]$ associated to the ${\bf ev}:\HOM(C,D)\otimes C\to D$.
This comorphism is couniversal in the following sense: for any coalgebra $E$ and any comorphism $f:E\to [C,D]$, there exists a unique map of coalgebras $g:E\to \HOM(C,D)$ such that $\Psi g=f$.
$$
\xymatrix{
 && \HOM(C, D)\ar[d]^{\Psi}\\
E \ar@{-->}[rru]^{g}  \ar[rr]^-{f}  && [C,D]
}$$
We shall say that the map $\Psi:\HOM(C,D)\to [C,D]$ is the {\em couniversal comorphism}. 
We shall prove in proposition \ref{strongforgetfulVectcoalg} that $\Psi$ is part of an enriched functor.
\end{defi}

\begin{rem}\label{remconstructionHOM}
The previous universal property in terms of comorphisms can be used to have a more constructive proof of the existence of the functor $\HOM$.
For $C$, $D$ and $E$ three coalgebras, a map $f:E\to [C,D]$ is a comorphism iff the following diagrams commute
$$
\vcenter{\xymatrix{
E\ar[rr]^-{\Delta_E} \ar[ddd]_f&& E\otimes E\ar[d]^{f\otimes f}\\
	&&[C,D]\otimes [C,D]\ar[d]^{can}\\
	&&[C\otimes C,D\otimes D]\ar[d]^{[\Delta_C,D\otimes D]}\\
[C,D]\ar[rr]^-{[C,\Delta_D]} &&[C,D\otimes D]\\
}}
\et
\vcenter{\xymatrix{
E \ar[r]^-f\ar[d]_{\epsilon_E} & [C,D]\ar[d]^{[C,\epsilon_D]}\\
\FF\ar[r]^-{\epsilon_C} & [C,\FF]
}}.$$

If $E=T^\vee([C,D])$, with $p:T^\vee([C,D])$ in stead of $f$, the diagrams do not commute, but each is defining a pair of parallel maps
in the category $\dgVect$,
$$\xymatrix{
T^\vee([C,D])\ar@<.6ex>[r]\ar@<-.6ex>[r]& [C,D\otimes D]
}
\et
\xymatrix{
T^\vee([C,D])\ar@<.6ex>[r]\ar@<-.6ex>[r]& D
}
$$
From these pairs, we obtain by coextension two other pairs of parallel maps in the category $\dgCoalg$,
$$\xymatrix{
v,v':T^\vee([C,D])\ar@<.6ex>[r] \ar@<-.6ex>[r]& T^\vee([C,D\otimes D])
}
\et
\xymatrix{
w,w':T^\vee([C,D])\ar@<.6ex>[r]\ar@<-.6ex>[r]& T^\vee(D)
}
$$
The coalgebra $\HOM(C,D)$ is the common equalizer of the pairs $(v,v')$ and $(w,w')$ in the category $\dgCoalg$ (see appendix \ref{commonequalizer}). 

This gives another proof of the existence of $\HOM(C,D)$, but, as limits in $\dgCoalg$ are difficult to construct, this construction stays somehow formal.
In particular, with this construction, it is not clear that $\HOM(C,D)$ is a subcoalgebra of $T^\vee([C,D])$.
We shall give in corollary \ref{HOMsubcoalg} another construction of $\HOM(C,D)$ using the comonadicity theorem \ref{comonadic1}.
\end{rem}

\medskip

\begin{prop} \label{closedmonoidcoalg4}  
For every coalgebra $C$, we have
\begin{enumerate}
\item $\HOM(\FF,C)=C$, in particular $\HOM(\FF,\FF)=\FF$,
\item $\HOM(C,\FF)=\FF$ and 
\item $\HOM(0,C)=\FF$.
\end{enumerate}
\end{prop}

\begin{proof}
The coalgebra $\FF$ is the unit object the monoidal category $\dgCoalg$. Hence we have $\HOM(\FF,C)=C$.
Then, the coalgebra $0$ and $\FF$ are respectively the initial and terminal objects of the category $\dgCoalg$,
since $\HOM$ is right adjoint, they are sent to the terminal object by $\HOM(-,C)$ and $\HOM(C,-)$ respectively.
\end{proof}

\medskip

\begin{cor} \label{corcounitHOM}  
Let $C$ and $D$ be two coalgebras, then the counit of $\HOM(C,D)$ is given by the map
$$\xymatrix{
\HOM(C,D)\ar[rr]^-{\HOM(C,\epsilon_D)}&& \HOM(C,\FF)=\FF
}$$
\end{cor}

\begin{proof}
$\FF$ is terminal in $\dgCoalg$ and that the unique map to $\FF$ is given by the counit.
Then the result follows from $\HOM(C,\FF)=\FF$.
\end{proof}

\bigskip

The composition law $\mathbf{c}:\HOM(D,E)\otimes \HOM(C,D) \to  \HOM(C,E)$ is defined to be the unique map of coalgebras $\mathbf{c}$
such that the following square commutes
$$\xymatrix{
\HOM(D,E)\otimes \HOM(C,D)\otimes C \ar@{-->}[rr]^-{\mathbf{c}\otimes C} \ar[d]_{\HOM(D,E) \otimes {\bf ev}} && \HOM(C,E)\otimes C \ar[d]^{{\bf ev}}\\
\HOM(D,E)\otimes D\ar[rr]^-{{\bf ev}} && E
}$$
Thus, $\mathbf{c}={\Lambda}^3({\bf ev}^2 )$, where ${\bf ev}^2={\bf ev}(\HOM(D,E) \otimes {\bf ev})$
is the strong evaluation iterated two times.
Similarly, the unit $u_C$ of $C$ is the unique map of coalgebras $u_C:\FF \to  \END(C)$ 
such that the following triangle commutes 
$$
\xymatrix{
&& \END(C) \otimes C  \ar[d]^-{{\bf ev}} \\
\FF \otimes C \ar@{-->}[urr]^-{u_C\otimes C}  \ar@{=}[rr]    && C
}$$
The element $u_C(1)=e_C$ is an atom of the coalgebra $ \END(C) $ by proposition \ref{atomsandmap}.
We shall identify $u_C$ with $e_C$ and write $e_C:\FF\to C$.
$\END(C)$ has the structure of a monoid in $\dgCoalg$ and is therefore a bialgebra.

\medskip
\begin{lemma}\label{lemmecommutcoalg}
The composition $\mathbf{c}:\HOM(D,E)\otimes \HOM(C,D) \to  \HOM(C,E)$ is the unique map of coalgebras $\mathbf{c}$ such that the following square commutes
$$\xymatrix{
\HOM(D,E)\otimes \HOM(C,D) \ar@{-->}[rr]^-{\mathbf{c}} \ar[d]_{\Psi\otimes \Psi} && \HOM(C,E) \ar[d]^{\Psi}\\
[D,E]\otimes [C,D]\ar[rr]^-{c} && [C,E]
}$$
where $[D,E]\otimes [C,D]\to [C,E]$ is the strong composition in $\dgVect$.
\end{lemma}
\begin{proof}
Recall that $ev(\Psi\otimes C)={\bf ev}$ by definition of $\Psi$. We deduce the commutative diagram
$$\xymatrix{
\HOM(D,E)\otimes \HOM(C,D)\otimes C \ar[rrr]^-{{\bf c}\otimes C}\ar[d]_{\HOM(D,E)\otimes\Psi\otimes C}\ar[rrd]^{\qquad\HOM(D,E)\otimes {\bf ev}}&&& \HOM(C,E)\otimes C \ar[dd]^{\bf ev}\\
\HOM(D,E)\otimes [C,D]\otimes C \ar[rr]_-{\HOM(D,E)\otimes ev}\ar[d]_{\Psi \otimes [C,D]\otimes C}&&\HOM(D,E)\otimes D\ar[d]^{\Psi\otimes D}\ar[rd]^{\bf ev}& \\
[D,E]\otimes [C,D]\otimes C \ar[rr]^-{[D,E]\otimes ev}&&[D,E]\otimes D\ar[r]^-{ev}& E.
}$$
The result is proven by transposing $C$ on the border of the diagram.
\end{proof}

\bigskip
Recall that if $\bf C$ is a category enriched over a monoidal category $({\bf V},\otimes,\un)$, 
the {\em underlying set} of morphisms between $X$ and $Y$ in $\bf C$ is defined as $\Hom_{\bf V}(\un,{\bf C}(X,Y))$ and
the {\em underlying category} of $\bf C$ is defined as the category with the same objects as $\bf C$ but with sets of morphisms $\Hom_{\bf V}(\un,{\bf C}(X,Y))$. 
Let $\dgCoalg^\$ $ be the category $\dgCoalg$ viewed as enriched over itself.
Recall from section \ref{atoms} that an atom of a coalgebra $E$ is the same thing as a coalgebra map $\FF\to E$.
The set of underlying elements of the coalgebra $\HOM(C,D)$ is thus the set of atoms of $\HOM(C,D)$.

\begin{lemma}\label{underlyingCoalglemma}
Let $C$ and $D$ be two coalgebras, the set of atoms of $\HOM(C,D)$ is in bijection with the set of coalgebras maps from $C$ to $D$.
\end{lemma}
\begin{proof}
By the universal property of $\HOM(C,D)$, coalgebra maps $C\to D$ are in bijection with coalgebra maps $\FF\to \HOM(C,D)$.
\end{proof}

\begin{prop}\label{underlyingCoalg}
The underlying category of $\dgCoalg^\$$ is the ordinary category $\dgCoalg$.
\end{prop}
\begin{proof}
This is a formal property of monoidal closed categories.
\end{proof}

\bigskip

Let $C$ be a (dg-)coalgebra and $X$ be a (dg-)vector space.
If $q: T^\vee([C,X])\to [C,X]$ is the cofree map, then the composite of the maps
$$\xymatrix{
T^\vee([C,X])\otimes C \ar[rr]^-{q\otimes C} && [C,X]\otimes C \ar[r]^-{ev} & X
}$$
can be coextended as a map of coalgebras $e: T^\vee([C,X])\otimes C \to T^\vee(X)$.

\begin{prop}\label{exemplehomcoalg}
Let $C$ be a (dg-)coalgebra and $X$
be a (dg-)vector space. For any coalgebra $E$ and any map of coalgebras
 $f:E\otimes C\to T^\vee(X)$, there exists a unique map of coalgebras
 $k:E\to  T^\vee([C,X])$ such that $e(k\otimes C)=f$. Thus, 
 $e$ is a strong evaluation ${\bf ev}:\HOM(C,T^\vee(X)) \otimes C\to  T^\vee([C,X])$
 and we have
$$
\HOM(C,T^\vee(X)) = T^\vee([C,X]).
$$
\end{prop}

\begin{proof} For any coalgebra $E$, we have a chain of natural bijections between
\begin{center}
\begin{tabular}{lc}
\rule[-2ex]{0pt}{4ex} coalgebra maps & $f:E\otimes C \to T^\vee(V)$, \\
\rule[-2ex]{0pt}{4ex} linear maps  & $E\otimes C \to V$, \\
\rule[-2ex]{0pt}{4ex} linear maps  & $E \to [C,V]$, \\
\rule[-2ex]{0pt}{4ex} and coalgebra maps  & $k:E \to T^\vee([C,V])$.
\end{tabular}
\end{center}
Hence the coalgebra $T^\vee([C,V])$ is representing the functor $E\mapsto \dgCoalg(E\otimes C, T^\vee(V))$.
Moerover, if $E=T^\vee([C,V])$ and $f=e :T^\vee([C,X])\otimes C \to T^\vee(X)$
then $k$ is the identity of  $T^\vee([C,V])$.
This means that $e$ is playing the role of the evaluation map ${\bf ev}:\HOM(C,T^\vee(X)) \otimes C \to T^\vee(X)$.
\end{proof}

\medskip

\begin{cor}
For any dg-vector spaces $X$ and $Y$, we have
$$
\HOM(T^\vee(X),T^\vee(Y))=T^\vee([T^\vee(X),Y]).
$$
Hence the hom object between cofree coalgebras is cofree.  
\end{cor}

\bigskip

For a general coalgebra $D$, it is possible to give a copresentation of $\HOM(C,D)$ in terms of a copresentation of $D$.
The comonadicity of the adjunction $U:\dgCoalg\rightleftarrows \dgVect:T^\vee$ implies that it is possible to present $D$ as a reflexive equalizer in $\dgCoalg$ of cofree coalgebras
$$\xymatrix{
D\ar[r]^-{\Delta} &T^\vee(D) \ar@<.6ex>[r]\ar@<-.6ex>[r]& T^\vee(T^\vee(D)).
}$$
Then, because $\HOM(C,-)$ commute with limits, the object $\HOM(C,D)$ is the reflexive equalizer in $\dgCoalg$
$$\xymatrix{
\HOM(C,D)\ar[r]^-{\Delta'} &\HOM(C,T^\vee(D)) \ar@<.6ex>[r]\ar@<-.6ex>[r]& \HOM(C,T^\vee(T^\vee(D)))
}$$
Using proposition \ref{exemplehomcoalg}, this is the same equalizer as
$$\xymatrix{
\HOM(C,D)\ar[r]^-{\Delta'} &T^\vee([C,D]) \ar@<.6ex>[r]\ar@<-.6ex>[r]& T^\vee([C,T^\vee(D)]).
}$$

\begin{cor}\label{HOMsubcoalg}
$\HOM(C,D)$ is naturally a subcoalgebra of $T^\vee([C,D])$.
Moreover, the map $\HOM(C,D)\to [C,D]$, obtained by composition with the natural projection $q:T^\vee([C,D])\to [C,D]$, is the couniversal comorphism map
$$\xymatrix{
\HOM(C,D)\ar[rrd]_{\Psi}\ar[rr]^-{\Delta'}&& T^\vee([C,D])\ar[d]^q\\
&&[C,D].
}$$
\end{cor}
\begin{proof}
The comonadicity theorem that we proved says that coreflexive equalizers are reflected by the forgetful functor $U:\dgCoalg\to \dgVect$.
Applied to the equalizer above we deduced that the map $\HOM(C,D)\to T^\vee([C,D])$ is injective.

To prove the second statement let us consider the commutative diagram
$$\xymatrix{
\HOM(C,D)\ar[r]^-{\Delta'}\ar[d]^{\Psi}&\HOM(C,T^\vee(D))\ar[d]^{\Psi}\\
[C,D]\ar[r]\ar@{=}[rd]&[C,T^\vee(D)]\ar[d]^f\\
&[C,D]
}$$
where the horizontal map are induced by the coalgebra map $\Delta:D\to T^\vee(D)$ and the map $f$ is induced by the natural projection $T^\vee(D)\to D$.
Then, to prove that $\Psi = q\Delta'$, it is enough to prove that $f\Psi = q$.
But, this is true by definition of the map $e: T^\vee([C,X])\otimes C \to T^\vee(X)$ and the fact that $\Psi=\lambda^2(e)$.
\end{proof}

Recall the notion of cogenerating and separating maps from definition \ref{coseparation}.

\begin{prop} \label{cogeneratingcomorphism}
If $C$ and $D$ are coalgebras, then 
the couniversal comorphism $\Psi:\HOM(C,D)\to [C,D]$ is cogenerating, hence separating.
\end{prop}

\begin{proof}
Trivial from corollary \ref{HOMsubcoalg}.
\end{proof}

\bigskip
Recall that if $D\to C$ is a coalgebra map, then $D$ has a canonical structure of a $C$-bicomodule.
This produces a functor $\dgCoalg/C\to \Bicomod(C)$. 
The following construction is an analog for coalgebras of the tensor algebra over a bimodule.

\begin{prop}\label{cofreecoalgebraforbicomodules}
The functor $\dgCoalg/C\to \Bicomod(C)$ has a right adjoint $T^\vee_C$.
\end{prop}
\begin{proof}
Dual to proposition \ref{tensoralgebraforbimodules}.
\end{proof}

The following result, which gives another example of $\HOM$, is a coalgebra analog of the proposition \ref{derham} for algebras.
\begin{prop}\label{coderham}
Let $\d_n$ be the Hopf primitive coalgebra of example \ref{primitiveHopf}, then $\HOM(\d_n,C)$
is $T^\vee_C(S^{-n}\Omega^C)$ the cofree $C$-coalgebra cogenerated by the bicomodule of codifferential $\Omega^C$.
\end{prop}
\begin{proof}
Using the universal bicomodule of coderivations $\Omega^C$ from propositions \ref{univcoderiv},
and \ref{cofreecoalgebraforbicomodules} and proposition \ref{primicoderiv} (which is proven independently),
we have natural bijections between:
\begin{center}
\begin{tabular}{lc}
\rule[-2ex]{0pt}{4ex} coalgebra maps & $E\otimes \HOM(\d_n,C)$, \\
\rule[-2ex]{0pt}{4ex} coalgebra maps & $d_n\otimes \HOM(E,C)$, \\
\rule[-2ex]{0pt}{4ex} pairs (coalgebra map, coderivation of degree $n$) & $(E\to C,E\rhup C)$, \\
\rule[-2ex]{0pt}{4ex} pairs (coalgebra map, $C$-bicomodule map) & $(E\to C,E\to S^{-n} \Omega^C)$, \\
\rule[-2ex]{0pt}{4ex} and coalgebra maps & $E\to T^\vee_C(S^{-n}\Omega^C)$.
\end{tabular}
\end{center}
Thus $\HOM(\d_n,C) =T^\vee_C(S^{-n}\Omega^C)$ by Yoneda's lemma.
\end{proof}

\bigskip
We finish on a lemma useful to study the strength of the functor $\HOM$ in proposition \ref{strengthHOM}.

\begin{lemma}\label{lemmacomposeparation}
If $\psi:D\to X$ is a separating map, then, for any coalgebra $C$, the composition
$$\xymatrix{
h:\HOM(C,D) \ar[r]^-\Psi & [C,D]\ar[rr]^-{\hom(C,\psi)}&& [C,X]
}$$
is separating.
\end{lemma}
\begin{proof}
Let $f:D\to T^\vee(X)$ be the monomorphism of coalgebra associated to $\psi$.
$\HOM(C,-)$ is right adjoint, hence for any coalgebra $C$, the map $g:\HOM(C,D)\to  \HOM(C,T^\vee(X))=T^\vee([C,X])$ is still a monomorphism.
Let us consider the diagram
$$\xymatrix{
\HOM(C,D)\ar[rr]^-{\HOM(C,f)}\ar[d]_\Psi&& \HOM(C,T^\vee(X))\ar@{=}[r]\ar[d]_\Psi& T^\vee([C,X])\ar[d]^p\\
[C,D]\ar[rr]^-{\hom(C,f)} && [C,T^\vee(X)] \ar[r]& [C,X].
}$$
The left square commutes by naturality of $\Psi$ and the right square commutes by the characterization of the strong evaluation ${\bf ev}=e$ in proposition \ref{exemplehomcoalg}.
The top composition $\HOM(C,D)\to [C,X]$ is separating because $g$ is a monomorphism.
The bottom composite is the map $h$, it is separating by commutation of the diagram.
\end{proof}

\subsection{Monoidal strength and lax structures}\label{monoidalstrengthcoalg}

We finish this chapter by a section to mention that both functors $\otimes$ and $\HOM$ are strong and that $\otimes$ is in fact a strong symmetric monoidal structure. We omit the proofs as these are formal consequences of $(\dgCoalg,\otimes, \HOM)$ being a symmetric monoidal closed category. Some details on definitions and proofs are given in Appendix \ref{enrichedcategoryapp}.

\bigskip

$\dgCoalg \times \dgCoalg$ is a symmetric monoidal category when equipped with the termwise tensor product 
$$
(C_1,C_2)\otimes (D_1,D_2) = (C_1\otimes D_1,C_2\otimes D_2).
$$
This monoidal structure $\dgCoalg \times \dgCoalg$ is closed with internal hom defined by
$$
\HOM((C_1,C_2), (D_1,D_2)) = (\HOM(C_1, D_1),\HOM(C_2, D_2)).
$$

Let us consider the isomorphisms $\sigma_{23}:(C_1\otimes C_2)\otimes (D_1\otimes D_2)\simeq (C_1\otimes D_1)\otimes (C_2\otimes D_2)$ and $\sigma_0:\FF\simeq \FF\otimes \FF$.
The pair $(\sigma_{23},\sigma_0)$ is a monoidal structure on the functor $\otimes:\dgCoalg\times \dgCoalg\to \dgCoalg$.
We can then transfer the natural enrichment of $\dgCoalg\times \dgCoalg$ over itself into an enrichment over $\dgCoalg$.
We shall call $(\dgCoalg\times \dgCoalg)^{\$\otimes \$}$ the corresponding category enriched over $\dgCoalg$.
Its hom coalgebras are
$$
\HOM((C_1,C_2),(D_1,D_2)) := \HOM(C_1,D_1)\otimes \HOM(C_2,D_2)
$$
for any four coalgebras $C_1$, $C_2$, $D_1$ and $D_2$.
Similary we shall note $(\dgCoalg^{op}\times \dgCoalg)^{\$\otimes \$}$ the enrichment of $\dgCoalg^{op}\times \dgCoalg$ over $\dgCoalg$. Its hom coalgebras are
$$
\HOM((C_1,C_2),(D_1,D_2)) := \HOM(D_1,C_1)\otimes \HOM(C_2,D_2)
$$
for any four coalgebras $C_1$, $C_2$, $D_1$ and $D_2$.

The functor $\otimes:\dgCoalg\times \dgCoalg \to \dgCoalg$ commutes with colimits and finite limits in both variables but commute only globally to sifted colimits and finite sifted limits (same proof as in proposition \ref{siftedtensor1}). In particular it does not have any adjoint and the transfer along $\otimes$ is not compatible with the tensor and cotensor operations (see appendix \ref{basechange}).
Hence, the categories $(\dgCoalg\times \dgCoalg)^{\$\otimes \$}$ and $(\dgCoalg^{op}\times \dgCoalg)^{\$\otimes \$}$ are neither tensored nor cotensored over $\dgCoalg$.

\medskip

The strength of $\otimes:\dgCoalg\times \dgCoalg\to \dgCoalg$ is the coalgebra map
$$\xymatrix{
\Theta_\otimes:=\Lambda^{24}({\bf ev \otimes ev}):
\HOM(C_1,D_1)\otimes \HOM(C_2,D_2) \ar[r]& \HOM(C_1\otimes D_1, C_2\otimes D_2)
}$$
where ${\bf ev \otimes ev}=\HOM(C_1,D_1)\otimes C_1 \otimes \HOM(C_2,D_2) \otimes C_2\to D_1\otimes D_2$.
Using comorphisms, $\Theta_\otimes$ can also be characterized as the unique coalgebra map such making the following square commute
\begin{align}
\tag{Strength $\otimes$ (coalg)}\label{strengthotimescoalg}
\vcenter{\xymatrix{
\HOM(C_1,D_1) \otimes \HOM(C_2,D_2) \ar[rr]^-{\Theta_\otimes}\ar[d]_{\Psi\otimes \Psi}&& \HOM(C_1\otimes C_2, D_1\otimes D_2)\ar[d]^\Psi\\
[C_1,D_1] \otimes [C_2,D_2] \ar[rr]^\theta&& [C_1\otimes C_2, D_1\otimes D_2]
}}
\end{align}
where $\theta$ is the strength of $\otimes$ in $\dgVect$.

\medskip

The strength of $\HOM:\dgCoalg^{op}\times \dgCoalg\to \dgCoalg$ is the map
$$\xymatrix{
\Theta_\HOM:=\Lambda^2({\bf c^2}):\HOM(D_1,C_1)\otimes \HOM(C_2,D_2) \ar[r]& \HOM(\HOM(C_1, C_2),\HOM(D_1, D_2))
}$$
where ${\bf c^2}= \HOM(D_1,C_1) \otimes \HOM(C_1, C_2)\otimes \HOM(C_2,D_2)\to \HOM(D_1, D_2)$.
The functor $\HOM$ also inherits a lax monoidal structure $(\alpha,\alpha_0)$ given by $\alpha=\Theta_\otimes$ and $\alpha_0=\HOM(\FF,\FF)\simeq \FF$.

\begin{prop}\label{strengthHOM}
The map $\Theta_\HOM$ is the unique coalgebra map such that the following diagram commutes.
\begin{align}
\tag{Strength $\HOM$}\label{strengthHOMdiag}
\vcenter{\xymatrix{
\HOM(D_1,C_1)\otimes \HOM(C_2,D_2) \ar[rr]^-{\Theta_\HOM}\ar[d]_{\Psi\otimes \Psi}&& \HOM(\HOM(C_1, C_2),\HOM(D_1, D_2))\ar[d]^{\Psi}\\
[D_1,C_1]\otimes [C_2,D_2] \ar[d]_-{\theta} && [\HOM(C_1, C_2),\HOM(D_1, D_2)]\ar[d]^{[\HOM(C_1, C_2),\Psi]}\\
[[C_1, C_2],[D_1, D_2]] \ar[rr]^-{[\Psi,[D_1, D_2]]}&& [\HOM(C_1, C_2),[D_1, D_2]]
}}
\end{align}
where $\theta$ is the strength of $[-,-]$ in $\dgVect$.
\end{prop}
\begin{proof}
The proof of the commutation of this diagram is a careful unravelling of the definition of $\Theta_\HOM=\Lambda^2(\bf c^2)$ left to the reader.
Then, to prove the assertion, it is enough to prove that the right side map $[\HOM(C_1, C_2),\Psi]\circ \Psi:\HOM(\HOM(C_1, C_2),\HOM(D_1, D_2))\to $ is separating $[\HOM(C_1, C_2),[D_1, D_2]]$. But this is lemma \ref{lemmacomposeparation}.
\end{proof}

\medskip

\begin{prop}\label{strength for otimes and hom of coalg}
\begin{enumerate}
\item The maps $\Theta_\otimes$ enhance $\otimes$ into a strong functor
$$\xymatrix{
\otimes :(\dgCoalg\times \dgCoalg)^{\$\otimes \$} \ar[r]& \dgCoalg^\$.
}$$
Moreover, the associativity and unital morphisms of $\otimes$ are strong natural transformation.
Also, the monoidal structure $(\sigma_{23},\sigma_0)$ of $\otimes$ (with respect to itself) is a strong monoidal structure.

\item The maps $\Theta_\HOM$ enhance $\HOM$ into a strong functor
$$\xymatrix{
\HOM :(\dgCoalg^{op}\times \dgCoalg)^{\$\otimes \$} \ar[r]& \dgCoalg^\$.
}$$
$\HOM$ is a strong adjoint to $\otimes$, there exists isomorphisms
$$
\HOM(C\otimes D,E )\simeq \HOM(C,\HOM(D,E))\simeq \HOM(D,\HOM(C,E)).
$$
Moreover the lax monoidal structure $(\alpha,\alpha_0)$ of $\HOM$ is a strong lax monoidal structure.
\end{enumerate}
\end{prop}

\subsection{Meta-morphisms}\label{meta-morphismcoalg}

\subsubsection{Reduction functor and meta-morphisms}\label{reducetmeta-morphismcoalg}

Recall that $\dgCoalg^\$$ is our notation for the enrichment of $\dgCoalg$ over itself. 
We can transfer this enrichment along the monoidal functor $U:\dgCoalg\to \dgVect$ into an enrichment over $\dgVect$.
We denote by $\dgCoalg^{U\$}$ the corresponding category. 
Recall also that $\dgVect^\$$ is our notation for $\dgVect$ viewed as enriched over itself.

\begin{prop}\label{strongforgetfulVectcoalg}
The reduction maps $\Psi:\HOM(C,D)\to [C,D]$ are the strengths of an enrichment over $\dgVect$ of the forgetful functor
$$\xymatrix{
U:\dgCoalg^{U\$} \ar[r]&\dgVect^\$
}$$
\end{prop}
\begin{proof}
The compatibility with compositions is the content of lemma \ref{lemmecommutcoalg}.
\end{proof}

\begin{rem}
It is also possible to prove that this functor is strong lax monoidal, the lax structure is given by by \eqref{strengthotimesalg}. We shall not use this, but this gives a nice meaning to some diagrams.
\end{rem}

\begin{rem}\label{remstrongnotadjointcoalg}
The category $\dgCoalg^{U\$}$ is nor tensored nor cotensor over $\dgVect$,
such a structure would require to have a left adjoint for the functor $U:\dgCoalg\to \dgVect$.
Also the strong functor $U$ do not have any strong adjoint anymore.
If $T^\vee$ were a strong right adjoint, the adjunction strength would be given by the cogenerating map 
$U\HOM(C,T^\vee X)\simeq T^\vee[UC,X]\to [UC,X]$ but this map is never an isomorphism.
We will see in theorem \ref{strongcomonadicitycoalg} that the whole adjunction $U\dashv T^\vee$ can be enriched provided we enriched the two categories over $\dgCoalg$ instead of $\dgVect$.
\end{rem}

\medskip

\begin{defi}
We define a {\it meta-morphism of coalgebras} $f:C\leadsto D$ to be an element $f\in U\HOM(C,D)$. The composite $gf$ of two meta-morphisms of coalgebras  $f:C\leadsto D$ and $g:D\leadsto E$ is defined by putting $gf=\mathbf{c}(g\otimes f)$, where $\mathbf{c}$ is the strong composition law
$$\xymatrix{
\mathbf{c}:\HOM(D,E)\otimes \HOM(C,D)\ar[r]&  \HOM(C,E)
}$$
of section \ref{closedmonoidalstructurecoalg}. We have $|gf|=|g|+|f|$.

By opposition, we define a {\em pro-morphism of coalgebras} $f:C\rhup D$ to be an element of $[C,D]$. The composition of pro-morphisms is defined through the strong composition of $\dgVect$. 
\end{defi}

\begin{rem}\label{remfactocoalg}
The strong functor $U:\dgCoalg^{U\$} \to \dgVect^\$$ factors as 
$$\xymatrix{
U:\dgCoalg^{U\$} \ar[r]& (U\dgCoalg)^{\$}\ar[r]& \dgVect^\$
}$$
where the objects of $(U\dgCoalg)^{\$}$ are the coalgebras but the hom between two coalgebras $C$ and $D$ is simply $[C,D]$.
The functor $(U\dgCoalg)^{\$}\to \dgVect^\$$ is thus strongly fully faithful.
If $\bC$ is a category enriched over $\dgVect$ let us call the elements of $\bC(X,Y)$ the {\em graded morphisms} from $X$ to $Y$.
Then the meta-morphisms are the graded morphisms of $\dgCoalg^{U\$}$ and the pro-morphisms are the graded morphisms of $(U\dgCoalg)^{\$}$.
\end{rem}

\medskip

If $f:C\leadsto D$ is a meta-morphism, the reduction map $\Psi:\HOM(C,D)\to [C,D]$ defines a pro-morphism $\Psi(f)$. We shall simplify notations and simply put $f^\flat :=\Psi(f)$.
If $g:D\leadsto E$ and $f:C\leadsto D$ are two meta-morphisms, lemma \ref{lemmecommutcoalg} says that
$$
(gf)^\flat = g^\flat f^\flat.
$$

\medskip

Recall the strong evaluation ${\bf ev}=\HOM(C,D)\otimes C\to D$ from section \ref{closedmonoidalstructurecoalg}.
We define the {\em evaluation} of a meta-morphism $f:C\leadsto D$ on an element $x\in C$ to be 
$$
f(x):={\bf ev}(f\otimes x)\in D.
$$
If $f:C\leadsto D$ and $g:D\leadsto E$ are meta-morphisms of coalgebras, we have $(gf)(x)=g(f(x))$ for every $x\in C$ by definition of the strong composition map from $\bf ev$.

From the definition of the reduction map $\flat=\Psi:\HOM(C,D)\to [C,D]$ as $\Psi=\Lambda^2{\bf ev}$, we have a commutative diagram
$$\xymatrix{
\HOM(C,D)\otimes C   \ar[d]_{\flat \otimes C}  \ar@/^1pc/[rrd]^-{\mathbf{ev}} &&  \\
[C,D]\otimes C \ar[rr]^-{ev} && D.
}$$
In terms of elements, this gives $f(x) = f^\flat(x)$; the strong evaluation of a meta-morphism coincides with the evaluation of the corresponding pro-morphism.

\medskip

A meta-morphism $f\in \HOM(C,D)$ is said to be {\em atomic} if it is an atom in $\HOM(C,D)$.
We have seen in lemma \ref{underlyingCoalglemma} that the atomic meta-morphisms in $\HOM(C,D)$ are the same thing as coalgebra maps $C\to D$. If $f:C\to D$ is a coalgebra map we shall denote by $f^\sharp$ the corresponding atom of $\HOM(C,D)$.
This provides a map (of graded sets) $\sharp:\dgCoalg(C,D)\to \HOM(C,D)$.

\begin{lemma}\label{musicallemmacoalg}
\begin{enumerate}
\item Let $f:C\to D$ be a coalgebra map, then we have 
$$
(f^\sharp)^\flat = f.
$$
In other terms, for any two coalgebras $C$ and $D$, we have a commutative diagram of graded sets
$$\xymatrix{
At(\HOM(C,D))\ar[d]^\flat \ar[r] & \HOM(C,D)\ar[d]^{\Psi = \flat}\\
\dgCoalg(C,D) \ar@/^1pc/[u]^-\sharp \ar[r] &[C,D].
}$$
where $At(\HOM(C,D))$ is the set of atoms of $\HOM(C,D)$ and
where the horizontal maps are the canonical inclusions.

\item The maps $\sharp$ and $\flat$ induce inverse bijections of sets
$$\xymatrix{
\dgCoalg(C,D) \ar@<.6ex>[rr]^-\sharp&& At(\HOM(C,D)).\ar@<.6ex>[ll]^-{\flat=\Psi}
}$$

\item For a coalgebra map $f:C\to D$ and $x\in C$, $f(x)$ can any way: $f(x)=(f^\sharp)^\flat(x)= f^\sharp(x)$.

\item If $f:C\to D$ and $g:D\to E$ are maps of coalgebras, we have
$$
g^\sharp f^\sharp = (gf)^\sharp.
$$
\end{enumerate}
\end{lemma}

\begin{proof}
\begin{enumerate}
\item Let $f:C\to D$ be a coalgebra map, then $f$ is the element corresponding to the map $\lambda(f):\FF\to \HOM(C,D)$.
The atom corresponding to $f$ is $f^\sharp = \Lambda(f)(1)$ where $\Lambda(f):\FF\to \HOM(C,D)$
is the unique coalgebra map lifting the comorphism $\lambda(f):\FF\to [C,D]$, in particular $\Psi\circ \Lambda(f) = \lambda(f)$.
Thus we can write
$$
(f^\sharp)^\flat =\Psi(f^\sharp) = \Psi(\Lambda(f)(1)) = \lambda(f)(1) = f.
$$
\item This a reformulation of 1. using lemma \ref{underlyingCoalglemma}.
\item Direct from 1.
\item For any coalgebra map $h:C\to D$, $h^\sharp$ is the unique element of $\HOM(C,D)$ such that $(f^\sharp)^\flat=f$.
Then the result is a consequence of $(g^\sharp f^\sharp)^\flat = (g^\sharp)^\flat( f^\sharp)^\flat = gf = (gf)^\sharp)^\flat$.
\end{enumerate}
\end{proof}

\begin{rem}
Recall the category $(U\dgCoalg)^{\$}$ from remark \ref{remfactocoalg}.
Lemma \ref{musicallemmacoalg} says that we have a commutative diagram of categories (enriched over graded sets)
$$\xymatrix{
 && \dgCoalg^{U\$}\ar[d]^{\flat}\\
\dgCoalg \ar[rru]^-\sharp \ar[rr]&&(U\dgCoalg)^{\$}.
}$$
\end{rem}

\subsubsection{Calculus of meta-morphisms}

\paragraph{Tensor product}

Recall from section \ref{monoidalstrengthcoalg} the strength of $\otimes :\dgCoalg\times \dgCoalg \to \dgCoalg$
$$\xymatrix{
\Theta_\otimes:\HOM(C_1,D_1)\otimes  \HOM(C_2, D_2) \ar[r]& \HOM(C_1\otimes C_2,  D_1 \otimes D_2).
}$$
If $f:C_1\leadsto D_1$ and $g:C_2\leadsto D_2$ are  meta-morphisms of coalgebras,
let us define the {\em tensor product of meta-morphisms} by  
$$\xymatrix{
f\otimes g:=\Theta_\otimes(f\otimes g):C_1 \otimes C_2 \ar@{~>}[r]& D_1 \otimes D_2.
}$$
In addition, let us put $C \otimes g: =1_C\otimes g$ and $f \otimes D:= f\otimes 1_D$ where $1_C$ and $1_D$ are the units of the bialgebras $\END(C)$ and $\END(D)$.

\medskip
The underlying functor of the strong functor $\otimes$ is the functor $\otimes$.
This implies, for $f:C_1\to D_1$ and $g:C_2\to D_2$ two coalgebra maps, the relation
$$
(f\otimes g)^\sharp = f^\sharp \otimes g^\sharp.
$$

\begin{prop}\label{tensormusicalcoalg}
For $f:C_1\leadsto D_1$ and $g:C_2\leadsto D_2$ two meta-morphisms of coalgebras, we have
$$
(f\otimes g)^\flat = f^\flat \otimes g^\flat.
$$
In particular, we can reconstruct $f\otimes g$ from $f^\sharp\otimes g^\sharp$ as $(f^\sharp\otimes g^\sharp)^\flat$.
\end{prop}
\begin{proof}
$(f\otimes g)^\flat=f^\flat\otimes g^\flat$ is a consequence of the commutation of \eqref{strengthotimescoalg}. 
The last assertion is a consequence of $(f^\sharp)^\flat=f$.
\end{proof}

\begin{prop}\label{metatensorcoalg1} 
If $u:C_1\leadsto D_1$,  $f:D_1\leadsto E_1$, $v:C_2\leadsto D_2$ and $g:D_2\leadsto E_2$
are meta-morphisms of coalgebras, then we have 
$$
(f \otimes g)(u \otimes v)=fu \otimes gv (-1)^{|u||g|}.
$$
In particular, we have
$$
f \otimes g = (f \otimes D_2 )(C_1\otimes g)=(D_1\otimes g) (f \otimes C_2 )  (-1)^{|f||g|}.
$$
\end{prop}

\begin{proof} The first identity follows from the functoriality of the strong functor $\otimes :\dgCoalg\times \dgCoalg \to  \dgCoalg$.
The second identity is a special case of the first.
\end{proof}

\begin{cor}\label{metatensorcoalg2} 
If $f:D_1\leadsto E_1$ and $g:D_2\leadsto E_2$ are meta-morphisms of coalgebras,
then we have $$(f\otimes g)(x\otimes y)=f(x)\otimes g(y) (-1)^{|g||x|}$$
for every $x\in D_1$ and $y\in D_2$. 
\end{cor}

\begin{proof} 
Recall from proposition \ref{closedmonoidcoalg4} that $\HOM(\FF,C)=C$, then apply proposition \ref{metatensorcoalg1} with $u=x:\FF\leadsto C_2$ and $v= y:\FF\leadsto D_2$.
\end{proof}

\bigskip

\paragraph{Internal hom}

The strength of the internal hom functor $\HOM:\dgCoalg^{op}\times \dgCoalg\to \dgCoalg$ is the map
$$\xymatrix{
\Theta_{\HOM}:\HOM(D_1,C_1) \otimes  \HOM(C_2,D_2) \ar[r]& \HOM\bigl(\HOM(C_1,C_2),  \HOM(D_1,D_2) \bigr)
}$$
defined in section \ref{monoidalstrengthcoalg}.
If $f:D_1\leadsto C_1$ and $g:C_2 \leadsto  D_2$ are meta-morphisms of coalgebras, let us define the {\em hom of meta-morphisms} by  
$$\xymatrix{
\HOM(f,g):=\Theta_{\HOM}(f\otimes g): \HOM(C_1,D_1) \ar@{~>}[r]& \HOM(C_2,D_2).
}$$
In addition, let us put $\HOM(C,g):=\HOM(1_C,g)$ and $\HOM(f,D):=\HOM(f,1_D)$ where $1_C$ and $1_D$ are the units of the bialgebras $\END(C)$ and $\END(D)$.

\bigskip

The underlying functor of the strong functor $\HOM$ is the functor $\HOM$.
This implies, for $f:D_1\to C_1$ and $g:C_2\to D_2$ two coalgebra maps, the relation
$$
\HOM(f,g)^\sharp = \HOM(f^\sharp,g^\sharp).
$$

\begin{prop}\label{redmorphismcoalg1} 
If $f:C_2\leadsto D_2$, $g:D_1\leadsto C_1$, $u:D_2\leadsto E_2$,  and $v:E_1\leadsto D_1$
are meta-morphisms of coalgebras, then we have 
$$
\HOM(gv,uf)=\HOM(v,u)\HOM(g,f)\ (-1)^{|g|(|v|+|u|)}.
$$
In particular, we have
$$
\HOM(g,f)=\HOM(g,C_2)\HOM(D_1,f)=\HOM(D_1,f)\HOM(g,C_1)\ (-1)^{|f||g|}.
$$
\end{prop}

\begin{proof}
The first identity follows from the functoriality of the strong functor $\HOM$.
The second identity is a special case of the first.
\end{proof}

\begin{prop}\label{redmorphismcoalg2} 
If $g:C_2\leadsto D_2$ and $f:D_1\leadsto C_1$ are  meta-morphisms of coalgebras,
then we have 
$$
\HOM(g,f)(h)=fhg\ (-1)^{|g|(|f|+|h|)}
$$
for every meta-morphism $h:C_1\leadsto D_1$.
In particular, the following square of graded morphisms commutes in $\dgVect$
$$\xymatrix{
\HOM(C_1,D_1) \ar[d]_{\Psi=\flat}\ar@^{>}[rrr]^-{\HOM(g,f)^\flat} &&&  \HOM(C_2,D_2)\ar[d]^{\Psi=\flat} \\
[C_1,D_1] \ar@^{>}[rrr]^-{\hom(g^\flat,f^\flat)}  &&& [C_2,D_2].
}$$
In other terms, we have $\left(\HOM(g,f)^\flat(h)\right)^\flat = \hom(g^\flat,f^\flat)(h^\flat)$ for any $h:C_1\leadsto D_1$.
\end{prop}

\begin{proof}
By definition, 
\begin{eqnarray*}
\HOM(g,f)(h)&=&\Theta_\HOM(g\otimes f)(h) \\
&=&{\Lambda}^2(\mathbf{c^2})\sigma(g\otimes f)(h)\\
&=&{\Lambda}^2(\mathbf{c^2})(f\otimes g)(h) (-1)^{|g||f|}\\
&=&\mathbf{c^2}(f\otimes h \otimes g) (-1)^{|g||f|+|f||h|}\\
&=& ghf(-1)^{|g|(|f|+|h|)}
\end{eqnarray*}
This proves the first assertion. The second is due to the fact that the formula $\left(\HOM(g,f)^\flat(h)\right)^\flat = \hom(g^\flat,f^\flat)(h^\flat)$ is a reformulation of the commutation of the diagram \eqref{strengthHOMdiag}.
\end{proof}

\subsubsection{Module-coalgebras}

\begin{defi}\label{Qcoalgebradefi}
A $Q$-module coalgebra $C$ is a comonoid in the category of $Q$-modules, \ie it is the data of a $Q$-module $C$ and maps
$$
\Delta_C:C\to C\otimes C \et \epsilon_C:C\to \FF
$$
that are $Q$-equivariant. We shall denote the category of $Q$-coalgebras by $Q\dgCoalg$.
\end{defi}

\begin{lemma}\label{modulecoalgebra}
With the previous notations, $\Delta_C$ and $\epsilon_C$ are $Q$-equivariant if and only if the map $a:Q\otimes C\to C$ is a map of coalgebras.
\end{lemma}

\begin{proof} 
The maps $\Delta_C$ and $\epsilon_C$ are equivariant if and only if the following square commutes
$$
\vcenter{\xymatrix{
Q\otimes C \ar[ddd]_{a}\ar[rr]_-{Q\otimes \Delta_C}&& Q\otimes C\otimes C\ar[d]^{\Delta_Q\otimes C\otimes C}\\
&& Q\otimes Q\otimes C\otimes C\ar@{=}[d]\\
&&Q\otimes C\otimes Q\otimes C\ar[d]^-{a\otimes a}\\
C\ar[rr]^-{\Delta_C}&&C\otimes C
}}
\et\vcenter{
\xymatrix{
Q\otimes C\ar[rr]\ar[d]_{Q\otimes \epsilon_C}&& C\ar[d]^{\epsilon_C}\\
Q\otimes \FF \ar[r]^-{\epsilon_Q\otimes \FF}& \FF\otimes \FF\ar@{=}[r]&\FF.
}}$$
But this is exactly means that $a$ is a coalgebra map.
\end{proof}

\begin{ex} Every coalgebra $C$ is a $Q$-module-coalgebra over the bialgebra $Q=\END(C)$.
The action of the bialgebra $\END(C)$ on $C$ is given by the evaluation map $\mathbf{ev}:\END(C) \otimes C \to C$.
\end{ex}

\begin{ex}  If $C$ and $D$ are algebras, then the coalgebra $\HOM(C,D)$ has the structure of a $Q$-module-coalgebra 
over the bialgebra $Q=\END(C)^o\otimes \END(D)$.
\end{ex}

\begin{ex}
A module-coalgebra over the bialgebra $\FF[\delta]$ of example \ref{shufflehopf} 
is a coalgebra $C$ equipped with a coderivation $\delta_C:C\to C$ of degree $|\delta|$.
\end{ex}

\medskip

\begin{defi}
Let $Q$ be a bialgebra and $C$ a coalgebra. 
A {\em meta-action} (or simply an {\em action}) of $Q$ on $C$ is a map of coalgebras $a:Q\otimes C\tto C$ such that the following squares are commutative
$$\vcenter{
\xymatrix{
Q\otimes Q\otimes C\ar[rr]^-{Q\otimes a} \ar[d]_{m_Q\otimes C} && Q\otimes C\ar[d]^a\\
Q\otimes C\ar[rr]^-{a} && C
}}
\et
\vcenter{
\xymatrix{
\FF\otimes C\ar@{=}[rrd]\ar[rr]^-{e_Q\otimes C}&&Q\otimes C\ar[d]^-{a}\\
&&C.
}}$$
\end{defi}

\begin{lemma}\label{lemmametaactioncoalg}
A meta-action is equivalent to the data of a comorphism $\pi=\lambda^2a:Q\to [C,C]$ which is an algebra map.
\end{lemma}
\begin{proof}
By forgetting the coalgebra structure a meta-action is an action in the category $\dgVect$, thus $\pi=\lambda^2a:Q\to [C,C]$ is an algebra map. The result then follows by definition of comorphisms.
\end{proof}

Recall from lemma \ref{lemmecommutcoalg} that, for any coalgebra $C$, the reduction mapping
$$
\Psi:\END(C)\to [C,C]
$$
is a map of algebras.

\begin{prop} \label{hombialg}
If $Q$ is a bialgebra, any comorphism $f:Q\to [C,C]$ which is an algebra map, lift to a unique bialgebra map $\phi:Q\to \END(C)$ such that $\Psi\phi=f$
$$\xymatrix{
&&\END(C)\ar[d]^\Psi\\
Q\ar[rr]_-f\ar[rru]^-\phi &&[C,C].
}$$
In particular composition with $\Psi$ provides a bijection between bialgebras maps $Q\to \END(C)$ and comorphisms $Q\to[C,C]$ that are algebra maps.
\end{prop}

\begin{proof}
Let $f:Q\to [C,C]$ be a comorphism and $\phi:Q\to \END(C)$ the unique coalgebra map such that $\Psi\phi=f$.
To prove that $\phi$ is an algebra map, let us consider the diagram
$$\xymatrix{
&\END(C)\otimes \END(C)\ar[rr]^-{\mathbf{c}}\ar[d]^{\Psi \otimes \Psi }&& \END(C)\ar[d]^{\Psi}\\
&[C,C]\otimes [C,C]\ar[rr]^-{c} && [C,C]\\
Q\otimes Q\ar[rr]^-(.6){m_Q}\ar[ru]_-{f\otimes f}\ar@/^1pc/[ruu]^{\phi\otimes \phi}&&Q\ar[ru]_-{f}\ar@/^1pc/[ruu]^\phi.
}$$
The bottom face is commutative because $f:Q\to [C,C]$ is a map of algebras,
the back face is commutative by lemma \ref{lemmecommutcoalg}
and each side triangle is commutative by construction of $\phi$, hence the top face is commutative.
This prove that $\phi$ is a map of algebras.
\end{proof}

\medskip

If $Q$ is a bialgebra, then the endo-functor $ Q\otimes (-)$ of the category $\dgCoalg$
has the structure of a monad. The multiplication of the monad $ Q\otimes (-)$
is given by the map $\mu \otimes  C:Q\otimes Q\otimes C\to Q\otimes C$
and the unit by the map $e\otimes C: C \to Q\otimes C$. 
The endo-functor $\HOM(Q,-)$ of the category 
$\dgCoalg$ is right adjoint to the endo-functor $Q\otimes (-)$ by theorem \ref{homcoalg}.
Hence the endo-functor $\HOM(Q,-)$ has the structure
of a comonad, since the  endo-functor $Q\otimes (-)$
has the structure of a monad.
The comultiplication of the comonad  $\HOM(Q,-)$ 
is is given by the map $\HOM(\mu,C):\HOM(Q,C)\to \HOM(Q\otimes Q, C)=\HOM(Q, \HOM(Q,C))$
and the counit by the map $\HOM(e,C):\HOM(Q,C)\to \HOM(\FF,C)=C$.

A map $a:Q\otimes C \to C$ is an action of the monad $Q\otimes (-)$
if and only if the map $\Lambda^1(a):C\to \HOM(Q,X)$ is a coaction
of the comonad $\HOM(Q,-)$.
The category of algebras over the monad $ Q\otimes (-)$ is equivalent to the
category of coalgebras over the comonad  $\HOM(Q,-)$.

\begin{prop}\label{caracgammacoalg}
If $C$ is a coalgebra, then the following data are equivalent:
\begin{enumerate}
\item a meta-action $Q\otimes C\to C$ of the monoid $Q$ on the object $C$;
\item an action $Q\otimes C\to C$ of the monad $Q\otimes (-)$;
\item a coaction $C\to \HOM(Q,C)$ of the comonad $\HOM(Q,-)$;
\item the structure of a $Q$-module coalgebra on $C$;
\item a map of algebras $\pi:Q\to [C,C]$ which is a comorphism;
\item a map of bialgebras $\pi:Q\to \END(C)$.
\end{enumerate}
\end{prop}

\begin{proof}
The equivalence between (1) and (2) is a general fact true in any closed category.
The equivalence between (2) and (3) is the remark above.
The equivalence between (1) and (4) is lemma \ref{modulecoalgebra}.
The equivalence between (5) and (6) is proposition \ref{hombialg}.
Finally, the equivalence between (4) and (5) is lemma \ref{lemmametaactioncoalg}.
\end{proof}

\begin{cor} \label{gammacoalgebraforget}
 The forgetful functor $Q\dgCoalg\to \dgCoalg$ has a left adjoint $C\to Q\otimes C$ 
 and a right adjoint $C\mapsto \HOM(Q,C)$.
In particular, limits and colimits exists in $Q\dgCoalg$ and can be computed in $\dgCoalg$.
\end{cor}

\begin{proof}
This follows from the general theory of monads and comonads.
\end{proof}

Let now $Q$ be a cocommutative Hopf algebra.
Then, the category $Q\dgCoalg$ is enriched, tensored and cotensored over the category $\dgCoalg$. 
The tensor product of a $Q$-module-coalgebra $C$ by a coalgebra $D$ is the coalgebra $D\otimes C$ and the cotensor is $\HOM(D,C)$. 
Moreover, the tensor product and internal hom of two $Q$-module coalgebras is again a $Q$-module algebra, we just have to copy the formulas from section \ref{moduleoverhopfalgebra}. Let $C$ and $D$ be two $Q$-module coalgebras, the action of $Q$ on the tensor product and internal hom are defined respectively by
$$\xymatrix{
Q \ar[r]^-{\Delta} & Q \otimes Q \ar[rr]^-{\pi_C\otimes \pi_D} && \END(C)\otimes \END(D) \ar[r]^-{\Theta_\otimes} &  \END(C\otimes D)
}$$
and by
$$\xymatrix{
Q \ar[r]^-{\Delta} & Q \otimes Q\ar[r]^-{S\otimes Q} & Q^o  \otimes Q \ar[rr]^-{\pi_C\otimes \pi_D} && \END(C)^o \otimes \END(D) \ar[r]^-{\Theta_\HOM} &  \END(\HOM(C,D))
}$$
where the $\pi$s are the meta-actions of $Q$ and $\Theta_\otimes$ and $\Theta_\HOM$ are the strength of $\otimes$ and $\HOM$ in $\dgCoalg$ (see section \ref{monoidalstrengthcoalg}).
In the case where the action of $Q$ on $D$ is the trivial action $\epsilon_Q\otimes D:Q\otimes D\to D$, these formulas give the tensor and cotensor of $Q\dgCoalg$ over $\dgCoalg$. Then, as in the end of section \ref{moduleoverhopfalgebra} the enrichment of $Q\dgCoalg$ over $\dgCoalg$ is defined as the equalizer in $\dgCoalg$
$$\xymatrix{
\HOM_Q(C,D) \ar[r]& \HOM(C,D)=\HOM(\FF,\HOM(C,D))\ar@<.6ex>[rrr]^-{\HOM(\epsilon_Q,\HOM(C,D))}\ar@<-.6ex>[rrr]_-{\Lambda^1a} &&& \HOM(Q,\HOM(C,D))
}$$
where $a:Q\otimes \HOM(C,D)\to \HOM(C,D)$ is the action. We leave the reader to check the details.

\medskip
The forgetful functor $U:Q\dgCoalg\to \dgCoalg$ is strong and the adjunctions $Q\otimes (-)\dashv U\dashv \HOM(Q,-)$ of proposition \ref{gammacoalgebraforget} are strong.
Moreover we have the stronger result, which is a generalisation of theorem \ref{homcoalg}.

\begin{thm}\label{Qhomcoalg}
The category $Q\dgCoalg$ is symmetric monoidal closed and the forgetful functor $U:Q\dgCoalg\to \dgCoalg$ is symmetric monoidal and preserves the internal hom.
\end{thm}
\begin{proof}
Let $C$ and $D$ be two $Q$-module coalgebras, the tensor product and internal hom are those in $\dgCoalg$ with the actions of $Q$ defined above. The closeness, is proven using the second version of the proof of proposition \ref{freeandcofreeQmodules2}.
The second statement is obvious by construction.
\end{proof}

It is interesting to note that the map
$$\xymatrix{
Q \ar[r]^-{\Delta} & Q \otimes Q \ar[rr]^-{\pi_C\otimes \pi_D} && \END(C)\otimes \END(D) \ar[r]^-{\Theta_\otimes} &  \END(C\otimes D)
}$$
is the unique bialgebra map lifting the algebra map 
$$\xymatrix{
Q \ar[r]^-{\Delta} & Q \otimes Q \ar[rr]^-{\pi_C\otimes \pi_D} && [C,C]\otimes [D,D] \ar[r]^-{\theta_\otimes} &  [C\otimes D,C\otimes D]
}$$
through $\Psi:\END(C\otimes D)\to [C\otimes D,C\otimes D]$.

\bigskip
It should be clear to the reader that the same result holds if we work with graded vector spaces, graded coalgebras etc. instead of their differential graded analogs. We can then apply theorem \ref{Qhomcoalg} to deduce the following important result.

\begin{thm}\label{dgtoghomcoalg}
The category $\gCoalg$ is symmetric monoidal closed.
The forgetful functor $U_d:\dgCoalg\to \gCoalg$ is symmetric monoidal,
it has left adjoint $\d \otimes-$ and right adjoint $\HOM(\d,-)$,
and preserves the internal hom.
\end{thm}
\begin{proof}
The proof of the first statement is analog of that of theorem \ref{homcoalg}.
The second statement can be deduce from 
corollary \ref{gammacoalgebraforget} and (the graded analog of) 
theorem \ref{Qhomcoalg} if, as in example \ref{fromgtodg}, we describe $\dgVect$ as $\Mod(Q)$ in $\gVect$ for the cocommutative Hopf algebra $Q=\d=\FF\delta_+$. But we need to prove that the tensor product and internal hom in $\dgCoalg$ constructed from those of $\gCoalg$ via theorem \ref{Qhomcoalg} coincide with those already constructed.
By adjunction, it is enough to prove that the two tensor products are the same, but this is easy to see that the differential in $C\otimes D$ is the same computed in $\dgCoalg$ or in $Q\gCoalg$ for $Q=\FF\delta_+$.
\end{proof}

This result says that to compute $\HOM(C,D)$ in the dg-context, we can first compute it in the graded context and there will be a unique differential on the graded coalgebra $\HOM(C,D)$, induced by that of $C$ and $D$, enhancing it into the dg-internal hom.
We shall detail how to compute this differential in section \ref{derivativesweedlercoalg}.

\subsubsection{Primitive meta-morphisms}

For a coalgebra map $f:D\to C$, recall that $\Coder(f)$ and $\Prim_f(\HOM(D,C))$ are respectively the dg-vector spaces of $f$-coderivations and $f$-primitive elements.

\begin{prop}\label{primicoderiv}
Let $f:D\to C$ be a coalgebra map, then there are natural bijections between
\begin{enumerate}
\item maps of vector spaces $p:X\to \Prim_f(\HOM(D,C))$;
\item maps of pointed coalgebras $k:T^c_{\bullet,1}(X)\to (\HOM(D,C),f)$;
\item maps of coalgebras $g:T^c_{\bullet,1}(X)\otimes D\to C$ such that $f=g(e\otimes D):D\to T^c_{\bullet,1}(X)\otimes D\to C$;
\item $f$-derivations $d:X\otimes D\to C$;
\item maps of $D$-bicomodules $X\otimes D\to \Omega^{C,f}$;
\item linear maps $h:X\to \Coder(f)=\hom_{D,D}(D,\Omega^{C,f})$.
\end{enumerate}
If $C=D$ and $f=id_C$, there are more natural bijections with 
\begin{enumerate}
\setcounter{enumi}{6}
\item the maps of bialgebras $T^{csh}(X)\to \END(C)$;
\item the meta-actions $T^{csh}(X)\otimes C\to C$.
\end{enumerate}
\end{prop}
\begin{proof}
If $k$ decomposes into $k_0+k_1$ with respect to the decomposition $T^c_{\bullet,1}(X)=\FF\oplus X$.
By assumption $k_0:\FF\to \HOM(D,C)$ is the atom corresponding to $f$.
The bijection $1\leftrightarrow 2$ is from proposition \ref{univpropofdeltanew} and identifies $k_1=p$.
The bijection $2\leftrightarrow 3$ is by adjunction, we have $k=\Lambda^2g$ and $g={\bf ev}(k\otimes D)$.
The coalgebra $T^c_{\bullet,1}(X)\otimes D=D\oplus X\otimes D$ is of the type $D\oplus N$ where $N$ is a $D$-bicomodule.
Hence, a coalgebra map $g:T^c_{\bullet,1}(X)\otimes D\to C$ decomposes into a coalgebra map $f: D\to C$ and a $f$-derivation $d:X\otimes D\to C$. We have $k_1=\Lambda^2d$ and $d={\bf ev}(k_1\otimes D)$. This proves the bijection $3\leftrightarrow 4$. 
The bijection $4\leftrightarrow 5$ is the definition of $\Omega^{C,f}$.
Recall that the category of bicomodules is tensored over $\dgVect$, hence $D$-bicomodules maps $D\otimes X\to \Omega^{C,f}$ are in bijection with linear maps $D$-bicomodules $X\to \Hom_{D,D}(D,\Omega^{C,f})=\Coder(D;C,f) = \Coder(f)$.
This proves the bijection $5\leftrightarrow 6$. Remark that the bijection $4\leftrightarrow 6$ is given by $h=\lambda^2 d$ and $d=ev(h\otimes D)$.
Finally, the bijections $2\leftrightarrow 7\leftrightarrow 8$ are from proposition \ref{univpropofdeltabialg} and proposition \ref{caracgammacoalg}.
\end{proof}

\begin{cor}\label{corprimicoderiv}
Let $f:D\to C$ be a coalgebra map, the reduction map $\Psi:\HOM(D,C)\to [D,C]$ induces an isomorphism in $\dgVect$
$$
\Prim_f(\HOM(D,C)) = \Coder(f).
$$
\end{cor}
\begin{proof}
By proposition \ref{primicoderiv}, both objects have the same functor of points. They are isomorphic by Yoneda's lemma.
Let us prove that the isomorphism is induced by $\Psi$.
According to the proof of proposition \ref{primicoderiv}, a map $p:X\to \Prim_f(\HOM(D,C))$ is send to the coderivation
$h=\lambda^2{\bf ev}(p\otimes D) = \Psi (p)$ by definition of $\Psi$.
\end{proof}

\begin{cor}
Let $f:D\to C$ be a coalgebra map, then there are natural isomorphism in $\dgVect$ between
$$
\Prim_f(\HOM(D,C)) = \Hom_{D,D}(D,\Omega^{C,f}) = \Hom_{C,C}(D,\Omega^C)
$$
\end{cor}
\begin{proof}
This comes from proposition \ref{primicoderiv} and $\Coder(f) = \Hom_{D,D}(D,\Omega^{C,f})= \Hom_{C,C}(D,\Omega^C)$ 
\end{proof}

\medskip

\begin{nota}\label{notationatomcoderiv}
Recall from lemma \ref{musicallemmacoalg} that if $D$ and $C$ are coalgebras, the map $\Psi: \HOM(D,C)  \to [D,C]$
induces a bijection between the atoms of the coalgebras $\HOM(D,C)$ and the maps of coalgebras $D\to C$.
If $f:D\to C$ is a map of coalgebras, we noted $f^\sharp \in \HOM(D,C)$  the unique atom such that $\Psi(f^\sharp)=f$.

Corollary \ref{corprimicoderiv} says that, if $d:D\rhup_n C$ is a $f$-coderivation, there exists a unique  element $b\in \HOM(D,C)_n$ primitive with respect to $f^\sharp$ such that $\Psi(b)=d$. We shall denote by $d^\sharp$ this element.

Recall the notation $\Psi(b)=b^\flat$ from section \ref{reducetmeta-morphismcoalg}. We have $(d^\sharp)^\flat=d$ by definition.
\end{nota}

With these notations we can explain proposition \ref{primicoderiv} by the commutative diagram in $\dgVect$
$$\xymatrix{
\Prim_f(\HOM(C,D))\ar[r]\ar[d]^\flat  & \HOM(C,D) \ar[d]^{\Psi=\flat}\\
\Coder(f) \ar[r] \ar@/^1pc/[u]^-\sharp & [C,D].
}$$
The maps $\sharp$ and $\flat$ induce inverse bijections $\Prim_f(\HOM(D,C)) = \Coder(f)$.

\bigskip

If $C$ is a coalgebra, then the coalgebra $\END(C)$ has the structure of a bialgebra.
$C$ is a module-coalgebra over the bialgebra $\END(C)$.
The action of this is given by the strong evaluation $\mathbf{ev}: \END(C) \otimes C\to C$.

\begin{thm}\label{primiLiecoderiv}
If $C$ is a coalgebra, we have a commutative diagram in $\dgLie$
$$\xymatrix{
\Prim(\END(C))\ar[r]\ar[d]^\flat  & \END(C) \ar[d]^{\Psi=\flat}\\
\Coder(C) \ar[r] \ar@/^1pc/[u]^-\sharp& [C,C].
}$$
In particular, the maps $\sharp$ and $\Psi(=\flat)$ induce inverse Lie algebra isomorphisms 
$$
\flat :\Prim(\END(C))\simeq\Coder(C):\sharp
$$
which preserve the square of odd elements.
\end{thm}

\begin{proof}
The result is proposition \ref{primicoderiv} but for the Lie structures.
$\END(C)$ and $[C,C]$ are algebras hence Lie algebras.
The inclusion $\Coder(C)\subset [C,C]$ is a Lie algebra map by definition of the Lie structure on $\Coder(C)$.
The map $\flat=\Psi:\END(C)\to [C,C]$ is an algebra map, hence a Lie algebra map which preserve the square of odd elements.
It is bijective when restricted to $\Prim(\END(C))$ by proposition \ref{primicoderiv}. 
$\sharp$ is an injective section of $\flat=\Psi':\Prim(\END(C))\to \Coder(C)$, hence it is an inverse and a Lie algebra morphism.
This proves the last point. Using this isomorphism the map $\sharp$ is the inclusion $\Prim(C)\to \END(C)$ is a Lie algebra morphism which preserve the square of odd elements.
\end{proof}

\bigskip

Recall from proposition \ref{caracgammacoalg} that a $Q$-module-coalgebra is a coalgebra $C$ 
equipped with a left $Q$-module structure defined by an action $a:Q\otimes C\to C$ which is a coalgebra map.

\begin{lemma} \label{primitiveactioncoalgebras} 
Let $C$ be a $Q$-module-coalgebra.
If $b\in Q$ is primitive, then the map $\pi(b):=b\cdot(-):C\to C$ is a coderivation of $C$.
Moreover the map $\pi:\Prim(Q)\to \Coder(C)$ so defined is a homomorphism of Lie algebras which preserves the square of odd elements.
\end{lemma}

\begin{proof} For every $x\in C$, we have
$$\Delta(b\cdot x)=(b\otimes 1+1\otimes b)\cdot (x^{(1)}\otimes x^{(2)})=(b\cdot x^{(1)}) \otimes x^{(2)} + x^{(1)} \otimes 
(b\cdot x^{(2)} ) (-1)^{|b|x^{(1)}||}.$$
This shows that  the map
$b\cdot(-):C\to C$
is a coderivation of degree $|b|$ of the coalgebra $C$.
We have have $\pi([b_1,b_2])=[\pi(b_1),\pi(b_2)]$
since $C$ is a left $Q$-module.
Moreover, we have  $\pi(b^2)=\pi(b)^2$, for every $b\in Q$, for the same reason.
\end{proof}

The following proposition says that the map $\pi:\Prim(Q)\to \Coder(C)$ can be computed either from the action map or from the meta-action map.

\begin{prop}\label{Qprim=Qcoderiv}
The commutative triangle of algebras
$$\xymatrix{
&& \END(C)\ar[d]^\Psi\\
Q\ar[rr]^-\beta\ar[rru]^-\alpha&& [C,C]
}$$
induces a commutative triangle of Lie dg-algebras morphisms preserving the squares of odd elements
$$\xymatrix{
&& \Prim(\END(C))\ar[d]^\Psi_\simeq\\
\Prim(Q)\ar[rr]^-{\beta'}\ar[rru]^-{\alpha'}&& \Coder(C)
}$$
\end{prop}
\begin{proof}
All maps $\alpha$, $\beta$ and $\Psi$ are algebra maps, so they induce Lie algebra maps preserving the square of odd elements.
The result will be proven if we show that these maps restricts to the subspaces of the second diagram.
$\beta$ restricts to a Lie algebra map $\beta':\Prim(Q)\to \Coder(C)$ by lemma \ref{primitiveactioncoalgebras}.
$\alpha$ restricts to a Lie algebra map $\alpha':\Prim(Q)\to \Prim(\END(C))$ by lemma \ref{functoprim}.
and $\Psi$ induces an isomorphism of Lie algebra by proposition \ref{primiLiecoderiv}.
\end{proof}

\bigskip

We finish this section with a characterization of primitive elements and coderivation of cofree coalgebras.

\begin{prop}
\begin{enumerate}
\item The atoms of $\HOM(C,T^\vee(X))$ are in bijection with the morphisms of dg-vector spaces $C\to X$.
\item If $e\in \HOM(C,T^\vee(X))$ is an atom, then 
$$
\Prim_e(\HOM(C,T^\vee(X))) = [C,X].
$$
\item In particular, if $C=T^\vee(X)$, we have
$$
\Coder(T^\vee(X))=[T^\vee(X),X].
$$
\item And if $C=T^c(X)$, we have
$$
\Coder(T^c(X))=[T^c(X),X].
$$
\end{enumerate}
\end{prop}
\begin{proof}
\begin{enumerate}
\item By proposition \ref{atomcofree} the atoms are un bijection with $Z_0([C,X])$ which is in bijection with map of dg-vector space $C\to X$.
\item This is corollary \ref{primcofree} applied to $\HOM(C,T^\vee(X))=T^\vee[C,X]$.
\item When $C=T^\vee(X)$, we deduce from 2. and theorem \ref{primiLiecoderiv} that
$$
\Coder(T^\vee(X)) = \Prim(\END(T^\vee(X))) = \Prim(T^\vee([T^\vee(X),X])) = [T^\vee(X),X].
$$
\item Let $\iota:T^c(X)\to T^\vee(X)$ be the canonical inclusion, the we have the isomorphisms
$$
\Coder(T^c(X)) = \Prim(\END(T^c(X))) = \Prim_\iota(\HOM(T^c(X),T^\vee(X))) = [T^c(X),X].
$$
\end{enumerate}
\end{proof}
The last two statements were already proven in lemmas \ref{coextensionofcoderlemma} and \ref{coextensioncoderconillemma}. 
The proofs given here are more conceptual.

\subsubsection{Derivative of Sweedler operations}\label{derivativesweedlercoalg}

Let $C$ and $D$ be two coalgebras.
Recall that the strengths of $\otimes$ and $\HOM$ give bialgebra maps
$$\xymatrix{
\Theta_\otimes:\END(C)\otimes  \END(D) \ar[r]& \END(C\otimes D)
}$$
$$\xymatrix{
\Theta_{\HOM}:\END(C)^o \otimes  \END(D) \ar[r]& \END(\HOM(C,D)).
}$$
By lemma \ref{derivationmorbialg}, we have the derivative maps between the Lie algebras of primitive elements
$$\xymatrix{
\Theta'_\otimes:\Prim(\END(C))\times  \Prim(\END(D)) \ar[r]& \Prim(\END(C\otimes D))
}$$
$$\xymatrix{
\Theta'_{\HOM}:\Prim(\END(C)^o) \times  \Prim(\END(D)) \ar[r]& \Prim(\END(\HOM(C,D))).
}$$
By theorem \ref{primiLiecoderiv}, these are equivalent to maps of Lie algebras
$$\xymatrix{
\Theta'_\otimes:\Coder(C)\times  \Coder(D) \ar[r]& \Coder(C\otimes D)
}$$
$$\xymatrix{
\Theta'_{\HOM}:\Coder(C) \times  \Coder(D) \ar[r]& \Coder\bigl(\HOM(C,D)\bigr).
}$$
Using the calculus of meta-morphisms, they are given respectively by 
\begin{eqnarray*}
(d_1, d_2) &\mto& \big(d_1^\sharp\otimes D + C\otimes d_2^\sharp\big)^\flat\\
(d_1, d_2) &\mto& \big(\HOM(C, d_2^\sharp) - \HOM(d_1^\sharp, D) \big)^\flat.
\end{eqnarray*}

\medskip
The calculus of meta-morphisms also tells us (proposition \ref{tensormusicalcoalg}) that
$$
\big(d_1^\sharp\otimes D + C\otimes d_2^\sharp\big)^\flat = (d_1^\sharp)^\flat\otimes D + C\otimes (d_2^\sharp)^\flat = d_1\otimes D + C\otimes d_2,
$$
\ie that the coderivation induced by $(d_1,d_2)$ through the strength of $\otimes$ is the classical coderivation constructed on a tensor product (proposition \ref{actioncodertensor}).
However for $\big(\HOM(C, d_2^\sharp) - \HOM(d_1^\sharp, D) \big)^\flat$, the situation is new. 
The following lemma will help us to understand these derivations.

\medskip
Recall from proposition \ref{redmorphismcoalg2}, that if $f\in \END(C)$ and $g\in \END(D)$, we have a commutative square of graded morphisms in $\dgVect$
$$\xymatrix{
\HOM(C,D) \ar[d]_{\Psi=\flat}\ar@^{>}[rrr]^-{\HOM(f,g)^\flat} &&&  \HOM(C,D)\ar[d]^{\Psi=\flat} \\
[C,D] \ar@^{>}[rrr]^-{\hom(f^\flat,g^\flat)}  &&& [C,D].
}$$
The graded morphism $\HOM(f,g)^\flat$ is not in general uniquely determined by $\hom(f^\flat,g^\flat)$. 
The following proposition proves that this is somehow the case for primitive elements.

\begin{prop}\label{uniqueextensioncoder}
If $d_1^\sharp$ and $d_2^\sharp$ are the primitive elements of $\END(C)$ and $\END(D)$ associated to coderivations $d_1$ and $d_2$ of $C$ and $D$ then
$\HOM(C,d_2^\sharp)^\flat$ and $\HOM(d_1^\sharp,D)^\flat$ are coderivations and 
$d=\HOM(C,d_2^\sharp)^\flat - \HOM(d_1^\sharp,D)^\flat$ is the unique coderivation such that the square
$$\xymatrix{
\HOM(C,D) \ar[d]_{\Psi=\flat}\ar@^{>}[rrrr]^-{d} &&&&  \HOM(C,D)\ar[d]^{\Psi=\flat} \\
[C,D] \ar@^{>}[rrrr]^-{\hom(C,d_2)-\hom(d_1,D)}  &&&& [C,D].
}$$
commutes.
Equivalently, $d$ is the unique coderivation such that
$$
d(h)^\flat = d_2h^\flat-h^\flat d_1 \ (-1)^{|h||d_1|}
$$
for any $h\in \HOM(C,T^\vee(X))$.
\end{prop}

\begin{proof}
It is possible to deduce the result from proposition \ref{strengthHOM} but we are going to give a more direct proof.
If $d_1$ and $d_2$ are coderivations, then so are $\HOM(C,d_2^\sharp)^\flat$ and $\HOM(d_1^\sharp,D)^\flat$ by theorem \ref{primiLiecoderiv}. The commutation of the diagram is from proposition \ref{redmorphismcoalg2} and the fact that $d_i=(d_i^\sharp)^\flat$.
The coderivation $\HOM(C,d_2)^\flat-\HOM(d_1,D)^\flat$ is equivalent to a coalgebra map $\HOM(C,D)\oplus S^{-n}\HOM(C,D)\to \HOM(C,D)$ by proposition \ref{fcoderiv} ($n$ is the degree of $d_1$ and $d_2$). Then, the unicity result follows by the separation property of $\Psi:\HOM(C,D)\to [C,D]$.
\end{proof}

If $D=T^\vee(X)$ is cofree, we have the following strengthening of the previous result.

\begin{prop}\label{uniqueextensioncodercofree}
If $d_1^\sharp$ and $d_2^\sharp$ are the primitive elements of $\END(C)$ and $\END(T^\vee(X))$ associated to coderivations $d_1$ and $d_2$, 
then the coderivation $d=\HOM(C, d_2^\sharp)^\flat - \HOM(d_1^\sharp, T^\vee(X))^\flat$ is the unique coderivation on $\HOM(C,T^\vee(X))=T^\vee([C,X])$ such the following square commutes
$$\xymatrix{
T^\vee([C,X]) \ar[d]_{\Psi}\ar@^{>}[rrrr]^-{d} &&&&  T^\vee([C,X])\ar[d]^{q} \\
[C,T^\vee(X)] \ar@^{>}[rrrr]^-{\hom(C, pd_2) - \hom(d_1, p)}  &&&& [C,X].
}$$
where $p:T^\vee(X)\to X$ and $q:T^\vee([C,X])\to [C,X]$ are the cogenerating maps.
Equivalently, $d$ is the unique coderivation such that
$$
q(d(h)) = pd_2h^\flat-ph^\flat d_1 \ (-1)^{|h||d_1|}
$$
for any $h\in \HOM(C,T^\vee(X))$.
\end{prop}

\begin{proof}
We have a commutative diagram
$$\xymatrix{
\HOM(C,T^\vee(X)) \ar[d]_{\Psi}\ar@^{>}[rrrr]^-{\HOM(C, d_2^\sharp)^\flat - \HOM(d_1^\sharp, T^\vee(X))^\flat} &&&&  \HOM(C,T^\vee(X))\ar[d]^{\Psi} \ar@{=}[r]&T^\vee([C,X])\ar[d]^q\\
[C,T^\vee(X)] \ar@^{>}[rrrr]^-{\hom(C, d_2) - \hom(d_1, D)}  &&&& [C,T^\vee(X)] \ar[r]^-{[C,p]}& [C,X].
}$$
Then the proof is the same as in lemma \ref{uniqueextensioncoder} but using the separating property of $q$ instead of that of $\Psi$.
\end{proof}

\bigskip

We have presented how to transport coderivations along the tensor product and internal hom of dg-coalgebras, but it should be clear for the reader that everything could be done the same way for graded coalgebras. There actually lies our main application of these formulas.

Recall from theorem \ref{dgtoghomcoalg} that the functor forgetting the differential $U_d:\dgCoalg\to \gCoalg$ preserve the tensor product and the internal hom. Let $(C,d_C)$ and $(D,d_D)$ be two dg-coalgebras viewed as graded coalgebra equipped with coderivation $d_C$ and $d_D$. It is classical that the tensor product $(C,d_C)\otimes (D,d_)D$ in dg-coalgebras can be described as the tensor product $C\otimes D$ of the underlying graded coalgebras equipped with the differential $d_C\otimes D+C\otimes d_D$.
This can be read a different way using the internal hom of coalgebras.
The pair $(d_C,d_D)$ of coderivation define a primitive element in the graded coalgebra $\END(C)\otimes \END(D)$, which can be transported as a primitive element in $\END(C\otimes D)$ using the strength of $\otimes$. The formula for this primitive element is given by the calculus of meta-morphisms as $d_C\otimes D+C\otimes d_D$. 
The differential enhancement of the graded coalgebras $C\otimes D$ is the same computed classically or using the enrichment of $\otimes$ over $\gCoalg$. 
This is an important feature of the theory which center is theorem \ref{primiLiecoderiv}.

The situation is more interesting for the internal hom. With the same notations, theorem \ref{dgtoghomcoalg} describes the dg-coalgebra $\HOM(C,d_C;D,d_D)$ as the graded coalgebra $\HOM(C,D)$ equipped with a differential. Again the calculus of meta-morphisms
tells us how to define the differential on $\HOM(C,D)$: it is the derivation $d:\HOM(C,d_D)-\HOM(d_C,D)$ (we drop the musical signs for simplicity).
One good thing about the calculus of meta-morphisms is that an odd derivation is of square zero iff the corresponding meta-morphism is of square zero, in particular they are preserved by any strong functor. We can compute that $d$ is indeed of square zero:
\begin{eqnarray*}
(\HOM(C,d_D)-\HOM(d_C,D))^2 &=& \HOM(C,d_D)\HOM(C,d_D) - \HOM(C,d_D)\HOM(d_C,D) \\
	&& - \HOM(d_C,D)\HOM(C,d_D) + \HOM(d_C,D)\HOM(d_C,D)\\
&=& \HOM(C,d^2_D) - \HOM(d_C,d_D)(-1) \\
	&&- \HOM(d_C,d_D) + \HOM(d_C^2,D)\\
&=& 0.
\end{eqnarray*}

\subsubsection{The enrichment of vector spaces over coalgebras}\label{enrichmentVectoverCoalg}

Recall that $\dgVect^\$$ is our notation for $\dgVect$ viewed as enriched over itself.
The functor $T^\vee:\dgVect \to \dgCoalg$ is left adjoint to the monoidal functor $U:\dgCoalg \to \dgVect$, hence it is a lax monoidal functor (see appendix \ref{laxmonoidaladj}). 
We can thus tranfer the enrichment of $\dgVect^\$$ along $T^\vee$ into an enrichment over $\dgCoalg$ which we denote $\dgVect^{T^\vee\$}$. The new hom object between two vector spaces $X$ and $Y$ is $T^\vee [X,Y]$.
The composition law is defined to be the composite 
$$\xymatrix{
{\bf c}:T^\vee [Y,Z]\otimes T^\vee[X,Y] \ar[r]^-\alpha & T^\vee([Y,Z] \otimes [X,Y]) \ar[rr]^-{ T^\vee(c)}  && T^\vee[X,Z],
}$$
where $\alpha$ is the lax structure on the functor $T^\vee$ and where $c$ is the strong composition law $[Y,Z] \otimes [X,Y] \to [X,Z]$ in the category $\dgVect^\$$.
The unit $e_X:\FF\to T^\vee([X,X])$ of $X\in \dgVect^{T^\vee\$}$ is defined to be the composite 
$$\xymatrix{
e_X:\FF\ar[r]^-{\alpha_0} & T^\vee(\FF)\ar[rr]^-{T^\vee(u_X)} && T^\vee[X,X]}
$$
where $\alpha_0$ is the unit of the lax structure on the functor $T^\vee$.

\begin{prop} \label{compositionandcofreemap}
The composition law  ${\bf c}:T^\vee [Y,Z]\otimes T^\vee[X,Y] \to T^\vee [X,Z]$
is the unique map of coalgebras for which following square commutes,
$$\xymatrix{
T^\vee [Y,Z]\otimes T^\vee[X,Y] \ar[rr]^-{{\bf c}} \ar[d]_{p\otimes p} &&  T^\vee [X,Z] \ar[d]^p\\
[Y,Z]\otimes [X,Y] \ar[rr]^-c  && [X,Z], \\
}
$$
where $p$ denotes the cofree maps and $c$ is the composition law in $\dgVect^\$$.  The unit of an object $X\in \dgVect^{T^\vee\$}$
is the unique atom $e_X\in T^\vee [X,X]$ such that $p(e_X)=1_X$.
\end{prop}

\begin{proof} Let us show that the square commutes.  
The triangle of the following diagram commutes by construction of $\alpha$ in appendix \ref{laxmonoidaladj}.

$$\xymatrix{
T^\vee [Y,Z]\otimes T^\vee[X,Y] \ar[rr]^\alpha \ar[drr]_-{p\otimes p}&& T^\vee([Y,Z] \otimes [X,Y]) \ar[rr]^(0.6){ T^\vee(c)} \ar[d]^p && T^\vee[X,Z] \ar[d]^p \\
&& [Y,Z]\otimes [X,Y] \ar[rr]^c  && [X,Z]. \\
}
$$
The square on the right hand side commutes by naturality of the cofree maps $p$.
It follows that the boundary diagram commutes. The uniqueness of $c'$
is clear, since the map $p: T^\vee [X,Z] \to  [X,Z] $ is cogenerating.
Let us show that $p(e_X)=1_X$.  
The triangle of the following diagram commutes by construction of $\alpha_0$ in appendix \ref{laxmonoidaladj},
$$\xymatrix{
\FF \ar[rr]^-{\alpha_0} \ar@{=}[drr] && T^\vee(\FF) \ar[rr]^(0.5){ T^\vee(u_X)} \ar[d]^p && T^\vee[X,X] \ar[d]^p \\
&&\FF \ar[rr]^-{u_X}  && [X,X]. \\
}
$$
Thus, $pe_X=u_X$ and this shows that $p(e_X)=1_X$, since $u_X(1)=1_X$.
The uniqueness of $e_X$ is clear, since the map $p: T^\vee [X,X] \to  [X,X] $ is cogenerating.
\end{proof}

\bigskip
\begin{prop}\label{vectbicompletecoalg} \label{TVectclosed}
The category $\dgVect^{T^\vee\$}$ is bicomplete and symmetric monoidal closed over $\dgCoalg$.
Moreover any ordinary (co)limit is automatically strong.
\end{prop}
\begin{proof}
Recall that an enriched category is bicomplete if it is tensored and cotensored and its underlying category is bicomplete. 
$\dgVect$ is bicomplete, so we need only to prove that $\dgVect^{T^\vee\$}$ is tensored and cotensored.
This a formal consequence of $T^\vee$ having a left adjoint (see proposition \ref{transferadjunction}).
Let $C$ be a coalgebra and $X$ and $Y$ be two vector spaces, we have natural bijection between
\begin{center}
\begin{tabular}{lc}
\rule[-2ex]{0pt}{4ex} coalgebra maps & $C\to T^\vee([X,Y])$\\
\rule[-2ex]{0pt}{4ex} linear maps & $UC\to [X,Y]$\\
\rule[-2ex]{0pt}{4ex} linear maps & $UC\otimes X\to Y$\\
\rule[-2ex]{0pt}{4ex} lienar maps & $X\to [UC,Y]$.
\end{tabular}
\end{center}
This proves that $UC\otimes X$ and $[UC,X]$ are the tensor and cotensor in $\dgVect^{T^\vee\$}$.
The proof of the symmetric closed structure is a consequence of proposition \ref{transferadjunction}.

Let us now prove that all limits are strong, the argument will be similar for colimits.
Let $X:I\to \dgVect$ be a diagram with ordinary limit $Y\in \dgVect$ viewed as the equalizer
$$\xymatrix{
Y \ar[r]&  \prod_{i\in I}X_i \ar@<.6ex>[r]\ar@<-.6ex>[r]& \prod_{i\to j\in I}X_i.
}$$
The universal property of $Y$ as a strong limit is, for any $E\in \dgVect$, the exactness in $\dgCoalg$ of the diagram
$$\xymatrix{
T^\vee([E,Y]) \ar[r]&  \prod_{i\in I}T^\vee([E,X_i]) \ar@<.6ex>[r]\ar@<-.6ex>[r]& \prod_{i\to j\in I}T^\vee([E,X_i]) 
}$$
which we deduce from the previous equalizer by the commutation of $T^\vee$ and $[E,-]$ with limits.
\end{proof}

\medskip

We generalized now the notion of meta-morphisms to vector spaces.

\begin{defi}
We shall say that an element $f\in T^\vee[X,Y]$ is a {\it meta-morphism of vector spaces} and we shall denote it by $f:X\leadsto Y$.
\end{defi}

The {\it composite} $gf=g\circ f$ of two meta-morphisms $f:X\leadsto Y$ and $g:Y\leadsto Z$ is defined by putting $gf=\mathbf{c}(g\otimes f)$. We have $|gf|=|g|+|f|$. where $\bf c$ is the composition of proposition \ref{compositionandcofreemap}.

The category $\dgVect^{T^\vee \$}$ is symmetric monoidal closed and bicomplete, the tensor product of $X$ by a coalgebra $C$ is $UC\otimes X$ and the cotensor is $[UC,X]$. The {\em strong evaluation} map ${\bf ev }:UT^\vee([X,Y])\otimes X\to Y$ is given by the composition 
$$\xymatrix{
UT^\vee([X,Y])\otimes X\ar[rr]^-{p\otimes X} &&[X,Y]\otimes X\ar[r]^-{ev} & Y
}$$
where $p:T^\vee([X,Y])\to [X,Y]$ is the cogenerating map and $ev$ is the evaluation in $\dgVect^\$$.

If $f:X\leadsto Y$ and $x\in X$, the {\em evaluation} of a meta-morphism $f:X\leadsto Y$ on $x$ is defined by $f(x):={\bf ev}(f\otimes x)=p(f)(x)$.

\medskip

\begin{prop}\label{reductionmappingandmeta} If $f:X\leadsto Y$
and $g:Y\leadsto Z$ are meta-morphisms of vector spaces, then we have
$(gf)(x)=g(f(x))$ for every $x\in X$.
\end{prop}

\begin{proof}
The commutative square of proposition \ref{compositionandcofreemap} shows that we have $p(gf)=p(g)p(f)$.
Thus $(gf)(x)=p(gf)(x)=p(g)p(f)(x)=g(f(x))$.
\end{proof}

A morphism $X\to Y$ in the category  $\dgVect^{T^\vee\$}$ is a map of coalgebras $u:\FF\to T^\vee[X,Y]$, equivalently
it is an atom $u(1)$ of the coalgebra  $T^\vee[X,Y]$.
Let us denote by $At(X,Y)$ the set of atoms of the coalgebra  $T^\vee [X,Y]$.

\begin{lemma}
The map $At(X,Y)\to [X,Y]_0$ induced by the cofree map $p:T^\vee [X,Y]\to [X,Y]$ induce a bijection
$$
At(X,Y) \simeq \dgVect(X,Y).
$$
The underlying category of $\dgVect^{T^\vee\$}$ is $\dgVect$.
\end{lemma}
\begin{proof}
We have the following bijections
\begin{center}
\begin{tabular}{lc}
\rule[-2ex]{0pt}{4ex} coalgebra maps & $\FF\to T^\vee([X,Y])$\\
\rule[-2ex]{0pt}{4ex} linear maps & $\FF\to [X,Y]$\\
\rule[-2ex]{0pt}{4ex} maps of dg-vector spaces & $X\to Y$.
\end{tabular}
\end{center}
This proves the first assertion. 
The second is a consequence of proposition \ref{transferadjunction}.
\end{proof}

\begin{nota}\label{notationatomvslinearmap}
If $f:X\to Y$ is a linear map, we shall denote by $f^\sharp$ the unique atom $e\in At(X,Y)$ such that $p(e)=f$.
If $f:X\leadsto Y$ is a meta-morphism, we shall denote by $f^\flat$ the linear map $p(f)$.
The previous lemma says that $(f^\sharp)^\flat=f$, $(fg)^\sharp = f^\sharp g^\sharp$ for $f$ and $g$ two linear maps 
and that $(fg)^\flat = f^\flat g^\flat$ for $f$ and $g$ two meta-morphisms. We have also $f(x)=f^\flat(x)$ by definition of the evaluation.
\end{nota}

\bigskip

By proposition \ref{transferadjunction}, the tensor product functor of vector spaces defines a strong functor
$$
\otimes :\dgVect^{T^\vee\$} \times \dgVect^{T^\vee\$} \to \dgVect^{T^\vee\$}
$$
Its strength $\Theta_\otimes$ is the composite 
$$\xymatrix{
T^\vee[X_1,X_2]\otimes  T^\vee[Y_1,Y_2] \ar[rr]^\psi  && T^\vee([X_1,X_2]\otimes [Y_1,Y_2]) \ar[rr]^-{T^\vee(\theta)} && T^\vee [X_1\otimes Y_1,X_2\otimes Y_2], \\
}$$
where $\psi$ is the lax structure on $T^\vee$ and $\theta$ is the strength of the tensor product functor in $\dgVect$.

\begin{prop} \label{strongtensorcofreemapvect}
The strength $\Theta_\otimes:T^\vee[X_1,X_2]\otimes  T^\vee[Y_1,Y_2]  \to T^\vee[X_1\otimes Y_1,X_2\otimes Y_2]$
is the unique map of coalgebras for which following square commutes,
$$\xymatrix{
T^\vee[X_1,X_2]\otimes  T^\vee[Y_1,Y_2]  \ar[rr]^-{\Theta_\otimes} \ar[d]_{p\otimes p} &&  T^\vee([X_1\otimes Y_1,X_2\otimes Y_2])\ar[d]^p\\
[X_1,X_2]\otimes  [Y_1,Y_2]  \ar[rr]^\theta  &&[X_1\otimes Y_1,X_2\otimes Y_2], 
}
$$
where $p$ denotes the cofree maps and  $\theta$ the strength of the tensor product in $\dgVect$.
\end{prop}

\begin{proof}
The triangle of the following diagram commutes by construction of $\alpha$ in appendix \ref{laxmonoidaladj},
$$\xymatrix{
T^\vee[X_1,X_2]\otimes  T^\vee[Y_1,Y_2]   \ar[rr]^\alpha \ar[drr]_-{p\otimes p}&& T^\vee([X_1,X_2]\otimes  [Y_1,Y_2] ) \ar[rr]^-{ T^\vee(\theta)} \ar[d]^p && T^\vee[X_1\otimes Y_1,X_2\otimes Y_2] \ar[d]^p \\
&& [X_1,X_2]\otimes  [Y_1,Y_2]  \ar[rr]^\theta  &&[X_1\otimes Y_1,X_2\otimes Y_2],
}
$$
The square on the right hand side commutes by naturality of the cofree maps $p$.
It follows that the boundary diagram commutes. The uniqueness of $\Theta$
is clear, since the map $p: T^\vee[X_1\otimes Y_1,X_2\otimes Y_2]   \to [X_1\otimes Y_1,X_2\otimes Y_2] $ is cogenerating.
\end{proof}

If $f:X_1\leadsto X_2$ and $g:Y_1\leadsto Y_2$ are  meta-morphisms of vector spaces
let us put 
$$f\otimes g:=\Theta(f\otimes g):X_1 \otimes Y_1 \leadsto X_2 \otimes Y_2$$
where  $\mu':T^\vee[X_1,X_2]\otimes  T^\vee[Y_1,Y_2]  \to T^\vee[X_1\otimes Y_1,X_2\otimes Y_2]$
is the  strength \ref{strongtensorcofreemapvect} of the tensor product in $\dgVect^{T^\vee\$}$.

\begin{prop}\label{metatensorvect1} 
If $f':X_1\leadsto X_2$,  $f:X_2\leadsto X_3$, $g':Y_1\leadsto Y_2$ and $g:Y_2\leadsto Y_3$
are meta-morphisms of vector spaces, then we have 
$$
(f \otimes g)(f' \otimes g')=ff' \otimes gg' (-1)^{|f'||g|}.
$$
\end{prop}

\begin{proof}
The identity follows from the functoriality of the strong functor $\otimes :\dgVect^{T^\vee\$} \times\dgVect^{T^\vee\$} \to \dgVect^{T^\vee\$}$.
\end{proof}

\begin{cor}\label{metatensorvect2} 
If $f:X_1\leadsto X_2$ and $g:Y_1\leadsto Y_2$ are meta-morphisms of vector spaces,
then we have $$(f\otimes g)(x\otimes y)=f(x)\otimes g(y) (-1)^{|g||x|}$$
for every $x\in X_1$ and $y\in Y_1$. 
\end{cor}

\bigskip

The hom functor 
$$
(\dgVect^{T^\vee\$})^{op} \times \dgVect^{T^\vee\$} \to \dgVect^{T^\vee\$}
$$
is strong with respect to the enrichment of the category $\dgVect^{T^\vee\$}$ over itself,
since the category $\dgVect^{T^\vee\$}$ is symmetric monoidal closed.
Its strength $\Theta_{[-,-]}$ is the composite 
$$\xymatrix{
T^\vee[X_2,X_1]\otimes  T^\vee[Y_1,Y_2] \ar[rr]^\psi  && T^\vee([X_2,X_1]\otimes [Y_1,Y_2]) \ar[rr]^-{T^\vee(\nu)} && T^\vee [[X_1,Y_1],[X_2,Y_2]] \\
}$$
where $\psi$ is the lax structure on $T^\vee$  and $\nu$ is the strength of the hom functor in $\dgVect$.

\begin{prop} \label{strongtensorcofreemap}
The strength $\Theta_{[-,-]}:T^\vee[X_2,X_1]\otimes  T^\vee[Y_1,Y_2]  \to T^\vee [[X_1,Y_1],[X_2,Y_2]] $
is the unique map of coalgebras for which following square commutes,
$$\xymatrix{
T^\vee[X_2,X_1]\otimes  T^\vee[Y_1,Y_2] \ar[rr]^-{\Theta_{[-,-]}} \ar[d]_{p\otimes p} && T^\vee [[X_1,Y_1],[X_2,Y_2]]  \ar[d]^p\\
[X_2,X_1]\otimes [Y_1,Y_2] \ar[rr]^\theta  &&[[X_1,Y_1],[X_2,Y_2]], 
}
$$
where $p$ denotes the cofree maps and  $\theta$ the strength of the hom functor in $\dgVect$.
\end{prop}

\begin{proof}
The triangle of the following diagram commutes by construction of $\alpha$ in appendix \ref{laxmonoidaladj},
$$\xymatrix{
T^\vee[X_2,X_1]\otimes  T^\vee[Y_1,Y_2]    \ar[rr]^-\alpha \ar[drr]_-{p\otimes p}&& T^\vee([X_1,X_2]\otimes  [Y_1,Y_2] ) \ar[rr]^(0.6){ T^\vee(\mu)} \ar[d]^p && T^\vee[[X_1,Y_1],[X_2,Y_2]]
 \ar[d]^p \\
&&[X_2,X_1]\otimes [Y_1,Y_2] \ar[rr]^\mu  &&[[X_1,Y_1],[X_2,Y_2]],
}
$$
The square on the right hand side commutes by naturality of the cofree maps $p$.
It follows that the boundary diagram commutes. The uniqueness of $\Theta_{[-,-]}$
is clear, since the map $p: T^\vee[[X_1,Y_1],[X_2,Y_2]]   \to [[X_1,Y_1],[X_2,Y_2]],$ is cogenerating.
\end{proof}

If $f:X_2\leadsto X_1$ and $g:Y_1 \leadsto  Y_2$ are meta-morphisms of vector spaces, let us put 
$$\xymatrix{
[g,f]:=\Theta_{[-,-]}(f\otimes g): [X_1,Y_1]  \ar@{~>}[r] &  [X_2,Y_2] 
}$$
where $\Theta_{[-,-]}:T^\vee[X_2,X_1]\otimes  T^\vee[Y_1,Y_2]  \to T^\vee [[X_1,Y_1],[X_2,Y_2]] $
is the strength of of the hom functor in $\dgVect^{T^\vee\$}$.

\begin{prop}\label{redmorphismvect1} 
If $f':X_3\leadsto X_2$,  $f:X_2\leadsto X_1$, $g:Y_1\leadsto Y_2$ and $g':Y_2\leadsto Y_3$ 
are  meta-morphisms of vector spaces, then we have 
 $$[f',g'][f,g]=[ff',g'g](-1)^{|f|(|f'|+|g'|)}.$$
\end{prop}

\begin{proof}
The identity follows from the functoriality of the strong hom functor  
$(\dgVect^{T^\vee\$})^{op} \times \dgVect^{T^\vee\$} \to \dgVect^{T^\vee\$}$.
\end{proof}

\begin{prop}\label{redmorphismcoalg3} 
If $f:X_2\leadsto X_1$ and $g:Y_1\leadsto Y_2$ are  meta-morphisms of vector spaces
then we have $$[f,g](h)=ghf(-1)^{|f|(|g|+|h|)}$$
for every meta-morphism $h:X_1\leadsto Y_1$.
\end{prop}

\begin{proof} This follows from \ref{strongtensorcofreemap}.
\end{proof}

\subsubsection{Strong comonadicty}

We finish this section by a strengthening of theorem \ref{comonadic1}.

\begin{thm}\label{strongcomonadicitycoalg}
The adjunction $U\dashv T^\vee$ enriches into a strong lax monoidal adjunction
$$\xymatrix{
U: \dgCoalg^\$ \ar@<.6ex>[r]& \dgVect^{T^\vee\$}:T^\vee. \ar@<.6ex>[l]
}$$
Moreover the adjunction is strong comonadic.
\end{thm}
\begin{proof}
Let us prove first that $U$ and $T^\vee$ are strong functors.
By proposition \ref{vectbicompletecoalg}, both categories are bicomplete so we can use proposition \ref{triplestrength} to describe the strength of $U$ and $T^\vee$ as (co)lax structures.
The lax structure of $U$ is given by the isomorphism
$$\xymatrix{
UC \otimes UD \ar[r]^-\simeq & U(C\otimes D),
}$$
the conditions of \ref{laxmodularstructure} are equivalent to the pentagon and the unit identities for $\otimes$ in $\dgVect$.
The colax structure of $T^\vee$ is given by the isomorphism
$$\xymatrix{
T^\vee [UC,X]\ar[r]^-\simeq & \HOM(C,T^\vee X)
}$$
of proposition \ref{exemplehomcoalg}.
The unit condition of \ref{colaxmodularstructure} is clear 
and the associativitiy condition is equivalent to the commutation of
$$\xymatrix{
\HOM(D\otimes C,T^\vee X) \ar[d] \ar@{=}[r] &  \HOM(D,\HOM(C,T^\vee X)) &\ar[l]_-\simeq \HOM(D,T^\vee [UC,X]) \ar[d]\\
[UD\otimes UC,X]  && \ar@{=}[ll]  [UD,[UC,X]]
}$$
which is an easy consequence of how the isomorphism $T^\vee [UC,X]\simeq \HOM(C,T^\vee X)$ is constructed.
Finally, the natural isomorphism $T^\vee [UC,X]\simeq \HOM(C,T^\vee X)$ is also the strength of the adjunction $U\dashv T^\vee$.

Let us prove now that $U$ is strong symmetric monoidal, the strong symmetric lax monoidal structure of $T^\vee$ will then be a formal consequence.
We need to prove the commutation of the diagram
$$\xymatrix{
\HOM(C_1,C_2) \otimes \HOM(D_1,D_2)\ar[d] \ar[r] & T^\vee([UC_1,UC_2])\otimes T^\vee([UC_1,UC_2])\ar[d]\\
\HOM(C_1\otimes D_1,C_2\otimes D_2)\ar[d] & T^\vee([UC_1,UC_2]\otimes [UC_1,UC_2])\ar[d]\\
T^\vee([U(C_1\otimes D_1),U(C_2\otimes D_2)]) \ar@{=}[r]&  T^\vee([UC_1\otimes UD_1,UC_2\otimes UD_2])
}$$
By adjunction this is equivalent to the commutation of
$$\xymatrix{
\HOM(C_1,C_2) \otimes \HOM(D_1,D_2)\ar[d]_{\Theta_\otimes} \ar[r] \ar[rd]_{\Psi\otimes \Psi}& T^\vee([UC_1,UC_2])\otimes T^\vee([UC_1,UC_2])\ar[d]^{p\otimes p}\\
\HOM(C_1\otimes D_1,C_2\otimes D_2)\ar[d]_\Psi & [UC_1,UC_2]\otimes [UC_1,UC_2]\ar[d]^{\theta}\\
[U(C_1\otimes D_1),U(C_2\otimes D_2)] \ar@{=}[r]&  [UC_1\otimes UD_1,UC_2\otimes UD_2]
}$$
But the upper triangle commute by construction of the strength of $U$
and the above square commute by definition of $\Theta_\otimes$.
We leave the proof of the unit and symmetry conditions to the reader.

Let us prove now the strong comonadicity statement. It will essentially be a consequence of the already proven (co)monadicty.
We are going to use the enriched monadicity result of \cite[thm. II.2.1]{Dubuc}.
We need to prove that $U$ detects, preserves and reflects $U$-split $\dgCoalg$-equalizers.
From the comonadicity theorem \ref{comonadic1}, and the classical monadicity them \cite{Barr-Wells} we know that $U$ detects, preserves and reflects ordinary $U$-split equalizers. First, recall that a strong equalizer is called $U$-split if it is $U$-split as an ordinary equalizer, then we need to compare strong and ordinary equalizers. By proposition \ref{stronglimcoalg} (resp. \ref{vectbicompletecoalg}) any limit in $\dgCoalg$ (resp. $\dgVect$) is automatically a strong limit in $\dgCoalg^{\$}$ (resp. $\dgVect^{T^\vee\$}$). Hence, equalizers and $\dgCoalg$-equalizers coincide in both $\dgCoalg$ and $\dgVect$, hence $U$ detects, preserves and reflects $U$-split strong equalizers.
\end{proof}

\begin{rem}
The functors $U:\dgCoalg^{U\$}\to \dgVect^\$$ and $U:\dgCoalg^{\$}\to \dgVect^{T^\vee\$}$ corresponds to each other through the adjunction of proposition \ref{transferadjunction}.

\end{rem}

\begin{rem} \label{calculHOM}\label{remreduction}
The strong comonadicity of $U\dashv T^\vee$ over $\dgCoalg$ has some interesting consequence for the computation of $\HOM(C,D)$.
Corollary \ref{HOMsubcoalg} proves that the maps $\HOM(C,D)\to T^\vee([UC,UD])$ defining the functor $U:\dgCoalg^{\$}\to \dgVect^{T^\vee\$}$ are injection of coalgebra, hence the functor is strongly faithful.
This property can be understood as subsuming our proofs by reduction (see section \ref{proofbyreduction}): morally, $\dgVect^{T^\vee\$}$ keeps all the structures of $\dgVect$ and the structure of $\dgCoalg^\$$ is inherited from that of $\dgVect^{T^\vee\$}$. In particular, the coalgebras $\HOM(C,D)$ can be computed via the (enriched version) of the usual computation of hom of coalgebras. In the same that that we have the common equalizer in $\Set$ (see appendix \ref{commonequalizer})
$$\xymatrix{
\dgCoalg(C,D)\ar[r]&\dgVect(C,D) \ar@<.6ex>[d]^-{\alpha'}\ar@<-.6ex>[d]_-{\beta'}\ar@<.6ex>[r]^-{\alpha}\ar@<-.6ex>[r]_-{\beta}&\dgVect(C,D\otimes D)\\
& \dgVect(C,\FF)
}$$
(where for $f:C\to D$, we define  $\alpha (f)=(f\otimes f)\Delta_C:C\to D\otimes D$, $\beta (f)=\Delta_Df:C\to D\otimes D$,
 $\alpha' (f)=\epsilon_Cf:C\to \FF$ and $\beta' (f)=\epsilon_Df:C\to \FF$),
we have the common equalizer in $\dgCoalg$
$$\xymatrix{
\HOM(C,D)\ar[r]^-{\Delta'} &T^\vee([C,D]) \ar@<.6ex>[d]^-{T^\vee(\alpha')}\ar@<-.6ex>[d]_-{T^\vee(\beta')}\ar@<.6ex>[r]^-{T^\vee(\alpha)}\ar@<-.6ex>[r]_-{T^\vee(\beta)}& T^\vee([C,D\otimes D])\\
& T^\vee([C,\FF]).
}$$
where $\Delta'$ is the map from corollary \ref{HOMsubcoalg}. 
In other words, the meta-morphisms of coalgebras are the meta-morphisms of vector spaces preserving the coalgebra structure.
\end{rem}

\newpage
\section{The category of algebras}\label{algebras}\label{SweedlerTheory}

This chapter will deal solely with dg-algebras and dg-vector spaces.
We shall use the same convention as in the previous chapter and simplify the langage by refering to dg-algebras and dg-vector spaces simply as {\em algebras} and {\em vector spaces} and by removing the "dg" prefix as often as possible. In particular a map of dg-vector spaces will be simply called a linear map or simply a map.

The chapter contains the following results.
\begin{itemize}

\item The category of (dg-)algebras $\dgAlg$ is locally presentable 
and the adjunction $T:\dgVect\leftrightarrows \dgAlg:U$ is monadic.

\item The category of (dg-)algebras $\dgAlg$ is enriched over the category $\dgCoalg$; the hom object between two (graded) algebras 
$A$ and $B$ is the (graded) measuring coalgebra $\{A, B\}$ introduced by Sweedler (section \ref{sweedlerhom}).

\item The enriched category  $\dgAlg$ is tensored over the category $\dgCoalg$.
The tensor product of an algebra $A$ by a coalgebra $C$ is an algebra $C\rhd A$ called the {\it Sweedler product}  of $A$ by $C$ (section \ref{sweedlerproduct}).

\item The enriched category  $\dgAlg$  is cotensored over the category $\dgCoalg$.
The cotensor product of an algebra $A$ by a coalgebra $C$ is the convolution algebra $[C,A]$ (section \ref{convolution}).

\item The Sweedler product functor $\rhd:\dgCoalg\times \dgAlg\to \dgAlg$ has a colax structure
and the Sweedler hom functor $\{-,-\}:\dgAlg^{op}\times \dgAlg\to \dgCoalg$ has a lax structure
(corollaries \ref{Sproductcolax} and \ref{Shomlax}).

\item The tensor product functor $\otimes: \dgAlg\times \dgAlg \to \dgAlg$ is strong with respect to the enrichment of the category $\dgAlg$ (theorem \ref{tensorenrichmentof}).

\item The operations $\rhd$, $[-,-]$, $\{-,-\}$ and $\otimes$ are endowed with natural lax or colax structure with respect to the tensor products of algebras and coalgebras (section \ref{laxstructures}).

\item The adjunction $T:\dgVect\rightleftarrows\dgAlg:T$ is in fact strongly monadic if both categories are enriched over $\dgCoalg$ (theorem \ref{strongmonadicityalg}).

\end{itemize}
\medskip

In the text, we shall talk of {\em Sweedler operations} to refer to the set of operations $\rhd$, $[-,-]$, $\{-,-\}$ and $\otimes$ on algebras.

\subsection{Presentability and monadicity}

The following result is classical but we give a quick proof for sake of completude.
Recall that the forgetful functor $U:\dgAlg\to dgVect$ has a left adjoint $T$ which send $X$ to the tensor algebra $TX$.

\begin{thm}\label{monadicalg}
The category $\dgAlg$ is locally presentable and monadic over $\dgVect$.
\end{thm}
\begin{proof}
The $\omega$-compact objects of $\dgAlg$ are easily characterized as the algebras of finite presentation.
In particular the algebras freely generated by a finite graded set are $\omega$-compact.
Every free algebra can be written as a colimit of finitely generated free algebras and every algebra can be written as a colimit of free algebras.
This proves the finite presentability of $\dgAlg$.
An explicit computation proves that an algebra over the monad $UT$ on $\dgVect$ is the same thing as an associative unital algebra.
This proves the monadicity.
\end{proof}

\begin{rem}\label{differencealgcoalg}
The notions of compact objects and presentable categories have a definition in terms of colimits rather than limits.
The proofs of the above results are easier for algebras than for coalgebras because algebras are easily seen as generated by colimits from free algebras. On the other side, coalgebras are naturally described as limits of cofree coalgebras.
This explain the a priori difficulty to characterized the $\omega$-compact objects in $\dgCoalg$ and to write every coalgebra as a colimit of such. This difficulty could be overcome if the tensor product were commuting with all limits (which would provide an explicit construction of the cofree coalgebra) but this is not the case.
\end{rem}

\subsection{Convolution}\label{convolution}

Recall that, if $A=(A,m,e)$ is a dg-algebra and $C=(C,\Delta,\epsilon)$ is a dg-coalgebra, 
then the dg-vector space $[C,A]$ has the structure of a dg-algebra with respect to the {\em convolution product} defined as the composition
$$\xymatrix{
\star : [C,A]\otimes [C,A] \ar[rr]^-{can} && [C\otimes C,A\otimes A] \ar[rr]^-{[\Delta,m]} && [C,A].
}$$
Explicitly, the convolution product of two graded morphisms $f,g\in [C,A]$ is a graded morphism of degree $|f|+|g|$ which evaluated on an element $x\in C$ is 
$$
(f\star g)(c)=(f\otimes g)(c^{(1)}\otimes c^{(2)})= f(c^{(1)})g(c^{(2)})\ (-1)^{|g||c^{(1)}|}.
$$
The unit element for the convolution product is the map $e\epsilon:C\to \FF \to A$.
We shall say that $([C,A],\star,e\epsilon)$ is the {\it convolution dg-algebra}.

The construction of the convolution make sense in the non-unital and non-differential cases.

\begin{ex}[Dual]\label{dualcoalg}
The graded dual $C^\star=[C,\FF ]$ of a dg-coalgebra $C$ has the structure of a dg-algebra.
The unit element of $C^\star$ is the counit of $C$. 
If $E$ is a finite algebra, then $E^\star$ is a finite coalgebra and $(E^\star)^\star =E$.
This can be proven by an explicit computation left to the reader.
\end{ex}

\begin{ex}[Product]\label{exconvolproduct}
If $\FF I$ is the diagonal coalgebra of a set $I$, then the convolution algebra $[\FF I,A]$ is isomorphic to the algebra $A^I=A\times \dots \times A$.
Indeed, on one side the product in $A^I$ of two families $(a_i)$ and $(b_i)$ is $(a_ib_i)$.
on the other side if $a,b:\FF I\to A$ are the functions such that $a(i)=a_i$ and $b(i)=b_i$, the convolution product $a\star b$ is defined by
$(a\star b)(i) = a(i)b(i) = a_ib_i$.
\end{ex} 

\begin{ex}[Base change]\label{convolbasechange}
If $E$ is a finite algebra, then $E^\star$ is a finite coalgebra and we have $(E^\star)^\star=E$ (see example \ref{dualfinitealgex}). 
Then, for any algebra $B$, the canonical map
$$
E\otimes B = [E^\star,\FF] \otimes [\FF,B] \to  [E^\star,B].
$$
given by the strength of $\otimes$ in $\dgVect$ is an algebra isomorphism.
It is clear that the map is a linear isomorphism, we need only to check that it is an algebra map.
This can be proven by an explicit computation left to the reader.
A more conceptual proof is to say that the strength of $\otimes$ is also the the lax structure of $[-,-]$.
\end{ex} 

\begin{ex}[Formal power series]\label{convolutionformalpowerseries}\label{dividedformalpower}
Let  $T^c(x)=\FF[x]$ be the tensor coalgebra on one generator $x$ of degree $-d$.
By definition $\Delta(x^n)=\sum_{i=0}^n x^i\otimes x^{n-i}$. 
Let $A[[t]]$ be the ring of formal power series in a variable $t$ of degree $d$.
The map $i:[T^c(x),A]\to A[[t]]$ defined by putting 
$$
i(\phi)=\sum_n \phi(x^n)t^n
$$
is an isomorphism of dg-algebras.

It is clearly an isomorphism of dg-vector spaces. Let us prove that it is an algebra map.
For $\phi, \psi\in [T^c(x),A]$, then we have
$$
i(\phi\star\psi) = \sum_n (\phi\star\psi)(x^n)t^n = \sum_n \sum_{i=0}^n \phi(x^i)\psi(x^{n-i})t^n = i(\phi)i(\psi).
$$

\end{ex}

\begin{ex}[Formal divided power series]\label{convolutionformaldividedpowerseries}
Let $T^{csh}(c)=\FF[x]$ be the {\it coshuffle coalgebra} on one generator $x$ of degree $0$.
By definition, $\Delta(x^n)=\sum_{i=0}^n {n\choose i} x^i\otimes x^{n-i}$. 
Let $A\{\!\{t\}\!\}$ be the dg-algebra of formal divided power series of example \ref{dividedpower}, then
the map $i:[T^{csh}(x),A]\to A\{\!\{t\}\!\}$ defined by putting 
$$
i(\phi)=\sum \phi(x^n)\Dfrac{t^n}{n!}
$$
is an isomorphism of dg-algebras. The proof is the same as in the example above.
\end{ex}

\medskip

\begin{prop} \label{smallconvolution}
If $C$ is a finite coalgebra and $A$ is an algebra,  then the canonical map $C^\star \otimes A\to [C,A]$
is an isomorphism of algebras.
 \end{prop}
 
\begin{proof}
By a straightforward calculation.
\end{proof}

\medskip

\begin{lemma}\label{dualofmapsofcoalgebras} If $C$ and $D$ are coalgebras,
then a linear map $g:C\to D$ is a map of coalgebras
if and only if $g^\star:D^\star \to C^\star$
is a map of algebras.
\end{lemma}

\begin{proof} If $g^\star$ is an algebra map, let us show
that $g$ is a coalgebra map. We have $g^\star(\epsilon)=\epsilon$,
since $g^\star$ preserves units, Thus, $\epsilon g=\epsilon$
and this shows that $g$ preserves the counits.
Let us show that $g$ preserves the coproducts.
We have to show that $\Delta g(x)=g(x^{(1)})\otimes g(x^{(2)})$
for every $x\in C$.
But for this it suffices to prove the equality
$$(\phi\otimes \psi)\Delta g(x)=(\phi\otimes \psi)\bigl(g(x^{(1)})\otimes g(x^{(2)})\bigr)$$
for every pair $(\phi,\psi)$ of linear forms in $D^\star$
(we are implicitly using the fact that the elements of $D^\star \otimes D^\star$
are separating the  elements of $D\otimes D$). We have $ g^\star(\phi \star \psi)=g^\star(\phi) \star g^\star(\psi)$,
since $g^\star$ preserves the convolution product.
Thus,
\begin{eqnarray*}
(\phi\otimes \psi)\Delta g(x) &=& (\phi \star \psi)(g(x))\\
 &=& g^\star(\phi \star \psi)(x)\\
  &=& \bigl(g^\star(\phi) \star g^\star(\psi)\bigr)(x)\\
    &=& g^\star(\phi)(x^{(1)}) g^\star(\psi)(x^{(2)}) (-1)^{|x^{(1)}||\psi|}\\
     &=& \phi (g(x^{(1)})) \psi (g(x^{(2)})) (-1)^{|x^{(1)}||\psi|}\\
  &=& (\phi\otimes \psi) \bigl(g(x^{(1)})\otimes g(x^{(2)})\bigr) \\
\end{eqnarray*}
\end{proof}

\medskip
We shall see in section \ref{sweedlerduality} that the functor $(-)^\star:\dgCoalg^{op}\to \dgAlg$ has a left adjoint. Meanwhile we have the following easy result.

\begin{prop}\label{dualfinitealg} 
The duality functor $(-)^\star:\dgVect^{op} \to \dgVect$ induces a contravariant equivalence between
\begin{enumerate}
\item the category of finite dg-coalgebras $\dgCoalg^{\sf fin}$ and the category of finite dg-algebras $\dgAlg^{\sf fin}$.
\item and also between the category of graded finite dg-coalgebras bounded below $\dgCoalg^{\sf gr.fin, b}$ (resp, bounded above $\dgCoalg^{\sf gr.fin, a}$) and the category of graded finite dg-algebras bounded above $\dgAlg^{\sf gr.fin, a}$ (resp, bounded below $\dgAlg^{\sf gr.fin, b}$) 
\end{enumerate}
\end{prop}
The result is also true for non-unital or non-differential (co)algebras.
 
\begin{proof}
The contravariant functor $(-)^\star:\dgVect \to \dgVect$ induces a contravariant self-equivalence of the monoidal category of finite dg-vector space $\dgVectfin$.
It thus induces a contravariant equivalence between the category of comonoid objects in $\dgVectfin$ and the category of monoid objects in $\dgVectfin$. This proves the first statement.
The second statement is proved similarly by using lemma \ref{lemmagradedfinite}.
\end{proof}

\begin{ex} \label{deltaepsilon}
Let $\FF\varepsilon$ be the dg-vector space generated by a variable of degree $d$ and let $\FF\delta=(\FF\varepsilon)^\star$ be the dual dg-vector space generated by the dual variable $\delta =\varepsilon^\star$. 
Recall the square zero dg-algebra $T_1(\varepsilon)$ and the primitive dg-coalgebra $\FF\delta_+ = T^c_1(\delta)$ then we have isomorphisms of dg-algebras $[T^c_1(\delta),\FF] = T_1(\varepsilon)$ and an isomorphism of dg-coalgebras $[T_1(\varepsilon),\FF] = T^c_1(\delta)$.
\end{ex}

\begin{prop}\label{dualtcvtvprop}
If a dg-vector space $X$ is graded finite and strictly positive or strictly negative, then there exists
\begin{enumerate}
\item an isomorphism of algebras between $T^c(X)^\star$ and $T(X^\star)$
\item and an isomorphism of coalgebras between $T(X)^\star$ and $T^c(X^\star)$.
\end{enumerate}
\end{prop}

\begin{proof} Let us prove the first statement. The canonical map $[X,\FF]^{\otimes n}\to  [X^{\otimes n},\FF] $
is an isomorphism for every $n\geq 0$ by lemma  \ref{lemmagradedfinite},
since $X$ is graded finite and bounded below. 
It follows 
that we have 
$$T(X)^\star =[T(X),\FF] =  \prod_{n\geq 0}  [X^{\otimes n} ,\FF]=\prod_{n\geq 0}  [X,\FF]^{\otimes n},$$
since the contravariant functor $(-)^\star$
takes direct sums to products.
But the family of graded vector spaces $(X^{\otimes n}:n\in \mathbb{Z})$
is locally finite, since the vector space $X$ strictly positive.
Hence also the family $([X,\FF]^{\otimes n}:n\in \mathbb{Z})$.
It then follows by lemma  \ref{convergence}
that the canonical map
 $$\bigoplus_{n\geq 0}  [X,\FF]^{\otimes n} \to \prod_{n\geq 0}  [X,\FF]^{\otimes n}$$
is an isomorphism.
This shows that the canonical map $T(X^\star)\to T(X)^\star $ is an isomorphism of vector spaces.
It is easy to verify that it is an isomorphism of algebras
$T(X^\star)\to  T^c(X)^\star$. The second statement
follows by duality if we use proposition \ref{dualfinitealg}.
\end{proof}

\medskip
The following lemma is the heart of the whole enrichment of $\dgAlg$ over $\dgCoalg$.

\begin{lemma}\label{assoconvolution}
If $C$ and $D$ are coalgebras and $A$ is an algebra, then the canonical isomorphism
$$
\lambda^{2}:[C\otimes D,A]\simeq [C,[D,A]]
$$
is a map of algebras.
\end{lemma}

\begin{proof}
A straightforward computation.
\end{proof}

\begin{thm} \label{actionconvolutionproduct}
The functor $[-,-]:\dgCoalg^{op}\times \dgAlg\to\dgAlg$ equips $\dgAlg$ with the structure of op-module over the monoidal category $\dgCoalg$ (see Appendix \ref{Vopmodule}).
The associativity constraint is given by 
$$
\lambda^{2}_{C,D;A}:[C\otimes D,A]\simeq [C,[D,A]]
$$ 
and the unit constraint by the isomorphism $[\FF,A]\simeq A$. 
\end{thm}

\begin{proof}

We need to prove that, for $C,D,E$ three coalgebras and $A$ an algebra, we want to verify the pentagon condition:
$$\xymatrix{
[(C\otimes  D)\otimes E, A] \ar[rr]^-{\simeq} \ar[d]_{\lambda^{2}_{C\otimes D,E;A}}&&[C\otimes  (D\otimes E), A] \ar[rr]^-{\lambda^{2}_{C,D\otimes E;A}}&& [C,[D\otimes E, A]\ar[d]^{[C, \lambda^{2}_{D,E;A}]}\\
[C\otimes  D, [E, A]] \ar[rrrr]^-{\lambda^{2}_{C,D;[E, A]}} &&&& [C, [D,[E, A]]]
}$$
But this condition is true in $\dgVect$ where it is a consequence of the pentagon equation for the graded tensor product.
Then the result follows from lemma \ref{assoconvolution}.
The proof of the unit condition is similar.
\end{proof}

In reference to this op-module structure, we shall call the convolution functor
$$\xymatrix{
[-,-]:\dgCoalg^{op}\times \dgAlg\ar[r]&\dgAlg
}$$
the {\em convolution product}.

\subsubsection{Convolution and (co)derivations}\label{convolderivation}

As in propositions \ref{actiondertensor} and \ref{actioncodertensor}, we establish the compatibility of the convolution product with (co)derivations.

\medskip
Recall from proposition \ref{Liebimodules2} that the formula for the Lie algebra map $\pi: gl(Y) \times  gl(X)  \to gl([X,Y])$
induced by the $(End(Y),End(X)^o)$-bimodule structure of $[X,Y]$ is given by $\pi(g,f)=[X,g]-[f,Y]$.

\begin{prop}\label{derconvolution3}
Let $A$ be an algebra and $C$ be a coalgebra. 
If $d_1:A\rhup A$ is a derivation 
and $d_2:C\rhup C$ is a coderivation, 
both of degree $n$, 
then the morphism $\varpi(d_1,d_2)=[C,d_1]-[d_2,A]$
is a derivation of degree $n$ of the convolution 
algebra $[C,A]$. 
Moreover, the map
$$\varpi:\Der(A) \times \Coder(C) \to \Der([C,A])$$
so defined is a homomorphism of Lie algebras.
The homomorphism preserves the square of odd elements.
\end{prop}

\begin{proof} If $d:A\rhup A$ is a derivation of degree
$n$, let us show that the morphism $[C,d]$ is a derivation of degree $n$.
For every $f,g\in [C,A]$, we have
\begin{eqnarray*}
[C,d](f\star g)&=&d\mu (f\otimes g)\Delta \\
&=&\mu(d \otimes C+C\otimes d) (f\otimes g)\Delta \\
&=& \mu (df\otimes g) \Delta +\mu (f\otimes dg) \Delta  (-1)^{n|f|}   \\
&=&df\star g+ f\star dg (-1)^{n|f|}   \\
&=&[C,d](f)\star g+f\star [C,d](g) (-1)^{n|f|}.  \\
\end{eqnarray*}
This proves that $[C,d]$ is a derivation of degree $n$.  
If $d:C\rhup C$ is a coderivation of degree $n$, let us show that $[d,A]$ is a derivation of degree $n$. 
For every $f,g\in [C,A]$, we have
\begin{eqnarray*}
[d,A](f\star g)&=&(f\star g)d(-1)^{n|f|+n|g|} \\
 &=&\mu (f\otimes g)\Delta d(-1)^{n|f|+n|g|} \\
&=&\mu (f\otimes g) (d \otimes C)\Delta (-1)^{n|f|+n|g|} +\mu (f\otimes g) (C\otimes d)\Delta (-1)^{n|f|+n|g|} \\
&=& \mu (fd\otimes g) \Delta (-1)^{n|f|} +\mu (f\otimes gd) \Delta  (-1)^{n|f|+n|g|}   \\
&=& \mu( [d,A](f)\otimes g)\Delta  +\mu (f\otimes  [d,A](g)) \Delta (-1)^{n|f|}   \\
&=& ([d,A](f)\star g)+(f\star [d,A](g)) (-1)^{n|f|}   \\
\end{eqnarray*}
This proves that $[d,A]$ is a derivation of degree $n$.  
Thus, id $d_1:A\rhup A$ is a derivation of degree $n$
and $d_2:C\rhup C$ is a coderivation of degree $n$, then the morphism $\varpi(d_1,d_2)=[C,d_1]-[d_2,A]$
is a derivation of degree $n$.
It then follows from proposition \ref{Liebimodules2} that the map
$\varpi :\Der(A) \times \Coder(C) \to \Der([C,A])$
is a homomorphism of Lie algebras preserving
the square of odd pairs.
\end{proof}

\subsection{The measuring functor}\label{measuring}

Recall that the transformation $\lambda^1:[X\otimes Y, Z]\to [Y,[X,Z]]$ is
taking a map $f:X\otimes Y \to Z$ in $\dgVect$ to the map $\lambda^1(f):Y \to [X,Z]$
defined by putting $\lambda^1(f)(y)(x)=f(x\otimes y)(-1)^{|x||y|}$,
and that the transformation $\lambda^2:[X\otimes Y, Z]\to [X,[Y,Z]]$ is
taking $f$ to the map $\lambda^2(f):X \to [Y,Z]$
defined by putting $\lambda^2(f)(x)(y)=f(x\otimes y)$.

\begin{defi}\label{defmeasuring}
Let $C$ be a coalgebra and $A$, $B$ be algebras.
We shall say that a linear map $f:C\otimes A\to B$ is a {\em measuring} if $\lambda^1(f):A\to [C,B]$ is an algebra map.
\end{defi}

\medskip

Let $C=(C,\Delta,\epsilon)$ be a coalgebra and let $A=(A,m_A,e_A)$ and $B=(B,m_B,e_B)$ be algebras.
Then a map $f:C\otimes A\to B$ is a measuring if and only if the following two diagrams commute
$$
\vcenter{
\xymatrix{
C\otimes A\otimes A \ar[d]_{\Delta \otimes A\otimes A} \ar[r]^-{C\otimes m_A}	& C\otimes A \ar[ddd]^f\\
C\otimes C\otimes A \otimes A \ar[d]_{can}			\\
C\otimes A\otimes C\otimes A \ar[d]_-{f\otimes f}						\\
B \otimes B\ar[r]^-{m_B} 									& B
}}
\et
\vcenter{
\xymatrix{
C \ar[d]_{\epsilon} \ar[rr]^-{C\otimes e_A} && C\otimes A \ar[d]^f\\
\FF  \ar[rr]^-{e_B}&  & B
}}.$$

In other words, a linear map $f:C\otimes A\to B$ is a measuring if and only if the following conditions are satisfied  
for every $a,b\in A$ and $c\in C$,
$$
f(c,ab)= f(c^{(1)},a)f(c^{(2)},b) \ (-1)^{|a||c^{(2)}|}
\et
f(c,e_A)=\epsilon(c)e_B.
$$
Of course the map $f:C\otimes A\to B$ is implicitely a map of complexes, this add a third condition
$$\xymatrix{
A\otimes C \ar@^{>}[rrr]^-{d_A\otimes C+A\otimes d_C} \ar[d]_{f} &&& A\otimes C \ar[d]^-{f} \\
B\ar@^{>}[rrr]^-{d_B} &&& B.
}$$
or with elements
$$
df(c,a) = f(dc,a)+f(c,da) (-1)^{|c|}.
$$

\begin{rem}
At this point we would like to advertise that the notion of measuring can be adapted to many other contexts (non-differential, non-unital, non graded, with an action of a Hopf algebra...) as soon as the convolution algebra make sense, there is a corresponding notion of measuring.
For non-differential (co)algebras, it suffices to remove the condition of compatibility with the differential.
For non-unital (co)algebras, it suffices to remove the unit condition.
For pointed (co)algebras, $f$ must satisfy the extra equation $\epsilon_B(f(c,a)) = \epsilon_C(c) \epsilon_A(a)$.
We shall come back to these variations in chapter \ref{Sweedlercontexts}.
\end{rem}

\begin{ex}
A measuring $\FF\otimes A\to B$ is just an algebra map.
\end{ex}

\begin{ex}
If $C$ is a coalgebra and $A$ is an algebra, then the reversed evaluation $rev: C\otimes [C,A]\to A$ is a measuring.
In particular, the reversed evaluation $ev:C\otimes C^\star \to \FF$ is a measuring.
\end{ex}

It is sometimes convenient to say that a map $f:A\otimes C\to B$ is a {\em right measuring} if $\lambda^2(f):A\to [C,B]$ is an algebra map. By opposition the previous notion of measuring can be called {\em left measuring}.
A map $f:A\otimes C\to B$ is a right measuring iff the composite $f\sigma:C\otimes A\simeq A\otimes C\to B$ is a left measuring. 
In the previous example, the evaluation $ev:[C,A]\otimes A\to A$ is a right measuring.
We shall seldom use the notion of right measuring.

\medskip

The following result is obvious:

\begin{lemma}\label{natcomes}\label{natmes}
Let $f:A'\to A$ and $g:B\to B'$ be two map of algebras, and $v:C'\to C$ be a map of coalgebras. 
 If $u:C\otimes A\to B$ is a measuring, then so is the composite
$$
\xymatrix{
C'\otimes A'\ar[rr]^-{v\otimes f}&& C\otimes A\ar[r]^-u &B \ar[r]^g &B'.
}$$
\end{lemma}

\medskip

\begin{defi}
We will denote the set of measurings $C\otimes A\to B$ by $\cM(C,A;B)$.
The lemma shows that these sets define a functor of three variables
$$\xymatrix{
\cM(-,-;-):\dgCoalg^{op}\times \dgAlg^{op}\times \dgAlg\ar[r]& \Set.
}$$
We shall call this functor the {\em measuring functor}. It will be the main protagonist of the enrichment of $\dgAlg$ over $\dgCoalg$.
\end{defi}

This functor is a triality in the sense of appendix \ref{triality}. By definition of measuring we have a natural isomorphism
$$
\cM(C,A;B) \simeq \dgAlg(A,[C,B])
$$
which say that the functor $\cM$ is representable in the second variable. We are going to prove that this functor is n fact representable in all its variables and that is is equivalent to an enrichment of $\dgAlg$ over $\dgCoalg$ which is tensored and cotensored:
the enrichment will be given by the operation representing $\cM$ in the first variable and called the {\em Sweedler hom},
the tensor product is the operation representing $\cM$ in the third variable and is called the {\em Sweedler product}
and the cotensor is simply given by the convolution product.

\subsection{Sweedler product}\label{sweedlerproduct}

Let $C$ be a (dg-)coalgebra and $A$ be a (dg-)algebra.
We shall say that a measuring $u:C\otimes A\to E$ is {\it universal} if the pair $(E,u)$ is representing the functor 
$$\xymatrix{
\cM(C,A;-): \dgAlg \ar[r]& \Set.
}$$
The universality means that for any algebra $B$ and any measuring $f:C\otimes A\to B$ there exists a unique map of algebras 
$g:E\to B$ such that $gu=f$. 
$$
\xymatrix{
  C\otimes A \ar[d]_u \ar[rr]^-{f} && B \\
		 E \ar@{-->}[urr]_-{g}  &&
}$$
The codomain of a universal measuring $u:C\otimes A\to E$ is well defined up to a unique isomorphism.
We shall put $C\rhd A=E$ and write $c\rhd a:=u(c\otimes a)$ for $c\in C$ and $a\in A$.
We shall say that the algebra $C\rhd A$ is  the (left) {\it Sweedler product} of the algebra $A$ by the coalgebra $C$. 
For any algebra $B$ and any measuring $f:C\otimes A\to B$
there exists a unique map of algebras $g:C\rhd A\to B$ such that $f(x\otimes y)=g(x\rhd y)$.

\begin{thm} \label{Sproduct}
The Sweedler product $C\rhd A$ exists for any algebra $A$ and any coalgebra $C$.
\end{thm}

\begin{proof} The algebra $C\rhd A$ is generated by symbols,
$c\rhd a$ for $c\in C$ and $a\in A$  on which the following 
relations are imposed,
\begin{itemize}
\item[(d)] the map $(c,a)\mapsto c\rhd a$ is a bilinear map of dg-vector spaces, $d(c\rhd a) = dc\rhd a + c\rhd da (-1)^{|c|}$;
\item[(m)] $c\rhd (ab) = (c^{(1)}\rhd a)(c^{(2)}\rhd b) (-1)^{|a||c^{(2)}|}$, for every $c\in C$ and $a,b\in A$;
\item[(u)] and $c\rhd 1 = \epsilon(c)$, for every $c\in C$.
\end{itemize}
In other terms, $C\rhd A$ is the quotient of the tensor dg-algebra $T(C\otimes A)$ by the relations (m) and (u).
\end{proof}  

\begin{rem}
The functor $C\rhd-$ is by definition left adjoint to the functor $[C,-]$.
It is also possible to prove theorem \ref{Sproduct} by using the continuity of $[C,-]$ together the $\omega$-presentability of $\dgAlg$, in analogy to theorems \ref{homcoalg} and \ref{Coalghomalg1}.
\end{rem}

\bigskip
The algebra $C\rhd A$ depends functorially on $C$ and $A$. We thus obtain a functor of two variables
$$\xymatrix{
\rhd : \dgCoalg\times \dgAlg\ar[r]& \dgAlg.
}$$
By definition, there is a natural bijection between
\begin{center}
\begin{tabular}{lc}
\rule[-2ex]{0pt}{4ex} the algebra maps & $C\rhd A\to B$, \\
\rule[-2ex]{0pt}{4ex} the measurings & $C\otimes A\to B$, \\
\rule[-2ex]{0pt}{4ex} and the algebra maps & $A\to [C,B]$.
\end{tabular}
\end{center}

\begin{prop} \label{adjunctionSproduct}
The functor $C\rhd(-) :\dgAlg \to \dgAlg$ is left adjoint to the functor $[C,-]$ for any coalgebra $C$.
\end{prop}

\begin{proof}
Obvious from the bijections above.
The unit of this adjunction is the map $\eta:A\to [C,C\rhd A]$ defined by putting $\eta(a)(c)=c\rhd a(-1)^{|a||c|}$.
The counit $\epsilon:C\rhd [C,B]\to B$ is defined by putting $\epsilon(c\rhd \phi)=\phi(c)(-1)^{|c||\phi|}$.
\end{proof}

\medskip

\begin{prop}\label{unitforSweedler}.
We have $C\rhd \FF = \FF$ for any coalgebra $C$, and we have $\FF\rhd A=A$ for any algebra $A$.
\end{prop}
\begin{proof} The functor $C\rhd(-) :\dgAlg \to \dgAlg$ preserves initial objects since it is
a left adjoint. It is easy to verify directly that the canonical isomorphism $\FF \otimes A\to A$ 
is a measuring and that it is universal.
\end{proof}

If $C$ is a coalgebra, then the counit $\epsilon:C\to \FF $ induces a map $\epsilon \rhd A:C \rhd A \to \FF \rhd A=A$.
Hence the algebra $C\rhd A$ is equipped with a canonical algebra map $\epsilon\rhd A:C\rhd A\to A$.
If $C$ is pointed with base point $e:\FF \to C$ then the algebra map $e\rhd A:A\to C\rhd A$ is a section of $\epsilon\rhd A$.

\begin{prop} \label{smallalgexp}
If $C$ is a finite coalgebra, then the functor $C\rhd (-) :\dgAlg\to \dgAlg$ is left adjoint to the functor $C^\star\otimes -$.
Equivalently, if $A$ is a finite algebra, the functor $A^\star\rhd (-) :\dgAlg\to \dgAlg$ is left adjoint to the functor $A\otimes -$.
 \end{prop}
\begin{proof}
The functor $C\rhd(-):\dgAlg\to \dgAlg$ is left adjoint to the functor $[C, -]$.
But the canonical map $C^\star\otimes B\to [C,B]$
is an isomorphism for any algebra $B$, since $C$ is finite.
 \end{proof}

\bigskip

We present here only a two examples of Sweedler products,
more are given in section \ref{examplesweedlerproduct}, including applications to matrix and jet algebras.

\begin{ex}[coproduct]
If  $\FF I$ is the diagonal coalgebra of a set $I$, then the algebra $\FF I\rhd A$ is isomorphic to the 
coproduct $I\cdot  A$ of $I$-copies of the algebra $A$. 
This can be seen by the explicit computation of $\FF I\rhd A$: it is generated by symbols $i\rhd a$ for each $i\in I$ and $a\in A$ and by relations $i\rhd ab=(i\rhd a)(i\rhd b)$. In particular any $i$ defines an algebra embedding $A\to \FF I\rhd A$.
It is clear that $\FF I\rhd A$ is generated freely from these embeddings.
\end{ex}

Recall that, for $A$ an algebra and $M$ an $A$-bimodule, the {\em tensor algebra of $M$ over $A$} is 
$$
T_A(M) = A\oplus M\oplus (M\otimes_AM) \oplus (M\otimes_AM\otimes_AM ) \oplus \dots
$$
The following result is the analog for algebra of proposition \ref{coderham} for coalgebras.

\begin{prop}\label{derham}
For $\delta$ a graded variable of degree $n$, there is an isomorphism of algebras $\FF[\delta]\rhd A\simeq T_A(S^n\Omega_A)$ where $\Omega_A$ is the bimodule of differentials of $A$.
\end{prop}
\begin{proof}
For the purpose of this proof, it is convenient to see $\Omega_A$ generated as a $A$-bimodule by symbols $dx$ for every $x\in A$ and the relations
\begin{eqnarray*}
&d1=0\\
&d(xy)=(dx)y +x(dy),
\end{eqnarray*}
the extra relation $d1=0$ is harmless as it is implied by $d(xy) = (dx)y+ xdy$.

In consequence, the tensor algebra $T_A(S^n\Omega_A)$ of $S^n\Omega_A$ over $A$ is generated as an algebra over $\FF$ by elements 
$x$ and $s^ndx$ for every $x\in A$ and the relations
\begin{eqnarray*}
&1.x = x\\
&s^nd1=0\\
&x.y = xy\\ 
&s^nd(xy)=(s^ndx) y +x(s^ndy) \ (-1)^{n|x|}.
\end{eqnarray*}
On the other side, the proof of theorem \ref{Sproduct} constructs $\FF[\delta]\rhd A$ as generated by symbols $1\rhd x$ and $\delta \rhd x$ and relations:
\begin{eqnarray*}
&1\rhd 1 = 1,\\
&\delta\rhd 1 = 0,\\
&1\rhd (xy)=(1\rhd x)(1\rhd y) \\
&\delta\rhd (xy)=(\delta \rhd x)(1\rhd y) + (1\rhd x)(\delta \rhd y)\ (-1)^{n|x|}.
\end{eqnarray*}
The identification $1\rhd x = x$ and $\delta\rhd x = s^n(dx)$ provides an obvious isomorphism.
\end{proof}

We shall call the algebra $T_A(S^n\Omega_A)$ the {\em differential algebra} of $A$.
In section \ref{derhamexampleadj} it will be proved to be a left adjoint functor.

\bigskip

\begin{prop}\label{Sweedleroftensor}
Let $C$ be a coalgebra, $X$ is a (dg-)vector space and $A$ be an algebra.
Then every linear map $f:C\otimes X\to A$ can be extended uniquely as a measuring $f' :C\otimes T(X)\to A$.
Moreover, if $i$ is the inclusion $C\otimes X \to T(C\otimes X)$, then the measuring $i':C\otimes T(X)\to T(C\otimes X)$ is universal.
Hence we have 
$$
C\rhd T(X)= T(C\otimes X).
$$
\end{prop}

\begin{proof} Let us prove the first statement. The linear map $f:C\otimes X\to A$
correponds to a linear map $g:X\to [C,A]$ which
can be extended uniquely as an algebra map $g':T(X)\to [C,A]$
which correponds to a measuring $f' :C\otimes T(X)\to A$.
Let us now show that the measuring $i':C\otimes T(X)\to T(C\otimes X)$ is universal.
If $A$ is an algebra, we have defined a  chain of natural bijections between 
\begin{center}
\begin{tabular}{lc}
\rule[-2ex]{0pt}{4ex} the algebra maps & $h:T(C\otimes X)\to A$,\\
\rule[-2ex]{0pt}{4ex} the linear maps & $C\otimes X\to A$,\\
\rule[-2ex]{0pt}{4ex} the linear maps & $X\to [C,A]$,\\
\rule[-2ex]{0pt}{4ex} the algebra maps & $T(X)\to [C,A]$,\\
\rule[-2ex]{0pt}{4ex} the left measurings & $k:C\otimes T(X) \to A$.
\end{tabular}
\end{center}
The bijections show that the functor $\cM(C,T(X);-)$ is represented by the algebra $T(C\otimes X)$.
Moreover, if $A=T(C\otimes X)$ and $h$ is the identity, then $k=i'$.
This proves that the left measuring $i':C\otimes T(X)\to T(C\otimes X)$ is universal,
and hence that we have $C\rhd T(X)= T(C\otimes X)$.
\end{proof}

\medskip

We study now the associative structure of the Sweedler product.

\begin{lemma}\label{compositionmeasuring}
Let $A$, $B$ and $E$ be algebras and let $C$ and $D$ be coalgebras.
If $f:C\otimes A\to B$ and $g:D\otimes B\to E$ are measurings, then so is the composite
$$
\xymatrix{
D\otimes C \otimes A \ar[rr]^-{D\otimes f} && D\otimes B \ar[r]^-g& E.
}
$$
\end{lemma}
This operation will be refered as the {\em composition} of measurings. 

\begin{proof} This follows from lemma \ref{assoconvolution} or by a straightforward computation.
\end{proof}

As a consequence of this lemma, we can compose the universal measurings $u_1:C\otimes A\to C\rhd A$ 
with the universal measuring
$u_2:D\otimes (C\rhd A) \to D\rhd (C \rhd A)$ to obtain a measuring $u_3=u_2(D\otimes u_1)$,
$$
\xymatrix{
D\otimes C \otimes A \ar[rr]^-{D\otimes u_1} && D\otimes(C\rhd A)  \ar[r]^-{u_2}& D\rhd(C\rhd A)
}
$$

By definition, $u_3(x,y,z)=x\rhd (y \rhd z)$.
There is then a unique map of algebras 
$$\alpha:(D\otimes C)\rhd A \to D\rhd (C\rhd A)$$
such that 
$\alpha ((x\otimes y)\rhd z ) = (x\rhd y)\rhd z$.

\begin{prop}\label{assocSproduct2} 
The measuring $u_3$ is universal and
the map $\alpha:(D\otimes C)\rhd A\to D\rhd(C\rhd A)$ is an isomorphism.
\end{prop}

\begin{proof} If $B$ is an algebra, 
let us show that the map $f\mapsto fu_3$ is a bijection between
the algebra maps $D\rhd(C\rhd A)$ and the measurings $(D\otimes C)\otimes A  \to B$.
For this it suffices to verify that the map $f\mapsto fu_3$ is the composite of
the following sequence of bijections
\begin{center}
\begin{tabular}{lc}
\rule[-2ex]{0pt}{4ex} algebra maps & $D\rhd(C\rhd A)\to B$.\\
\rule[-2ex]{0pt}{4ex} algebra maps & $C\rhd A \to [D,B]$,\\
\rule[-2ex]{0pt}{4ex} algebra maps & $A \to [C,[D,B]]$,\\
\rule[-2ex]{0pt}{4ex}  algebra maps  & $A \to [C\otimes D,B]$,\\
\rule[-2ex]{0pt}{4ex} measurings & $(C\otimes D)\otimes A\to B$,
\end{tabular}
\end{center}
where the first and second bijections depends on the adjunction in proposition \ref{adjunctionSproduct},
where the third bijection depends on lemma \ref{assoconvolution},
and where the last bijection depends on definition \ref{defmeasuring}.
The universality of $u_3$ is proved.
It follows that $\alpha$ is an isomorphism.
\end{proof}

Recall that a measuring $\FF \otimes A\to B$ is just a map of algebras, thus the map $l:A\to \FF \rhd A$ defined by putting $l(x)=1\rhd x$ is an isomorphism.

\begin{thm} \label{actionsweedlerproduct}
The functor $\rhd:\dgCoalg\times \dgAlg\to\dgAlg$ equips $\dgAlg$ with the structure of module over the monoidal category $\dgCoalg$ (see Appendix \ref{Vmodule}).
The associativity constraint is given by 
$$
\alpha = \alpha_{C,D;A}:(C\otimes D)\rhd A \simeq C\rhd (D\rhd A)
$$ 
and the unit constraint $\FF \rhd A \simeq A$ is the map $1\rhd x\mapsto x$. 
\end{thm}

\begin{proof}
The proof of the associativity of the enrichment over $\dgCoalg$ is formal consequence of theorem \ref{actionconvolutionproduct} and of the definition via the measuring functor: if any operation representing the measuring functor is associative so are the other ones.

We shall nonetheless give a proof using the notion of universal measuring.
For $C,D,E$ three coalgebras and $A$ an algebra, we want to verify the pentagon condition:
$$\xymatrix{
((C\otimes  D)\otimes E)\rhd A \ar[rr]^-{\simeq} \ar[d]_{\alpha_{C\otimes D,E;A}}&&(C\otimes  (D\otimes E))\rhd A \ar[rr]^-{\alpha_{C,D\otimes E;A}}&& C\rhd ((D\otimes E)\rhd A)\ar[d]^{C\rhd \alpha_{D,E;A}}\\
(C\otimes  D)\rhd (E\rhd A) \ar[rrrr]^-{\alpha_{C,D;E\rhd A}} &&&& C\rhd  (D\rhd (E\rhd A))
}$$
By property of universal measuring, it is enough to check the commutativity of the diagram on the generators $(c\otimes d\otimes e)\rhd a$ of $(C\otimes  D\otimes E)\rhd A$. Then the commutativity is obvious.
The proof of the unit condition is similar.
\end{proof}

\subsection{Sweedler Hom and comeasurings}\label{sweedlerhom}

Let $A$ and $B$ be two (dg-)algebras.
We shall say that a measuring $u:E\otimes A\to B$ is {\it couniversal} if the pair $(E,u)$ represents the functor 
$$
\cM(-,A;B): \dgCoalg^{op}\to \Set.
$$
The coalgebra $E$ of a couniversal measuring $v:E\otimes A\to B$ is well defined up to a unique isomorphism 
and we shall denote it by $\{A,B\}$. We shall denote the couniversal measuring as a strong evaluation $\mathbf{ev}:\{A,B\}\otimes A \to B$.
By definition, for any coalgebra $C$ and any measuring $f:C\otimes A\to B$, there is a unique coalgebra map 
$g:C\to \{A,B\}$ such that ${\bf ev}(g\otimes A)=f$.
$$
\xymatrix{
&& \{A,B\}\otimes A	\ar[d]^-{\mathbf{ev}} \\
C\otimes A	\ar@{-->}[urr]^-{g\otimes A} \ar[rr]^-{f} && B. 
}$$
Couniversal measurings are constructed by Sweedler in \cite{Sw}. Our approach will use the dual notion of comeasuring. This notion is analogous to that of comorphism of coalgebras (definition \ref{remreductioncoalg}).

\begin{defi}\label{defcomeasuring}
Let $A=(A,m_A,e_A)$ and $B=(B,m_B,e_B)$ be algebras and $C=(C,\Delta,\epsilon)$ be a coalgebra.
We say that a linear map $g:C\to [A,B]$ is a  {\em comeasuring}  if the map $ev(g\otimes A):C\otimes A\to B$ is measuring.
\end{defi}

Hence a map $g:C\to [A,B]$ is a comeasuring if and only if we have
$$
g(c)(ab)= g(c^{(1)})(a)g(c^{(2)})(b) \ (-1)^{|a||c^{(2)}|}
\et
k(c)(1)=\epsilon(c)1
$$
for every $a,b\in A$ and $c\in C$. These conditions also means that the 
 following two diagrams commute
$$
\vcenter{
\xymatrix{
C \ar[d]_\Delta \ar[rr]^k				&& [A,B] \ar[ddd]^{[m_A,B]}\\
C\otimes C \ar[d]_{k\otimes k}			\\
[A,B]\otimes [A,B] \ar[d]_{can}			\\
[A\otimes A,B\otimes B]\ar[rr]_-{[A,m_B]}	&& [A\otimes A,B]
}}
\et
\vcenter{
\xymatrix{
C\ar[d]_{\epsilon} \ar[rr]^-k && [A,B] \ar[d]^{[e_A,B]}\\
\FF  \ar[rr]^-{e_B}&  & B.
}}$$

There are canonical bijections between
\begin{center}
\begin{tabular}{lc}
\rule[-2ex]{0pt}{4ex} algebra maps & $A\to [C,B]$, \\
\rule[-2ex]{0pt}{4ex} measurings & $C\otimes A\to B$, \\
\rule[-2ex]{0pt}{4ex} and comeasurings & $C\to [A,B]$.
\end{tabular}
\end{center}

We shall say that a comeasuring  $k:E\to [A,B]$ is {\it couniversal} if the corresponding measuring $ev( k\otimes A):E\otimes A\to B$
is couniversal. The couniversality of $k$
means concretely that for any coalgebra $C$ and any comeasuring $g:C\to [A,B]$ there exists a unique coalgebra map 
$f:C\to E$ such that $g=k f$.
$$
\xymatrix{
	&& E\ar[d]^k \\
C\ar[rr]_-{g}\ar@{-->}[urr]^f && [A,B] 
}$$

The domain of a couniversal comeasuring $u:E\to [A,B]$ is well defined up to a unique isomorphism and shall be denoted it $\{A,B\}$
The couniversal comeasuring shall be noted as $\Psi:\{A,B\}\to [A,B]$.

\medskip

We leave to the reader to check that the two definitions of $\{A,B\}$ agree.
The strong evaluation $\bf ev$ can be defined from $\Psi$ by 
the following triangle
$$\xymatrix{
\{A,B\}\otimes A\ar[d]_{\Psi \otimes A}	 \ar[rr]^-{\mathbf{ev}} && B \\
[A,B]\otimes A. \ar[rru]_{ev} && 
}$$
Reciprocally $\Psi$ can be constructed from $\bf ev$ by 
$$\xymatrix{
\{A,B\}\ar[rrd]_{\Psi}\ar[rr]^-{coev} && [A,\{A,B\}\otimes A]\ar[d]^{[A,\mathbf{ev}]} \\
&& [A,B]
}$$
where the map $coev$ is the unit of the adjunction $-\otimes A \dashv [A,-]$.

\medskip

The coalgebra $\{A,B\}$ is the {\it measuring coalgebra} $M(A,B)$ introduced by Sweedler \cite[ch. VII]{Sw}. 
We shall say that it is the {\it Sweedler hom} between the algebras $A$ and $B$. 

\medskip

With these notions, we have natural bijections between
\begin{center}
\begin{tabular}{lc}
\rule[-2ex]{0pt}{4ex} the measurings & $f:C\otimes A \to B$, \\
\rule[-2ex]{0pt}{4ex} the algebra maps & $g:C\rhd A\to B$, \\
\rule[-2ex]{0pt}{4ex} the algebra maps & $h:A \to [C,B]$, \\
\rule[-2ex]{0pt}{4ex} the comeasurings & $k:C \to [A,B]$, \\
\rule[-2ex]{0pt}{4ex} and the coalgebra maps & $l:C\to \{A,B\}$.
\end{tabular}
\end{center}
We shall put 
$\Lambda^1(g):=\lambda^1(f)=h$ and
$\Lambda^2(f)=\Lambda^2(g):=l$.

\medskip

\begin{thm}[\cite{Sw}] \label{Coalghomalg1}
There exists a couniversal comeasuring $\Psi:\{A,B\} \to [A,B]$  for any pair of algebras $A$ and $B$.
Equivalently, there exists a couniversal measuring $\mathbf{ev}: \{A,B\} \otimes A\to B$.
 \end{thm}

\begin{proof} 
By corollary \ref{catcoalgcomplete}, it suffices to show that the functor
$\cM(-,A;B): \dgCoalg^{op}\to \Set$
is continuous. But the functor $\cM( -,A;B)$
is by definition isomorphic to the functor $\dgAlg(A,[-,B])$.
The functor $\dgAlg(A,-):\dgAlg \to  \Set$ 
is continuous, since it is representable.
Hence it suffices to show that the functor 
$[-,B]:\dgCoalg^{op}\to   \dgAlg$ is continuous.
The forgetful functor $U_1:  \dgAlg\to   \dgVect$
preserves and reflects limits, since it is continuous and reflect isomorphisms.
Hence it suffices to show that the composite functor 
$U_1[-,B]:\dgCoalg^{op} \to   \dgVect$ is continuous.
But we have
$U_1[C,B]=[U_2(C),U_1(B)]$ for every
coalgebra $C$, where $U_2:  \dgAlg\to   \dgVect$
is the forgetful functor. 
The functor $[-,U_1(B)]: \dgVect^{op}\to   \dgVect$
is continuous since it is a right adjoint.
The functor $U_2$ is cocontinuous by proposition \ref{coalgcocomp}. 
Hence the composite $[U_2(-),U_1(B)]$ is continuous.
  
\end{proof}

\begin{rem}\label{remconstructionSweedlerHOM}

As in remark \ref{remconstructionHOM}, we can use the universal property of the comeasuring $\Psi$ to construct $\{A,B\}$.
For $A$ and $B$ two algebras, the diagrams
$$
\vcenter{
\xymatrix{
T^\vee([A,B]) \ar[d]_\Delta \ar[rr]				&& [A,B] \ar[ddd]^{[m_A,B]}\\
T^\vee([A,B])\otimes T^\vee([A,B]) \ar[d]_{p\otimes p}	\\
[A,B]\otimes [A,B] \ar[d]_{can}					\\
[A\otimes A,B\otimes B]\ar[rr]_-{[A,m_B]}			&& [A\otimes A,B]
}}
\et
\vcenter{
\xymatrix{
T^\vee([A,B])\ar[d]_{\epsilon} \ar[rr]^-g && [A,B] \ar[d]^{[e_A,B]}\\
\FF  \ar[rr]^-{e_B}&  & B
}}$$
where $p:T^\vee([A,B])\to [A,B]$ is the cogenerating map, do not commute, but they define a pair of parallel maps
in the category $\dgVect$,
$$\xymatrix{
T^\vee([A,B])\ar@<.6ex>[r]\ar@<-.6ex>[r]& [A\otimes A,B]
}
\et
\xymatrix{
T^\vee([A,B])\ar@<.6ex>[r]\ar@<-.6ex>[r]& B
}
$$
From these pairs, we obtain by coextension two other pairs of parallel maps in the category $\dgCoalg$,
$$\xymatrix{
v,v':T^\vee([A,B])\ar@<.6ex>[r] \ar@<-.6ex>[r]& T^\vee([A\otimes A,B])
}
\et
\xymatrix{
w,w':T^\vee([A,B])\ar@<.6ex>[r]\ar@<-.6ex>[r]& T^\vee(B)
}
$$
The coalgebra $\{A,B\}$ is the common equalizer of the pairs $(v,v')$ and $(w,w')$ in the category $\dgCoalg$ (see appendix \ref{commonequalizer}).
With this construction, the couniversal comeasuring $\Psi$ is defined to be the composite $\{A,B\}\to T^\vee([A,B]) \to [A,B]$ and
the couniversal measuring ${\bf ev}:\{A,B\}\otimes A \to B$ is obtained by putting ${\bf ev} =ev(\Psi\otimes A)$.

This construction gives another proof of the existence of $\{A,B\}$ but as remarked in \ref{remconstructionHOM}, limits of coalgebras are difficult to compute.
Corollary \ref{betterShom} gives another construction of $\{A,B\}$ using the monadicity theorem \ref{monadicalg}.
\end{rem}

Recall the notion of cogenerating and separating maps from definition \ref{coseparation}.

\begin{lemma} \label{cogeneratingcomeasuring}
The comeasuring $\Psi:\{A,B\}\to [A,B]$ is separating.
\end{lemma}

\begin{proof}
This follows from the uniqueness condition in the couniversal property of $\Psi$.
\end{proof}

\bigskip

The Sweedler hom defines a functor of two variables
$$\xymatrix{
\{-,-\}:\dgAlg^{op}\times \dgAlg \ar[r]& \dgCoalg
}$$
which we are going to prove to be an enrichment of $\dgAlg$ over $\dgCoalg$.

\begin{lemma} \label{compocomes}
If $A$, $B$ and $E$ are algebras, there is then a unique coalgebra map
$\mathbf{c}: \{B,E\}\otimes \{A,B\}\to   \{A,E\}$
such that the following square commutes
$$\xymatrix{
\{B,E\}\otimes \{A,B\} \otimes A \ar@{-->}[rr]^-{\mathbf{c}\otimes A}\ar[d]_(0.6){\{B,E\}\otimes {\bf ev}}&& \{A,E\} \otimes A \ar[d]^{{\bf ev}} \\
 \{B,E\}\otimes B  \ar[rr]^-{{\bf ev}} && E
}$$
And equivalently, such that the following square commutes
$$\xymatrix{
\{B,E\}\otimes \{A,B\}\ar@{-->}[rr]^-{\mathbf{c}}\ar[d]_{\Psi\otimes \Psi}&& \{A,E\}\ar[d]^{\Psi}\\
[B,E]\otimes [A,B]\ar[rr]^c && [A,E].
}$$
\end{lemma} 

\begin{proof}
The composite
$$\xymatrix{
\{B,E\}\otimes \{A,B\}\otimes A \ar[rr]^-{\{B,E\}\otimes {\bf ev}} &&\{B,E\}\otimes B
\ar[r]^-{{\bf ev}}& E
}$$
is a measuring by lemma \ref{compositionmeasuring}. 
 Hence there is a unique coalgebra map $\mathbf{c}: \{B,E\}\otimes \{A,B\}\to   \{A,E\}$
 such that the first square commutes. It is easy to see that the 
 second square commutes iff the first commutes, using a diagram as in lemma \ref{lemmecommutcoalg}.
 \end{proof}

\begin{lemma} \label{unitcomes}
For any algebra $A$, there is also a unique coalgebra map $e_A:\FF \to \{A, A\}$ such that $\Psi e_A=1_A$,
$$
\xymatrix{
	&& \{A,A\} \ar[d]^\Psi \\
\FF \ar[rr]_-{1_A}\ar@{-->}[urr]^-{e_A} && [A,A] 
}$$
\end{lemma}

\begin{proof}
The unit $1_A:\FF \to [A,A]$ is a comeasuring, since
the canonical isomorphism $\FF \otimes A\to A$ is a measuring.
 \end{proof}

\begin{thm}\label{enrichmentalgcoal}
The map $\mathbf{c}: \{B,E\}\otimes \{A,B\} \to  \{A,E\}$
defined above is the composition law for an enrichment of the category $\dgAlg$ over the closed monoidal  category $\dgCoalg$. 
The unit of the composition law is the map $e_A:\FF \to \{A, A\}$.
The enriched category $\dgAlg$ is bicomplete over $\dgCoalg$;
the tensor product of an algebra $A$ by a coalgebra $C$ is the algebra $C\rhd A$ and the cotensor 
of $A$ by $C$ is the convolution algebra $[C,A]$.
Hence there are natural isomorphisms of coalgebras
$$\{C\rhd A,B\} \simeq \HOM(C,\{A,B \})\simeq \{A,[C,B] \}$$
for a coalgebra $C$ and for algebras  $A$ and $B$.
\end{thm}

\begin{proof}
As for theorem \ref{actionsweedlerproduct}, the proof is formal consequence of theorem \ref{actionconvolutionproduct} and of the definition via the measuring functor: if any operation representing the measuring functor is associative so are the other ones.

We shall nonetheless detail the proof as it is an interesting manipulation of the notion of couniversal comeasuring.

Let us prove the associativity of the composition law $\mathbf{c}: \{B,E\}\otimes \{A,B\} \to  \{A,E\}$.
For this we have to show that the top face of the following cube commutes
for any quadruple of algebras $A,B,E,F$,
$$\xymatrix{
&\{E,F\}\otimes \{B,E\}\otimes \{A,B\}  \ar[dl]_-{\mathbf{c}\otimes \{A,B\}\qquad} \ar[rr]^-{\{E,F\}\otimes \mathbf{c}}\ar'[d]^{\Psi\otimes \Psi \otimes \Psi}[dd]
	&& \{E,F\}\otimes \{A,E\} \ar[dd]^{\Psi\otimes \Psi}  \ar[dl]_-{\mathbf{c}}   \\
\{B,F\}\otimes \{A,B\}\ar[rr]^(0.4){\mathbf{c}} \ar[dd]_{\Psi\otimes \Psi} 
	&&  \{A,F\} \ar[dd]^(0.3)\Psi \\
&[E,F]\otimes [B,E]\otimes [A,B]  \ar'[r][rr] \ar[dl]_-{c\otimes [A,B]\qquad}
	&& [E,F]\otimes [A,E]  \ar[dl]_c \\
[B,F]\otimes [A,B] \ar[rr]^c 
	&& [A,F]
 }$$
The vertical faces of the cube commute by lemma \ref{compocomes}, and the bottom face commutes since
composition is associative in $\dgVect$. Hence the top face commutes 
after post-composition with the map $\Psi:\{A,F\} \to [A,F]$.
It follows that the top face commutes, since $\Psi$ is separating by lemma \ref{cogeneratingcomeasuring}.
Let us now verify that the map $e_A:\FF \to \{A, A\}$ is a right 
unit for the composition law. For this we have to show that the top triangle
of the following diagram commutes,
$$\xymatrix{
\{A,B\} \ar@{=}[dr] \ar[rr]^-{\{A,B\}\otimes e_A} \ar[dd]_{\Psi}
	&& \{A,B\}\otimes \{A,A\} \ar[dd]^{\Psi\otimes \Psi}  \ar[dl]_-{\mathbf{c}}    \\
&  \{A,B\} \ar[dd]_(0.3)\Psi \\
[A,B] \ar'[r][rr]^-{[A,B]\otimes 1_A} \ar@{=}[dr] 
	&& [A,B]\otimes [A,A]  \ar[dl]_c\\
& [A,B]
}$$
The vertical faces of the diagram commute by definition of $e_A$ and by lemma \ref{compocomes} .
The bottom face commutes since the map $1_A:\FF \to [A, A]$ is a unit
for the composition law in the category $\dgVect$. 
Hence the top face commutes 
after composition with the map $\Psi:\{A,B\} \to [A,B]$.
It follows that the top face commutes, since $\Psi$ is separating by lemma \ref{cogeneratingcomeasuring}.
The proof that the map $e_A:\FF \to \{A, A\}$ is a left
unit for the composition law is similar. 

Let us now prove that the convolution algebra $[C,B]$ is the cotensor product of the algebra $B$
by the coalgebra $C$. If $A$ is an algebra and $D$ is a coalgebra, we have a chain of natural bijections between
\begin{center}
\begin{tabular}{lc}
\rule[-2ex]{0pt}{4ex} coalgebra maps & $D\to  \{A,[C,B]\}$,\\
\rule[-2ex]{0pt}{4ex} measurings & $ D\otimes A \to [C,B]$,\\
\rule[-2ex]{0pt}{4ex} algebra maps & $A \to  [D,[C,B]]$,\\
\rule[-2ex]{0pt}{4ex} algebra maps & $A \to  [D\otimes C,B]$, \\
\rule[-2ex]{0pt}{4ex} measurings & $D\otimes C \otimes A \to  B$, \\
\rule[-2ex]{0pt}{4ex} coalgebra maps & $D\otimes C \to  \{A,B\}$, \\
\rule[-2ex]{0pt}{4ex} coalgebra maps & $D\to \HOM(C, \{A,B\})$, \\
\end{tabular}
\end{center}
where the third bijection is given by lemma \ref{assoconvolution}.
This shows by Yoneda lemma that we have a natural isomorphism
$$\{A,[C,B]\}\simeq  \HOM(C, \{A,B\})$$
and hence that $[C,B]$ is the cotensor of $B$ by $C$.
It remains to prove that 
the algebra $C\rhd A$ is the tensor product of the algebra $A$
by the coalgebra $C$. 
 If $B$ is an algebra and $D$ is a coalgebra,
we have a chain of natural bijections between
\begin{center}
\begin{tabular}{lc}
\rule[-2ex]{0pt}{4ex} coalgebra maps & $D\to  \{C\rhd A,B\}$,\\
\rule[-2ex]{0pt}{4ex} measurings & $D\otimes (C\rhd A) \to B$,\\
\rule[-2ex]{0pt}{4ex} algebra maps & $D\rhd (C\rhd A) \to B$, \\
\rule[-2ex]{0pt}{4ex} algebra maps & $(D\otimes C)\rhd A \to B$, \\
\rule[-2ex]{0pt}{4ex} coalgebra maps & $D\otimes C \to  \{A,B\}$, \\
\rule[-2ex]{0pt}{4ex} coalgebra maps & $D\to \HOM(C, \{A,B\})$, \\
\end{tabular}
\end{center}
where the third bijection is given by proposition \ref{assocSproduct2}.
This shows by Yoneda lemma that we have a natural isomorphism
$$ \{C\rhd A,B\} \simeq  \HOM(C, \{A,B\})$$
and hence that $C\rhd A$ is the tensor of $A$ by $C$.
\end{proof}

In particular, theorem \ref{enrichmentalgcoal} says that, for any algebra $A$, $\{A,A\}$ has the structure of a monoid in $\dgCoalg$ and is therefore a bialgebra. Moreover the second square in lemma \ref{compocomes} ensure that the map $\Psi:\{A,A\}\to [A,A]$ is a map of algebras.

\medskip

\begin{cor} \label{adjunctiontypes}
Let $C$ be a coalgebra and $A$ be an algebra.
We have the following strong adjunctions
$$\xymatrix{
C\rhd(-):\dgAlg  \ar@<.6ex>[r]& \dgAlg:[C,-]\ar@<.6ex>[l]
}$$
$$\xymatrix{
[-,A]: \dgCoalg^{op}  \ar@<.6ex>[r]& \dgAlg:\{-,A \}\ar@<.6ex>[l]
}$$
$$\xymatrix{
(-)\rhd A: \dgCoalg  \ar@<.6ex>[r]& \dgAlg:\{A,- \}\ar@<.6ex>[l].
}$$
\end{cor}

The second adjunction expresses the fact that the contravariant functors $\{-,A \}$ and $[-,A] $ are right adjoints.
Chapter \ref{adjunctions} is dedicated to example of these three types of adjunctions.

\medskip
The notion of strong limit is given in appendix \ref{stronglimit}.
\begin{prop}\label{stronglimalg}
All ordinary (co)limits in $\dgAlg$ are strong.
\end{prop}
\begin{proof}
This is a formal consequence of $\dgAlg$ begin bicomplete.
We shall prove only the result for colimits, the proof for limits is similar.
Let $A:I\to \dgAlg$ be a diagram with colimit $B$, then $B$ can be writtent as the coequalizer in $\dgAlg$
$$\xymatrix{
\coprod_{j\to i}C_i\ar@<.6ex>[r]\ar@<-.6ex>[r] & \coprod_{i}C_i\ar[r]& B.
}$$
$B$ is a strong colimit, if for any algebra $E$, we have an equalizer in $\dgCoalg$
$$\xymatrix{
\{B,E\}\ar[r]&\prod_i\{A_i,E\}\ar@<.6ex>[r]\ar@<-.6ex>[r] & \coprod_{j\to i}\{C_i,E\}.
}$$
But this can be deduced form the previous equalizer and the fact that $\{-,E\}$ being a contravariant right adjoint to $[-,E]$, it sends coequalizers to equalizers.
\end{proof}

\bigskip

\begin{prop}\label{Shomunit}
For every coalgebra $A$, we have $\{\FF,A\}=\FF$, in particular $\{\FF,\FF\}=\FF$.
\begin{enumerate}
\item $\{\FF,A\}=\FF$, in particular $\{\FF,\FF\}=\FF$
\item and $\{A,0\}=\FF$.
\end{enumerate}
\end{prop}
\begin{proof}
The functors $\{-,A\}:\dgAlg^{op} \to \dgCoalg$ and $\{A,-\}:\dgAlg \to \dgCoalg$ are right adjoint to by corollary \ref{adjunctiontypes}, so they send the terminal object to a terminal object.
The result follows, since $\FF$ and $0$ are respectively the initial and terminal objects of $\dgAlg$ and $\FF$ is the terminal object of $\dgCoalg$.
\end{proof}

\bigskip

Recall that if $\bf C$ is a category enriched over a monoidal category $({\bf V},\otimes,\un)$, 
the {\em underlying set} of morphisms between $X$ and $Y$ in $\bf C$ is defined as $\Hom_{\bf V}(\un,{\bf C}(X,Y))$ and
the {\em underlying category} of $\bf C$ is defined as the category with the same objects as $\bf C$ but with sets of morphisms $\Hom_{\bf V}(\un,{\bf C}(X,Y))$. 
Let $\EAlg$ be the category $\dgAlg$ viewed as enriched over $\dgCoalg$.
In this case, the set of underlying elements of $\{A,B\}$ is the set of atoms of $\{A,B\}$.

\begin{lemma}\label{atomShom}
Let $A$ and $B$ be two algebras, the set of atoms of $\{A,B\}$ is in bijection with the set of algebras maps from $A$ to $B$.
\end{lemma}
\begin{proof}
By the universal property of $\{A,B\}$, algebra maps $A\to B$ are in bijection with coalgebra maps $\FF\to \{A,B\}$.
\end{proof}

\begin{prop}\label{underlyingAlg}
The underlying category of $\dgAlg^\$$ is the ordinary category $\dgAlg$.
\end{prop}
\begin{proof}
This is essentially the previous proposition. We leave to the reader the proof that the composition is the good one.
\end{proof}

\bigskip

Let $A$ be an algebra and $X$ be a vector space.
If $p: T^\vee([X,A])\to [X,A]$
is the cofree map, then the composite of the maps
$$
\xymatrix{
 T^\vee([X,A])\otimes X \ar[rr]^-{p\otimes X} && [X,A]\otimes X \ar[r]^-{ev} & A
}
$$
can be extended uniquely as a measuring $p':T^\vee([X,A])\otimes T(X)\to A$ by proposition \ref{Sweedleroftensor}.
In particular, the map $\lambda^2(p):T^\vee([X,A])\to [X,A]$ is the composition of $\lambda^2(p'):T^\vee([X,A])\to [T(X),A]$ with the 'restriction to the generators' map $[T(X),A]\to [X,A]$.

\begin{prop} \label{exSweedlerhom}
If $A$ is an algebra and $X$ is a vector space, 
then the measuring $p':T^\vee([X,A])\otimes T(X)\to A$
defined above is  couniversal.
Hence we have 
$$\{T(X),A\}= T^\vee([X,A]).$$
\end{prop}

\begin{proof} 
 If $C$ is a coalgebra, there is a chain of natural bijections between
\begin{center}
\begin{tabular}{lc}
\rule[-2ex]{0pt}{4ex} the measurings & $f:C\otimes T(X)\to A$,\\
\rule[-2ex]{0pt}{4ex} the algebra maps & $T(X)\to [C,A]$,\\
\rule[-2ex]{0pt}{4ex} the linear maps &$X\to [C,A]$,\\
\rule[-2ex]{0pt}{4ex} the linear maps &  $C\otimes X\to A$, \\
\rule[-2ex]{0pt}{4ex} the linear maps &  $C\to [X,A]$, \\
\rule[-2ex]{0pt}{4ex} and the coalgebra maps & $s:C\to  T^\vee([X,A])$.
\end{tabular}
\end{center}
Hence the functor $\cM(-,T(X);A):\dgCoalg \to \Set$ is represented by the coalgebra $T^\vee([X,A])$.
If $C= T^\vee([X,A])$ and $s=id$, then $f=p'$.
This shows that the measuring $p'$ is couniversal.
\end{proof}

\begin{cor}
For any dg-vector spaces $X$ and $Y$, we have
$$\{T(X),T(Y)\}=T^\vee([(X,T(Y)]).$$
Hence the Swedler hom object between free algebras is cofree.  
\end{cor}

\bigskip

For a general algebra $A$, it is possible to give a copresentation of $\{A,B\}$ in terms of a presentation of $A$.
From the monadicity of the adjunction $T:\dgVect \rightleftarrows \dgAlg:U$ implies that it is possible to present $A$ as a reflexive coequalizer in $\dgAlg$ of free algebras
$$\xymatrix{
T(T(A)) \ar@<.6ex>[r]\ar@<-.6ex>[r]& T(A)\ar[r]^-{m} & A.
}$$
$\{-,B\}$ is a contravariant right adjoint and send colimits to limits, the object $\{A,B\}$ is then the reflexive equalizer in $\dgCoalg$
$$\xymatrix{
\{A,B\}\ar[r]^-{m'} &\{T(A),B\} \ar@<.6ex>[r]\ar@<-.6ex>[r]& \{T(T(A)),B\}
}$$
Using proposition \ref{exSweedlerhom}, this is the same equalizer as
$$\xymatrix{
\{A,B\}\ar[r]^-{m'} &T^\vee([A,B]) \ar@<.6ex>[r]\ar@<-.6ex>[r]& T^\vee([T(A),B]).
}$$

\begin{cor}\label{betterShom}
$\{A,B\}$ is naturally a subcoalgebra of $T^\vee([A,B])$.
Moreover, the map $\{A,B\}\to [A,B]$, obtained by composition with the natural projection $q:T^\vee([A,B])\to [A,B]$, is the couniversal comeasuring $\Psi$
$$\xymatrix{
\{A,B\}\ar[rrd]_{\Psi}\ar[rr]^-{m'}&& T^\vee([A,B])\ar[d]^q\\
&&[A,B].
}$$
\end{cor}
\begin{proof}
The monadicity theorem that we proved says that coreflexive equalizers are reflected by the forgetful functor $U:\dgCoalg\to \dgVect$.
Applied to the equalizer above we deduced that the map $\{A,B\}\to T^\vee([A,B])$ is injective.

To prove the second statement let us consider the commutative diagram
$$\xymatrix{
\{A,B\}\ar[r]^-{m'}\ar[d]^{\Psi}&\{T(A),B\}\ar[d]^{\Psi}\\
[A,B]\ar[r]\ar@{=}[rd]&[T(A),B]\ar[d]^f\\
&[A,B]
}$$
where the horizontal map are induced by the algebra map $T(A)\to A$ and the map $f$ is induced by the natural inclusion $A\to T(A)$.
Then, to prove that $\Psi = qm'$, it is enough to prove that $f\Psi = q$.
But, by definition of $p$ and $p'$ in proposition \ref{exSweedlerhom}, we have $q=\lambda^2(p):T^\vee([X,A])\to [X,A]$ and $\lambda^2(p)$ is the composition of the universal comeasuring $\Psi=\lambda^2(p'):T^\vee([X,A])\to [T(X),A]$ with the 'restriction to the generators' map $f:[T(X),A]\to [X,A]$.
\end{proof}

Recall the notion of cogenerating and separating maps from definition \ref{coseparation}.

\begin{prop}
The comeasuring $\Psi:\{A,B\}\to [A,B]$ is cogenerating, hence separating.
\end{prop}

\begin{proof}
Trivial from corollary \ref{betterShom}.
\end{proof}

\medskip
We will need the following lemma to study the strength of the functor $\rhd$ in proposition \ref{strengthSproduct}.

\begin{lemma}\label{lemmacomposeparationalg}
If $\phi:X\to A$ is a separating map for the algebra $A$, then, for any algebra $B$, the composition
$$\xymatrix{
h:\{A,B\} \ar[r]^-\Psi & [A,B]\ar[rr]^-{\hom(\phi,B)}&& [X,B]
}$$
is separating.
\end{lemma}
\begin{proof}
The proof is dual to that of lemma \ref{lemmacomposeparation} using proposition \ref{exSweedlerhom} and the diagram
$$\xymatrix{
\{A,B\}\ar[rr]^-{\{f,B\}}\ar[d]_\Psi&& \{T(X),B\}\ar@{=}[r]\ar[d]_\Psi& T^\vee([X,B])\ar[d]^p\\
[A,B]\ar[rr]^-{\hom(f,B)} && [T(X),B] \ar[r]& [X,B].
}$$
where $f:T(X)\to A$ is the algebra map corresponding to $\phi:X\to A$.
\end{proof}

\subsubsection{Reduction maps and proofs by reduction}\label{proofbyreduction}\label{reductionseparation}

We would like to emphasize in this section a technique of proof that we have used several times already and that we will use again.

\medskip
We have shown the existence of distinguished maps
\begin{center}
\begin{tabular}{lc}
\rule[-2ex]{0pt}{4ex} the universal comorphism & $\Psi:\HOM(C,D)\to [C,D]$\\
\rule[-2ex]{0pt}{4ex} the universal measuring & $\Phi:C\otimes A\to C\rhd A$,\\
\rule[-2ex]{0pt}{4ex} and the couniversal comeasuring & $\Psi:\{A,B\}\to [A,B]$.
\end{tabular}
\end{center}
We shall call the $\Psi$ maps {\em reduction maps} and the $\Phi$ map the {\em coreduction map}. 
They all have a universal property and a separating property (definitions \ref{separation} and \ref{coseparation}).

\medskip

We shall call {\em reduction} a technique to prove the commutation of certain diagrams of coalgebras or algebras such as diagrams 
as 
the pentagon in the proof of theorem \ref{actionsweedlerproduct},
the top faces of the cube and prism in the proof of theorem \ref{enrichmentalgcoal},
or the top face of the cube in proposition \ref{functorialityofcolax}.

Let us look at the cube of theorem \ref{enrichmentalgcoal}. The diagram of interest is a square of coalgebras which are tensor products of Sweedler homs $\{A,B\}$. To study it, we introduce an analogous square where the coalgebras $\{A,B\}$ are replaced by $[A,B]$ (the bottom face of the cube). Between these two squares, we have vertical maps which are given by tensor products of reduction maps $\Psi$. The other example of the prism in the same proof is based on a triangle rather thant a square but the vertical maps are still reduction maps. We shall call informally such diagrams {\em reduction diagrams}.
The top and bottom faces of a reduction diagram have the same shape : they have both an initial (named $I_{top}$ and $I_{bot}$) and a terminal object (named $T_{top}$ and $T_{bot}$) and they have two sides whose composition we want to prove equal.
In practice, the bottom face of a reduction diagram is easily proven to be commutative so we will assume this. 
Also we will assume that the map $T_{top}\to T_{bot}$ is a single reduction maps $\Psi$.

Let us picture such a diagram as
$$\xymatrix{
I_{top}\ar@/^1pc/@{-->}[rr]^a\ar@/_1pc/@{-->}[rr]_b\ar[d]_{\Psi\otimes \dots \otimes\Psi}&& T_{top}\ar[d]^{\Psi}\\
I_{bot}\ar@{-->}[rr]^c&& T_{bot}
}$$
where the arrows $a$ and $b$ represent the total composite of the two sides of the top diagram
and where the arrow $c$ represents the total composite of the bottom diagram.

Our purpose is to prove $a=b$. Because the right map $\Psi$ is separating, this is equivalent to $\Psi a=\Psi b$.
So the proof will be finished if we know that the lateral faces of the diagram commute.

\medskip
The maps of the top face of a reduction diagram will always be defined using the universal property of some $\Psi$s (or as tensor product of such) and to each of these maps will be associated a commutative square. For example, in the cube of theorem \ref{enrichmentalgcoal}, the map ${\bf c}:\{B,E\}\otimes \{A,B\} \to \{A,E\}$ is associated to the commutative square
$$\xymatrix{
\{B,E\}\otimes \{A,B\} \ar[r]^-{\bf c}\ar[d]_{\Psi\otimes \Psi}& \{A,E\}\ar[d]^{\Psi}\\
[B,E]\otimes [A,B] \ar[r]& [A,E].
}$$
Let us call informally these squares {\em reduction squares}.
In practice, the lateral faces of a reduction diagram decomposes into reduction squares so their commutativity is never difficult to prove.

\bigskip
This technique of {\em proof by reduction} will be used again in section \ref{monoidalalg} to prove that the tensor product of algebras is compatible with the enrichment over $\dgCoalg$. Results of section \ref{laxstructures} about (co)lax structures of Sweedler operations can also be proven by reduction, but we have chosen to present them using a different approach.

\medskip
There are dual notions of {\em coreduction diagram} and {\em proof by coreduction} using the coreduction maps instead of reduction maps.
A coreduction diagram (based on a pentagon) is used implicitely in the proof of theorem \ref{actionsweedlerproduct}.

\begin{rem}
A conceptual explanation of these proof by reduction is sketched in remarks \ref{remreduction} and \ref{remreduction2}.
\end{rem}

\subsection{Monoidal strength}\label{monoidalalg}\label{monoidalstrengthalg}

In this section, we prove that the tensor product of algebras 
$$\xymatrix{
\otimes : \dgAlg \times \dgAlg \ar[r]& \dgAlg
}$$
can be enriched as a symmetric monoidal structure over $\dgCoalg$.
And that all Sweedler operations on algebras are strong functors.

\medskip

$\dgAlg \times \dgAlg$ is naturally enriched bicomplete over $\dgCoalg \times \dgCoalg$ with all operations defined termwise.
For the tensor product of algebras to be an enriched functor, we need to enrich $\dgAlg \times \dgAlg$ over $\dgCoalg$.
We can do this by using a transfer along the monoidal functor $\otimes:\dgCoalg\times \dgCoalg \to \dgCoalg$. 
We shall call $(\dgAlg\times \dgAlg)^{\$\otimes \$}$ the resulting enriched category, 
for any pairs of algebras $(A_1,A_2)$ and $(B_1,B_2)$ its hom coalgebra are
$$
\{(A_1,A_2),(B_1,B_2)\} = \{A_1,B_1\}\otimes \{A_2,B_2\}.
$$
As in section \ref{monoidalstrengthcoalg}, this enrichment is neither tensored nor cotensored over $\dgCoalg$.

The main result of this section is that $\otimes$ is a strong symmetric monoidal structure on $\dgAlg$ (theorem \ref{tensorenrichmentof}).
We shall see in section \ref{laxstructures} that $\otimes$ can also be viewed as an enriched functor in another way (proposition \ref{laxmoduleotimes}).

\bigskip

If $A_1$, $A_2$, $B_1$ and $B_2$ are algebras, then there is a unique map of coalgebras
$$\xymatrix{
\Theta_\otimes:\{A_1,B_1\}\otimes \{A_2,B_2\} \ar[r]&  \{A_1\otimes A_2, B_1\otimes B_2\}
}$$
such that the following square commutes,
\begin{align}
\tag{Strength $\otimes$ (alg)}\label{strengthotimesalg}
\vcenter{\xymatrix{
\{A_1,B_1\}\otimes \{A_2,B_2\} \ar[d]_{\Psi\otimes \Psi} \ar[rr]^-{\Theta_\otimes} &&\{A_1\otimes A_2, B_1\otimes B_2\}\ar[d]^{\Psi}\\
[A_1,B_1]\otimes [A_2,B_2] \ar[rr]^-{\gamma} &&[A_1\otimes A_2, B_1\otimes B_2],
}}
\end{align}
where $\gamma$ is the canonical map and the $\Psi$s are the couniversal comeasurings.

\begin{prop}\label{functorialityofcolax}
The maps $\Theta_\otimes$ enhanced $\otimes$ into a strong functor 
$$\xymatrix{
\otimes :(\dgAlg \times \dgAlg)^{\$\otimes \$} \to \dgAlg^\$.
}$$
\end{prop}

\begin{proof} 
We need to verify that the following diagram commutes
$$\xymatrix{
\{(B_1,B_2),(E_1,E_2)\} \otimes \{(A_1,A_2),(B_1,B_2)\} \ar[d]_{\Theta\otimes \Theta} \ar[rr]^-{\mathbf{c}} && \{(A_1,A_2),(E_1,E_2)\}
\ar[d]^\Theta \\
\{B_1\otimes B_2, E_1\otimes E_2\}  \otimes  \{A_1\otimes A_2, B_1\otimes B_2\}  \ar[rr]^-{\mathbf{c}} &&
\{A_1\otimes A_2, E_1\otimes E_2\}.
}$$
The proof is a reduction using the cube 
{\tiny
$$
\xymatrix{
&\{(B_1,B_2),(E_1,E_2)\} \otimes \{(A_1,A_2),(B_1,B_2)\} \ar[ld]_{\Theta\otimes \Theta} \ar[rr]^-{\mathbf{c}}\ar'[d]^{\Psi\otimes \Psi}[dd]
	&& \{(A_1,A_2),(E_1,E_2)\}\ar[ld]_\Theta \ar[dd]^{\Psi}\\
\{B_1\otimes B_2, E_1\otimes E_2\}  \otimes  \{A_1\otimes A_2, B_1\otimes B_2\}  \ar[rr]^-(0.7){\mathbf{c}} \ar[dd]^{\Psi\otimes \Psi}
	&& \{A_1\otimes A_2, E_1\otimes E_2\}\ar[dd]^(0.3){\Psi} \\
&[B_1,E_1]\otimes [B_2,E_2] \otimes [A_1,B_1]\otimes [A_2,B_2] \ar[ld]_{\gamma \otimes \gamma} \ar'[r]^-{c}[rr] && 
[A_1,E_1]\otimes [A_2,E_2] \ar[ld]_\gamma \\
[B_1\otimes B_2, E_1\otimes E_2]  \otimes  [A_1\otimes A_2, B_1\otimes B_2]  \ar[rr]^-{c} &&
[A_1\otimes A_2, E_1\otimes E_2].
}$$}
The bottom commutes because the tensor product is a strong functor in $\dgVect$ 
and the lateral faces commutes by definition of $\bf c$ and $\Theta$.
It follows that the top face commutes after composition with the map
$\Psi:  \{A_1\otimes A_2, E_1\otimes E_2\}  \to [A_1\otimes A_2, E_1\otimes E_2].$
This proves that the top face commutes, since $\Psi$ is separating by lemma \ref{cogeneratingcomeasuring}. 

We leave to the reader the verification that the functor $\otimes :\dgAlg \times \dgAlg \to \dgAlg$ preserves the identities.
\end{proof}

\bigskip

As in the case of two variables, we can define an enrichment $(\dgCoalg^{\times 3})^{\$\otimes\$\otimes\$}$ of $\dgCoalg^{\times 3}$ over $\dgCoalg$ and prove that the functors $-\otimes (-\otimes -)$ and $(-\otimes -)\otimes -$ are strong.
Moreover, the associator and symmetry of the monoidal structure of $\dgAlg$ are strong natural transformations. This is the meaning of the following theorem.

\begin{thm}\label{tensorenrichmentof}
The category $\dgAlg$ is symmetric monoidal as a category enriched over the category $\dgCoalg$. 
\end{thm}

\begin{proof} 
Let us show that the associativity isomorphism
$$\xymatrix{
as =as(A_1,A_2,A_3):(A_1\otimes A_2)\otimes A_3\ar[r]& A_1\otimes (A_2\otimes A_3)
}$$
is defining a strong natural transformation. For this we have 
to show that the following square commutes for six-tuples of algebras
$(A_1, A_2, A_3, B_1,B_2,B_3)$.
$$
\xymatrix{
  \{A_1, B_1\}\otimes  \{A_2, B_2\} \otimes  \{A_3, B_3\}  \ar[dd]_{\alpha(\{A_1, B_1\} \otimes \alpha)} \ar[rrr]^-{\alpha(\alpha \otimes\{A_3,B_3\})} && &
   \{(A_1\otimes A_2)\otimes A_3, (B_1\otimes B_2)\otimes B_3\} \ar[dd]^{ \{(A_1\otimes A_2)\otimes A_3, as \} }  \\
  & &&\\
  \{A_1\otimes (A_2\otimes A_3), B_1\otimes (B_2\otimes B_3) \}
  \ar[rrr]^-{  \{as, B_1\otimes (B_2\otimes B_3)\} }  &&&  \{(A_1\otimes A_2)\otimes A_3, B_1\otimes (B_2\otimes B_3) \} 
}$$
We will prove it by reduction using the cube 
{\tiny
$$\xymatrix{
&\{A_1, B_1\}\otimes  \{A_2, B_2\} \otimes  \{A_3, B_3\}  \ar[ld]_-{\alpha(\{A_1, B_1\} \otimes \alpha)\qquad} \ar[rr]^-{\alpha(\alpha \otimes\{A_3,B_3\})} \ar'[d]_-{\Psi\otimes \Psi\otimes \Psi}[dd]
	&& \{(A_1\otimes A_2)\otimes A_3, (B_1\otimes B_2)\otimes B_3\} \ar[ld]_-{ \{(A_1\otimes A_2)\otimes A_3, as \} \qquad} \ar[dd]^\Psi \\
\{A_1\otimes (A_2\otimes A_3), B_1\otimes (B_2\otimes B_3) \} \ar[rr]^-(.75){\{as, B_1\otimes (B_2\otimes B_3)\} }  \ar[dd]_\Psi
	&& \{(A_1\otimes A_2)\otimes A_3, B_1\otimes (B_2\otimes B_3)\} \ar[dd]^-(0.3)\Psi \\
&[A_1, B_1]\otimes  [A_2, B_2] \otimes  [A_3, B_3]  \ar[ld]_{\gamma ([A_1, B_1] \otimes \gamma)\qquad} \ar'[r]^-{\gamma (\gamma \otimes[A_3,B_3])}[rr]
	&& [(A_1\otimes A_2)\otimes A_3, (B_1\otimes B_2)\otimes B_3] \ar[ld]_{[(A_1\otimes A_2)\otimes A_3, as ] \qquad}  \\
[A_1\otimes (A_2\otimes A_3), B_1\otimes (B_2\otimes B_3)] \ar[rr]^-{ [as, B_1\otimes (B_2\otimes B_3)] }
	&&  [(A_1\otimes A_2)\otimes A_3, B_1\otimes (B_2\otimes B_3)] 
}$$}

\noindent The bottom face commutes, since the category $\dgVect$ is symmetric monoidal as a category enriched over itself.
The back and left faces commute essentially by definition of $\alpha$.
Finally the right and front faces commute because $\Psi$ is compatible with composition by lemma \ref{compocomes}.

It follows that the top face commutes after composition with the map
$\Psi: \{(A_1\otimes A_2)\otimes A_3, B_1\otimes (B_2\otimes B_3)  \to [(A_1\otimes A_2)\otimes A_3, B_1\otimes (B_2\otimes B_3)].$
This proves that the top face commutes, since $\Psi$ is separating by lemma \ref{cogeneratingcomeasuring}.
This shows that the enriched tensor product is equipped with an associativity isomorphism.
It is similarly equipped with strong left and right unit isomorphisms.
Mac Lane's coherence conditions are satisfied, since they are true for un-enriched tensor
product of algebras. We can prove similarly that the symmetry isomorphism
$\sigma(A,B):A\otimes B\to B\otimes A$ is strong, and that Mac Lane's coherence conditions 
for a symmetric monoidal structure are satisfied.
\end{proof}

\subsection{Strength and lax structures}\label{laxstructures}

In this section, we prove that the four Sweedler operations on algebras are equipped with natural strong lax or colax structures.

\bigskip

We already now that the tensor product of algebras is a strong functor, we prove the result for the other operations.
Let us consider the following maps.
$$\xymatrix{
\Theta_{\{-,-\}}:=\Lambda^2({\bf c^2}):\{B_1,A_1\}\otimes \{A_2,B_2\} \ar[r]& \HOM(\{A_1,A_2\},\{B_1,B_2\})
}$$
where ${\bf c^2}= \{B_1,A_1\}\otimes \{A_1,A_2\}\otimes \{A_2,B_2\}  \to \{B_1,B_2\}$,
$$\xymatrix{
\Theta_\rhd:= \Lambda^3({\bf ev}\rhd {\bf ev}):\HOM(C,D)\otimes \{A,B\} \ar[r]& \{C\rhd A,D\rhd B\}
}$$
where
${\bf ev}\rhd {\bf ev}=\HOM(C,D)\otimes \{A,B\}\otimes C\rhd A \simeq (\HOM(C,D)\otimes C)\rhd (\{A,B\} \rhd A) \to D\rhd B$,
and
$$\xymatrix{
\Theta_{[-,-]}:=\Lambda^2(\lambda^4({\bf ev}^3)):\HOM(D,C)\otimes \{A,B\} \ar[r]& \{[C,A],[D,B]\}
}$$
where ${\bf ev}^3$ is the composition
$$\xymatrix{
\{A,B\} \otimes [C,A]\otimes \HOM(D,C)\otimes D  \ar[d]_{\{A,B\} \otimes [C,A]\otimes {\bf ev}}\\
\{A,B\} \otimes [C,A]\otimes C  \ar[rr]^-{\{A,B\} \otimes ev} &&
\{A,B\} \otimes A  \ar[r]^-{\bf ev} & B.
}$$

\begin{prop}\label{strengthSHOM}
The map $\Theta_{\{-,-\}}$ is the unique coalgebra map such that the following diagram commutes.
\begin{align}
\tag{Strength $\{-,-\}$}\label{strengthSHOMdiag}
\vcenter{\xymatrix{
\{B_1,A_1\}\otimes \{A_2,B_2\} \ar[rr]^-{\Theta_{\{-,-\}}}\ar[d]_{\Psi\otimes \Psi}&&  \HOM(\{A_1,A_2\},\{B_1,B_2\}) \ar[d]^{\Psi}\\
[B_1,A_1]\otimes [A_2,B_2] \ar[d]_-{\theta} && [\{A_1,A_2\},\{B_1,B_2\}]\ar[d]^{[\{A_1,A_2\},\Psi]}\\
[[A_1, A_2],[B_1, B_2]] \ar[rr]^-{[\Psi,[B_1, B_2]]}&& [\{A_1,A_2\},[B_1, B_2]]
}}
\end{align}
where $\theta$ is the strength of $[-,-]$ in $\dgVect$.
\end{prop}
\begin{proof}
The proof of the commutation of this diagram is a careful unravelling of the definition of $\Theta_{\{-,-\}}=\Lambda^2(\bf c^2)$ left to the reader.
Then, the assertion will be proven if we show that the right side map 
$[\{A_1,A_2\},\Psi]\circ \Psi:\HOM(\{A_1,A_2\},\{B_1,B_2\})\to [\{A_1,A_2\},[B_1, B_2]]$
is separating. This is a consequence of lemma \ref{lemmacomposeparation}.
\end{proof}

\begin{prop}\label{strengthSproduct}
The map $\Theta_\rhd$ is the unique algebra map such that the following diagram commutes.
\begin{align}
\tag{Strength $\rhd$}\label{strengthSproductdiag}
\vcenter{\xymatrix{
\HOM(C,D)\otimes \{A,B\} \ar[rr]^-{\Theta_\rhd}\ar[d]_{\Psi\otimes \Psi}&&   \{C\rhd A,D\rhd B\} \ar[d]^{\Psi}\\
[C,D]\otimes [A,B] \ar[d]_-{\theta} && [C\rhd A,D\rhd B]\ar[d]^{[\Phi,D\rhd B]}\\
[C\otimes A,D\otimes B] \ar[rr]^-{[C\otimes A,\Phi]}&& [C\otimes A,D\rhd B]
}}
\end{align}
where $\theta$ is the strength of $\otimes$ in $\dgVect$.
\end{prop}
\begin{proof}
The proof of the commutation of this diagram is a careful unravelling of the definition of $\Theta_{\rhd}$ left to the reader.
The unicity of $\Theta_\rhd$ is proven by lemma \ref{lemmacomposeparationalg}.
\end{proof}

\begin{prop}\label{strengthconvol}
The map $\Theta_{[-,-]}$ is the unique coalgebra map such that the following diagram commutes.
\begin{align}
\tag{Strength $[-,-]$}\label{strengthconvoldiag}
\vcenter{\xymatrix{
\HOM(D,C)\otimes \{A,B\}\ar[rr]^-{\Theta_{[-,-]}}\ar[d]_{\Psi\otimes \Psi}&& \{[C,A],[D,B]\} \ar[d]^{\Psi}\\
[D,C]\otimes [A,B] \ar[rr]_-{\theta} && [[C,A],[D,B]]
}}
\end{align}
where $\theta$ is the strength of $[-,-]$ in $\dgVect$.
\end{prop}
\begin{proof}
The proof of the commutation of this diagram is a careful unravelling of the definition of $\Theta_{[-,-]}$ left to the reader.
Then, the assertion is proven by the separation property of $\Psi:\{[C,A],[D,B]\}\to [[C,A],[D,B]]$.
\end{proof}

\medskip

\begin{prop}\label{strongSP}\label{strongconvol}
\begin{enumerate}
\item The maps $\Theta_{\{-,-\}}$ enhance the Sweedler hom $\{-,-\}$ into a strong functor
$$\xymatrix{
\{-,-\}:(\dgAlg^{op}\times \dgAlg)^{\$\otimes \$} \ar[r] &\dgCoalg^\$.
}$$

\item The maps $\Theta_\rhd$ enhance the Sweedler product $\rhd$ into a strong functor
$$\xymatrix{
\rhd:(\dgCoalg\times \dgAlg)^{\$\otimes \$} \ar[r] &\dgAlg^\$.
}$$

\item The maps $\Theta_{[-,-]}$ enhance the convolution product $[-,-]$ into a strong functor
$$\xymatrix{
[-,-]:(\dgCoalg^{op}\times \dgAlg)^{\$\otimes \$} \ar[r] &\dgAlg^\$.
}$$

\end{enumerate}
\end{prop}

\begin{proof}
This is always the case when an enriched category is tensored and cotensored (see appendix \ref{enrichedcategoryapp}).
\end{proof}

Notice that because of the variance of the Sweedler hom and convolution product, their strength give morphisms
$$
\{A,A\}^o\otimes \{B,B\} \to \HOM(\{A,B\},\{A,B\})
\et
\END(C)^o\otimes \{A,A\} \to \{[C,A],[C,A]\}
$$
where $\{A,A\}^o$ and $\END(C)^o$ are the opposite bialgebras of $\{A,A\}$ and $\END(C)$ (definition \ref{defoppositebialgebra}).

\bigskip
We now turn to the proof of the lax structures.
We have shown that the measuring functor
$$\xymatrix{
\cM:\dgCoalg^{op}\times \dgAlg^{op}\times \dgAlg\ar[r]& \Set.
}$$
is representable in all its variables and interpreted this as a bicomplete enrichment of $\dgAlg$ over $\dgCoalg$.
In this section we prove that the measuring functor is compatible with the monoidal structure of $\dgCoalg$ and $\dgAlg$ in the sense that it is a lax monoidal functor. 
As a consequence, all the operations representing $\cM$ will inherit a lax or colax structure.

\bigskip

Recall from \ref{trialitydgVect} the lax monoidal triality associated to the symmetric monoidal closed structure of $\dgVect$
$$\xymatrix{
\cT:\dgVect^{op}\times \dgVect^{op}\times \dgVect \ar[r]& \Set
}$$
where $\cT(X,Y;Z) = [X\otimes Y,Z]_0$.
When compose with the lax monoidal forgetful functor
$$\xymatrix{
U\times U\times U:\dgCoalg^{op}\times\dgAlg^{op}\times\dgAlg\ar[r]&\dgVect^{op}\times\dgVect^{op}\times\dgVect
}$$
it gives a lax monoidal functor
$$\xymatrix{
\cU:\dgCoalg^{op}\times\dgAlg^{op}\times\dgAlg\ar[r]&\Set.
}$$
The set $\cM(C,A;B)=\dgAlg(A,[C,B])$ is naturally a subset of the set of linear maps $[A,[C,B]]_0$ and this gives an inclusion of functors
$$
\cM \subset \cU.
$$
We are going to prove that the lax structure of $\cU$ restrict to $\cM$.
The proof is based on the following lemma.

\begin{lemma}\label{laxlemma}
Let $C_1$ and $C_2$ be two coalgebras and $A_1$ and $A_2$ be two algebras.
The canonical map
$$\xymatrix{
\alpha:[C_1,A_1]\otimes [C_2,A_2]\ar[r]&  [C_1\otimes C_2,A_1\otimes A_2]
}$$
is an algebra map.
\end{lemma}
\begin{proof}
A straightforward computation.
\end{proof}

\begin{prop}\label{laxmeasuring}
The measuring functor 
$$\xymatrix{
\cM:\dgCoalg^{op}\times \dgAlg^{op}\times \dgAlg\ar[r]& \Set
}$$
is a symmetric lax monoidal functor.
\end{prop}
\begin{proof}
We need only to prove that the lax structure of $\cU$ restrict to $\cT$.
The lax structure of $\cU$ is given by the surjection
$\beta_0:[\FF,[\FF,\FF]]_0\simeq \FF \to \{*\}$ and 
the map
$$\xymatrix{
\beta:Z_0[A_1,[C_1;B_1]]\times Z_0[A_2,[C_2;B_2]]\ar[r]& Z_0[A_1\otimes A_2,[C_1\otimes C_2;B_1\otimes B_2]].
}$$
The map $\beta_0$ restricts to $\cM(\FF,\FF;\FF)\subset [\FF,[\FF,\FF]]_0$ to give a bijection $\cM(\FF,\FF;\FF)\simeq \{*\}$
and $\beta$ restricts to $\cM$ iff the following square commutes
$$\xymatrix{
\dgAlg(A_1,[C_1,B_1])\times \dgAlg(A_1,[C_2,B_2])\ar[r]\ar[d] &  \dgAlg(A_1\otimes A_2,[C_1\otimes C_2,B_1\otimes B_2])\ar[d]\\
Z_0[A_1,[C_1,B_1]]\times Z_0[A_2,[C_2,B_2]]\ar[r]&  Z_0[A_1\otimes A_2,[C_1\otimes C_2,B_1\otimes B_2]]
}$$
where the vertical maps are the inclusions and the horizontal maps are essentially the composition with the map $\alpha$ of lemma \ref{laxlemma}. 
The lower composite map sends $(f,g)$ to $\alpha\circ (f\otimes g)$. But the tensor product of two algebra maps is an algebra map, and by lemma \ref{laxlemma} the composition with $\alpha$ is still an algebra map.

\end{proof}

For the purposes of proofs to come, let us introduce notations for the associativity condition on the lax structure of $\cM$.
For $i=1,2,3$, 
let $X_i=(C_i,A_i;B_i)$ be a triplet of a coalgebra and two algebras and let 
let $X_i\otimes X_j=(C_i\otimes C_j,A_i\otimes A_j;B_i\otimes B_j)$,
the condition on $\cM$ is then
\begin{align}
\tag{Lax $\cM$}\label{laxM}
\vcenter{\xymatrix{
\cM(X_1)\times \cM(X_2)\times \cM(X_3)\ar[d]_{\beta\times id}\ar[r]^-{id\times \beta}& \cM(X_1)\times \cM(X_2 \otimes X_3)\ar[d]^{\beta}\\
\cM(X_1 \otimes X_2)\times \cM(X_3)\ar[r]_-{\beta}& \cM(X_1\otimes X_2\otimes X_3)
}}
\end{align}

\begin{cor}\label{convolutionlax}
If $C_1$ and $C_2$ are coalgebras and $A_1$ and $A_2$ are algebras,  then
the natural map
$$\xymatrix{
\alpha_{[-,-]}:=\alpha: [C_1,A_1]\otimes [C_2, A_2] \ar[r]& [C_1\otimes C_2, A_1\otimes A_2]
}$$
together with the isomorphism $\alpha_0: [\FF,\FF]=\FF$ define a strong symmetric lax structure $(\alpha_{[-,-]},\alpha_0)$ on the functor $[-,-]:\dgCoalg^{op} \times \dgAlg \to \dgAlg$.
\end{cor}

\begin{proof}
The associativity condition is the commutation of the square
$$\xymatrix{
[C_1,A_1]\otimes [C_2, A_2]\otimes [C_3, A_3] \ar[r]\ar[d] &[C_1,A_1]\otimes [C_2\otimes C_3, A_2\otimes A_3]\ar[d]\\
[C_1\otimes C_2, A_1\otimes A_2]\otimes [C_3, A_3] \ar[r]& [C_1\otimes C_2\otimes C_3, A_1\otimes A_2\otimes A_3].
}$$
for any triple of coalgebras $(A_1,A_2,A_3)$ and triple of coalgebras $(C_1,C_2,C_3)$.
Let $f$ and $g$ be the two side of this square.
This is a consequence of the condition \eqref{laxM} applied to $X_i=(C_i,[C_i,A_i],A_i)$: in the square
$$\xymatrix{
\cM(X_1)\times \cM(X_2)\times \cM(X_3)\ar[d]_{\beta\times id}\ar[r]^-{id\times \beta}& \cM(X_1)\times \cM(X_2 \otimes X_3)\ar[d]^{\beta}\\
\cM(X_1 \otimes X_2)\times \cM(X_3)\ar[r]_-{\beta}& \cM(X_1\otimes X_2\otimes X_3)
}$$
the two images of the triplet $(id_1,id_2,id_3)$ where $id_i$ is the identity of $[C_i,A_i]$ are exactly the two morphisms $f$ and $g$.
The unitality and symmetry conditions are proven similarly.

Let us now prove that $\alpha_{[-,-]}$ is a strong natural transformation.
We need to prove the commutativity of
$$\xymatrix{
\HOM(D_1,C_1)\otimes \HOM(D_2,C_2)\otimes \{A_1,B_1\}\otimes \{A_2,B_2\} \ar[r]\ar[d]& \{[C_1, A_1]\otimes [C_2, A_2],[D_1, B_1]\otimes [D_2, B_2]\}\ar[d]\\
\{[C_1\otimes C_2,A_1\otimes A_2],[D_1\otimes D_2,B_1\otimes B_2]\} \ar[r] &\{[C_1\otimes C_2,A_1\otimes A_2],[D_1, B_1]\otimes [D_2, B_2]\}.
}$$
This can be done by reduction using a cube whose bottom face
$$\xymatrix{
[D_1,C_1]\otimes [D_2,C_2]\otimes [A_1,B_1]\otimes [A_2,B_2] \ar[r]\ar[d]& [[C_1, A_1]\otimes [C_2, A_2],[D_1, B_1]\otimes [D_2, B_2]]\ar[d]\\
[[C_1\otimes C_2,A_1\otimes A_2],[D_1\otimes D_2,B_1\otimes B_2]] \ar[r] &[[C_1\otimes C_2,A_1\otimes A_2],[D_1, B_1]\otimes [D_2, B_2]].
}$$
We leave the details to the reader.
\end{proof}

\bigskip

Let $C_1$, $C_2$ be two coalgebras and $A_1$, $A_2$ be two algebras, we define a map
$$\xymatrix{
\alpha_\rhd:(C_1\otimes C_2) \rhd (A_1\otimes A_2) \ar[r]& (C_1\rhd A_1)\otimes (C_2\rhd A_2)
}$$
as the unique map map of algebras such that we have
$$
\alpha_\rhd\bigl((c_1\otimes c_2)\rhd (a_1\otimes a_2)\bigr)=(c_1\rhd a_1)\otimes (c_2\rhd a_2)(-1)^{|c_2||a_1|}
$$
for every $c_1\in C_1, c_2\in C_2$, $a_1\in A_1$ and $a_2\in A_2$.

\begin{cor}\label{Sproductcolax}
If $C_1$ and $C_2$ are coalgebras and $A_1$ and $A_2$ are algebras, then the map
$$\xymatrix{
\alpha_\rhd:(C_1\otimes C_2) \rhd (A_1\otimes A_2) \ar[r]& (C_1\rhd A_1)\otimes (C_2\rhd A_2)
}$$
together with the isomorphism $\alpha_0: [\FF,\FF]=\FF$ define a strong symmetric colax structure $(\alpha_\rhd,\alpha_0)$ 
on the functor $\rhd:\dgCoalg \times \dgAlg \to \dgAlg$.
\end{cor}

\begin{proof}
We can present more conceptually the construction of $\alpha$, if we consider the image $\alpha'$ of the pair $(f_1,f_2)$ of the universal measurings $f_i:C_i\otimes A_i\to C_i\rhd A_i$ by 
$$
\beta:\cM(C_1,A_1,C_1\rhd A_1)\times \cM(C_2,A_2,C_2\rhd A_2) \to \cM(C_1\otimes C_2,A_1\otimes A_2,(C_1\rhd A_1)\otimes (C_2\rhd A_2)),
$$
the map $\alpha$ is the image of $\alpha'$ under the bijection
$$
\cM(C_1\otimes C_2,A_1\otimes A_2,(C_1\rhd A_1)\otimes (C_2\rhd A_2)) \simeq \dgAlg((C_1\otimes C_2) \rhd (A_1\otimes A_2),(C_1\rhd A_1)\otimes (C_2\rhd A_2)).
$$

The associativity condition is the commutation of the square
$$\xymatrix{
(C_1\otimes C_2\otimes C_3) \rhd (A_1\otimes A_2 \otimes A_3)  \ar[r] \ar[d] &
(C_1\rhd A_1)\otimes \bigl((C_2\otimes C_3) \rhd (A_2\otimes A_3)\bigr)  \ar[d] \\
\bigl((C_1\otimes C_2)\rhd (A_1\otimes A_2) \bigr)\otimes (C_3\rhd A_3)\ar[r]& (C_1\rhd A_1)\otimes (C_2\rhd A_2)\otimes(C_3\rhd A_3) 
}$$
for any triple of coalgebras $(A_1,A_2,A_3)$ and triple of coalgebras $(C_1,C_2,C_3)$.
As in corollary \ref{convolutionlax}, this is a consequence of \eqref{laxM} applied to $X_i=(C_i,A_i;C_i\rhd A_i)$.
The unitality and symmetry conditions are proven similarly.

Let us now prove that $\alpha_\rhd$ is a strong natural transformation.
We need to prove the commutativity of
$$\xymatrix{
\HOM(C_1,D_1)\otimes \HOM(C_2,D_2)\otimes \{A_1,B_1\}\otimes \{A_2,B_2\} \ar[r]\ar[d]& \{(C_1\rhd A_1)\otimes (C_2\rhd A_2),(D_1\rhd B_1)\otimes (D_2\rhd B_2)\}\ar[d]\\
\{(C_1\otimes C_2) \rhd (A_1\otimes A_2),(D_1\otimes D_2) \rhd (B_1\otimes B_2)\} \ar[r] &\{(C_1\otimes C_2) \rhd (A_1\otimes A_2),(D_1\rhd B_1)\otimes (D_2\rhd B_2)\}.
}$$
This can be done by reduction using a cube whose bottom face
$$\xymatrix{
\HOM(C_1,D_1)\otimes \HOM(C_2,D_2) \otimes \{A_1,B_1\}\otimes \{A_2,B_2\} \ar[r]\ar[d]& [(C_1\rhd A_1)\otimes (C_2\rhd A_2),(D_1\rhd B_1)\otimes (D_2\rhd B_2)]\ar[d]\\
[(C_1\otimes C_2) \rhd (A_1\otimes A_2),(D_1\otimes D_2) \rhd (B_1\otimes B_2)] \ar[r] &\{(C_1\otimes C_2) \rhd (A_1\otimes A_2),(D_1\rhd B_1)\otimes (D_2\rhd B_2)\}.
}$$
We leave the details to the reader.
\end{proof}

\bigskip

If $A_1$, $A_2$, $B_1$ and $B_2$ are algebras, recall the map of coalgebras
$$\xymatrix{
\alpha_{\{-,-\}}:=\Theta_{\otimes}:\{A_1,B_1\}\otimes \{A_2,B_2\} \ar[r]&  \{A_1\otimes A_2, B_1\otimes B_2\}
}$$
of proposition \ref{functorialityofcolax}.
It is the unique map of coalgebras such that the following square commutes
$$\xymatrix{
\{A_1,B_1\}\otimes \{A_2,B_2\} \ar[d]_{\Psi\otimes \Psi} \ar[rr]^-{\alpha_{\{-,-\}}} &&\{A_1\otimes A_2, B_1\otimes B_2\} \ar[d]^{\Psi}\\
[A_1,B_1]\otimes [A_2,B_2] \ar[rr]^-{\gamma} &&[A_1\otimes A_2, B_1\otimes B_2],
}
$$
where $\gamma$ is the canonical map and the $\Psi$s are the couniversal comeasurings.

\begin{cor}\label{Shomlax}
If $A_1$, $A_2$, $B_1$ and $B_2$ are algebras, then the map
$$\xymatrix{
\alpha_{\{-,-\}}:\{A_1,B_1\}\otimes \{A_2,B_2\} \ar[r]&  \{A_1\otimes A_2, B_1\otimes B_2\}
}$$
together with the isomorphism $\alpha_0:\FF \to \{\FF,\FF\}$, 
defines a strong symmetric lax structure $(\alpha_{\{-,-\}},\alpha_0)$ on the functor $\{-,-\}:\dgAlg^{op} \times \dgAlg \to \dgCoalg$.
\end{cor}

\begin{proof}
$\alpha_{\{-,-\}}$ can be constructed from the measuring functor $\cM$. 
Let $\alpha'$ be the image of the pair $(g_1,g_2)$ of the the couniversal comeasurings $g_i:\{A_i,B_i\}\to [A_i,B_i]$ by 
$$\xymatrix{
\beta:\cM(\{A_1,B_1\},A_1,B_1)\times \cM(\{A_2,B_2\},A_2,B_2) \ar[r]& \cM(\{A_1,B_1\}\otimes \{A_2,B_2\},A_1\otimes A_2,B_1\otimes B_2),
}$$
then, the map $\alpha_{\{-,-\}}$ is the image of $\alpha'$ under the bijection
$$
\cM(\{A_1,B_1\}\otimes \{A_2,B_2\},A_1\otimes A_2,B_1\otimes B_2) \simeq \dgCoalg(\{A_1,B_1\}\otimes \{A_2,B_2\},\{A_1\otimes A_2, B_1\otimes B_2\}).
$$

The associativity condition is the commutation of the square
$$
\xymatrix{
\{A_1, B_1\}\otimes  \{A_2, B_2\} \otimes  \{A_3, B_3\}  \ar[d] \ar[rr] && 
\{A_1, B_1\}\otimes \{A_2\otimes A_3, B_2\otimes B_3 \} \ar[d] \\
\{A_1\otimes A_2, B_1\otimes B_2\} \otimes  \{A_3, B_3\} \ar[rr] &&  \{A_1\otimes A_2\otimes A_3, B_1\otimes B_2\otimes B_3 \} 
}$$
for every six-tuples of algebras $(A_1, A_2, A_3, B_1,B_2,B_3)$.
As in corollary \ref{convolutionlax}, it is a consequence of \eqref{laxM} applied to $X_i=(\{A_i,B_i\},A_i;B_i)$.
The unitality and symmetry conditions are proven similarly.

Let us prove that $\alpha_{\{-,-\}}$ is a strong natural transformation. We need to prove the commutativity of
$$\xymatrix{
\{A_1,E_1\}\otimes \{A_2,E_2\}\otimes \{B_1,F_1\}\otimes \{B_1,F_2\} \ar[r] \ar[d] 
	& \HOM(\{A_1\otimes A_2,B_1\otimes B_2\},\{E_1\otimes E_2,F_1\otimes F_2\}) \ar[d]\\
\HOM(\{A_1,B_1\}\otimes \{A_2,B_2\},\{E_1,F_1\}\otimes \{E_2,F_2\}) \ar[r]
	&\HOM(\{A_1,B_1\}\otimes \{A_2,B_2\},\{E_1\otimes E_2,F_1\otimes F_2\}) 
}$$
for any algebras $A_1,A_2,E_1,E_2, B_1, B_2, F_1, F_2$.
This can be done by reduction using a cube whose bottom face is
$$\xymatrix{
\{A_1,E_1\}\otimes \{A_2,E_2\}\otimes \{B_1,F_1\}\otimes \{B_1,F_2\} \ar[r] \ar[d] 
	& [\{A_1\otimes A_2,B_1\otimes B_2\},\{E_1\otimes E_2,F_1\otimes F_2\}] \ar[d]\\
[\{A_1,B_1\}\otimes \{A_2,B_2\},\{E_1,F_1\}\otimes \{E_2,F_2\}] \ar[r]
	&[\{A_1,B_1\}\otimes \{A_2,B_2\},\{E_1\otimes E_2,F_1\otimes F_2\}].
}$$
We leave the details to the reader.
\end{proof}

\medskip

\bigskip
The three maps giving the (co)lax structures
$$\xymatrix{
\alpha_{[-,-]}:[C_1,A_1]\otimes [C_2,A_2]\ar[r]&  [C_1\otimes C_2,A_1\otimes A_2]
}$$
$$\xymatrix{
\alpha_\rhd:(C_1\otimes C_2) \rhd (A_1\otimes A_2) \ar[r]& (C_1\rhd A_1)\otimes (C_2\rhd A_2)
}$$
$$\xymatrix{
\alpha_{\{-,-\}}:\{A_1,B_1\}\otimes \{A_2,B_2\} \ar[r]&  \{A_1\otimes A_2, B_1\otimes B_2\}
}$$
have a second interpretation.
We have already seen that the map $\alpha_{\{-,-\}}$ was the strength of the functor $\otimes$ in theorem \ref{tensorenrichmentof}. Let us give a corresponding interpretation for the other maps.
\medskip

The Sweedler product endows $\dgAlg$ with the structure of a module over the monoidal category $(\dgCoalg,\otimes)$ (appendix \ref{Vmodule}). This structure can be restricted along the monoidal functor $\otimes:\dgCoalg\times \dgCoalg\to \dgCoalg$ to give a module structure over $\dgCoalg\times \dgCoalg$ defined by
$$
(C_1,C_2)\rhd A := (C_1\otimes C_2)\rhd A.
$$
Similarly, the convolution product endows $\dgAlg$ with an opmodule structure over $(\dgCoalg,\otimes)$ (appendix \ref{Vmodule}) and this structure can also be restricted to give an opmodule structure over $\dgCoalg\times \dgCoalg$ defined by
$$
[(C_1,C_2),A] := [C_1\otimes C_2,A]
$$

The definition of (co)lax morphisms of (op)modules is given in appendix \ref{Vmodule}.

\begin{prop}\label{laxmoduleotimes}
\begin{enumerate}
\item The map $\alpha_{[-,-,]}$ is a colax modular structure on $\otimes:\dgAlg \times \dgAlg \to \dgAlg$ for the opmodule structure over $\dgCoalg\times \dgCoalg$.
\item The map $\alpha_\rhd$ is a lax modular structure on $\otimes:\dgAlg \times \dgAlg \to \dgAlg$ for the module structure over $\dgCoalg\times \dgCoalg$.
\end{enumerate}
\end{prop}
\begin{proof}
The diagrams to check are written in appendix \ref{Vmodule}.
We leave the proof to the reader: 1) is a straightforward computation and 2) can be proven by reduction.
\end{proof}

We deduce from the proposition that the two functors $(-\otimes -)\otimes -,-\otimes (-\otimes -):(\dgAlg)^3 \to \dgAlg$ have a natural (co)lax modular structure if $(\dgAlg)^3$ and $\dgAlg$ are viewed as $(\dgCoalg)^3$-(op)modules.
Moreover we can prove that the associativity, unital and symmetry isomorphisms of $\otimes$ are strong natural transformations of (op)modules. We leave the proof of this to the reader.

\subsection{Consequences on bialgebras}\label{applicbialg}

The symmetric (co)lax monoidal structures that we have proven to exist on Sweedler operations in section \ref{laxstructures} imply that these operations induce some functors between the categories of (co)algebras in $\dgCoalg$ and $\dgAlg$ over any symmetric operad in $\dgCoalg$ (\ie Hopf operad). We shall detail only two applications to (co)commutative (co)algebras.

\medskip

We introduce the following categories:
$\dgAlg^{\sf com}$ is the category of commutative algebras,
$\dgCoalg^{\sf coc}$ be the category of cocommutative coalgebras,
$\dgBialg$ is the category of bialgebras,
$\dgBialg^{\sf com}$ is the category of commutative bialgebras
and $\dgBialg^{\sf coc}$ is the category of cocommutative bialgebras.

For a monoidal category $(\bC,\otimes)$, let $\Mon(\bC)$ and $\coMon(\bC)$ be respectively the categories of monoids and comonoids in $\bC$. All of the previous categories of (co/bi)algebras inherit the symmetric monoidal structure of $\dgVect$ and,
according to the Eckman-Hilton argument, we have the following canonical equivalences of categories:
\begin{center}
\begin{tabular}{l}
\rule[-2ex]{0pt}{4ex} $\dgAlg^{\sf com}=\Mon(\dgAlg)$, \\
\rule[-2ex]{0pt}{4ex} $\dgCoalg^{\sf coc}=\coMon(\dgCoalg)$, \\
\rule[-2ex]{0pt}{4ex} $\dgBialg=\Mon(\dgCoalg)=\coMon(\dgAlg)$, \\
\rule[-2ex]{0pt}{4ex} $\dgBialg^{\sf com}=\Mon(\dgBialg)=\coMon(\dgAlg^{\sf com})$, \\
\rule[-2ex]{0pt}{4ex} $\dgBialg^{\sf coc}=\coMon(\dgBialg)=\Mon(\dgCoalg^{\sf coc})$.
\end{tabular}
\end{center}

\begin{prop}\label{corcolaxSproduct}
The Sweedler product $\rhd$ induces functors
$$\xymatrix{
\rhd : \dgCoalg^{\sf coc} \times \dgBialg \ar[r] & \dgBialg
}$$
$$\xymatrix{
\rhd : \dgCoalg^{\sf coc} \times \dgBialg^{\sf coc} \ar[r] & \dgBialg^{\sf coc}.
}$$
\end{prop}
\begin{proof}
Direct from the symmetric colax monoidal structure of corollary \ref{Sproductcolax}.
The first result is obtain for comonoids, the second for cocommutative comonoids.
\end{proof}

If $C$ is a cocommutative coalgebra and $H$ a bialgebra, $C\rhd H$ is a bialgebra whose coproduct $\Delta:C\rhd H\to (C\rhd H)\otimes (C\rhd H)$ is defined to be the
composite 
$$\xymatrix{
C\rhd H \ar[rr]^-{\Delta\rhd \Delta} && (C\otimes C)\rhd (H\otimes H)\ar[r]^\alpha &(C\rhd H)\otimes (C\rhd H)
}$$
where $\alpha $ is the colax structure (we need $C$ to be cocommutative for the coproduct $\Delta:C\to C\otimes C$ to be an coalgebra map).
Explicitely on elements, this gives
$$
\Delta(c\rhd h)=\alpha\bigl( (c^{(1)}\otimes c^{(2)}) \rhd (h^{(1)} \otimes h^{(2)})  \bigr)=(c^{(1)}\rhd h^{(1)})\otimes (c^{(2)}\rhd h^{(2)}) (-1)^{|c^{(2)}||h^{(1)}|}.
$$
The counit of $C\rhd H$ is the map $\epsilon \rhd \epsilon :C\rhd H\to \FF\rhd \FF\simeq \FF$.

\begin{cor}\label{iterationcobar}
If $H$ is a cocommutative bialgebra, it defines an endofunctor
$$\xymatrix{
-\rhd H : \dgCoalg^{\sf coc}\ar[r] & \dgBialg^{\sf coc}\ar[r]^-{U} & \dgCoalg^{\sf coc}.
}$$
where the functor $U:\dgBialg^{\sf coc}\to \dgCoalg^{\sf coc}$ is the functor forgetting the algebra structure.
\end{cor}

\medskip

\begin{prop}\label{corlaxShom}
The Sweedler hom $\{-,-\}$ induces functors
$$\xymatrix{
\{-,-\} : \dgBialg^{op} \times \dgAlg^{\sf com} \ar[r] & \dgBialg
}$$
$$\xymatrix{
\{-,-\} : \left(\dgBialg^{\sf coc}\right)^{op} \times \dgAlg^{\sf com} \ar[r] & \dgBialg^{\sf com}.
}$$
\end{prop}
\begin{proof}
Direct from the symmetric lax monoidal structure of corollary \ref{Shomlax}.
The first result is obtain for monoids, the second for commutative monoids.
\end{proof}

If $H$ is a cocommutative bialgebra and $A$ a commutative algebra, $\{H,A\}$ is a bialgebra whose product $\mu:\{H,A\}\otimes \{H,A\}\to \{H,A\}$ is a convolution defined as the composite 
$$\xymatrix{
\{H,A\}\otimes \{H,A\} \ar[r]^-{\alpha} & \{H\otimes H,A\otimes A\} \ar[rr]^-{\{\Delta, \mu\}} && \{H,A\}
}$$
where $\alpha $ is the lax structure (we need $A$ to be commutative for the product $\mu:A\otimes A\to A$ to be an algebra map).
The unit is the map $\{\epsilon,e\}:\FF\simeq \{\FF,\FF\} \to \{H,A\}$.

\begin{cor}\label{iterationbar}
If $H$ is a cocommutative bialgebra, it defines an endofunctor
$$\xymatrix{
\{H,-\} : \dgAlg^{\sf com} \ar[r] & \dgBialg^{\sf com}\ar[r]^-{U}& \dgAlg^{\sf com}
}$$
where the functor $U:\dgBialg^{\sf com}\to \dgAlg^{\sf com}$ is the functor forgetting the coalgebra structure.
\end{cor}

\medskip
We shall apply both these results to iteratation of the bar and cobar constructions in section \ref{iteratedbarcobar}.

\subsection{Meta-morphisms}

\subsubsection{Reduction functor and meta-morphisms}\label{reducetmeta-morphismalg}

Recall that $\dgAlg^\$$ is our notation for the enrichment of $\dgAlg$ over $\dgCoalg$. 
We can transfer this enrichment along the monoidal functor $U:\dgCoalg\to \dgVect$ into an enrichment over $\dgVect$.
We denote by $\dgAlg^{U\$}$ the corresponding category. 
Recall also that $\dgVect^\$$ is our notation for $\dgVect$ viewed as enriched over itself.

\begin{prop}\label{strongforgetfulVectalg}
The reduction maps $\Psi:\{A,B\}\to [A,B]$ are the strengths of an enrichment over $\dgVect$ of the forgetful functor
$$\xymatrix{
U:\dgAlg^{U\$} \ar[r]&\dgVect^\$
}$$
\end{prop}
\begin{proof}
The compatibility with compositions is the content of lemma \ref{compocomes}.
\end{proof}

\begin{rem}
As for coalgebras, this functor is strong lax monoidal, the lax structure being given by \eqref{strengthotimesalg}. We shall not develop this.

Also $\dgAlg^{U\$}$ is nor tensored nor cotensor over $\dgVect$ and the strong functor $U$ do not have any strong adjoint anymore.
We will see in theorem \ref{strongmonadicityalg} that the adjunction $T\dashv U$ can be enriched provided we enriched the two categories over $\dgCoalg$ instead of $\dgVect$.
\end{rem}

\medskip

\begin{defi}
We define a {\it meta-morphism of algebras} $f:A\leadsto B$ to be an element $f\in U\{A,B\}$. The composite $gf$ of two meta-morphisms of coalgebras  $f:A\leadsto B$ and $g:B\leadsto E$ is defined by putting $gf=\mathbf{c}(g\otimes f)$, where $\mathbf{c}$ is the strong composition law
$$\xymatrix{
\mathbf{c}:\{B,E\}\otimes \{A,B\}\ar[r]& \{A,E\}
}$$
of section \ref{sweedlerhom}. We have $|gf|=|g|+|f|$.

By opposition, we define a {\em pro-morphism of algebras} $f:A\rhup B$ to be an element of $[A,B]$. The composition of pro-morphisms is defined through the strong composition of $\dgVect$. 
\end{defi}

\begin{rem}\label{remfactoalg}
The strong functor $U:\dgAlg^{U\$} \to \dgVect^\$$ factors as 
$$\xymatrix{
U:\dgAlg^{U\$} \ar[r]& (U\dgAlg)^{\$}\ar[r]& \dgVect^\$
}$$
where the objects of $(U\dgAlg)^{\$}$ are the algebras but the hom between two algebras $A$ and $B$ is simply $[A,B]$.
The functor $(U\dgAlg)^{\$}\to \dgVect^\$$ is strongly fully faithful. With the vocabulary of remark \ref{remfactocoalg}, we can say that the meta-morphisms are the graded morphisms of $\dgCoalg^{U\$}$ and the pro-morphisms are the graded morphisms of $(U\dgCoalg)^{\$}$.
\end{rem}

\medskip

If $f:A\leadsto B$ is a meta-morphism, the reduction map $\Psi:\{A,B\}\to [A,B]$ defines a pro-morphism $\Psi(f)$. We shall simplify notations and simply put $f^\flat :=\Psi(f)$.
If $g:B\leadsto E$ and $f:A\leadsto B$ are two meta-morphisms, proposition \ref{strongforgetfulVectalg} says that
$$
(gf)^\flat = g^\flat f^\flat.
$$

\medskip

Recall the strong evaluation ${\bf ev}=\{A,B\}\otimes A\to B$ from section \ref{sweedlerhom}.
We define the {\em evaluation} of a meta-morphism $f:A\leadsto B$ on an element $x\in A$ to be 
$$
f(x):={\bf ev}(f\otimes x)\in B.
$$
If $f:A\leadsto B$ and $g:B\leadsto E$ are meta-morphisms of coalgebras, we have $(gf)(x)=g(f(x))$ for every $x\in A$ by definition of the strong composition map from $\bf ev$.

From the definition of the reduction map $\flat=\Psi:\{A,B\}\to [A,B]$ as $\Psi=\Lambda^2{\bf ev}$, we have a commutative diagram
$$\xymatrix{
\{A,B\}\otimes A   \ar[d]_{\flat \otimes A}  \ar@/^1pc/[rrd]^-{\mathbf{ev}} &&  \\
[A,B]\otimes A \ar[rr]^-{ev} && B.
}$$
In terms of elements, this gives $f(x) = f^\flat(x)$; the strong evaluation of a meta-morphism coincides with the evaluation of the corresponding pro-morphism.

\medskip

\begin{nota}
A meta-morphism $f\in \{A,B\}$ is said to be {\em atomic} if it is an atom in $\{A,B\}$.
We have seen in lemma \ref{atomShom} that the atomic meta-morphisms in $\{A,B\}$ are the same thing as algebra maps $A\to B$. If $f:A\to B$ is a coalgebra map we shall denote by $f^\sharp$ the corresponding atom of $\{A,B\}$.
This provides a map (of graded sets) $\sharp:\dgAlg(A,B)\to \{A,B\}$.
\end{nota}

\begin{lemma}\label{musicallemmaalg}
\begin{enumerate}
\item Let $f:A\to B$ be a map of algebras, then we have 
$$
(f^\sharp)^\flat = f.
$$
In other terms, for any two algebras $A$ and $B$, we have a commutative diagram of graded sets
$$\xymatrix{
At(\{A,B\}) \ar[r]\ar[d]^{\flat}& \{A,B\}\ar[d]^{\Psi = \flat}\\
\dgAlg(A,B) \ar@/^1pc/[u]^-\sharp \ar[r]&[A,B].
}$$
where $At(\{A,B\})$ is the set of atoms of $\{A,B\}$
and the horizontal maps are the canonical inclusions.

\item The maps $\sharp$ and $\flat$ induce inverse bijections of sets
$$\xymatrix{
\dgAlg(A,B) \ar@<.6ex>[rr]^-\sharp&& At(\{A,B\}).\ar@<.6ex>[ll]^-{\flat=\Psi}
}$$

\item For a coalgebra map $f:A\to B$ and $x\in A$, $f(x)$ can any way: $f(x)=(f^\sharp)^\flat(x)= f^\sharp(x)$.

\item If $f:A\to B$ and $g:B\to A$ are maps of coalgebras, we have
$$
g^\sharp f^\sharp = (gf)^\sharp.
$$
\end{enumerate}
\end{lemma}

\begin{proof}
Same as lemma \ref{musicallemmacoalg}.
\end{proof}

\begin{rem}
Recall the category $(U\dgAlg)^{\$}$ from remark \ref{remfactoalg}.
Lemma \ref{musicallemmaalg} says that we have a commutative diagram of categories (enriched over graded sets)
$$\xymatrix{
&& \dgAlg^{U\$}\ar[d]^{\flat}\\
\dgAlg \ar[rru]^-\sharp \ar[rr] &&(U\dgAlg)^{\$}.
}$$
\end{rem}

\subsubsection{Calculus of meta-morphisms}

\paragraph{Tensor product}

Recall from section \ref{monoidalstrengthalg} the strength of $\otimes :\dgAlg\times \dgAlg \to \dgAlg$
$$\xymatrix{
\Theta_\otimes:\{A_1,B_1\}\otimes \{A_2,B_2\} \ar[r]& \{A_1\otimes A_2,B_1\otimes A_2\}
}$$
If $f:A_1\leadsto B_1$ and $g:A_2\leadsto B_2$ are  meta-morphisms of algebras,
let us define the {\em tensor product of meta-morphisms} by  
$$\xymatrix{
f\otimes g:=\Theta_\otimes(f\otimes g):A_1 \otimes A_2 \ar@{~>}[r]& B_1 \otimes B_2.
}$$
In addition, let us put $A \otimes g: =1_A\otimes g$ and $f \otimes B:= f\otimes 1_B$ where $1_A$ and $1_B$ are the units of the bialgebras $\{A,A\}$ and $\{B,B\}$.

\bigskip

The underlying functor of the strong functor $\otimes$ is the functor $\otimes$.
This implies, for $f:C_1\to D_1$ and $g:C_2\to D_2$ two algebra maps, the relation
$$
(f\otimes g)^\sharp = f^\sharp \otimes g^\sharp.
$$

\begin{prop}\label{tensormusicalalg}
For $f:A_1\leadsto B_1$ and $g:A_2\leadsto B_2$ two meta-morphisms of algebras, we have
$$
(f\otimes g)^\flat = f^\flat \otimes g^\flat.
$$
In particular, we can reconstruct $f\otimes g$ from $f^\sharp\otimes g^\sharp$ as $(f^\sharp\otimes g^\sharp)^\flat$.
\end{prop}
\begin{proof}
$(f\otimes g)^\flat=f^\flat\otimes g^\flat$ is a consequence of the commutation of \eqref{strengthotimesalg}. 
The last computation is a consequence of $(f^\sharp)^\flat=f$.
\end{proof}

\begin{prop}\label{metatensoralg1} 
If $u:A_1\leadsto B_1$,  $f:B_1\leadsto E_1$, $v:A_2\leadsto B_2$ and $g:D_2\leadsto E_2$
are meta-morphisms of algebras, then we have 
$$
(f \otimes g)(u \otimes v)=fu \otimes gv (-1)^{|u||g|}.
$$
In particular, we have
$$
f \otimes g = (f \otimes B_2 )(A_1\otimes g)=(B_1\otimes g) (f \otimes A_2 )  (-1)^{|f||g|}.
$$
\end{prop}

\begin{proof} 
The first identity follows from the functoriality of the strong functor $\otimes :\dgAlg\times \dgAlg \to  \dgAlg$.
The second identity is a special case of the first.
\end{proof}

\begin{cor}\label{metatensoralg2} 
If $f:B_1\leadsto E_1$ and $g:B_2\leadsto E_2$ are meta-morphisms of algebras,
then we have $$(f\otimes g)(x\otimes y)=f(x)\otimes g(y) (-1)^{|g||x|}$$
for every $x\in B_1$ and $y\in B_2$. 
\end{cor}

\begin{proof}
Using proposition \ref{tensormusicalalg} we have
\begin{eqnarray*}
(f\otimes g)(x\otimes y) & = & (f\otimes g)^\flat(x\otimes y)\\
	& = & f^\flat(x) \otimes g^\flat(y) \ (-1)^{|g||x|}\\
	& = & f(x) \otimes g(y) \ (-1)^{|g||x|}.
\end{eqnarray*}
\end{proof}

\paragraph{Sweedler hom}

The strength of the Sweedler hom functor $\{-,-\}:\dgAlg^{op}\times \dgAlg\to \dgCoalg$ is the map
$$\xymatrix{
\Theta_{\{-,-\}}:\{B_1,A_1\} \otimes  \{A_2,B_2\} \ar[r]& \HOM\bigl(\{A_1,A_2\},\{B_1,B_2\} \bigr)
}$$
defined in section \ref{laxstructures}.
If $f:B_1\leadsto A_1$ and $g:A_2 \leadsto B_2$ are meta-morphisms of algebras, let us define the {\em Sweedler hom of meta-morphisms} by  
$$\xymatrix{
\{f,g\}:=\Theta_{{\{-,-\}}}(f\otimes g): \{A_1,A_2\} \ar@{~>}[r]& \{B_1,B_2\}.
}$$
In addition, let us put $\{A,g\}:=\{1_A,g\}$ and $\{f,B\}:=\{f,1_B\}$ where $1_A$ and $1_B$ are the units of the bialgebras $\{A,A\}$ and $\{B,B\}$.

\bigskip

The underlying functor of the strong functor $\{-,-\}$ is the functor $\{-,-\}$.
This implies, for $f:B_1\to A_1$ and $g:A_2\to B_2$ two algebra maps, the relation
$$
\{f,g\}^\sharp = \{f^\sharp,g^\sharp\}.
$$

\begin{prop}\label{redmorphismalg1} 
If $u:A_2\leadsto B_2$, $v:B_1\leadsto A_1$, $f:B_2\leadsto E_2$,  and $g:E_1\leadsto B_1$
are meta-morphisms of algebras, then we have 
$$
\{vg,fu\}=\{g,f\}\{v,u\}\ (-1)^{|v|(|g|+|f|)}.
$$
In particular, we have
$$
\{v,u\}=\{(v,A_2\}\{B_1,u\}=\{B_1,u\}\{v,A_1\}\ (-1)^{|u||v|}.
$$
\end{prop}

\begin{proof}
The first identity follows from the functoriality of the strong functor  $\{-,-\}$.
The second identity is a special case of the first.
\end{proof}

\begin{prop}\label{redmorphismalg2} 
If $f:A_2\leadsto A_1$ and $g:B_1\leadsto B_2$ are meta-morphisms of algebras, then we have 
$$
\{f,g\}(h)=ghf\ (-1)^{|f|(|g|+|h|)}
$$
for every meta-morphism $h:A_1\leadsto B_1$.
In particular, the following square of graded morphisms commutes in $\dgVect$
$$\xymatrix{
\{A_1,B_1\} \ar[d]_{\Psi=\flat}\ar@^{>}[rrr]^-{\{f,g\}^\flat} &&&  \{A_2,B_2\}\ar[d]^{\Psi=\flat} \\
[A_1,B_1] \ar@^{>}[rrr]^-{[f^\flat,g^\flat]}  &&& [A_2,B_2].
}$$
In other terms, we have $\left(\{f,g\}^\flat(h)\right)^\flat = [f^\flat,g^\flat](h^\flat)$ for any $h:A_1\leadsto B_1$.
\end{prop}

\begin{proof} 
Same as proposition \ref{redmorphismcoalg2} using the diagram \eqref{strengthSHOMdiag}.
\end{proof}

\paragraph{Sweedler product}

The strength of the Sweedler product functor $\rhd:\dgCoalg\times \dgAlg\to \dgAlg$ is the map
$$\xymatrix{
\Theta_{\rhd}:\HOM(C_1,C_2)\otimes \{A_1,A_2\}\ar[r]& \{C_1\rhd A_1, C_2\rhd A_2\}
}$$
defined in section \ref{laxstructures}.
If $f:C_1\leadsto C_2$ and $g:A_1 \leadsto A_2$ are meta-morphisms of algebras, let us define the {\em Sweedler product of meta-morphisms} by  
$$\xymatrix{
f\rhd g:=\Theta_{\rhd}(f\otimes g): C_1\rhd A_1 \ar@{~>}[r]& C_2\rhd A_2.
}$$
In addition, let us put $C\rhd g:=1_C\rhd g$ and $f\rhd A:=f\rhd 1_A$ where $1_A$ and $1_C$ are the units of the bialgebras $\{A,A\}$ and $\END(C)$.

\bigskip

The underlying functor of the strong functor $\rhd$ is the functor $\rhd$.
This implies, for every meta-morphism $f:C_1\leadsto C_2$ and $g:A_1\leadsto A_2$, the relation
$$
(f\rhd g)^\sharp = f^\sharp\rhd g^\sharp.
$$

\begin{prop}\label{Sproductmeta1} 
For every meta-morphisms of coalgebras $v:C_1\leadsto C_2$ and  $w:C_2\leadsto C_3$ 
and meta-morphisms of algebras  $f:A_1\leadsto A_2$ and $g:A_2\leadsto A_3$ we have 
 $$(wv \rhd gf)=(w \rhd g)(v \rhd f)  (-1)^{|g||v|}.$$
In particular, we have
$$v \rhd f= (v \rhd A_2 )(C_1\rhd f)=(C_2\rhd f) (v \rhd A_1 )  (-1)^{|v||f|}.$$
\end{prop}

\begin{proof}
The first identity follows from the functoriality of the strong functor $\rhd$.
The second identity is a special case of the first.
\end{proof}

\begin{prop}\label{Sproductmeta2}  
If $g:A_1\leadsto A_2$ is a meta-morphism of algebras and $f:C_1\leadsto C_2$ is a meta-morphism of coalgebras,
then we have 
$$
(f\rhd g)(x\rhd y)=f(x)\rhd g(y)(-1)^{|g||x|}
$$
for every $x\in C_1$ and $y\in A_2$.
In particular, the following square of graded morphisms commutes in $\dgVect$
$$\xymatrix{
C_1\otimes  A_1\ar[d]_{\Phi}\ar@^{>}[rr]^-{f^\flat\otimes g^\flat} &&C_2\otimes A_2 \ar[d]^{\Phi} \\
C_1\rhd A_1 \ar@^{>}[rr]^-{(f\rhd g)^\flat}  && C_2\rhd A_2,
}$$
where the $\Phi$ maps are the universal measurings.
With a slight abuse of notation, we will symbolically write this property as
$$
f^\flat\rhd g^\flat = (f\rhd g)^\flat.
$$
\end{prop}

\begin{proof} 
By definition, 
\begin{eqnarray*}
(f\rhd g)(x\rhd y)&=&{\bf ev}((f\rhd g)\otimes (x\rhd y)) \\
&=&\beta(f\otimes g\otimes (x\rhd y)) \\
&=&{\bf ev}(f\otimes x)\rhd {\bf ev}(g\otimes y)(-1)^{|g||x|} \\
&=&f(x)\rhd g(y) (-1)^{|g||x|}
\end{eqnarray*}
The second assertion is a refomulation of the commutation of \eqref{strengthSproductdiag}.
\end{proof}

\paragraph{Convolution product}

The strength of the convolution product functor $[-,-]:\dgCoalg^{op}\times \dgAlg\to \dgAlg$ is the map
$$\xymatrix{
\Theta_{[-,-]}:\HOM(C_2,C_1)\otimes \{A_1,A_2\}\ar[r]& \{C_1\rhd A_1, C_2\rhd A_2\}
}$$
defined in section \ref{laxstructures}.
If $f:C_2\leadsto C_1$ and $g:A_1 \leadsto A_2$ are meta-morphisms of algebras, let us define the {\em Sweedler product of meta-morphisms} by  
$$\xymatrix{
[f, g]:=\Theta_{\rhd}(f\otimes g): [C_1,A_1] \ar@{~>}[r]& [C_2, A_2].
}$$
In addition, let us put $[C,]g:=[1_C, g]$ and $[f, A]:=[f, 1_A]$ where $1_A$ and $1_C$ are the units of the bialgebras $\{A,A\}$ and $\END(C)$.

\bigskip

The underlying functor of the strong functor $[-,-]$ is the functor $[-,-]$.
This implies, for every meta-morphism $f:C_2\leadsto C_1$ and $g:A_1\leadsto A_2$, the relation
$$
[f,g]^\sharp = [f^\sharp, g^\sharp].
$$

\begin{prop}\label{convolutionmeta1} 
If $v:C_2\leadsto C_1$ and $w:C_3\leadsto C_2$ 
are meta-morphisms of coalgebras 
and $f:A_1\leadsto A_2$ and $g:A_2\leadsto A_3$ are meta-morphisms of algebras,
then we have
$$
[vw,gf]= [w,g]  [v,f](-1)^{|v|(|w|+|g|)}.
$$
In particular, we have
$$
[v,f] = [v,A_2] [C_1,f]=[C_2,f] [v,A_1] (-1)^{|v||f|}.
$$
\end{prop}

\begin{proof}
The first identity follows from the functoriality of the strong functor  $[-,-]$.
The second identity is a special case of the first.
\end{proof}

\begin{prop}\label{convolutionmeta2} 
If $v:C_2 \leadsto C _1$ is a meta-morphism of coalgebras and $f:A_1\leadsto A_2$ is a meta-morphism of algebras, 
then for every pro-morphism $h:C_1\rhup A_1$ we have 
$$
[v,f](h)=f^\flat hv^\flat(-1)^{|v|(|f|+|h|)}.
$$
Hence we have 
$$
[v,f]^\flat=[v^\flat,f^\flat].
$$
\end{prop}

\begin{proof}
By definition, for every $x\in C_2$ we have
\begin{eqnarray*}
[v,f](h)(x)&=&\psi(v\otimes f)(h)(x)\\
&=& {\Lambda}^2(\beta)\sigma(v\otimes f)(h)(x)\\
&=& {\Lambda}^2(\beta)(f\otimes v)(h)(x) (-1)^{|v||f|}\\
&=&\beta(f\otimes h\otimes v)(x)(-1)^{|v||f|+|h||v|}\\
&=& \alpha(f\otimes h\otimes v \otimes x)(-1)^{|v||f|+|h||v|}\\
&=& f(h(v(x))) (-1)^{|v||f|+|v||h|} \\
&=& \Psi(f)h\Psi(v)(x)(-1)^{|v||f|+|v||h|} \
\end{eqnarray*}
Thus, 
\begin{eqnarray*}
[v,f](h)&=&\Psi(f)h\Psi(v)(-1)^{|v||h|+|v||f|} \\
&=&[\Psi(v),\Psi(f)](h)
\end{eqnarray*}
and this shows that $\Psi([v,f])=[\Psi(v),\Psi(f)]$.
The second assertion is also a reformulation of the commutation of \eqref{strengthconvoldiag}.
\end{proof}

\subsubsection{Module-algebras}

\begin{defi}\label{Qalgebradefi}
A $Q$-module algebra $A$ is a monoid in the category of $Q$-modules, \ie it is the data of a $Q$-module $A$ and maps
$$
m_A:A\otimes A\to A \et e_A:\FF\to A
$$
that are $Q$-equivariant. We shall denote the category of $Q$-algebras by $Q\dgAlg$.
\end{defi}

\begin{lemma}\label{modulealgebrameasuring}
With the previous notations, $m_A$ and $e_A$ are $Q$-equivariant if and only if the map $a:Q\otimes A\to A$ is a measuring with respect to the coalgebra structure of $Q$.
\end{lemma}

\begin{proof} The maps $m_A$ and $e_A$ are $Q$-module maps iff the following diagrams commute
$$
\vcenter{\xymatrix{
Q \otimes A\otimes A \ar[d]_{\Delta \otimes A\otimes A} \ar[r]^-{Q \otimes m} & Q \otimes A \ar[ddd]^{a}\\
Q \otimes Q \otimes A\otimes A \ar[d]_{\simeq} &\\
Q \otimes A \otimes  Q \otimes A \ar[d]_{a\otimes a} &\\
A \otimes  A \ar[r]^-{m} & A
}}
\et
\vcenter{\xymatrix{
Q  \ar[rr]^-{Q\otimes e }
\ar[d]_{\epsilon} && Q\otimes A \ar[d]^{a}\\
\mathbb{F}\ar[rr]^-{e_A} && A. 
}}$$
But this exactly means that $a:Q\otimes A\to A$ is a measuring.
\end{proof}

\begin{ex}
Every algebra $A$ is a module-algebra over the bialgebra $\{A, A\}$.
The action of $\{A,A\}$ on $A$ is provided by the evaluation map $\mathbf{ev}: \{A,A\}\otimes A\to A$.
\end{ex}

\begin{ex}
If $A$ is an algebra and $C$ is a coalgebra, then 
the algebra $C\rhd A$ has the structure of a module-algebra over the algebra $\END(C)\otimes \{A,A\}$.
and the algebra $[C,A]$ has the structure of a module-algebra over the algebra $\END(C)^o\otimes \{A,A\}$.
The actions are given by the strength of the functor $\rhd$ and $[-,-]$.
\end{ex}

\begin{ex}
If $A$ and $B$ are two algebras, then 
the algebra $A\otimes B$ has the structure of a module-algebra over the algebra $\{A,A\}\otimes \{B,B\}$.
The action is given by the strength of the functor $\otimes$.
\end{ex}

\begin{ex}
A module-algebra over the bialgebra $\FF[\delta]$ of example \ref{shufflehopf} 
is an algebra $A$ equipped with a derivation $\delta_A:A\to A$ of degree $|\delta|$.
\end{ex}

\medskip

\begin{defi}
Let $Q$ be a bialgebra and $A$ an algebra. 
A {\em meta-action} (or simply an {\em action}) of $Q$ on $A$ is a map of algebras $a:Q\rhd A\tto A$
such that the following squares are commutative
$$\vcenter{
\xymatrix{
(Q\otimes Q)\rhd A \ar[d]_{m_Q\otimes C} \ar@{=}[r] & Q\rhd (Q\rhd A) \ar[rr]^-{Q\rhd a}&& Q\otimes C\ar[d]^a\\
Q\otimes C\ar[rrr]^-{a} &&& C
}}
\et
\vcenter{
\xymatrix{
\FF\rhd A \ar@{=}[rrd]\ar[rr]^-{e_Q\rhd A} &&Q\rhd A\ar[d]^a\\
&& A.
}}$$
\end{defi}

Recall from lemma \ref{compocomes} that, for any coalgebra $C$, the reduction mapping
$$
\Psi:\END(C)\to [C,C]
$$
is a map of algebras.

\begin{prop} \label{hombialgalg}
For any coalgebra $A$, the reduction mapping
$$
\Psi:\{A,A\}\to [A,A]
$$
is a map of algebras.
Moreover, if $Q$ is a bialgebra, any comorphism $f:Q\to [A,A]$ which is an algebra map, lift to a unique bialgebra map $\phi:Q\to \{A,A\}$ such that $\Psi\phi=f$
$$\xymatrix{
&&\{A,A\}\ar[d]^\Psi\\
Q\ar[rr]_-f\ar[rru]^-\phi &&[A,A].
}$$
In particular composition with $\Psi$ provides a bijection between bialgebras maps $Q\to \{A,A\}$ and comeasurings $Q\to[A,A]$ that are algebra maps.
\end{prop}

\begin{proof}
Same as corollary \ref{hombialg} using  lemma \ref{compocomes}.
\end{proof}

\medskip

If $Q$ is a bialgebra, then the endo-functor $ Q\rhd (-)$ of the category $\dgAlg$
has the structure of a monad. The multiplication of the monad $ Q\rhd (-)$
is given by the map $\mu \rhd A:Q\rhd (Q\rhd A)=(Q\otimes Q)\rhd A\to Q\rhd A$
and the unit by the map $e\rhd A:A=\FF\rhd A\to Q\rhd A$. 
The endo-functor $[Q,-]$ of the category $\dgAlg$ is right adjoint to the
endo-functor $ Q\rhd (-)$ by theorem \ref{enrichmentalgcoal}.  
Hence the endo-functor $[Q,-]$ has the structure of a comonad, since the endo-functor 
$ Q\rhd (-)$ has the structure of a monad. The comultiplication of the comonad is given
by the map $[m_Q,A]:[Q,A]\to [Q\otimes Q, A]=[Q,[Q,A]]$
and its counit by the map $[e,A]:[Q,A]\to [\FF,A]=A$.

A map $a:Q\rhd A\to A$ is an action of the monad $ Q\rhd (-)$
iff the map $\Lambda^1(a):A\to [Q, A]$ is a coaction of the comonad $[Q,-]$.
The category of algebras over the monad $ Q\rhd (-)$  is equivalent to the
category of coalgebras over the comonad $[Q,-]$.

\begin{prop}\label{caracgammaalg}
If $A$ is a graded algebra, then the following data are equivalent:
\begin{enumerate}
\item a meta-action $Q\rhd A\to A$ of $Q$ on $A$;
\item an action $Q\rhd A\to A$ of the monad $Q\rhd (-)$;
\item a coaction $A\to [Q,A]$ of the comonad $[Q,-]$;
\item the structure of a $Q$-module algebra on $A$;
\item a map of algebras $\pi:Q\to [A,A]$ which is a comeasuring;
\item a map of bialgebras $\pi:Q\to \{A,A\}$.
\end{enumerate}
\end{prop}

\begin{proof}
The equivalence between (1) and (2) is a general fact true in any closed category.
The equivalence between (2) and (3) is the remark above.
The equivalence between (1) and (4) is lemma \ref{modulealgebrameasuring}.
The equivalence between (5) and (6) is proposition \ref{hombialgalg}.
Finally, the equivalence between (4) and (5) is the fact that a $Q$-module structure is equivalent to an algebra map $Q\to [A,A]$ and the definition of comeasurings.
\end{proof}

\begin{prop}\label{modulealgebramonad}
The forgetful functor $Q\dgAlg\to \dgAlg$ has a left adjoint $A\to Q\rhd A$ 
and a right adjoint $A\to [Q, A]$.
In particular, limits and colimits exist in $Q\dgAlg$ and can be computed in $\dgAlg$.
\end{prop}

\begin{proof} This follows from the general theory of monads and comonads.
\end{proof}

\bigskip

Let now $Q$ be a cocommutative Hopf algebra.
Then, the category $Q\dgAlg$ is enriched, tensored and cotensored over the category $\dgCoalg$. 
Moreover, we shall prove that $Q\dgAlg$ is in fact enriched, tensored and cotensored over the category $Q\dgCoalg$. 
Let $C$ be a $Q$-module coalgebra and $A$ a $Q$-module algebra, 
the tensor product of $A$ by $C$ is the algebra $C\rhd A$ and the cotensor is $[C,A]$.

The actions of $Q$ on $C\rhd A$ and $[C,A]$ are defined by 
$$\xymatrix{
Q \ar[r]^-{\Delta} & Q \otimes Q \ar[rr]^-{\pi_C\otimes \pi_A} && \END(C)\otimes \{A,A\} \ar[r]^-{\Theta_\rhd} &  \{C\rhd A, C\rhd A\}
}$$
and
$$\xymatrix{
Q \ar[r]^-{\Delta} & Q \otimes Q\ar[r]^-{S\otimes Q} & Q^o  \otimes Q \ar[rr]^-{\pi_C\otimes \pi_A} && \END(C)^o \otimes \{A,A\} \ar[r]^-{\Theta_{[-,-]}} &  \{[C,A],[C,A]\}
}$$
where the $\pi$s are the meta-actions of $Q$ and $\Theta_\otimes$ and $\Theta_\HOM$ are the strength of $\otimes$ and $\HOM$ in $\dgCoalg$ (see section \ref{laxstructures}).
For $A$ and $B$ two $Q$-module algebras, the enrichment of $Q\dgAlg$ over itself is given by $\{A,B\}$ with the action
$$\xymatrix{
Q \ar[r]^-{\Delta} & Q \otimes Q\ar[r]^-{S\otimes Q} & Q^o  \otimes Q \ar[rr]^-{\pi_A\otimes \pi_B} && \{A,A\}^o \otimes \{B,B\} \ar[r]^-{\Theta_{\{-,-\}}} &  \END(\{A,A\})
}$$
where $\Theta_{\{-,-\}}$ is the strength of $\{-,-\}$ defined in section \ref{laxstructures}.

In the case where the action of $Q$ on $C$ is the trivial action $\epsilon_Q\otimes C:Q\otimes C\to C$, these formulas specializes to give the tensor and cotensor of $Q\dgCoalg$ over $\dgCoalg$. Then, as in the end of section \ref{moduleoverhopfalgebra} the enrichment of $Q\dgCoalg$ over $\dgCoalg$ is defined as the equalizer in $\dgCoalg$
$$\xymatrix{
\{A,B\}_Q \ar[r]& \{A,B\}=\HOM(\FF,\{A,B\})\ar@<.6ex>[rrr]^-{\HOM(\epsilon_Q,\{A,B\})}\ar@<-.6ex>[rrr]_-{\Lambda^1a} &&& \HOM(Q,\{A,B\})
}$$
where $a:Q\otimes \{A,B\}\to \{A,B\}$ is the above action.

\medskip
The forgetful functor $U:Q\dgAlg\to \dgAlg$ is strong and  the adjunctions $Q\rhd (-)\dashv U\dashv [Q,-]$ of proposition \ref{modulealgebramonad} are strong. We leave the details to the reader, they are similar to what was done in section \ref{moduleoverhopfalgebra}.
Moreover we have the promised stronger result, which is a generalisation of theorem \ref{enrichmentalgcoal}.

\begin{thm}\label{Qenrichmentalgcoal}
The category $Q\dgAlg$ is enriched, bicomplete and monoidal over $Q\dgCoalg$.
The forgetful functor $U:Q\dgAlg\to \dgAlg$ is symmetric monoidal and preserves all Sweedler operations.
\end{thm}
\begin{proof}
To prove the first statement, we need to prove that the strong natural transformations defining the adjunctions between $\rhd$, $[-,-]$ and $\{-,-\}$ are $Q$-equivariant and that the strength of the monoidal structure is $Q$-equivariant.

Let us prove that the strong natural transformation $\Lambda=\Lambda^2:\{C\rhd A,B\}\simeq \HOM(C,\{A,B\})$ is $Q$-equivariant.
We consider the square 
$$\xymatrix{
\END(C)^o \otimes \{A,A\}^o \otimes \{B,B\} \ar[r]\ar[d] &\END(C)^o \otimes \END(\{A,B\})\ar[d]\\
\{C\rhd A,C\rhd A]^o  \otimes \{B,B\}\ar[r]& \{C\rhd A,B\}\simeq \HOM(C,\{A,B\}).
}$$
stating that $\Lambda$ is a strong natural transformation. Then, as in the end of section \ref{moduleoverhopfalgebra}, the proof that $\Lambda$ is $Q$-equivariant is nothing but the commutation of the diagram where we have precompose this square with the map $(S\otimes S\otimes Q)\circ \Delta^{(3)}:Q\to Q^o \otimes Q^o \otimes Q$.
The other equivariances are proven the same way from the corresponding strong natural transformations. We leave the details to the reader.

The second statement is obvious by construction.
\end{proof}

\bigskip
As for theorem \ref{Qhomcoalg} and theorem \ref{dgtoghomcoalg}, it should be clear to the reader that the same result holds if we work in the graded setting instead of the differential graded setting. We can then apply theorem \ref{Qenrichmentalgcoal} to deduce the following important result.

\begin{thm}\label{dgtoghomalg}
The category $\gAlg$ is enriched and bicomplete over $\gCoalg$.
The forgetful functor $U_d:\dgAlg\to \gAlg$ is symmetric monoidal, 
it has left adjoint $\d \rhd-$ and right adjoint $[\d,-]$,
it preserves the enrichment and commute with tensors and cotensors.
\end{thm}
\begin{proof}
The proof of the first statement is analog of that of theorem \ref{enrichmentalgcoal}.
Then we proceed as in theorem \ref{dgtoghomcoalg}. We describe $\dgVect$ as $\Mod(Q)$ in $\gVect$ for the cocommutative Hopf algebra $Q=\d=\FF\delta_+$ and apply proposition \ref{modulealgebramonad} and (the graded analog of) 
theorem \ref{Qenrichmentalgcoal}. We only have to check that the enrichment, tensor and cotensors of $\dgAlg$ constructed from those of $\gAlg$ from theorem \ref{Qenrichmentalgcoal} coincide with those of theorem \ref{enrichmentalgcoal}. 
By adjunction it is enough to check the coincidence of the cotensor only, but it is easy to see that the differential in $[C,A]$ is the same computed in $\dgAlg$ or in $Q\gAlg$ for $Q=\FF\delta_+$. A similar computation proves the result about the monoidal structure.
\end{proof}

This theorem says in particular that to compute $C\rhd A$ (or $\{A,B\}$) in the dg-context, we can first compute them in the graded context and there will be a unique differential induced by that of $C$ and $A$ (or $A$ and $B$) enhancing it into the dg-Sweedler product (or the dg-Sweedler hom).
We shall detail how to compute these differentials in the next section.

\subsubsection{Primitive meta-morphisms}

For an algebra map $f:A\to B$, recall that $\Der(f)$ and $\Prim_f(\{A,B\})$ are respectively the dg-vector spaces of $f$-derivations and $f$-primitive elements.

\begin{prop}\label{primideriv}
Let $f:A\to B$ be an algebra map, then there are natural bijections between
\begin{enumerate}
\item maps of vector spaces $p:X\to \Prim_f(\{A,B\})$;
\item maps of pointed coalgebras $k:T^c_{\bullet,1}(X)\to (\{A,B\},f)$;
\item maps of algebras $g:T^c_{\bullet,1}(X)\rhd A\to B$ such that $f=g(e\rhd A):A\to T^c_{\bullet,1}(X)\rhd A\to B$;
\item maps of algebras $j:A\to [T^c_{\bullet,1}(X),B]$ such that $f=[e,B]j:A\to [T^c_{\bullet,1}(X),B]\to [\FF,B]=B$;
\item $f$-derivations $d:A\to [X,B]$;
\item maps of $B$-bimodules $\Omega_{A,f}\to [X,B]$;
\item linear maps $h:X\to \Der(f)=\hom_{B,B}(\Omega_{A,f},B)$.
\end{enumerate}
If $A=B$ and $f=id_A$, there are more natural bijections with 
\begin{enumerate}
\setcounter{enumi}{7}
\item the maps of bialgebras $T^{csh}(X)\to \{A,A\}$;
\item the meta-actions $T^{csh}(X)\rhd A\to A$.
\end{enumerate}
\end{prop}
\begin{proof}
If $k$ decomposes into $k_0+k_1$ with respect to the decomposition $T^c_{\bullet,1}(X)=\FF\oplus X$.
By assumption $k_0:\FF\to \{A,B\}$ is the atom corresponding to $f$.
The bijection $1\leftrightarrow 2$ is from proposition \ref{univpropofdeltanew} and identifies $k_1=p$.
The bijection $2\leftrightarrow 3$ is by adjunction, we have $k=\Lambda^2g$ and $g={\bf ev}(k\rhd A)$.
The bijection $3\leftrightarrow 4$ is also by adjunction, we have $j=\lambda^2gu$ where $u:T^c_{\bullet,1}(X)\otimes A\to T^c_{\bullet,1}(X)\rhd A$ is the universal measuring.
The algebra $[T^c_{\bullet,1}(X),B]=B\oplus [X,B]$ is of the type $B\oplus M$ where $M$ is a $B$-bicomodule.
Hence, an algebra map $g:A\to [T^c_{\bullet,1}(X),B]$ decomposes into an algebra map $f: A\to B$ and a $f$-derivation $d:A\to [X,B]$. We have $k_1=\Lambda^2d$ and $d={\bf ev}(k_1\otimes D)$. This proves the bijection $4\leftrightarrow 5$. 
The bijection $5\leftrightarrow 6$ is the definition of $\Omega_{A,f}$.
Recall that the category of bimodules is cotensored over $\dgVect$, hence $B$-bimodules $\Omega_{A,f}\to [X,B]$ are in bijection with linear maps $B$-bimodules $X\to \hom_{B,B}(\Omega_{A,f},B)=\Der(A,f;B) = \Der(f)$.
This proves the bijection $6\leftrightarrow 7$. Remark that the bijection $4\leftrightarrow 6$ is given by $h=\lambda^2 d$ and $d=ev(h\otimes D)$.
Finally, the bijections $2\leftrightarrow 8\leftrightarrow 9$ are from proposition \ref{univpropofdeltabialg} and proposition \ref{caracgammacoalg}.
\end{proof}

\begin{cor}\label{corprimideriv}
Let $f:A\to B$ be an algebra map, the reduction map $\Psi:\{A,B\}\to [A,B]$ induces an isomorphism in $\dgVect$
$$
\Prim_f(\{A,B\}) = \Der(f).
$$
\end{cor}
\begin{proof}
By proposition \ref{primideriv}, both objects have the same functor of points. They are isomorphic by Yoneda's lemma.
Let us prove that the isomorphism is induced by $\Psi$.
According to the proof of proposition \ref{primideriv}, a map $p:X\to \Prim_f(\{A,B\})$ is send to the coderivation
$h=\lambda^2{\bf ev}u(p\otimes A) = \Psi (p)$ by definition of $\Psi$.
\end{proof}

\medskip

\begin{nota}\label{notationatomderiv}
Recall from lemma \ref{musicallemmaalg} that if $A$ and $B$ are algebras, the map $\Psi: \{A,B\}  \to [A,B]$
induces a bijection between the atoms of the coalgebras $\{A,B\}$ and the maps of coalgebras $A\to B$.
If $f:A\to B$ is a map of algebras, we noted $f^\sharp \in \HOM(D,C)$ the unique atom such that $\Psi(f^\sharp)=f$.

Corollary \ref{corprimideriv} says that, if $d:A\rhup_n B$ is a $f$-derivation, there exists a unique  element $b\in \{A,B\}_n$ primitive with respect to $f^\sharp$ such that $\Psi(b)=d$. We shall denote by $d^\sharp$ this element.

Recall the notation $\Psi(b)=b^\flat$ from section \ref{reducetmeta-morphismalg}. We have $(d^\sharp)^\flat=d$ by definition.
\end{nota}

With these notations we can explain proposition \ref{primideriv} by the commutative diagram in $\dgVect$
$$\xymatrix{
\Prim_f(\{A,B\})\ar[r]\ar[d]^\flat  & \{A,B\} \ar[d]^{\Psi=\flat}\\
\Der(f) \ar[r] \ar@/^1pc/[u]^-\sharp & [A,B].
}$$
The maps $\sharp$ and $\flat$ induce inverse bijections $\Prim_f(\{A,B\}) = \Der(f)$.

\bigskip

If $C$ is a coalgebra, then the coalgebra $\END(C)$ has the structure of a bialgebra.
$C$ is a module-coalgebra over the bialgebra $\END(C)$.
The action of this is given by the strong evaluation $\mathbf{ev}: \END(C) \otimes C\to C$.

\begin{thm}\label{primiLiederiv}
If $A$ is an algebra, we have a commutative diagram in $\dgLie$
$$\xymatrix{
\Prim(\{A,A\})\ar[r]\ar[d]^\flat  & \{A,A\} \ar[d]^{\Psi=\flat}\\
\Der(A) \ar[r] \ar@/^1pc/[u]^-\sharp& [A,A].
}$$
In particular, the maps $\sharp$ and $\Psi(=\flat)$ induce inverse Lie algebra isomorphisms 
$$
\flat :\Prim(\{A,A\})\simeq\Der(A):\sharp
$$
which preserve the square of odd elements.
\end{thm}

\begin{proof}
Same as in theorem \ref{primiLiecoderiv}.
\end{proof}

\bigskip

Recall from proposition \ref{caracgammaalg} that a $Q$-module-algebra is an algebra $A$ 
equipped with a left $Q$-module structure defined by an action $a:Q\otimes A\to A$ which is a measuring.

\begin{lemma} \label{primitiveactionalgebras} 
Let $A$ be a $Q$-module-algebra.
If $b\in Q$ is primitive, then the map $\pi(b):=b\cdot(-):A\to A$ is a derivation of $A$.
Moreover the map $\pi:\Prim(Q)\to \Der(A)$ so defined is a homomorphism of Lie algebras which preserves the square of odd elements.
\end{lemma}

\begin{proof} 
Same as in lemma \ref{primitiveactioncoalgebras}.
\end{proof}

The following proposition says that the map $\pi:\Prim(Q)\to \Der(A)$ can be computed either from the action map or from the meta-action map.

\begin{prop}\label{Qprim=Qderiv}
The commutative triangle of algebras
$$\xymatrix{
&& \{A,A\}\ar[d]^\Psi\\
Q\ar[rr]^-\beta\ar[rru]^-\alpha&& [A,A]
}$$
induces a commutative triangle of Lie dg-algebras morphisms preserving the squares of odd elements
$$\xymatrix{
&& \Prim(\{A,A\})\ar[d]^\Psi_\simeq\\
\Prim(Q)\ar[rr]^-{\beta'}\ar[rru]^-{\alpha'}&& \Der(A)
}$$
\end{prop}
\begin{proof}
Same as in proposition \ref{Qprim=Qcoderiv}.
\end{proof}

\subsubsection{Derivative of Sweedler operations}\label{derivativesweedleralg}

Let $A$ and $B$ be two algebras and $C$ be a coalgebra.
Recall that the strengths of the Sweedler operations $\otimes$, $\{-,-\}$, $\rhd$ and $[-,-]$ give bialgebra maps
$$\xymatrix{
\Theta_\otimes:\{A,A\}\otimes  \{B,B\} \ar[r]& \{A\otimes B,A\otimes B\},
}$$
$$\xymatrix{
\Theta_{\{-,-\}}:\{A,A\}^o \otimes  \{B,B\} \ar[r]& \END(\{A,B\}),
}$$
$$\xymatrix{
\Theta_{\rhd}:\END(C)\otimes \{A,A\} \ar[r]& \{C\rhd A,C\rhd A\}),
}$$
$$\xymatrix{
\Theta_{[-,-]}:\END(C)^o\otimes \{A,A\} \ar[r]& \{[C, A],[C,A]\}.
}$$
By lemma \ref{derivationmorbialg}, we have the derivative maps between the Lie algebras of primitive elements
$$\xymatrix{
\Theta'_\otimes:\Prim(\{A,A\})\times  \Prim(\{B,B\}) \ar[r]& \Prim(\{A\otimes B,A\otimes B\}),
}$$
$$\xymatrix{
\Theta'_{\{-,-\}}:\Prim(\{A,A\}^o) \times  \Prim(\{B,B\}) \ar[r]& \Prim(\END(\{A,B\})),
}$$
$$\xymatrix{
\Theta'_{\rhd}:\Prim(\END(C))\times \Prim(\{A,A\}) \ar[r]& \Prim(\{C\rhd A,C\rhd A\})),
}$$
$$\xymatrix{
\Theta'_{[-,-]}:\Prim(\END(C)^o)\times \Prim(\{A,A\}) \ar[r]& \Prim(\{[C, A],[C,A]\}).
}$$
By theorems \ref{primiLiecoderiv} and \ref{primiLiederiv}, these are equivalent to maps of Lie algebras
$$\xymatrix{
\Theta'_\otimes:\Der(A)\times  \Der(B) \ar[r]& \Der(A\otimes B),
}$$
$$\xymatrix{
\Theta'_{\{-,-\}}:\Der(A) \times  \Der(B) \ar[r]& \Coder(\{A,B\}),
}$$
$$\xymatrix{
\Theta'_{\rhd}:\Coder(C)\times \Prim(A) \ar[r]& \Der(C\rhd A),
}$$
$$\xymatrix{
\Theta'_{[-,-]}:\Coder(C)\times \Prim(A) \ar[r]& \Der([C, A]).
}$$
Using the calculus of meta-morphisms, they are given respectively by 
\begin{eqnarray*}
(d_1, d_2) &\mto& \big(d_1^\sharp\otimes B + A\otimes d_2^\sharp\big)^\flat,\\
(d_1, d_2) &\mto& \big(\{A, d_2^\sharp\} - \{d_1^\sharp, B\} \big)^\flat,\\
(d_1, d_2) &\mto& \big(d_1^\sharp\rhd A + C\rhd d_2^\sharp\big)^\flat,\\
(d_1, d_2) &\mto& \big(\hom(C, d_2^\sharp) - \hom(d_1^\sharp, A) \big)^\flat.
\end{eqnarray*}

The calculus of meta-morphisms also tells us that (proposition \ref{tensormusicalalg})
$$
\big(d_1^\sharp\otimes B + A\otimes d_2^\sharp\big)^\flat = (d_1^\sharp)^\flat\otimes B + A\otimes (d_2^\sharp)^\flat = d_1\otimes B + A\otimes d_2,
$$
and that (proposition \ref{convolutionmeta2})
$$
\big([C,d_1^\sharp]-[d_2^\sharp,A]\big)^\flat = [C,(d_1^\sharp)^\flat]-[(d_2^\sharp)^\flat,A] = [C,d_1]-[d_2,A].
$$
In other words, the derivations induced through the strength of the tensor and convolution products $\otimes$ and $[-,-]$ are the classical derivations constructed on a tensor product (proposition \ref{actiondertensor}) and on a hom complex (proposition \ref{derconvolution3}).
However for the Sweedler hom and Sweedler product operations, the computation is new. 
The following results detail how to deal with them.

\paragraph{Sweedler hom}
\medskip
Recall from proposition \ref{redmorphismalg2}, that if $f\in \{A,A\}$ and $g\in \{B,B\}$, we have a commutative square of graded morphisms in $\dgVect$
$$\xymatrix{
\{A,B\} \ar[d]_{\Psi=\flat}\ar@^{>}[rrr]^-{\{f,g\}^\flat} &&&  \{A,B\}\ar[d]^{\Psi=\flat} \\
[A,B] \ar@^{>}[rrr]^-{\hom(f^\flat,g^\flat)}  &&& [A,B].
}$$
The graded morphism $\HOM(f,g)^\flat$ is not in general uniquely determined by $\hom(f^\flat,g^\flat)$. 
The following proposition proves that this is somehow the case for primitive elements.

\begin{prop}\label{uniqueextensionderSHOM}
If $d_1^\sharp$ and $d_2^\sharp$ are the primitive elements of $\{A,A\}$ and $\{B,B\}$ associated to derivations $d_1$ and $d_2$ of $A$ and $B$ then
$\{A,d_2^\sharp\}^\flat$ and $\{d_1^\sharp,B\}^\flat$ are coderivations and 
$d=\{A,d_2^\sharp\}^\flat - \{d_1^\sharp,B\}^\flat$ is the unique coderivation such that the square
$$\xymatrix{
\{A,B\} \ar[d]_{\Psi=\flat}\ar@^{>}[rrrr]^-{d} &&&&  \{A,B\}\ar[d]^{\Psi=\flat} \\
[A,B] \ar@^{>}[rrrr]^-{\hom(A,d_2)-\hom(d_1,B)}  &&&& [A,B].
}$$
commutes.
Equivalently, $d$ is the unique coderivation such that
$$
d(h)^\flat = d_2h^\flat-h^\flat d_1 \ (-1)^{|h||d_1|}
$$
for any $h\in \{A,B\}$.
\end{prop}

\begin{proof}
If $d_1$ and $d_2$ are derivations, then $\{A,d_2^\sharp\}^\flat$ and $\{d_1^\sharp,B\}^\flat$ are coderivation by theorem \ref{primiLiederiv}. 
The commutation of the diagram is from proposition \ref{redmorphismalg2} and the fact that $d_i=(d_i^\sharp)^\flat$.
The coderivation $\{A,d_2\}^\flat-\{d_1,B\}^\flat$ is equivalent to a coalgebra map $\{A,B\}\oplus S^{-n}\{A,B\}\to \{A,B\}$ by proposition \ref{fderiv} ($n$ is the degree of $d_1$ and $d_2$). Then, the unicity result follows by the separation property of $\Psi:\{A,B\}\to [A,B]$.
\end{proof}

If $A=T(X)$ is free, we have the following strengthening of the previous result.

\begin{prop}\label{uniqueextensionderSHOMfree}
If $d_1^\sharp$ and $d_2^\sharp$ are the primitive elements of $\{T(X),T(X)\}$ and $\{B,B\}$ associated to derivations $d_1$ and $d_2$ of $T(X)$ and $B$ then
the coderivation $d=\{T(X), d_2^\sharp\}^\flat - \{d_1^\sharp, B\}^\flat$ is the unique coderivation on $\{T(X),B\}=T^\vee([X,B])$ such the following square commutes
$$\xymatrix{
T^\vee([X,B]) \ar[d]_{\Psi}\ar@^{>}[rrrr]^-{d} &&&&  T^\vee([X,B])\ar[d]^{q} \\
[T(X),B] \ar@^{>}[rrrr]^-{\hom(i, d_2) - \hom(d_1i, B)}  &&&& [X,B].
}$$
where $i:X\to T(X)$ is the generating map and $q:T^\vee([X,B])\to [X,B]$ is the cogenerating map.
Equivalently, $d$ is the unique coderivation such that
$$
q(d(h)) = d_2h^\flat i-h^\flat d_1i \ (-1)^{|h||d_1|}
$$
for any $h\in \{T(X),B\}$.
\end{prop}

\begin{proof}
We have a commutative diagram
$$\xymatrix{
\{T(X),B\} \ar[d]_{\Psi}\ar@^{>}[rrrr]^-{\{T(X), d_2^\sharp\}^\flat - \{d_1^\sharp, B\}^\flat} &&&&  \{T(X),B\}\ar[d]^{\Psi} \ar@{=}[r]&T^\vee([X,B])\ar[d]^q\\
[T(X),B] \ar@^{>}[rrrr]^-{\hom(T(X), d_2) - \hom(d_1, B)}  &&&& [T(X),B] \ar[r]^-{[C,p]}& [X,B].
}$$
Then the proof is the same as in lemma \ref{uniqueextensionderSHOM} but using the separating property of $q$ instead of that of $\Psi$.
\end{proof}

\paragraph{Sweedler product}
If $g:A\leadsto B$ is a meta-morphism of algebras and $f:C\leadsto D$ is a meta-morphism of coalgebras,
recall from proposition \ref{Sproductmeta2} that we have a commutative square
$$\xymatrix{
C\otimes  A\ar[d]_{\Phi}\ar@^{>}[rr]^-{f^\flat\otimes g^\flat} && D \otimes B  \ar[d]^{\Phi} \\
C\rhd A \ar@^{>}[rr]^-{(f\rhd g)^\flat}  && D\rhd B.
}$$

\begin{prop}\label{uniqueextensionderSproduct}
If $d_1^\sharp$ and $d_2^\sharp$ are the primitive elements of $\END(C)$ and $\{A,A\}$ associated to a coderivation $d_1$ of $C$ and a derivation $d_2$ of $A$ then
$(d_1^\sharp\rhd A)^\flat$ and $(C\rhd d_2^\sharp)^\flat$ are derivations of $C\rhd A$ and 
$d=(d_1^\sharp\rhd A)^\flat + (C\rhd d_2^\sharp)^\flat$ is the unique derivation of $C\rhd A$ such that the square
$$\xymatrix{
C\otimes  A\ar[d]_{\Phi}\ar@^{>}[rrr]^-{d_1\otimes A+ C\otimes d_2} &&& C \otimes A  \ar[d]^{\Phi} \\
C\rhd A \ar@^{>}[rrr]^-{d}  &&& C\rhd A.
}$$
commutes.
Equivalently, $d$ is the unique coderivation such that
$$
d(c\rhd a) = (d_1c)\rhd a + c\rhd (d_2a)\ (-1)^{|c||d_2|}
$$
for any $c\otimes a\in C\otimes A$.
\end{prop}

\begin{proof}
If $d_1$ and $d_2$ are as in the statement, then 
$(d_1^\sharp\rhd A)^\flat$ and $(C\rhd d_2^\sharp)^\flat$ are derivations of $C\rhd A$ by theorem \ref{primiLiederiv}.
The commutation of the diagram is from proposition \ref{Sproductmeta2} and the fact that $d_i=(d_i^\sharp)^\flat$.
The derivation 
$(d_1^\sharp\rhd A)^\flat + (C\rhd d_2^\sharp)^\flat$
is equivalent to an algebra map $C\rhd A\oplus S^{-n}(C\rhd A)\to C\rhd A$ by proposition \ref{fderiv} ($n$ is the degree of $d_1$ and $d_2$). Then, the unicity result follows by the separation property of the universal measuring $\Phi:C\otimes A\to C\rhd A$.
\end{proof}

In $A=T(X)$ is a free algebra, we have the following strengthening of the previous result.

\begin{prop}\label{uniqueextensionderSproductfree}
If $d_1^\sharp$ and $d_2^\sharp$ are the primitive elements of $\END(C)$ and $\{T(X),T(X)\}$ associated to a coderivation $d_1$ of $C$ and a derivation $d_2$ of $T(X)$ then
$d=(d_1^\sharp\rhd T(X))^\flat + (C\rhd d_2^\sharp)^\flat$ is the unique derivation of $C\rhd T(X)=T(C\otimes X)$ such that the square
$$\xymatrix{
C\otimes  X\ar[d]_{j}\ar@^{>}[rrr]^-{d_1\otimes i+ C\otimes d_2i} &&& C \otimes T(X)  \ar[d]^{\Phi} \\
T(C\otimes X) \ar@^{>}[rrr]^-{d}  &&& T(C\otimes X).
}$$
commutes ($j:C\otimes X\to T(C\otimes X)$ is the generating map).
Equivalently, $d$ is the unique coderivation such that
$$
d(c\rhd x) = d_1(c)\rhd x + c\rhd d_2(x)\ (-1)^{|c||d_2|}
$$
for any $c\otimes x\in C\otimes X$.
\end{prop}

\begin{proof}
We have a commutative diagram
$$\xymatrix{
C\otimes X\ar[d]_j\ar[r]^-{C\otimes i}&C\otimes  T(X)\ar[d]_{\Phi}\ar@^{>}[rrr]^-{d_1\otimes T(X)+ C\otimes d_2} &&& C \otimes T(X)  \ar[d]^{\Phi} \\
T(C\otimes X) \ar@{=}[r] & C\rhd T(X) \ar@^{>}[rrr]^-{d}  &&& C\rhd T(X).
}$$
Then the proof is the same as in lemma \ref{uniqueextensionderSproduct} but using the separating property of $j$ instead of that of $\Phi$.
\end{proof}

\subsubsection{Strong monadicity}

We finish this section by a strengthening of theorem \ref{monadicalg}.

\begin{thm}\label{strongmonadicityalg}
The adjunction $T\dashv U$ enriches into a strong colax monoidal adjunction
$$\xymatrix{
T: \dgVect^{T^\vee \$} \ar@<.6ex>[r]& \dgAlg^\$:U. \ar@<.6ex>[l]
}$$
Moreover the adjunction is strongly monadic.
\end{thm}

\begin{proof}

Let us prove first that $T$ and $U$ are strong functors.
Both categories are bicomplete because of proposition \ref{vectbicompletecoalg}, so we can use proposition \ref{triplestrength} to describe the strength of $T$ and $U$ as (co)lax structures.
The colax structure of $U$ is given by the isomorphism
$$\xymatrix{
U[C,A]\ar[r]^-\simeq & [UC,UA]
}$$
The conditions of \ref{colaxmodularstructure} are equivalent to the pentagon and the unit identities for $[-,-]$ in $\dgVect$.

The lax structure of $T$ is given by the isomorphism
$$\xymatrix{
C\rhd TX \ar[r]^-\simeq & T(UC\otimes X)
}$$
of proposition \ref{Sweedleroftensor}.
The unit condition of \ref{laxmodularstructure} is clear and the associativitiy condition is equivalent to the commutation of
$$\xymatrix{
(D\otimes C)\rhd TX \ar@{=}[r] &  D\rhd (C\rhd TX) \ar[r]^-\simeq &D\rhd T(UC\otimes X) \\
(UD\otimes UC)\otimes X\ar[u]  && \ar@{=}[ll]  UD\otimes (UC\otimes X)\ar[u]
}$$
which is an easy consequence of how the isomorphism $T(UC\otimes X)\simeq C\rhd T(X)$ is constructed.
Finally, the natural isomorphism $\{TX,A\}\simeq T^\vee [X,UA]$ of proposition \ref{exSweedlerhom} is the strength of the adjunction $T\dashv U$.

The argument to prove that $U$ is strong symmetric monoidal (and thus that $T$ is strong symmetric colax)
is similar to that of theorem \ref{strongcomonadicitycoalg} and we do not reproduce it.

The proof of the strong monadicity is analogous to the proof of the strong comonadicity in theorem \ref{strongcomonadicitycoalg}.
We need only to replace proposition \ref{stronglimcoalg} by proposition \ref{stronglimalg}.
\end{proof}

\begin{rem}\label{remreduction2}
The functors $U:\dgAlg^{U\$}\to \dgVect^\$$ and $U:\dgAlg^{\$}\to \dgVect^{T^\vee\$}$ corresponds to each other through the adjunction of proposition \ref{transferadjunction}.
\end{rem}

\begin{rem}\label{calculSHOM}
As for coalgebras in remark \ref{calculHOM}, the strong monadicity of $T\dashv U$ over $\dgCoalg$ leads to a new computation of the coalgebras $\{A,B\}$.
The functor $U:\dgAlg^\$\to \dgVect^{T^\vee\$}$ is strongly faithful and again this subsumes our proofs by reduction (see section \ref{proofbyreduction}).
The coalgebras $\{A,B\}$ can be computed via the (enriched version) of the usual computation of hom of algebras. 
In the same way that we have the common equalizer in $\Set$ (see appendix \ref{commonequalizer})

$$\xymatrix{
\dgAlg(A,B)\ar[r]&\dgVect(A,B) \ar@<.6ex>[d]^-{\alpha'}\ar@<-.6ex>[d]_-{\beta'}\ar@<.6ex>[r]^-{\alpha}\ar@<-.6ex>[r]_-{\beta}& \dgVect(A\otimes A,B)\\
& \dgVect(\FF,B)
}$$
(where for $f:A\to B$, we define $\alpha (f)=m_B(f\otimes f):A\otimes A\to B$, $\beta (f)=fm_A:A\otimes A\to B$,
$\alpha' (f)=1_B:\FF\to B$ and $\beta' (f)=f(1_A):\FF\to B$)
we have the equalizer in $\dgCoalg$
$$\xymatrix{
\{A,B\}\ar[r]^-{m'} &T^\vee([A,B]) \ar@<.6ex>[d]^-{T^\vee(\alpha')}\ar@<-.6ex>[d]_-{T^\vee(\beta')}\ar@<.6ex>[r]^-{T^\vee(\alpha)}\ar@<-.6ex>[r]_-{T^\vee(\beta)}& T^\vee([A\otimes A,B])\\
& T^\vee([\FF,B]).
}$$
where $m'$ is the map from corollary \ref{betterShom}. 
In other words, the meta-morphisms of algebras are the meta-morphisms of vector spaces preserving the algebra structure.
\end{rem}

\newpage
\section{Other Sweedler contexts}\label{Sweedlercontexts}

In the last two sections, we have detailled what we called the Sweedler theory of dg-coalgebras and dg-algebras.
This consists in a set of six functors:
\begin{center}
\begin{tabular}{lrl}
\rule[-2ex]{0pt}{4ex} the tensor product of non-counital coalgebras & $\otimes$&$\!\!\!\!: \dgCoalg\times \dgCoalg \to \dgCoalg$, \\
\rule[-2ex]{0pt}{4ex} the non-counital coalgebra internal hom & $\HOM$&$\!\!\!\!:\dgCoalg^{op}\times \dgCoalg \to \dgCoalg$, \\
\rule[-2ex]{0pt}{4ex} the non-unital Sweedler hom & $\{-,-\}$&$\!\!\!\!:\dgAlg^{op}\times \dgAlg \to \dgCoalg$, \\
\rule[-2ex]{0pt}{4ex} the non-unital Sweedler product & $\rhd$&$\!\!\!\!:\dgCoalg\times \dgAlg \to \dgAlg$, \\
\rule[-2ex]{0pt}{4ex} the non-unital convolution product & $[-,-]$&$\!\!\!\!:\dgCoalg^{op}\times \dgAlg \to \dgAlg$, \\
\rule[-2ex]{0pt}{4ex} and the tensor product of non-unital algebras & $\otimes$&$\!\!\!\!:\dgAlg\times \dgAlg \to \dgAlg$,
\end{tabular}
\end{center}
such that
\begin{itemize}
\item the category $(\dgCoalg,\otimes, \HOM)$ is locally presentable symmetric monoidal closed and comonadic over $\dgVect$
\item and the category $(\dgAlg,\{-,-\},\rhd, [-,-],\otimes)$ is locally presentable, enriched, bicomplete and symmetric monoidal over $\dgCoalg$, and  monadic over $\dgVect$.
\end{itemize}
Together with this, we have distinguished some isomorphisms
\begin{center}
\begin{tabular}{cl}
\rule[-2ex]{0pt}{4ex} $T^\vee [C,X]\simeq \HOM(C,T^\vee (X))$ & (proposition \ref{exemplehomcoalg})\\
\rule[-2ex]{0pt}{4ex} $\{TX,A\}\simeq T^\vee([X,A])$ &  (proposition \ref{exSweedlerhom})\\
\rule[-2ex]{0pt}{4ex} $C\rhd T(X) \simeq T(C\otimes X)$ &  (proposition \ref{Sweedleroftensor})\\
\end{tabular}
\end{center}
which are interpreted as strengthening the adjunctions $U\dashv T^\vee$ and $T\dashv U$ over $\dgCoalg$:
\begin{itemize}
\item the adjunction $U:\dgCoalg \rightleftarrows \dgVect:T^\vee $ is strongly comonadic over $\dgCoalg$
\item and the adjunction $T:\dgVect \rightleftarrows \dgAlg:U $ is strongly monadic over $\dgCoalg$.
\end{itemize}

\bigskip

We are going to consider in this chapter other Sweedler theories.
We will fully develop the theory of non-(co)unital (co)algebras and pointed (co)algebras
and sketched a few other in the last section.

It should be clear enough that our proofs in the unital dg-context can be adapted to these new contexts 
and the sections of this chapter gives the statements without proofs.

\subsection{The non-unital context}\label{nonunitalsweedlertheory}

In this section we construct six functors:
\begin{center}
\begin{tabular}{lrl}
\rule[-2ex]{0pt}{4ex} the tensor product of coalgebras & $\otimes$&$\!\!\!\!: \dgCoalg_\circ\times \dgCoalg_\circ \to \dgCoalg_\circ$, \\
\rule[-2ex]{0pt}{4ex} the coalgebra internal hom & $\HOM_\circ$&$\!\!\!\!:(\dgCoalg_\circ)^{op}\times \dgCoalg_\circ \to \dgCoalg_\circ$, \\
\rule[-2ex]{0pt}{4ex} the Sweedler hom & $\{-,-\}_\circ$&$\!\!\!\!:(\dgAlg_\circ)^{op}\times \dgAlg_\circ \to \dgCoalg_\circ$, \\
\rule[-2ex]{0pt}{4ex} the Sweedler product & $\rhd_\circ$&$\!\!\!\!:\dgCoalg_\circ\times \dgAlg_\circ \to \dgAlg_\circ$, \\
\rule[-2ex]{0pt}{4ex} the convolution product & $[-,-]$&$\!\!\!\!:(\dgCoalg_\circ)^{op}\times \dgAlg_\circ \to \dgAlg_\circ$, \\
\rule[-2ex]{0pt}{4ex} and the tensor product of algebras & $\otimes$&$\!\!\!\!:\dgAlg_\circ\times \dgAlg_\circ \to \dgAlg_\circ$.
\end{tabular}
\end{center}
such that
\begin{itemize}
\item the category $(\dgCoalg_\circ,\otimes, \HOM_\circ)$ is locally presentable symmetric monoidal closed and comonadic over $\dgVect$
\item and the category $(\dgAlg_\circ,\{-,-\}_\circ,\rhd_\circ, [-,-],\otimes)$ is locally presentable, enriched, bicomplete and monoidal over $\dgCoalg_\circ$, and monadic over $\dgVect$.
\end{itemize}
We will also distinguish the following isomorphism:
\begin{center}
\begin{tabular}{c}
\rule[-2ex]{0pt}{4ex} $T_\circ^\vee [C,X]\simeq \HOM_\circ(C,T_\circ^\vee (X))$ \\
\rule[-2ex]{0pt}{4ex} $\{T_\circ X,A\}\simeq T_\circ^\vee([X,A])$ \\
\rule[-2ex]{0pt}{4ex} $C\rhd_\circ T_\circ(X) \simeq T_\circ(C\otimes X)$
\end{tabular}
\end{center}
which are the data to strengthen the adjunctions $U\dashv T_\circ^\vee$ and $T_\circ\dashv U$ over $\dgCoalg_\circ$:
\begin{itemize}
\item the adjunction $U:\dgCoalg_\circ \rightleftarrows \dgVect:T_\circ^\vee $ is strongly comonadic over $\dgCoalg_\circ$
\item and the adjunction $T_\circ:\dgVect \rightleftarrows \dgAlg_\circ:U $ is strongly monadic over $\dgCoalg_\circ$.
\end{itemize}

Most of the result of the section are analogs of those of sections \ref{coalgebras} and \ref{algebras} so we will not repeat the proofs.
However, we draw the attention of the reader on some comparison result between the (co)unital and the non-(co)unital Sweedler operations
(they are at the end of each section).

\subsubsection{Presentability, comonadicity and cofree non-unital coalgebra}

\begin{thm}\label{noncounitalpresentability}
The category $\dgCoalg_\circ$ is finitary presentable. 
The $\omega$-compact objects are the finite dimensional coalgebras.
\end{thm}
\begin{proof}
As for theorem \ref{finprescoalg}.
\end{proof}

\begin{thm}\label{noncounitalcomonadicity}
The forgetful functor $U:\dgCoalg_\circ\to \dgVect$ has a right adjoint $T^\vee_\circ$ and the adjunction $U\dashv T^\vee_\circ$ is comonadic.
\end{thm}
\begin{proof}
As for theorem \ref{cofreecoalg}.
\end{proof}

\begin{prop}\label{isocofreeconil}
There exists a counital dg-coalgebra isomorphism $(T^\vee_\circ(X))_+ \simeq T^\vee(X)$.
\end{prop}
\begin{proof}
Any dg-vector space $X$ is naturally pointed by $0$, thus we have canonical unital coalgebra maps $T^\vee(0)\to T^\vee(X)\to T^\vee(0)$.
The free dgcoalgebra $T^\vee(0)$ is $\FF$ so the previous maps enhanced $T^\vee(X)$ in to a pointed dg-coalgebra.
Moreover since the maps $0\to X\to 0$ are natural in $X$, this enhancement of $T^\vee$ defines a functor $T^\vee_\bullet:\dgVect\to \dgCoalg_\bullet$.
It is easy to see that $T^\vee_\bullet$ is right adjoint to the forgetful functor $U:\dgCoalg_\bullet \to \dgVect$.
The result then follows using the equivalence $(-)_-:\dgCoalg_\bullet \simeq \dgCoalg_\circ:(-)_+$.
\end{proof}

\medskip

Recall from example \ref{ncutensorcoalg} the non-unital tensor coalgebra $T^c_\circ(X)$ on a dg-vector space $X$
and from proposition \ref{radical} the notion of radical $R^c$ of a non-unital coalgebra.

\begin{prop}\label{nilradicalcofree}
We have $T^c_\circ(X)=R^cT^\vee_\circ(X)$ for any vector space $X$.
\end{prop}

\begin{proof}
It follows from proposition \ref{structureconil} that
the functor $T^c_\circ: \dgVect \to \dgCoalgnilcirc$ is right adjoint to the forgetful functor $\dgCoalgnilcirc\to \dgVect$.
But the functor $R^cT^\vee_\circ $ is also right adjoint to the forgetful functor $\dgCoalgnilcirc\to \dgVect$,
since the functor $R^c:\dgCoalg \to \dgCoalgnilcirc$ is right adjoint to the inclusion $\dgCoalgnilcirc\subseteq \dgCoalg$
and the functor $T^\vee_\circ : \dgVect \to \dgCoalg$  is right adjoint to the forgetful functor  $\dgCoalg\to \dgVect$.
It follows that $T^c_\circ(X)=R^cT^\vee_\circ(X)$.
\end{proof}

\begin{prop}\label{comonadicTc}
The adjunction $U:\dgCoalgnilcirc \rightleftarrows \dgVect: T_\circ^c$ is comonadic.
\end{prop}
\begin{proof}
This adjunction decomposes as
$$\xymatrix{
\dgCoalgnilcirc \ar@<.6ex>[r]^-U & \dgCoalg_\circ \ar@<.6ex>[l]^-{R^c}\ar@<.6ex>[r]^-U & \dgVect \ar@<.6ex>[l]^-{T^\vee}
}$$
The adjunction $U\dashv T^\vee$ is comonadic by theorem \ref{noncounitalcomonadicity}
and the adjunction $U\dashv R^c$ is comonadic because it is a coreflexion.
Therefore the composite of the two forgetful functor satisfies the hypothesis of the comonadicity theorem.
\end{proof}

\begin{rem}
In particular $\dgCoalgnilcirc$ is finitely presentable. Because of the equivalence of example \ref{dualfinitealgex}, the category $\dgCoalgnilcirc^{op}$ is equivalent to the category $\mathsf{Pro}(\nfAlg_\circ)$ of (strict) pro-finite nilpotent non-unital algebras.
\end{rem}

\subsubsection{Internal hom}

The tensor product of non-counital coalgebras is defined in section \ref{monoidalstructurecoalgebras}.

\begin{thm}\label{closedmonoidunitalcoalg}
The category $(\dgCoalg_\circ,\otimes,\FF)$ is symmetric monoidal closed. 
\end{thm}
\begin{proof}
Same as for theorem \ref{homcoalg}.
\end{proof}

We shall denote the hom object between two non-counital coalgebras $C$ and $D$ by $\HOM_{\circ}(C,D)$. 
The counit of the adjunction $C\otimes(-)\dashv \HOM_{\circ}(C,-)$ is the {\it evaluation map} $\mathbf{ev} :\HOM_{\circ}(C,D)\otimes C\to D$.
For any non-counital coalgebra $E$ and any map of non-counital coalgebras $f:E\otimes C\to D$ there exists a unique map of non-counital coalgebras $g:E\to \HOM_\circ(C,D)$ such that ${\bf ev}(g\otimes C)=f$.
$$\xymatrix{
&& \HOM_\circ(C, D) \otimes C  \ar[d]^-{{\bf ev}} \\
E \otimes C \ar@{-->}[urr]^-{g\otimes C}  \ar[rr]^-{f}    && D
}$$
We shall put ${\Lambda}^2(f):=g$.
In addition, we shall put ${\Lambda}^1(f):={\Lambda}^2(f\sigma):C\to \HOM_\circ (E,D)$, where $\sigma:C\otimes E\to E\otimes C$ is the symmetry. More generally, we shall use the same notation $\Lambda^{i}$ and $\Lambda^{i,j}$ as in section \ref{closedmonoidalstructurecoalg} for multivariable lambda-transforms.

\begin{rem}
The coalgebras $\END_\circ(C)$ are examples of non-counital bialgebras in the sense of section \ref{pointed bialgebras}.
\end{rem}

\medskip
As in definition \ref{remreductioncoalg}, we can define a notion of {\em non-unital comorphism} and prove that the canonical map 
$\Psi =\lambda^2{\bf ev}:\HOM_\circ(C,D)\to [C,D]$ is the couniversal comorphism.
We leave the details to the reader.

\begin{prop} \label{closedmonoidcoalg4circ}  
For every non-unital coalgebra $C$, we have
\begin{enumerate}
\item $\HOM_\circ(\FF,C)=C$, in particular $\HOM_\circ(\FF,\FF)=\FF$,
\item $\HOM_\circ(C,0)=0$ and 
\item $\HOM_\circ(0,C)=0$.
\end{enumerate}
\end{prop}

\begin{proof}
Same as in proposition \ref{closedmonoidcoalg4}, using that $0$ is both initial and terminal in $\dgCoalg_\circ$.
\end{proof}

$\HOM_\circ (C,D)$ can be described as in proposition \ref{HOMsubcoalg} via the functor $T^\vee_\circ$.
In case $D$ is cofree we have also the following result, analogous to proposition \ref{exemplehomcoalg}.

\medskip
Let $C$ be a non-counital coalgebra and $X$ be a dg-vector space.
If $q: T^\vee_\circ([C,X])\to [C,X]$ is the cogenerating map, then the composite of the maps
$$\xymatrix{
T^\vee_\circ([C,X])\otimes C \ar[rr]^-{q\otimes C} && [C,X]\otimes C \ar[r]^-{ev} & X
}$$
can be coextended as a map of non-counital coalgebras $e: T^\vee_\circ ([C,X])\otimes C \to T^\vee_\circ (X)$.

\begin{prop}\label{exemplehomcoalgnu}
With the previous notations, for any non-counital coalgebra $E$ and any map of non-counital coalgebras
$f:E\otimes C\to T^\vee_\circ (X)$, there exists a unique map of non-counital coalgebras
$k:E\to  T^\vee_\circ([C,X])$ such that $e(k\otimes C)=f$. 
Thus, $e$ is a strong evaluation $\mathbf{ev}:\HOM_\circ(C,T^\vee_\circ(X)) \otimes C\to  T^\vee_\circ([C,X])$ and we have
$$
\HOM_\circ(C,T^\vee_\circ(X)) = T^\vee_\circ ([C,X]).
$$
\end{prop}
\begin{proof}
Similar to the proof of proposition \ref{exemplehomcoalg}.
\end{proof}

\medskip

The strong composition law $\mathbf{c}:\HOM_\circ(D,E)\otimes \HOM_\circ(C,D) \to  \HOM_\circ(C,E)$ is defined as $\Lambda^3({\bf ev}^2)$ where ${\bf ev}^2$ is the map
$$\xymatrix{
\HOM_\circ(D,E)\otimes \HOM_\circ(C,D)\otimes C \ar[rrr]^-{\HOM_\circ(D,E)\otimes {\bf ev}}&&&\HOM_\circ(D,E)\otimes D \ar[r]^-{\bf ev}& E.
}$$
Equivalently it can be proven to be the unique map of non-counital coalgebras $\mathbf{c}$ such that the following square commutes
$$\xymatrix{
\HOM_\circ(D,E)\otimes \HOM_\circ(C,D)\ar[rr]^-{\mathbf{c}}\ar[d]_{\Psi \otimes \Psi }&& \HOM_\circ(C,E)\ar[d]^{\Psi}\\
[D,E]\otimes [C,D]\ar[rr]^-{c} && [C,E].
}$$
We shall denote $\dgCoalg_\circ^\$$ the category $\dgCoalg_\circ$ viewed as enriched over itself.

\begin{prop}\label{atomhomnucoalg}
The atoms of $\HOM_\circ(C,D)$ are in bijection with the maps of non-counital coalgebras.
The underlying category of $\dgCoalg_\circ^\$$ is $\dgCoalg_\circ$.
\end{prop}
\begin{proof}
Similar to lemma \ref{underlyingCoalglemma} and proposition \ref{underlyingCoalg}.
\end{proof}

\medskip
\begin{rem}
The enriched structure of $\dgCoalg_\circ^\$$ restricts to defined an enriched category $(\dgCoalgnilcirc)^\$$ of conilpotent coalgebras, but the coalgebra of morphisms between two non-zero conilpotent coalgebras is never conilpotent because it always has at least two atoms, the identity and the zero map. 
This situation look like a construction of B. Keller in \cite{Keller} where he defines a right adjoint to the tensor product of conilpotent coalgebras as a dg-cocategory. A comparison of this dg-cocategory with the hom dg-coalgebra would be nice.
\end{rem}

\bigskip

We finish this section by comparing the functors $\HOM$ and $\HOM_\circ$.

\begin{prop}\label{HOMvsHOMnu}
Let $C$ and $D$ be counital coalgebra, then we have a pullback square of non-counital coalgebras
$$\xymatrix{
\HOM(C,D) \ar[rr]^-\iota \ar[d]_-{\epsilon} && \HOM_\circ(C,D)\ar[d]^{\HOM_\circ(C,\epsilon_D)}\\
\FF \ar[rr]^-{\epsilon_C} && \HOM_\circ(C,\FF)
}$$ 
such that the map $\epsilon:\HOM(C,D)\to \FF$ is the counit of $\HOM(C,D)$.
In particular, $\HOM(C,D)$ is a non-counital sub-coalgebra of $\HOM_\circ(C,D)$.
\end{prop}

\begin{proof}
Let us prove that the square is commutative. 
The naturality of the universal counital comorphism gives a square
$$\xymatrix{
\HOM(C,D) \ar[rr]^-{\Psi} \ar[d]_{\HOM(C,\epsilon_D)} && [C,D] \ar[d]^{[C,\epsilon_D]}\\
\HOM(C,\FF)\ar[rr]^-{\Psi} && [C,\FF].
}$$
A counital comorphism is in particular a non-counital comorphism so the square factors into
$$\xymatrix{
\HOM(C,D) \ar[rr]^-\iota\ar[d]_{\HOM(C,\epsilon_D)} &&\HOM_\circ (C,D) \ar[rr]^-{\Psi} \ar[d]_{\HOM_\circ(C,\epsilon_D)} && [C,D] \ar[d]^{[C,\epsilon_D]}\\
\HOM(C,\FF)\ar[rr]&&\HOM_\circ(C,\FF)\ar[rr]^-{\Psi} && [C,\FF].
}$$
where we have noted $\Psi_\circ$ the non-counital universal comorphism to distinguish it from the counital one.
Using corollary \ref{corcounitHOM}, the left vertical map is the counit $\epsilon:\HOM(C,D)\to \FF$.
Moreover the atom generating $\FF=\HOM(C,\FF)$ is $\epsilon_C^\sharp$.
By lemma \ref{musicallemmacoalg}, it is send to $\epsilon_C\in [C,\FF]$ by $\Psi:\HOM(C,\FF)\to [C,\FF]$.
Then by universal property of $\Psi_\circ:\HOM_\circ(C,\FF)\to [C,\FF]$, the map $\FF=\HOM(C,\FF)\to \HOM_\circ(C,D)$ lifting $\FF\to [C,\FF]$ is the non-counital atom $\FF\to \HOM_\circ(C,D)$ corresponding to $\epsilon_C$.
This proves that the left square of the diagram identifies with the square of the proposition, hence its commutativity.

Let us prove now that the square is cartesian.
We consider a commutative square 
$$\xymatrix{
E \ar[d]_{\epsilon} \ar[rr]^-\iota && \HOM_\circ(C,D)\ar[d]^{\HOM_\circ(C,\epsilon_D)} \\
\FF \ar[rr]^-{\epsilon_C} && \HOM_\circ(C,\FF)
}$$ 
we need to construct a map $E\to \HOM(C,D)$ which commutes with the $\epsilon$s and the $\iota$s

The composite $\psi_\circ\iota:E\to \HOM_\circ(C,D)\to [C,D]$ is a non-counital comorphism.
$E\to [C,D]$ is a counital comorphism iff there is a commutative square
$$\xymatrix{
E\ar[d]\ar[rr] && [C,D] \ar[d]^{[C,\epsilon_D]}\\
\FF\ar[rr]^-{\epsilon_C} && [C,\FF]
}$$
But by naturality of $\Psi_\circ$ there is a commutative diagram
$$\xymatrix{
E\ar[d]_\epsilon \ar[rr]^-\iota && \HOM_\circ(C,D)\ar[rr]^-{\Psi_\circ} \ar[d]^{\HOM_\circ(C,\epsilon_D)} && [C,D] \ar[d]^{[C,\epsilon_D]} \\
\FF\ar[rr]^-{\HOM(\epsilon_C,\FF)} && \HOM_\circ(C,\FF)\ar[rr]^-{\Psi_\circ} && [C,\FF]
}$$
This proves that $E\to [C,D]$ is a counital morphism and we deduce a counital coalgebra map $E\to \HOM(C,D)$.
The same reasoning apply if $E=\FF=D$ and we have a commutative diagram of non-counital coalgebras
$$\xymatrix{
& \HOM(C,D)\ar[rd]^-\iota\ar'[d][dd]_-(.6)\epsilon\\
E\ar[ru]\ar[dd]_\epsilon \ar[rr]^-(.6)\iota && \HOM_\circ(C,D)\ar[dd]^{\HOM_\circ(C,\epsilon_D)} \\
& \HOM(C,\FF)\ar[rd]^-\iota\\
\FF\ar@{=}[ru]\ar[rr]^-{\HOM(\epsilon_C,\FF)} && \HOM_\circ(C,\FF)
}$$
This proves that the square of the proposition is cartesian.

The last statement is a consequence of the following lemma.
\end{proof}

\begin{rem}\label{remHOMvsHOMnu}
In the proof of proposition \ref{HOMvsHOMnu}, we have constructed the inclusion $\HOM(C,D)\to \HOM_\circ(C,D)$ as the unique non-counital coalgebra map such that the triangle
$$\xymatrix{
\HOM(C,D)\ar[rr]\ar[rrd]_\Psi&&\HOM_\circ(C,D) \ar[d]^{\Psi_\circ}\\
&&[C,D]
}$$
commutes.
\end{rem}

\begin{lemma} \label{directandinverseimage5}
If a map of coalgebras $u:C\to D$ is injective then a commutative square of  coalgebras 
$$\xymatrix{
C'\ar[r]^f\ar[d]_{u'} & C \ar[d]^u\\
D'\ar[r]^g & D 
}$$
is cartesian if and only if it is cartesian in the category of vector spaces.
\end{lemma}
\begin{proof} We may suppose that the map $u:C\to D$
is defined by an inclusion $C\subseteq D$.
We then have a pullback square of coalgebras
 $$\xymatrix{
g^{-1}(C) \ar[r]^h\ar[d]_{i} & C \ar[d]^u \\
D'\ar[r]^g & D 
}$$
where $h$ is induced by $g$ and $i$ is the inclusion $g^{-1}(C)\subseteq D$.
Hence there is a unique morphism of coalgebras $v:C'\to g^{-1}(C)$ such that $iv=u'$ and $hv'=f$.
The square $S$ of coalgebras in the statement of the proposition is a pulback if and only if the morphism $v$ is an isomorphism, if and only if the linear map $U(v)$ is invertible if and only if the image of $S$ is a pullback in the category of vector spaces.
\end{proof}

\medskip

\begin{prop}  \label{radical3} 
If $C$ is a nilpotent non-unital coalgebra, then we have
$$
\HOM_{\circ}(C,R^cD)\simeq \HOM_{\circ}(C,D)
$$
for any non-unital coalgebra $D$.
In other words, the radical adjunction 
$$\xymatrix{
\iota:\dgCoalgnilcirc \ar@<.6ex>[r] &\dgCoalg_\circ:R^c \ar@<.6ex>[l]
}$$
is strong.
\end{prop}
\begin{proof} 
By lemma \ref{radical2}, the non-unital coalgebra $E\otimes C$ is conilpotent for any non-unital coalgebra $E$, since $C$ is conilpotent by hypothesis.
It follows that every map of non-unital coalgebras $E\otimes C\to D$ factors through the inclusion $R^cD\subseteq D$.
This defines a natural bijection between the maps of non-unital coalgebras $E\otimes C\to D$
and the maps of non-unital coalgebras $E\otimes C\to R^cD$.
Thus, we obtain a natural bijection between the maps of non-unital coalgebras $E\to \HOM_{\circ}(C,D)$
and the maps of non-unital coalgebras $E\to \HOM_{\circ}(C,R^cD)$.
It follows by Yoneda lemma that the inclusion $R^cD\subseteq D$ induces a natural isomorphism $\HOM_{\circ}(C,R^cD)\simeq \HOM_{\circ}(C,D)$.
\end{proof}

\begin{prop}\label{HomCplus} 
If $C$ is a coalgebra and $D$ is a non-counital coalgebras, then we have a natural
isomorphism of coalgebras
$$\HOM(C,D_{+}) \simeq \HOM_{\circ}(C,D)_{+}.$$
\end{prop}

\begin{proof} 
 If $E$ is a coalgebra,  there is a chain of natural bijections, 
between the 
\begin{center}
\begin{tabular}{lc}
\rule[-2ex]{0pt}{4ex}  coalgebra maps & $E \to \HOM(C,D_{+}) $,\\
\rule[-2ex]{0pt}{4ex}  coalgebra maps  & $E\otimes C \to D_{+}$,\\
\rule[-2ex]{0pt}{4ex}  non-counital coalgebra maps &$ E\otimes C \to D$,\\
\rule[-2ex]{0pt}{4ex} non-counital coalgebra maps & $E \to \HOM_\circ (C,D) $, \\
\rule[-2ex]{0pt}{4ex}  coalgebra maps &  $E \to \HOM_\circ (C,D)_{+}$. \\
\end{tabular}
\end{center}
It then follows by Yoneda lemma that the coalgebra $\HOM(C,D_{+}) $
is isomorphic to the algebra $\HOM_\circ (C,D)_{+}$.
\end{proof}

\subsubsection{Convolution and measuring}

Let $C$ be a non-unital coalgebra and $A$ and $B$ be two non-unital algebras. 
The dg-vector space $[C,A]$ is a non-unital algebra called the {\em non-unital convolution algebra}.

\medskip
A map $f:C\otimes A\to B$ is called a {\em non-unital measuring} if the map $\lambda^1f:A\to [C,B]$ is a non-unital algebra map.
They can be characterized as in the unital case, we just have to remove the unit condition:
a linear map $f:C\otimes A\to B$ is a measuring if and only if, for every $a,b\in A$ and $c\in C$,
$$
f(c,ab)= f(c^{(1)},a)f(c^{(2)},b) \ (-1)^{|a||c^{(2)}|}.
$$
Let $\cM_\circ(C,A;B)$ be the set of non-unital measurings $C\otimes A\to B$, theses sets define a functor
$$\xymatrix{
\cM_\circ(-,-;-):(\dgCoalg_\circ)^{op}\times (\dgAlg_\circ)^{op}\times \dgAlg_\circ\ar[r]& \Set.
}$$
By definition, it is representable in the second variable by the {\em non-unital convolution product} functor
$$\xymatrix{
[-,-]:(\dgCoalg_\circ)^{op}\times \dgAlg \ar[r]& \dgAlg_\circ.
}$$

\subsubsection{Sweedler product}

A non-unital measuring $u:C\otimes A\to E$ is said to be {\it universal} if the pair $(E,u)$ is representing the functor 
$$\xymatrix{
\cM_\circ(C,A;-): \dgAlg \ar[r]& \Set.
}$$
The object $E$ of a universal measuring is well defined up to a unique isomorphism. 
We shall denote it by $C\rhd_\circ A$ and write $c\rhd_\circ a:=u(c\otimes a)$ for $c\in C$ and $a\in A$.
We shall say that the non-unital algebra $C\rhd_\circ A$ is the {\it Sweedler product} of $A$ by $C$.

\begin{thm} \label{Sproductcirc}
The Sweedler product $C\rhd_\circ A$ exists for any non-unital algebra $A$ and any non-counital coalgebra $C$.
The functor $\rhd_\circ:\dgCoalg_\circ\times \dgAlg_\circ\to\dgAlg_\circ$ has the structure of a left action of the monoidal category $\dgCoalg_\circ$ on the category $\dgAlg_\circ$.

\end{thm} 
\begin{proof}
As in theorems \ref{Sproduct} and \ref{actionsweedlerproduct}.
\end{proof}

The non-unital algebra $C\rhd_{\circ} A$ can be constructed explicitely as follows.
It is generated by symbols, $c\rhd_\circ a$ for $c\in C$ and $a\in A$ on which the following relations are imposed,
\begin{itemize}
\item[(d)] the map $(c,a)\mapsto c\rhd a$ is a bilinear map of dg-vector spaces, $d(c\rhd a) = dc\rhd a + c\rhd da (-1)^{|c|}$;
\item[(m)] $c\rhd_{\circ} (ab) = (c^{(1)}\rhd_{\circ} a)(c^{(2)}\rhd_{\circ} b) (-1)^{|a||c^{(2)}|}$, for every $c\in C$ and $a,b\in A$;
\end{itemize}
In other terms, $C\rhd_\circ A$ is the quotient of the non-unital tensor dg-algebra $T_\circ(C\otimes A)$ by the relations (m).

\medskip

The canonical isomorphism $\FF\otimes A\to A$
is a non-unital measuring and it is universal.  It follows that we have $\FF\rhd_{\circ} A\simeq A$.

\medskip

\begin{prop}\label{Sweedleroftensornu}
Let $C$ be a non-counital coalgebra, $X$ a graded vector space and $A$ be a non-unital algebra.
Then every linear map $f:C\otimes X\to A$ can be extended uniquely as a non-unital measuring $f' :C\otimes T_\circ (X)\to A$.
Moreover, if $i$ is the inclusion $C\otimes X \to T_\circ(C\otimes X)$, then the non-unital measuring $i':C\otimes T_\circ(X)\to T_\circ (C\otimes X)$ is universal.
Hence we have 
$$
C\rhd_\circ T_\circ (X)= T_\circ(C\otimes X).
$$
\end{prop}

\begin{proof} Same proof as proposition \ref{Sweedleroftensor}
\end{proof}

\bigskip

We finish this section by comparing the functors $\rhd$ and $\rhd_\circ$.

\begin{prop} \label{SproductvsSproductnu}
If $C$ is a counital coalgebra and $A$ is a unital algebra, then we have a pushout in the category of non-unital algebras
$$\xymatrix{
C\rhd_\circ \FF\ar[rr]^-{C\rhd_\circ e_A}\ar[d]_{\epsilon_C\rhd_\circ \FF}&& C\rhd_\circ A\ar[d]\\
\FF\ar[rr]^-{e}&& C\rhd A.
}$$
such that the map $e:\FF\to C\rhd A$ is the unit of $C\rhd A$.
In particular, $C\rhd A$ is a quotient of the algebra $C\rhd_\circ A$.
\end{prop} 

\begin{proof}
This can be proven with a scheme dual to that of proposition \ref{HOMvsHOMnu}, where the universal unital and non-unital measurings $\Phi:C\otimes A\to C\rhd A$ and $\Phi:C\otimes A\to C\rhd_\circ A$ replace their comorphisms analogs. We leave it to the reader.
\end{proof}

\begin{prop} \label{Sproductcoalgnualg}
If $C$ is a coalgebra and $A$ is a non-unital algebra, then we have $C\rhd A_{+}=(C\rhd_{\circ} A)_{+}$.
\end{prop} 

\begin{proof} 
 If $B$ is an algebra,
  there is a chain of natural bijections, 
between the 
\begin{center}
\begin{tabular}{lc}
\rule[-2ex]{0pt}{4ex}  algebra maps & $C\rhd A_{+} \to B$,\\
\rule[-2ex]{0pt}{4ex}  algebra maps  & $A_{+}\to [C,B]$,\\
\rule[-2ex]{0pt}{4ex}  non-unital algebra maps &$A \to [C,B]$,\\
\rule[-2ex]{0pt}{4ex}  non-unital algebra maps & $ C\rhd_\circ A \to B$, \\
\rule[-2ex]{0pt}{4ex}  algebra maps &  $(C\rhd_\circ A)_{+} \to B$. \\
\end{tabular}
\end{center}
It then follows by Yoneda lemma that the algebra $C\rhd A_{+}$
is isomorphic to the algebra $(C\rhd_{\circ} A)_{+}$.
\end{proof}

\subsubsection{Sweedler hom}

Let $A$ and $B$ be non-unital algebras and $E$ be a non-counital coalgebra.
We shall say that a non-unital measuring $v:E\otimes A\to B$ is {\it couniversal} if the pair $(E,v)$ is representing the functor 
$$
\cM_\circ(-,A;B): \Coalg^{op}_\circ \to \Set.
$$
The non-counital coalgebra $E$ of a couniversal non-unital measuring $v:E\otimes A\to B$ is well defined up to a unique isomorphism 
and we shall denote it by $\{A,B\}_\circ$ and call it the {\em non-counital Sweedler Hom}. 
We shall denote the couniversal non-unital measuring as a {\it strong evaluation} 
$\mathbf{ev}:\{A,B\}_\circ \otimes A \to B$.
By definition, for any non-counital coalgebra $C$ and any non-unital measuring $f:C\otimes A\to B$, there is a unique map of non-counital coalgebras
$g:C\to \{A,B\}_\circ $ such that ${\bf ev}(g\otimes A)=f$.
$$
\xymatrix{
			&& \{A,B\}_\circ \otimes A	\ar[d]^-{\mathbf{ev}} \\
C\otimes A	\ar@{-->}[urr]^-{g\otimes A} \ar[rr]^-{f} && B. 
}$$

As in the counital case, it is useful in some proofs to use the dual notion of non-unital comeasuring.

\begin{defi}
Let $A$ and $B$ be non-unital algebras and $C$ be a non-counital coalgebra.
We shall say that a linear map $g:C\to [A,B]$ is a  {\em non-unital comeasuring} if the map $ev(g\otimes A):C\otimes A\to B$
is a non-unital measuring.
\end{defi}

A linear map $g:C\to [A,B]$ is a non-unital comeasuring if and only if we have 
$$
g(c)(ab)= g(c^{(1)})(a)g(c^{(2)})(b) \ (-1)^{|a||c^{(2)}|}
$$
for every $a,b\in A$ and $c\in C$.
There are canonical bijections between
\begin{center}
\begin{tabular}{lc}
\rule[-2ex]{0pt}{4ex} non-unital algebra maps & $A\to [C,B]$, \\
\rule[-2ex]{0pt}{4ex} non-unital measurings & $C\otimes A\to B$, \\
\rule[-2ex]{0pt}{4ex} non-unital comeasurings & $C\to [A,B]$.
\end{tabular}
\end{center}

We shall say that a non-unital comeasuring  $k:E\to [A,B]$ is {\it couniversal} if the corresponding measuring $ev( k\otimes A):E\otimes A\to B$
is couniversal. The couniversality of $k$
means concretely that for any non-counital coalgebra $C$ and any non-unital comeasuring $g:C\to [A,B]$ there exists a unique map of non-counital coalgebras
$f:C\to E$ such that $g=k f$.
$$
\xymatrix{
	&& E\ar[d]^k \\
C\ar[rr]_-{g}\ar@{-->}[urr]^f && [A,B] 
}$$

\begin{thm} \label{SemiCoalghomalg1}
There exists a couniversal non-unital measuring $\mathbf{ev}: \{A,B\}_\circ \otimes A\to B$
any pair of non-unital algebras $A$ and $B$. Equivalently, there exists a couniversal 
non-unital comeasuring $\Psi:\{A,B\}_\circ \to [A,B]$.
\end{thm}
\begin{proof} Similar to the proof of theorem \ref{Coalghomalg1}.
\end{proof}

As in the unital case, explicit constructions of $\{A,B\}_\circ$ rely on the cofree non-counital coalgebra $T^\vee_\circ$.
In particular, the comeasuring $\Psi:\{A,B\}_\circ\to [A,B]$ is cogenerating and there is a non-unital analog of corollary \ref{betterShom}.
We leave its statement to the reader, we shall state only the analog of proposition \ref{exSweedlerhom}.

\medskip

Let $A$ be a non-unital algebra and $X$ be a vector space.
If $p: T^\vee_\circ([X,A])\to [X,A]$ is the cofree map, then the composite of the maps
$$\xymatrix{
T^\vee_\circ([X,A])\otimes X \ar[rr]^-{p\otimes X} && [X,A]\otimes X \ar[r]^-{ev} & A
}$$
can be extended as a non-unital measuring $p':T^\vee_\circ([X,A])\otimes T_\circ (X)\to A$
by proposition \ref{Sweedleroftensornu}.

\begin{prop} \label{exSweedlerhomnu}
If $A$ is a non-unital algebra and $X$ is a vector space, 
then the non-unital measuring $p':T^\vee_\circ ([X,A])\otimes T_\circ(X)\to A$
defined above is couniversal.
Hence we have 
$$
\{T_\circ(X),A\}_\circ= T^\vee_\circ([X,A]).
$$
\end{prop}

\begin{proof} Same proof as proposition \ref{exSweedlerhom} . 
\end{proof}

\begin{cor} \label{corexSweedlerhomnu}
If $A$ is a non-unital algebra and $X$ is a vector space, 
there is a natural isomorphism
$$
R^c\{T_\circ(X),A\}_\circ= T^c_\circ([X,A]).
$$
and the non-counital measuring corresponding to the inclusion $R^c\{T_\circ(X),A\}_\circ\to \{T_\circ(X),A\}_\circ$ is the 
unique non-counital measuring extending 
$$\xymatrix{
T^c_\circ([X,A])\otimes X \ar[rr]^-{p\otimes X} && [X,A]\otimes X \ar[r]^-{ev} & A
}$$
where $p:T^c_\circ([X,A])\to [X,A]$ is the cogenerating map.
\end{cor}

\begin{proof}
The first assertion is a direct consequence from $R^cT^\vee_\circ=T^c_\circ$. 
The second is by composing the universal measuring with the coalgeba map $T^c_\circ([X,A])\to T^\vee_\circ([X,A])$
and proposition \ref{Sweedleroftensornu}.
\end{proof}

\bigskip

If $A$, $B$ and $E$ are non-unital algebras, there is then a unique map of non-counital coalgebras
$\mathbf{c}: \{B,E\}_\circ \otimes \{A,B\}_\circ \to   \{A,E\}_\circ$
such that the following square commutes
$$\xymatrix{
\{B,E\}_\circ \otimes \{A,B\}_\circ \ar@{-->}[rr]^-{\mathbf{c}}\ar[d]_{\Psi\otimes \Psi}&& \{A,E\}_\circ \ar[d]^{\Psi}\\
[B,E]\otimes [A,B]\ar[rr]^-c && [A,E].
}$$
For any non-unital algebra $A$, there is also a unique map of non-counital coalgebra $e_A:\FF \to \{A, A\}_\circ$ such that $\Psi e_A=1_A$,
$$
\xymatrix{
	&& \{A,A\}_\circ \ar[d]^\Psi \\
\FF \ar[rr]_-{1_A}\ar@{-->}[urr]^-{e_A} && [A,A] 
}$$

\begin{thm}\label{enrichmentalgcoalpartial}
The coalgebras $\{A,B\}_\circ$ define an enrichment of the category $\dgAlg_\circ$ over the closed monoidal category $\dgCoalg_\circ$. 
The unit of the composition law is the map $e_A:\FF \to \{A, A\}$.
The enriched category $\dgAlg_{\circ}$ is bicomplete;
the tensor product of a non-unital algebra $A$ by a non-counital coalgebra $C$ is the non-unital algebra $C\rhd_{\circ} A$ and 
the cotensor product of $A$ by $C$ is the convolution non-unital algebra $[C,A]$.
Hence there are natural isomorphisms of non-counital coalgebras
$$
\{C\rhd_{\circ} A,B\}_{\circ} \simeq \HOM_{\circ}(C,\{A,B \}_{\circ})\simeq \{A,[C,B] \}_{\circ}
$$
for a non-counital coalgebra $C$ and non-unital algebras $A$ and $B$.
\end{thm}

\begin{proof} Same proof as theorem \ref{enrichmentalgcoal}.
\end{proof}

\begin{rem}
The coalgebras $\{A,A\}_\circ$ are examples of a non-counital bialgebras in the sense of section \ref{pointed bialgebras}.
\end{rem}

Let $\dgAlg_\circ^\$$ be the category $\dgAlg_\circ$ viewed as enriched over $\dgCoalg_\circ$.

\begin{prop}\label{nuatomShom} 
Let $A$ and $B$ be two non-unital algebras, the set of atoms of $\{A,B\}_\circ$ is in bijection with the set of maps of non-unital  algebras from $A$ to $B$. The underlying category of $\dgAlg_\circ^\$$ is the ordinary category $\dgAlg_\circ$.
\end{prop}
\begin{proof}
As in lemma \ref{atomShom} and proposition \ref{underlyingAlg}
\end{proof}

\bigskip

\begin{cor} \label{nuadjunctiontypes}
Let $C$ be a non-counital coalgebra and $A$ be a non-unital algebra.
We have the following strong adjunctions
$$\xymatrix{
C\rhd_\circ(-):\dgAlg_\circ  \ar@<.6ex>[r]& \dgAlg_\circ:[C,-]\ar@<.6ex>[l]
}$$
$$\xymatrix{
[-,A]: (\dgCoalg_\circ)^{op}  \ar@<.6ex>[r]& \dgAlg_\circ:\{-,A \}_\circ\ar@<.6ex>[l]
}$$
$$\xymatrix{
(-)\rhd_\circ A: \dgCoalg_\circ  \ar@<.6ex>[r]& \dgAlg_\circ:\{A,- \}_\circ\ar@<.6ex>[l].
}$$
\end{cor}

\begin{cor}
Let $A$ be a non-unital algebra, we have the following strong adjunctions
$$\xymatrix{
[-,A]: \dgCoalgnilcirc^{op}  \ar@<.6ex>[r]& \dgAlg_\circ:R^c\{-,A \}_\circ\ar@<.6ex>[l]
}$$
$$\xymatrix{
(-)\rhd_\circ A: \dgCoalgnilcirc  \ar@<.6ex>[r]& \dgAlg_\circ:R^c\{A,- \}_\circ\ar@<.6ex>[l].
}$$

\end{cor}
\begin{proof}
Compose the previous result with the radical adjunction $\iota:\dgCoalgnilcirc\rightleftarrows \dgCoalg_\circ:R^c$.
\end{proof}

\begin{cor}
Let $X$ be a (dg-)vector space, we have the following strong adjunctions
$$\xymatrix{
T_\circ(X\otimes -): \dgCoalg_\circ  \ar@<.6ex>[r]& \dgAlg_\circ:T^\vee_\circ([X,-])\ar@<.6ex>[l].
}$$
$$\xymatrix{
T_\circ(X\otimes -): \dgCoalgnilcirc  \ar@<.6ex>[r]& \dgAlg_\circ:T^c_\circ([X,-])\ar@<.6ex>[l].
}$$
\end{cor}
\begin{proof}
By propositions \ref{Sweedleroftensornu} and proposition \ref{exSweedlerhomnu}
and corollary \ref{corexSweedlerhomnu}, 
we have $T_\circ(X\otimes -)=(-)\rhd_\circ T_\circ X$, 
$T^\vee_\circ([X,-])=\{T_\circ X,- \}_\circ$
and $T^c_\circ([X,-])=R^c\{T_\circ X,- \}_\circ$.
\end{proof}

\bigskip

We finish this section by comparing the functors $\{-,-\}$ and $\{-,-\}_\circ$.

\begin{prop} \label{SHOMvsSHOMnu}
If $A$ and $B$ are two unital algebras, then we have a pushout in the category of non-counital coalgebras
$$\xymatrix{
\{A,B\} \ar[rr]^-\iota \ar[d]_-{\epsilon} && \{A,B\}_\circ\ar[d]^{\{e_A,B\}_\circ}\\
\FF \ar[rr]^-{e_B} && \{\FF,B\}
}$$ 
such that the map $\epsilon:\{A,B\}\to \FF$ is the counit of $\{A,B\}$.
In particular $\{A,B\}$ is a non-counital sub-coalgebra of $\{A,B\}_\circ$.
\end{prop} 

\begin{proof}
The proof is the same as for proposition \ref{HOMvsHOMnu} with the universal unital and non-unital comeasurings $\Psi:\{A,B\}\to [A,B]$ and $\Psi:\{A,B\}_\circ\to [A,B]$ instead of their comorphisms analogs.
We leave it to the reader.
\end{proof}

\begin{prop} \label{SHOMnualgalg}
If $A$ is a non-unital algebra and $B$ is an algebra, then we have $\{A_{+},B\}=\{A,B\}_{\circ +}$.
\end{prop} 
\begin{proof}
 If $C$ is a coalgebra, there is a chain of natural bijections, between the 
\begin{center}
\begin{tabular}{lc}
\rule[-2ex]{0pt}{4ex} the coalgebra maps & $C\to \{A_{+},B\}$,\\
\rule[-2ex]{0pt}{4ex}  the measurings  & $C\otimes A_{+}\to B$,\\
\rule[-2ex]{0pt}{4ex} the algebra maps &$A_{+}\to [C,B]$,\\
\rule[-2ex]{0pt}{4ex} the non-unital algebra maps & $ A \to [C,B]$, \\
\rule[-2ex]{0pt}{4ex}  the non-unital measurings &  $C\otimes A \to B$, \\
\rule[-2ex]{0pt}{4ex} the non-counital coalgebra maps & $C\to \{A,B\}_{\circ}$,\\
\rule[-2ex]{0pt}{4ex} the coalgebra maps & $C\to \{A,B\}_{\circ +}$.
\end{tabular}
\end{center}
It then follows by Yoneda lemma that the coalgebra $\{A_{+},B\}$
is isomorphic to the coalgebra $(\{A,B\}_\circ)_{+}$.
\end{proof}

\subsubsection{Monoidal strength and lax structures}

$\dgCoalg_\circ$ and $\dgAlg_\circ$ are enriched over $\dgCoalg$ and so are their opposite categories.
A cartesian product $\bC\times \bD$ where $\bC$ and $\bD$ are any of the categories $\dgCoalg_\circ$, $\dgAlg_\circ$ or their opposites, is naturally enriched over $\dgCoalg\times \dgCoalg$. We can transfer this into an enrichment over $\dgCoalg_\circ$ along the functor
$\otimes : \dgCoalg_\circ\times \dgCoalg_\circ\to \dgCoalg_\circ$. We shall note $(\bC\times \bD)^{\$\otimes\$}$ this enrichment.

For example, we have categories $(\dgCoalg_\circ\times \dgCoalg_\circ)^{\$\otimes \$}$ where the coalgebra of morphisms from $(C_1,C_2)$ to $(D_1,D_2)$ is
$$
\HOM_\circ(C_1,C_2) \otimes \HOM_\circ(D_1,D_2)
$$
and we have $((\dgCoalg_\circ)^{op}\times \dgAlg_\circ)^{\$\otimes \$}$ where the coalgebra of morphisms from $(C,A)$ to $(D,B)$ is 
$$
\HOM_\circ(D,C)\otimes \{A,B\}_\circ.
$$

\begin{prop}\label{numonoidalstrength}
The functors 
$$\xymatrix{
\otimes:\dgCoalg_\circ\times \dgCoalg_\circ\ar[r]& \dgCoalg_\circ,
}$$ 
$$\xymatrix{
\otimes:\dgAlg_\circ\times \dgAlg_\circ\ar[r]& \dgAlg_\circ
}$$ 
enhance into strong symmetric monoidal structures.
\end{prop}
\begin{proof}
Same as in section \ref{monoidalstrengthcoalg} and theorem \ref{tensorenrichmentof}.
\end{proof}

\begin{prop}\label{nulaxmonoidalstrength}
\begin{enumerate}
\item The tensor products of algebras and coalgebras enhance to strong symmetric monoidal functors
$$\xymatrix{
\otimes:(\dgCoalg_\circ\times \dgCoalg_\circ)^{\$\otimes \$}\ar[r]& \dgCoalg_\circ^{\$}
}$$ 
$$\xymatrix{
\otimes:(\dgAlg_\circ\times \dgAlg_\circ)^{\$\otimes \$}\ar[r]& \dgAlg_\circ^{\$}.
}$$

\item The functors $\HOM_\circ$, $\{-,-\}_\circ$ and $[-,-]$ enhance to strong symmetric lax monoidal functors
$$\xymatrix{
\HOM_\circ:((\dgCoalg_\circ)^{op}\times \dgCoalg_\circ)^{\$\otimes \$}\ar[r]& \dgCoalg_\circ^\$
}$$ 
$$\xymatrix{
\{-,-\}_\circ:((\dgAlg_\circ)^{op}\times \dgAlg_\circ)^{\$\otimes \$}\ar[r]& \dgCoalg_\circ^\$
}$$ 
$$\xymatrix{
[-,-]:((\dgCoalg_\circ)^{op}\times \dgAlg_\circ)^{\$\otimes \$}\ar[r]& \dgAlg_\circ^\$
}$$

\item The Sweedler product enhance to strong symmetric colax monoidal functor
$$\xymatrix{
\rhd_\circ:(\dgCoalg_\circ\times \dgAlg_\circ)^{\$\otimes \$}\ar[r]& \dgAlg_\circ^\$.
}$$ 
\end{enumerate}
\end{prop}
\begin{proof}
As in section \ref{laxstructures}.
\end{proof}

In particular, these functors induces functors on the corresponding categories of (co)monoids and (co)commutative (co)monoids. We leave all the statements to the reader. We only detail the construction of the (co)lax structures that we will use.

\medskip
The lax structure of $\{-,-\}_\circ$ is given by $(\alpha,\alpha_0)$ where $\alpha_0:\FF\simeq \{\FF,\FF\}_\circ$ and 
$\alpha$ is the unique map of non-counital coalgebras such that the following square commute
$$\xymatrix{
\{A_1,B_1\}_\circ\otimes \{A_2,B_2\}_\circ \ar[rr]^-{\alpha} \ar[d]_{\Psi\otimes \Psi}&& \{A_1\otimes A_2,B_1\otimes B_2\}_\circ\ar[d]^{\Psi}\\
[A_1,B_1]\otimes [A_2,B_2] \ar[rr]^-\theta && [A_1\otimes A_2,B_1\otimes B_2]
}$$
where $\theta$ is the lax structure of $[-,-]$ in $\dgVect$.

The colax structure of $\rhd_\circ$ is given by $(\alpha,\alpha_0)$ where $\alpha_0:\FF\simeq \FF\rhd_\circ \FF$ and 
$\alpha$ is the unique map of non-unital algebras such that the following square commute
$$\xymatrix{
C_1\otimes C_2\otimes A_1\otimes A_2 \ar[rr]^-{\sigma_{23}}\ar[d]_{\Phi\otimes \Phi}&& (C_1\rhd_\circ A_1)\otimes (C_2\rhd_\circ A_2)\ar[d]^{\Phi}\\
(C_1\otimes C_2)\rhd_\circ (A_1\otimes A_2) \ar[rr]^-{\alpha}&& (C_1\rhd_\circ A_1)\otimes (C_2\rhd_\circ A_2)
}$$
where the $\Phi$s are the universal non-unital measurings.
In terms of elements, $\alpha$ is the unique map such that $\alpha((c_1\otimes c_2)\rhd_\circ(a_1\otimes a_2)) = (c_1\rhd_\circ a_1)\otimes (c_2\rhd_\circ a_2)$.

\subsubsection{Reduction functor}

Let $\dgVect^\$$ and $\dgCoalg_\circ^\$$ be the categories $\dgVect$ and $\dgCoalg_\circ$ viewed as enriched over themselves.
Let also $\dgAlg_\circ^\$$ be $\dgAlg_\circ$ viewed as enriched over $\dgCoalg_\circ$.
We can transfer the enrichment of $\dgCoalg_\circ^\$$ and $\dgAlg_\circ^\$$ along the lax monoidal functor $U:\dgCoalg_\circ\to \dgVect$.
Let $\dgCoalg_\circ^{U\$}$ and $\dgAlg_\circ^{U\$}$ be the resulting categories enriched over $\dgVect$.

\begin{prop}\label{strongforgetfulVectnu}
The reduction maps
$$
\HOM_\circ(C,D) \to [C,D]
\et
\{A,B\}_\circ \to [A,B]
$$
are the strengths of enriched functors over $\dgVect$
$$\xymatrix{
U:\dgCoalg_\circ^{U\$} \ar[r] & \dgVect^\$
}
\et
\xymatrix{
U:\dgAlg_\circ^{U\$} \ar[r] & \dgVect^\$
}$$
\end{prop}
\begin{proof}
Same as in propositions \ref{strongforgetfulVectcoalg} and \ref{strongforgetfulVectalg}.
\end{proof}

\subsubsection{Meta-morphisms and (co)derivations}

We define {\em non-counital meta-morphisms} to be the elements of the coalgebra $\HOM_\circ(C,D)$ and $\{A,B\}_\circ$.
By opposition we shall say that elements of $\HOM(C,D)$ and $\{A,B\}$ are {\em counital meta-morphisms}.
Non-counital meta-morphisms can be composed, evaluted and passed through Sweedler operations the same way as the counital ones. 
We leave all details to the reader. 

\medskip

Recall the canonical inclusions $\HOM(C,D)\subset \HOM_\circ(C,D)$ and $\{A,B\}\subset \{A,B\}_\circ$ 
of propositions \ref{HOMvsHOMnu} and \ref{SHOMvsSHOMnu}.
We have the following concrete characterization of counital meta-morphisms.

\begin{prop}\label{propmetanutou}
\begin{enumerate}
\item A non-counital meta-morphism $f:C\leadsto D$ is counital iff $\epsilon_D(f(c))=\epsilon(f)\epsilon_C(c)$ for every $c\in C$.
\item A non-unital meta-morphism $f:A\leadsto B$ is unital iff $f(e_A)=\epsilon(f)e_B$.
\end{enumerate}
\end{prop}
\begin{proof}
This is a reformulation of propositions \ref{HOMvsHOMnu} and \ref{SHOMvsSHOMnu} using the calculus of meta-morphisms.
\end{proof}

\medskip
We use the same musical notation as for the counital meta-morphisms.

\begin{prop}\label{musicallemmacoalgnu}
Let $C$ and $D$ be two non-counital coalgebras, then the maps $\sharp$ and $\flat=\Psi:\HOM_\circ(C,D)\to [C,D]$ induce:
\begin{enumerate}
\item inverse bijections of sets
$$
\dgCoalg_\circ(C,D)\simeq At_\circ(\HOM_\circ(C,D))
$$
where $At_\circ(\HOM_\circ(C,D))$ is the set of non-counital atoms of $\HOM_\circ(C,D)$;

\item inverse isomorphisms in $\dgVect$
$$
\Prim_f(\HOM_\circ(D,C)) \simeq \Coder(f);
$$

\item and, if $C=D$ inverse Lie algebra isomorphisms 
$$
\Prim(\END_\circ(C))\simeq\Coder(C)
$$
which preserve the square of odd elements.
\end{enumerate}
\end{prop}

\begin{proof}
Same as lemma \ref{musicallemmacoalg}, corollary \ref{corprimicoderiv} and theorem \ref{primiLiecoderiv}.
\end{proof}

\begin{prop}\label{musicallemmaalgnu}
Let $A$ and $B$ be two non-unital algebras, then the maps $\sharp$ and $\flat=\Psi:\{A,B\}_\circ\to [A,B]$ induce:
\begin{enumerate}
\item inverse bijections of sets
$$
\dgAlg_\circ(A,B)\simeq At_\circ(\{A,B\}_\circ)
$$
where $At_\circ(\{A,B\}_\circ)$ is the set of non-counital atoms of $\{A,B\}_\circ$.

\item inverse isomorphisms in $\dgVect$
$$
\Prim_f(\{A,B\}_\circ) = \Der(f);
$$

\item and, if $A=B$ inverse Lie algebra isomorphisms 
$$
\Prim(\{A,A\}_\circ)\simeq\Der(A)
$$
which preserve the square of odd elements.
\end{enumerate}
\end{prop}

\begin{proof}
Same as lemma \ref{musicallemmaalg}, corollary \ref{corprimideriv} and theorem \ref{primiLiederiv}.
\end{proof}

\bigskip
We give now the results of transport of coderivations by the Sweedler operations $\HOM_\circ$, $\{-,-\}_\circ$ and $\rhd_\circ$ (the result for the other operations are the usual ones).

\begin{prop}\label{uniqueextensioncodernu}
If $d_1^\sharp$ and $d_2^\sharp$ are the primitive elements of $\END_\circ(C)$ and $\END_\circ(D)$ associated to coderivations $d_1$ and $d_2$ of $C$ and $D$ then
$\HOM_\circ(C,d_2^\sharp)^\flat$ and $\HOM_\circ(d_1^\sharp,D)^\flat$ are coderivations and 
$d=\HOM_\circ(C,d_2^\sharp)^\flat - \HOM_\circ(d_1^\sharp,D)^\flat$ is the unique coderivation such that the square
$$\xymatrix{
\HOM_\circ(C,D) \ar[d]_{\Psi=\flat}\ar@^{>}[rrrr]^-{d} &&&&  \HOM_\circ(C,D)\ar[d]^{\Psi=\flat} \\
[C,D] \ar@^{>}[rrrr]^-{\hom(C,d_2)-\hom(d_1,D)}  &&&& [C,D].
}$$
commutes.
Equivalently, $d$ is the unique coderivation such that
$$
d(h)^\flat = d_2h^\flat-h^\flat d_1 \ (-1)^{|h||d_1|}
$$
for any $h\in \HOM_\circ(C,T^\vee(X))$.
\end{prop}
\begin{proof}
Same as proposition \ref{uniqueextensioncoder}.
\end{proof}

We leave to the reader the cofree and conilpotent variations of this result.

\medskip
\begin{prop}\label{uniqueextensionderSHOMnu}
If $d_1^\sharp$ and $d_2^\sharp$ are the primitive elements of $\{A,A\}_\circ$ and $\{B,B\}_\circ$ associated to derivations $d_1$ and $d_2$ of $A$ and $B$ then
$\{A,d_2^\sharp\}_\circ^\flat$ and $\{d_1^\sharp,B\}_\circ^\flat$ are coderivations and 
$d=\{A,d_2^\sharp\}_\circ^\flat - \{d_1^\sharp,B\}_\circ^\flat$ is the unique coderivation such that the square
$$\xymatrix{
\{A,B\}_\circ \ar[d]_{\Psi=\flat}\ar@^{>}[rrrr]^-{d} &&&&  \{A,B\}_\circ\ar[d]^{\Psi=\flat} \\
[A,B] \ar@^{>}[rrrr]^-{\hom(A,d_2)-\hom(d_1,B)}  &&&& [A,B].
}$$
commutes.
Equivalently, $d$ is the unique coderivation such that
$$
d(h)^\flat = d_2h^\flat-h^\flat d_1 \ (-1)^{|h||d_1|}
$$
for any $h\in \{A,B\}_\circ$.
\end{prop}
\begin{proof}
Same as proposition \ref{uniqueextensionderSHOM}.
\end{proof}

If $A=T_\circ(X)$ is free, we have the following strengthening of the previous result.

\begin{prop}\label{uniqueextensionderSHOMnufree}
If $d_1^\sharp$ and $d_2^\sharp$ are the primitive elements of $\{T_\circ(X),T_\circ(X)\}_\circ$ and $\{B,B\}_\circ$ associated to derivations $d_1$ and $d_2$ of $A$ and $B$ then
the coderivation $d=\{T_\circ(X), d_2^\sharp\}_\circ^\flat - \{d_1^\sharp, B\}_\circ^\flat$ is the unique coderivation on $\{T_\circ(X),B\}_\circ=T^\vee_\circ([X,B])$ such the following square commutes
$$\xymatrix{
T^\vee_\circ([X,B]) \ar[d]_{\Psi}\ar@^{>}[rrrr]^-{d} &&&&  T^\vee_\circ([X,B])\ar[d]^{q} \\
[T_\circ(X),B] \ar@^{>}[rrrr]^-{\hom(i, d_2) - \hom(d_1i, B)}  &&&& [X,B].
}$$
where $i:X\to T_\circ(X)$ is the generating map and $q:T^\vee_\circ([X,B])\to [X,B]$ is the cogenerating map.
Equivalently, $d$ is the unique coderivation such that
$$
q(d(h)) = d_2h^\flat i-h^\flat d_1i \ (-1)^{|h||d_1|}
$$
for any $h\in \{T_\circ(X),B\}_\circ=T^\vee_\circ([X,B])$.
\end{prop}
\begin{proof}
Same as proposition \ref{uniqueextensionderSHOMfree}.
\end{proof}

Recall from corollary \ref{nilradicalcodernu} that a coderivation of a non-counital coalgebra preserves the radical.
We deduce the conilpotent form of the previous result.

\begin{prop}\label{uniqueextensionderSHOMnufreeconil}
If $d_1^\sharp$ and $d_2^\sharp$ are the primitive elements of $\{T_\circ(X),T_\circ(X)\}_\circ$ and $\{B,B\}_\circ$ associated to derivations $d_1$ and $d_2$ of $A$ and $B$ then
the coderivation $d=R^c\{T_\circ(X), d_2^\sharp\}_\circ^\flat - R^c\{d_1^\sharp, B\}_\circ^\flat$ is the unique coderivation on $R^c\{T_\circ(X),B\}_\circ=T^c_\circ([X,B])$ such the following square commutes
$$\xymatrix{
T^c_\circ([X,B]) \ar[d]_{\Psi}\ar@^{>}[rrrr]^-{d} &&&&  T^c_\circ([X,B])\ar[d]^{q} \\
[T_\circ(X),B] \ar@^{>}[rrrr]^-{\hom(i, d_2) - \hom(d_1i, B)}  &&&& [X,B].
}$$
where $i:X\to T_\circ(X)$ is the generating map and $q:T^c_\circ([X,B])\to [X,B]$ is the cogenerating map.
Equivalently, $d$ is the unique coderivation such that
$$
q(d(h)) = d_2h^\flat i-h^\flat d_1i \ (-1)^{|h||d_1|}
$$
for any $h\in R^c\{T_\circ(X),B\}_\circ=T^c_\circ([X,B])$.
\end{prop}

\medskip
\begin{prop}\label{uniqueextensionderSproductnu}
If $d_1^\sharp$ and $d_2^\sharp$ are the primitive elements of $\END_\circ(C)$ and $\{A,A\}_\circ$ associated to a coderivation $d_1$ of $C$ and a derivation $d_2$ of $A$ then
$(d_1^\sharp\rhd_\circ A)^\flat$ and $(C\rhd_\circ d_2^\sharp)^\flat$ are derivations of $C\rhd_\circ A$ and 
$d=(d_1^\sharp\rhd_\circ A)^\flat + (C\rhd_\circ d_2^\sharp)^\flat$ is the unique derivation of $C\rhd_\circ A$ such that the square
$$\xymatrix{
C\otimes  A\ar[d]_{\Phi}\ar@^{>}[rrr]^-{d_1\otimes A+ C\otimes d_2} &&& C \otimes A  \ar[d]^{\Phi} \\
C\rhd_\circ A \ar@^{>}[rrr]^-{d}  &&& C\rhd_\circ A.
}$$
commutes.
Equivalently, $d$ is the unique coderivation such that
$$
d(c\rhd_\circ a) = (d_1c)\rhd_\circ a + c\rhd_\circ (d_2a)\ (-1)^{|c||d_2|}
$$
for any $c\otimes a\in C\otimes A$.
\end{prop}
\begin{proof}
Same as proposition \ref{uniqueextensionderSproduct}.
\end{proof}

If $A=T_\circ(X)$ we have the following strengthening.

\begin{prop}\label{uniqueextensionderSproductfreenu}
If $d_1^\sharp$ and $d_2^\sharp$ are the primitive elements of $\END_\circ(C)$ and $\{T_\circ(X),T_\circ(X)\}_\circ$ associated to a coderivation $d_1$ of $C$ and a derivation $d_2$ of $T(X)$ then
$d=(d_1^\sharp\rhd_\circ T_\circ(X))^\flat + (C\rhd_\circ d_2^\sharp)^\flat$ is the unique derivation of $C\rhd_\circ T_\circ(X)=T_\circ(C\otimes X)$ such that the square
$$\xymatrix{
C\otimes  X\ar[d]_{j}\ar@^{>}[rrr]^-{d_1\otimes i+ C\otimes d_2i} &&& C \otimes T_\circ(X)  \ar[d]^{\Phi} \\
T_\circ(C\otimes X) \ar@^{>}[rrr]^-{d}  &&& T_\circ(C\otimes X).
}$$
commutes ($j:C\otimes X\to T_\circ(C\otimes X)$ is the generating map).
Equivalently, $d$ is the unique coderivation such that
$$
d(c\rhd_\circ x) = d_1(c)\rhd_\circ x + c\rhd_\circ d_2(x)\ (-1)^{|c||d_2|}
$$
for any $c\otimes x\in C\otimes X$.
\end{prop}
\begin{proof}
Same as proposition \ref{uniqueextensionderSproductfree}.
\end{proof}

\subsubsection{Strong (co)monadicity}

We have a non-(co)unital analog of theorems \ref{strongcomonadicitycoalg} and \ref{strongmonadicityalg}.

\begin{thm}\label{strongmonadicitynu}
The adjunction $U\dashv T^\vee_\circ$ enriches into a strong lax monoidal comonadic adjunction
$$\xymatrix{
U: \dgCoalg_\circ^\$ \ar@<.6ex>[r]& \dgVect^{T^\vee_\circ\$}:T^\vee_\circ. \ar@<.6ex>[l]
}$$
The adjunction $T\dashv U$ enriches into a strong colax monoidal monadic adjunction
$$\xymatrix{
T_\circ: \dgVect^{T^\vee_\circ \$} \ar@<.6ex>[r]& \dgAlg_\circ^\$:U. \ar@<.6ex>[l]
}$$
\end{thm}

In consequence, we can construct the hom coalgebras as equalizers in $\dgCoalg$
$$\xymatrix{
\HOM_\circ(C,D)\ar[r]^-{\Delta'} &T^\vee_\circ([C,D]) \ar@<.6ex>[rr]^-{T^\vee_\circ(\alpha)}\ar@<-.6ex>[rr]_-{T^\vee_\circ(\beta)}&& T^\vee_\circ([C,D\otimes D])
}$$
where for $f:C\to D$, we put $\alpha (f)=(f\otimes f)\Delta_C:C\to D\otimes D$ and $\beta (f)=\Delta_Df:C\to D\otimes D$,
and 
$$\xymatrix{
\{A,B\}_\circ\ar[r]^-{m'} &T^\vee_\circ([A,B]) \ar@<.6ex>[rr]^-{T^\vee_\circ(\alpha)}\ar@<-.6ex>[rr]_-{T^\vee_\circ(\beta)}&& T^\vee_\circ([A\otimes A,B])
}$$
where for $f:A\to B$, we put $\alpha (f)=m_B(f\otimes f):A\otimes A\to B$ and $\beta (f)=fm_A:A\otimes A\to B$.

\bigskip

Moreover, we have the following result, which make sense of more distinguished isomorphisms.
Recall from corollaries \ref{adjunctionunitalnonunital} and \ref{adjunctioncounitalnoncounital} the adjunctions
$U_\epsilon :\dgCoalg \rightleftarrows \dgCoalg_\circ {:(-)_+}$ and 
$(-)_+ :\dgAlg_\circ \rightleftarrows \dgAlg {:U_e}$.
Both functors $U$ are monoidal, in particular we can transfer the enrichement of each categories along the lax monoidal right adjoint
$(-)_+:\dgCoalg_\circ\to \dgCoalg$.

\begin{thm}\label{relativestrongmonadicity}
The adjunction $U_\epsilon\dashv (-)_+$ enriches into a strong lax monoidal adjunction
$$\xymatrix{
U_\epsilon: \dgCoalg^\$ \ar@<.6ex>[r]& \dgCoalg_\circ^{\$_+}:(-)_+ \ar@<.6ex>[l]
}$$
where the functor $U_\epsilon$ is strongly faithful.
The adjunction $(-)_+\dashv U_e$ enriches into a strong colax monoidal adjunction
$$\xymatrix{
(-)_+: \dgAlg_\circ^{\$_+} \ar@<.6ex>[r]& \dgAlg^\$:U_e \ar@<.6ex>[l]
}$$
where the functor $U_e$ is strongly faithful.
\end{thm}
\begin{proof}
As for theorems \ref{strongcomonadicitycoalg} and \ref{strongmonadicityalg} but 
the strength of $U_\epsilon: \dgCoalg\to \dgCoalg_\circ$ is given by the lax modular structure
$$
U_\epsilon C\otimes U_\epsilon D= U_\epsilon(C\otimes D);
$$
the strength of $(-)_+: \dgCoalg_\circ\to \dgCoalg$ is given by the colax modular structure
$$
\HOM(C,D_{+}) = \HOM_{\circ}(U_\epsilon C,D)_{+}
$$
of proposition \ref{HomCplus}
which is also the strength of the adjunction $U_\epsilon\dashv (-)_+$;
the strength of $(-)_+: \dgAlg \to \dgAlg_\circ$ is given by the lax modular structure
$$
C\rhd A_{+}=(U_\epsilon C\rhd_{\circ} A)_{+}
$$
of proposition \ref{Sproductcoalgnualg};
the strength of $U_e: \dgAlg_\circ \to  \dgAlg$ is given by the colax modular structure
$$
U_e[C,A]=[U_\epsilon C,U_eA]
$$
and the strength of the adjunction $ (-)_+\dashv U_e$ is given by the isomorphism
$$
\{A_{+},B\}=\{A,B\}_{\circ +}
$$
of proposition \ref{SHOMnualgalg}.

Finally, the faithfulness of functors $U_\epsilon$ and $U_e$ is a consequence of the injections of coalgebras
$$
\xymatrix{\HOM(C,D)\ar[r]& \HOM_\circ(U_\epsilon C,U_\epsilon D)} \et 
\xymatrix{\{A,B\}\ar[r]& \{U_eA,U_eB\}_\circ}
$$
of propositions \ref{HOMvsHOMnu} and \ref{SHOMvsSHOMnu}.
\end{proof}

\begin{rem}
It is in fact possible to prove that these adjunctions are strongly comonadic and monadic over $\dgCoalg$.
\end{rem}

\subsection{The pointed context}\label{pointedsweedlertheory}

In this section we are going to transpose all the results of the previous section through the equivalences
$$\xymatrix{
(-)_-:\dgCoalg_\bullet \ar@<.6ex>[r] & \dgCoalg_\circ \ar@<.6ex>[l]{\ :(-)_+}
}$$
$$\xymatrix{
(-)_-:\dgAlg_\bullet \ar@<.6ex>[r] & \dgAlg_\circ \ar@<.6ex>[l]{\ :(-)_+.}
}$$
Most of this section is logically equivalent to the previous one, but it is convenient in practice to have both languages of pointed an non-(co)untial algebras, so we repeat the statements. Also some changes are made: we adapt the notions of comorphism and of (co)measuring to avoid dealing with smash products, and the comparison of pointed meta-morphisms with unpointed ones is not the one coming from the equivalence with non-counital meta-morphisms.

\medskip
Notice first that the forgetful adjunctions
$$\xymatrix{
U:\dgCoalg_\circ \ar@<.6ex>[r] & \dgVect \ar@<.6ex>[l]{\ :T^\vee_\circ}
}
\et
\xymatrix{
T_\circ:\dgVect \ar@<.6ex>[r] & \dgAlg_\circ \ar@<.6ex>[l]{\ :U}
}$$
are replaced by the adjunctions
$$\xymatrix{
U_-:= U(-)_-:\dgCoalg_\bullet \ar@<.6ex>[r] & \dgVect \ar@<.6ex>[l]{\ :(T^\vee_\circ-)_+ = T^\vee_\bullet}
}$$
$$\xymatrix{
T_\bullet = (T_\circ-)_+ :\dgVect \ar@<.6ex>[r] & \dgAlg_\bullet \ar@<.6ex>[l]{\ :U(-)_- =: U_-.}
}
$$
We will then construct six functors:
\begin{center}
\begin{tabular}{lrl}
\rule[-2ex]{0pt}{4ex} the smash product of pointed coalgebras & $\wedge$&$\!\!\!\!: \dgCoalg_\bullet\times \dgCoalg_\bullet \to \dgCoalg_\bullet$, \\
\rule[-2ex]{0pt}{4ex} the pointed coalgebra internal hom & $\HOM_\bullet$&$\!\!\!\!:(\dgCoalg_\bullet)^{op}\times \dgCoalg_\bullet \to \dgCoalg_\bullet$, \\
\rule[-2ex]{0pt}{4ex} the pointed Sweedler hom & $\{-,-\}_\bullet$&$\!\!\!\!:\dgAlg_\bullet^{op}\times \dgAlg_\bullet \to \dgCoalg_\bullet$, \\
\rule[-2ex]{0pt}{4ex} the pointed Sweedler product & $\rhd_\bullet$&$\!\!\!\!:\dgCoalg_\bullet\times \dgAlg_\bullet \to \dgAlg_\bullet$, \\
\rule[-2ex]{0pt}{4ex} the pointed convolution product & $[-,-]_\bullet$&$\!\!\!\!:(\dgCoalg_\bullet)^{op}\times \dgAlg_\bullet \to \dgAlg_\bullet$, \\
\rule[-2ex]{0pt}{4ex} and the smash product of pointed algebras & $\wedge$&$\!\!\!\!:\dgAlg_\bullet\times \dgAlg_\bullet \to \dgAlg_\bullet$.
\end{tabular}
\end{center}
such that
\begin{itemize}
\item the category $(\dgCoalg_\bullet,\wedge, \HOM_\bullet)$ is locally presentable symmetric monoidal closed and comonadic over $\dgVect$
\item and the category $(\dgAlg_\bullet,\{-,-\}_\bullet,\rhd_\bullet, [-,-]_\bullet,\wedge)$ is locally presentable, enriched, bicomplete and symmetric monoidal over $\dgCoalg_\bullet$, and monadic over $\dgVect$.
\end{itemize}
We will also distinguish the following isomorphism:
\begin{center}
\begin{tabular}{c}
\rule[-2ex]{0pt}{4ex} $ T_\bullet^\vee [U_-C,X]\simeq \HOM_\bullet(C,T_\bullet^\vee (X))$ \\
\rule[-2ex]{0pt}{4ex} $ \{T_\bullet X,A\}\simeq T_\bullet^\vee([X,U_-A])$\\
\rule[-2ex]{0pt}{4ex} $C\rhd_\bullet T_\bullet(X) \simeq T_\bullet(U_-C\otimes X)$
\end{tabular}
\end{center}
which are the data to strengthen the adjunctions $U_-\dashv T_\bullet^\vee$ and $T_\bullet\dashv U_-$ over $\dgCoalg_\bullet$:
\begin{itemize}
\item the adjunction $U_-:\dgCoalg_\bullet \rightleftarrows \dgVect:T_\bullet^\vee $ is strongly comonadic over $\dgCoalg_\bullet$
\item and the adjunction $T_\bullet:\dgVect \rightleftarrows \dgAlg_\bullet:U_-$ is strongly monadic over $\dgCoalg_\bullet$.
\end{itemize}

\begin{rem}
If we use the equivalence of categories $(-)_+:\dgVect\simeq \dgVect_\bullet:(-)_-$ of section \ref{pointedvectorspaces}, there is a commutative diagram of adjunctions
$$\xymatrix{
\dgCoalg_\bullet \ar@<.6ex>[rr]^-{(-)_-} \ar@<-.6ex>[dd]_{U}&& \dgCoalg_\circ \ar@<.6ex>[ll]^-{(-)_+}\ar@<-.6ex>[dd]_{U}\\
&&\\
\dgVect_\bullet \ar@<.6ex>[rr]^-{(-)_-} \ar@<-.6ex>[uu]_{T^\vee_\bullet(-)_-}&& \dgVect \ar@<.6ex>[ll]^-{(-)_+} \ar@<-.6ex>[uu]_{T^\vee_\circ}
}$$
where the horizontal ones are equivalences.
The forgetful functor $U_-:\dgCoalg_\bullet \to \dgVect$ is then replaced by the functor $U_\bullet:\dgCoalg_\bullet \to \dgVect_\bullet$ forgetting the coassociative coproduct but not the counit and coaugmentation.
Some formulas becomes nicer with $U_\bullet$ but we have chosen to work mainly with $U_-$ because it is closest to what already exist in the litterature when manipulating pointed (co)algebras.
\end{rem}

\subsubsection{Presentability, comonadicity and cofree pointed coalgebra}

\begin{thm}
The category $\dgCoalg_\bullet$ is finitary presentable. 
The $\omega$-compact objects are the finite dimensional coalgebras.
\end{thm}
\begin{proof}
Direct from theorem \ref{noncounitalpresentability}.
\end{proof}

\begin{thm}
The forgetful functor $U_-:\dgCoalg_\bullet\to \dgVect:T^\vee_\bullet$ is comonadic.
\end{thm}
\begin{proof}
Direct from theorem \ref{noncounitalcomonadicity}.
\end{proof}

\begin{prop}
For an (dg-)vector space $X$, there exists a counital dg-coalgebra isomorphism $T^\vee_\bullet(X) \simeq T^\vee(X)$.
\end{prop}
\begin{proof}
This is actually what is proven in proposition \ref{isocofreeconil}.
The canonical coaugmentation of $\FF=T^\vee(0)\to T^\vee(X)$ is the image by $T^\vee$ of the zero map $0\to X$.
\end{proof}

\medskip
Recall from example \ref{ptensorcoalg} the pointed tensor coalgebra $T^c_\bullet(X)$ on a dg-vector space $X$
and from proposition \ref{radical} the notion of radical $R^c$ of a non-unital coalgebra.

\begin{prop}\label{nilradicalcofreepointed}
We have $T^c_\bullet(X)=R^cT^\vee_\bullet(X)$ for any vector space $X$.
\end{prop}

\begin{proof}
Direct from proposition \ref{nilradicalcofree}.
\end{proof}

\begin{prop}
The adjunction $U:\dgCoalgnilbullet \rightleftarrows \dgVect: T_\bullet^c$ is comonadic.
\end{prop}
\begin{proof}
Direct from proposition \ref{comonadicTc}.
\end{proof}

\subsubsection{Internal hom}

The {\it pointed hom object} between two pointed coalgebras $C$ and $D$ is the pointed coalgebra
$\HOM_\bullet(C,D)$  defined by putting 
$$
\HOM_\bullet(C,D)=\HOM_{\circ}(C_{-},D_{-})_{+}.
$$
Recall from section \ref{monoidalstructurecoalgebras} that the smash product of pointed coalgebras is defined by $C\wedge D := (C_-\otimes D_-)_+$ and correspond to $\otimes$ through the equivalence $\dgCoalg_\bullet\simeq \dgCoalg_\circ$.
$\HOM_\bullet$ correspond to $\HOM_\circ$ through the same equivalence and we deduce that $\HOM_\bullet$ is right adjoint to $\wedge$, \ie that we have canonical bijections between
\begin{center}
\begin{tabular}{lc}
\rule[-2ex]{0pt}{4ex} pointed coalgebra maps & $C\wedge D\to E$,\\
\rule[-2ex]{0pt}{4ex} pointed coalgebra maps & $ D\to \HOM_\bullet(C,E)$\\
\rule[-2ex]{0pt}{4ex} and pointed coalgebra maps & $ C\to \HOM_\bullet(D,E)$.
\end{tabular}
\end{center}

We proved in proposition \ref{smashproductofcoalgebras} that the smash product could be written as a pushout. Correspondingly, the pointed coalgebra $\HOM_\bullet(C,D)$ can be also described as a pull-back.

\begin{prop}\label{pHOMvsHOM}
We have a pullback square of coalgebras
$$\xymatrix{
\HOM_\bullet(C,D) \ar[d]_-\epsilon \ar[rr] && \HOM(C,D) \ar[d]^-{ \HOM(e_C,D)} \\
\FF \ar[rr]^-{e_D}  && D=\HOM(\FF,D).
}$$ 
In particular, $\HOM_\bullet(C,D)$ is a sub-coalgebra of $\HOM(C,D)$.

The coalgebra $\HOM_\bullet(C,D)$ is pointed by the map $\FF\to \HOM_\bullet(C,D)$
lifting the map $\HOM(\epsilon_C,e_D):\FF=\HOM(\FF,\FF)\to \HOM(C,D)$.
\end{prop}

\begin{proof}
If $(E,e)$ is a pointed  coalgebra, then it follows from proposition
\ref{smashproductofcoalgebras} that there is a bijection
between the maps of pointed coalgebras $f:E\wedge C\to D$ 
and the maps of coalgebras $g:E\otimes C\to D$ 
fitting in the commutative square
$$\xymatrix{
E \oplus C \ar[d]_{(\epsilon_E,\epsilon_C)} \ar[rr]^-{(E \otimes e_C, e \otimes  C)} && E\otimes C \ar[d]^g  \\
 \ar[rr]^-{e_D}  \FF && D
}$$
But the square commutes if and only if  the following diagram commutes,
$$\xymatrix{
C  \ar[d]_{\epsilon_C} \ar[rr]^-{e \otimes C} &&E\otimes C \ar[d]_g  && \ar[ll]_-{E\otimes e_C} E \ar[d]^{\epsilon_E}\\
\FF \ar[rr]^-{e_D} && D && \FF \ar[ll]_-{e_D}
}$$ 
iff the following diagram commutes,
$$\xymatrix{
\FF \ar[rrr]^-{e} \ar[d]_{e_D}  &&&  E \ar[d]_h   \ar[rrr]^-{\epsilon_E}&&& \FF \ar[d]^{e_D}  \\
D  \ar[rrr]^(0.4){ \HOM(\epsilon_C,D)}    &&& \HOM(C,D)  \ar[rrr]^(0.6){ \HOM(e_C,D)} &&& D
}$$ 
where $h=\lambda^2(g)$. Let us put $K=\HOM(C,D)\times_D \FF$.
Composition with the first projection $p_1:K\to \HOM(C,D)$
induces bijection between the maps of coalgebras
$k:E\to  K$ and the maps of coalgebras $h=p_1k:E\to \HOM(C,D)$ for which
the square on the right hand side commutes. 
There a unique map of coalgebras $e_K:\FF\to K$ such that
$p_1e_K= \HOM(\epsilon_C,D)e_D$ and $p_2e_K=1_\FF$.
Moreover, we have $ke_K=e$ iff the square on the left hand side commutes.
Hence the composite $f\mapsto g\mapsto  h\mapsto k$  is a bijection between the 
maps of pointed coalgebras $f:E\wedge C\to D$ and the maps
of pointed coalgebras $k:E\to K$. The bijection is natural and it follows by Yoneda lemma
that $K=\HOM_\bullet(C,D)$.

Then $\HOM_\bullet(C,D)$ is a sub-coalgebra of $\HOM(C,D)$ by lemma \ref{directandinverseimage5}.
\end{proof}

\begin{rem}
Recall that $\FF$ is the terminal object in $\dgCoalg$, hence the unit for the cartesian product.
By proposition \ref{closedmonoidcoalg4} we have $\HOM(C,\FF)=\FF$, hence we can replace $\HOM(\FF,D)$ by $\HOM(\FF,D)\times \HOM(C,\FF)$ in the square of proposition \ref{pHOMvsHOM}. This way it becomes analog to the cartesian square of lemma \ref{pointedhomcart}.
\end{rem}

\begin{rem}\label{rempHOMvsHOM}
In the proof of proposition \ref{pHOMvsHOM}, we have constructed the inclusion $\HOM_\bullet (C,D)\to \HOM(C,D)$ as the unique coalgebra map such that the triangle
$$\xymatrix{
\HOM_\bullet(C,D)\ar[rr]\ar[rrd]_\Psi&&\HOM(C,D) \ar[d]^{\Psi}\\
&&[C,D]
}$$
commutes.
\end{rem}

\begin{thm} \label{closedmonoidpointedcoalg} 
The symmetric monoidal category $(\dgCoalg_\bullet,\wedge, \FF_{+})$ is closed and the hom object is the pointed coalgebras $\HOM_\bullet(C,D)$.
\end{thm}
\begin{proof}
This follows from theorem \ref{closedmonoidunitalcoalg}.
\end{proof}

The counit of the adjunction $C\wedge(-)\dashv \HOM_{\bullet}(C,-)$ is the {\it evaluation map} $\mathbf{ev} :\HOM_\bullet(C,D)\wedge C\to D$.
For any pointed coalgebra $E$ and any map of pointed coalgebras $f:E\wedge C\to D$ there exists a unique map of non-counital coalgebras $g:E\to \HOM_\bullet(C,D)$ such that ${\bf ev}(g\wedge C)=f$.
$$\xymatrix{
&& \HOM_\bullet(C, D) \wedge C  \ar[d]^-{{\bf ev}} \\
E \wedge C \ar@{-->}[urr]^-{g\wedge C}  \ar[rr]^-{f}    && D
}$$
We shall put ${\Lambda}^2(f):=g$.
In addition, we shall put ${\Lambda}^1(f):={\Lambda}^2(f\sigma):C\to \HOM_\bullet (E,D)$, where $\sigma:C\wedge E\to E\wedge C$ is the symmetry of $\wedge$.

\begin{rem}
The coalgebras $\END(C)_\bullet$ are examples of pointed bialgebras in the sense of section \ref{pointed bialgebras}.
Moreover, using the equivalence of proposition \ref{equivBilagpointed}, we have an isomorphism of pointed bialgebras
$$
HOM_{\bullet}(C_{+},C_{+})=\FF\times HOM_{\circ}(C,C).
$$
\end{rem}

\bigskip

Let us now turn to the notion of comorphism. 
Recall that maps $E\wedge C\to D$ of pointed vector spaces are in bijection with map of pointed vector spaces $E\to [C,D]_\bullet$ and linear maps $E_-\to [C_-,D_-]$. We shall say that maps $E\to [C,D]_\bullet$ and $E_-\to [C_-,D_-]$ corresponding to pointed coalgebra map $E\wedge C\to D$ are respectively a {\em pointed comorphism in $\dgVect_\bullet$} and a a {\em pointed comorphism in $\dgVect$}.
We leave to the reader the definition of couniversal pointed comorphism in both contexts.
The couniversal comorphism are given by the maps $\Psi_\bullet:\HOM(C,D)\to [C,D]_\bullet$ in $\dgVect_\bullet$ and $\HOM_\bullet(C,D)_-\to [C_-,D_-]$ in $\dgVect$.
We shall prove in proposition \ref{pointedreduction2} that $\Psi_\bullet$ is part of an enriched functor.

We will also introduce an equivalent notion of {\em expanded pointed comorphism} such that the couniversal comorphism is simply given by a map $\HOM_\bullet(C,D)\to [C,D]$. 
Recall that map of pointed coalgebras $f:E\wedge C\to D$ is equivalent to a map of coalgebras $g:E\otimes C\to D$ entering a commutative square 
\begin{center}
\hfill $\vcenter{\xymatrix{
E\oplus C \ar[r]\ar[d]_{E\otimes e_C\oplus e_E\otimes C} & \FF\ar[d]^{e_D}\\
E\otimes C \ar[r]^-g &D.
}}$
\hfill
(P)
\end{center}
We shall say that a map $g:E\otimes C\to D$ such that $P$ commutes is an {\em expanded pointed coalgebra map}.
This is equivalent to say that is is a map of coalgebras and an expanded pointed map of pointed vector space in the sense of section \ref{pointedvectorspaces}.

The commutation of (P) is equivalent to the commutation of
\begin{center}
\hfill $\vcenter{\xymatrix{
\FF	\ar[rr]^-{e_E}\ar[d]_{e_D} && E	\ar[d]^{\lambda g}\ar[rr]^-{\epsilon_E} && \FF \ar[d]^{e_D}\\
D \ar[rr]^-{[\epsilon_C,D]} && [C,D] \ar[rr]^-{[e_C,D]} && D.
}}$
\hfill
(P')
\end{center}
We shall say that a map $E\to [C,D]$ satisfying condition (P') is a {\em expanded pointed comorphism}.

For three pointed coalgebras, $C$, $D$ and $E$, we then have natural bijections between
\begin{center}
\begin{tabular}{lrl}
\rule[-2ex]{0pt}{4ex} pointed coalgebra maps & $E\wedge C\to D$\\
\rule[-2ex]{0pt}{4ex} pointed comorphisms in $\dgVect_\bullet$ & $E\to [C,D]_\bullet$\\
\rule[-2ex]{0pt}{4ex} pointed comorphisms in $\dgVect$ & $E_-\to [C_-,D_-]$\\
\rule[-2ex]{0pt}{4ex} expanded pointed coalgebra maps & $E\otimes C\to D$\\
\rule[-2ex]{0pt}{4ex} and expanded pointed comorphisms & $E\to [C,D]$.
\end{tabular}
\end{center}
In particular the pointed coalgebra map ${\bf ev}:\HOM_\bullet(C,D)\wedge C\to D$ corresponds to a map $\Psi:\HOM_\bullet(C,D)\to [C,D]$
which is just the composition $\HOM_\bullet(C,D)\to [C,D]_\bullet\to [C,D]$.
We leave to the reader the definition of a couniversal expanded pointed comorphism in this context and the proof that $\Psi$ is couniversal.
We shall prove in section \ref{pointedreductionsection} that $\Psi$ is also part of an enriched functor.

The following lemmas help compare the notions of pointed, co-unital and non-counital comorphisms.

Recall that if $C$ and $D$ are pointed, we have a canonical decomposition $[C,D]=\FF\oplus [C_-,\FF]\oplus D_-\oplus [C_-,D_-]$.
Then, the equation $\epsilon_Df=\epsilon_C$ of proposition \ref{propmetanutou} implies that the couniversal (counital) comorphism $\Psi:\HOM(C,D)\to [C,D]$ factors through $[C,D]_{\sf half\bullet} = \FF\oplus D_-\oplus [C_-,D_-]$.

\begin{lemma}\label{comparisonallmeta-morphismscoalg}
We have a commutative diagram
$$\xymatrix{
\HOM_\bullet(C,D) \ar[rr]\ar[d]_{\Psi_\bullet} && \HOM(C,D)\ar[d]\ar[rrd]^{\Psi}\ar[rr]&& \HOM_\circ(C,D)\ar[d]^{\Psi_\circ}\\
[C,D]_\bullet \ar[rr] && [C,D]_{\sf half\bullet} \ar[rr] && [C,D].
}$$
where the vertical maps are the canonical inclusions.
In particular, the couniversal expanded pointed comorphism is the total diagonal
$\HOM_\bullet(C,D)\to \HOM(C,D)\to [C,D]$.
\end{lemma}
\begin{proof}
$\Psi:\HOM_\bullet(C,D)\to [C,D]_\bullet\to [C,D]$ is in particular a measuring, hence the left commutative square.
$\Psi:\HOM(C,D)\to [C,D]_{\sf half\bullet}\to [C,D]$ is in particular a non-counital measuring, hence the right commutative square.
Then we need to prove that the top horizontal maps are the inclusions of propositions \ref{HOMvsHOMnu} and \ref{pHOMvsHOM}.
This comes from remarks \ref{remHOMvsHOMnu} and \ref{rempHOMvsHOM}.
\end{proof}

\medskip

\begin{prop} \label{closedmonoidcoalg4bullet}  
For every non-unital coalgebra $C$, we have
\begin{enumerate}
\item $\HOM_\bullet(\FF_+,C)=C$,
\item $\HOM_\bullet(C,\FF)=\FF$ and 
\item $\HOM_\bullet(\FF,C)=\FF$.
\end{enumerate}
\end{prop}
\begin{proof}
Same as in proposition \ref{closedmonoidcoalg4}, using that $\FF$ is both initial and terminal in $\dgCoalg_\bullet$.
\end{proof}

$\HOM_\bullet (C,D)$ can be described as in proposition \ref{HOMsubcoalg} using the functor $T^\vee_\bullet$.
In case $D$ is pointed cofree we have also the following result, analogous to proposition \ref{exemplehomcoalg}.

\medskip
Let $C$ be a pointed coalgebra and $X$ be a dg-vector space.
Let us simplify the forgetful functor $U_-C$ as $C_-$.
If $q: T^\vee_\bullet([C_-,X])\to [C_-,X]$ is the cogenerating map, then the composite of the maps
$$\xymatrix{
T^\vee_\bullet([C_-,X])\wedge C \simeq T^\vee_\circ([C_-,X])\otimes C_- \ar[rr]^-{q\otimes C} && [C_-,X]\otimes C_- \ar[r]^-{ev} & X
}$$
can be coextended as a map of pointed coalgebras $e: T^\vee_\bullet ([C,X])\wedge C \to T^\vee_\bullet (X)$.

\begin{prop}\label{exemplehomsemicoalgp}
With the previous notations, for any pointed coalgebra $E$ and any map of pointed coalgebras
$f:E\wedge C\to T^\vee_\bullet (X)$, there exists a unique map of pointed coalgebras
$k:E\to  T^\vee_\bullet([C_-,X])$ such that $e(k\wedge C)=f$. 
Thus, $e$ is a strong evaluation $\mathbf{ev}:\HOM_\bullet(C,T^\vee_\bullet(X)) \wedge C\to  T^\vee_\bullet([C_-,X])$ and we have
$$
\HOM_\bullet(C,T^\vee_\bullet(X)) = T^\vee_\bullet ([C_-,X]).
$$
\end{prop}
\begin{proof}
Direct from proposition \ref{exemplehomcoalgnu}.
\end{proof}

\medskip

The strong composition law $\mathbf{c}:\HOM_\bullet(D,E)\wedge \HOM_\bullet(C,D) \to  \HOM_\bullet(C,E)$ is derived from the strong composition of non-counital coalgebras. It can be defined as $\Lambda^3({\bf ev}^2)$ where ${\bf ev}^2$ is the map
$$\xymatrix{
\HOM_\bullet(D,E)\wedge \HOM_\bullet(C,D)\wedge C \ar[rrr]^-{\HOM_\bullet(D,E)\wedge {\bf ev}}&&&\HOM_\bullet(D,E)\wedge D \ar[r]^-{\bf ev}& E.
}$$
It can also be defined through the notion of pointed comorphism but it is indirect.
Let us consider the commutative square 
$$\xymatrix{
\HOM_\bullet(D,E)\otimes \HOM_\bullet(C,D)\ar[rr]^-{\mathbf{c}}\ar[d]_{\Psi \otimes \Psi }&& \HOM_\bullet(C,E)\ar[d]^{\Psi}\\
[D,E]\otimes [C,D]\ar[rr]^-{c} && [C,E]
}$$
where $c$ is the strong composition in $\dgVect$ and the $\Psi$ are the couniversal pointed comorphisms. 
We leave the reader to check that $c(\Psi\otimes \Psi)$ is a comorphism, the map $\bf c$ is then constructed by the universal property of $\Psi:\HOM_\bullet(C,E)\to [C,E]$. We leave also to the reader the proof that $\bf c$ factors through $\HOM_\bullet(D,E)\wedge \HOM_\bullet(C,D)$.

\medskip
We shall denote $\dgCoalg_\bullet^\$$ the category $\dgCoalg_\bullet$ viewed as enriched over itself.

\begin{prop}
The atoms of $\HOM_\bullet(C,D)$ are in bijection with the maps of pointed coalgebras.
The underlying category of $\dgCoalg_\bullet^\$$ is $\dgCoalg_\bullet$.
\end{prop}
\begin{proof}
Direct from proposition \ref{atomhomnucoalg}.
\end{proof}

We have also the analog of proposition \ref{radical3}.

\begin{prop}  \label{radical3pointed} 
If $C$ is a nilpotent pointed coalgebra, then we have
$$
\HOM_{\bullet}(C,R^cD)\simeq \HOM_{\bullet}(C,D)
$$
for any pointed coalgebra $D$.
In other words, the radical adjunction 
$$\xymatrix{
\iota:\dgCoalgnilbullet \ar@<.6ex>[r] &\dgCoalg_\bullet:R^c \ar@<.6ex>[l]
}$$
is strong.
\end{prop}

\subsubsection{Convolution and measurings}
If $A=(A,\epsilon_A)$ is a pointed algebra and $C=(C,e_C)$ is a pointed coalgebra then the {\it pointed convolution algebra} $[C,A]_\bullet$ is defined by putting 
$$
[C,A]_\bullet=[C_{-},A_{-}]_+.
$$
By equivalence with $\dgAlg_\circ$, it is an algebra for the smash product of pointed vector spaces.

\begin{prop}\label{convolutionpointed}
We have a pullback square of dg-algebras
$$\xymatrix{
[C,A]_\bullet \ar[rr]\ar[d]_\epsilon && [C,A] \ar[d]^{([C,\epsilon_A],[e_C,A])}\\
\FF \ar[rr]^(0.4){(\epsilon_C,e_A)} &&[C,\FF]\times [\FF,A].
}$$
Hence the augmentation ideal of the algebra $[C,A]_\bullet$ is the intersection of the kernel of the maps $[C,\epsilon_A]$ and $[e_C,A]$.
In particular $[C,A]_\bullet$ is a sub-algebra of $[C,A]$.
\end{prop}

\begin{proof} 
We know from lemma \ref{pointedhomcart} that the square is cartesian. We need only to prove that the maps are algebras maps, which is a straightforward computation left to the reader.
\end{proof}

\bigskip

\begin{defi}\label{pointedmeasuring} 
For $C$ a pointed coalgebra and $A$ and $B$ two pointed algebras, we shall say that a map $C\wedge A\to B$ is a {\em pointed measuring} if the corresponding maps $A\to [C,B]_\bullet$ is a map of pointed algebras. As before we will define an expanded version of this notion.

We shall say that a linear map $f:C\otimes A\to B$ is {\it expanded pointed measuring} 
if it is an expanded map of pointed vector space and if the corresponding map $C\wedge A\to B$ is a pointed measuring.
Using the characterization of expanded pointed map of section \ref{pointedvectorspaces}, a map $C\otimes A\to B$ is an expanded pointed measuring iff 
$$
f(c,ab)= f(c^{(1)},a)f(c^{(2)},b) \ (-1)^{|a||c^{(2)}|} \ ,
\qquad
f(c,e_A)=\epsilon(c)e_B\ ,
$$
(which is the measuring condition) and
$$
\epsilon_B(f(c,a)) = \epsilon_C(c)\epsilon_A(a)
\et
f(e_C,a) = \epsilon_A(a)e_B
$$
for every $a\in A$ and $c\in C$. 
\end{defi}

\begin{lemma}\label{pointedmeasuringdetermined}
Let $C$ be a pointed coalgebra and $A$ and $B$ be a pointed algebras.
Then every non-unital measuring $f:C_{-}\otimes A_{-}\to B_-$ can be extended uniquely as a expanded pointed measuring $f':C\otimes A\to B$.
\end{lemma}

\begin{proof} 
Let us write $A= A_-\oplus \FF e_A$, $B= B_-\oplus \FF e_B$ and $C= C_-\oplus \FF e_C$.
We need to define the extension of $f$ to $\FF e_C\otimes A_- \oplus C_-\otimes \FF e_A$ we use the formulas
$$
f'(e_C,a) := \epsilon_A(a)e_B
\et 
f'(c,e_A) := \epsilon_C(c)e_B
$$
for every $a\in A$ and $c\in C$.
Finally we need to check that $\epsilon_B(f(c,a)) = \epsilon_C(c)\epsilon_A(a)$.
If $c\in C_-$ and $a\in A_-$ the $f'(c,a)=f(c,a)\in B_-$ so both terms are zero.
Then form the above relations we have
$$
\epsilon_B(f'(e_C,a)) = \epsilon_A(a)\epsilon_B(e_B)
\et 
\epsilon_B(f'(c,e_A)) = \epsilon_C(c)\epsilon_B(e_B)
$$
which finishes to prove the formula.
\end{proof}

We shall denote the set of expanded pointed measurings $C\otimes A\to B$ by $\cM_{\bullet}(C,A;B)$. 
This defines a functor of three variables,
$$\xymatrix{
\cM_{\bullet}(-,-;-):\dgCoalg_{\bullet}^{op}\times \dgAlg_{\bullet}^{op}\times \dgAlg_{\bullet}\ar[r]& \Set.
}$$
By definition of expanded pointed measurings, this functor is representable in the second variable:
$$
\cM_{\bullet}(C,A;B) = \dgAlg_\bullet(A,[C,B]_\bullet).
$$

Under the equivalence of non-(co)unital (co)algebras and pointed (co)algebras, we have natural bijections
$$
\cM_{\bullet}(C,A;B) = \cM_{\circ}(C_-,A_-;B_-).
$$
In particular we can deduce the representability of $\cM_{\bullet}$ in all variables from that of $\cM_{\circ}$.

\subsubsection{Sweedler product}

We shall say that an expanded pointed measuring $u:C\otimes A\to B$ is {\it universal} if the pair $(B,u)$ is representing the functor $\cM_\bullet(C,A;-): \dgAlg_\bullet\to \Set$.
The {\it pointed Sweedler product} $C\rhd_\bullet  A$ of $A=(A,\epsilon_A)$ by $C=(C,e_C)$ is defined by
putting
$$
C\rhd_\bullet  A =(C_{-}\rhd_{\circ} A_{-})_{+}.
$$
The pointed algebra $C\rhd_\bullet  A$ is the codomain of a universal expanded pointed measuring $u:C\otimes A\to C\rhd_\bullet A$.
Let us put $x\rhd_\bullet y=u(x\otimes y)$ for every $x\in C$ and $y\in A$.

\begin{prop}\label{Sproductpointed}
We have a pushout square of algebras
$$\xymatrix{
A\ar[rr]^-{e_C \rhd A} \ar[d]_{\epsilon_A} && \ar[d] C\rhd A \ar[d] \\
\FF \ar[rr] && C\rhd_\bullet A.
}$$ 
The augmentation $\epsilon:C\rhd_\bullet A\to \FF$ is induced by the map $\epsilon_C \rhd \epsilon_A:C\rhd A\to \FF \rhd \FF=\FF$.
\end{prop}

\begin{proof} 
Let us denote by $E$ the algebra defined by pushout square of algebras
$$\xymatrix{
A\ar[rr]^-{e_C \rhd A} \ar[d]_{\epsilon_A} && \ar[d] C\rhd A \ar[d]^q \\
\FF \ar[rr]^e && E.
}$$ 
Let us show first that $E$ is pointed.
By evaluating the diagram on $e_A\in A$, we find that $e$ is the identity of $E$.
Then by considering the square 
$$\xymatrix{
A\ar[rr]^-{e_C \rhd A} \ar[d]_{\epsilon_A} && \ar[d] C\rhd A \ar[d]^{\epsilon_C\rhd \epsilon_A} \\
\FF \ar[rr]^-{id} && \FF.
}$$ 
we contruct an augmentation $\epsilon:E\to \FF$. 

The algebra $C\rhd A$ is naturally pointed by $(e_C\rhd e_A, \epsilon_C\rhd \epsilon_A)$ and the map $q:C\rhd A\to E$ is a map of pointed algebras. All the other maps of the square are also pointed. This proves that the square is also cartesian in $\dgAlg_\bullet$.

The proof will be finished if we prove that a commutative square
$$\xymatrix{
A\ar[rr]^-{e_C \rhd A} \ar[d]_{\epsilon_A} && \ar[d] C\rhd A \ar[d]^{g} \\
\FF \ar[rr]^-{e_B} && B
}$$
in $\dgAlg_\bullet$ is equivalent to an expanded pointed measuring.
Remark that the condition to be a commutative square in $\dgAlg_\bullet$ reduces to the condition that $g$ preserves the counit.

We need to show that the composition $m:C\otimes A\to C\rhd A\to B$ is an expanded pointed measuring.
It is a measuring because post composition of measurings with algebra maps stay measurings.
It is a pointed map because $m(e_c\otimes e_A) = e_B$ and $\epsilon_B(m(c,a)) = \epsilon_C(c)\epsilon_A(a)$ by hypothesis on $g$.
Let us prove that $e_C\otimes A$ is send to the $e_B$ component of $B$:
we have $m(e_C,a) = g(e_C\rhd a) = \epsilon_A(a)e_B$ by hypothesis on $B$
Let us prove that $C\otimes e_A$ is send to the $e_B$ component of $B$:
we have $m(c,e_A) = \epsilon_C(c)e$ by the measuring property.
This finishes to prove that $m$ is an expanded pointed measuring.
Then the universal property of $E$ says that it represents the functor of expanded pointed measuring.
\end{proof}

In terms of elements, $C\rhd_\bullet A$ is the algebra generated by symbols $c\rhd_\bullet a$ for $c\in C$ and $a\in A$ with the following 
relations:
\begin{itemize}
\item[(d)] the map $(c,a)\mapsto c\rhd_\bullet a$ is a bilinear map of dg-vector spaces, $d(c\rhd_\bullet a) = dc\rhd_\bullet a + c\rhd_\bullet da (-1)^{|c|}$;
\item[(m)] $c\rhd_\bullet (ab) = (c^{(1)}\rhd_\bullet a)(c^{(2)}\rhd_\bullet b) (-1)^{|a||c^{(2)}|}$, for every $c\in C$ and $a,b\in A$;
\item[(u)] $c\rhd_\bullet e_A = \epsilon(c)$, for every $c\in C$;
\item[(a)] and $e_C\rhd_\bullet a = \epsilon(a)$, for every $a\in A$.
\end{itemize}
The augmentation is then given by $\epsilon(c\rhd_\bullet a) = \epsilon(c)\epsilon(a)$.

Equivalently $C\rhd_\bullet A$ is the quotient of the tensor dg-algebra $T(C\otimes A)$ by the relations (m), (u) and (a).
It can also be described as the quotient of the pointed tensor dg-algebra $T_\bullet(C_-\otimes A_-)$ by the relations (m).

\medskip

\begin{lemma} \label{extensionpointedmeasuring} 
Let $A$ be a pointed algebra,  $X$ be a vector space and $i:X\to T_\bullet (X)$ be the inclusion. If $C$ is a pointed coalgebra,
then any linear map $f:C_{-}\otimes X \to A_-$ can be extended uniquely as an expanded pointed measuring $f':C\otimes T_\bullet (X) \to A$.
\end{lemma}

\begin{proof} 
The map $\lambda^1(f):X\to [C_{-},A_-]$ corresponding to 
$f:C_{-}\otimes X\to A_-$ can be extended uniquely as a map of non-unital algebras  $g:T_\circ (X)\to [C_{-},A_-]$.
The non-unital measuring $h:T_\circ (X)\otimes C_{-} \to A$ defined by $g$
can be extended uniquely as an expanded pointed measuring 
$f':T_\bullet (X)\otimes C \to A$ by lemma \ref{pointedmeasuringdetermined}.
\end{proof}

\begin{prop}\label{Sweedleroftensorpointed}
If $C$ is a pointed coalgebra,  $X$ be a vector space and $i$ is the inclusion $C_{-}\otimes X \to T(C_{-}\otimes X)$, 
then the expanded pointed measuring $i':C\otimes T_\bullet(X)\to T_\bullet (C_{-}\otimes X)$ defined above is universal.
Hence we have 
$$
C\rhd_\bullet T_\bullet(X)= T_\bullet(C_{-}\otimes X).
$$
\end{prop}

\begin{proof}
If $A$ is a pointed algebra, then we have a  chain of natural bijections between 
\begin{center}
\begin{tabular}{lc}
\rule[-2ex]{0pt}{4ex} the maps of pointed algebras maps & $h:T_\bullet (C_{-}\otimes X)\to A$,\\
\rule[-2ex]{0pt}{4ex} the linear maps & $C_{-}\otimes X\to A_{-}$,\\
\rule[-2ex]{0pt}{4ex} the linear maps & $X\to [C_{-},A_{-}]$,\\
\rule[-2ex]{0pt}{4ex} the maps of pointed algebras & $T_\bullet(X)\to [C,A]_\bullet $,\\
\rule[-2ex]{0pt}{4ex} the expanded pointed measurings & $k:C\otimes T_\bullet(X)  \to A$.
\end{tabular}
\end{center}
The bijections show that the functor $\cM_\bullet(C,T_\bullet(X);-)$ is represented by the algebra $T_\bullet(C\otimes X)$.
Moreover, if $A=T(C\otimes X)$ and $h$ is the identity, then $k=i'$.
This proves that the pointed measuring $i':C\otimes T_\bullet (X)\to T_\bullet (C_{-}\otimes X)$ is universal,
and hence that we have $C\rhd_\bullet T_\bullet(X)= T_\bullet(C_{-}\otimes X)$.
\end{proof}

\subsubsection{Sweedler hom}

If $A$ and $B$ are pointed algebras, we define
$$
\{A,B\}_\bullet =(\{A_{-},B_{-}\}_{\circ})_{+}.
$$
The pointed coalgebra $\{A,B\}_\bullet $ is also described by the following result.

\begin{prop}\label{pSHOMvsSHOM}
We have a pullback square of coalgebras
$$\xymatrix{
\{A,B\}_\bullet  \ar[d]_-\epsilon \ar[rr] &&  \{A,B\} \ar[d]^{ \{A,\epsilon_B\}}  \\
\FF \ar[rr]^-{\epsilon_A}  && \{A,\FF\}.
}$$ 
In particular, $\{A,B\}_\bullet$ is a sub-coalgebra of $\{A,B\}$.

The coalgebra $\{A,B\}_\bullet$ is pointed by the map $\FF\to \{A,B\}_\bullet  $ lifting the  map $\{ \epsilon_A,e_B\}:\FF=\{\FF,\FF\} \to \{A,B\}$.
\end{prop}

\begin{proof}
If $(C,e)$ is a pointed  coalgebra, then there is a bijection
between the maps of pointed algebras $f:C\rhd_\bullet A\to B$ 
and the maps of algebras $g:C\rhd A\to B$ 
fitting in a commutative diagram
$$\xymatrix{
 A  \ar[d]_{\epsilon_A} \ar[rr]^-{e_C \rhd A} &&C\rhd A \ar[d]_g   \ar[rr]^-{\epsilon_C\rhd A} && A \ar[d]^{\epsilon_A}\\
\FF \ar[rr]^-{e_B} && B \ar[rr]^-{\epsilon_B}&& \FF
}$$ 
since the augmentation $\epsilon: C\rhd_\bullet A \to \FF$ is induced by the map $\epsilon_C \rhd \epsilon_A:C\rhd A\to \FF \rhd \FF=\FF$.
By adjointness, the diagram commutes iff the   following diagram commutes,
$$\xymatrix{
\FF \ar@{=}[d]\ar[rrr]^-{e_C}  &&&  C \ar[d]_h   \ar[rrr]^-{\epsilon_C}&&& \FF \ar[d]^{\lceil\epsilon_A\rceil}  \\
 \{\FF,B\}  \ar[rrr]^(0.4){\{\epsilon_A,B\}}    &&& \{A,B\} \ar[rrr]^(0.6){ \{A,\epsilon_B\}} &&& \{A,\FF\}.
}$$ 
where $h=\lambda^1(g)$. Let us put $K= \{A,B\} \times_{\{A,\FF\}} \FF$.
Composition with the first projection $p_1:K\to \{A,\FF\}$
induces bijection between the maps of coalgebras
$k:C\to  K$ and the maps of coalgebras $h:C\to \{A, B\}$ for which
the square on the right hand side commutes. 
There a unique map of coalgebras $e_K:\FF\to K$ such that
$p_1e_K= \{\epsilon_A,B\}e_B$ and $p_2e_K=1_\FF$.
We have $ke_K=e$ iff the square on the left hand side commutes.
Hence the composite $f\mapsto g\mapsto  h\mapsto k$  is a bijection between the 
maps of pointed coalgebras $g:C\rhd A\to B$  and the maps
of pointed coalgebras $k:C\to K$. The bijection is natural and it follows by Yoneda lemma
that $K=\{A,B\}_\bullet $.

Then $\{A,B\}_\bullet$ is a sub-coalgebra of $\{A,B\}$ by lemma \ref{directandinverseimage5}.
\end{proof}

Notice that the canonical map $\{A,B\}_\bullet \to\{A,B\}$ in the square of proposition \ref{pSHOMvsSHOM} is injective, 
since the square is a pullback and $\epsilon_A$ is injective (it is a split monomorphism).

\begin{rem}
The coalgebras $\{A,A\}_\bullet$ are examples of pointed bialgebras in the sense of section \ref{pointed bialgebras}.
Moreover, using the equivalence of proposition \ref{equivBilagpointed}, we have an isomorphism of pointed bialgebras
$$
\{A_{+},A_{+}\}_\bullet=\FF\times \{A,A\}_\circ.
$$
\end{rem}

\medskip

For $A$ and $B$ two pointed algebras and $C$ a pointed coalgebra, we shall say that a map $C\to [A,B]_\bullet$ is a {\em pointed comeasuring} if the corresponding map $C\wedge A\to B$ is a pointed measuring. Let us define the expanded analog.
We shall say that a map $C\to [A,B]$ is an {\em expanded pointed comeasuring} if the corresponding map $C\otimes A\to B$ is an expanded pointed measuring. We leave to the reader the definition of a couniversal (expanded or not) pointed measuring and the proof that the pointed comeasuring $\Psi_\bullet:\{A,B\}_\bullet\to [A,B]_\bullet$ and the expanded pointed comeasuring $\Psi:\{A,B\}_\bullet\to \{A,B\}\to [A,B]$ are couniversal. These can be used to give a direct proof of the following theorem.

To compare the notions of pointed, co-unital and non-counital comorphisms, there is an analog of lemma \ref{comparisonallmeta-morphismscoalg} that can be proven with a similar argument.
Recall that if $A$ and $B$ are pointed, we have a canonical decomposition $[A,B]=\FF\oplus [A_-,\FF]\oplus B_-\oplus [A_-,B_-]$.
Then, the equation $f(e_A)=e_B$ of proposition \ref{propmetanutou} implies that the couniversal (counital) comeasuring $\Psi:\{A,B\}\to [A,B]$ factors through $[A,B]_{\sf \bullet half} = \FF\oplus [A_-,\FF]\oplus [A_-,B_-]$.

\begin{lemma}\label{comparisonallmeta-morphismsalg}
We have a commutative diagram
$$\xymatrix{
\{A,B\}_\bullet \ar[rr]\ar[d]_{\Psi_\bullet} && \{A,B\} \ar[d]\ar[rrd]^{\Psi}\ar[rr]&& \{A,B\}_\circ\ar[d]^{\Psi_\circ}\\
[A,B]_\bullet \ar[rr] && [A,B]_{\sf \bullet half} \ar[rr] && [A,B].
}$$
where the vertical maps are the canonical inclusions.
The couniversal expanded pointed comeasuring is the total diagonal.
\end{lemma}

\bigskip

\begin{thm} \label{enrichmentpointedalgcoal}
The functor
$$
\{-,-\}_\bullet :\dgAlg_\bullet^{op}\times \dgAlg_\bullet\to \dgCoalg_\bullet
$$
defined above is an enrichment of the category $\dgAlg_\bullet$ over the monoidal category $(\dgCoalg_\bullet,\wedge, \FF_+)$.
The resulting enriched category is bicomplete. 
The tensor product of a pointed algebra $A$ by a pointed coalgebra $C$ is the pointed Sweedler product $C\rhd_\bullet A$ and the cotensor product 
of $A$ by $C$ is the pointed convolution algebra $[C,A]_\bullet$.
Hence there are natural isomorphisms of pointed coalgebras
$$\{C\rhd_{\bullet} A,B\}_{\bullet} \simeq \HOM_{\bullet}(C,\{A,B \}_{\bullet})\simeq \{A,[C,B]_\bullet \}_{\bullet}$$
for a pointed coalgebra $C$ and  pointed algebras $A$ and $B$.
\end{thm}

\begin{proof}
This follows from theorem \ref{enrichmentalgcoalpartial} and propositions \ref{equivnonunitalpointedalg} and \ref{equivnonunitalpointedcoalg}.
\end{proof}

If $A$ is a pointed algebra, $X$ is a vector space and $p: T^\vee_\circ([X,A_{-}])\to [X,A_{-}]$
is the canonical map, then the composite
of the map
$$
\xymatrix{
 T^\vee_\circ ([X,A_{-}])\otimes X \ar[rr]^-{p\otimes X} && [X,A_{-}]\otimes X \ar[r]^(0.6){ev} & A_{-}
}
$$
can be extended as an expanded pointed measuring $p':T^\vee_\bullet([X,A_{-}])\otimes T_\bullet (X)\to A$
by lemma \ref{extensionpointedmeasuring}.

\begin{prop} \label{exSweedlerhompointed}
If $A$ is a pointed algebra and $X$ is a vector space, 
then the pointed measuring $p'$
defined above is couniversal.
Hence we have 
$$
\{T_\bullet(X),A\}_\bullet= T^\vee_\bullet([X,A_{-}]).
$$
\end{prop}
\begin{proof}
By equivalence with proposition \ref{exSweedlerhomnu}.
\end{proof}

\begin{cor} \label{corexSweedlerhompointed}
If $A$ is a pointed algebra and $X$ is a vector space, 
we have a natural isomorphism
$$
R^c\{T_\bullet(X),A\}_\bullet= T^c_\bullet([X,A_{-}]).
$$
\end{cor}

\begin{proof}
Direct from $R^cT^\vee_\bullet=T^c_\bullet$. 
\end{proof}

\bigskip
The following result is a direct consequence of the enrichment structure of $\dgAlg_\bullet$.

\begin{prop} \label{padjunctiontypes}
Let $C$ be a pointed coalgebra and $A$ be a pointed algebra.
We have the following strong adjunctions
$$\xymatrix{
C\rhd_\bullet(-):\dgAlg_\bullet  \ar@<.6ex>[r]& \dgAlg_\bullet:[C,-]\ar@<.6ex>[l]
}$$
$$\xymatrix{
[-,A]_\bullet: (\dgCoalg_\bullet)^{op}  \ar@<.6ex>[r]& \dgAlg_\bullet:\{-,A \}_\bullet\ar@<.6ex>[l]
}$$
$$\xymatrix{
(-)\rhd_\bullet A: \dgCoalg_\bullet  \ar@<.6ex>[r]& \dgAlg_\bullet:\{A,- \}_\bullet\ar@<.6ex>[l].
}$$
\end{prop}

\begin{cor}
Let $A$ be a pointed algebra, we have the following strong adjunctions
$$\xymatrix{
[-,A]_\bullet: \dgCoalgnilbullet^{op}  \ar@<.6ex>[r]& \dgAlg_\bullet:R^c\{-,A \}_\bullet\ar@<.6ex>[l]
}$$
$$\xymatrix{
(-)\rhd_\bullet A: \dgCoalgnilbullet  \ar@<.6ex>[r]& \dgAlg_\bullet:R^c\{A,- \}_\bullet\ar@<.6ex>[l].
}$$
\end{cor}
\begin{proof}
Compose the previous result with the radical adjunction $\iota:\dgCoalgnilbullet\rightleftarrows \dgCoalg_\bullet:R^c$.
\end{proof}

\begin{cor}\label{exgenbarcobar}
Let $X$ be a (dg-)vector space, we have the following strong adjunction
$$\xymatrix{
T_\bullet(X\otimes (-)_-): \dgCoalg_\bullet  \ar@<.6ex>[r]& \dgAlg_\bullet:T^\vee_\bullet([X,(-)_-])\ar@<.6ex>[l].
}$$
$$\xymatrix{
T_\bullet(X\otimes (-)_-): \dgCoalgnilbullet  \ar@<.6ex>[r]& \dgAlg_\bullet:T^c_\bullet([X,(-)_-])\ar@<.6ex>[l].
}$$
\end{cor}
\begin{proof}
By propositions \ref{Sweedleroftensorpointed} and proposition \ref{exSweedlerhompointed}
and corollary \ref{corexSweedlerhompointed}, 
we have $T_\bullet(X\otimes (-)_-)=(-)\rhd_\bullet T_\bullet X$,
$T^\vee_\bullet([X,(-)_-]) =\{T_\bullet X,- \}_\bullet$
and $T^c_\bullet([X,(-)_-]) =R^c\{T_\bullet X,- \}_\bullet$.
\end{proof}

\subsubsection{Monoidal strengths and lax structures}

Let us view $\dgCoalg_\bullet\times \dgCoalg_\bullet$ and $\dgAlg_\bullet\times \dgAlg_\bullet$ enriched over $\dgCoalg_\bullet$ with respective hom
$$
\HOM_\bullet((C_1,C_2),(D_1,D_2)) := \HOM_\bullet(C_1,C_2) \wedge \HOM_\bullet(D_1,D_2),
$$
$$
\{(C_1,C_2),(D_1,D_2)\}_\bullet := \{C_1,C_2\}_\bullet\wedge \{D_1,D_2\}_\bullet.
$$

\begin{prop}
The functors 
$$\xymatrix{
\wedge:\dgCoalg_\bullet\times \dgCoalg_\bullet\ar[r]& \dgCoalg_\bullet,
}$$ 
$$\xymatrix{
\wedge:\dgAlg_\bullet\times \dgAlg_\bullet\ar[r]& \dgAlg_\bullet
}$$ 
are strong symmetric monoidal structures.
\end{prop}
\begin{proof}
Direct from proposition \ref{numonoidalstrength}.
\end{proof}

\begin{prop}\label{plaxmonoidalstrength}
\begin{enumerate}
\item The functors 
$$\xymatrix{
\wedge:\dgCoalg_\bullet\times \dgCoalg_\bullet\ar[r]& \dgCoalg_\bullet
}$$ 
$$\xymatrix{
\wedge:\dgAlg_\bullet\times \dgAlg_\bullet\ar[r]& \dgAlg_\bullet
}$$ 
are strong symmetric monoidal functors.

\item The functors
$$\xymatrix{
\HOM_\bullet:(\dgCoalg_\bullet)^{op}\times \dgCoalg_\bullet\ar[r]& \dgCoalg_\bullet
}$$ 
$$\xymatrix{
\{-,-\}_\bullet:(\dgAlg_\bullet)^{op}\times \dgAlg_\bullet\ar[r]& \dgCoalg_\bullet
}$$ 
$$\xymatrix{
[-,-]_\bullet:(\dgAlg_\bullet)^{op}\times \dgAlg_\bullet\ar[r]& \dgAlg_\bullet
}$$ 
are strong symmetric lax monoidal functors.

\item The functor
$$\xymatrix{
\rhd_\bullet:\dgAlg_\bullet\times \dgAlg_\bullet\ar[r]& \dgAlg_\bullet
}$$ 
is a strong symmetric colax monoidal functor.
\end{enumerate}
\end{prop}
\begin{proof}
By equivalence with proposition \ref{laxstructures}.
\end{proof}

In particular, these functors induces functors on the corresponding categories of (co)monoids and (co)commutative (co)monoids. We are going to state explicitely two of them and leave all the others to the reader.
Let us first describe some of the structure maps of these (co)lax structures.

\medskip
The lax structure of $\{-,-\}_\bullet$ can be constructed from that of $\{-,-\}_\circ$, it is given by $(\alpha,\alpha_0)$ where $\alpha_0:\FF_+\simeq \{\FF_+,\FF_+\}_\bullet$ and 
where
$$\xymatrix{
\alpha: \{A_1,B_1\}_\bullet\wedge \{A_2,B_2\}_\bullet \ar[r] & \{A_1\wedge A_2,B_1\wedge B_2\}_\bullet
}$$
is such that the corresponding map $\alpha':\{A_1,B_1\}_\bullet\otimes \{A_2,B_2\}_\bullet \to \{A_1\wedge A_2,B_1\wedge B_2\}_\bullet$
can also be decribed as the unique expanded map of pointed coalgebras making the following square commute
$$\xymatrix{
\{A_1,B_1\}_\bullet\otimes \{A_2,B_2\}_\bullet \ar[rr]^-{\alpha'} \ar[d]_{\Psi\otimes \Psi}&& \{A_1\otimes A_2,B_1\otimes B_2\}_\bullet\ar[d]^{\Psi}\\
[A_1,B_1]\otimes [A_2,B_2] \ar[rr]^-\theta && [A_1\otimes A_2,B_1\otimes B_2]
}$$
where $\theta$ is the lax structure of $[-,-]$ in $\dgVect$ and the $\Psi$s are the couniversal expanded comeasurings.

\medskip
The colax structure of $\rhd_\bullet$ is given by $(\alpha,\alpha_0)$ where $\alpha_0:\FF_+\simeq \FF_+\rhd_\bullet \FF_+$ and 
where $\alpha$ is the unique map of pointed algebras such that the following square commute
$$\xymatrix{
C_1\otimes C_2\otimes A_1\otimes A_2 \ar[rr]^-{\sigma_{23}}\ar[d]_{\Phi\otimes \Phi}&& (C_1\rhd A_1)\otimes (C_2\rhd A_2)\ar[d]^{\Phi}\\
(C_1\otimes C_2)\rhd_\bullet (A_1\otimes A_2) \ar[rr]^-{\alpha}&& (C_1\rhd_\bullet A_1)\otimes (C_2\rhd_\bullet A_2)
}$$
where the $\Phi$s are the universal expanded pointed measurings.
In terms of elements, $\alpha$ is the unique map such that $\alpha((c_1\otimes c_2)\rhd_\bullet(a_1\otimes a_2)) = (c_1\rhd_\bullet a_1)\otimes (c_2\rhd_\bullet a_2)$.

\bigskip

We introduce the following categories:
$\dgAlg_\bullet^{\sf com}$ is the category of pointed commutative algebras,
$\dgCoalg_\bullet^{\sf coc}$ be the category of pointed cocommutative coalgebras.
Recall that a bialgebra is said to be commutative or cocommutative if its underlying algebra or coalgebra is.
Let $\dgBialg^{\sf com}_\bullet$ be the category of commutative pointed bialgebras
and $\dgBialg^{\sf coc}_\bullet$ be the category of cocommutative pointed bialgebras.

For a monoidal category $(\bC,\otimes)$, let $\Mon_{\not e}(\bC)$ and $\coMon_{\!\not\,\epsilon}(\bC)$ be respectively the categories of non-unital monoids and non-counital comonoids in $\bC$. All of the previous categories of (co)algebras inherit the symmetric monoidal structure of $\dgVect_\bullet$ and, according to section \ref{nonbiunitalbialgebras} and to the Eckman-Hilton argument, we have the following canonical equivalences of categories:
\begin{center}
\begin{tabular}{l}
\rule[-2ex]{0pt}{4ex} $\dgAlg_\bullet^{\sf com}=\Mon_{\not e}(\dgAlg_\bullet,\wedge)$, \\
\rule[-2ex]{0pt}{4ex} $\dgCoalg_\bullet^{\sf coc}=\coMon_{\,\!\not\,\epsilon}(\dgCoalg_\bullet,\wedge)$, \\
\rule[-2ex]{0pt}{4ex} $\dgBialg =\Mon_{\not e}(\dgCoalg_\bullet,\wedge)=\coMon_{\!\not\,\epsilon}(\dgAlg_\bullet,\wedge)$, \\
\rule[-2ex]{0pt}{4ex} $\dgBialg^{\sf com}=\Mon_{\not e}(\dgBialg_\bullet,\wedge)=\coMon_{\!\not\,\epsilon}(\dgAlg^{\sf com}_\bullet,\wedge)$, \\
\rule[-2ex]{0pt}{4ex} $\dgBialg^{\sf coc}=\coMon_{\!\not\,\epsilon}(\dgBialg,\wedge)=\Mon_{\not e}(\dgCoalg^{\sf coc}_\bullet,\wedge)$.
\end{tabular}
\end{center}

\begin{prop}\label{corcolaxSproductpointed}
The Sweedler product $\rhd_\bullet$ induces functors
$$\xymatrix{
\rhd_\bullet : \dgCoalg^{\sf coc}_\bullet \times \dgBialg \ar[r] & \dgBialg
}$$
$$\xymatrix{
\rhd_\bullet : \dgCoalg^{\sf coc}_\bullet \times \dgBialg^{\sf coc} \ar[r] & \dgBialg^{\sf coc}.
}$$
\end{prop}
\begin{proof}
Direct from the symmetric colax monoidal structure of proposition \ref{plaxmonoidalstrength}.
The first result is obtain for non-counital comonoids, the second for cocommutative non-counital comonoids.
\end{proof}

If $C$ is a cocommutative pointed coalgebra and $H$ a bialgebra, let us describe the coproduct $\Delta$ of $(C\rhd_\bullet H)_-=C_-\rhd_\circ H_-$. 
It is given by composite 
$$\xymatrix{
\Delta : C_-\rhd_\circ H_- \ar[rr]^-{\Delta\rhd_\circ \Delta} && (C_-\otimes C_-)\rhd_\circ (H_-\otimes H_-)\ar[r]^\alpha &(C_-\rhd_\circ H_-)\otimes (C_-\rhd_\circ H_-)
}$$
where $\alpha $ is the colax structure.
Using the formula for $\alpha$, this gives on elements
$$
\Delta(c\rhd_\circ h)=\alpha\bigl( (c^{(1)}\otimes c^{(2)}) \rhd_\circ (h^{(1)} \otimes h^{(2)})  \bigr)=(c^{(1)}\rhd_\circ h^{(1)})\otimes (c^{(2)}\rhd_\circ h^{(2)}) (-1)^{|c^{(2)}||h^{(1)}|}.
$$

\begin{cor}\label{iterationcobarpointed}
If $H$ is a cocommutative bialgebra, it defines an endofunctor
$$\xymatrix{
-\rhd H : \dgCoalg^{\sf coc}_\bullet\ar[r] & \dgBialg^{\sf coc} \ar[r]^-{U_m} & \dgCoalg^{\sf coc}_\bullet.
}$$
where the functor $U_m:\dgBialg^{\sf coc}\to \dgCoalg^{\sf coc}$ is the functor forgetting the algebra structure.
\end{cor}

\medskip

\begin{prop}\label{corlaxShompointed}
The pointed Sweedler hom $\{-,-\}_\bullet$ induces functors
$$\xymatrix{
\{-,-\}_\bullet : \dgBialg^{op} \times \dgAlg^{\sf com}_\bullet \ar[r] & \dgBialg
}$$
$$\xymatrix{
\{-,-\}_\bullet : \left(\dgBialg^{\sf coc}\right)^{op} \times \dgAlg^{\sf com}_\bullet \ar[r] & \dgBialg^{\sf com}.
}$$
\end{prop}
\begin{proof}
Direct from the symmetric lax monoidal structure of proposition \ref{plaxmonoidalstrength}.
The first result is obtain for non-unital monoids, the second for commutative non-unital monoids.
\end{proof}

If $H$ is a cocommutative bialgebra and $A$ a commutative pointed algebra, $(\{H,A\}_\bullet)_-=\{H_-,A_-\}_\circ$ is a non-biunital bialgebra whose product a convolution defined as the composite 
$$\xymatrix{
\mu:\{H_-,A_-\}_\circ\otimes \{H_-,A_-\}_\circ \ar[r]^-{\alpha} & \{H_-\otimes H_-,A_-\otimes A_-\}_\circ \ar[rr]^-{\{\Delta_H, m_A\}} && \{H_-,A_-\}_\circ
}$$
where $\alpha $ is the lax structure.

\begin{cor}\label{iterationbarpointed}
If $H$ is a cocommutative bialgebra, it defines an endofunctor
$$\xymatrix{
\{H,-\}_\bullet : \dgAlg^{\sf com}_\bullet \ar[r] & \dgBialg^{\sf com}\ar[r]^-{U_\Delta}& \dgAlg^{\sf com}_\bullet
}$$
where the functor $U_\Delta:\dgBialg^{\sf com}\to \dgAlg^{\sf com}$ is the functor forgetting the coalgebra structure.
\end{cor}

\medskip
We shall apply both these results to iteration of the bar and cobar constructions in section \ref{iteratedbarcobar}.

\subsubsection{Reduction functors}\label{pointedreductionsection}

We interpret the couniversal pointed comorphism and pointed comeasuring (expanded or not) as the strength of enriched functors.
\medskip

Let $\dgVect_\bullet^\$$ and $\dgCoalg_\bullet^\$$ be the categories $\dgVect_\bullet$ and $\dgCoalg_\bullet$ viewed as enriched over themselves. Let also $\dgAlg_\bullet^\$$ be $\dgAlg_\bullet$ viewed as enriched over $\dgCoalg_\bullet$.
We can transfer the enrichment of $\dgCoalg_\bullet$ and $\dgAlg_\bullet$ along the lax monoidal functor $U_\bullet:\dgCoalg_\bullet\to \dgVect_\bullet$.
Let $\dgCoalg_\bullet^{U_\bullet\$}$ and $\dgAlg_\bullet^{U_\bullet\$}$ be the resulting categories enriched over $\dgVect$.

\begin{prop}\label{pointedreduction2}
The maps $\Psi_\bullet:\HOM_\bullet(C,D) \to [C,D]_\bullet$ and $\Psi_\bullet:\{A,B\}_\bullet \to [A,B]_\bullet$
are the strengths of strong monoidal functors over $\dgVect_\bullet$
$$\xymatrix{
U_\bullet:\dgCoalg_\bullet^{U\$} \ar[r] & \dgVect_\bullet^\$
}
\et
\xymatrix{
U_\bullet:\dgAlg_\bullet^{U\$}\ar[r] &\dgVect_\bullet^\$
}$$
\end{prop}
\begin{proof}
By equivalence with proposition \ref{strongforgetfulVectnu}.
\end{proof}

Let $\dgVect^\$$ be the category $\dgVect$ viewed as enriched over itself.
We can transfer the enrichment of $\dgCoalg_\bullet^\$$ and $\dgAlg_\bullet^\$$ along the lax monoidal functor $U:\dgCoalg_\bullet\to \dgVect$.
Let $\dgCoalg_\bullet^{U\$}$ and $\dgAlg_\bullet^{U\$}$ be the resulting categories enriched over $\dgVect$.

\begin{prop}\label{pointedreduction1}
The maps $\Psi:\HOM_\bullet(C,D) \to [C,D]$ and $\Psi:\{A,B\}_\bullet \to [A,B]$ are the strengths of enriched functors over $\dgVect$
$$\xymatrix{
U:\dgCoalg_\bullet^{U\$} \ar[r] & \dgVect^\$
}
\et
\xymatrix{
U:\dgAlg_\bullet^{U\$}\ar[r] &\dgVect^\$
}$$
\end{prop}
\begin{proof}
Essentially the same proof as for proposition \ref{strongforgetfulVectnu}.
\end{proof}

\subsubsection{Pointed meta-morphisms and pointed (co)derivations}

Recall the canonical inclusion $\HOM_\bullet(C,D)\subset \HOM(C,D)$ and $\{A,B\}_\bullet\subset \{A,B\}$ 
of propositions \ref{pHOMvsHOM} and \ref{pSHOMvsSHOM}.
We shall say that a meta-morphism $f:C\leadsto D$ is {\it pointed} if $f\in \HOM_\bullet(C,D)$
and that a meta-morphism $f:A\leadsto B$ is {\it pointed} if $f\in \{A,B\}_\bullet$.

In particular, we can apply on pointed meta-morphisms the calculus of (unpointed) meta-morphisms.
There is also a pointed calculus for pointed meta-morphisms, using all the pointed operations.
We leave all details to the reader. 

\medskip
We have the following concrete characterization of pointed meta-morphisms.

\begin{prop}
\begin{enumerate}
\item A meta-morphism $f:C\leadsto D$ is pointed iff $f(e_C)=\epsilon(f)e_D$.
\item A meta-morphism $f:A\leadsto B$ is pointed iff $\epsilon_B(f(a))=\epsilon(f)\epsilon_A(a)$ for every $a\in A$.
\end{enumerate}
\end{prop}

\begin{proof}
This is a reformulation of propositions \ref{pHOMvsHOM} and \ref{pSHOMvsSHOM} using the calculus of meta-morphisms.
\end{proof}

\medskip
We use the same musical notation as for the counital and non-counital meta-morphisms.
Recall that the notion of atom for pointed coalgebras coincides with that of atom for the underlying coalgebra.

\begin{prop}\label{musicallemmacoalgp}
Let $C$ and $D$ be two pointed coalgebras, then the maps $\sharp$ and $\flat=\Psi:\HOM_\bullet(C,D)\to [C,D]_\bullet$ induce:
\begin{enumerate}
\item inverse bijections of sets
$$
\dgCoalg_\bullet(C,D)\simeq At(\HOM_\bullet(C,D))
$$
where $At(\HOM_\bullet(C,D))$ is the set of atoms of $\HOM_\bullet(C,D)$;

\item inverse isomorphisms in $\dgVect$
$$
\Prim_f(\HOM_\bullet(D,C)) \simeq \Coder_\bullet(f);
$$

\item and, if $C=D$ inverse Lie algebra isomorphisms 
$$
\Prim(\END_\bullet(C))\simeq\Coder_\bullet(C)
$$
which preserve the square of odd elements.
\end{enumerate}
\end{prop}

\begin{proof}
Same as lemma \ref{musicallemmacoalg}, corollary \ref{corprimicoderiv} and theorem \ref{primiLiecoderiv}.
\end{proof}

\begin{prop}\label{musicallemmaalgp}
Let $A$ and $B$ be two pointed algebras, then the maps $\sharp$ and $\flat=\Psi:\{A,B\}_\bullet\to [A,B]_\bullet$ induce:
\begin{enumerate}
\item inverse bijections of sets
$$
\dgAlg_\bullet(A,B)\simeq At(\{A,B\}_\bullet)
$$
where $At(\{A,B\}_\bullet)$ is the set of atoms of $\{A,B\}_\bullet$.

\item inverse isomorphisms in $\dgVect$
$$
\Prim_f(\{A,B\}_\bullet) = \Der_\bullet(f);
$$

\item and, if $A=B$ inverse Lie algebra isomorphisms 
$$
\Prim(\{A,A\}_\bullet)\simeq\Der_\bullet(A)
$$
which preserve the square of odd elements.
\end{enumerate}
\end{prop}

\begin{proof}
Same as lemma \ref{musicallemmaalg}, corollary \ref{corprimideriv} and theorem \ref{primiLiederiv}.
\end{proof}

\bigskip
We leave the reader to translate the results of transport of non-(co)unital (co)derivations by the Sweedler operations to the pointed case.
Because of the functors $[-,-]_\bullet$ and $\wedge$, these results are more conveniently handled in the non-unital context.

\subsubsection{Strong (co)monadicity}

By equivalence with theorem \ref{strongmonadicitynu}, we have a pointed analog of theorems \ref{strongcomonadicitycoalg} and \ref{strongmonadicityalg}.

\begin{thm}\label{strongmonadicityp}
The adjunction $U_-\dashv T^\vee_\bullet$ enriches into a strong lax monoidal comonadic adjunction
$$\xymatrix{
U_-: \dgCoalg_\bullet^\$ \ar@<.6ex>[r]& \dgVect^{T^\vee_\bullet\$}:T^\vee_\bullet. \ar@<.6ex>[l]
}$$
The adjunction $T\dashv U$ enriches into a strong colax monoidal monadic adjunction
$$\xymatrix{
T_\bullet: \dgVect^{T^\vee_\bullet \$} \ar@<.6ex>[r]& \dgAlg_\bullet^\$:U_-. \ar@<.6ex>[l]
}$$
\end{thm}

In consequence, we can construct the hom coalgebras as equalizers in $\dgCoalg_\bullet$
$$\xymatrix{
\HOM_\bullet(C,D)\ar[r]^-{\Delta'} &T^\vee_\bullet([C_-,D_-]) \ar@<.6ex>[rr]^-{T^\vee_\bullet(\alpha)}\ar@<-.6ex>[rr]_-{T^\vee_\bullet(\beta)}&& T^\vee_\bullet([C_-,D_-\otimes D_-])
}$$
where for $f:C_-\to D_-$, we put $\alpha (f)=(f\otimes f)\Delta_{C_-}:C_-\to D_-\otimes D_-$ and $\beta (f)=\Delta_{D_-}f:C_-\to D_-\otimes D_-$,
and 
$$\xymatrix{
\{A,B\}_\bullet\ar[r]^-{m'} &T^\vee_\bullet([A_-,B_-]) \ar@<.6ex>[rr]^-{T^\vee_\bullet(\alpha)}\ar@<-.6ex>[rr]_-{T^\vee_\bullet(\beta)}&& T^\vee_\bullet([A_-\otimes A_-,B_-])
}$$
where for $f:A_-\to B_-$, we put $\alpha (f)=m_{B_-}(f\otimes f):A_-\otimes A_-\to B_-$ and $\beta (f)=fm_{A_-}:A_-\otimes A_-\to B_-$.

\subsection{Other contexts}

In this section, we mentionned briefly som other contexts where a Sweedler theory exists and may be meaningful.

\subsubsection{The Hopf context}\label{hopfsweedler}

Let $Q$ be a cocommutative Hopf algebra, we claim we have the following strengthening of theorems \ref{Qhomcoalg} and \ref{Qenrichmentalgcoal}..

\begin{thm}\label{QSweedlertheory}
\begin{enumerate}
\item The category $Q\dgCoalg$ is symmetric monoidal closed and the forgetful functor $U:Q\dgCoalg\to \dgCoalg$ is symmetric monoidal and preserves the internal hom. 
Moreover the adjunction $Q\otimes - \dashv U$ is enriched over $\dgCoalg$ and is strongly monadic.

\item The category $Q\dgAlg$ is enriched, bicomplete and monoidal over $Q\dgCoalg$.
The forgetful functor $U:Q\dgAlg\to \dgAlg$ is symmetric monoidal and preserves all Sweedler operations.
Moreover the adjunction $Q\rhd- \dashv U$ is enriched over $\dgCoalg$ and is strongly monadic.
\end{enumerate}
\end{thm}

This theorem says in particular that there is a canonical action of $Q$ on the Sweedler dg-constructions computed on $Q$-module (co)algebras. This action can be precised using the calculus of meta-morphisms.
Moreover we claim that {\sl all} the statements of the Sweedler theory of dg-(co)algebras stay true in this context.
All the functors (like the free algebra and cofree coalgebra functors) are compatible with the action of $Q$, 
in particular the distinguished isomorphisms
$$
T^\vee [C,X]\simeq \HOM(C,T^\vee (X))
\quad,\qquad
\{TX,A\}\simeq T^\vee([X,A])\et 
C\rhd T(X) \simeq T(C\otimes X)
$$
are $Q$-equivariant.

\bigskip
We have already given an application of this to pass from graded to differential graded Sweedler theory, we will mentionned another one after theorem \ref{commQSweedlertheory}.

\subsubsection{The commutative context}
We claim that all our result stay true if associative (co)algebras (unital or not, differential or not) are replaced by commutative (co)algebras.
We give details on some theorems.

\medskip
Let $\dgAlg^{\sf com}$ be the category of commutative (unital associative) algebras and $\dgCoalg^{\sf cocom}$ the category of cocommutative (counital coassociative) coalgebras.
For a commutative algebra $A$, it is easy to see the product map $A\otimes A\to A$ is a map of algebras.
The tensor product of two commutative algebras is again commutative and using the previous remark, $A\otimes B$ becomes the sum in the category $\dgAlg^{\sf com}$.
Dually, the coproduct map $C\to C\otimes C$ is a map of coalgebras, the tensor product of two cocommutative coalgebras is again cocommutative and $C\otimes D$ becomes the cartesian product in the category $\dgCoalg^{\sf cocom}$.

Recall that the canonical inclusion $\dgAlg^{\sf com}\to \dgAlg$ has a left adjoint called the {\em abelianization}, it sends an algebra $A$ to the commutative algebra $A_{ab}$ defined by the coequalizer
$$\xymatrix{
A\otimes A \ar@<.6ex>[rr]^-{m_A}\ar@<-.6ex>[rr]_-{m_A\sigma}&&A\ar[r]&A_{ab}
}$$
where $\sigma:A\otimes A\to A\otimes A$ is the symmetry of the tensor product. 

Similarly the canonical inclusion $\dgCoalg^{\sf cocom}\to \dgCoalg$ has a right adjoint called the {\em coabelianization}, it sends a coalgebra $C$ to the commutative coalgebra $C^{ab}$ defined by the equalizer
$$\xymatrix{
C^{ab}\ar[r]& C\ar@<.6ex>[rr]^-{\Delta_C}\ar@<-.6ex>[rr]_-{\sigma \Delta_C}&&C\otimes C
}$$
where $\sigma$ is again the symmetry of the tensor product.

\medskip

We can then recover some results of \cite{Barr}.
\begin{thm}
The cofree cocommutative coalgebra on a space $X$ exists and is isomorphic to $T^\vee(X)^{ab}$.
\end{thm}
\begin{proof}
For $C$ a cocommutative coalgebra, we have natural bijections between
\begin{center}
\begin{tabular}{lc}
\rule[-2ex]{0pt}{4ex} the linear map & $C\to X$\\
\rule[-2ex]{0pt}{4ex} the coalgebra maps & $C\to T^\vee(X)$\\
\rule[-2ex]{0pt}{4ex} the cocommutative coalgebra maps & $C\to (T^\vee(X))^{ab}$.
\end{tabular}
\end{center}
\end{proof}

\begin{thm}
\begin{enumerate}
\item The category $\dgCoalg^{\sf cocom}$ is cartesian closed and the functor $(-)^{ab}:\dgCoalg\to \dgCoalg^{\sf cocom}$ is monoidal and preserves the internal hom. In other terms, the internal hom of $\dgCoalg^{\sf cocom}$ is given by $\HOM(C,D)^{ab}$.

\item The category $\dgAlg^{\sf com}$ is enriched and bicomplete over $\dgCoalg^{\sf cocom}$.
Moreover the functor $(-)_{ab}:\dgAlg\to \dgAlg^{\sf com}$ preverses all Sweedler operations.
In other terms, 
the hom coalgebra is $\{A,B\}^{ab}$,
the tensor is $(C\rhd A)_{ab}$,
the cotensor is $[C,A]_{ab}=[C,A]$
and we have also $(A\otimes B)_{ab} = A_{ab} \otimes B_{ab}$.
\end{enumerate}
\end{thm}

\begin{proof}
The monoidal structure of the functors $(-)^{ab}$ and $(-)_{ab}$ is essentially the isomorphisms $(C\otimes D)^{ab} = C^{ab}\otimes D^{ab}$ and $(A\otimes B)_{ab} = A_{ab}\otimes B_{ab}$ that we leave to the reader to prove. 
Then, to prove the first assertion, it is sufficient to prove that $\HOM(C,D)^{ab}$ is an internal hom for $\dgCoalg^{\sf cocom}$.
Let $C$, $D$ and $E$ be three cocommutative coalgebras, the result is a consequence of the natural bijections between
\begin{center}
\begin{tabular}{lc}
\rule[-2ex]{0pt}{4ex} the cocommutative coalgebra maps & $C\otimes D\to E$\\
\rule[-2ex]{0pt}{4ex} the coalgebra maps & $C\to \HOM(D,E)$\\
\rule[-2ex]{0pt}{4ex} the cocommutative coalgebra maps & $C\to \HOM(D,E)^{ab}$.
\end{tabular}
\end{center}

Let $A$ and $B$ be two commutative algebras and $C$ be a cocommutative coalgebra, we have natural bijection between
\begin{center}
\begin{tabular}{lc}
\rule[-2ex]{0pt}{4ex} the coalgebra maps & $C\to \{A,B\}$\\
\rule[-2ex]{0pt}{4ex} the cocommutative coalgebra maps & $C\to \{A,B\}^{ab}$\\
\rule[-2ex]{0pt}{4ex} the algebra maps & $C\rhd A \to B$\\
\rule[-2ex]{0pt}{4ex} the commutative algebra maps & $(C\rhd A)_{ab} \to B$\\
\rule[-2ex]{0pt}{4ex} and the algebra maps & $A \to [C,B]$.
\end{tabular}
\end{center}
With the remark that $[C,B]$ is commutative when $B$ and $C$ are, this proves that $\dgAlg^{\sf cocom}$ is enriched tensored and cotensored over $\dgCoalg^{\sf cocom}$ and that $(-)_{ab}:\dgAlg\to \dgAlg^{\sf cocom}$ preserves all operations.
\end{proof}

If $Q$ is a cocommutative Hopf algebra, we claim that have a commutative analog of theorem \ref{QSweedlertheory}.

\begin{thm}\label{commQSweedlertheory}
\begin{enumerate}
\item The category $Q\dgCoalg^{\sf cocom}$ is symmetric monoidal closed and the forgetful functor $U:Q\dgCoalg^{\sf cocom}\to \dgCoalg^{\sf cocom}$ is monoidal and preserves the internal hom. 
Moreover the adjunction $Q\otimes- \dashv U$ is enriched over $\dgCoalg^{\sf cocom}$ and is strongly monadic.

\item The category $Q\dgAlg^{\sf com}$ is enriched, bicomplete and symmetric monoidal over $Q\dgAlg^{\sf cocom}$.
The forgetful functor $U:Q\dgAlg^{\sf com}\to \dgAlg^{\sf com}$ preverses all Sweedler operations.
Moreover the adjunction $Q\rhd- \dashv U$ is enriched over $\dgCoalg^{\sf cocom}$ and is strongly monadic.
\end{enumerate}
\end{thm}

\bigskip
Let us mentionned two applications of this theorem for $Q=\d_{-1}$ and $Q=\d_1$ (with the notation of example \ref{primitiveHopf}).
Recall that for $Q=\d_{-1}$, the category of $\d_1$-modules is the category of dg-graded vector spaces equipped with an extra differential commuting with the other one. A $Q$-module dg-(co)algebras is a dg-(co)algebras, with a distinguished (co)derivation of (homological) degree $-1$, of square zero and commutting with the structural (co)derivation of degree $-1$.
For $Q=\d_1$ of example \ref{primitiveHopf}, the category of $\d_1$-modules is the category of dg-graded vector spaces equipped with an extra differential of degree $1$ commuting with the other one. Such objects are called {\em mixed complexes} and are related to cyclic (co)homology \cite{Kassel}. A $Q$-module dg-(co)algebras are mixed dg-(co)algebras, \ie a dg-(co)algebras with a distinguished (co)derivation of (homological) degree $1$, of square zero and commutting with the structural (co)derivation of degree $-1$.

\begin{lemma}
For $n$ an odd integer and $A$ a commutative algebra, the commutative algebra $(\d_n\rhd A)_{ab}$ is $S_A(S^n\Omega_A^{\sf com})$ the symmetric $A$-algebra generated by the $n$-th suspension of the module $\Omega_A^{\sf com}$ of commutative differentials of $A$.
\end{lemma}
\begin{proof}
We give the sketch of a proof.
Recall that for a commutative ring $A$ a bimodule $M$ is said to be {\em symmetric} if the two left and right actions coincide through the isomorphism $A\otimes M\simeq M\otimes A$.
The full subcategory of symmetric bimodules is reflexive, the left adjoint send a $A$-bimodule $M$ to $M/[A,M]$ where $[A,M]$ is the sub-bimodule generated by the elements $ax-xa\ (-1)^{|a||x|}$ for $a\in A$ and $x\in M$.
Let $\Omega_A$ be the module of non-commutative differential of $A$, then $\Omega_A^{\sf com}=\Omega_A/[A,\Omega_A]$.

From the definition of Sweeder product, $(\d_n\rhd A)_{ab}$ is generated as an $A$-algebra $\delta\rhd x$ for any $x\in A$ and the relations
\begin{eqnarray*}
&\delta\rhd 1 = 0,\\
&\delta\rhd (xy)=(\delta \rhd x)(1\rhd y) + (1\rhd x)(\delta \rhd y)\ (-1)^{n|x|}\\
&x(\delta\rhd y)=(\delta \rhd y)x \ (-1)^{|x|(n+|y|)}.
\end{eqnarray*}
But this coincides with the presentation of $S_A(\Omega_A^{\sf com}$ given from the isomorphism $\Omega_A^{\sf com}=\Omega_A/[A,\Omega_A]$.
\end{proof}

In consequence, the underlying graded objects of the dg-algebra $\d_{-1}\rhd A$ and of the mixed dg-algebra $\d_1\rhd A$ are the two version of the de Rham algebra (where the differential are in degree $\pm 1$).

\medskip
The first adjunctions  
$$\xymatrix{
U:\d_{-1}\dgCoalg^{\sf cocom}\ar@<.6ex>[r]&\ar@<.6ex>[l] \dgCoalg^{\sf cocom}:\HOM(\d_{-1},-)
}$$$$\xymatrix{
\d_{-1}\rhd - :\dgAlg^{\sf cocom}\ar@<.6ex>[r]&\ar@<.6ex>[l] \d_{-1}\dgAlg^{\sf cocom}:U
}$$ 
can be called the {\em (classical) de Rham adjunctions}
And the second adjunctions  
$$\xymatrix{
U :\d_{1}\dgCoalg^{\sf cocom}\ar@<.6ex>[r]&\ar@<.6ex>[l] \dgCoalg^{\sf cocom}:\HOM(\d_{1},-)
}$$$$\xymatrix{
\d_{1}\rhd - :\dgAlg^{\sf cocom}\ar@<.6ex>[r]&\ar@<.6ex>[l] \d_{1}\dgAlg^{\sf cocom}:U
}$$ 
are the {\em mixed (or cyclic, or $S^1$-equivariant) de Rham adjunctions} (see \cite{TVderham} for the algebra side).

\subsubsection{The general context}

After the last sections, the reader should begin to suspect that Sweedler theory is a fairly general fact concerning algebras and coalgebras. 
In fact we claim that the following is true.

\begin{thm}
If $(\cV,\otimes,\un)$ is a locally presentable symmetric monoidal closed category then 
\begin{itemize}
\item the category $(\coMon(\cV),\otimes)$ is locally presentable symmetric monoidal closed and comonadic over $\cV$ 
\item and the category $(\Mon(\cV),\otimes)$ is locally presentable, enriched, bicomplete and symmetric monoidal over $\coMon(\cV)$, and monadic over $\cV$.
\end{itemize}
Moreover if we enriched $\bV$ over $\coMon(\bV)$ (as in section \ref{enrichmentVectoverCoalg}) all the structures are strong.
\end{thm}
The first part of this theorem is proven in \cite{Porst} as well as the presentability and symmetric monoidal structure of $\Mon(\bV)$.

\bigskip
This result applies in particular to $(\cV,\otimes,\un)=(\Set,\times,1)$ to give the elementary fact that $\Mon(\Set)$ is enriched over $\coMon(\Set)=\Set$.

\newpage
\section{Adjunctions between algebras and coalgebras}\label{adjunctions}\label{applications}

In this section, we use the Sweedler operations to construct adjunctions between the categories of algebras and coalgebras.
If $C$ is a coalgebra and $A$ an algebra, we have from corollary \ref{adjunctiontypes} three types of adjunctions
\begin{center}
\hfill
\begin{tabular}{c}
\rule[-2ex]{0pt}{5ex} $\xymatrix{C\rhd - :\dgAlg \ar@<.6ex>[r]&\ar@<.6ex>[l] \dgAlg :[C,-]}$\\
\rule[-2ex]{0pt}{5ex} $\xymatrix{[-,A] :\dgCoalg \ar@<.6ex>[r]&\ar@<.6ex>[l] \dgAlg^{op} :\{-,A\}}$\\
\rule[-2ex]{0pt}{5ex} $\xymatrix{-\rhd A :\dgCoalg \ar@<.6ex>[r]&\ar@<.6ex>[l] \dgAlg :\{A,-\}}$
\end{tabular}
\hfill
\begin{tabular}{c}
\rule[-2ex]{0pt}{5ex} (I)\\
\rule[-2ex]{0pt}{5ex} (II)\\
\rule[-2ex]{0pt}{5ex} (III)
\end{tabular}
\end{center}
This section details examples of these adjunction types. We are mainly going to work in the unital and unpointed case. Only for type III adjunctions where the most important example is the classical Bar-Cobar constructions we will need the pointed context.

\subsection{Type I - Examples}\label{examplesweedlerproduct}

\subsubsection{Products and coproducts}\label{productcoproduct}

Let $I$ be a set and $\FF I$ be the dg-vector space with zero differential generated by $I$.
Let $e_i$ be the canonical basis of $\FF I$ indexed by the elements of $I$.
As observed in example \ref{diagonalcoalg}, $\FF I$ is a coalgebra with coproduct defined by $\Delta(e_i)=e_i\otimes e_i$.

Let $B$ be an algebra, we saw in example \ref{exconvolproduct} that the convolution algebra $[\FF I,B]$ is the product $B^I$ of $I$ copies of $B$
Let $A$ be another algebra, we have bijection between
\begin{center}
\begin{tabular}{lc}
\rule[-2ex]{0pt}{4ex} the algebra maps & $A\to [\FF I,B] = B^I$\\
\rule[-2ex]{0pt}{4ex} and the algebra maps & $\FF I\rhd A\to B$.
\end{tabular}
\end{center}
We deduce that the algebra $\FF I\rhd A$ is isomorphic to the sum of $I$ copies of $A$.

\medskip
This conclusion can also be reached by an explicit computation of $\FF I\rhd A$: it is generated by symbols $i\rhd a$ for each $i\in I$ and $a\in A$ and by relations $i\rhd ab=(i\rhd a)(i\rhd b)$. This proves that $\FF I\rhd A$ is the free product of $I$ copies of $A$.
In particular any $i$ defines an embedding $A\to \FF I\rhd A$.

\medskip
Therefore, the adjunction
$$\xymatrix{
\FF I\rhd - :\dgAlg \ar@<.6ex>[r]&\ar@<.6ex>[l] \dgAlg :[\FF I,-]
}$$
is the $I$-indexed sum-product adjunction.

\bigskip
More generally, and conformally to the philosophy of strongly bicomplete categories,
if $C$ is a coalgebra $C\rhd A$ is the $C$-indexed sum of copies of $A$ and $[C,B]$ is the $C$-indexed product of copies of $B$.
In other words, the adjunction
$$\xymatrix{
C\rhd - :\dgAlg \ar@<.6ex>[r]&\ar@<.6ex>[l] \dgAlg :[C,-]
}$$
is the $C$-indexed sum-product adjunction.

\subsubsection{Weil restriction}\label{Weilrestriction}

Let $E$ be a graded finite algebra bounded above or below, then $E^\star$ is a graded finite coalgebra (bounded below or above) and we have $(E^\star)^\star=E$ (proposition \ref{equivgrfinalgcoalg}).
Let $B$ be another algebra, we have a canonical isomorphism of algebras $E\otimes B = [E^\star,B]$ where $[E^\star,B]$ is the convolution algebra (example \ref{convolbasechange}). 
Therefore the Sweedler product $E^\star\rhd -$ gives us a left adjoint to the base change functor $E\otimes -$ (proposition \ref{smallalgexp}). 
$$\xymatrix{
E^\star \rhd - :\dgAlg \ar@<.6ex>[r]&\ar@<.6ex>[l] \dgAlg :[E^\star,-] = E\otimes -
}$$
Because of the analogy with Weil restriction for commutative algebras, this adjunction can be call the non-commutative Weil restriction.

In other words, for $A$ and $B$ two algebras, we have bijection between
\begin{center}
\begin{tabular}{lc}
\rule[-2ex]{0pt}{4ex} the algebra maps & $A\to E\otimes B$\\
\rule[-2ex]{0pt}{4ex} and the algebra maps & $E^\star\rhd A\to B$.
\end{tabular}
\end{center}

\begin{rem}
The construction of $C\rhd A$ by generators and relations of theorem \ref{Sproduct} gives an explicit construction of Weil restriction.
Let us mentionned that this presentation also work for the construction of Weil restriction in commutative algebra, we need only to replace the free algebras by the free commutative algebra, but the generators and relations stay the same.
\end{rem}

The next two examples study in detail the case where $E$ is a matrix algebra and a dual number algebra.

\subsubsection{Matrix (co)algebras}\label{matrixexampleadj}

Let us denote by $Mat(n,A)$ the algebra of $n\times n$ matrices with coefficients in an algebra $A$.
We have $Mat(n,A)=Mat(n,\FF )\otimes A$.
By Weil restriction we have an adjunction 
$$\xymatrix{
Mat(n,\FF )^\star \rhd - :\dgAlg \ar@<.6ex>[r]&\ar@<.6ex>[l] \dgAlg :Mat(n,- ).
}$$
Where $Mat(n,\FF )^\star$ is the endomorphism coalgebra of example \ref{coendomorphisms}.

\medskip
For any algebra $A$, let us put $A^{[n]}:=Mat(n,\FF )^\star\rhd A $. Let us detail the structure of $A^{[n]}$.

A map of algebra $A\to Mat(n,B)$ is equivalent to the data of $n^2$ maps $(-)_{ij}:A\to B$ satisfying the relations
$$
(ab)_{ij} = \sum a_{ik}b_{kj}.
$$
Let us call these maps {\em representative functions of size $n$} of $A$ with values in $B$.
The algebra $A^{[n]}$ is then generated by symbols $a_{ij}$ corresponding to the values of representative functions of size $n$ on every element $a\in A$.

The unit of the adjunction is a map  $\eta:A\to Mat(n,A^{[n]})$ sends an element $a$ to the matrix $[a_{ij}]$ of its representative values.
It is universal in the following sense:
for any algebra $B$ and any morphism of algebras $f:A\to Mat(n,B)$,
there exists a unique morphism of algebras $g:A^{[n]}\to B$ such that $f=Mat(n,g)\circ \eta$,
$$\xymatrix{
A \ar[rr]^-{f}\ar[d]_{\eta} && Mat(n,B)  \\
Mat(n,A^{[n]}) \ar@{-->}[urr]_-{\qquad Mat(n,g)}  &&
}$$

\subsubsection{Differentials and de Rham algebra}\label{derhamexampleadj}

Let $\FF\delta_{+}=\FF \oplus \FF \delta=T^c_1(\delta)$ be the primitive coalgebra (see example \ref{primitivecoalgebra}) generated by an element $\delta$ of degree $n$. By definition, we have $\Delta(1)=1\otimes 1$ and $\Delta(\delta)= \delta\otimes 1+1\otimes \delta$.
The dual algebra $(\FF\delta_{+})^\star =\FF \oplus \FF\varepsilon= \FF\varepsilon_{+}$ is generated by an element $\varepsilon=\delta^\star$ of degree $-n$ and square 0. We have $\FF\delta_{+} = (\FF\varepsilon_{+})^\star$.

Weil restriction gives an adjunction
$$\xymatrix{
\FF\delta_{+} \rhd - :\dgAlg \ar@<.6ex>[r]&\ar@<.6ex>[l] \dgAlg : \FF\varepsilon_{+}\otimes -
}$$

For any algebra $B$, we have $\FF\varepsilon_{+}\otimes B=B\oplus  \varepsilon B= B[\varepsilon]$.
If $A$ is an algebra, then every map of algebras $A\to B\oplus  \varepsilon B$ is of the form $a\mapsto f(a)+\varepsilon D(a)$, where $f:A\to B$ is a map of algebras, and $D:A\rhup B$ is a graded morphism of degree $n$ satisfying
$$
D(ab)=D(a)f(b)+f(a)D(b) \ (-1)^{n|a|}
$$
\ie $D$ is a $f$-derivation of $A$ in $B$ of degree $n$.
(Beware that if we use the notation $B[\varepsilon]$, the map $A\to B[\varepsilon]$ corresponding to $D$ is 
$a\mapsto f(a)+ D(a)\varepsilon (-1)^{|\varepsilon|(|D|+|a|)}$.)

\medskip
Recall from section \ref{derivationalgebra} the bimodule $\Omega_A$ of differentials of the algebra $A$ (proposition \ref{univderiv}).
The bimodule $\Omega_A$ is the target of a universal derivation $d:A\to \Omega_A$.
Let $S^n\Omega$ bethe $n$-th suspension of $\Omega$. 
Recall from proposition \ref{derham} the differential algebra $T_A(S^n\Omega_A)$ defined as the tensor algebra over $A$ of bimodule $\Omega_A$, we have $T_A(S^n\Omega_A)=\FF\delta_+\rhd A$.
The universal measuring $\FF\delta_{+}\otimes A\to T_A(S^n\Omega_A)$
takes an element $a+\delta b\in A\oplus  \delta A= \FF\delta_{+}\otimes A$
to the element $a+s^ndb \in T_A(S^n\Omega_A)$.

\medskip
For $\delta$ a graded symbol of degree $n$ and $\varepsilon$ its dual, we have natural bijection between
\begin{center}
\begin{tabular}{lc}
\rule[-2ex]{0pt}{4ex} algebra maps & $A\to \FF\varepsilon_+\otimes B= B[\varepsilon]$,\\
\rule[-2ex]{0pt}{4ex} pairs (algebra maps, $A$-bimodule maps) & $(A\to B, \Omega_A\to \varepsilon B)$,\\
\rule[-2ex]{0pt}{4ex} pairs (algebra maps, $A$-bimodule maps) & $(A\to B, S^n\Omega_A\to B)$,\\
\rule[-2ex]{0pt}{4ex} and algebra maps & $T_A(S^n\Omega_A)=\FF\delta_+\rhd A\to B$.
\end{tabular}
\end{center}
And the previous adjunction $\FF\delta_+\rhd - \dashv \FF\varepsilon_+\otimes -$ rewrites in a more common form
$$\xymatrix{
T_{(-)}(S^n\Omega_{(-)}) :\dgAlg \ar@<.6ex>[r]&\ar@<.6ex>[l] \dgAlg : (-)[\varepsilon]
}$$
where $n=-|\varepsilon|$.

\bigskip

The differential algebra $T_A(S^n\Omega_A)$ is a non-commutative analog of the de Rham algebra. 
Let us now explain how it is canonically equipped with an analog of the de Rham differential when $\delta$ is of odd degree.

\medskip
If $\delta$ is of odd degree $n$ then the coalgebra $\d_n =\FF\delta_+$ is a cocommutative Hopf algebra by example \ref{primitiveHopf}.
An action of $\d_n $ on a dg-vector space $X$ is the data of a graded endomorphism $X\rhup X$ of degree $n$, of square zero and which commutes with the differential of $X$. Similarly a $\d_n $-module algebras, is a dg-algebra $A$ equipped with a graded derivation $A\rhup A$ of degree $n$, of square zero and which commutes with the differential of $A$.
Let $\d_n \dgAlg$ be the category of $\d_n $-module algebras, 
we have proven in theorem \ref{Qenrichmentalgcoal} that the functor $\d_n \rhd-$ is left adjoint to the forgetful functor
$U:\d_n \dgAlg \to \dgAlg$.
This says that not only $\d_n \rhd A$ is equipped with an action of $\d_n $, \ie a square zero derivation of degree $n$, but it is the free $\d_n $-module algebra generated by $A$.

Let us explain the action of $\delta\in \d_n $ on $\d_n \rhd A$.
According to proposition \ref{derham}, $\d_n \rhd A$ is generated as an algebra by elements of $A$ and by elements  $\delta\rhd a$ for $a\in A$.
It is sufficient to explain the action of $\delta\in \d_n $ on this generators, it is given by 
$$
\delta\cdot a = \delta \rhd a
\et
\delta\cdot(\delta \rhd a) = (\delta^2 \rhd a) = 0.
$$
The formulas are analogous to that of the classical de Rham differential.

\begin{rem}
The adjunction $\d_n\otimes-: \dgCoalg \rightleftarrows \dgCoalg:\HOM(\d_n,-)$ can be though as an analog for coalgebra of the 
Indeed by proposition \ref{coderham}, $\HOM(\d_n,C) = T^\vee_C(S^{-n}\Omega^C)$ which looks like the de Rham algebra $T_A(S^n\Omega_A)$.
\end{rem}

\bigskip
The following examples study the higher order generalisations of de Rham algebra (jet algebras).

\subsubsection{Jet algebras}\label{jetexample1adj}

Let  $T^c(x)=\FF[x]$ be the tensor coalgebra on the vector space $\FF x$ generated by an element $x$ of degree 0.
By definition, we have $\Delta(x^n)=\sum_{i=0}^n x^i\otimes x^{n-i}$. 
We want to understand the adjunction
$$\xymatrix{
T^c(x) \rhd - :\dgAlg \ar@<.6ex>[r]&\ar@<.6ex>[l] \dgAlg : [T^c(x),-].
}$$

\medskip
If $A$ is an algebra, then the map $i:[T^c(x),A]\to A[[t]]$ defined by putting $i(\phi)=\sum \phi(x^i)t^i$ is an isomorphism of algebras by example \ref{formalpower}.

Let us put $Jet(A)=T^c(x)\rhd A$. It is generated by symbols $x^n\rhd a$ and relations
$$
x^n\rhd (ab) = \sum_{i=0}^n (x^i\rhd a)(x^{n-1}\rhd b).
$$
which look like the Leibniz rule for divided higher derivations.
We shall call $x^n\rhd a$ is the {\it divided $n$-fold differential} of the element $a$ and denote it by $\dfrac{D^n}{n!}(a)$.
The algebra $Jet(A)$ is call the {\it jet algebra} of $A$.

In the previous example, the differential algebra of $A$ was generated by symbols $a$ and $da$ for any $a\in A$, where $da$ was a formal differential for $a$. Similarly, the jet algebra of $A$ is generated by symbols $x^n\rhd a=\,\dfrac{D^n}{n!}(a)$ which capture the higher differentials of $a$.
The unit of the adjunction $Jet(-)\dashv [T^c(x),-]$ is given by the map $\eta:A\to Jet(A)[[t]]$
which associates to an element $a\in A$ its {\em Taylor power series}
$$
\eta(a)=\sum_{n\geq 0} (x^n\rhd a) t^n = \sum_{n\geq 0} \dfrac{D^n}{n!}(a) t^n.
$$

\medskip
Finally the adjunction take the familiar form
$$\xymatrix{
Jet(-) :\dgAlg \ar@<.6ex>[r]&\ar@<.6ex>[l] \dgAlg : (-)[[t]].
}$$

\bigskip
As in the example \ref{derhamexampleadj}, $T^c(x)$ is a cocommutative Hopf algebra (example \ref{shufflehopf}).
We deduce from theorem \ref{Qenrichmentalgcoal} that the functor $Jet(A)=T^c(x)\rhd A$ is the free $T^c(x)$-module algebra generated by $A$. The action of $x\in T^c(x)$ generalizes the de Rham differential, but it is no longer an operator of square zero.

\subsubsection{Divided powers jet algebras}\label{jetexample2adj}

We consider now the coshuffle coalgebra $T^{csh}(x)=\FF[x]$ on one generator $x$ (of degree 0).
By definition $\Delta(x^n)=\sum_{i=0}^n {n\choose i} x^i\otimes x^{n-i}$. 
We want to understand the adjunction
$$\xymatrix{
T^{csh}(x) \rhd - :\dgAlg \ar@<.6ex>[r]&\ar@<.6ex>[l] \dgAlg : [T^{csh}(x),-].
}$$

\medskip
If $A$ is an algebra, recall from example \ref{Formal divided power series algebra} that $A\{\!\{t\}\!\}$ is the formal divided power series algebra on $A$.
Then the map $i:[T^{csh}(x),A]\to A\{\!\{t\}\!\}$ defined by putting $i(\phi)=\sum \phi(x^i)\dfrac{t^i}{i!}$ is an isomorphism of algebras by example \ref{dividedformalpower}.

Let us put $Jet_{div}(A)=T^{csh}(x)\rhd A$. It is generated by symbols $x^n\rhd a$ and relations
$$
x^n\rhd (ab) = \sum_{i=0}^n {n\choose i}(x^i\rhd a)(x^{n-1}\rhd b).
$$
which look like the Leibniz rule for higher derivations.
We shall call $x^n\rhd a$ is the {\it $n$-fold differential} of the element $a$ and denote it by $D^n(a)$.
The algebra $Jet_{div}(A)$ is call the {\it divided jet algebra} of $A$.

The unit of the adjunction $Jet_{div}(-)\dashv [T^{csh}(x),-]$ is given by the map $\eta:A\to J(A)[[t]]$
which associates to an element $a\in A$ its {\em Taylor divided power series}
$$
\eta(a)=\sum_{n\geq 0} (x^n\rhd a)\Dfrac{t^i}{i!} \ = \sum_{n\geq 0} D^n(a)\Dfrac{t^i}{i!}.
$$

\medskip
Finally the adjunction take the familiar form
$$\xymatrix{
Jet_{div}(-) :\dgAlg \ar@<.6ex>[r]&\ar@<.6ex>[l] \dgAlg : (-)\{\!\{t\}\!\}.
}$$

\subsection{Type II - Sweedler duality}\label{sweedlerduality}

The {\em Sweedler dual} of a dg-algebra $A$ is defined to be the dg-coalgebra  $A^\vee=\{A,\FF\}$.
($A^\vee$ was defined in \cite[Ch. VI]{Sw} with the notation $A°$.)
By the general theory it is part of a contravariant adjunction
$$\xymatrix{
[-,\FF] :\dgCoalg \ar@<.6ex>[r]&\ar@<.6ex>[l] \dgAlg^{op} :\{-,\FF\} = (-)^\vee.
}$$

\medskip
The main results of this section are the following
\begin{enumerate}
\item The Sweedler dual of an algebra $A$  is the coalgebra $A^\vee=\{A,\FF\}$; the functor $A\mapsto A^\vee$ is a strong right adjoint to the functor $C\mapsto C^\star$ (proposition \ref{Sweedlerdual}). We have $A^\vee =A^\star$ when the graded algebra $A$ is graded finite and bounded below or above (proposition \ref{Sweedlerdualgrinitefinite}).

\item For every graded finite, or stricly positive, or strictly negative dg-vector space we have $T^\vee(V)=T^c(V)$ (theorem \ref{tc=tv}).
\end{enumerate}

\bigskip

\begin{prop} \label{Sweedlerdual} 
The contravariant functors 
$(-)^\vee:\dgAlg \to \dgCoalg$ and  $(-)^\star:\dgCoalg \to \dgAlg$ are strong and mutually right adjoints. 
Moreover, if $A$ is an algebra and $C$ is a coalgebra,
then we have two natural isomorphisms
$$
 \HOM(C,A^\vee)\simeq(C\rhd A)^\vee   \simeq  \{A,C^\star\}.
$$
\end{prop}

\begin{proof} 
We have $\HOM(C,\{A,\FF  \}) \simeq \{C\rhd A,\FF \} \simeq  \{A,[C,\FF ] \}$ by theorem \ref{enrichmentalgcoal}.
\end{proof}

\medskip

\begin{prop} \label{Sweedlerduallaxandcolax}
The Sweedler dual of a bialgebra $A$ is a bialgebra $A^\vee$.
The bialgebra $A^\vee$ is commutative when $A$ is cocommutative.
\end{prop}

\begin{proof} The functor $F=\{-,\FF\}:\dgAlg^{op} \to \dgCoalg$ is symmetric lax by corollary \ref{Shomlax}, since
the algebra $\FF$ is commutative. It follows that the image by $F$ of a monoid
object $M\in \dgAlg^{op}$ is a monoid object $FM\in  \dgCoalg$. 
Moreover, $FM$ of a commutative if $M$ is commutative, since the
lax structure is symmetric.
\end{proof}

\medskip
The coalgebra $A^\vee$ is equipped with a couniversal measuring ${\bf ev}:A^\vee \otimes A\to \FF$.
Equivalently, it is equipped with a couniversal comeasuring $\Psi:A^\vee \to A^\star$.

\medskip
If $C$ is a coalgebra, then the evaluation $ev:C^\star \otimes C \to \FF$
is a right measuring, since the map $\lambda^2(ev):C^\star \to C^\star$
is an algebra map (it is the identity map).
Hence the canonical map $i:C\to [C^\star, \FF]$ is a comeasuring.

\begin{prop}\label{dualgradedfinitecoalg}
If $C$ is a graded finite dg-coalgebra, then the comeasuring $i:C\to [C^\star, \FF]$
is couniversal. Hence we have $C^{\star\vee} =C$.
\end{prop}

\begin{proof} 
 If $E$ is a coalgebra and $f:E\to [C^\star, \FF]$ is a comeasuring,
let us show that there is a unique map of coalgebras $g:E\to C$ such
that $ig=f$. The canonical map $i:C\to [C^\star, \FF]$ is invertible, since $C$
is graded finite. The uniqueness of $g$ follows.
Let us show that $g=i^{-1}f$ is a map of coalgebras. 
By lemma  \ref{dualofmapsofcoalgebras} it suffices to show
that $g^\star : C^\star \to E^\star$ is an algebra map.
By definition, for every $\phi\in C^\star$  and $x\in E$ we have $f(x)(\phi)=i(g(x))(\phi)=\phi(g(x))(-1)^{|x||\phi|}$.
The map $h={ev}(f\otimes C^\star):E\otimes C^\star \to \FF$
is a measuring, since $f$ is a comeasuring. Hence
 $k=\lambda^1(h):C^\star \to E^\star $ is an algebra map.
 For every $\phi\in C^\star $ and every $x\in E$ we have
 $$k(\phi)(x)=h(x\otimes \phi)(-1)^{|x||\phi|}=f(x)(\phi)(-1)^{|x||\phi|}=\phi(g(x)).$$
 Thus, $k(\phi)=g^\star(\phi)$ and it follows that $k=g^\star$.
 This shows that  $g^\star$ is an algebra map,
 since $k$ is an algebra map.
 \end{proof}

\begin{prop}\label{Sweedlerdualgrinitefinite}
If $A$ is a graded finite dg-algebra bounded above or below, then $A^\star=[A,\FF]$ has the structure
of a dg-coalgebra and the identity map $A^\star \to [A,\FF]$ is a couniversal comeasuring.
Hence we have $A^\vee=A^\star$.
\end{prop}

\begin{proof}
The dual $C= [A,\FF]$ is a graded finite coalgebra by proposition \ref{equivgrfinalgcoalg}.
The canonical map $i_A:A\to [[A,\FF], \FF]=[C,\FF]$ is an isomorphism
by proposition \ref{dualfinitealg}. We have $[i_A,\FF] i_C=1_C$ by an adjunction identity,
since the contravariant functor $[-,\FF]:\dgVect \to \dgVect$ is right adjoint to itself.
Hence the following square commutes
$$\xymatrix{
C\ar@{=}[r]  \ar[d]_{i_C}& A^\star  \ar@{=}[d] \\
[[C,\FF], \FF] \ar[r]^-{[i_A,\FF]}&  [A, \FF]
}$$
But the map $i_C$ is a couniversal comeasuring by proposition \ref{dualgradedfinitecoalg}.
It follows that the identity map $A^\star \to [A,\FF]$ is a couniversal comeasuring.
\end{proof}

\medskip

\begin{prop} \label{hommingtofinite}
If $B$ is a graded finite dg-algebra bounded below or above, then we have two natural isomorphisms
$$
\{A,B\}\simeq \HOM(B^\star ,A^\vee)\simeq (B^\star \rhd A)^\vee 
$$
for any algebra $A$.
Dually, if $D$ is a graded finite coalgebra bounded below or above, then we have two natural isomorphisms
$$
\HOM(C ,D) \simeq \{D^\star ,C^\star \}\simeq (C \rhd D^\star)^\vee
$$
for any coalgebra $C$.
\end{prop}
\begin{proof}
We have $B=B^{\star\star}$ by proposition \ref{dualfinitealg}, since $B$ is graded finite. Thus, 
$$
\{A,B\}\simeq \{A,B^{\star\star}\}\simeq \HOM(B^\star,A^\vee)\simeq (B^\star\rhd A)^\vee
$$
for any algebra $A$ by proposition \ref{Sweedlerdual}.
Dually, we have $D=D^{\star \vee}$ by proposition \ref{dualgradedfinitecoalg}, since the coalgebra $D$ is graded finite and bounded above or below. Thus, 
$$
\HOM(C,D)\simeq  \HOM(C,D^{\star\vee})\simeq  \{D^\star,C^\star\}\simeq (C\rhd D^\star)^\vee
$$
for any coalgebra $C$ by proposition \ref{Sweedlerdual}.
\end{proof}

\begin{cor} \label{Sweedlerdualfinite2} 
The contravariant adjunction
$$\xymatrix{
(-)^\vee:\dgAlg\ar@<.6ex>[r]& \dgCoalg:(-)^\star \ar@<.6ex>[l]
}$$
induces a strong contravariant equivalence between the following categories:
\begin{enumerate}
\item the category of graded finite algebras bounded below and the category of graded finite coalgebras bounded above
\item the category of graded finite algebras bounded above and the category of graded finite coalgebras bounded below
\item the category of finite algebras and the category of finite coalgebras.
\end{enumerate}
\end{cor}
\begin{proof}
The first two equivalences follows from propositions \ref{Sweedlerdualgrinitefinite} and \ref{dualgradedfinitecoalg}.
The last one is a consequence of the first two.
\end{proof}

The contravariant adjunction
$\xymatrix{
(-)^\vee:\dgAlg\ar@<.6ex>[r]& \dgCoalg:(-)^\star \ar@<.6ex>[l]
}$
induces a contravariant adjunction
$$\xymatrix{
F:(A\backslash \dgAlg)\ar@<.6ex>[r]& \dgCoalg/A^\vee:G \ar@<.6ex>[l]
}$$
for any algebra $A$. By construction, the functor $F$ takes a map $f:A\to B$ to the map $f^\vee:B^\vee \to A^\vee$ and the functor $G$ takes a map $g:C\to A^\vee$ to the map $g^\star u_A:A\to (A^\vee)^\star\to C^\star$, where $u_A:A\to  (A^\vee)^\star$ is the canonical map.
It then follows from \ref{Sweedlerdualfinite2} that the contravariant adjunction $(F,G)$ induces a contravariant equivalence of categories 
$$\xymatrix{
F':(A\backslash \fAlg)\ar@<.6ex>[r]& \fCoalg/A^\vee:G \ar@<.6ex>[l]
}$$
where $\fAlg$ is the category of finite algebras and $\fCoalg$ is the category of finite coalgebras (not to be confused with the categories $\gfAlg$ and $\gfCoalg$ of graded finite algebras and graded finite coalgebras). 

\medskip

\begin{prop} \label{Sweedlerdual3} 
The coalgebra $A^\vee$ is the colimit of the coalgebras $(A/J)^\star$, where $J$ runs in the poset of two-sided graded ideals of finite codimension in $A$.
It is also the colimit of coalgebras $(A/K)^\star$ where $K$ runs in the poset of two-sided graded ideal $K$ for which the quotient $A/K$ is graded finite and bounded below (resp. above).\end{prop}
\begin{proof} 
Recall that the codimension of a graded ideal $I$ of a graded algebra $A$ is the dimension of the total space of $A/I$.
Every coalgebra $C$ is the colimit of the canonical diagram of finite coalgebras $E\to C$, since the category $ \dgCoalg$ is finitary presentable and the $\omega$-compact coalgebras are the finite coalgebras by theorem \ref{finprescoalg}. 
More precisely, every coalgebra $C$ is the colimit of the forgetful functor $U:\fCoalg/C \to  \dgCoalg$.
In particular, the coalgebra $A^\vee$ is the colimit of the forgetful functor $U:  \fCoalg/A^\vee \to  \dgCoalg$.
The coalgebra $A^\vee$ is thus also the colimit of the  functor 
$$\xymatrix{
U'=UF': (A\backslash \fAlg)^{op} \ar[r]&  \dgCoalg.
}$$
since the functor $F'$ is a contravariant equivalence of categories. It is easy verify that the full subcategory $\cal F$ of 
$A\backslash \fAlg$ spanned by the surjections $p:A\to B$
is coreflexive in $A\backslash \fAlg$.
It follows that the coalgebra $A^\vee$ is the colimit of the functor $U'\mid {\cal F}^{op}$.
This proves the result, since the category ${\cal F}$ is equivalent to the poset of two-sided ideals of finite codimension in $A$.
The proof of the second statement is similar, but it uses lemma \ref{Sweedlerdualfinite2}.
\end{proof}

It follows from proposition \ref{Sweedlerdual3} that the algebra $A^{\vee \star}=(A^\vee)^\star$ is the limit of the algebras $A/J$,
where $J$ runs in the poset of two-sided ideals of finite codimension in $A$. 
It is thus the {\it completion} $A^\wedge$ of $A$ with respect to the linear topology defined by the two-sided ideals of  finite codimension in $A$. The canonical map $i_A:A\to (A^\vee)^\star$ is the canonical map $A\to A^\wedge$.
 
If $(A,\epsilon)$ is a pointed algebra, $A^\vee=\{A,\FF\}$ is pointed by the atom corresponding to $\epsilon:A\to \FF$.
The radical $R^c_\epsilon A^\vee$ is related to the completion $A^{\wedge,\epsilon}$ of $A$ along the maximal ideal $\ker \epsilon$ by the formula
$$
(R^cA^\vee)^\star = A^{\wedge,\epsilon}.
$$

\bigskip

If $X$ is a vector space and $p:T^\vee(X^\star)\to X^\star$ is the cofree map,
 then the composite
$$\xymatrix{
T^\vee(X^\star)\otimes X\ar[rr]^-{p\otimes X} && X^\star \otimes X \ar[r]^-{ev} & \FF
}$$
can be extended uniquely as a measuring $\mu:T^\vee(X^\star)\otimes T(X)\to \FF$
by proposition \ref{Sweedleroftensor}.

\begin{prop}\label{Sdualtensor} The measuring $\mu:T^\vee(X^\star)\otimes T(X)\to \FF$
defined above is couniversal. Hence we have $$T(X)^\vee=T^\vee(X^\star).$$
\end{prop}

\begin{proof} This is a special case of proposition \ref{exSweedlerhom}, since $T(X)^\vee= \{T(X),\FF\}$
\end{proof}

\begin{lemma}\label{dualtcvtv}
Let $X$ be a graded finite vector space and $i:X\to T(X)$ be the inclusion.
If $X$ is strictly positive (resp. strictly negative), then 
$T(X)^\star$ has the structure of a coalgebra and the map $i^\star:T(X)^\star\to X^\star$
is cofree. Hence we have $T^\vee(X^\star)=T(X)^\star$.
\end{lemma}

\begin{proof} The algebra $T(X)$ is graded
finite and bounded below, since  $X$ is graded finite and strictly positive.
Hence the dual $T(X)^\star$ has the structure of
a coalgebra and the map ${ev}:T(X)^\star \otimes T(X)\to \FF$
is a couniversal measuring by proposition \ref{dualfinitealg}.
But the measuring $\mu:T^\vee(X^\star)\otimes T(X)\to \FF$
of proposition \ref{Sdualtensor} is also couniversal.
Hence the map $\lambda^2(\mu):T^\vee(X^\star)\to T(X)^\star$
is an isomorphism of coalgebras, since we have $ev(\lambda^2(\mu)\otimes T(X)=\mu$.
If $p:T^\vee(X^\star) \to X^\star$ is the cofree map,
then the following square commutes by construction of $\mu$,
$$\xymatrix{
T^\vee(X^\star)\otimes X  \ar[d]_{T^\vee(X^\star) \otimes i} \ar[rr]^-{p\otimes X}&& X^\star \otimes X \ar[d]^{ev} \\
T^\vee(X^\star)\otimes T(X) \ar[rr]^\mu &&  \FF
}$$  
It follows that the following triangle commutes,
$$\xymatrix{
T^\vee(X^\star)  \ar[rrd]_{p} \ar[rr]^-{\lambda^2(\mu)}&& T(X)^\star  \ar[d]^{i^\star} \\
 && X^\star
}$$  
This shows that $i^\star$ is a cofree map, since $p$ is a cofree map.
\end{proof}

\begin{thm}\label{tc=tv}
If $X$ is a strictly positive (resp. negative) vector space, then 
$T^\vee(X)=T^c(X)$.
\end{thm}

\begin{proof} Let us first suppose that $X$ is finite dimensional and strictly positive.
Then we have $T^\vee(X^\star)=T(X)^\star$
and the map $i^\star:T(X)^\star\to X^\star$
is cofree by lemma \ref{dualtcvtv}, where $i:X\to T(X)$ is the inclusion.
 But the coalgebra $T(X)^\star$ is isomorphic
to the coalgebra $T^c(X^\star)$ by proposition \ref{dualtcvtvprop}.
Moreover, the map $i^\star$ is isomorphic to the projection
$T^c(X^\star)\to X^\star$.
This shows that the coalgebra $T^c(X^\star)$ is cofreely cogenerated
by the projection $T^c(X^\star)\to X^\star$.
It follows by duality that that the coalgebra $T^c(X)$ is cofreely cogenerated
by the projection $T^c(X)\to X$.
Hence the canonical map $T^c(X)\to T^\vee(X)$
is an isomorphism.
Let us now remove the finiteness condition on $X$.
The vector space $X$ is the directed
 union of its finite subspaces. 
 Every subspace of $X$ is strictly positive, since $X$ is strictly positive.
 The functor $T^\vee$ preserves directed colimits by theorem \ref{cofreecoalg}.
 Obviously, the functor  $T^c$
preserves directed colimits. 
Hence the canonical map $T^c(X)\to T^\vee(X)$ is the directed
colimit of the canonical maps $T^c(E)\to T^\vee(E)$, 
when $E$ runs in the poset a finite subspace of $X$.
The result follows, since the  canonical maps $T^c(E)\to T^\vee(E)$
is an isomorphism by the first part of the proof.
\end{proof}

\subsection{Type III - Bar-Cobar}

In this section, we apply the third kind of adjunction to an algebra $\mc$ that we called the Maurer-Cartan algebra. We prove that the corresponding adjunctions (in the unpointed and pointed contexts) are related the classical bar and cobar adjunction,
$$\xymatrix{
-\rhd \mc :\dgCoalg \ar@<.6ex>[r]&\ar@<.6ex>[l] \dgAlg :\{\mc,-\}
}$$
$$\xymatrix{
-\rhd_\bullet \mc :\dgCoalg_\bullet \ar@<.6ex>[r]&\ar@<.6ex>[l] \dgAlg_\bullet :\{\mc,-\}_\bullet.
}$$

\subsubsection{The Maurer-Cartan algebra}

\begin{defi}\label{mcelt}
If $A$ is a dg-algebra (unital or not), we shall say that an element $a\in A$ de degree $-1$ is a {\it Maurer-Cartan element} if it satisfies the {\it Maurer-Cartan equation}:
$$
da+aa =0.
$$
We shall denote $MC(A)$ the set of Maurer-Cartan elements of $A$. It is never empty as 0 is always a Maurer-Cartan element.
The image of a Maurer-Cartan element by a morphism of dg-algebras (unital or not) is again a Maurer-Cartan element and the sets $MC(A)$ define functors
$$\xymatrix{
MC: \dgAlg \ar[r] & \Set.
}\et
\xymatrix{
MC_\circ: \dgAlg_\circ \ar[r] & \Set.
}$$
We shall prove that these functors are representable.
\end{defi}

\medskip

Let $u$ be a graded variable of degree $-1$ and let $T(u)=\bigoplus_n \FF u^n$ be the tensor algebra on $u$.
Recall that $T(u)$ is pointed by the map $T(u)\to \FF$ sending $u$ to zero.
By proposition \ref{extensionderivationpointed}, the graded map of degree -1
\begin{eqnarray*}
\FF u &\underset{-1}{\lrhup} & T_\circ(u)\\
u &\mto& -u^2
\end{eqnarray*}
extend to a unique pointed derivation $d_\mc$ of degree $-1$ of $T(u)$.
By Leibniz' rule, 
\begin{eqnarray*}
d_\mc(u^n) &=& \sum_{i+j+1=n} u^i(du)u^j \ (-1)^{i|u|}\\
	&=& \sum_{i+j+1=n} u^i(-u^2)u^j \ (-1)^{i}\\
	&=& u^{n+1} \sum_{1\leq i \leq n} (-1)^{i+1}\\
	&=& 
	\left\{\begin{array}{ll}
	\rule[-2ex]{0pt}{4ex} 0 & \textrm{if $n$ is even}\\
	-u^{n+1} & \textrm{if $n$ is odd}.
	\end{array}\right.
\end{eqnarray*}
It is then clear that $d_\mc$ is of square zero, hence $(T(u),d_\mc)$ is a dg-algebra. 

\begin{defi}
We shall denote $\mc$ the dg-algebra $(T(u),d_\mc)$ and call it the {\em Maurer-Cartan algebra}.
By definition of $d_\mc$, $d_\mc u+u^2=0$, \ie the element $u$ is a Maurer-Cartan element of $\mc$. 

$\mc$ is pointed by the map $\mc \to \FF$ sending $u$ to 0.
We shall call the algebra $\mc_-=(T_\circ(u),d_\mc)$ and call it the {\em non-unital Maurer-Cartan algebra}.
\end{defi}

It can be helpful to picture $\mc$ as the complex
$$\xymatrix{
\FF\ar[r]^-0& \FF u\ar[r]^-{-u}&\FF u^2\ar[r]^-0&\FF u^3\ar[r]^-{-u}& \FF u^4\ar[r]^-0& \dots
}$$
In particular, it is clear that the homology of $\mc$ is only $\FF$ in degree 0. 
In fact $\FF$ is a retract of $\mc$ since there exists an augmentation $\epsilon:T(u)\to \FF$ sending $u$ to zero.

\begin{rem}\label{remudual}
Let $u^\star$ be the dual variable of $u$, we can think of $u^\star$ as an element of the dual graded space $\mc^\star$, it is the linear map whose value on $u$ is 1 and values on the other powers of $u$ is 0.
From Koszul's sign rule, $(u^\star)^{\otimes n}$ is the dual variable of $u^{\otimes n}$ up to a sign:
$$
(u^\star)^{\otimes n}(u^{\otimes n}) = (-1)^{\sum_{i=0}^{n-1} i} = (-1)^{\binom{n-1}{2}}.
$$
Viewed as an element of $\mc^\star$, $(u^\star)^{\otimes n}$ is the linear map whose value on $u^n$ is $(-1)^{\binom{n-1}{2}}$ and values on the other powers of $u$ is 0.
In particular, $(u^\star \otimes u^\star)(u) = 0$ and $(u^\star \otimes u^\star)(u^2) = -1$.
\end{rem}

\medskip

We shall say that a dg-algebra $A$ (unital or not) is {\em generated by a universal Maurer-Cartan element $u\in A$}, if it represents the functor $MC$. That is if, for any dg-algebra $B$ and any Maurer-Cartan element $b\in B$, there exists a unique map of dg-algebras $A\to B$ sending $u$ to $b$.

\begin{prop}\label{mcalg}
The Maurer-Cartan algebra $\mc$ is generated by a universal Maurer-Cartan element $u$. 
In other terms, we have a natural isomorphism
$$
MC(A) = \dgAlg(\mc,A).
$$
Similarly the non-unital Maurer-Cartan algebra $\mc_-$ represents the functor $MC_\circ$.
\end{prop}
\begin{proof}
Let $(A,d)$ be a dg-algebra and $a\in MC(A)$. 
There exists a unique map of graded algebras $T(u)\to A$ sending $u$ to $a$ since $|u|=|a|$. 
Let us show that this map commute with the differential: 
we have to prove that $d(a^n) = 0$ if $n$ is even and $d(a^n) = -a^{n+1}$ if $n$ is odd, but the computation is the same as for $u$.
\end{proof}

If $\alpha\in A$ is a Maurer-Cartan element, we shall note $\lceil \alpha \rceil:\mc\to A$ the corresponding algebra map and call it the {\em classifying map of $\alpha$}.

\medskip
The following result shows that $MC(\mc)$ has two elements.

\begin{lemma}
The algebra $\mc$ has only two Maurer-Cartan elements: $0$ and $u$.
\end{lemma}
\begin{proof}
The elements of degree $-1$ are all of the form $\lambda u$ for some $\lambda\in \FF$.
The Maurer-Cartan equation is $\lambda du +\lambda^2u^2 = (\lambda^2-\lambda )u^2= 0$.
Solutions are given by $\lambda=0$ or $\lambda=1$.
\end{proof}

\bigskip

We study now the Hopf structure of $\mc$.

\begin{lemma}\label{sommeEMC}
Let $A$ and $B$ be two dg-algebras. 
Si $a\in A$ et $b\in B$ are Maurer-Cartan elements, then the sum $a\otimes 1+1\otimes b$ is a Maurer-Cartan element of $A\otimes B$.
\end{lemma}
\begin{proof} 
Let us put $c=a\otimes 1+1\otimes b$. 
Then we have
\begin{eqnarray*}
cc &=&(a\otimes 1+1\otimes b)(a\otimes 1+1\otimes b)\\
	&=& aa\otimes 1+a\otimes b+(-1)^{|a||b|} a\otimes b+1\otimes bb\\
	&=& aa\otimes 1+1\otimes bb
\end{eqnarray*}
since $|a|=|b|=-1$.
Thus,
\begin{eqnarray*}
dc+cc &=& (da)\otimes 1+1\otimes db+aa\otimes 1+1\otimes bb \\
	&=& (da+aa)\otimes 1+1\otimes (db+bb)\\
	&=& 0.
\end{eqnarray*}
 \end{proof}

\medskip

If $A$ is an algebra and $A^o$ its opposite algebra, recall that, if $a\in A$ is an odd element $a^oa^o = (aa)^o \ (-1)^{|a||a|} = -(aa)^o$.

\begin{lemma} \label{mcop} 
If $a\in A$ is a Maurer-Cartan element in a dg-algebra $A$, then $-a^o$ is a Maurer-Cartan element in the opposite algebra $A^o$. 
\end{lemma} 
\begin{proof}
\begin{eqnarray*}
d(-a^o)+(-a^o)(-a^o) &=&- d(a)^o+a^oa^o\\
	&=&-d(a)-(aa)^o\\
	&=&0.
\end{eqnarray*}
\end{proof}

If $u\in \mc$ is the universal Maurer-Cartan element, then $u\otimes 1+1\otimes u $ is a Maurer-Cartan element of the algebra $\mc\otimes \mc$ by lemma \ref{sommeEMC}; hence 
there exists a unique morphism of dg-algebras $\Delta:\mc\to \mc\otimes \mc$ such that $\Delta(u)=u\otimes 1+1\otimes u$.
The element $0\in \mathbb{F}$ is Maurer-Cartan; hence there exists a unique morphism of dg-algebras $\epsilon:\mc\to \mathbb{F}$ such that $\epsilon(u)=0$.
The element $-u^o\in \mc^o$ is Maurer-Cartan by lemma \ref{mcop}; hence there exists a unique anti-homomorphism of dg-algebras $S:\mc\to \mc$ such that $S(u)=-u$.

\begin{prop}\label{mchopf}
The dg-algebra $\mc$ has the structure of a cocommutative conilpotent Hopf algebra with the coproduct $\Delta:\mc\to \mc\otimes \mc$, the counit $\epsilon:\mc\to \mathbb{F}$ and the antipode $S:\mc\to \mc$ defined above.
\end{prop}
\begin{proof}
But for the differential, the Hopf structure is that of the coshuffle Hopf algebra $T^{csh}(u)$ on an odd variable.
In particular, $\Delta (u^n) = \sum_{k=0}^n \oddbinom{n}{k}u^k\otimes u^{n-k}$ where the $\oddbinom{n}{k}$ are the odd binomial coefficients of example \ref{coshufflecoalgebra} and the coproduct is conilpotent.
We are left to check that $d_\mc$ is a coderivation. By construction of the coshuffle coproduct, it is sufficient to check that the map $\FF u\to \FF[u]\otimes \FF[u]:u\mapsto u\otimes 1+1\otimes u$ is compatible with the differential but this is obvious.
\end{proof}

\bigskip

We finish by a characterisation of modules over $\mc$.
Recall our convention (section \ref{dgvectorspaces}) that if $X$ is a dg-vector space, $|X|$ is the underlying graded vector space.
In particular, we can write $X=(|X|,d_X)$.

\begin{lemma} \label{mcmod}
Let $A$ be a dg-algebra and $(M,d)$ be a left $A$-module.
If $a\in A$ is a Maurer-Cartan element, then the map $d_a:M\to M$ defined by putting 
$$
d_a(x)=dx+ax
$$
is a differential, \ie $(d_a)^2=0$.
\end{lemma}
\begin{proof} 
For every $x\in M$ we have
\begin{eqnarray*}
d_a d_a(x)  &  = &  d(d(x)+ax)+ a(d(x)+ax)\\
	&=&dd(x)+   d (ax)+ ad (x)+aax \\
	&=&   d(a)x-  ad(x)+  ad (x) +aax \\
	&= &  (d (a)+aa)x=0.
\end{eqnarray*}
\end{proof} 

\begin{lemma}\label{mcmod2} 
If $X =(|X|,d)$ is a dg-vector space, then the map of graded vector spaces
\begin{eqnarray*}
[X,X] &\tto & [X,X]\\
\alpha &\mto& d_\alpha
\end{eqnarray*}
induces a bijection between the Maurer-Cartan elements of the dg-algebra $[X,X]$ and the differentials $d':|X|\to |X|$ on the underlying graded vector space of $X$.
\end{lemma} 
\begin{proof}
From lemma \ref{mcmod}, if $\alpha \in [X,X]$ is a Maurer-Cartan element, then the map $d_\alpha :|X|\to |X|$ is a differential. 
Conversely, if $ d_\alpha^2=0$, let us show that $\alpha$ is a Maurer-Cartan element of the algebra $ [X,X]$.
But we have 
$$
d_\alpha ^2=(d +\alpha )(d +\alpha )= \alpha d+d\alpha +\alpha \circ \alpha =d(\alpha )+\alpha \circ \alpha .
$$
Hence the condition $d_\alpha^2=0$ is equivalent to the Maurer-Cartan equation on $\alpha$.
\end{proof}

We shall say that a graded vector space $X$ equipped with two differentials $d',d'':X\to X$ (with no commutation conditions) is a {\it dicomplex}.
A map of dicomplexes $f:X\to Y$ is a map respecting the two differentials.
We shall denote the category of dicomplexes by $\ddgVect$.

\begin{prop}
The category $\Mod(\mc)$ of $\mc$-modules in $\dgVect$ is equivalent to the category $\ddgVect$ of dicomplexes.
\end{prop}
\begin{proof} 
For $(X,d)$ a dg-vector space, a dg-algebra map $\mc\to [X,X]$ indeuces a second differential by lemma \ref{mcmod2}.
Conversely, if $(X,d',d'')$ is a dicomplex, the complex $(X,d')$ has the structure of a left $\mc$-module by lemma \ref{mcmod2} if we put $u\cdot x =d''x-d'x$ for every $x\in X$.
\end{proof}

\subsubsection{Representation of twisting cochains}\label{twistingcochainsection}

\begin{defi} \label{cochainenetordantedef} 
Let $C$ be a dg-coalgebra and $A$ be a dg-algebra.
We shall say that a differential graded morphism $C\rhup_{-1} A$ is a {\it twisting cochain} if it is a Maurer-Cartan element of the the convolution algebra $[C,A]$.
More explicitly, a linear map $\alpha:C\to A$ of degree $-1$ is a twisting cochain if it satisfies the Maurer-Cartan equation:
$$
d_A\alpha+\alpha d_C+\alpha \star \alpha=0.
$$
where $\star$ is the convolution product.

If $C$ and $A$ are non-(co)unital, we shall define a {\em non-unital twisting cochain} $\alpha:C\rhup_{-1} A$ as a Maurer-Cartan element of the the convolution algebra $[C,A]$.

If the coalgebra $(C,e_C)$ and the algebra $(A,\epsilon_A)$ are pointed, we shall say that a twisting cochain $\alpha:C\rhup_{-1} A$ is {\it pointed}, or {\it admissible}, if $\alpha(e_C)=0=\epsilon_A \alpha$. Equivalently a pointed twisting cochain is a Maurer-Cartan element of the pointed convolution algebra $[C,A]_\bullet$.
\end{defi}

Recall that $\mc$ has a natural augmentation.

\begin{lemma}\label{lemmaTCmes}\label{mcalgpointed}
\begin{enumerate}
\item Twisting cochains $\alpha:C\rhup_{-1} A$ are in bijection with 
	\begin{itemize}
	\item algebras maps $\lceil \alpha\rceil:\mc\to [C,A]$
	\item and measuring $\alpha_\mc:C\otimes \mc\to A$.
	\end{itemize}

\item Pointed twisting cochains $\alpha:C\rhup_{-1} A$ are in bijection with 
	\begin{itemize}
	\item pointed algebras maps $\lceil \alpha\rceil:\mc\to [C,A]_\bullet$,
	\item expanded pointed measuring $\alpha_\mc:C\otimes \mc\to A$,
	\item non-unital twisting cochains  $\alpha_-:C_-\rhup A_-$,
	\item non-unital algebras maps $\lceil \alpha_-\rceil:\mc_-\to [C_-,A_-]$
	\item and non-unital measuring $\alpha_{\mc_-}:C_-\otimes \mc_-\to A_-$.
	\end{itemize}
\end{enumerate}
\end{lemma}
\begin{proof}
\begin{enumerate}

\item Clear by definition of twisting cochains and of measurings.
The bijection $\alpha \leftrightarrow \lceil \alpha\rceil$ is given by $\lceil \alpha\rceil(u)=\alpha$
and the bijection $\lceil \alpha\rceil\leftrightarrow \alpha_\mc$ is given by $\alpha_\mc(c,u) = \alpha(c)\ (-1)^{|c|}$.

\item The first bijections $\alpha \leftrightarrow \lceil \alpha\rceil \leftrightarrow \alpha_\mc$ are by definition of pointed twisting cochains and expanded pointed measurings. They are defined as in 1.
The bijection $\lceil \alpha\rceil \leftrightarrow \lceil \alpha_-\rceil$ is given by the equivalence between pointed and non-unital algebras.
The bijections $\alpha \leftrightarrow \lceil \alpha\rceil \leftrightarrow \alpha_\mc$ are by definition of non-unital twisting cochains and measurings.
\end{enumerate}
\end{proof}

We shall denote the set of twisting cochains $C\rhup A$ by $Tw(C,A)$. These sets define a functor of two variables
$$\xymatrix{
Tw:\dgCoalg^{op}\times \dgAlg\ar[r]& \Set.
}$$
If $A$ and $C$ are pointed, we shall denote the set of pointed twisting cochains $C\rhup A$ by $Tw_\bullet(C,A)$. 
This defines a functor of two variables
$$\xymatrix{
Tw_\bullet:(\dgCoalg_\bullet)^{op}\times \dgAlg_\bullet \ar[r]&\Set.
}$$
Restricted to conilpotent pointed dg-coalgebras, this defines also a functor
$$\xymatrix{
Tw_\bullet^c:(\dgCoalgnilbullet)^{op}\times \dgAlg_\bullet \ar[r]&\Set.
}$$
We have also the non-(co)unital analog functors $Tw_\circ(C,A)$ and $Tw_\circ^c(C,A)$.

\medskip

We now apply Sweedler theory to have an elementary proof of the representability of the functors $Tw$. 
The definition of a binary relator representable by an adjunction is recalled in appendix \ref{2rel}.

\begin{thm}\label{represtw}
\begin{enumerate}
\item For $C$ a dg-coalgebra and $A$ a dg-algebra, there exist natural isomorphisms 
$$
\dgAlg(C\rhd \mc, A) = Tw(C,A) = \dgCoalg(C,\{\mc,A\})
$$
In other words, the binary relator $Tw$ is representable by the adjunction $(-)\rhd\mc\dashv \{\mc,-\}$.
\item For $C$ a pointed dg-coalgebra and $A$ a pointed dg-algebra, there exist natural isomorphisms 
$$
\dgAlg_\bullet(C\rhd_\bullet \mc, A) = Tw_\bullet(C,A) = \dgCoalg_\bullet(C,\{\mc,A\}_\bullet)
$$
In other words, the binary relator $Tw_\bullet$ is representable by the adjunction $(-)\rhd_\bullet\mc\dashv \{\mc,-\}_\bullet$.
\item For $C$ a conilpotent pointed dg-coalgebra and $A$ a pointed dg-algebra, there exist natural isomorphisms 
$$
\dgAlg_\bullet(C\rhd_\bullet \mc, A) = Tw_\bullet^c(C,A) = \dgCoalg_\bullet(C,\{\mc,A\}_\bullet^c)
$$
where $\{\mc,A\}_\bullet^c:=R^c\{\mc,A\}_\bullet$.
In other words, the binary relator $Tw_\bullet^c$ is representable by the adjunction $(-)\rhd_\bullet\mc\dashv \{\mc,-\}_\bullet^c$.
\end{enumerate}
\end{thm}
\begin{proof}
\begin{enumerate}
\item By proposition \ref{mcalg}, we have
$$
Tw(C,A) = \dgAlg(\mc,[C,A]),
$$
the other isomorphisms follows from the properties of Sweedler operations.

\item By lemma \ref{mcalgpointed}, we have
$$
Tw_\bullet(C,A) = \dgAlg_\bullet(\mc,[C,A]),
$$
the other isomorphisms follows from the properties of pointed Sweedler operations.

\item This follows from the previous result and from $\dgCoalg_\bullet(C,\{\mc,A\}_\bullet)=\dgCoalgnilbullet(C,R^c\{\mc,A\}_\bullet)$ when $C$ is conilpotent (proposition \ref{radical}).
\end{enumerate}
\end{proof}

\bigskip
It is a classical result that $Tw_\bullet^c$ is representable by the cobar-bar adjunction $\Omega\dashv \Beta $ (see theorem \ref{classicbarcobar}), in consequence there must exist natural isomorphisms 
$$
\Omega C\simeq C\rhd_\bullet\mc
\et
\Beta A\simeq \{\mc,A\}_\bullet^c
$$
that we will unravel.

\medskip
We start by fixing some notations.
Let $u$ be the variable of degree $-1$ generating $\mc$, and $\FF u$ be the graded vector space generated by it.
Let $u^\star$ the dual variable, \ie the generator $(\FF u)^\star$ defeind by $u^\star(u)=1$.
We shall see $\FF u$ and $\FF u^\star$ as a graded or a dg-vector space and put 
$X\otimes u$ or $Xu$ instead of $X\otimes \FF u$ 
as well as $[u,X]$, or $X\otimes u^\star$, or $Xu^\star$ instead of $[\FF u,X]= X\otimes \FF u^\star$.
Let $s$ be a variable of degree $1$ generating the vector space $S=\FF s$, 
and $s^{-1}$ be another variable of degree $-1$ (unrelated to $s$) generating the vector space $S^{-1}=\FF s^{-1}$.
For $X$ a graded vector space or a dg-vector space $X$, we shall put $sX$ and $s^{-1}X$ instead of $SX=S\otimes X$ and $S^{-1}X=S^{-1}\otimes X$.

\medskip
We investigate now the isomorphism $\Omega C=C\rhd_\bullet \mc$.
For $C$ a pointed coalgebra, we have an isomorphism $s^{-1}C_- = C_-\otimes u$ given by $s^{-1}c= c\otimes u\ (-1)^{|c|}$.
Recall from theorems \ref{dgtoghomcoalg} and \ref{dgtoghomalg} that the functor $X\mapsto |X|$ sending a dg-(co)algebra to its underlying graded (co)algebra preserves all Sweedler operations.
From proposition \ref{Sweedleroftensorpointed} we have the isomorphisms of graded algebras
$$
|C\rhd_\bullet \mc| =  |C|\rhd_\bullet |\mc| = |C|\rhd_\bullet T(u) =T_\bullet(|C_-| u) = T_\bullet(s^{-1}|C_-|).
$$
This gives the classical description of $|\Omega C|$.

To described the differential we are going to use the equivalent language of non-(co)unital (co)algebras which is more convenient.
Let $d$ be the differential of $C$, from proposition \ref{uniqueextensionderSproductfreenu}, 
the differential $d \rhd_\bullet \mc + C\rhd_\bullet d_\mc$ induced by $d$ and $d_\mc$ is the unique pointed derivation extending
\begin{eqnarray*}
|C_-| u &\underset{-1}{\lrhup} & T_\circ(|C_-| u) \\
c\otimes u &\mto &  (dc)\rhd_\circ u + c\rhd_\circ (d_\mc u) \ (-1)^{|c|}\\
	&& = (dc)\rhd_\circ u + c\rhd_\circ (- u^2) \ (-1)^{|c|}\\
	&& = (dc)\otimes u - (c^{(1)}\otimes u)(c^{(2)}\otimes u) \ (-1)^{|c^{(1)}|}
\end{eqnarray*}
where we have used that, for $c\in C_-$, $c\rhd_\circ u = c\otimes u\in T_\circ(|C_-| u)$.

Transporting the structure along the isomorphism $T_\circ(|C_-| u) = T_\circ(s^{-1}|C_-|)$, we found that the differential of $C\rhd_\bullet \mc = T_\bullet(s^{-1}C_-)$ is the unique pointed derivation extending
\begin{eqnarray*}
s^{-1}|C_-| &\underset{-1}{\lrhup} & T_\circ(s^{-1}|C_-|) \\
s^{-1}c &\mto &  \underbrace{-s^{-1}dc}_{d^{int}(s^{-1}c)}   \underbrace{-(s^{-1}c^{(1)})(s^{-1}c^{(2)}) \ (-1)^{|c^{(1)}|}}_{d^{ext}(s^{-1}c)}.
\end{eqnarray*}

\begin{rem}\label{remdmcdextcobar}
This is exactly the definition of the differential of the cobar construction (see appendix \ref{cobarconstruction}).
The identification $s^{-1}=u$ leads to the correspondance:
\begin{eqnarray*}
d^{int} &=& d_C\rhd_\bullet \mc,\\
d^{ext} &=& C \rhd_\bullet d_\mc,\\
d^{tot} &=& d^{int} + d^{ext} = d_C\rhd_\bullet \mc + C \rhd_\bullet d_\mc,
\end{eqnarray*}
\ie the internal and external differentials corresponds exactly to the differential induced by the differential of $C$ and the Maurer-Cartan differential. This correspondance justifies our convention for the definition of $d^{ext}$ in \ref{cobarconstruction}.
\end{rem}

\bigskip
Let us now turn to the isomorphism $\Beta A=\{\mc,A\}_\bullet^c$.
For $A$ a pointed algebra, we have an isomorphism $s|A_-| = [u,A_-] = (A_-)u^\star$ given by $s a = a\otimes u^\star\ (-1)^{|a|}$.
From corollary \ref{corexSweedlerhompointed}, we have the isomorphisms of graded coalgebras
$$
|\{\mc,A\}_\bullet^c| = \{|\mc|,|A|\}_\bullet^c= \{T(u),|A|\}_\bullet^c = T^c_\bullet([u,|A_-|]) = T^c_\bullet(|A_-| u^\star) = T^c_\bullet(s|A_-|).
$$
which extract the classical description of $|\Beta  A|$.

As for the cobar construction, we are going to describe the differential in the more convenient language of non-(co)unital (co)algebras.
Let $d$ be the differential of $A$, proposition \ref{uniqueextensionderSHOMnufreeconil}, 
the differential $\{\mc,d\}_\bullet-\{A,d_\mc\}_\bullet$ induced by $d$ and $d_\mc$ is the unique pointed coderivation coextending
\begin{eqnarray*}
b:T_\circ^c(|A_-| u^\star) &\underset{-1}{\lrhup} & |A_-| u^\star \\
h &\mto & dh^\flat i-h^\flat d_\mc i \ (-1)^{|h||d_\mc|}
\end{eqnarray*}
Let us simplify this expression. We have 
\begin{eqnarray*}
\big(dh^\flat i-h^\flat d_\mc i \ (-1)^{|h||d_1|}\big)(u) &=& d(h^\flat(u)) - h^\flat( d_\mc(u)) \ (-1)^{|h||d_\mc|}\\
 &=& d(h^\flat(u)) - h^\flat( -u^2) \ (-1)^{|h|}\\
 &=& d(h^\flat(u)) + h^\flat(u^2) \ (-1)^{|h|}.
\end{eqnarray*}
An homogeneous element $h\in T^c(|A_-| u^\star)$ is of the type $h=a_1u^\star \otimes \dots \otimes a_n u^\star$,
and $h^\flat=\pm (a_1\dots a_n)(u^\star)^{\otimes n}$. 
Then,
$h^\flat(u)$ is non zero iff $h=au^\star$ in which case $d(h^\flat(u)) = da$
and 
$h^\flat(u^2)$ is non zero iff $h=a u^\star\otimes bu^\star$ in which case 
$$
h^\flat(u^2) = ab (u^\star\otimes u^\star)(u^2)\ (-1)^{|b|} = ab (- u^\star(u) u^\star(u)) \ (-1)^{|b|} = - ab \ (-1)^{|b|}
$$
(we have used the computation of remark \ref{remudual}).

This computation characterizes the differential as the only pointed coderivation coextending
\begin{eqnarray*}
b:T_\circ^c(|A_-| u^\star) &\underset{-1}{\lrhup} & |A_-| u^\star \\
a u^\star &\mto & (da) u^\star\\
au^\star\otimes bu^\star &\mto & - (ab) u^\star \ (-1)^{|a|}\\
\textrm{other terms} & \mto & 0.
\end{eqnarray*}

Recall the isomorphism $s|A_-| = |A_-|u^\star$ given by $s a= a\otimes u^\star\ (-1)^{|a|}$.
Transporting the structure along this isomorphism, we found that the differential of $\{\mc, A\}^c_\bullet = T(sA_-)$ is the unique differential extending
$$
\begin{array}{rcll}
T_\circ^c(s|A_-|) &\underset{-1}{\lrhup} & s|A_-|\\
\rule[-2ex]{0pt}{4ex}	sa			&\mto &	-s(da) 			& = d^{int}(sa)\\
\rule[-2ex]{0pt}{4ex}	sa\otimes s b	&\mto &	-s(ab)\ (-1)^{|a|}	& = - d^{ext}(sa\otimes sb)\\
\rule[-2ex]{0pt}{4ex}	\textrm{other terms} & \mto & 0.
\end{array}
$$

\begin{rem}\label{remdmcdextbar}
We recognize the formulas from the bar construction (see appendix \ref{barconstruction}).
The identification $s=u^\star$ leads to the correspondance:
\begin{eqnarray*}
d^{int} &=& \{\mc,d_A\}_\bullet,\\
d^{ext} &=& \{d_\mc,A\}_\bullet,\\
d^{tot} &=& d^{int} - d^{ext} = \{\mc,d_A\}_\bullet - \{d_\mc,A\}_\bullet,
\end{eqnarray*}
\ie the internal and external differentials corresponds exactly to the differential induced by the differential of $C$ and the Maurer-Cartan differential. This nice correspondance justifies our convention for the definition of $d^{ext}$ in \ref{barconstruction}.
\end{rem}

\bigskip
We now turn to the study of universal twisting cochains. Although the definitions make sense in any of the unpointed, pointed and conilpotent context, we will study only the conilpotent case.

\begin{defi}
If $C$ is a conilpotent dg-coalgebra, we shall say that a pointed twisting cochain $\alpha:C\rhup_{-1} A$ is {\em universal} if the pair $(A,\alpha)$ represents the functor 
$$\xymatrix{
Tw_\bullet^c(C,-) :\dgAlg_\bullet \ar[r]&\Set.
}$$
If $A$ is a dg-algebra, we shall say that a twisting cochain $\alpha:C\rhup_{-1} A$ is {\em couniversal} if the pair $(C,\alpha)$ represents the functor 
$$\xymatrix{
Tw_\bullet^c(-,A) :(\dgCoalgnilbullet)^{op}\ar[r]&\Set.
}$$
\end{defi}

\medskip
Because of the representability of $Tw_\bullet(-,A)$ and $Tw_\bullet(C,-)$
by $\Beta A$ and $\Omega C$, the (co)universal twisting morphisms can be described as a graded morphisms $\beta:\Beta A\rhup A$ and $\omega:C\rhup \Omega C$ that we are going to compute.

\begin{prop}\label{twistingcochain}
\begin{enumerate}
\item The couniversal twisting cochain is the pointed graded morphism $\beta:\Beta  A\rhup A$ of degree $-1$
	\begin{eqnarray*}
	\beta:\Beta  A=T^c_\bullet(sA_-) &\underset{-1}{\lrhup}& A\\
	1 &\mto & e_A\\
	s a &\mto & -a\\
	\textrm{\em other terms} &\mto & 0.
	\end{eqnarray*}
\item The universal twisting cochain $\omega:C \rhup_{-1}\Omega C$ is the graded morphism of degree $-1$
	\begin{eqnarray*}
	\omega:C&\underset{-1}{\lrhup} & \Omega C = T_\bullet(s^{-1}C_-)\\
	e_c &\mto & 1\\
	c\in C_- &\mto & s^{-1}c
	\end{eqnarray*}
\end{enumerate}
\end{prop}

\begin{proof}

\begin{enumerate}

\item To compute $\beta$, we are going to use the non-counital langage, which is more convenient.
By lemma \ref{mcalgpointed}, we have the following correspondance
\begin{center}
\begin{tabular}{lc}
\rule[-2ex]{0pt}{4ex} the identity map & $id:\Beta  A\to \Beta  A$,\\
\rule[-2ex]{0pt}{4ex} the canonical inclusion (by property of the radical) & $\iota:\Beta  A\to \Betaext A=\{\mc,A\}_\bullet$,\\
\rule[-2ex]{0pt}{4ex} the canonical inclusion & $\iota_-:(\Beta  A)_-\to \{\mc_-,A_-\}_\circ$,\\
\rule[-2ex]{0pt}{4ex} a certain non-unital measuring & $\beta_{\mc_-}:(\Beta  A)_-\otimes \mc_-\to A_-$,\\
\rule[-2ex]{0pt}{4ex} a non-unital algebra map & $\lceil \beta_-\rceil: \mc_-\to [(\Beta  A)_-,A_-]$,\\
\rule[-2ex]{0pt}{4ex} the couniversal non-unital twisting cochain & $\beta_-:(\Beta  A)_-\rhup A_-$.
\end{tabular}
\end{center}
Forgetting the differentials we have $|(\Beta  A)_-\otimes \mc_-| = T^c_\circ(|A_-|u^\star) \otimes T_\circ(u)$
and by the graded analog of proposition \ref{corexSweedlerhomnu} we found that the measuring $\beta_{\mc_-}$ is the unique non-unital measuring extending the map 
\begin{eqnarray*}
\beta':T_\circ^c(|A_-|u^\star)\otimes u &\tto & |A_-|\\
(a_1u^\star\otimes \dots \otimes a_nu^\star)\otimes u &\mto& 
\left\{\begin{array}{cl}
a_1& \textrm{ if $n=1$}\\
0 & \textrm{ if $n\not=1$}.
\end{array}\right.
\end{eqnarray*}
From lemma \ref{lemmaTCmes}, the twisting cochain $\beta_-:(\Beta  A)_-=T^c_\circ(sA_-)\rhup A_-$ corresponding to $\beta'$ is computed by the formula
$$
\beta_-(c) =\beta_{\mc_-}(c,u) \ (-1)^{|c|} = \beta'(c,u) \ (-1)^{|c|}.
$$
Hence, $\beta_-$ is the map
\begin{eqnarray*}
\beta_-:T^c_\circ(|A_-|u^\star ) &\underset{-1}{\lrhup}& A_-\\
au^\star&\mto& a\ (-1)^{|a|+1}\\
\textrm{other terms} &\mto & 0.
\end{eqnarray*}
Then, using the isomorphism $sA_- = A_- u^\star $ given by $sa = au^\star \ (-1)^{|a|}$ we can write the couniversal twisting cochain as
\begin{eqnarray*}
\beta_-:T^c_\circ(s|A_-|) &\underset{-1}{\lrhup}& A_-\\
sa&\mto& -a\\
\textrm{other terms} &\mto& 0.
\end{eqnarray*}

\item We proceed similarly to compute $\omega$.
By lemma \ref{mcalgpointed}, we have the following correspondance
\begin{center}
\begin{tabular}{lc}
\rule[-2ex]{0pt}{4ex} the identity map & $id:\Omega C=C\rhd_\bullet \mc \to \Omega C$,\\
\rule[-2ex]{0pt}{4ex} the identity map & $id:C_-\rhd_\circ \mc_- \to (\Omega C)_-$,\\
\rule[-2ex]{0pt}{4ex} the universal non-unital measuring & $\omega_{\mc_-}:C_-\otimes \mc_- \to (\Omega C)_-$,\\
\rule[-2ex]{0pt}{4ex} a non-unital algebra map & $\lceil \omega_-\rceil: \mc_-\to [C_-,(\Omega C)_-]$,\\
\rule[-2ex]{0pt}{4ex} the universal non-unital twisting cochain & $\omega_-:C_-\rhup(\Omega C)_-$.
\end{tabular}
\end{center}
Forgetting the differentials we have $|(\Omega C)_-| = T_\circ(|C_-|u)$
and by the graded analog of proposition \ref{Sweedleroftensornu} we found 
that the measuring $\omega_{\mc_-}$ is the unique non-unital measuring extending the canonical inclusion
$\omega':|C_-|u \to T_\circ(|C_-|u)$.
From lemma \ref{lemmaTCmes}, the couniversal twisting cochain $\omega_-:C_-\to (\Omega C)_-$ corresponding to $\omega'$ is computed by the formula
$$
\omega_-(c) =\omega_{\mc_-}(c,u) \ (-1)^{|c|} = \omega'(c,u) \ (-1)^{|c|}.
$$
Hence, $\omega_-$ is the map
\begin{eqnarray*}
\omega_-:|C_-| &\underset{-1}{\lrhup}& T_\circ(|C_-|u)\\
c &\mto& c \otimes u \ (-1)^{|c|}.
\end{eqnarray*}
Then, using the isomorphism $s^{-1}C_- = C_- u$ given by $s^{-1}c = cu \ (-1)^{|c|}$, the universal twisting cochain is 
\begin{eqnarray*}
\omega_-:|C_-| &\underset{-1}{\lrhup}& T_\circ(s^{-1}|C_-|)\\
c &\mto& s^{-1} c.
\end{eqnarray*}

\end{enumerate}
\end{proof}

\begin{rem}
Note that the minus sign in the definition of the universal twisting cochain $\beta:\Beta A\rhup A$ appeared because of the formula for the lambda transform in the first variable.
\end{rem}

\subsubsection{Consequences of Sweedler formalism}\label{csqsweedlerformalism}

In this section, we use the enrichment of Sweedler operations to deduce some facts about the bar and cobar constructions.

\begin{defi}
By analogy with the classical notations for bar and cobar constructions, we shall use the following notations:
\begin{center}
\begin{tabular}{clcl}
\rule[-2ex]{0pt}{4ex} $\bOmega C = C\rhd\mc $ & and & $\bB A = \{\mc,A\}$ & in the unpointed case,\\
\rule[-2ex]{0pt}{4ex} $\Omega C = C\rhd_\bullet\mc $ & and & $\Betaext A = \{\mc,A\}_\bullet$ & in the pointed case, and\\
\rule[-2ex]{0pt}{4ex} $\Omega C = C\rhd_\bullet\mc $ & and & $\Beta  A = R^c\{\mc,A\}_\bullet$ & in the conilpotent (\ie classical) case.
\end{tabular}
\end{center}
\end{defi}

By definition we have the bijections of sets
\begin{center}
\begin{tabular}{cl}
\rule[-2ex]{0pt}{4ex} $\dgAlg(\bOmega C, A) = Tw(C,A) = \dgCoalg(C,\bB A)$ & in the unpointed case,\\
\rule[-2ex]{0pt}{4ex} $\dgAlg_\bullet(\Omega C, A) = Tw_\bullet(C,A) = \dgCoalg_\bullet(C,\Betaext A)$ & in the pointed case, and\\
\rule[-2ex]{0pt}{4ex} $\dgAlg_\bullet(\Omega C, A) = Tw_\bullet^c(C,A) = \dgCoalg_\bullet(C,\Beta  A)$ & in the conilpotent case.
\end{tabular}
\end{center}
The fact that Sweedler operations are all enriched over $\dgCoalg$ strenghten these bijections into coalgebra isomorphisms.
In particular, this leads to an interpretation of the bar construction of the convolution algebra.

\begin{prop}\label{isobarcobar}
There exists canonical isomorphisms of dg-coalgebras:
\begin{center}
\begin{tabular}{cl}
\rule[-2ex]{0pt}{4ex} $\{\bOmega C,A\} = \bB[C,A] = \HOM(C,\bB A)$ & in the unpointed case,\\
\rule[-2ex]{0pt}{4ex} $\{\Omega C,A\}_\bullet = \Betaext [C,A]_\bullet = \HOM_\bullet(C,\Betaext A)$ & in the pointed case, and\\
\rule[-2ex]{0pt}{4ex} $\{\Omega C,A\}_\bullet = \Betaext [C,A]_\bullet = \HOM_\bullet(C,\Beta  A)$ & in the conilpotent case.
\end{tabular}
\end{center}
\end{prop}
\begin{proof}
The enrichment of the bar cobar adjunction gives isomorphisms of pointed dg-coalgebras
$$
\{C\rhd_\bullet  \mc ,A\}_\bullet = \{\mc ,[C,A]_\bullet\}_\bullet = \HOM_\bullet(C,\{\mc,A\}_\bullet).
$$
The proof is similar in the unpointed case.
For the conilpotent case we have used proposition \ref{radical3pointed} to have
$$
\HOM_\bullet(C,\{\mc,A\}_\bullet) = \HOM_\bullet(C,R^c\{\mc,A\}_\bullet).
$$
\end{proof}

\begin{cor}\label{cormcdual}
For any dg-coalgebra $C$,
$$
\bB(C^\star) = (\bOmega C)^\vee.
$$
\end{cor}
\begin{proof} 
This is proposition \ref{isobarcobar} applied to $A=\FF$.
\end{proof}

The next result says in particular that the fundamental object $\mc$ can be defined from the $(\bB,\bOmega)$ adjunction.

\begin{cor}
$$
\mc = \bOmega \FF\ ,\qquad
\mc^\vee = \bB \FF
\et 
\mc = \Omega \FF_+
$$
\end{cor}
\begin{proof}
The first formula is true because $\FF\rhd \mc = \mc$.
The second formula is deduced from corollary \ref{cormcdual}.
The third formula is $\FF_+\rhd\bullet \mc = \mc$.
\end{proof}

The following corollary states in a sense that the bar and cobar constructions are "enriched over the bar construction".
These isomorphisms bear a strong similarity to those of \cite[ch. 5.7]{Keller} and can be use as a basis to construct an "hom $A_\infty$-algebras" between $A_\infty$-(co)algebras.

\begin{cor}\label{hombar}
There exists canonical isomorphisms of dg-coalgebras
$$
\HOM(\bB A,\bB A') = \bB [\bB A,A']
\ , \qquad
\{\bOmega C',\bOmega C\} = \bB [C',\bOmega C]\ ,
$$
$$
\HOM_\bullet(\Betaext A,\Betaext A') = \Betaext [\Betaext A,A']_\bullet
\et 
\{\Omega C',\Omega C\}_\bullet = \Betaext [C',\Omega C]_\bullet.
$$
We have also
$$
\HOM_\bullet(\Beta  A,\Beta  A') = \Betaext [\Beta  A,A']_\bullet.
$$
In particular, $\bB [\bB A,A]$, $\bB [C,\bOmega C]$, $\Betaext [\Betaext A,A]_\bullet$, $\Betaext [C,\Omega C]_\bullet$ and $\Betaext [\Beta A,A]_\bullet$ have a canonical bialgebras structure. 
In the bialgebras $\Betaext [C,\Omega C]_\bullet$ and $\Betaext [\Beta A,A]_\bullet$ the identity atom correspond to the (co)universal twisting cochains.
\end{cor}
\begin{proof} 
This is direct from proposition \ref{isobarcobar}.
\end{proof}

\bigskip

\begin{prop}
\begin{enumerate}
\item For any dg-algebra $A$, there exists a canonical bijection between
\begin{enumerate}
\item Maurer-Cartan elements of $A$
\item and coaugmentations of $\bB A$.
\end{enumerate}
\item For any dg-coalgebra $C$, there exists a canonical bijection between 
\begin{enumerate}
\item Maurer-Cartan elements of $C^\star$, 
\item coaugmentations of $\bB C^\star=(\bOmega C)^\vee$ 
\item and augmentations of $\bOmega C$.
\end{enumerate}
\end{enumerate}
\end{prop}
\begin{proof} 
\begin{enumerate}
\item We have $\dgAlg(\mc, A) = \dgCoalg(\FF,\{\mc, A\})$.
\item We have $\dgAlg(C\rhd \mc ,\FF) = \dgAlg(\mc ,[C,\FF]) = \dgCoalg(\FF,\{\mc, [C,\FF]\})$.
\end{enumerate}
\end{proof}

\begin{rem}[A moduli interpretation]
We present here a moduli interpretation of the bar construction inspired by the previous result, this is merely a new vocabulary, but it can be useful to think and echoes to constructions of \cite{Hinich, Getzler} and others.
If we call a map $C\rhd \mc \to A$ a {\em family of Maurer-Cartan elements of $A$ parametrized by $C$}, 
the {\em moduli functor of Maurer-Cartan elements of $A$} is
\begin{eqnarray*}
MC(A):\dgCoalg^{op} &\tto &\Set\\
C & \mto & \dgAlg(C\rhd \mc, A).
\end{eqnarray*}
It is obvious that the coalgebra $\{\mc, A\}=\bB A$ represents this functor.
In particular, the set of rational points of this functor, \ie the value at $\FF$, is the set of Maurer-Cartan elements of $A$.
\end{rem}

\subsubsection{Iterated bar constructions}\label{iteratedbarcobar}

From the (co)lax structure of Sweedler operations it is easy to deduce the existenc of (co)shuffle products on the bar and cobar constructions applied to (co)commutative (co)algebras \cite{EM,HMS}.
\begin{thm}\label{barcobarcommutative}
\begin{enumerate}

\item If $C$ is a cocommutative coalgebra (resp. a cocommutative pointed coalgebra), the cobar construction $\bOmega C$ (resp. $\Omega C$) is a cocommutative bialgebra. Moreover, in both cases, the cocommutative coproduct is the coshuffle coproduct.

\item If $A$ is a commutative algebra (resp. a commutative pointed algebra), the bar construction $\bB A$ (resp. $\Betaext A$ and $\Beta A$) is a commutative bialgebra. Moreover, in all cases, the commutative product is the shuffle product.

\end{enumerate}
\end{thm}

\begin{proof}
Recall from proposition \ref{mchopf} that $\mc$ is a cocommutative Hopf algebra.
Then, the fact tha the image have a bialgebra structure follows from 
propositions \ref{corcolaxSproduct} and \ref{corlaxShom} and their pointed analogs 
propositions \ref{corcolaxSproductpointed} and \ref{corlaxShompointed}.
For $\Beta A$ we deduce the result from the structure of $\Betaext A$ and proposition \ref{laxcorad}.

Let us now prove that the (co)products are given by the (co)shuffle.
By the computation after proposition \ref{corcolaxSproduct}, the coproduct on $C\rhd \mc$ is the composition
$$\xymatrix{
C\rhd \mc \ar[rr]^-{\Delta\rhd \Delta} && (C\otimes C)\rhd (\mc\otimes \mc)\ar[r]^\alpha &(C\rhd \mc)\otimes (C\rhd \mc)
}$$
where $\alpha $ is the colax structure of $\rhd$.
Explicitely on elements, this gives
\begin{eqnarray*}
\Delta(c\rhd u)&=& (c^{(1)}\rhd u)\otimes (c^{(2)}\rhd 1) (-1)^{|c^{(2)}|} + (c^{(1)}\rhd 1)\otimes (c^{(2)}\rhd u) \\
	&=& (c^{(1)}\rhd u)\otimes \epsilon(c^{(2)}) + \epsilon(c^{(1)})\otimes (c^{(2)}\rhd u) \\
	&=& (c\rhd u)\otimes 1 + 1\otimes (c\rhd u).
\end{eqnarray*}
Using $|\bOmega(C)| = |C\rhd \mc| = T(Cu)$ we can recognize the coshuffle coproduct formula of example \ref{coshuffle bialgebra}.

The computation is analog for $\Omega C=C\rhd_\bullet \mc$, we leave it the the reader.

\medskip

Similarly, the product on $\{\mc,A\}$ is defined as the composite 
$$\xymatrix{
\{\mc,A\}\otimes \{\mc,A\} \ar[r]^-{\alpha} & \{\mc\otimes \mc,A\otimes A\} \ar[rr]^-{\{\Delta, \mu\}} && \{\mc,A\}
}$$
where $\alpha $ is the lax structure of $\{-,-\}$, which is also the strength of the $\otimes$. Using this last remark and the calculus of meta-morphisms, the product of $f,g\in \{\mc,A\}$ can be written as $\mu(f\otimes g)\Delta$.

Passing to the underlying graded objects the product is a map $T^\vee(Au^\star)\otimes T^\vee(Au^\star) \to T^\vee(Au^\star)$.
The composition 
$$\xymatrix{
\phi:T^\vee(Au^\star)\otimes T^\vee(Au^\star) \ar[r] &T^\vee(Au^\star) \ar[r]^-p & Au^\star = [u,A],
}$$
where $p$ is the cogenerating map, sends $(f,g)$ to the function 
$$
\phi(f,g):u\mto (\mu(f\otimes g)\Delta)(u) = f(u)g(1) \ (-1)^{|g||u|} + f(1)g(u)
$$
where $1$ is the unit of $\mc$.
By the caracterisation of unital meta-morphisms in proposition \ref{propmetanutou} we have $f(1) = \epsilon(f).1\in A$ for any $f\in \{\mc,A\}$ and $\phi(f,g)$ is the function
$$
u\mto f(u)\epsilon(g) + \epsilon(f)g(u)
$$
which is exactly the definition of the shuffle product on $T^\vee(Au^\star)$ of corollary \ref{shuffleproductTv}.

The computations are similar for $\Betaext A$ and $\Beta A$, we leave them the the reader.
\end{proof}

In particular, we find that we can iterate the bar and cobar constructions (pointed, conilpotent or not) defined on (co)commutative (co)algebras. The abstract form of this result is given by the corollaries \ref{iterationcobar} and \ref{iterationbar} and their pointed analogs corollaries \ref{iterationcobarpointed} and \ref{iterationbarpointed}.

\subsubsection{Generalized bar-cobar adjunctions}\label{variationbarcobar}
The bar-cobar adjunction generalizes in adjunctions of the type $-\rhd K \dashv \{K,-\}$ (pointed or not) where $K$ is any dg-algebra.
The relations (BCB) are generalized in
$$
\dgAlg(C\rhd K, A) = \dgAlg(K, [C,A]) = \dgCoalg(C, \{K,A\})
$$
and this says that the functors $-\rhd K \dashv \{K,-\}$ represents the set of certain elements of the convolution algebra classified by maps $K\to [C,A]$.

\medskip
A noticeable example is when $K=T(x)$ is the free graded algebra on one variable $x$ of degree $-1$ with the zero differential, which corresponds to the replacement of the Maurer-Cartan equation by the equation
$$
dx=0.
$$
Then a map $T(x)\to [C,A]$ is simply a map of dg-vector spaces $C\to A$.
We leave the reader to check that the adjunction $-\rhd T(x) \dashv \{T(x),-\}$ is the composition
$$\xymatrix{
-\rhd T(x):\dgCoalg \ar@<.6ex>[rr]^-U && \dgVect \ar@<.6ex>[ll]^-{T^\vee}\ar@<.6ex>[rr]^-T && \dgAlg:\{T(x),-\}. \ar@<.6ex>[ll]^-{U}
}$$
If $|x|=n$, some suspensions are involved.

\medskip
Another example is when $K$ is the free differential graded algebra on one variable $x$ of degree $-1$. As a graded algebra it is the free tensor algebra $T(x,dx)$ where $dx$ is in degree $-2$. This corresponds to replace the Maurer-Cartan equation by no equation at all!

A map $K\to [C,A]$ is a map of graded vector spaces $C\to A$.
We leave the reader to check that 
$$
C\rhd K = T(Cx\oplus Cdx) \et \{K,A\} = T^\vee(Ax^\star \oplus A(dx)^\star)
$$
These constructions bear a similarity with the Weil algebra of a Lie algebra \cite{Cartan}.
Sweedler theory seems to provides a context where to develop a non-commutative analog of Weil theory.

\medskip
Finally corollary \ref{exgenbarcobar} give also some generalizations of the barcobar adjunctions.
They corresponds to the case where $K=T(X)$ for some dg-vector space $X$.

\newpage
\appendix

\section{Classical bar and cobar constructions}\label{barcobar}

This appendix contains some recollections on the bar and cobar construction after \cite{Ad2, EB, EM, LV, Proute}. 
Its main purpose is to establish the good signs in the formulas for the differential.
We shall focus only on the differential graded approach and say nothing about the underlying simplicial structures.

\medskip

Recall that, for $X$ a graded vector space, the non-unital graded coalgebra $T_\circ^c(X)$ is $\oplus_{n>0}X^{\otimes n}$, we shall denote $p_n:T_\circ^c(X)\to X^{\otimes n}$ the projection to the $n$-th factor.
Similarly, the non-unital graded algebra $T_\circ(X)$ is $\oplus_{n>0}X^{\otimes n}$ and we shall denote $i_n:X^{\otimes n}\to T_\circ(X)$ the inclusion of the $n$-th factor.
Recall also that if $X$ is a dg-vector space, $|X|$ is its underlying graded object.

\subsection{The bar construction}\label{barconstruction}

Recall that, for $(A,d)$ an augmented dg-algebra, $A_-$ is the associated non-unital dg-algebra.
We consider the non-unital graded coalgebra $T_\circ^c(s|A_-|)$.
We consider the two morphisms of degree $-1$:

$$\xymatrix{
b^{int}: T_\circ^c(s|A_-|)\ar[r]^-{p}& s|A_-|\ar@{-^{>}}[rr]^-{s\otimes d}_>>>{-1}&& s|A_-|,
}$$
sending $sa$ to $-s(da)$ and other elements to 0; 
and
$$\xymatrix{
b^{ext}: T_\circ^c(s|A_-|)\ar[r]^-{p_2} & s|A_-|\otimes s|A_-|\ar[r]^-{\simeq}& s^2(|A_-|\otimes |A_-|)\ar@{-^{>}}[rr]^-{(s^2\mapsto s)}_>>>{-1}&&s(|A_-|\otimes |A_-|) \ar[rr]^-{s\otimes m} && s|A_-|
}$$
sending $sa\otimes sb$ to $s(ab)\ (-1)^{|a|}$ and other elements to 0.

\medskip

By proposition \ref{coextensioncoderconilpointed}, $b^{int}$ and $b^{ext}$ extend to two pointed coderivations of degree $-1$ of $T_\circ^c(s|A_-|)$, respectively noted $d^{int}$ and $d^{ext}$, given the formulas
\begin{eqnarray*}
d^{int}(sa_1\otimes \dots \otimes sa_n) &=& \sum_{1\leq i\leq n} sa_1\otimes \dots \otimes -sd(a_i)\otimes \dots \otimes  sa_n\ (-1)^{|sa_1|+\dots +|sa_{i-1}|}\\
	&=& \sum_{1\leq i\leq n} sa_1\otimes \dots \otimes sd(a_i)\otimes \dots \otimes  sa_n\ (-1)^{i+|a_1|+\dots +|a_{i-1}|}\\
\end{eqnarray*}
and
\begin{eqnarray*}
d^{ext}(sa_1\otimes \dots \otimes sa_n) &=& \sum_{1\leq i<n} sa_1\otimes \dots \otimes s(a_ia_{i+1})\otimes \dots \otimes  sa_n\ (-1)^{|a_i|+|sa_1|+\dots +|sa_{i-1}|}\\
	&=& \sum_{1\leq i<n} sa_1\otimes \dots \otimes s(a_ia_{i+1})\otimes \dots \otimes  sa_n\ (-1)^{i-1 + |a_1|+\dots +|a_i|}.
\end{eqnarray*}

\begin{lemma}
The coderivations $d^{int}$ and $d^{ext}$ commute.
Moreover, $d^{int}$, $d^{ext}$ and $d^{int}\pm d^{ext}$ are coderivations of square zero.
\end{lemma}
\begin{proof}
Let us prove the commutation of $d^{int}$ and $d^{ext}$. Because $d^{int}$ and $d^{ext}$ are of odd degrees, their graded commutator is $d^{ext}d^{int}+d^{int}d^{ext}$.
It is again a coderivation by stability of coderivations by the commutator.

Thus, by corollary \ref{caraccoderconil}, $d^{ext}d^{int}+d^{int}d^{ext}=0$ if and only if $b^{ext} d^{int}+b^{int} d^{ext}=0$.
It is easy to see that $b^{ext} d^{int}+\alpha d^{ext}=0$ when applied to an element $sa_1\otimes \dots \otimes sa_n$ for $n\not=2$.
Then, if $n=2$, we have on one side
\begin{eqnarray*}
b^{ext} d^{int}(sa_1\otimes sa_2) &=& b^{ext} \left(-sd(a_1)\otimes sa_2 + sa_1\otimes -sd(a_2)\ (-1)^{1+|a_1|}\right)\\
	&=& b^{ext} \left(s^2d(a_1)\otimes a_2\ (-1)^{1+(|a_1|-1)} + s^2a_1\otimes d(a_2)\ (-1)^{|a_1|+|a_1|})\right)\\
	&=& s(da_1)a_2\ (-1)^{|a_1|} + sa_1(da_2)\\
	&=& sd(a_1a_2)\ (-1)^{|a_1|}
\end{eqnarray*}
and on the other side
\begin{eqnarray*}
b^{int} d^{ext}(sa_1\otimes sa_2) &=& b^{int} (s(a_1a_2))\ (-1)^{|a_1|}\\
	&=& -sd(a_1a_2)\ (-1)^{|a_1|}
\end{eqnarray*}
hence $b^{ext} d^{int}+b^{int} d^{ext}=0$.

Let us prove now that $d^{int}$ and $d^{ext}$ are of square zero.
Using again corollary \ref{caraccoderconil}, $d^{int}d^{int}=0$ if and only if $b^{int} d^{int}=0$.
It is easy to see that $b^{int} d^{int}=0$ when applied to an element $sa_1\otimes \dots \otimes sa_n$ for $n\not=1$, and, if $n=1$, $b^{int} d^{int}(sa_1) = -b^{int} sda_1 = s d(da_1) = 0$.
Similarly, $d^{ext}d^{ext}=0$ if and only if $b^{ext} d^{ext}=0$.
It is easy to see that $\beta d^{ext}=0$ when applied to an element $sa_1\otimes \dots \otimes sa_n$ for $n\not=3$, and, if $n=3$,
\begin{eqnarray*}
b^{ext} d^{ext}(sa_1\otimes sa_2\otimes sa_3) &=& b^{ext}\left( s(a_1a_2)\otimes sa_3\ (-1)^{|a_1|} + sa_1\otimes s(a_2a_3)\ (-1)^{|a_2|+(|a_1|+1)}\right)\\
	&=& s((a_1a_2)a_3)\ (-1)^{|a_1|+|a_2|+|a_1|} + s(a_1(a_2a_3))\ (-1)^{|a_1|+|a_2|+|a_1|+1}\\
	&=& 0
\end{eqnarray*}
by associativity.

The last statement follows from $(d^{int}\pm d^{ext})^2 = (d^{int})^2 \pm (d^{int}d^{ext} + d^{ext}d^{int}) + (d^{ext})^2= (d^{int})^2 + (d^{ext})^2$ .
\end{proof}

We shall call $d^{int}$ the {\em internal differential} and $d^{ext}$ the {\em external differential}. 
Recall that the non-unital graded coalgebra $T_\circ^c(s|A_-|)$ is equivalent to the pointed graded coalgebra $T^c_\bullet(s|A_-|)$ and that coderivations of $T_\circ^c(s|A_-|)$ corresponds bijectively to pointed coderivations of $T^c_\bullet(s|A_-|)$. We shall again call $d^{int}$ and $d^{ext}$ the internal and external differentials of $T^c_\bullet(s|A_-|)$ induced by that of $T_\circ^c(s|A_-|)$. Remark that the conilpotent pointed dg-coalgebra $(T^c_\bullet (s|A_-|),d^{int})$ is simply $T^c_\bullet(sA_-)$.

\medskip

The {\em bar construction} of a pointed dg-algebra $A$ is denoted $\Beta  A$ and is defined as the conilpotent pointed dg-coalgebra 
$$
(T^c_\bullet(s|A_-|),d^{int} - d^{ext}).
$$
The bar construction $\Beta $ defines a functor
$$\xymatrix{
\Beta :\dgAlg_\bullet \ar[r]&\dgCoalgnilbullet.
}$$

\begin{rem}\label{signbar}
The sign in the definition of the total differential of $\Beta A$ depends on the authors, in \cite{EM}, the differential is $d^{int} - d^{ext}$, but in \cite{EB, LV,Proute}, the differential is $d^{int} + d^{ext}$. Our convention is chosen to agree with our construction of $\Omega$ via Sweedler operations (see also section \ref{signtw} for a discussion of signs).
The following lemma proves that the two definitions of the bar construction are isomorphic.
\end{rem}

Let $X$ be a graded vector space, the graded space $T(X)=\bigoplus_{n\geq 0}X^{\otimes n}$ has an automorphism $\pi$ which send $x\in X^{\otimes n}$ to $x\ (-1)^n \in X^{\otimes n}$. Moreover $\pi$ respects the algebra and coalgebra structures of $T(X)$.

\begin{lemma}
The map $\pi:T^c_\bullet(s|A_-|)\to T^c_\bullet(s|A_-|)$ induces an isomorphism between the two bar constructions
$$
(T^c_\bullet(s|A_-|),d^{int} - d^{ext}) \et  (T^c_\bullet(s|A_-|),d^{int} + d^{ext}).
$$
\end{lemma}
\begin{proof}
We want to prove $\pi^{-1}(d^{int} + d^{ext})\pi = d^{int} - d^{ext}$.
Let $X = s|A_-|$, $d^{int}$ send $X^{\otimes n}$ to itself so $\pi^{-1}(d^{int})\pi = d^{int}$ and
$d^{ext}$ sends $X^{\otimes n}$ to $X^{\otimes n-1}$ so $\pi^{-1}(d^{ext})\pi = -d^{ext}$.
\end{proof}

\subsection{The cobar construction}\label{cobarconstruction}

Recall that, for $(C,d)$ an augmented dg-coalgebra, $C_-$ is the associated non-unital dg-coalgebra.
We consider the non-unital graded algebra $T_\circ(s^{-1}|C_-|)$.
We consider the two morphisms of degree $-1$:
$$\xymatrix{
b^{int}: s^{-1}|C_-|\ar@{-^{>}}[rr]^-{s^{-1}\otimes d}_>>>{-1}&& s^{-1}|C_-|\ar[r]^-{i_1}& T_\circ(s^{-1}|C_-|)
}$$
sending $s^{-1}c$ to $-s^{-1}dc$;
and
$$\xymatrix{
b^{ext}: s^{-1}|C_-|\ar[rr]^-{s^{-1}\otimes \Delta}&& s^{-1}|C_-|\otimes |C_-| \ar@{-^{>}}[rr]^-{s^{-1}\mapsto -s^{-2}}_>>>{-1} && s^{-2}|C_-|\otimes |C_-| \ar[r]^-\simeq & s^{-1}|C_-|\otimes s^{-1}|C_-| \ar[r]^-{i_2}&T_\circ(s^{-1}|C_-|)
}$$
sending $s^{-1}c$ to $-s^{-1}c^{(1)}\otimes s^{-1}c^{(2)}\ (-1)^{|c^{(1)}|}$.

Some explanation on the minus sign in $b^{ext}$ will be given below, see remark \ref{signcobar}.

\medskip

By proposition \ref{extensionderivationpointed}, $b^{int}$ and $b^{ext}$ extend to two pointed derivations of degree $-1$ of $T_\circ(s^{-1}|C_-|)$, respectively noted $d^{int}$ and $d^{ext}$ given by the formulas 
\begin{eqnarray*}
d^{int}(s^{-1}c_1\otimes \dots \otimes s^{-1}c_n) &=& \sum_{1\leq i\leq n} s^{-1}c_1\otimes \dots \otimes -s^{-1}d(c_i)\otimes \dots \otimes  s^{-1}c_n\ (-1)^{|s^{-1}c_1|+\dots +|s^{-1}c_{i-1}|}\\
	&=& \sum_{1\leq i\leq n} s^{-1}c_1\otimes \dots \otimes s^{-1}d(c_i)\otimes \dots \otimes  s^{-1}c_n\ (-1)^{i+|c_1|+\dots +|c_{i-1}|}\\
\end{eqnarray*}
and
\begin{eqnarray*}
d^{ext}(s^{-1}c_1\otimes \dots \otimes s^{-1}c_n) &=& \sum_{1\leq i<n} s^{-1}c_1\otimes \dots \otimes (-s^{-1}c_i^{(1)}\otimes s^{-1}c_i^{(2)}) \otimes \dots \otimes  s^{-1}c_n\ (-1)^{|c_i^{(1)}|+|s^{-1}c_1|+\dots +|s^{-1}c_{i-1}|}\\
	&=& \sum_{1\leq i<n} s^{-1}c_1\otimes \dots \otimes s^{-1}c_i^{(1)}\otimes s^{-1}c_i^{(2)} \otimes \dots \otimes  s^{-1}c_n\ (-1)^{i + |c_1|+\dots +|c_{i-1}|+|c_i^{(1)}|}.
\end{eqnarray*}

\begin{lemma}
The derivations $d^{int}$ and $d^{ext}$ commute.
Moreover, $d^{int}$, $d^{ext}$ and $d^{int}\pm d^{ext}$ are derivations of square zero.
\end{lemma}
\begin{proof}
Let us prove the commutation of $d^{int}$ and $d^{ext}$. Because $d^{int}$ and $d^{ext}$ are of odd degrees, their graded commutator is $d^{ext}d^{int}+d^{int}d^{ext}$. It is again a derivation.
Thus, by corollary \ref{caracderfree}, $d^{ext}d^{int}+d^{int}d^{ext}=0$ if and only if $d^{int}b^{ext}+d^{ext}b^{int}=0$.
\begin{eqnarray*}
d^{int}b^{ext}(s^{-1}c) &=& d^{int}(-s^{-1}c^{(1)}\otimes s^{-1}c^{(2)}\ (-1)^{|c^{(1)}|})\\
	&=& s^{-1}d(c^{(1)})\otimes s^{-1}c^{(2)} \ (-1)^{|c^{(1)}|} - s^{-1}c^{(1)}\otimes \left(-s^{-1}d(c^{(2)})\right)\ (-1)^{(|c^{(1)}|-1)+|c^{(1)}|}\\
	&=& s^{-1}d(c^{(1)})\otimes s^{-1}c^{(2)} \ (-1)^{|c^{(1)}|} - s^{-1}c^{(1)}\otimes s^{-1}d(c^{(2)})
\end{eqnarray*}
and on the other side
\begin{eqnarray*}
d^{ext}b^{int}(s^{-1}c) &=& d^{ext}(-s^{-1}dc)\\
	&=& s^{-2}(dc)^{(1)}\otimes (dc)^{(2)}\\
	&=& s^{-2}\left(d(c^{(1)})\otimes c^{(2)} + c^{(1)}\otimes d(c^{(2)}) \ (-1)^{|c^{(1)}|}\right)\\
	&=& s^{-1}d(c^{(1)})\otimes s^{-1}c^{(2)} \ (-1)^{|c^{(1)}|-1} + s^{-1}c^{(1)}\otimes s^{-1}d(c^{(2)})
\end{eqnarray*}
hence $d^{int}b^{ext}+d^{ext}b^{int}=0$.

Let us prove now that $d^{int}$ and $d^{ext}$ are of square zero.
Using again corollary \ref{caracderfree}, $d^{int}d^{int}=0$ if and only if $d^{int}b^{int}=0$, but $d^{int}b^{int}(s^{-1}c) = d^{int}-s^{-1}dc = s^{-1}d^2c = 0$.
Similarly, $d^{ext}d^{ext}=0$ if and only if $d^{ext}b^{ext}=0$.
Using the coassociativity,
\begin{eqnarray*}
d^{ext}b^{ext} (s^{-1}c) &=& d^{ext}(-s^{-1}c^{(1)}\otimes s^{-1}c^{(2)}\ (-1)^{|c^{(1)}|})\\
	&=& -(-s^{-1}c^{(1)}\otimes s^{-1}c^{(2)})\otimes  s^{-1}c^{(3)} \ (-1)^{|c^{(1)}|+(|c^{(1)}|+|c^{(2)}|)}\\
	&& - s^{-1}c^{(1)}\otimes (-s^{-1}c^{(2)}\otimes s^{-1}c^{(3)})\ (-1)^{(|c^{(1)}|-1)+|c^{(2)}|+|c^{(1)}|}\\
	&=& s^{-3}c^{(1)}\otimes c^{(2)} \otimes  c^{(3)} - s^{-3}c^{(1)}\otimes c^{(2)}\otimes c^{(3)}\\
	&=&0.
\end{eqnarray*}

The last statement follows from $(d^{int}\pm d^{ext})^2 = (d^{int})^2 \pm (d^{int}d^{ext} + d^{ext}d^{int}) + (d^{ext})^2= (d^{int})^2 + (d^{ext})^2$ .
\end{proof}

We shall call $d^{int}$ the {\em internal differential} and $d^{ext}$ the {\em external differential}. 
Recall that the non-unital graded algebra $T_\circ(s^{-1}|C_-|)$ is equivalent to the pointed graded algebra $T_\bullet(s^{-1}|C_-|)$ and that derivations of $T_\circ(s^{-1}|C_-|)$ corresponds bijectively to pointed derivations of $T_\bullet(s^{-1}|C_-|)$. We shall again call $d^{int}$ and $d^{ext}$ the internal and external differentials of $T_\bullet(s^{-1}|C_-|)$ induced by that of $T_\circ(s^{-1}|C_-|)$.
Remark that the pointed dg-algebra $(T_\bullet(s^{-1}|C_-|),d^{int})$ is simply $T_\bullet(s^{-1}C_-)$.

\medskip 

The {\em cobar construction} of a pointed dg-coalgebra $C$ is noted $\Omega C$ and is defined as the pointed dg-algebra 
$$
(T_\bullet(s^{-1}|C_-|),d^{int} + d^{ext}).
$$
The cobar construction $\Omega$ defines a functor
$$\xymatrix{
\Omega:\dgCoalg_\bullet \ar[r]&\dgAlg_\bullet.
}$$

\begin{rem}\label{signcobar}
The sign in the definition of the total differential of $\Beta A$ depends on the authors, in \cite{Ad,EB}, the differential is $d^{int} - d^{ext}$, but in \cite{LV,Proute}, it is $d^{int} + d^{ext}$. 
Our convention is chosen to agree with our construction of $\Beta $ via Sweedler operations (see also section \ref{signtw} for a discussion of signs).
The following lemma prove that the two definitions of the bar construction are isomorphic.
\end{rem}

Recall the automorphism $\pi$ from after remark \ref{signbar}.

\begin{lemma}
The map $\pi:T_\bullet(s^{-1}|C_-|)\to T_\bullet(s^{-1}|C_-|)$ induces an isomorphism between the two cobar constructions
$$
(T_\bullet(s^{-1}|C_-|),d^{int} - d^{ext}) \et  (T_\bullet(s^{-1}|C_-|),d^{int} + d^{ext}).
$$
\end{lemma}
\begin{proof}
We want to prove $\pi^{-1}(d^{int} + d^{ext})\pi = d^{int} - d^{ext}$.
Let $X = s^{-1}|C_-|$, $d^{int}$ send $X^{\otimes n}$ to itself so $\pi^{-1}(d^{int})\pi = d^{int}$ and
$d^{ext}$ sends $X^{\otimes n}$ to $X^{\otimes n+1}$ so $\pi^{-1}(d^{ext})\pi = -d^{ext}$.
\end{proof}

\subsection{Universal twisting cochains}

Recall that a twisting cochain from a dg-coalgebra $C$ to a dg-algebra $A$ is an element $\alpha\in [C,A]$ of degree $-1$ which satisfies the Maurer-Cartan equation
$$
d_A\alpha+\alpha d_C + \alpha\star \alpha = 0.
$$
When $(C,e)$ and $(A,\epsilon)$ are pointed, a {\em pointed twisting cochain} is a twisting cochain such that $\epsilon \alpha = 0$ and $\alpha e=0$.
Let us denote $Tw_\bullet^c(C,A)$ the set of pointed twisting cochains from $C$ to $A$, this defines a functor
$$\xymatrix{
Tw_\bullet^c:(\dgCoalgnilbullet)^{op}\times \dgAlg_\bullet \ar[r]&\Set.
}$$

If $C$ is fixed, we shall say that a twisting cochain $C\to A$ is {\em universal} if it represents the functor $Tw_\bullet^c(C,-)$.
If $A$ is fixed, we shall say that a twisting cochain $C\to A$ is {\em couniversal} if it represents the functor $Tw_\bullet^c(-,A)$.

Let us consider the map
\begin{eqnarray*}
\beta :T_\circ^c(sA_-) &\underset{-1}{\lrhup} & A_-\\
sa_1\otimes \dots \otimes sa_n &\mto &
	\left\{\begin{array}{ll}
	\rule[-2ex]{0pt}{4ex} -a_1 & \textrm{if $n=1$}\\
	0 & \textrm{if $n\not =1$ is even}.
	\end{array}\right.
\end{eqnarray*}

For the justification of the sign see appendix \ref{signtw}.

\begin{lemma}
$\beta $ is a pointed twisting cochain.
\end{lemma}
\begin{proof}
Recall that for any coalgebra $C$, $[C,A]_\bullet = [C_-,A_-]_+$ and that a Maurer-Cartan element of $[C,A]_\bullet$ is pointed if and only if it is in the augmentation ideal of $[C,A]_\bullet$, \ie in $[C_-,A_-]$.
It is clear by definition that $\beta $ is equivalent to a map $T_\circ(s|A_-|)\to |A_-|$, hence $\beta $ defines an element in the kernel of the augmentation of $[T_\circ^c(sA_-),A_-]_+=[\Beta  A,A]_\bullet$. To prove that $\beta $ is a pointed Maurer-Cartan element, it remains to prove that it is a Maurer-Cartan element.
The Maurer-Cartan equation in $[T_\circ(sA_-),A_-]$ is 
$$
d_A\alpha + \alpha d_{BA} + \alpha\star \alpha =  d_A\alpha + \alpha d^{int} - \alpha d^{ext} + \alpha\star \alpha = 0.
$$
We are going to check that $\beta $ is a Maurer-Cartan element by evaluating this equation on elements $sa_1\otimes \dots \otimes sa_n$.
If $n=1$,
$$
\begin{array}{rcccl}
(d_A\beta  + \beta d^{int} -\beta d^{ext} + \beta \star \beta )(sa_1) &=& -d_A(a_1) + d_A(a_1) - 0 + 0 &=& 0\ ;
\end{array}
$$
if $n=2$,
$$
\begin{array}{rcccl}
(d_A\beta  + \beta d^{int} -\beta d^{ext} + \beta \star \beta )(sa_1\otimes sa_2) &=& 0 + 0 + a_1a_2 \ (-1)^{|a_1|} + a_1a_2\ (-1)^{|a_1|+1} &=& 0\ ;
\end{array}
$$
and, if $n\geq 3$,
\begin{eqnarray*}
(d_A\beta  + \beta d^{int} -\beta d^{ext} + \beta \star \beta )(sa_1\otimes \dots \otimes sa_n) &=& 0 + 0 + 0 + 0.
\end{eqnarray*}
\end{proof}

\begin{rem}\label{remmcbar}
The necessity of the sign in the formula of $\beta$ depends on our choice $d^{tot} = d^{int}-d^{ext}$ for the differential of $\Beta  A$ (the $\beta$ sign is necessary to have 0 in the case $n=2$). 
If the differential is chosen to be $d^{int}+d^{ext}$ the universal pointed twisting cochain would be $-\beta:sa \mapsto a$. This the choice of many authors, see section \ref{signtw}.
\end{rem}

\bigskip

Let us consider the map
\begin{eqnarray*}
\omega :C_- &\underset{-1}{\lrhup} & T_\circ(s^{-1}C_-)\\
c &\mto & s^{-1}c
\end{eqnarray*}

\begin{lemma}
$\omega $ is a pointed twisting cochain.
\end{lemma}
\begin{proof}
It is clear by definition that $\omega $ is equivalent to a map $|C_-|\to T_\circ(s^{-1}|C_-|)$, hence $\beta $ defines an element in the kernel of the augmentation of $[C_-,T_\circ(s^{-1}|C_-|)]_+=[C,\Omega C]_\bullet$. To prove that $\omega $ is a pointed Maurer-Cartan element it is sufficient to prove that it is a Maurer-Cartan element.
The Maurer-Cartan equation in $[C_-,T_\circ(s^{-1}|C_-|)]$ is 
$$
d_{\Omega C}\alpha + \alpha d_C + \alpha\star \alpha =  d^{int}\alpha +d^{ext}\alpha +\alpha d_C + \alpha\star \alpha = 0.
$$
We are going to check that $\omega $ is a Maurer-Cartan element by evaluating this equation on an element $c\in C$
\begin{eqnarray*}
(d^{int}\omega  +d^{ext}\omega  +\omega d_C + \omega \star \omega )(c) &=& -s^{-1}d_Cc - s^{-1}c^{(1)}\otimes s^{-1}c^{(1)}\ (-1)^{|c^{(1)}|} \\
&& + s^{-1}d_Cc + s^{-1}c^{(1)}\otimes s^{-1}c^{(1)}\ (-1)^{|c^{(1)}|}\\
&=& 0.
\end{eqnarray*}
\end{proof}

\begin{rem}\label{remmccobar}
Again, the choice of the twisting cochain $\omega$ is conditionned by that of the differential $d_{\Omega C}=d^{int}+d^{ext}$ for of $\Omega C$.
If the differential is chosen to be $d^{int}-d^{ext}$ the universal pointed twisting cochain would be $-\omega:c \mapsto -s^{-1}c$.
\end{rem}

\medskip

\begin{thm}\label{classicbarcobar}
$\beta $ is a universal pointed twisting cochain and $\omega $ is a couniversal twisting cochain. 
In other words, the functors 
$$\xymatrix{
\Omega:\dgCoalgnilbullet \ar[r]&\dgAlg_\bullet
}\et
\xymatrix{
\Beta :\dgAlg_\bullet \ar[r]&\dgCoalgnilbullet
}$$
are adjoint and represent the binary relator $Tw_\bullet^c$:
$$
\dgAlg_\bullet(\Omega C,A) = Tw_\bullet^c(C,A) = \dgCoalgnilbullet(C,\Beta  A).
$$
\end{thm}
\begin{proof}
The proof is taken from \cite[thm 2.2.9]{LV} (see also \cite[thm 2.10 \& 2.11]{Proute}).

Because $\Omega C$ is free as a graded algebra, a pointed graded algebra map $f:\Omega C\to A$ is equivalent to a map $\phi:s^{-1}C_-\to A_-$, and to a degree $-1$ morphism $\alpha:C_-\rhup_{-1} A_-$ such that $\alpha(c) = \phi (s^{-1}c)$.
The commutation with the differential is equivalent to 
\begin{eqnarray*}
\rule[-2ex]{0pt}{4ex} &&f(d^{int}(c)) + d^{ext}(c)) = d_A(f(c)) \\
\rule[-2ex]{0pt}{4ex} & \iff & \phi (-s^{-1}d_Cc) + m_A(\phi\otimes \phi) (-s^{-1}c^{(1)}\otimes s^{-1}c^{(2)})\ (-1)^{|c^{(1)}|}= d_A(\phi(c))\\
\rule[-2ex]{0pt}{4ex} & \iff & - \alpha (d_Cc) - \alpha(c^{(1)})\alpha(c^{(1)}) \ (-1)^{|c^{(1)}|}= d_A(\alpha(c))\\
\rule[-2ex]{0pt}{4ex} & \iff & d_A(\alpha(c)) + \alpha (d_Cc) + \alpha(c^{(1)}) \alpha(c^{(1)}) \ (-1)^{|c^{(1)}|}= 0
\end{eqnarray*}
which is the Maurer-Cartan equation. This produces the first bijection.
The second one is analog but there is a sign issue.

Because $\Beta  A$ is cofree as a graded conilpotent coalgebra, a pointed graded coalgebra map $f:C\to \Beta  A$ is equivalent to a map $\phi:C_-\to sA_-$.
Then by proposition \ref{structureconil} $f(c) = \sum_n\otimes_{i=1}^n\phi(c^{(i)})$. 
Let us write $\phi=s\alpha$ for some degree $-1$ morphism $\alpha:C_-\rhup_{-1} A_-$.
The commutation with the differential  $(d^{int} - d^{ext})(f(c)) = f(d_Cc) $ is equivalent to $p(d^{int} - d^{ext})(f(c)) = \phi (d_Cc)$ where $p$ is the projection on cogenerators.
Using the formula $f(c) = \sum_n\otimes_{i=1}^n\phi(c^{(i)})$, most of the terms disappear and this is equivalent to
\begin{eqnarray*}
\rule[-2ex]{0pt}{4ex} &&d^{int}(\phi(c)) - pd^{ext}(\phi(c^{(1)})\phi(c^{(2)})) =\phi(d_Cc)\\
\rule[-2ex]{0pt}{4ex} &\iff& d^{int}(s\alpha(c)) - pd^{ext}(s\alpha(c^{(1)})s\alpha(c^{(2)})) =s\alpha(d_Cc) \\
\rule[-2ex]{0pt}{4ex} &\iff& -sd_A(s\alpha(c)) - s(\alpha(c^{(1)})\alpha(c^{(2)}))\ (-1)^{|\alpha(c^{(1)})|} = s\alpha(d_Cc) \\
\rule[-2ex]{0pt}{4ex} &\iff& d_A(\alpha(c)) + \alpha(d_Cc) - \alpha(c^{(1)})\alpha(c^{(2)})\ (-1)^{|c^{(1)}|} = 0
\end{eqnarray*}
where we have used $|\alpha(c^{(1)})| = |c^{(1)}|-1$.
This is not the Maurer-Cartan equation because of the minus sign. 
However $-\alpha$ does satisfies the Maurer-Cartan equation and the second bijection is $f \mapsto -\alpha$.
\end{proof}

\begin{rem}
Note that we would not have had the sign issue in the previous proof if we had chosen $d^{int}+d^{ext}$ for the differential of $\Beta A$.
\end{rem}

\subsection{Signs issues}\label{signtw}

As already noted in previous remarks, different authors have had different conventions of signs for the external differential in the total differential of the bar and cobar constructions as well as for the universal twisting cochains. 
We list here the conventions chosen in our main references according to {\sl our} definitions of $d^{ext}$ (an X means that the work does not consider this construction).

\begin{center}
\begin{tabular}{|c|cc|cc|}
\hline
\rule[-2ex]{0pt}{5ex} authors $\backslash$ sign of & $\quad d_{\Betaext A}=d^{int}\pm d^{ext}\quad$ & $\quad \beta\quad  $ & $\quad  d_{\Omega C}=d^{int}\pm d^{ext}\quad $ & $\quad \omega\quad $\\
\hline
\rule[-2ex]{0pt}{5ex} us & $-$ & $-$ & $+$ & $+$ \\
\rule[-2ex]{0pt}{5ex} \cite{EM} & $-$ & X & X & X\\
\rule[-2ex]{0pt}{5ex} \cite{Ad2} & X & X & $-$ & X\\
\rule[-2ex]{0pt}{5ex} \cite{EB} & $+$ & X & $-$ & X\\
\rule[-2ex]{0pt}{5ex} \cite{Proute} & $+$ & $+$ & $+$ & $+$\\
\rule[-2ex]{0pt}{5ex} \cite{LV}  & $+$ & $+$ & $+$ & $+$\\
\hline
\end{tabular}
\end{center}

In \cite{EM}, the construction of the external differential is based on the alternate sum of faces operators of (co)simplicial sets, this approach produces a natural minus sign for the external differential. 
The external differential of the bar construction in \cite{EB} and the cobar construction in \cite{Ad2, EB} are defined with no signs in front of the product or the coproduct, so they correspond to our $-d^{ext}$.
These articles do not compute the universal twisting cochains and the use of no sign gives the simplests formulas.
The notion of universal twisting cochains seems to have appear explicitely in \cite{Proute} where the author actually starts by fixing the universal twisting cochains and then deduces the formula for differential from the Maurer-Cartan equation. 
In \cite{LV}, the authors follows the conventions of \cite{Proute} but introduce some duality between the definitions of the bar and cobar constructions (the maps $\prod_s$ and $\Delta_s$ in their notations).

\medskip

All these conventions can be justified by different reasons (mostly simplicity of notations) but signs cannot totally disappear. The explanation is that there are {\sl three} different sources of minus signs in the formulas; some of them can compensate each other but we have found simpler to make them all of them explicit in our definitions:
\begin{enumerate}
\item the sign in $d_\mc u=-u^2$ is responsible for the sign in our definition of $d^{ext}_{\Omega C}$,
\item the contravariance of $\{\mc,-\}_\bullet$ in $\mc$ is responsible for the sign in $d^{tot}_{\Beta  A} = d^{int} - d^{ext}$ (see remark \ref{remdmcdextbar}),
\item and the lambda transform $\lambda^1$ is responsible for the sign in the universal twisting cochain $\beta$.
\end{enumerate}
A few remarks can be added to justify further our conventions.
\begin{enumerate}
\setcounter{enumi}{3}
\item The approach of the bar and cobar constructions using the universal Maurer-Cartan element $u$ is without arbitrary choice, not like the approach using the symbol $s$, \ie the {\sl bar notation}, where different "natural" conventions exists.
\item The relation between $s$ and $u$ exists but is not canonical. Using our definitions of the two $d^{ext}$ and the identifications $s^{-1}=u, s=u^\star$ (cf. the remark about the notations $s$ and $s^{-1}$ in section \ref{suspension}), we have the nice correspondance that $d^{ext}$ is exactly the action of the Maurer-Cartan differential.

\item Our conventions are in accordance with Loday and Vallette's philosophy of duality between $\Beta $ and $\Omega$. First, it should be remarked that our definitions of the two $d^{ext}$ coincides with their definitions of the two $d_2$ \cite[2.2.1, 2.2.5]{LV}, but our definitions of the total differentials coincide only for the cobar construction.
Then, to expand on the point 1. above, with $s^{-1}= u$, the Loday-Vallette map $\Delta_s:s^{-1}\mapsto -s^{-1}\otimes s^{-1}$ is secretly the Maurer-Cartan differential $u\mapsto -u\otimes u$, hence the sign in the definition of $d^{ext}_{\Omega C}$.
Dually, with $s = u^\star$, the Loday-Vallette map $\prod_s:s\otimes s\mapsto s$ is the dual of the Maurer-Cartan differential $u^\star \otimes u^\star \mapsto u^\star$, hence the definition of $d^{ext}_{\Beta  A}$ without sign.
These correspondances are explicit when computing the derivations $C\rhd_\bullet d_\mc$ and $\{d_\mc,A\}_\bullet$ using meta-morphisms
(see the computation after theorem \ref{represtw}).
\end{enumerate}

\newpage
\section{Appendix on category theory}\label{AppendixCattheory}

\subsection{Relators and Gray trialities}\label{triality}

In this appendix, we define the structure of Gray triality which is underlying Sweedler's operations. 
This structure generalizes the notion of an adjunction between two categories when more than two categories are involved.
Similar notions have been considered in \cite{Gray,CGR} but our approach articulate the situation around a functor that we call a {\em relator}.
The notion of triality is also a categorification of the notion of triality for vector spaces \cite{Ad}, this was our motivation for the name.

\bigskip

Let $\bC_1$, ..., $\bC_n$ be some categories, we shall call a {\em relator between $\bC_1$, ..., $\bC_n$} (or simply an {\em $n$-ary relator}) any functor
$$\xymatrix{
R:\bC_1\times \dots \times \bC_n\ar[r]& \Set
}$$
For objects $C_i\in \bC_i$, we shall say that an element in $R(C_1,\dots , C_n)$ is a {\em relation between $C_1$, ..., $C_n$}.

\medskip

For $1\leq i\leq n$, we shall say that a relator is {\em representable in variable $i$} if the functor
$$\xymatrix{
R_i(C_1,\dots, C_{i-1},C_{i+1},\dots, C_n) := R(C_1,\dots, C_{i-1}, - , C_{i+1},\dots, C_n) :\bC_i\ar[r] &\Set
}$$
is representable for any family $(C_1,\dots, C_{i-1},C_{i+1},\dots, C_n)\in \prod_{j\not= i}C_j$.

\medskip

We shall say that a relator is {\em representable} if it is representable in all its variables.

\paragraph{Unary relators}

A unary relator is just a functor $R:\bC\to \Set$; it is representable if it is of the type $R=\bC(A,-)$ for some object $A$.

\paragraph{Binary relators}\label{2rel}

A binary relator is a functor $R:\bC_1\times \bC_2 \to \Set$; it is representable if there exists two functors 
$$\xymatrix{
u:\bC_1^{op}\ar[r]& \bC_2
}
\et
\xymatrix{
v:\bC_2\ar[r]& \bC_1^{op}
}$$
together with natural isomorphisms
$$
\bC_2(u(A),B) = R(A,B) = \bC_1^{op}(A,v(B)) = \bC_1(v(B),A)
$$
In other words, the data of a representable relator is equivalent to that of an adjunction
$$\xymatrix{
u:\bC_1^{op}\ar@<.6ex>[r]& \bC_2:v\ar@<.6ex>[l]
}$$

In this case, we shall say that the relator $\bC_1^{op}\times \bC_2\to \Set$ is representable by the adjunction $u\dashv v$.

\paragraph{Ternary relators}\label{3rel}\label{graytrialitydef}

A ternary relator is a functor $R:\bC_1\times \bC_2\times \bC_3 \to \Set$; it is representable if there exists three functors 
$$\xymatrix{
u_1:\bC_2\times \bC_3\ar[r] &\bC_1^{op}
}$$
$$\xymatrix{
u_2:\bC_1\times \bC_3\ar[r] &\bC_2^{op}
}$$
$$\xymatrix{
u_3:\bC_1\times \bC_2\ar[r] &\bC_3^{op}
}$$
together with natural isomorphisms
$$
R(A,B,C) = \bC_1(u_1(B,C),A) = \bC_2(u_2(A,C),B) = \bC_3(u_3(A,B),C) 
$$
In particular, every objet $A\in \bC_1$ defines an adjunction
$$\xymatrix{
u_3(A,-): \bC_2 \ar@<.6ex>[r] &\ar@<.6ex>[l] \bC_3^{op} : u_2(A,-)
}$$
and similarly for objects in $\bC_2$ and $\bC_3$.

\medskip
We shall call a representable ternary relator a {\em Gray triality}.

\medskip
In other words, a Gray triality can be thought as the data of a certain game of adjunctions between three categories.
This can be illustrated by a diagram
$$
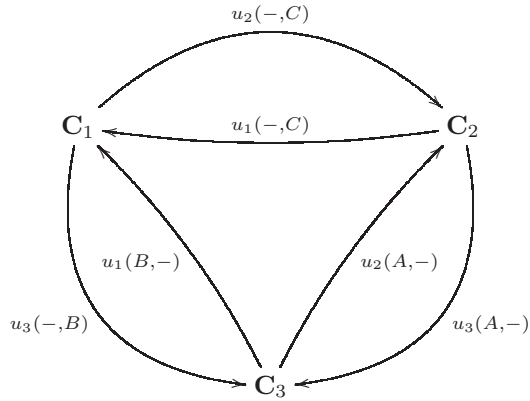
\xymatrix{
\bC_1   \ar@/^3pc/[rrrr]^-{u_2(-,C)}\ar@/_3pc/[dddrr]_-{u_3(-,B)}  &&&&\ar@/^0,5pc/[llll]_-{u_1(-,C)}   \bC_2  \ar@/^3pc/[dddll]^-{u_3(A,-)}\\
&&&& \\
&&&& \\
&& \bC_3 \ar@/_0,5pc/[uuull]^-{u_1(B,-)} \ar@/^0,5pc/[uuurr]_-{u_2(A,-)} && 
}.$$

\begin{ex} If ${\cal P}(X)$ denote the power set of a set $X$, the binary relator
$$\xymatrix{
R:\Set^{op}\times \Set^{op}\ar[r]& \Set
}$$
defined by $R(E,F) = {\cal P}(E\times F)$ is representable. 
Indeed, for any two sets $E$ and $F$, the functors
$$\xymatrix{
R(E,-):\Set^{op}\ar[r]& \Set
}\et
\xymatrix{
R(-,F):\Set^{op}\ar[r]& \Set
}$$
are respectively representable by ${\cal P}(E)$ and  ${\cal P}(F)$.
\end{ex}

\begin{ex}
If $\Vect$ is the category of (non-graded) vector spaces, for every triple of objects $(A,B,C)\in \Vect$ let $R(A,B,C)$ be the set of maps $A\times B\to C$ that are bilinear.
This defines a ternary relator
$$\xymatrix{
R:\Vect^{op}\times \Vect^{op}\times \Vect \ar[r]& \Set.
}$$
The classical definition of the tensor product $A\otimes B$ of two vector spaces $A$ and $B$ can be seen to be the representation property of this relator in its third variable.
This relator is also representable in the first variable by the vector spaces of linear map $[A,C]$ and in its second variable by $[B,C]$.

\end{ex}

\begin{ex}
If $\bV=(\bV,\otimes )$ is a monoidal category, let us put $R(A,B,C)=\bV(A\otimes B,C)$ for every triple of objects $(A,B,C)\in \bV^3$.
This defines a ternary relator
$$\xymatrix{
R:\bV^{op}\times \bV^{op}\times \bV \ar[r]& \Set.
}$$
By definition, this relator is representable in its third variable.
It is representable in the second variable if and only if the functor $A\otimes (-): \bV\to\bV$ has a right adjoint for every object $A\in \bV$, that is, if and only if the functor $\otimes$ is divisible on the left. Dually, the relator representable in the first variable if and only if the functor $(-)\otimes A: \bV\to \bV$ has a right adjoint for every object $A\in  \bV$, that is, if and only if the functor $\otimes$ is divisible on the right(\footnote{The fact that divisible monoidal categories produce Gray trialities is analog to the relation between division algebras and triality of vector spaces.}).

\medskip

If $(\bV,\otimes)$ is symmetric monoidal, $R$ is representable, \ie is a Gray triality, if and only if $(\bV,\otimes)$ is closed.
If $[A,B]$ is the internal hom in $\bV$, we can picture the triality as
$$
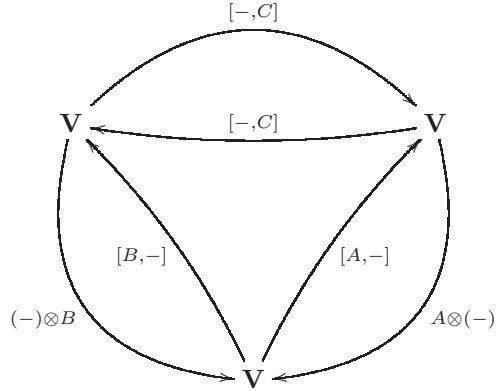
\xymatrix{
\bV \ar@/^3pc/[rrrr]^-{[-,C]}\ar@/_3pc/[dddrr]_-{(-)\otimes B}	&&&&\bV \ar@/^0.5pc/[llll]_-{[-,C]} \ar@/^3pc/[dddll]^-{A\otimes(-)}\\
&&&& \\
&&&& \\
&&	\bV \ar@/_0.5pc/[uuull]^-{[B,-]} \ar@/^0.5pc/[uuurr]_-{[A,-]} && 
}. $$

\end{ex}

\begin{ex}
If $\bV=(\bV,\otimes )$ is a monoidal category, we shall call category $\bC$ {\em tensored over $\bV$} if it is equipped with a functor $\rhd :\bV\times \bC\to \bC$ which defines an action of the monoid $\bV$ on $\bC$.
Let us define a ternary relator 
$$\xymatrix{
R:\bV^{op}\times \bC^{op}\times \bC \ar[r]& \Set.
}$$
by 
$$
R(A,X,Y)=\bC(A\rhd X,Y)
$$
for every triple of objects $(A,X,Y)\in \bV\times \bC\times \bC$.

This relator is by definition representable in its third position.
It is representable in its first position, if and only if $\bC$ is {\em enriched over $\bV$}.
And it is representable in its second position, if and only if $\bC$ is {\em cotensored over $\bV$}.

Thus $R$ is Gray triality if and only if and only if $\bC$ is enriched, tensored and cotensored over $\bV$.

\medskip

For such a category $\bC$, let $\{X,Y\}$ be the notation for enrichment of $\bC$ over $\bV$ and $[A,X]$ be the notation for the cotensor,
we can picture the triality as
$$
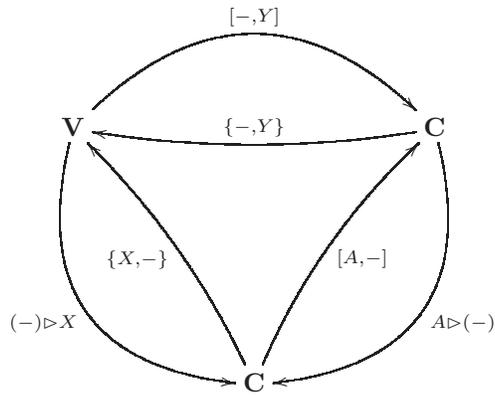
\xymatrix{
{\bV} \ar@/^3pc/[rrrr]^-{[-,Y]}\ar@/_3pc/[dddrr]_-{(-)\rhd X}	&&&&{\bC} \ar@/^0.5pc/[llll]_-{\{-,Y\}} \ar@/^3pc/[dddll]^-{A\rhd (-)}\\
&&&& \\
&&&& \\
&&	{\bC} \ar@/_0.5pc/[uuull]^-{\{X,-\}} \ar@/^0.6pc/[uuurr]_-{[A,-]} && 
} $$

\end{ex}

\subsection{Monoidal categories and functors}\label{Catmonoidales}

Recall \cite{MacLane} that a {\bf monoidal structure} 
on a category  $\bC$ is five-tuple
$(\otimes, \alpha, I, \lambda, \rho)$ consisting of

\begin{itemize}

\item{} a functor $\otimes : \bC\times \bC\to \bC$
called the {\it tensor product functor}
\item{} a natural isomorphism 
$$\alpha(A,B,C):(A\otimes B)\otimes C\to A\otimes (B\otimes C)$$
called the {\it associativity constraint}

\item{} an object $I\in \bC$ called the  {\it unit object}.

\item{}  two natural isomorphisms 
$$\lambda(A):I\otimes A\to A \quad \quad \rho(A):A\otimes I\to A,$$
 respectively called the {\it left and right unit contraints}.
\end{itemize}

\noindent
The fivetuple $(\otimes, \mathbf{a}, I, \lambda, \rho)$ must satisfy the following conditions:

\begin{itemize}

\item{} {\bf Coherence of the associativity constraint}

The following pentagon commutes for every quadruple of objects $(A,B,C,D)$,
where $AB=A\otimes B$.
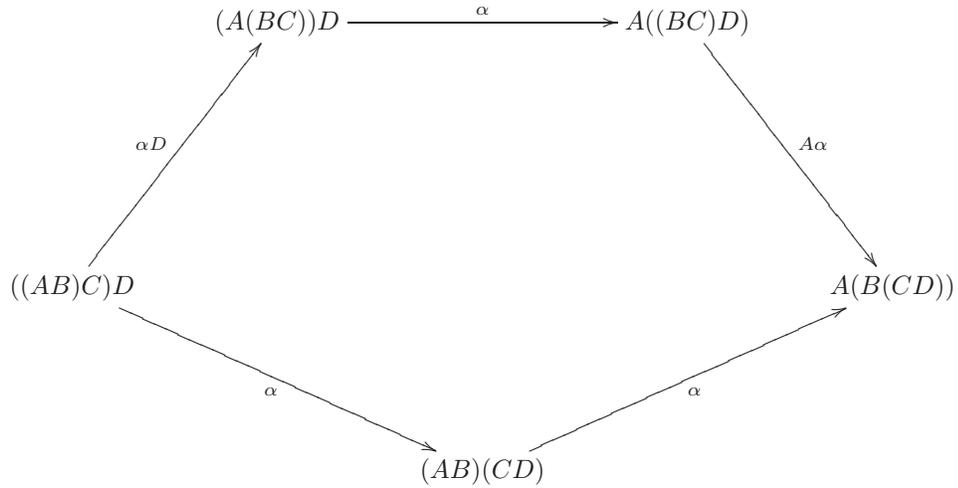
$$\xymatrix{
 & (A(BC))D  \ar[rr]^-{\alpha}  && A((BC)D)  \ar[dddr]^-{A\alpha} & \\
&& &&\\
&& &&\\
((AB)C)D \ar[ruuu]^{\alpha D}\ar[rrdd]_{\alpha} && &&   A(B(CD)) \\
&& &&\\
&& (AB)(CD) \ar[rruu]_{\alpha}  &&
}$$

\item{} {\bf Coherence of the unit constraints} The following triangle commutes for
every couple of objects $(A,B)$,
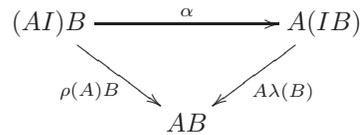
$$\xymatrix{
 (AI) B \ar[rr]^-{\alpha}  \ar[rd]_{\rho(A)B}  && A(I B)  \ar[dl]^-{A\lambda(B)} \\
&AB &
}$$
\end{itemize}

A {\bf monoidal category} is a category $ \bC$ equipped
with a monoidal structure $(\otimes, \mathbf{a}, I, \lambda, \rho)$.
We shall often denote a monoidal category $( \bC, \otimes, \mathbf{a}, I, \lambda, \rho)$
by $( \bC, \otimes,I)$.

\medskip

The coherence theorem of Mac Lane \cite{MacLane} asserts that the product of a sequence  
of objects $(A_1,\ldots ,A_n)$ in a monoidal category is unambiguously defined
up to a unique canonical isomorphism,  that is, any two products,
obtained by bracketing the objects in the sequence are canonically isomorphic.
We shall often denote a canonical isomorphism between by an equality sign.
For example, the associativity and unit constraints:
$$\xymatrix{
(A\otimes B)\otimes C \ar@{=}[r]& A\otimes (B\otimes C)
}\quad \quad
\xymatrix{
I\otimes A  \ar@{=}[r]& A  \ar@{=}[r] & A\otimes I
}$$

\medskip

Recall that if $(\bC, \otimes, I)$ is a monoidal category 
then a natural isomorphism
$$\sigma=\sigma(A,B):A\otimes B \to B\otimes A$$ 
is said to be a {\it braiding} if the following diagrams commute,
where $AB:=A\otimes B$,
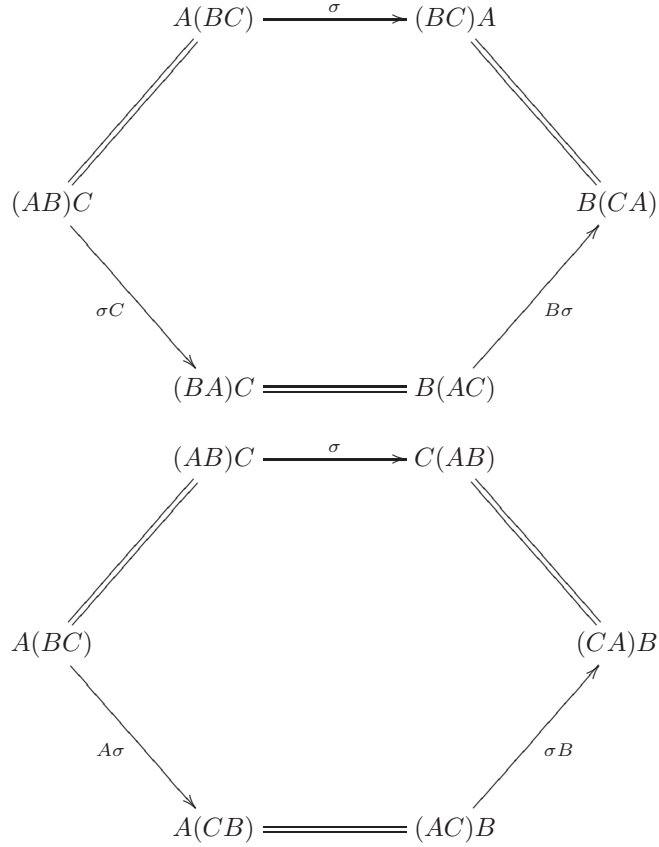
$$\xymatrix{
 &  \ar@{=}[ddl]    A(BC)     \ar[rr]^-{\sigma}  &&  (BC)A \ar@{=}[ddr]  &  \\
  && &  \\
  \ar[ddr]_-{\sigma C} (AB)C  &  & &&  B(CA) \\
  & & & \\
 &  (BA)C  \ar@{=}[rr] && B(AC) \ar[uur]_-{B \sigma} &  
}$$
$$\xymatrix{
 &  \ar@{=}[ddl]    (AB)C    \ar[rr]^-{\sigma}  &&  C(AB) \ar@{=}[ddr]  &  \\
  && &  \\
  \ar[ddr]_-{A\sigma} A(BC)  &  & &&  (CA)B \\
  & & & \\
 &  A(CB) \ar@{=}[rr] && (AC)B \ar[uur]_-{\sigma B} &  
}$$
A {\it braided monoidal category} is a monoidal
category $(\bC, \otimes,I)$ equipped with a braiding $\sigma$.
The braiding $\sigma$ is said to be a {\it symmetry} 
if the composite 
$$\xymatrix{A\otimes B \ar[r]^-{\sigma(A,B)} & B\otimes A \ar[r]^-{\sigma(B,A)} & A\otimes B
}$$
is the identity,
A {\it symmetric monoidal category} is a monoidal
category $(\bC, \otimes,I)$ equipped with a symmetry $\sigma$.

\medskip

A {\it monoidal functor} $F:\bC\to \bD$
between two monoidal categories
is a lax functor $F=(F,\alpha,\alpha_0)$ in which the product $\alpha$ and
the unit $\alpha_0$ are invertible. It can also be defined
as colax functor $F=(F,\alpha,\alpha_0)$ in which the coproduct $\alpha$ and
the counit $\alpha_0$ are invertible. 
A {\it homomorphism} $\theta:F\to G$ of monoidal functors $F,G: \bC\to \bD$
 is a homomorphism of lax functors. 
The composite of two 
 monoidal functors $F:\bC\to \bD$
and $G:\bD\to \bE$ has the structure of a monoidal functor
$G\circ F:\bC\to \bE$. The monoidal categories, monoidal functors 
and homomorphisms
form a 2-category.

\medskip

Similarly, 
a {\it symmetric monoidal functor} $F:\bC\to \bD$
between two symmetric monoidal categories
is a symmetric lax functor $F=(F,\alpha,\alpha_0)$ in which the product $\alpha$ and
the unit $\alpha_0$ are invertible.
A {\it homomorphism} $\theta:F\to G$ of symmetric monoidal functors $F,G: \bC\to \bD$
 is a homomorphism of lax functors. 
 The composite of two 
symmetric monoidal functors $F:\bC\to \bD$
and $G:\bD\to \bE$ has the structure of a symmetric monoidal functor
$G\circ F:\bC\to \bE$. The symmetric monoidal categories, symmetric monoidal functors 
and homomorphisms
form a 2-category.

\bigskip

\subsection{Closed categories}\label{appclosedcategories}

Let $\bV=( \bV ,\otimes, I)$ be a monoidal category
We shall
say that a pair of morphisms 
$\eta: I \to B\otimes A$ and $\epsilon:A \otimes B\to I$ 
is an {\bf adjunction} between two objects $A$ and $B$ if the
following triangles commute,
$$\xymatrix{
B \ar[r]^-{\eta \otimes B}  \ar@{=}[dr] & B\otimes A \otimes B \ar[d]^-{B\otimes \epsilon} \\
 & B
}\quad \quad
\xymatrix{
A \ar[r]^-{A\otimes  \eta }  \ar@{=}[dr] &A \otimes B \otimes A \ar[d]^-{ \epsilon \otimes A} \\
 &X
}
$$
In which case we shall write $(\eta,\epsilon):A\dashv B$.
The object $A$ is said to be the {\bf left adjoint} and the object $B$ the {\bf right adjoint}.
The morphism $\eta$ is the {\bf unit} of the adjunction and 
the morphism $\epsilon$ the {\bf counit}.
The functor $A\otimes(-) :\bV \to \bV$
is then left adjoint to the functor $B\otimes (-)$, 
and the functor $(-)\otimes B$ left adjoint to the
functor $(-)\otimes A $.
The right (left) adjoint of an object is unique up to unique isomorphism.
When the category $\bV$ is symmetric monoidal,
an adjunction $(\eta,\epsilon):A\dashv B$
is equivalent to an adjunction $(\sigma \eta, \epsilon \sigma):B\dashv A$,
where $\sigma:B\otimes A\to A\otimes B$ is the symmetry.
We shall say that an object $A$ of a symmetric monoidal
category is {\bf perfect} if it admits it admits a right (resp. right) adjoint.

\medskip

We shall say that a monoidal category $\bV=( \bV ,\otimes, I)$ is {\it left closed}
(resp. {\it right closed}, resp. {\it closed})
if the tensor product functor $\otimes :\bV \times\bV \to \bV$
is divisible on the left (resp. on the right, resp. on both sides).

\medskip

Let $\bV=( \bV ,\otimes, I)$ be a right closed monoidal category.
By definition, the functor $(-)\otimes X: \bV\to  \bV$
has a right adjoint $(-)/X: \bV\to  \bV$
for every object $X\in \bV$. We shall denote the
the counit of the adjunction $(-)\otimes X\dashv (-)/X$ 
as an evalutaion map $ev=ev(X,Y):Y/X\otimes X\to Y$.
 By definition,  for every object $A\in  \bV $
 and every map $f:A\otimes X\to Y$, there exist
 a unique map $g:A\to Y/X$
 such that $ev(g \otimes X)=f$. 
We shall write $\lambda^2(f):=g$.
  In particular, for any triple of objects $X,Y,Z\in \bV $, there exists a unique map 
 $c=c(Z,Y,X): Z/Y\otimes Y/X\to Z/X$ such that 
 the following square commutes
 $$
\xymatrix{
Z/Y\otimes Y/X \otimes X    \ar[rr]^-{c \otimes X} \ar[d]_{Z/Y\otimes ev} && Z/X\otimes X  \ar[d]^{ev}  \\
 Z/Y \otimes Y  \ar[rr]^-{ev}&& Z
}
$$
The map $c(Z,Y,X)$ is the composition law of an enrichment
 of the monoidal category $ \bV $ over itself if we put
 $Hom(X,Y)=Y/X$ for $X,Y\in \bV $.

\medskip

Let $\bV=( \bV ,\otimes, I)$ be a left closed monoidal category.
By definition, the functor $X\otimes (-): \bV\to  \bV$
has a right adjoint $X\backslash(-): \bV\to  \bV$
for every object $X\in \bV$. We shall denote the
the counit of the adjunction $X\otimes (-)\dashv X\backslash(-)$ 
as a reversed evalutaion map $rev=rev(X,Y):X \otimes X\backslash Y\to Y$.
 By definition,  for every object $A\in  \bV $
 and every map $f:X\otimes A\to Y$, there exist
 a unique map $g:A\to X\backslash Y$
 such that $rev(X \otimes g)=f$. 
We shall write $\lambda^1(f):=g$.
  In particular, for any triple of objects $X,Y,Z\in \bV $, there exists a unique map 
 $rc=rc(X,Y,Z): X\backslash Y\otimes Y\backslash Z \to X\backslash Z$ such that 
 the following square commutes
 $$
\xymatrix{
X\otimes X\backslash Y\otimes Y\backslash Z   \ar[rr]^-{X \otimes rc} \ar[d]_{rev \otimes (Y\backslash Z) } &&  X\otimes X\backslash Z  \ar[d]^{rev}  \\
Y \otimes Y\backslash Z  \ar[rr]^-{rev}&& Z
}
$$
The map $rc(X,Y,Z)$ is the composition law of an enrichment
of the reversed monoidal category $ \bV^{rev} $ over itself if we put $Hom(X,Y)=X\backslash Y$ for $X,Y\in \bV $.

\medskip
 
A symmeric monoidal category is  left closed if and only if it is right closed
if and only if it is closed. 
Let $\bE=(\bE,\otimes, I, \sigma) $ be a symmetric monoidal closed category
and let us put $Hom(X,Y)=Y/X$. The  reversed evaluation $rev:X\otimes Hom(X,Y)  \to Y$
is defined to be the composite $ev \circ \sigma :X\otimes Hom(X,Y)  \to Hom(X,Y) \otimes X \to Y$.
For every object $A\in  \bV $
 and every map $f:X\otimes A\to Y$, there exist
 a unique map $g=\lambda^1(f):A\to Hom(X,Y)$
 such that $rev(X \otimes g)=f$. 
 By construction, $\lambda^1(f)=\lambda^2(f\sigma)$, where $f\sigma:A\otimes X\to Y$.
 More generally, to every morphism $f:A_1\otimes \cdots \otimes A_n \to B$ 
we can associate a morphism
$${\lambda}^k(f):A_1\otimes \cdots \otimes \Hat{A_{k}} \otimes \cdots \otimes A_n \to Hom(A_k,B)$$
 for each $1\leq k\leq n$. We can also associate a morphism
 $${\lambda}^{k,r}(f)= A_1\otimes \cdots \otimes \Hat{A_{k}} \otimes \cdots  \otimes \Hat{A_{r}} \otimes \cdots  \otimes A_n \to 
 Hom(A_k\otimes A_r,B)$$
for each $1\leq k<r\leq n$. Etc.

\medskip

The tensor product  functor $\otimes:\bE\times \bE\to \bE$
of a symmetric monoidal closed category $\bE$ is strong. 
Its strength is defined by the canonical map
$$\psi=\lambda^{3,4}(\alpha) :Hom(A_1,A_2)\otimes Hom(B_1,B_2) \to  Hom(A_1\otimes B_1, A_2\otimes B_2),$$
where $\alpha$ is the composite
$$\xymatrix{
Hom(A_1,A_2)\otimes Hom(B_1,B_2)\otimes A_1\otimes B_1\ar@/_1pc/[rrrd]_\alpha
\ar[rrr]^-{id\otimes \sigma\otimes id}  &&& Hom(A_1,A_2)\otimes A_1\otimes  Hom(B_1,B_2)\otimes B_1 \ar[d]^{ev\otimes ev}   \\
&&&  A_2\otimes B_2.
}$$

\medskip

Similarly, the hom functor $Hom: \bE^{op}\times \bE \to \bE$
of a symmetric monoidal closed category $\bE$ is strong. 
Its strength is given by the canonical map
$$\psi={\lambda}^2(c^3)\sigma: Hom(A_2,A_1) \otimes  Hom(B_1,B_2) \to 
Hom(Hom(A_1,B_1),  Hom(A_2,B_2) )$$
where $c^3$ is the three-fold  composition law,
$$\xymatrix{
Hom(B_1,B_2)\otimes  Hom(A_1,B_1) \otimes Hom(A_2,A_1) \ar[rr] &&Hom(A_2,B_2) .
}$$

\medskip

For a symmetric monoidal closed category, the coherence constraint can be equivalently reformulated as condition on the internal hom, the explicit diagrams are the same as those defining an opmodule in appendix \ref{Vopmodule}.

\subsection{Lax and colax functors}\label{Structureslax,colaxetbilax}

Here we recall the notions of lax and  colax functors  \cite{AM}.

\begin{defi}\label{defstructurelax}

If $\bC=(\bC, \star, J)$ and  $\bD=(\bD, \diamond, I)$
are two monoidal categories, then a \emph{lax structure} on a functor
$F:\bC\to \bD$ is a pair $(\alpha, \alpha_0)$, where $\alpha$ is a natural
transformation
$$\alpha(A,B): FA \diamond FB \to F(A\star B)$$
and $\alpha_0$ is a  morphisme $:I\to FJ$ 
satisfying the following conditions,

\medskip

\noindent
{\bf associativity}: the following square commutes:
$$\xymatrix{
FA\diamond FB\diamond FC \ar[rr]^-{FA\diamond  \alpha(B,C)} \ar[dd]_{\alpha(A,B) \diamond FC}  && FA\diamond  F(B\star C)\ar[dd]^{\alpha(A,B\star C)} \\
&&  \\
F(A\star B)\diamond FC \ar[rr]^-{\alpha(A\star B,C)} && F(A\star B \star C). 
}$$

\medskip

\noindent
{\bf unitality}:  the following squares commute:

$$\xymatrix{
I\diamond FA  \ar@{=}[rr] \ar[d]_{\alpha_0\diamond FA} && FA   \ar@{=}[d]  \\
FJ\diamond FA \ar[rr]^-{\alpha(J,A)}  && F(J \star A) }
\quad \quad \quad 
 \xymatrix{
FA \diamond I  \ar@{=}[rr] \ar[d]_{FA\diamond \alpha_0} && FA   \ar@{=}[d]  \\
FA \diamond FJ \ar[rr]^-{\alpha(A,J)}  && F(A \star J) }.$$
We shall say that $\alpha$
is the {\bf product} and that $\alpha_0$ is the {\bf unit} of the lax structure.
The lax structure  is said to be {\bf normal} if $\alpha_0$ is an isomorphism. 
The lax structure  is said to be {\bf symmetric} if the monoidal categories $\bC$ and $\bD$ are symmetric and the following square commutes, 
$$ \xymatrix{
FA \diamond  FB \ar[rr]^-{\alpha}\ar[d]_{\sigma} && F(A\star B)  \ar[d]^{F\sigma}  \\
FB \diamond FA \ar[rr]^-{\alpha}  && F(B \star A)
 }$$
where $\sigma$ denotes the symmetries in both categories.
A  functor  $F:\bC\to \bD$ equipped with a lax structure $(\alpha, \alpha_0)$
is said to be {\bf lax}; a  functor  equipped with a normal  lax structure
is said to be {\bf normal lax}; a functor  equipped with a symmetric lax structure
is said to be {\bf symmetric lax}. 
\end{defi}

If $(\alpha, \alpha_0)$ is a lax structure on a functor $F:\bC\to \bD$,
then the map $\alpha(J,J): FJ \diamond FJ \to F(J\star J)=F(J)$
gives the object $FJ$ a monoid structure with unit $\alpha_0:I\to FJ$.
It follows that the unit of a lax structure $(\alpha, \alpha_0)$ is determined
by the product $\alpha$,
We may therefore describe the lax structure as a pair $(F,\alpha)$,

\begin{ex}
Let $\mathbb{I}$ be the monoidal category with a single object $I$
and a single morphism (the identity if $I$).
If $\bC=(\bC, \star, J)$ is a monoidal category,
then a lax functor $F:\mathbb{I}\to \bC$ is the same
thing as monoid object $C=F(I)$ in $\bC$.
The lax functor $F$ is symmetric if and only if the monoid $F(I)$
is commutative.
\end{ex}

\begin{ex}
Given a symmetric monoidal category $(\bC,\star,J)$, 
there is a canonical symmetric lax structure on the functor $\star : \bC \times \bC \to \bC$ given by the map
$$
\sigma_{23}:(A_1\otimes A_2)\otimes (B_1\otimes B_2)\simeq (A_1\otimes B_1)\otimes (A_2\otimes B_2).
$$
All conditions are consequence of the coherence of the symmetric structure of $\star$.
\end{ex}

\begin{ex}
Given a symmetric monoidal category $(\bC,\star,J)$, 
the hom functor $\hom:\bC^{op}\times \bC\to \Set$ has a canonical lax monoidal structure (for the cartesian product in $\Set$) given by the functoriality of $\otimes$
$$
\hom(A_1,A_2) \times \hom(B_1,B_2) \to \hom(A_1\otimes B_1,A_2\otimes B_2).
$$
All conditions are consequence of the coherence of the symmetric structure of $\star$.

If $(\bC,\star,J)$ is closed there is also a similar symmetric lax monoidal structure on the internal hom.
\end{ex}

We shall say that a natural transformation $\theta:F\to G: \bC\to \bD$ between
 lax functors $F=(F,\alpha,\alpha_0)$ and $G= G=(G,\beta,\beta_0))$
is a {\it homomorphism} if the following diagrams commutes
$$ \xymatrix{
FA \diamond  FB \ar[rr]^-{\alpha}\ar[dd]_{\theta(A)\diamond \theta(B)} && F(A\star B)  \ar[dd]^{\theta(A\star B)}  \\
&& \\
GA \diamond GB  \ar[rr]^-{\beta}  && G(A \star B)
 }\quad \quad \quad \quad
 \xymatrix{
  & F(J)  \ar[dd]^{\theta(J)}\\
I \ar[ru]^{\alpha_0}\ar[rd]_{\beta_0} &  \\
 & G(J).
 }
$$

\begin{prop}\label{compositionoflax}
The composite of two lax functors $F:\bC\to \bD$
and $G:\bD\to \bE$ has the structure of a lax functor
$GF:\bC\to \bE$. Monoidal categories, lax functors 
and homomorphisms form a 2-category.
\end{prop}

Similarly, the composite of two symmetric lax functors $F:\bC\to \bD$
and $G:\bD\to \bE$ has the structure of a symmetric lax functor
$GF:\bC\to \bE$. Symmetric monoidal categories, symmetric lax functors 
and homomorphisms form a 2-category.

\medskip

\begin{ex} \label{imageofmonoidsbylax} It follows from \ref{compositionoflax} that
the image $FM$ of a monoid $M$ by a lax functor $F:\bC \to \bD$
has the structure of a monoid. Moreover, the monoid $FM$ is commutative
if $M$ is commutative and the lax structure of $F$ is symmetric.
\end{ex}

\bigskip

\begin{defi}\label{defstructurecolax}
If $\bC=(\bC, \star, J)$ and  $\bD=(\bD, \diamond, I)$
are two monoidal categories, then a {\bf colax structure} on a functor
$F:\bC\to \bD$ is a pair $(\alpha, \alpha_0)$, where $\alpha$ is a natural
transformation
$$\alpha(A,B): F(A\star B) \to FA \diamond FB $$
and $\alpha_0$ is a  morphisme $:FJ\to I$ 
satisfying the following conditions,
\medskip

\noindent
{\bf associativity}: the following square commutes:
$$\xymatrix{
 F(A\star B \star C) \ar[rr]^-{\alpha(A\star B, C)} \ar[dd]_{\alpha(A,B \star C)}  && F(A\star B)\diamond FC\ar[dd]^{\alpha(A,B)\diamond FC} \\
&&  \\
FA\diamond F(B\star C) \ar[rr]^-{F\!A\diamond  \alpha(B,C)} && FA\diamond FB \diamond FC
}$$

\medskip

\noindent
{\bf unitality}:  the following squares commute:
$$\xymatrix{
F( J \star A)  \ar@{=}[rr] \ar[d]_{\alpha(J,A) } && FA   \ar@{=}[d]  \\
FJ\diamond FA \ar[rr]^-{\alpha_0\diamond F\!A}  && I \diamond FA }
\quad \quad \quad 
 \xymatrix{
F(A\star J)  \ar@{=}[rr] \ar[d]_{\alpha(A,J) } && FA   \ar@{=}[d]  \\
FA \diamond FJ \ar[rr]^-{F\!A \diamond \alpha_0}  &&  FA \diamond I}$$
We shall say that $\alpha$
is the {\bf coproduct} and that $\alpha_0$ is the {\bf counit} of the lax structure.
The colax structure  is said to be {\bf normal} if $\alpha_0$
est un isomorphisme. The colax structure  is said to be {\bf symmetric} 
if the monoidal categories $\bC$ and
$\bD$ are symmetric and the following square commutes, 
$$ \xymatrix{
F(A\star B)  \ar[rr]^-{\alpha}\ar[d]_{\sigma} &&  FA \diamond  FB  \ar[d]^{F\sigma}  \\
 F(B \star A) \ar[rr]^-{\alpha}  && FB \diamond FA
 }$$
where $\sigma$ denotes the symmetry in both categories.
A  functor  $F:\bC\to \bD$ equipped with a colax structure $(\beta, \beta_0)$
is said to be {\bf colax}; a functor  equipped with a normal colax structure 
is said to be {\bf normal colax}; a functor equipped with a symmetric colax structure 
is said to be {\bf symmetric colax}. 
\end{defi}

If $(\alpha, \alpha_0)$ is a colax structure on a functor $F:\bC\to \bD$,
then the map $\alpha(J,J): F(J)=F(J\star J) \to FJ \diamond FJ$
gives the object $FJ$ a comonoid structure with counit $\alpha_0:FJ\to I$.
It follows that the counit of a colax structure $(\alpha, \alpha_0)$ is uniquely determined 
by the coproduct $\alpha$. We may therefore describe the colax structure as a pair $(F,\alpha)$,

\begin{ex}
Let $\mathbb{I}$ be the monoidal category with a single object $I$
and a single morphism (the identity if $I$).
If $\bC=(\bC, \star, J)$ is a monoidal category,
then a colax functor $F:\mathbb{I}\to \bC$ is the same
thing as comonoid object $C=F(I)$ in  $\bC$.
The colax functor $F$ is symmetric if and only if the comonoid $F(I)$
is cocommutative.
\end{ex}

We shall say that a natural transformation $\theta:F\to G: \bC\to \bD$ between
colax functors $F=(F,\alpha,\alpha_0)$ and $G= G=(G,\beta,\beta_0))$
is a {\it homomorphism} of colax functors if the following diagrams commutes
$$ \xymatrix{
 F(A\star B) \ar[dd]_{\theta(A\star B)} \ar[rr]^-{\alpha} &&   FA \diamond  FB \ar[dd]^{\theta(A)\diamond \theta(B)}   \\
&& \\
 G(A \star B) \ar[rr]^-{\beta}  && GA \diamond GB  
 }\quad \quad \quad \quad
 \xymatrix{
  F(J)  \ar[dd]_{\theta(J)} \ar[rd]^{\alpha_0} &\\
 & I\\
 G(J)\ar[ru]_{\beta_0} . &
 }
$$

\begin{prop}\label{compositionofcolax}
The composite of two colax functors $F:\bC\to \bD$
and $G:\bD\to \bE$ has the structure of a colax functor
$GF:\bC\to \bE$. Monoidal categories, colax functors 
and homomorphisms form a 2-category.
\end{prop}

Similarly, the composite of two symmetric colax functors $F:\bC\to \bD$
and $G:\bD\to \bE$ has the structure of a symmetric colax functor
$G\circ F:\bC\to \bE$. Symmetric monoidal categories, symmetric colax functors 
and homomorphisms form a 2-category.

\begin{ex} \label{imageofcomonoidsbylax} It follows from \ref{compositionofcolax} that
the image $FC$ of a comonoid $C$ by a colax functor $F:\bC \to \bD$
has the structure of a comonoid. Moreover, the comonoid $FC$ is cocommutative
if $C$ is cocommutative and the colax structure of $F$ is symmetric.
\end{ex}

\medskip

\begin{prop}\label{rightadjointofcolax} The right adjoint 
of a colax functor (resp. symmetric colax functor) between monoidal categories (resp. symmetric monoidal categories)  has the structure of a lax functor (resp. symmetric lax functor).
Dually, the left adjoint  lax functor (resp. symmetric lax functor) between monoidal categories (resp. symmetric monoidal categories)  has the structure of a colax functo( resp. symmetric colax functor).
\end{prop}

\subsection{Enriched categories and strong functors}\label{enrichedcategoryapp}

Let $\bV=(\bV,\otimes, I)$ be a fixed monoidal category.
Recall that a category {\bf enriched
over}  $\bV$,  also called a $\bV$-{\bf category},
is a quadruple $\bC=(Ob(\bC), Hom, c, u)$ where,
\begin{itemize}
\item{} $Ob(\bC)$ is the class of {\it objects} of $\bC$);
\item{} $Hom$ is a map which associates
to every pair of objects $A,B\in Ob\bC$ an object $Hom(A,B)\in\bV$;
\item{} $c$ is a map
which associates to every triple $(A,B,C)$ of object of $\bC$ a morphism
$$c(A,B,C):Hom(B,C)\otimes Hom(A,B)\to Hom(A,C)$$
called the {\it composition law};
\item{} $u$ is a map 
which associates to every object $A\in\bC$ a morphism
$1_A:=u(A)=:I \to  Hom(A,A)$ called the {\it identity}, or the {\it unit} of $A$;
\item{} The composition law is required to be  associative
and the units to be two-sided with respect to the composition law.
\end{itemize}
We shall often put $\bC(A,B):=Hom(A,B)$ for every pair of objects $A,B\in \bC$.

\medskip

A {\bf morphism} $A\to B$ between two objects of a $\bV$-category $\bC$ is defined to be a map $I\to \bC(A,B)$
in the category $\bV$.  The composite of two morphisms $f:A\to B$ and $g:B\to C$
is the morphism $gf:A\to C$ obtained by composing the maps 
$$
\xymatrix{
I \ar[rr]^-{g\otimes f}&&  \bC(B,C)\otimes  \bC(A,B) \ar[rr]^-{c}  &&   \bC(A,C).
}$$
This defines an ordinary category $\bC_0$ (enriched over the category of sets) whose hom sets are $\bC_0(A,B)=\hom_V(I,\bV(A,B))$, called the {\bf underlying category} of $\bC$,

\begin{ex}
If $(\bV,\otimes,I)$ is a symmetric closed category, the internal hom defines an enrichment $\bV^\$$ of $\bV$ over itself.
The underlying category of $\bV^\$$ is $\bV$.
\end{ex}

\medskip

Recall that a {\bf strong functor}  $F:\bC \to\bE$
between two  $\bV$-categories, is a pair $(F_0,F)$,
where $F_0$ is a map $Ob(\bC)\to Ob(\bE)$ and
$F=(F_{AB})$ is a family of morphism
$F_{AB}:\bC (A,B)\to\bE (F_0A, F_0B)$
indexed by pairs $(A,B)$ of object of $\bC $ satisfying the following two conditions: (1)
the following diagram commutes
$$
\xymatrix{
\bC(B,C) \otimes \bC (A,B) \ar[d]_{F\otimes F}  \ar[rr]^c && \bC (A,C) \ar[d]^{F}  \\
\bE(F_0B,F_0C)\otimes \bE(F_0A,F_0B)  \ar[rr]^-{c}   & &\bE(F_0A,F_0C)  }
$$
for every triple  $(A,B,C)$ of objects of $\bC$; (2) $F(1_A)=1_{F_0A}$
for every object $A\in  \bC$. The morphism
$F_{AB}:\bC (A,B)\to\bE (F_0A, F_0B)$
is called the {\bf structure map} or the {\bf strength} of the functor. 
For simplicity, we shall write $FA$ instead of $F_0A$.

\medskip

Recall that a natural transformation $\alpha:F \to G$
between strong functors $F,G:\bC  \to\bE$
is said to be {\bf strong}
if the following square commutes
$$
\xymatrix{
\bC(A,B)  \ar[d]_F \ar[rr]^G  &&\bE(GA,GB) \ar[d]^{\bE(\alpha_A,GB)}  \\
\bE(FA,FB)  \ar[rr]^-{\bE(FA,\alpha_B) }   &  &\bE(FA,GB) .    }
$$
for every pair of objects $(A,B)$ of $ \bC$.
Strong natural transformations can be composed.
We thus obtain a category whose objects are the functors 
$\bC  \to\bE$ and whose morphisms are the
the strong natural transformations.
A  strong natural transformation $\alpha:F\to G$
is invertible if and only if the morphism $\alpha_A:FA\to GA$
is invertible for every object $A\in \bC $.

\medskip

The $\bV$-categories, strong functors and strong natural transformations
are forming a 2-category $\bV\bCat$.

\medskip

Let us now suppose that the category  $\bV$ is symmetric monoidal.
In which case  the (tensor) product $\bC\times \bE$ of
two  $\bV$-categories can be defined. An object
of $\bC\times \bE$ is a pair of objects $(C,E)\in Ob\bC\times Ob\bE$
and
$$Hom((C_1,E_1),(C_2,E_2)):=\bC(C_1,C_2)\otimes \bE(E_1,E_2).$$
The composition law 
$$Hom((C_2,E_2),(C_3,E_3))\otimes Hom((C_1,E_1),(C_2,E_2))\to 
Hom((C_1,E_1),(C_3,E_3))$$
is defined to be the composite
$$\xymatrix{
\bC(C _2,C_3)\otimes \bE(E_2,E_3)\otimes \bC(C_1,C_2)\otimes \bE(E_1,E_2)
\ar[d]_{id\otimes \sigma\otimes id}  &&  \\
\bC(C _2,C_3)\otimes \bC(C_1,C_2)\otimes \bE(E_2,E_3) \otimes \bE(E_1,E_2)
\ar[rr]^-{c\otimes c}  && \bC(C_1,C_3)\otimes \bE(E_1,E_3) 
}$$

\medskip

The operation of tensor product is a 2-functor
$$\times : \bV\bCat\times  \bV\bCat \to  \bV\bCat$$
and it defines a symmetric monoidal structure on the 2-category $ \bV\bCat$
(the notion of symmetric monoidal structure on a 2-category can be found in the litterature).

\medskip

The {\it opposite} $\bC^{op}$ of a $\bV$-category $\bC$
has the same objects as $\bC$, but the hom object is defined
by putting $\bC^{op}(A,B):=\bC(A,B)$.
The {\it reversed law composition}  
$$c^o:\bC^{op}(B,C)\otimes \bC^{op}(A,B)\to \bC^{op}(A,C)$$
is defined by composing of the maps
$$\xymatrix{
\bC(C,B)\otimes \bC(B,A)\ar[r]^-{ \sigma} & \bC(B,A) \otimes  \bC(C,B) \ar[r]^-{c} &
 \bC(C,A)
}$$
A {\bf strong contravariant functor} $F:\bC \to\bD$
between two  $\bV$-categories, is defined to be a strong (covariant)
functor $F:\bC^{op} \to\bD$.

\medskip

The tensor product  functor $\otimes:\bV\times \bV\to \bV$
of a symmetric monoidal closed category $\bV$
is strong. The strength is defined by the canonical map
$$\psi=\lambda^{3,4}(\alpha) :\bV(A_1,A_2)\otimes \bV(B_1,B_2) \to  \bV(A_1\otimes B_1, A_2\otimes B_2),$$
where $\alpha$ is the composite
$$\xymatrix{
\bV(A_1,A_2)\otimes \bV(B_1,B_2)\otimes A_1\otimes B_1\ar@/_1pc/[rrrd]_\alpha
\ar[rrr]^-{id\otimes \sigma\otimes id}  &&& \bV(A_1,A_2)\otimes A_1\otimes  \bV(B_1,B_2)\otimes B_1 \ar[d]^{ev\otimes ev}   \\
&&&  A_2\otimes B_2
}$$

\medskip

The hom functor  $\bV(-,-):\bV^{op}\times \bV\to \bV$
 of a symmetric monoidal closed category $\bV$ is strong.
 The strength is given by the canonical map
$$\psi={\lambda}^2(c^3)\sigma:\bV(A_2,A_1) \otimes  \bV(B_1,B_2) \to 
\bV(\bV(A_1,B_1),  \bV(A_2,B_2) )$$
where $c^3$ is the three-fold  composition law,
$$\xymatrix{
\bV(B_1,B_2)\otimes  \bV(A_1,B_1) \otimes \bV(A_2,A_1) \ar[rr] &&\bV(A_2,B_2) .
}$$
More generally, the hom functor  $\bC(-,-):\bC^{op}\times \bC\to \bV$
 of a $\bV$-category $ \bC$ is strong. The strength is given by the canonical map
$$\psi={\lambda}^2(c^3)\sigma:\bC(A_2,A_1) \otimes  \bC(B_1,B_2) \to 
\bV(\bC(A_1,B_1),  \bC(A_2,B_2) )$$
where $c^3$ is the three-fold  composition law,
$$\xymatrix{
\bC(B_1,B_2)\otimes  \bC(A_1,B_1) \otimes \bC(A_2,A_1) \ar[rr] &&\bC(A_2,B_2) .
}$$

\medskip

A {\bf strong  monoidal structure} on a  $\bV$-category $ \bC$
is a five-tuple $(\otimes, \alpha, I, \lambda, \rho)$ consisting of
a strong functor $\otimes : \bC\times \bC\to \bC$,
a strong  associativity isomorphism
 $\alpha(A,B,C):(A\otimes B)\otimes C\to A\otimes (B\otimes C)$,
 a unit object $I$, and two strong unit isomorphisms
$\lambda(A):I\otimes A\to A$ and $\rho(A):A\otimes I\to A$
satifying the coherence axioms of an ordinary monoidal category.
The strength of the tensor product functor is a natural map
$$\psi: \bC(A_1,A_2)\otimes \bC(B_1,B_2) \to  \bC(A_1\otimes B_1, A_2\otimes B_2)$$
for a quadruple of objects $(A_1,A_2, B_1, B_2)$ in the category $\bC$.
The map $\psi$ is defining a lax structure on the hom functor $\bC(-,-):\bC^{op} \times \bC\to \bV$.
A  {\bf monoidal $\bV$-category} 
 is defined to be a $\bV$-category $\bC$  equipped with a 
strong monoidal structure $(\otimes, I, \alpha, \lambda, \rho)$.
The notion of {\bf symmetric} monoidal  $\bV$-category is defined similarly.

\medskip

Let $\bC$ be a category enriched over a symmetric monoidal closed category $\bV$.
If $F:\bC\to \bV$ is a strong functor and
$A$ is an object of $\bC$, then the map $\alpha\mapsto \alpha(1_A)$
is a bijection between the strong natural transformations $\alpha:\bC(A,-)\to F$
and the maps $a:I\to F(A)$ (Yoneda lemma). The functor $F$ is said
to be {\it represented} by the pair $(A,a)$
if the corresponding natural transformation $\alpha:\bC(A,-)\to F$
is invertible.  The functor $F$ is said
to be {\bf representable} if it can be represented by
a pair $(A,a)$. Dually, 
if $F:\bC\to \bV$ is a strong contravariant functor and
$A$ is an object of $\bC$, then the map $\alpha\mapsto \alpha(1_A)$
is a bijection between the strong natural transformations $\alpha:\bC(-,A)\to F$
and the maps $a:I\to F(A)$ (Yoneda lemma). The functor $F$ is said
to be {\it represented} by the pair $(A,a)$
if the corresponding  natural transformation $\alpha:\bC(-,A)\to F$
is invertible.  The contravariant functor $F$ is said
to be {\bf representable} if it can be represented by
a pair $(A,a)$.

\medskip

Recall that a {\bf strong adjuntion} $\theta: F\dashv G$ between two strong functors $F:\bC\to \bE$
and $G:\bE\to \bC$ is defined to be a strong natural isomorphism
$$\theta=\theta_{AB}: \bC(A,GB)\simeq \bE(FA,B).$$
The {\bf unit} of the adjunction is the strong natural transformation
$\eta_A:A \to GFA$ obtained by putting $\eta_A=\theta^{-1}(1_{FA})$
for every object $A\in  \bC$ 
and its {\bf counit}  is the strong natural transformation
$\epsilon_B:FGB \to B$ obtained by putting $\epsilon_B=\theta(1_{GB})$
for every object $B\in  \bE$.
A strong functor $F:\bC\to \bE$
has a (strong) right adjoint if an only if its composite with
the functor $\bE(B,-): \bE \to \bV$
is representable for every object $B\in \bE $.
Dually, a strong functor $G:\bE\to \bC$
has a (strong) left adjoint if an only if its composite with
the contravariant functor $\bE(-,A): \bE \to \bV$
is representable for every object $A\in \bC $.

\medskip

A category $\bC$ enriched over a monoidal category $\bV$
 is said to be {\bf tensored} by $\bV$
 if the (strong) functor $\bC(A,-):\bC\to \bV$
 has a (strong) left adjoint $(-)\diamond A: \bV\to \bC$
 for each object $A\in \bC$.
 The object $X\diamond A \in  \bC$ is
equipped with a map $\eta: X\to \bC(A, X\diamond A)$ which is universal
in the following sense: 
for every object $B\in  \bC$
 and every map $f:X\to  \bC(A, B)$, there exist
 a unique map $g:X\diamond A \to B$
 such that $\bC(A, g)\eta =f $.
  Dually,  a  $\bV$-category $\bC$
 is said to be {\bf cotensored} by $\bV$
 if the (strong) contravariant functor $\bC(-,B):\bC\to \bV$
 has a (strong) right adjoint $[-,B]: \bV\to \bC$
  for each object $B\in \bC$. 
 The object $[X,B] \in  \bC$ is
equipped with a map $\eta: X\to \bC([X,B],B)$ having the following couniversal property:
for every object $A\in  \bC$
 and every map $f:X\to \bC(A,B)$, there exist
 a unique map $g:A \to [X,B]$
 such that $\bC(g, A)\eta =f $.

 \medskip

The $\bV$-category $\bC$ is said to be {\bf cocomplete}
if it is tensored over $\bV$ and cocomplete
as an ordinary category.
 Dually, $\bC$ is said to be {\bf complete}
if it is cotensored over $\bV$ and complete
as an ordinary category.
A $\bV$-category $\bC$ is said to be {\bf bicomplete}
if it is complete and cocomplete 
as a $\bV$-category.

\label{stronglimit}

For $\bC$ be a bicomplete $\bV$-category, $I$ an ordinary category and $X:I\to \bC_0$ a diagram with limit (resp. colimit) $L$.
We shall say that $L$ is a {\bf strong limit} (resp. a {\bf strong colimit}), if for any $A\in \bC$, 
$\bC(A,L)$ is the limit in $\bV_0$ of the diagram $\bC(A,X_i)$ (resp. $\bC(L,A)$ is the limit in $\bV_0$ of the diagram $\bC(X_i,A)$).
If $\bC$ is bicomplete, 
$\bC(A,-):\bC\to \bV$ is right adjoint to $-\otimes A:\bV\to \bC$ and 
$\bC(-,A):\bC^{op}\to \bV$ is right adjoint to $A^{(-)}:\bV^{op}\to \bC$.
In particular, $\bC(A,-)$ sends limits in $\bC_0$ to limits in $\bV_0$ and $\bC(-,A)$ sends colimits in $\bC_0$ to limits in $\bV_0$.
Hence every limit or colimit in $\bC$ is strong.

\medskip
We shall use the following notations.
$\bV\bCat^{tens}$ and $\bV\bCat^{cotens}$ are the full sub-2-categories of $\bV\bCat$ generated respectfully by the tensored, cotensored $\bV$-categories. $\bV\bCat^{bitens}$ is the full sub-2-category of $\bV\bCat$ generated by the $\bV$-categories that are both tensored and cotensored.
$\bV\bCat^{comp}$, $\bV\bCat^{cocomp}$ and $\bV\bCat^{bicomp}$ are the full sub-2-categories of $\bV\bCat$ generated respectfully by the complete, cocomplete and bicomplete $\bV$-categories.
The natural inclusions between them are summarized at the end of the next section.

\subsection{$\bV$-modules and $\bV$-opmodules}\label{Vopmodule}\label{Vmodule}

The category of endofunctors of an arbitrary
category $\bC$ has the structure of a monoidal category $End( \bC)$ 
in which the product
is composition of endofunctors.
If $\bV=(\bV,\otimes, I)$ is a monoidal category,
then a monoidal functor
$\pi: \bV\to End( \bC)$
is a {\bf representation} of the monoidal category $ \bV$ in the category of endomorphims of $\bC$.
To every representation $\pi: \bV\to End( \bC)$
we can associate an {\bf action} $\diamond: \bV\times  \bC\to  \bC$
of the category $ \bV$ on  $\bC$ by putting 
$X\diamond A=\pi(X)(A)$ for $X\in  \bV$
and  $A \in  \bC$.
From the natural isomorphism $ \pi(X\otimes Y) \simeq \pi(X)\circ \pi(Y)$,
we obtain a natural isomorphism
 $\mathbf{a}(X,Y,A):(X\otimes Y)\diamond A \simeq X\diamond (Y\diamond A)$,
and from the natural isomorphism $\pi(I)\simeq I_{\bC}$
we obtain a natural isomorphism $\lambda(A):I\diamond A\to A$.
The pair $(\mathbf{a},\lambda)$ satisfies the following coherence conditions:
 
 \begin{itemize}
\item{}  the following pentagon commutes for every quadruple of objects $(X,Y,Z,A)\in \bV^3\times \bC$,
where $XY:=X\otimes Y$.
$$\xymatrix{
 & (X(YZ))\diamond A  \ar[rr]^-{\mathbf{a}}  && X\diamond ((YZ)\diamond A)  \ar[dddr]^-{X\diamond \mathbf{a}} & \\
&& &&\\
&& &&\\
((XY)Z)\diamond A \ar[ruuu]^{\mathbf{a}\diamond A}\ar[rrdd]_{\mathbf{a}} && &&   X\diamond (Y\diamond (Z\diamond A)) \\
&& &&\\
&& (XY)\diamond (Z\diamond A) \ar[rruu]_{\mathbf{a}}  &&
}$$
\item{} the following two triangles commute for
every couple of objects $(X,A)\in \bV\times \bC$,
$$
 \xymatrix{
 (IX) \diamond A \ar[rr]^-{\mathbf{a}}  \ar[rd]_{\lambda(X)\diamond A}  && I\diamond (X \diamond A)  \ar[dl]^-{\lambda(X\diamond A)} \\
&X\diamond A &
}
\quad \quad \quad
\xymatrix{
 (XI) \diamond A \ar[rr]^-{\mathbf{a}}  \ar[rd]_{\rho(X)\diamond A}  && X\diamond (I \diamond A)  \ar[dl]^-{X\diamond \lambda(A)} \\
&X\diamond A &
}$$
\end{itemize}

We shall say that  a triple  $(\diamond,a,\lambda)$ satisfying these conditions 
 is a {\it left action} of the monoidal category $\bV$ on the category $\bC$.
There is a bijection beween the
 monoidal functors $\bV \to {End}( \bC)$
and the left actions of  $\bV$ on $\bC$.
 We shall say that a category $\bC$ equipped with a left
 action of $\bV$ is a  {\bf left $\bV$-module}.

\medskip

A $\bV$-category $\bC$ which is tensored by  $\bV$
has the structure of a left  $\bV$-module,
where $X\diamond A$ is the tensor product of $A$ by $X$.
The functor $(X,A)\mapsto X\diamond A$
is divisible on the right and we have $\bV(A,B)=B/A$
for every $A,B\in \bV$.

\medskip

Conversely, we shall say that a left $\bV$-module $\bC$
is {\bf divisible on the right} if the action functor
$\bV \times\bC \to\bC$
is divisible on the right. In which case the functor 
 $(-)\diamond A:\bV \to  \bC$
 has a right adjoint $(-)/A:  \bC \to   \bV$
for every object $A\in  \bC$. If $B\in  \bC$, then the object $B/A\in  \bV$
is equipped with a map $ev=ev(A,B):(B/A)\diamond A\to B$ 
having the following couniversal property:
for every object $C\in  \bV$
 and every map $f:C\diamond A\to B$, there exist
 a unique map $g:C\to B/A$
 such that $ev(g \diamond X)=f$.
  In particular, for any triple of objects $A,B,C\in \bC$, there exists a unique map 
 $c=c(A,B,C):(C/B)\otimes (B/A)\to C/A$ such that 
 the following diagram commutes
 $$
\xymatrix{
(C/B\otimes B/A)\diamond A  \ar@{=}[r]^a  \ar[d]_{c\diamond A} &C/B\diamond ((B/A)\diamond A) \ar[rr]^-{(C/B)\diamond ev} && (C/B) \diamond B\ar[d]^{ev} \\
(C/A)\diamond A \ar[rrr]^-{ev}&&& C
}
$$
 The map $c(A,B,C)$ is the composition law of an enrichment
 of  the category $\bC$ over $\bV$ if we put
 $Hom(A,B)=B/A$ for $A,B\in \bC$.
 The unit $u_A:I\to Hom(A,A)$ is the unique map
 $u_A:I\to A/A$ such that $ev(A,A)(u_A\diamond A)$ is the isomorphism $\lambda(A):I\diamond A\simeq A$.
 The resulting $ \bV$-category $ \bC$ is tensored by $ \bV$
 and the tensor product of an object $A\in \bC$
 by an object $X\in  \bV$ is the product $X\diamond A$.

\begin{prop}
If $\bV$ is a monoidal category, then the notion of a  $\bV$-category tensored by $\bV$
is equivalent to the notion of a left $\bV$-module divisible on the right.
\end{prop}

Let us now suppose that $\bC$ and $\bD$ are $\bV$-modules.
A {\bf lax modular structure}\label{laxmodularstructure} of a functor $\bC \to \bD$ is 
 a natural transformation $\psi:X\diamond F(A)\to F(X\diamond A)$
 for which the following
two diagrams commute,
$$\xymatrix{
(X\otimes Y)\diamond F(A)  \ar@{=}[r]^a  \ar[d]_{\psi} &  X\diamond (Y\diamond F(A)) \ar[rr]^-{X\diamond \psi} && X \diamond F(Y\diamond A) \ar[d]^{\psi} \\
 F( (X\otimes Y)\diamond A)   \ar[rrr]^-{F(a)}&&& F( X \diamond(Y\diamond A))
}\quad \quad 
\xymatrix{ 
I\diamond F(A)\ar[d]_{\psi} \ar[rd]^{\lambda(F(A))}&\\
F(I\diamond A) \ar[r]_-{F(\lambda_A)} &  F(A) \\
}
$$
We shall say that the pair $(F,\psi)$ is a  {\bf lax modular functor} $\bC \to \bD$.
We shall say that a natural transformation $\alpha:F\to G$ between
lax modular functors $(F,\psi)$ and $(G,\psi)$ is {\bf strong} if the following
diagram commutes,
$$\xymatrix{
X \diamond FA  \ar[d]_{\psi }  \ar[rr]^-{X\diamond \alpha } &&  X\diamond GA \ar[d]^{\psi}  \\
F(X\diamond A)  \ar[rr]^-{\alpha }&&  G(X\diamond A)
}$$
The composite of two strong transformations $\alpha:F\to G$ and $\beta:  G\to H$
is a stronf transformation $\beta \alpha:F\to H$.
There is then a category $\MMod(\bV)(\bC,\bD)$
whose objects are the lax modular functors $\bC \to \bD$  
and whose morphisms are the strong transformations.
Lax modular functors $F:\bC \to \bD$ and $F:\bC \to \bD$
can be composed. There is then a 2-category $\MMod(\bV)$ 
whose objects are the (left) $\bV$-modules
 and whose morphisms are lax modular functors
 and whose 2-cells are the strong transformations.
We shall denote by $\MMod(\bV)^{rdiv}$ the full sub-2-category of $\MMod(\bV)$  generated by $\bV$-modules that are divisible on the right.

\medskip

Let us now suppose that  $(F,\psi):\bC \to \bD$ 
is a lax modular functor between two $\bV$-modules
divisible on the right. For every pair of objects $A,B\in \bC$
there is then a unique map $\phi(A,B):\bC(A,B)\to \bD(FA,FB)$
for which the following square commutes,
$$
\xymatrix{
\bC(A,B) \diamond FA  \ar[rr]^\psi  \ar[d]_{\phi\diamond FA} &&  F( \bC(A,B) \diamond A) \ar[d]^{F(ev)}  \\
 \bD(FA,FB)\diamond FA  \ar[rr]^-{ev}&&  FB
}
$$
The map $\phi(A,B)$
is the strength of a $\bV$-functor $(F,\phi):\bC \to \bD$.
This defines a bijection between the lax modular functors $(F,\psi):\bC \to \bD$ 
and the $\bV$-functors $(F,\phi):\bC \to \bD$.
The inverse bijectition takes a $\bV$-functor $(F,\phi):\bC \to \bD$,
to the lax modular structure $\psi:X\diamond F(A)\to F(X\diamond A)$
obtained by composing the maps
$$
\xymatrix{
X \diamond FA  \ar[rr]^-{\eta \diamond FA } && \bC(A, X\diamond A)\diamond FA    \ar[rr]^-{\phi \diamond FA} && 
 \bD(FA,F(X\diamond A)) \diamond FA\ar[d]^{ev}  \\
& &&& F(X\diamond A)
}
$$
Moreover,  $(F,\psi)$ and $(G,\psi)$ are lax modular functors
$\bC \to \bD$, then a natural transformation $\alpha:F\to G$
is strong if and only if it is strong as a natural transformation of $\bV$-functors   $(F,\phi)\to (G,\phi)$.
Hence the category $\MMod(\bV)(\bC,\bD)$ is isomorphic
to the category $\bV\bCat(\bC,\bD)$

\bigskip

There is a dual notion of a  {\it right action} of a monoidal category $\bV$ on 
 a category $\bC$.  We shall say that a category $\bC$ equipped with a right
 action of $\bV$ is {\bf right $\bV$-module}.
 The {\it reverse} tensor product $X\otimes^r Y$
of two objects $X,Y\in \bV$ is defined by
putting $X\otimes^r Y=Y\otimes X$.
This defines a new monoidal category $\bV^r=({\cal V},\otimes^r ,I)$.
A right action of $\bV$ on a category
$\bC$ is the same thing as a left action of $\bV^r$
on $\bC$. The monoidal categories $\bV^r$ and $\bV$
are equivalent when $\bV$ is symmetric,
In which case the  left and right actions of  $\bV$
on a category $\bC$ coincide.

 \medskip

We shall say that the monoidal category $(\bV^{op},\otimes^r ,I)$
is the {\bf reverse opposite} of $\bV$
and we shall denote it by $\bV^{rop}$.
We shall say that a left module over $\bV^{rop}$
 is an {\bf opmodule} over $\bV$.
An opmodule structure on a category $\bC$ is defined
by a functor $[-,-]:\bV^{op} \times \bC  \to  \bC$
and two natural isomorphism  $\mathbf{a}:[Y\otimes X,A] \simeq [X,[Y,A]] $
satisfying the following conditions:

\begin{itemize}
\item{}  the following pentagon commutes for every quadruple of objects $(X,Y,Z,A)\in \bV^3\times \bC$,
where $XY:=X\otimes Y$.
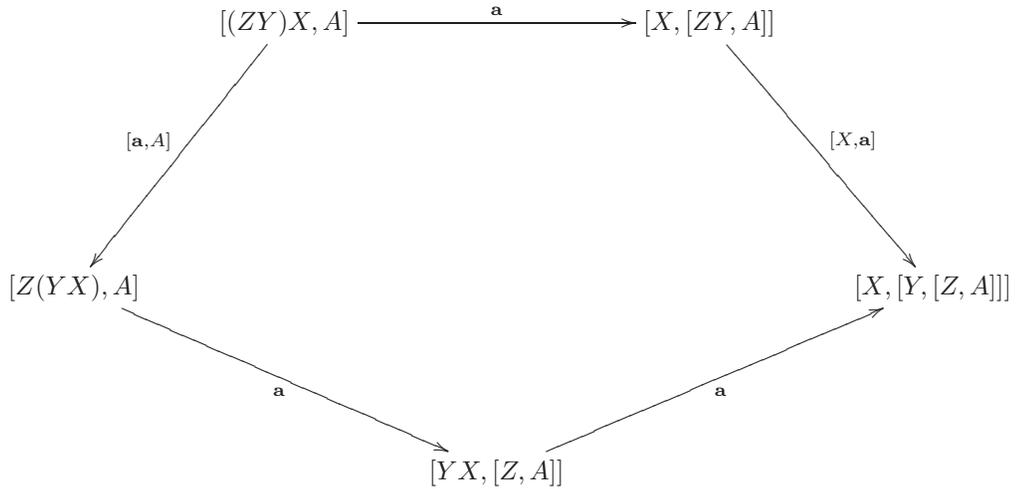
$$\xymatrix{
 & [(ZY)X,A]  \ar[rr]^-{\mathbf{a}}  \ar[lddd]_{[\mathbf{a}, A]}  && [X,[ZY,A]]  \ar[dddr]^-{[X,\mathbf{a}]} & \\
&& &&\\
&& &&\\
[Z(YX), A] \ar[rrdd]_{\mathbf{a}} && &&   [X, [Y,[Z,A]]] \\
&& &&\\
&& [YX, [Z,A]] \ar[rruu]_{\mathbf{a}}  &&
}$$
\item{} the following two triangles commute for
every couple of objects $(X,A)\in \bV\times \bC$,
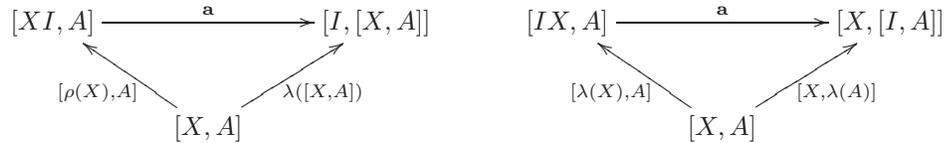
$$\xymatrix{
 [XI,A] \ar[rr]^-{\mathbf{a}} && [I,[X,A]] \\
&[X, A]  \ar[lu]^{[\rho(X),A]}  \ar[ur]_-{ \lambda([X,A])}   &
} \quad \quad \quad
\xymatrix{
 [IX,A] \ar[rr]^-{\mathbf{a}}  && [X,[I,A]]  \\
&[X,A]  \ar[lu]^{[\lambda(X),A]}  \ar[ur]_-{[X, \lambda(A)]} &
}$$
\end{itemize}

\medskip

Conversely, we shall say that a $\bV$-opmodule $\bC$
is {\bf codivisible on the right} if the action functor
$\bV^{op} \times\bC \to\bC$
is codivisible on the right. In which case the contravariant functor 
 $[-,A]:\bV \to  \bC$
 has a contravariant right adjoint $(-)\backslash A:  \bC \to   \bV$
for every object $A\in  \bC$. For each object $B\in  \bC$,
the object $A/B\in  \bV$
is equipped with a map $cv: B\to [A/B, A]$ having the following couniversal
property: for every object $C\in  \bV$
 and every map $f:B\to [C,A]$, there exist
 a unique map $g:C\to B/A$
 such that $[g,A](cv)=f$.
  In particular, for any triple of objects $A,B,C\in \bC$, there exists a unique map 
 $c=c(A,B,C): C/B\otimes B/A  \to C/A$ such that 
 the following diagram commutes
 $$
\xymatrix{
[C/B\otimes B/A,C]  \ar@{=}[r]^a   &[B/A, [C/B, C]] && \ar[ll]_(0.4){ [B/A, cv]} [B/A, B]  \\
 [C/A, C] \ar[u]^{[c,C]}&&& \ar[lll]_-{cv}  A \ar[u]_{cv}
}
$$
 The map $c(A,B,C)$ is the composition law of an enrichment
 of  the category $\bC$ over $\bV$ if we put
 $Hom(A,B)=B/A$ for $A,B\in \bC$.
The resulting $ \bV$-category $ \bC$ is cotensored by $ \bV$
 and the cotensor product of an object $A\in \bC$
 by an object $X\in  \bV$ is the object $[X,A]$.

\begin{prop}  If $\bV$ is a monoidal category, then the notion of a  $\bV$-category cotensored by $\bV$
is equivalent to the notion of $\bV$-opmodule 
codivisible on the right.
\end{prop}

Let us now suppose that $\bC$ and $\bD$ are $\bV$-opmodules.
A {\bf colax modular structure}\label{colaxmodularstructure} of a functor $\bC \to \bD$ is a natural transformation $\psi: F[X,A]\to [X,FA]$ for which the following two diagrams commute,
$$\xymatrix{
[Y\otimes X,FA]  \ar@{=}[r]^a &  [X, [Y, FA]] &&\ar[ll]_-{[X,\psi] }  [X, F[Y,A]] \\
F[Y\otimes X,A] \ar[u]^{\psi}   &&& \ar@{=}[lll]_-{F(a)}  \ar[u]_{\psi} F[X,[Y,A]]
}\quad \quad 
\xymatrix{ 
[I,FA] &\\
F[I,A] \ar[u]^{\psi} & \ar[l]^-{F\lambda(A)}  F(A)\ar[ul]_-{\lambda(F\! A)} 
} 
$$
We shall say that the pair $(F,\psi)$ is a  {\bf colax modular functor} $\bC \to \bD$.
We shall say that a natural transformation $\alpha:F\to G$ between
two colax modular functors $(F,\psi)$ and $(G,\psi)$ is {\bf strong} if the following
diagram commutes,
$$
\xymatrix{
[X,FA]  \ar[rr]^-{[X, \alpha] } &&  [X,GA]   \\
F[X,A]  \ar[u]_{\psi }  \ar[rr]^-{\alpha }&&  G[X,A] \ar[u]^{\psi}
}
$$
The composite of two strong transformations $\alpha:F\to G$ and $\beta:  G\to H$
is a strong transformation $\beta \alpha:F\to H$.
There is then a category $\OpMod(\bV)(\bC,\bD)$
whose objects are the colax modular functors $\bC \to \bD$  
and whose morphisms are the strong transformations.
Colax modular functors $F:\bC \to \bD$ and $F:\bC \to \bD$
can be composed. There is then a 2-category $\OpMod(\bV)$ 
whose objects are the $\bV$-opmodules
 and whose morphisms are  colax modular functors
 and whose 2-cells are strong transformations.
We shall denote by $\OpMod(\bV)^{rcodiv}$ the full sub-2-category of $\MMod(\bV)$  generated by $\bV$-modules that are divisible on the right.

\medskip

Let us now suppose that  $(F,\psi):\bC \to \bD$ 
is a colax modular functor between two $\bV$-opmodules
divisible on the right. For every pair of objects $A,B\in \bC$
there is then a unique map $\phi(A,B):\bC(A,B)\to \bD(FA,FB)$
for which the following square commutes,
$$
\xymatrix{
[\bC(A,B),FA]  &&  \ar[ll]_-\psi  F[\bC(A,B),A]  \\
[ \bD(FA,FB),FA]  \ar[u]^{[\phi,FA] }  && \ar[ll]_-{cv}  FB  \ar[u]_{F(cv)}
}
$$
The map $\phi(A,B)$
is the strength of a $\bV$-functor $(F,\phi):\bC \to \bD$.
This defines a bijection between the colax modular functors $(F,\psi):\bC \to \bD$ 
and the $\bV$-functors $(F,\phi):\bC \to \bD$.
The inverse bijection takes a $\bV$-functor $(F,\phi):\bC \to \bD$,
to the natural transformation $\psi: F[X,A] \to [X,FA]$
obtained by composing the maps
$$
\xymatrix{
[X,FA]   &&\ar[ll]_-{[\eta,FA] } [\bC(A, [X,A]),FA]  &&   \ar[ll]_-{[\phi , FA] } 
 [\bD(FA,F[X,A]),FA] \\
& &&& F[X,A]\ar[u]_{cv}.
}
$$
Moreover, if $(F,\psi)$ and $(G,\psi)$ are colax modular functors
$\bC \to \bD$, then a natural transformation $\alpha:F\to G$
is strong if and only if it is strong as a natural transformation of $\bV$-functors  $(F,\phi)\to (G,\phi)$.
Hence the category $\OpMod(\bV)(\bC,\bD)$ is isomorphic
to the category $\bV\bCat(\bC,\bD)$.

\bigskip

\begin{prop}\label{triplestrength}
There are natural equivalences of 2-categories 
$$
\bV\bCat^{tens}=\MMod(\bV)^{rdiv}
\et 
\bV\bCat^{cotens}=\OpMod(\bV)^{rcodiv}.
$$
In particular, the notion of strong functor between two bicompletes $\bV$-categories is equivalent to that of lax modular functor and to that of colax modular functor.
\end{prop}

\medskip

The following diagram of inclusions of 2-categories compares all the notions of the sections \ref{enrichedcategoryapp} and \ref{Vmodule}.
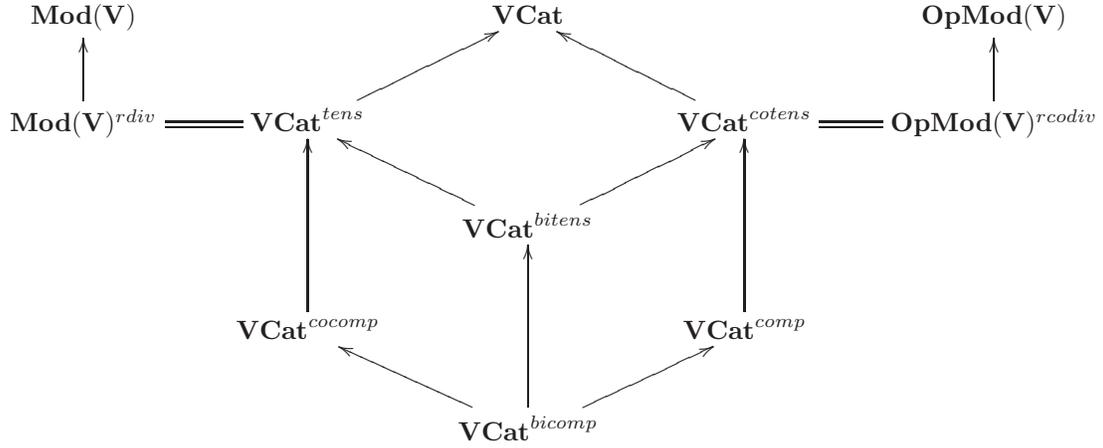
$$\xymatrix{
\MMod(\bV)&&\bV\bCat&&\OpMod(\bV)\\
\MMod(\bV)^{rdiv}\ar[u]&\bV\bCat^{tens}\ar[ru] \ar@{=}[l]&&\bV\bCat^{cotens}\ar[lu]\ar@{=}[r]&\OpMod(\bV)^{rcodiv}\ar[u]\\
&&\bV\bCat^{bitens}\ar[ru]\ar[lu]\\
&\bV\bCat^{cocomp}\ar[uu]&&\bV\bCat^{comp}\ar[uu]\\
&&\bV\bCat^{bicomp}\ar[ru]\ar[lu]\ar[uu]\\
}$$

\subsection{Base change}\label{basechange}

\subsubsection{Transfer of enrichments}\label{transferofenrichment}

The enrichments of a category over a monoidal category $\bV$ can be moved
to another base $\bV'$ along a lax monoidal functor $\bV\to \bV'$.

\medskip

Let $(U,\psi,\psi_0)$ be a lax functor  $U:\bV\to \bV'$ 
between two monoidal categories $\bV=(\bV,\otimes, I)$
and $\bV'=(\bV', \otimes, I)$. If $\bC$ is a $\bV$-category
we shall denote by $\bC^U$ the $\bV'$-category having 
same objects as $\bC$, but whose
hom object 
is defined by putting $\bC^U(A,B)=U\bC(A,B)$
and whose composition law $\bC^U(B,C)\otimes \bC^U(A,B)\to  \bC^U(A,C)$
is defined to be the composite 
$$\xymatrix{
U\bC(B,C)\otimes U\bC(A,B) \ar[r]^-{\psi} & U(\bC(B,C)\otimes \bC(A,B))\ar[r]^-{Uc}  & U\bC(A,C) 
}
$$
where $c$ is the composition law 
$\bC(B,C)\otimes \bC(A,B)\to  \bC(A,C)$
of the category $\bC$. The unit morphism $I\to  \bC^U(A,A)$
of an object $A\in  \bC^U$
is defined to be the composite $\psi_0U(1_A):I\to  UI \to U\bC(A,A)$
where $1_A:I\to \bC(A,A)$ is the unit morphism of the object $A\in \bC$.
We shall say that the category  $\bC^U$ is obtained by {\bf transfering} along $U$
the enrichment of category $\bC$.
A strong functor $F:\bC\to \bD$ between $\bV$-categories 
induces a strong functor $F^U:\bC^U\to \bD^U$.
Similarly, a strong natural transformation $\alpha:F\to G:\bC\to \bD$
induces a strong natural transformation $\alpha^U:F^U\to G^U:\bC^U\to \bD^U$
It is easy to see that this defines a 2-functor 
$$(-)^U:\bV\bCat \to \bV'\bCat.$$

\medskip

Let us now suppose that the monoidal categories $\bV$ and $\bV'$
are symmetric, and that the lax functor  $U:\bV\to \bV'$ 
is symmetric. In which case  the 2-functor
$(-)^U:\bV\bCat \to \bV'\bCat$
takes a monoidal $\bV$-category $\bC=( \bC,\star, J)$ 
to a monoidal $\bV'$-category $\bC^U=( \bC^U,\star, J)$.
The strength of the functor $\star:\bC^U\times \bC^U\to \bC^U$

$$\xymatrix{
U\bC(A_1,A_2) \otimes U\bC(B_1,B_2)\ar[r]^-\alpha& U(\bC(A_1,A_2) \otimes \bC(B_1,B_2)) \ar[rr]^-{U(\theta)} && U(\bC(A_1\otimes B_1,A_2\otimes B_2))
}$$
where $\alpha$ is the lax structure of $U$ and $\theta$ is the strength of the functor $\star:\bC\times \bC\to \bC$.
The strength of the associativity isomorphism $a:(A\otimes B)\otimes C\simeq A\otimes (B\otimes C)$ in $\bC$ is 
$$\xymatrix{
\bC(A_1,A_2)\otimes \bC(B_1,B_2)\otimes \bC(C_1,C_2) \ar[r]\ar[d] & \bC(A_1,A_2)\otimes \bC(B_1\otimes C_1,B_2\otimes C_2)\ar[dd] \\
\bC(A_1\otimes B_1,A_2\otimes B_2)\otimes \bC(C_1,C_2) \ar[d] \\
\bC((A_1\otimes B_1)\otimes C_1,(A_2\otimes B_2)\otimes C_1) \ar@{=}[r] & \bC(A_1\otimes (B_1\otimes C_1),A_2\otimes (B_2\otimes C_1))
}$$
Composing with $U$ and using its the lax structure of $U$ we deduce that $a$ is also strong in $(\bC)^R$.
The proof that the unit isomorphisms are strong is similar. 
The monoidal $\bV'$-category $\bC^U$ is symmetric when the  $\bV'$-category  $\bC$ is symmetric.
If $[B,C]$ is the internal hom in $\bC$, 
the isomorphism $U\bC(A\otimes B,C)=U\bC(A,[B,C])$ proves that the strong functor $(A,B,C)\mto U\bC(A\otimes B,C)$ is representable.
This proves that the monoidal $\bV'$-category $\bC^U$ is closed when the  $\bV'$-category  $\bC$ is closed.

If a $\bV$-category $\bA$ is enriched over $\bC$, then the $\bV'$-category $\bA^U$  is enriched over $\bC^U$;
if $\bA$ is tensored (resp.  cotensored) over $\bC$, then  $\bA^U$  is tensored (resp.  cotensored) over $\bC^U$.

\medskip

In particular, if the category $\bV$ is symmetric monoidal closed, then the $\bV'$-category
$\bV^U$ is symmetric monoidal closed.  If a category  $\bC$
 is enriched over $\bV$, then the $\bV'$-category $\bC^U$
 is enriched over $\bV^U$.

\bigskip
Recall that a {\bf symmetric lax monoidal adjunction}\label{laxmonoidaladj} is an adjunction $U\dashv R$ where $U:\bV\to \bV'$ is a symmetric monoidal functor.
Then the functor $R$ is automatically lax monoidal and symmetric. 
Let us describe the lax structure $(\alpha,\alpha_0)$ of the functor $R$.
The adjunction $U\dashv R$ is defined by a unit-counit pair $(\eta,\epsilon)$.
Then the map $\alpha:RX\otimes RY\to R(X\otimes Y)$ 
corresponds to the map $\epsilon_X\otimes \epsilon_Y:U(RX\otimes RY)=URX\otimes URY\to X\otimes Y$
by the adjunction  $U\dashv R$ and we have $\alpha_0=\eta_I:I\to RUI=RI$.
Let us denote by $\bV^R$ the category enriched over $\bC$ 
obtained by transfering the enrichment of $\bV$ over itself along the functor $R$. 
An object of  $\bV^R$ is an object of $\bV$, 
the hom object is defined by putting $\bV^R(X,Y)=R\bV(X,Y)$ and
the composition law $\bV^R(Y,Z)\otimes \bV^R(X,Y)\to \bV^R(X,Z)$
is the composite 
$$\xymatrix{
R\bV(Y,Z)\otimes R\bV(X,Y) \ar[r]^\alpha & R(\bV(Y,Z)\otimes \bV(X,Y))\ar[r]^-{R(c)}  & R\bV(X,Z). 
}$$
The category $\bV^R$ is symmetric monoidal closed,
since the category $\bV$ is symmetric monoidal closed.

Dually, there is the notion of {\bf symmetric colax monoidal adjunction}\label{colaxmonoidaladj} $L:\bV\rightleftarrows \bV':U$ where the right adjoint $U$ is monoidal. Then the left adjoint $L$ is automatically colax and symmetric.

\bigskip
If $U\dashv R$ is a lax monoidal adunction, we have two transfer functors $(-)^U:\bV\bCat \to \bV'\bCat$ and $(-)^R:\bV'\bCat \to \bV\bCat$.
\begin{prop}\label{transferadjunction}
There exists an adjunction
$$\xymatrix{
(-)^U:\bV\bCat \ar@<.6ex>[r]& \bV'\bCat:(-)^R. \ar@<.6ex>[l]
}$$
Both functors preserves symmetric monoidal closed categories.
The functor $(-)^R$ send tensored (or cotensored) $\bV'$-categories to tensored (or cotensored) $\bV$-categories
and for a $\bV'$-category $\bC'$, the underlying category of $(\bC')^R$ is isomorphic to the underlying category of $\bC'$.
\end{prop}
\begin{proof}
Let $\bC$ be a $\bV$-category and $\bC'$ be a $\bV'$-category
Let $F$ be a functor $F:\bC^U\to \bC'$, 
the strength of $F$ is a map $\Psi:U\bC(A,B)\to \bC'(FA,FB)$, by adjunction it give a map $\Theta:\bC(A,B)\to R\bC'(FA,FB)$ which defines a functor $G:\bC\to (\bC')^R$. 
We need to prove that the functor condition on $F$ is equivalent to the functor condition on $G$.

The preservation of composition by $\Theta$ is the commutative diagram
$$\xymatrix{
\bC(B,C)\otimes \bC(A,B)\ar[rr]^-{\mathbf{c}}\ar[d]_{\Theta \otimes \Theta }&& \bC(A,C)\ar[d]^{\Theta} \\
(\bC')^R(FB,UC)\otimes (\bC')^R(FA,FB)\ar[rr]^-{c^R} &&(\bC')^R(FA,UC)
}$$
for any triple of objects $A,B,C\in \bC$, where $\mathbf{c}$ and $c^R$ are denoting the composition laws 
in $\bC$ and $(\bC')^R$ respectively. 
Its commutativity is equivalent to that of 
$$\xymatrix{
\bC(B,C)\otimes \bC(A,B)\ar[rr]^-{\mathbf{c}}\ar[d]_{\Theta \otimes \Theta }&& \bC(A,C)\ar[dd]^{\Theta}\\
R\bC'(FB,FC)\otimes R\bC'(FA,FB)\ar[d]_{\alpha} &&\\
R(\bC'(FB,FC)\otimes \bC'(FA,FB)) \ar[rr]^-{R(c)} &&R\bC'(FA,FC),
}$$
where $\alpha$ is the lax structure or $R$ and $c$ is the composition laws in $\bC'$. 
By the adjunction $U\dashv R$, this is also equivalent to showing that the following square commutes 
$$\xymatrix{
U(\bC(B,C)\otimes \bC(A,B))\ar[rr]^-{U(\mathbf{c})}\ar[d]_{U(\Theta \otimes \Theta) }&& U\bC(A,C)\ar[ddd]^{\Psi}\\
U(R\bC'(FB,FC)\otimes R\bC'(FA,FB))\ar@{=}[d] &&\\
UR\bC'(FB,FC)\otimes UR\bC'(FA,FB))\ar[d]^{\epsilon\otimes \epsilon} &&\\
\bC'(FB,FC)\otimes \bC'(FA,FB) \ar[rr]^-{c} && \bC'(FA,FC),
}$$
But the composite of the vertical maps on the left hand side of this diagram is equal to $\Psi\otimes \Psi$.
This proves that $\Theta$ is compatible with composition iff $\Psi$ is.
We leave the verification that $\Theta$ preserves units to the reader iff $\Psi$ does to the reader.
We have proven the existence of a bijection between functors $F:\bC^U\to \bC'$ and $G:\bC\to (\bC')^R$, hence the adjunction.

Let $\bC'$ be a tensored (or cotensored) $\bV'$-category, for $X\in \bV$ and $A,B\in \bC'$ we have bijection between
\begin{center}
\begin{tabular}{lc}
\rule[-2ex]{0pt}{4ex} $UX\otimes A\to B$ (or $A\to [UX,B]$),\\
\rule[-2ex]{0pt}{4ex} $UX\to \bC(A,B)$,\\
\rule[-2ex]{0pt}{4ex} $X\to R\bC(A,B)$
\end{tabular}
\end{center}
which says $(\bC')^R$ is tensored (or cotensored) as a $\bV$-category.
This proves that the functor $(-)^R$ respect the full subcategories of tensored and cotensored categories.

\medskip
The proof that both functor preserves monoidal closed categories was given before the proposition.

\medskip
Let us prove that $(-)^R$ preserves the underlying category.
Let $\bC'$ be a $\bV'$-category and $A$ and $B$ two objects of $\bC'$, Let also $I$ and $J$ be the two unit objects of $\bV$ and $\bV'$,
the set of morphisms from $A$ to $B$ in $\bC'$ is the set of maps $I\to \bC'(A,B)$ in $\bV$, the set of morphisms from $A$ to $B$ in $(\bC')^R$ is the set of maps $J\to R\bC'(A,B)$ in $\bV'$. The bijection between the two is given by $J\to R\bC'(A,B) \iff UJ\to \bC'(A,B)$ and the unit $UJ=I$ of the monoidal structure of $U$.
\end{proof}

\medskip
For example, if $(\bV,\otimes,I)$ is a symmetric monoidal closed category with sums, there is a monoidal adjunction
$I\otimes -:\Set \rightleftarrows \bV : \hom(I,-)$ where $\hom(I,-)$ is the underlying category functor and, for a category $\bC$, $I\otimes \bC$ is the $\bV$-category with hom $I\otimes Hom(A,B) = \coprod_{Hom(A,B)}I$.
In particular, the underlyiing category functor preserves tensored and cotensored categories.

\medskip
Also, if $U:\bV\rightleftarrows \bV':R$ is a monoidal adjunction, we have two functors $\Psi:\bV^U\to \bV'$ and $\Theta:\bV\to (\bV')^R$
which coincides with the functor $U$ on objects.

\subsubsection{Base change for modules}\label{basechangemodule}

Let $U:\bV\to \bV'$ be a monoidal functor and let $\mathbf C$ be a $\mathbf V'$-module with the action given by a monoidal functor $\bV'\to End(\bV)$. Then we obtain a structure of $\mathbf V$-module on $\mathbf C$ by the composition
$\bV\to \bV'\to End(\bV)$.
If we write the $\bV'$-action on $C$ as $\diamond': \bV'\times  \bC\to  \bC$, the $\bV$-action is given by $X\diamond A := U(X)\diamond' A$.
This produces a {\bf restriction functor}
$$
Restr_U:\MMod(\bV')\to \MMod(\bV).
$$

Let $\bC'$ is a right divisible module or a right codivisible module, then, using the equivalence with (co)tensored categories $Restr_U\bC'$ coincides with $(\bC')^R$.
This is a consequence of the bijection between $UX\otimes A\to B$ and $X\to R\bC'(A,B)$ compute above.

\subsection{Common (co)equalizers}\label{commonequalizer}

The following notions are useful in some of our proofs.
For four maps $f_1,g_1:X\to Y_j$, $f_2,g_2:X\to Y_2$, we define their {\em common equalizer} as the limit of the diagram
$$\xymatrix{
X\ar@<.6ex>[r]^-{f_1}\ar@<-.6ex>[r]_-{g_1} \ar@<.6ex>[d]^-{f_2}\ar@<-.6ex>[d]_-{g_2}&Y_1\\
Y_2
}.$$
The common equalizer of the previous diagram is thus an object $Z$ and a map $h:Z\to X$ such that $f_ih=g_ih$ for $i=1,2$:
$$\xymatrix{
Z\ar@{-->}[rd]^-h\ar@{-->}@/^1pc/[rrd]^-{\quad f_1h=g_1h}\ar@{-->}@/_1pc/[rdd]_-{f_2h=g_2h}\\
&X\ar@<.6ex>[r]^-{f_1}\ar@<-.6ex>[r]_-{g_1} \ar@<.6ex>[d]^-{f_2}\ar@<-.6ex>[d]_-{g_2}&Y_1\\
&Y_2
}.$$
For $i=1,2$, if $Z_i$ is the equalizer of 
$\xymatrix{
X\ar@<.6ex>[r]^-{f_i}\ar@<-.6ex>[r]_-{g_i}&Y_i,
}$
then $Z=Z_1\cap Z_2$, where the intersection is taken in subobjects of $X$, is a common equalizer of $(f_1,g_1,f_2,g_2)$.

\bigskip
Dually, for four maps $f_1,g_1:Y_j\to X$, $f_2,g_2:Y_2\to X$, we define their {\em common coequalizer} as the colimit of the diagram
$$\xymatrix{
&Y_2\ar@<.6ex>[d]^-{f_2}\ar@<-.6ex>[d]_-{g_2}\\
Y_1\ar@<.6ex>[r]^-{f_1}\ar@<-.6ex>[r]_-{g_1} & X
}.$$
The common coequalizer of the previous diagram is thus an object $Z$ and a map $h:Z\to X$ such that $f_ih=g_ih$ for $i=1,2$:
$$
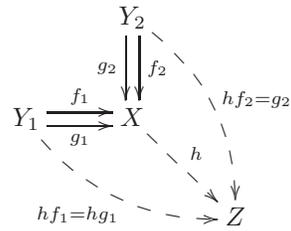
\xymatrix{
&Y_2\ar@<.6ex>[d]^-{f_2}\ar@<-.6ex>[d]_-{g_2}\ar@{-->}@/^1pc/[rdd]^-{ hf_2=g_2}\\
Y_1\ar@{-->}@/_1pc/[rrd]_-{hf_1=hg_1\quad}\ar@<.6ex>[r]^-{f_1}\ar@<-.6ex>[r]_-{g_1} & X\ar@{-->}[rd]^-h\\
&&Z
}.$$
For $i=1,2$, if $Z_i$ is the equalizer of 
$\xymatrix{
Y_i\ar@<.6ex>[r]^-{f_i}\ar@<-.6ex>[r]_-{g_i}&X,
}$
then the amalgamated sum $Z=Z_1\cup_X Z_2$ is a common coequalizer of $(f_1,g_1,f_2,g_2)$.

\end{document}